\documentclass[11pt,a4paper,reqno]{article}
\usepackage[OT1]{fontenc}

\usepackage[margin=25.4mm]{geometry}
\usepackage{amsmath,amssymb,amsthm,mathtools,microtype}
\usepackage[hidelinks]{hyperref}
\allowdisplaybreaks[1]
\AtBeginDocument{%
 \setlength{\abovedisplayskip}{6pt plus 2pt minus 2pt}%
 \setlength{\belowdisplayskip}{6pt plus 2pt minus 2pt}%
 \setlength{\abovedisplayshortskip}{0pt plus 2pt}%
 \setlength{\belowdisplayshortskip}{4pt plus 2pt minus 2pt}%
}
\newtheorem{lemma}{Lemma}
\newtheorem{proposition}[lemma]{Proposition}
\newtheorem{theorem}[lemma]{Theorem}
\newtheorem{corollary}[lemma]{Corollary}
\theoremstyle{remark}\theoremstyle{plain}
\newcommand{\PP}{\mathbb P}
\newcommand{\EE}{\mathbb E}
\newcommand{\ZZ}{\mathbb Z}
\newcommand{\ind}{\mathbf1}
\newcommand{\one}{\mathbf1}
\newcommand{\cL}{\mathcal L}
\newcommand{\cK}{\mathcal K}
\newcommand{\Qw}{\mathsf Q_N^I}
\newcommand{\Qaug}{\overline{\mathsf Q}_N^I}

\newcommand{\Ttwo}{\mathbb T^2}

\DeclareMathOperator{\Cov}{Cov}
\DeclareMathOperator{\Var}{Var}

\DeclareMathOperator{\supp}{supp}
\DeclareMathOperator{\diam}{diam}
\DeclareMathOperator{\dist}{dist}
\DeclareMathOperator{\Pack}{Pack}
\DeclareMathOperator{\Leb}{Leb}
\hypersetup{bookmarksnumbered=true,bookmarksdepth=4}
\let\mergedcontentsline\addcontentsline
\newcommand{\MergedChapter}{0}
\newif\ifMergedMain
\MergedMaintrue
\newcommand{\AppendixEquationPrefix}{}
\newcommand{\appendixchapter}[2]{%
 \section*{Appendix #1. #2}%
 \hypertarget{proof.appendix.#1}{}%
 \mergedcontentsline{toc}{subsection}{Appendix #1. #2}%
 \setcounter{lemma}{0}\setcounter{equation}{0}\setcounter{subsection}{0}%
 \renewcommand{\thelemma}{#1.\arabic{lemma}}%
 \renewcommand{\theequation}{\AppendixEquationPrefix#1.\arabic{equation}}%
 \renewcommand{\thesubsection}{#1.\arabic{subsection}}%
 \renewcommand{\theHequation}{\MergedChapter.appendix.#1.\arabic{equation}}%
 \renewcommand{\theHlemma}{\MergedChapter.appendix.#1.\arabic{lemma}}%
 \renewcommand{\theHsubsection}{\MergedChapter.appendix.#1.\arabic{subsection}}%
}
\newcommand{\mergedchaptersetup}[2]{%
 \renewcommand{\MergedChapter}{#1}%
 \renewcommand{\AppendixEquationPrefix}{}%
 \setcounter{lemma}{0}\setcounter{equation}{0}%
 \setcounter{section}{0}\setcounter{subsection}{0}\setcounter{subsubsection}{0}%
 \setcounter{footnote}{0}%
 \renewcommand{\thesection}{\arabic{section}}%
 \renewcommand{\thesubsection}{\thesection.\arabic{subsection}}%
 \renewcommand{\thesubsubsection}{\thesubsection.\arabic{subsubsection}}%
 \renewcommand{\thelemma}{#1.\arabic{lemma}}%
 \renewcommand{\theequation}{#1.\arabic{equation}}%
 \renewcommand{\theHlemma}{#1.main.\arabic{lemma}}%
 \renewcommand{\theHequation}{#1.main.\arabic{equation}}%
 \renewcommand{\theHsection}{#1.main.\arabic{section}}%
 \renewcommand{\theHsubsection}{#1.main.\arabic{section}.\arabic{subsection}}%
 \renewcommand{\theHsubsubsection}{#1.main.\arabic{section}.\arabic{subsection}.\arabic{subsubsection}}%
 \ifMergedMain
  \phantomsection
  \mergedcontentsline{toc}{section}{Chapter #1. #2}%
 \fi
 \renewcommand{\addcontentsline}[3]{}%
}
\let\MergedOriginalSubsection\subsection
\let\MergedOriginalSubsubsection\subsubsection
\def\MergedMainEntry#1. #2\MergedEntryEnd{%
 \mergedcontentsline{toc}{subsection}{\protect\numberline{#1}#2}%
}
\RenewDocumentCommand{\subsection}{s o m}{%
 \IfBooleanTF{#1}{%
  \MergedOriginalSubsection*{#3}%
  \ifMergedMain\MergedMainEntry#3\MergedEntryEnd\fi
 }{%
  \IfNoValueTF{#2}{\MergedOriginalSubsection{#3}}{\MergedOriginalSubsection[#2]{#3}}%
  \ifMergedMain\else
   \pdfbookmark[3]{\thesubsection\space #3}{nav.\MergedChapter.\theHsubsection}%
  \fi
 }%
}
\RenewDocumentCommand{\subsubsection}{s o m}{%
 \IfBooleanTF{#1}{\MergedOriginalSubsubsection*{#3}}{%
  \IfNoValueTF{#2}{\MergedOriginalSubsubsection{#3}}{\MergedOriginalSubsubsection[#2]{#3}}%
  \ifMergedMain\else
   \pdfbookmark[4]{\thesubsubsection\space #3}{nav.\MergedChapter.\theHsubsection.\arabic{subsubsection}}%
  \fi
 }%
}

\makeatletter
\newcommand{\CompactAppendixHeadings}{%
 \renewcommand{\section}{\@startsection{section}{1}{\z@}%
  {-2.5ex \@plus -1ex \@minus -.2ex}%
  {1ex \@plus .2ex}{\normalfont\large\bfseries}}%
 \renewcommand{\MergedOriginalSubsection}{\@startsection{subsection}{2}{\z@}%
  {-2ex \@plus -.5ex \@minus -.2ex}%
  {.7ex \@plus .2ex}{\normalfont\normalsize\bfseries}}%
 \renewcommand{\MergedOriginalSubsubsection}{\@startsection{subsubsection}{3}{\z@}%
  {-1.5ex \@plus -.5ex \@minus -.2ex}%
  {.5ex \@plus .2ex}{\normalfont\normalsize\itshape}}%
}
\makeatother

\hypersetup{
 pdfauthor={Yechi Zhou},
 pdftitle={A limit law for the cover time of the two-dimensional discrete torus},
 pdfsubject={Cover-time fluctuations on the two-dimensional discrete torus}
}
\begin{document}
\pagenumbering{roman}
\begin{center}
{\LARGE A limit law for the cover time of the\par
two-dimensional discrete torus\par}
\vspace{9pt}{\normalsize Yechi Zhou\par}
\vspace{4pt}{\small School of Mathematical Sciences, Fudan University, Shanghai, China\par}
\vspace{3pt}{\small\href{mailto:25210180114@m.fudan.edu.cn}{\texttt{25210180114@m.fudan.edu.cn}}\par}
\end{center}
\begin{abstract}
We determine the limiting distribution of the cover time of simple random walk on the two-dimensional discrete torus. For the continuous-time walk with total jump rate one on $(\mathbb Z/N\mathbb Z)^2$, let $T_N$ denote its cover time. We prove that
\[
\frac{T_N}{(2/\pi)N^2\log N}-2\log N+\log\log N
\ \Longrightarrow\ G+\log(\kappa Z),
\]
where $G$ is a standard Gumbel random variable, $Z$ is the total mass of the critical Gaussian multiplicative chaos associated with the zero-average Gaussian free field on the unit torus, $G$ and $Z$ are independent, and $\kappa>0$ is deterministic. This answers the limit-law question suggested by Aldous and recorded by Dembo, Peres, Rosen and Zeitouni. The proof identifies the random fluctuations in the number of small, well-separated unvisited components at a deterministic time before coverage, and then estimates the time needed to visit the remaining components.
\end{abstract}

\begingroup
\par
\endgroup

\clearpage
\begingroup
\small
\tableofcontents
\par\endgroup
\clearpage
\pagenumbering{arabic}
\begingroup
\renewcommand{\theequation}{Intro.\arabic{equation}}
\renewcommand{\theHequation}{intro.\arabic{equation}}
\setcounter{equation}{0}
\newcommand{\FrontSection}[1]{%
  \MergedOriginalSubsection*{#1}%
  \mergedcontentsline{toc}{subsection}{#1}%
}
\section*{Introduction}
\mergedcontentsline{toc}{section}{Introduction}

How long does a random walk take to visit every vertex of a large graph? On the two-dimensional torus, this question links a global exploration problem to the extreme values of a strongly correlated random field. Most vertices are visited long before coverage is complete. The final stage is governed by a sparse collection of regions that have received unusually little occupation across many spatial scales. Their number and arrangement retain information about the walk on the scale of the whole torus.

The leading and second-order terms of the cover time are known. We identify the remaining fluctuations on the scale $N^2\log N$: the limit is a Gumbel variable shifted by the logarithm of the total mass of a periodic critical Gaussian measure. This answers the limit-law question recorded in \cite[Section 9, item 4]{DPRZ04} and resolves the tightness and nondegeneracy conjecture of Bramson and Zeitouni \cite[Conjecture 5.1]{BZ09}.

\FrontSection{Model and main result}

Let $V_N=(\mathbb Z/N\mathbb Z)^2$, and let $X=(X_t)_{t\ge0}$ be continuous-time simple random walk on $V_N$ with generator
\[
\mathcal L_N f(x)=\frac14\sum_{e\in\{\pm e_1,\pm e_2\}}
\bigl(f(x+e)-f(x)\bigr).
\]
Thus the holding times have mean one. Write $\mathbb P_x$ for the law of the walk started at $x$, and define
\[
T_N=\inf\{t\ge0:\{X_s:0\le s\le t\}=V_N\}.
\]
Throughout the paper, logarithms are natural. Put
\begin{equation}\label{intro:eq:normalization}
g=\frac2\pi,\qquad b_N=gN^2\log N,\qquad
a_N=b_N(2\log N-\log\log N).
\end{equation}

The random shift is described by the critical Gaussian multiplicative chaos of the zero-average Gaussian free field on $\mathbb T^2=\mathbb R^2/\mathbb Z^2$. We denote this measure by $\mathcal Z$ and its total mass by $Z=\mathcal Z(\mathbb T^2)$. It satisfies $0<Z<\infty$ almost surely; its normalization is specified just below the theorem.

\medskip
\noindent\textbf{Main theorem (Theorem~\ref{chVIII:thm:main}).}
\emph{There is a deterministic constant $\kappa\in(0,\infty)$ such that}
\begin{equation}\label{intro:eq:main}
\lim_{N\to\infty}\sup_{x\in V_N}\sup_{z\in\mathbb R}
\left|\mathbb P_x\!\left(\frac{T_N-a_N}{b_N}\le z\right)
-\mathbb E\exp\{-\kappa Ze^{-z}\}\right|=0.
\end{equation}
\emph{Equivalently, $(T_N-a_N)/b_N$ converges in distribution to $G+\log(\kappa Z)$, where $\mathbb P(G\le z)=\exp\{-e^{-z}\}$ and $G$ is independent of $Z$.}
\medskip

To fix the normalization, let $h$ be the zero-average continuum GFF on $\mathbb T^2$, with covariance $4(-\Delta)^{-1}$ on functions of integral zero, where $\Delta=\partial_1^2+\partial_2^2$. Lebesgue measure on $\mathbb T^2$ has total mass one. Set
\[
h_\varepsilon=e^{\varepsilon^2\Delta}h,\qquad
V_\varepsilon=\mathbb E h_\varepsilon(x)^2
=g\log(1/\varepsilon)+O(1),\qquad \alpha=2/\sqrt g.
\]
The variance is independent of $x$. The critical measure is the weak limit in probability
\begin{equation}\label{intro:eq:chaos}
\mathcal Z(dx)=\lim_{\varepsilon\downarrow0}
\sqrt{\log(1/\varepsilon)}
\exp\!\left\{\alpha h_\varepsilon(x)-\frac{\alpha^2}{2}V_\varepsilon\right\}dx.
\end{equation}
Chapter II constructs this periodic limit using critical-chaos results of Duplantier, Rhodes, Sheffield and Vargas \cite{DRSV14a,chII:DRSV}, and relates it to the discrete field. The measure is almost surely finite, nonatomic, and positive on every nonempty open set. Under the convention \eqref{intro:eq:chaos}, the proof characterizes $\kappa$ as a product $\beta\Gamma$ of two positive comparison coefficients. Rescaling $\mathcal Z$ by a deterministic constant rescales $\kappa$ inversely.

The theorem specifies fluctuations of order $N^2\log N$ around
\[
a_N=\frac4\pi N^2(\log N)^2
-\frac2\pi N^2\log N\log\log N.
\]
Its square-root form is
\begin{equation}\label{intro:eq:squareroot}
\sqrt{T_N/N^2}-\frac2{\sqrt\pi}\log N
+\frac1{2\sqrt\pi}\log\log N
\ \Longrightarrow\ \frac{G+\log(\kappa Z)}{2\sqrt\pi}.
\end{equation}
The same limit holds for the discrete-time cover time: the rate-one continuous-time walk has the same embedded jump chain, and the independent exponential holding times contribute a clock error of order $O_{\mathbb P}(N\log N)=o_{\mathbb P}(b_N)$. Since the limiting distribution in \eqref{intro:eq:squareroot} is continuous, strictly increasing, and nondegenerate, centering by the median gives the formulation in \cite[Conjecture 5.1]{BZ09}. Translation invariance makes the cover-time law independent of the initial vertex.

\FrontSection{From first-order asymptotics to a limit law}

Dembo, Peres, Rosen and Zeitouni \cite{DPRZ04} proved that
\[
\frac{T_N}{N^2(\log N)^2}\longrightarrow\frac4\pi
\quad\text{in probability},
\]
settling the first-order conjecture of Aldous. Their work also established leading-order asymptotics for Brownian covering on two-dimensional manifolds. In Section 9, item 4, they recorded a further question suggested by Aldous: to find a limit law for an appropriately normalized cover time. Motivated by branching-random-walk analogies, they proposed that the square root of the cover time, divided by $N$ and centered by its median, should have a nondegenerate limit. Bramson and Zeitouni \cite[Conjecture 5.1]{BZ09} subsequently formulated the corresponding tightness and nondegeneracy conjecture.

Ding \cite{Ding12} controlled the centered square root of the cover time to order $\log\log N$. Belius and Kistler \cite{BK17} determined the second-order correction for small-radius Brownian covering on the continuous two-dimensional torus, connecting it to multiscale barrier constraints and entropic repulsion. Abe \cite{Abe21} obtained the corresponding correction for discrete simple random walk. More precisely, for some $c\in(0,1)$, his result gives
\[
\mathbb P\!\left(\left|\frac{T_N-a_N}{b_N}\right|
\le(\log\log N)^c\right)\longrightarrow1.
\]
The question left by this estimate is the behavior on a bounded fluctuation scale. It is raised explicitly in \cite[Remark 1.3]{Abe21}. The main theorem establishes convergence on this scale and identifies the limit.

The geometry underlying the last stages of coverage was developed in \cite{DPRZ06}, which describes the abundance, clustering and fractal structure of late points on the discrete torus. Comets, Popov and Vachkovskaia \cite{CPV16} studied microscopic neighborhoods of late points through two-dimensional random interlacements, while Comets, Gallesco, Popov and Vachkovskaia \cite{CGPV13} studied cover-time large deviations. Together, these works describe the spatial organization of rare unvisited points and their local environment.

\FrontSection{Related models and the periodic setting}

For Brownian motion, Belius, Rosen and Zeitouni \cite{BRZ20,BRZ20c} proved tightness of the cover time on the two-dimensional sphere. A randomly shifted Gumbel limit for small-radius Brownian covering of the two-dimensional torus was announced by Dembo, in joint work with Rosen and Zeitouni, in a 2022 seminar \cite{Dembo22}. That announcement concerns the continuum torus, whereas the present paper concerns coverage of all vertices by a discrete random walk.

The closest discrete results are those of Louidor and Saglietti. They established tightness for lattice approximations of planar domains with wired boundary conditions \cite{chIV:LS24}, and subsequently proved a joint limit for the cover time and the last visited vertex \cite{chVI:LS26}. The cover-time marginal is a randomly shifted Gumbel law, described through critical Liouville measure and the associated Gaussian field. Section 1.2.3 of \cite{chVI:LS26} identifies periodic boundary conditions as an outstanding case, with even tightness still open.

On the torus, random-walk behavior in disks much smaller than the system size is planar, whereas occupation profiles across the whole torus depend on the periodic Green kernel. In wired domains, inverse local time at the distinguished boundary vertex provides a useful parametrization for the Ray--Knight comparison. The cover-time distribution requires occupation estimates at deterministic physical times. The proof must relate the Gaussian comparison to these observation times on the $N^2\log N$ fluctuation scale.

Our argument builds on the low-local-time geometry and annular analysis of \cite{chIV:LS24,chVI:LS26}. The key step for the periodic model is to determine the random number of small, well-separated unvisited components shortly before coverage. Local estimates describe low occupation within individual regions; to obtain the cover-time law, we also need to identify the fluctuations in the total number of these components.

For the stationary walk, we prove that, at a suitable deterministic observation time, the number of geometrically selected unvisited components, divided by $\sqrt{\log N}$, converges in law to a deterministic multiple of $Z$. We first study a nonnegative, normalized local-time statistic whose asymptotics can be calculated using the Gaussian field. The Gaussian comparison initially determines a time average of this statistic. We recover its limiting law at the deterministic observation time and prove that the normalized component count differs from a deterministic multiple of the statistic by a term converging to zero in probability.

We then estimate the probabilities that representatives of the remaining components stay unvisited over the final time interval. We control the additional time needed to cover every vertex in each component and the probability that at least one excluded vertex remains unvisited. These estimates yield the cover-time law. The periodic critical mass determines the random coefficient in the component count, while the final visits give the Gumbel distribution.

\FrontSection{Gaussian fields and the random shift}

The connection between cover times and Gaussian free fields provides a route from occupation statistics to extreme-value theory. Ding, Lee and Peres \cite{DLP12} proved that, on a finite connected graph, the expected cover time is comparable, up to universal constants, to the number of edges times the square of the expected maximum of the associated pinned GFF. At the level of local times, the generalized second Ray--Knight theorem of Eisenbaum, Kaspi, Marcus, Rosen and Shi \cite{EKMRS00} relates occupation fields to squares of Gaussian fields.

For the two-dimensional Dirichlet GFF, Bramson, Ding and Zeitouni \cite{BDZ16} proved convergence in law of the centered maximum. Biskup and Louidor identified the limiting process of extreme local maxima \cite{BL16}, its full cluster structure \cite{chII:BL}, and the critical Liouville measure governing its random intensity \cite{BL20}. Conditional on this measure, the cluster maxima have Poisson statistics with an exponential height intensity, giving a randomly shifted Gumbel law for the maximum. The critical-chaos construction \cite{DRSV14a,chII:DRSV} provides a measure-valued description of the randomness retained at macroscopic scales. Background accounts are \cite{Biskup20,RV14}.

A related picture emerges for exceptional local times \cite{Abe15}. In wired planar domains, Abe and Biskup \cite{AB22} identified scaling limits for thick, thin and avoided points at times comparable with the cover time, using the boundary-local-time parametrization. Abe, Biskup and Lee \cite{ABL23} treated discrete-time walks and described the zero-average fields and associated factors in the limiting measures. These results connect low occupation to Liouville measures before complete coverage. The present argument identifies the periodic critical mass governing the components that enter the final covering window.

The contrast with higher dimensions is instructive \cite{BHT22}. For discrete tori of dimension $d\ge3$, Belius \cite{Belius13} proved convergence to an ordinary Gumbel law. In two dimensions, long-range correlations survive in the limiting random shift.

\FrontSection{Structure of the proof}

Write $L=\log N$ and $\ell_t(x)=\int_0^t\mathbf1_{\{X_s=x\}}\,ds$. The proof observes the walk at a deterministic time $s_N$ satisfying
\[
s_N=b_N\left(2L-\frac32\log L\right)
+O\!\left(N^2(\log L)^2\right).
\]
The remaining time to $a_N+zb_N$ is therefore
\begin{equation}\label{intro:eq:window}
a_N+zb_N-s_N=b_N\left(\frac12\log L+z+o(1)\right).
\end{equation}
At $s_N$, the relevant low-local-time vertices form small components separated by much larger distances. Their number is of order $\sqrt L$. The argument first identifies the random coefficient in this count, then converts it into the cover-time law.

\paragraph*{Low-local-time geometry and Gaussian activity.}
We remove an exceptional set of low-local-time vertices with neighbors at intermediate distances. The retained set has a canonical decomposition into small, well-separated components. Each component is assigned a Gaussian contact activity: a nonnegative weight measuring the scalar Gaussian shifts that create a small square-field value in that component. Summing these weights and dividing by $\sqrt L$ gives an observable $A_N$. Multiscale estimates control the exceptional set, establish tightness of $A_N$, and bound its variation over short time intervals. The annular analysis uses planar random-walk potential theory \cite{chVI:LL} and the estimates of \cite{chIV:LS24,chVI:LS26}.

Activity is the intermediate quantity linking Gaussian asymptotics to exact zeros of local time. The finite-graph coarea identity compares it with Gaussian configurations constrained by total energy. Their asymptotics are determined by the periodic critical measure, using Dirichlet extremal-cluster results \cite{chII:BL} and ballot estimates \cite{chII:GL}. The comparison initially yields a time average: the squared Gaussian root coordinate contributes a random shift of the observation time. Inverting the resulting one-sided Gamma average, and using short-window comparisons to return to deterministic time, gives
\begin{equation}\label{intro:eq:activity}
A_N\Longrightarrow\Gamma Z\qquad\text{for a deterministic }\Gamma>0.
\end{equation}

\paragraph*{From activity to unvisited components.}
The presence of an unvisited vertex is an indicator of the local occupation profile, whereas activity is a Gaussian-weighted function of the same profile. We compare these two observables under common conditional local laws, obtained by resampling interior excursions while retaining their exterior data. Aggregating the local comparisons identifies a deterministic proportionality coefficient. If $K_N$ is the number of canonical unvisited components at $s_N$ and $C_N=K_N/\sqrt L$, then
\begin{equation}\label{intro:eq:count}
C_N-\beta A_N\longrightarrow0\quad\text{in probability},
\qquad C_N\Longrightarrow\kappa Z,\qquad\kappa=\beta\Gamma>0.
\end{equation}
Here the intermediate limits are under the stationary path law. A packing and local-clearing argument establishes positivity of $\beta$. Thus the critical mass is identified before the final covering step.

\paragraph*{The final covering window.}
For a fixed collection of $k$ well-separated representatives, spectral and capacity estimates give a joint survival probability asymptotic to $L^{-k/2}e^{-kz}$ over the interval in \eqref{intro:eq:window}. Factorial-moment estimates and Bonferroni inequalities turn these probabilities into the probability that no representative survives. Local clearing estimates control the additional time needed to cover each component, and a separate estimate handles the exceptional set. The result is
\begin{equation}\label{intro:eq:cover}
\mathbb P(T_N\le a_N+zb_N)
=\mathbb E\exp\{-C_Ne^{-z}\}+o(1)
\longrightarrow\mathbb E\exp\{-\kappa Ze^{-z}\}.
\end{equation}
The random intensity comes from the walk before $s_N$; the subsequent disappearance of the remaining components produces the Gumbel factor.

\FrontSection{Organization}

Chapter I develops periodic potential theory and the spectral estimates for hitting times. Chapter II constructs the periodic Gaussian field and its critical measure, and establishes the Gaussian asymptotics used later. Chapter III introduces local-time activity, the finite-graph identity and the resampling laws. Chapter IV proves the geometric and time-comparison estimates, and Chapter V identifies the activity limit. Chapter VI analyzes conditional local profiles; Chapter VII deduces the limiting number of unvisited components. Chapter VIII completes the covering argument and proves Theorem~\ref{chVIII:thm:main}. Appendices A--AB contain the detailed supporting estimates.
A Lean~4 formalization of Theorem~\ref{chVIII:thm:main} is available in the
repository linked in the \hyperref[sec:availability-disclosure]{code availability statement}.
\par\endgroup

\clearpage
\MergedMaintrue
\begingroup
\mergedchaptersetup{I}{Random walk and potential theory}
\begin{center}
{\large\scshape Chapter I\par}\vspace{6pt}
{\LARGE Random walk and potential theory\par}
\end{center}
\vspace{6pt}
\subsection*{1.1. The walk and its Fourier representation}

Let $N\ge3$ be an integer and put $V_N=(\ZZ/N\ZZ)^2$,
$n=|V_N|=N^2$, and $L=\log N$. All limits below are as
$N\to\infty$, unless another limiting parameter is specified.
We consider continuous-time
simple random walk $X=(X_t)_{t\ge0}$ on $V_N$, with total jump rate one
and generator
\[
 \cL_Nf(x)=\frac14\sum_{e\in\{\pm e_1,\pm e_2\}}
                  \bigl(f(x+e)-f(x)\bigr).
\]
Here $e_1,e_2$ are the coordinate unit vectors; a unit lattice
direction means an element of $\{\pm e_1,\pm e_2\}$. We use $|\cdot|$
for the Euclidean norm of a spatial vector and $|\cdot|_\infty$
for its maximum norm. For a finite set, $|\cdot|$ denotes cardinality.
The law and expectation of the walk started at $x$ are $\PP_x$ and $\EE_x$;
its natural filtration is $\mathcal F_t=\sigma(X_s:0\le s\le t)$.
For an initial probability measure $\mu$ on $V_N$, write
$\PP_\mu=\sum_x\mu(x)\PP_x$, with expectation $\EE_\mu$.
The uniform measure $\pi_N(x)=N^{-2}$ is reversible, so
$\PP_{\pi_N}$ is the stationary law.
On $L^2(\pi_N)$ we use the inner product
\[
 \langle f,h\rangle_{\pi_N}=N^{-2}\sum_{x\in V_N}f(x)\overline{h(x)}.
\]
For vectors indexed by a finite set $D$, the unqualified norm is
$\|f\|_2^2=\sum_{x\in D}|f(x)|^2$, with counting measure.
In particular $\|f\|_{L^2(\pi_N)}^2=n^{-1}\|f\|_2^2$ on $V_N$.
We write $\ind_A$ for the indicator of a set or event $A$, and
$\ind$ for the constant function one. For real $a,b$, put
$a_+=\max\{a,0\}$, $a\vee b=\max\{a,b\}$ and
$a\wedge b=\min\{a,b\}$.
Write
\[
 p_t(x,y)=\PP_x(X_t=y),\qquad
 \ell_t(y)=\int_0^t\ind_{\{X_u=y\}}\,du.
\]
For $B\subseteq V_N$, put
\[
 H_y=\inf\{t\ge0:X_t=y\},\qquad
 H_B=\min_{y\in B}H_y,\qquad T_B=\max_{y\in B}H_y,
\]
with $H_\varnothing=\infty$ and $T_\varnothing=0$.
Thus $H_B$ is the entrance time into $B$, while $T_B$ is its cover
time; the cover time of the torus is $T_N=T_{V_N}$.
The normalization used throughout is
\[
 g=\frac2\pi,\qquad b_N=gN^2\log N,\qquad
 a_N=b_N(2\log N-\log\log N),\qquad t_N(z)=a_N+zb_N.
\]

Choose the frequency representatives
\[
 \cK_N=\bigl(\ZZ\cap[-N/2,N/2)\bigr)^2,\qquad
 \chi_k(x)=e^{2\pi i k\cdot x/N},
\]
and set
\[
 \lambda_N(k)=\sin^2\frac{\pi k_1}{N}
              +\sin^2\frac{\pi k_2}{N},\qquad
 \gamma_N=\sin^2\frac\pi N.
\]

\begin{lemma}\label{chI:lem:fourier}
The eigenvalues of $-\cL_N$ are $\lambda_N(k)$, $k\in\cK_N$, and
\begin{equation}\label{chI:eq:spectrum}
 \frac{4|k|^2}{N^2}\le\lambda_N(k)\le\frac{\pi^2|k|^2}{N^2},
 \qquad \min_{k\ne0}\lambda_N(k)=\gamma_N
 =\frac{\pi^2}{N^2}\bigl(1+O(N^{-2})\bigr).
\end{equation}
For $t\ge0$ and $x,y\in V_N$,
\begin{equation}\label{chI:eq:heat}
 N^2p_t(x,y)=1+\sum_{k\in\cK_N\setminus\{0\}}
                    e^{-t\lambda_N(k)}\chi_k(x-y).
\end{equation}
The integrals
\[
 G_N(x,y)=\int_0^\infty\bigl(p_u(x,y)-N^{-2}\bigr)\,du,
 \qquad R_t(x,y)=\int_t^\infty\bigl(p_u(x,y)-N^{-2}\bigr)\,du
\]
converge absolutely, and
\begin{equation}\label{chI:eq:green}
 \begin{aligned}
 G_N(x,y)&=\frac1{N^2}\sum_{k\ne0}
                  \frac{\chi_k(x-y)}{\lambda_N(k)},\\
 R_t(x,y)&=\frac1{N^2}\sum_{k\ne0}
                  \frac{e^{-t\lambda_N(k)}\chi_k(x-y)}{\lambda_N(k)}.
 \end{aligned}
\end{equation}
Here and below, these sums are over $\cK_N\setminus\{0\}$. Moreover,
\begin{equation}\label{chI:eq:occupation}
 \EE_x\ell_t(y)=\frac{t}{N^2}+G_N(x,y)-R_t(x,y).
\end{equation}
\end{lemma}

\begin{proof}
The character basis of $L^2(\pi_N)$ satisfies
\[
 \cL_N\chi_k=-\lambda_N(k)\chi_k,\qquad
 P_t\chi_k=e^{-t\lambda_N(k)}\chi_k.
\]
The bounds $2|j|/N\le|\sin(\pi j/N)|\le\pi|j|/N$ for $|j|\le N/2$
and the expansion of $\sin(\pi/N)$ give \eqref{chI:eq:spectrum}.
Spectral expansion gives \eqref{chI:eq:heat}; integrating its finite
sum gives \eqref{chI:eq:green}, with absolute convergence since
$\lambda_N(k)>0$ for $k\ne0$. Finally,
$\EE_x\ell_t(y)=\int_0^tp_u(x,y)\,du$ gives \eqref{chI:eq:occupation}.
\end{proof}

We regard $G_N$ as an operator acting by counting-measure summation:
\[
 (G_Nf)(x)=\sum_{y\in V_N}G_N(x,y)f(y).
\]
In particular, $G_N$ is symmetric and has zero row sums. On the
nonconstant Fourier modes its eigenvalues are $\lambda_N(k)^{-1}$;
hence it is positive semidefinite, with kernel consisting of the
constant functions, and
\[
 (-\cL_N)G_N=I-\Pi_N,\qquad
 \Pi_N f=\left(N^{-2}\sum_xf(x)\right)\ind.
\]
These identities fix the normalization of the Green operator used
in the capacity estimates.

\subsection*{1.2. Mixing and Green-function bounds}

For $x,y\in V_N$, let
\[
 d(x,y)=\min_{v\equiv x-y\ ({\rm mod}\ N)}|v|,\qquad
 d_1(x,y)=\min_{v\equiv x-y\ ({\rm mod}\ N)}(|v_1|+|v_2|).
\]
Thus $d_1$ is nearest-neighbor graph distance and
$d\le d_1\le\sqrt2d$. For a nonempty set $B\subseteq V_N$,
write $d(x,B)=\min_{y\in B}d(x,y)$ and
$\operatorname{diam}_d B=\max_{x,y\in B}d(x,y)$; the analogous
notation with subscript $d_1$ uses graph distance.

In estimates, $C,c>0$ denote constants whose values may change
between occurrences. Constants in this chapter are absolute;
later, subscripts or the accompanying statement specify their
dependence on fixed parameters. For a positive deterministic scale,
the notation $O(\cdot)$ and $o(\cdot)$ refers respectively to a bounded
ratio and a ratio tending to zero, uniformly in any parameters over
which the statement is uniform.
The mixing time used below is much shorter than $b_N$.

\begin{lemma}\label{chI:lem:mixing}
Put $h_N=4N^2\log L$. For all sufficiently large $N$, uniformly in
$t\ge h_N$,
\begin{equation}\label{chI:eq:mixing}
 \sup_{x,y}|N^2p_t(x,y)-1|\le L^{-8},\qquad
 \sup_{x,y}|R_t(x,y)|\le L^{-8}.
\end{equation}
Moreover, $h_N/b_N\to0$.
\end{lemma}

\begin{proof}
Equations \eqref{chI:eq:spectrum}--\eqref{chI:eq:green} give, for $t\ge N^2/4$,
\[
 \sup_{x,y}|N^2p_t(x,y)-1|+\sup_{x,y}|R_t(x,y)|
       \le Ce^{-4t/N^2}.
\]
For $t\ge h_N$, the right-hand side is at most $CL^{-16}\le L^{-8}$
for large $N$, and $h_N/b_N=(4/g)(\log L)/L\to0$.
\end{proof}

The spatial estimates require cancellation in the Fourier sum.
We isolate the required summation facts; their proofs are in
\hyperlink{proof.appendix.A}{Appendix~A}. For $z\in V_N$, set
\[
 \Sigma_N(z)=\sum_{k\in\cK_N\setminus\{0\}}
                      \frac{\chi_k(z)}{|k|^2}.
\]

\begin{lemma}\label{chI:lem:fourier-sums}
Uniformly in $N\ge3$,
\begin{equation}\label{chI:eq:quadratic-sums}
 \begin{aligned}
 \Sigma_N(0)&=2\pi\log N+O(1),\\
 |\Sigma_N(z)|&\le C\left(1+\log\frac N{d(z,0)}\right)
                       &&(z\ne0),\\
 \sum_{k\ne0}\min\{u^{-2},|k|^{-2}\}
   &\le C\bigl(1+\log\max\{N/u,1\}\bigr)
                       &&(u\ge1).
 \end{aligned}
\end{equation}
\end{lemma}

\begin{proposition}\label{chI:prop:green-bounds}
For all $N\ge3$ and $x,y\in V_N$,
\begin{equation}\label{chI:eq:green-bounds}
 \begin{aligned}
 G_N(x,x)&=G_N(0,0)=g\log N+O(1),\\
 |G_N(x,y)|&\le C\left(1+\log\frac N{\max\{d(x,y),1\}}\right),\\
 0\le G_N(x,x)-G_N(x,y)&\le C\log(2+d(x,y)).
 \end{aligned}
\end{equation}
In particular, if $d(x,y)\ge Ne^{-L^{1/8}}/4$, then
$|G_N(x,y)|\le CL^{1/8}$.
\end{proposition}

\begin{proof}
For $k\ne0$, the elementary inequality
$0\le u^2-\sin^2u\le u^4/3$ on $|u|\le\pi/2$ gives
\[
 0\le\frac{\pi^2|k|^2}{N^2}-\lambda_N(k)\le\frac{\pi^4|k|^4}{3N^4}.
\]
Together with the lower bound in \eqref{chI:eq:spectrum}, this implies
\begin{equation}\label{chI:eq:quadratic-error}
 0\le\frac1{N^2\lambda_N(k)}-\frac1{\pi^2|k|^2}
       \le\frac{\pi^2}{12N^2}.
\end{equation}
Summing the error yields
\[
 \left|G_N(x,y)-\pi^{-2}\operatorname{Re}\Sigma_N(x-y)\right|
 \le\pi^2/12.
\]
The first two assertions follow from Lemma~\ref{chI:lem:fourier-sums},
translation invariance, and $2\pi/\pi^2=g$.

Taking real parts in \eqref{chI:eq:green} gives the nonnegative sum
\[
 G_N(x,x)-G_N(x,y)
 =\frac1{N^2}\sum_{k\ne0}
       \frac{1-\cos(2\pi k\cdot(x-y)/N)}{\lambda_N(k)}.
\]
Choose a shortest representative of $x-y$, of length $d(x,y)$.
The inequalities $1-\cos v\le\min\{v^2/2,2\}$ show that each
summand is at most $C\min\{d(x,y)^2/N^2,|k|^{-2}\}$.
If $x\ne y$, then $1\le d(x,y)\le N$; apply the last bound in
\eqref{chI:eq:quadratic-sums} with $u=N/d(x,y)$. This gives
$C(1+\log d(x,y))\le C\log(2+d(x,y))$. For $x=y$ the increment
is zero. Finally, substitute
$\log(N/d(x,y))\le L^{1/8}+\log4$ in \eqref{chI:eq:green-bounds}.
\end{proof}

\subsection*{1.3. Spatial differences of the Green function}

The bounded error in \eqref{chI:eq:quadratic-error} is sufficient for
\eqref{chI:eq:green-bounds}, but does not give a decaying bound for
spatial differences. We therefore estimate the differentiated
Fourier multiplier itself.

\begin{lemma}\label{chI:lem:unit-differences}
For every $z\ne0$ and unit lattice directions $e,e'$,
\begin{equation}\label{chI:eq:unit-differences}
 \begin{aligned}
 |G_N(z+e,0)-G_N(z,0)|&\le\frac C{d(z,0)},\\
 |G_N(z+e+e',0)-G_N(z+e,0)&-G_N(z+e',0)+G_N(z,0)|
       \le\frac C{d(z,0)^2}.
 \end{aligned}
\end{equation}
\end{lemma}

\hyperlink{proof.appendix.B}{Appendix~B} proves this lemma by a smooth partition of the periodic
frequency multiplier. On a block of frequency radius $R$, the
absolute sum is of order $(R/N)^m$ for $m$ spatial differences;
three summations by parts supply the additional factor
$\min\{1,(Rd(z,0)/N)^{-3}\}$. Summing the blocks gives
$d(z,0)^{-m}$ for $m=1,2$.

\begin{proposition}\label{chI:prop:green-differences}
Let $D>0$ and suppose that
\[
 d(x,x'),d(y,y')\le D/4,\qquad
 \min\{d(x,y),d(x,y'),d(x',y),d(x',y')\}\ge D.
\]
Then
\begin{equation}\label{chI:eq:green-differences}
 \begin{aligned}
 |G_N(x,y)-G_N(x',y)|&\le C\frac{d(x,x')}{D},\\
 |G_N(x,y)-G_N(x',y)-G_N(x,y')+G_N(x',y')|
   &\le C\frac{d(x,x')d(y,y')}{D^2}.
 \end{aligned}
\end{equation}
Also,
\begin{equation}\label{chI:eq:spectral-square}
 \sum_{k\ne0}\frac1{N^4\lambda_N(k)^2}\le C.
\end{equation}
\end{proposition}

\begin{proof}
Take coordinate paths $(x_j)_{j=0}^{a}$ from $x$ to $x'$ and
$(y_k)_{k=0}^{b}$ from $y$ to $y'$, using shortest coordinate
displacements. Then
\[
 \begin{gathered}
 a\le\sqrt2d(x,x'),\qquad b\le\sqrt2d(y,y'),\\
 d(x_j,x)\le d(x,x'),\qquad d(y_k,y)\le d(y,y'),\\
 d(x_j,y_k)\ge D-d(x,x')-d(y,y')\ge D/2.
 \end{gathered}
\]
By \eqref{chI:eq:unit-differences},
\[
 \begin{aligned}
 |G_N(x',y)-G_N(x,y)|
 &\le\sum_{j=0}^{a-1}|G_N(x_{j+1},y)-G_N(x_j,y)|\\
 &\le 2Ca/D\le 2\sqrt2C\,d(x,x')/D.
 \end{aligned}
\]
For the mixed difference, the exact identity is
\[
\begin{aligned}
 &G_N(x,y)-G_N(x',y)-G_N(x,y')+G_N(x',y')\\
 &=\sum_{j=0}^{a-1}\sum_{k=0}^{b-1}
 \bigl[G_N(x_{j+1},y_{k+1})-G_N(x_{j+1},y_k)\\
 &\hspace{42mm}-G_N(x_j,y_{k+1})+G_N(x_j,y_k)\bigr].
\end{aligned}
\]
Translation invariance turns each bracket into a second difference
with directions $x_{j+1}-x_j$ and $-(y_{k+1}-y_k)$.
Its absolute value is at most $4C/D^2$; hence the sum is bounded by
$4Cab/D^2\le8C\,d(x,x')d(y,y')/D^2$. Empty sums are zero, and all path
coordinates are interpreted modulo $N$.

For \eqref{chI:eq:spectral-square}, use
\[
 \sum_{k\ne0}\frac1{N^4\lambda_N(k)^2}
 \le\frac1{16}\sum_{k\in\ZZ^2\setminus\{0\}}|k|^{-4}<\infty.
\]
Indeed, the square shell $|k|_\infty=j$ has $8j$ points, and its
contribution is at most $8j^{-3}$.
\end{proof}

\subsection*{1.4. Capacity and the killed generator}

Fix $N$. For real functions on $V_N$, write
\[
 \pi_N(f)=N^{-2}\sum_xf(x),\qquad
 \mathcal E_N(f,h)=-\sum_xf(x)\cL_Nh(x),\qquad
 \mathcal E_N(f)=\mathcal E_N(f,f).
\]
In this normalization the energy has no factor $N^{-2}$:
\begin{equation}\label{chI:eq:energy}
 \mathcal E_N(f)=\frac18\sum_x\sum_{e\in\{\pm e_1,\pm e_2\}}
                       (f(x+e)-f(x))^2.
\end{equation}
The Fourier representation implies Poincar\'e's inequality
\begin{equation}\label{chI:eq:poincare}
 \mathcal E_N(f)\ge N^2\gamma_N
                      \bigl(\pi_N(f^2)-\pi_N(f)^2\bigr).
\end{equation}

Let $\varnothing\ne B\subsetneq V_N$. The matrix
$G_N[B]=(G_N(x,y))_{x,y\in B}$ is positive definite: extending a
vector by zero on $B^c$, a zero quadratic form would make it a
constant vector on $V_N$, hence zero. We may therefore define
\[
 q_B=G_N[B]^{-1}\ind_B,
\]
where $\ind_B$ is the all-one vector indexed by $B$. Extend $q_B$
by zero on $B^c$ and put
\[
 c_B=\sum_{x\in B}q_B(x),\qquad \psi_B=G_Nq_B.
\]
We call $c_B$ the capacity
of $B$ for the mean-zero Green kernel. The restriction to proper
subsets is part of this definition.

\begin{proposition}\label{chI:prop:capacity}
The equilibrium potential and charges satisfy
\begin{equation}\label{chI:eq:equilibrium}
 \begin{gathered}
 c_B>0,\qquad q_B(x)\ge c_B/N^2\quad(x\in B),\\
 \psi_B|_B=1,\qquad\pi_N(\psi_B)=0,\qquad\psi_B\le1,\\
 -\cL_N\psi_B=q_B-c_B/N^2.
 \end{gathered}
\end{equation}
Moreover,
\begin{equation}\label{chI:eq:capacity-variation}
 c_B=\min\{\mathcal E_N(f):\pi_N(f)=0,\ f|_B=1\},
\end{equation}
and $\psi_B$ is the unique minimizer.
\end{proposition}

\begin{proof}
The definition and the Green operator identities give
$\psi_B|_B=1$, $\pi_N(\psi_B)=0$, and the last identity in
\eqref{chI:eq:equilibrium}. Thus
\[
 \mathcal E_N(\psi_B)
 =\sum_x\psi_B(x)\bigl(q_B(x)-c_B/N^2\bigr)=c_B>0.
\]
Strict positivity follows from \eqref{chI:eq:poincare}, since
$\psi_B$ is nonzero and has mean zero. On $B^c$,
$\cL_N\psi_B=c_B/N^2>0$, so a global maximum cannot occur there.
It follows that $\psi_B\le1$. For $x\in B$, we consequently have
$-\cL_N\psi_B(x)\ge0$, which proves the charge bound.

If $f$ is admissible in \eqref{chI:eq:capacity-variation}, put
$h=f-\psi_B$. Then $h|_B=0$ and $\pi_N(h)=0$, so
\[
 \mathcal E_N(h,\psi_B)
 =\sum_x h(x)(q_B(x)-c_B/N^2)=0.
\]
Hence $\mathcal E_N(f)=c_B+\mathcal E_N(h)\ge c_B$.
Equality forces $h$ to be constant by \eqref{chI:eq:energy}, and its
mean is zero. Thus $h=0$.
\end{proof}

Define the positive killed operator on functions on $B^c$ by
\[
 K_B=(-\cL_N)_{B^c,B^c},\qquad
 (K_Bf)(x)=-(\cL_N\widetilde f)(x)\quad(x\in B^c),
\]
where $\widetilde f$ extends $f$ by zero on $B$. Thus $K_B$ is
the negative of the generator killed upon entering $B$.
We use this zero extension whenever a function on $B^c$ occurs
in a sum or inner product over $V_N$.
Let $\lambda_B$ be its least eigenvalue. Its eigenvalues are
positive, since zero energy would force a constant function
vanishing on a nonempty set. A first eigenfunction can be chosen
nonnegative: replacing a Rayleigh minimizer by its absolute value
does not increase the energy. Normalize it by
$\pi_N(\varphi_B)=1$ and put
\[
 r_B=\frac{c_B}{N^2\gamma_N},\qquad
 v_B=\pi_N(\varphi_B^2)-1.
\]

\begin{proposition}\label{chI:prop:killed}
If $r_B<1$, the first eigenspace is one-dimensional and
\begin{equation}\label{chI:eq:killed-bounds}
 (1-r_B)\frac{c_B}{N^2}\le\lambda_B\le\frac{c_B}{N^2},
 \qquad 0\le v_B\le\frac{r_B}{1-r_B}.
\end{equation}
For every $f$ vanishing on $B$ with
$\sum_x f(x)\varphi_B(x)=0$,
\begin{equation}\label{chI:eq:killed-gap}
 \mathcal E_N(f)\ge\gamma_N\sum_x f(x)^2.
\end{equation}
In particular, the normalized nonnegative first eigenfunction is
unique.
\end{proposition}

\begin{proof}
The trial function $1-\psi_B$ vanishes on $B$ and has squared norm
$N^2(1+\pi_N(\psi_B^2))$. Therefore
\[
 \lambda_B\le\frac{c_B}{N^2(1+\pi_N(\psi_B^2))}
             \le\frac{c_B}{N^2}<\gamma_N.
\]
Poincar\'e's inequality applied to $\varphi_B$ gives
\[
 N^2\gamma_Nv_B\le\mathcal E_N(\varphi_B)
 =N^2\lambda_B(1+v_B)\le c_B(1+v_B).
\]
This proves the bound on $v_B$. The function $1-\varphi_B$ is
admissible for \eqref{chI:eq:capacity-variation}, whence
\[
 c_B\le N^2\lambda_B(1+v_B),\qquad
 \lambda_B\ge\frac{c_B}{N^2(1+v_B)}
              \ge(1-r_B)\frac{c_B}{N^2}.
\]

If $f$ is another first eigenfunction, then
$h=f-\pi_N(f)\varphi_B$ is a first eigenfunction of mean zero.
Thus
\[
 \lambda_B\sum_xh(x)^2=\mathcal E_N(h)
                      \ge\gamma_N\sum_xh(x)^2.
\]
Since $\lambda_B<\gamma_N$, we have $h=0$, proving simplicity.

For \eqref{chI:eq:killed-gap}, suppose $f\ne0$,
$\sum_xf(x)\varphi_B(x)=0$, and
$\mathcal E_N(f)<\gamma_N\|f\|_2^2$. Put
$h=f-\pi_N(f)\varphi_B$. Then
\[
 \pi_N(h)=0,\qquad
 \mathcal E_N(f,\varphi_B)=\lambda_B\sum_xf(x)\varphi_B(x)=0,
\]
and
\[
 \|h\|_2^2=\|f\|_2^2+\pi_N(f)^2\|\varphi_B\|_2^2>0,\qquad
 \mathcal E_N(h)
 =\mathcal E_N(f)+\pi_N(f)^2\lambda_B\|\varphi_B\|_2^2
 <\gamma_N\|h\|_2^2.
\]
This contradicts \eqref{chI:eq:poincare}.

\end{proof}

The capacity therefore controls the exponential rate of avoidance.
The following expansion also controls its prefactor and remainder.

\begin{proposition}\label{chI:prop:stationary}
Under the assumptions of Proposition~\ref{chI:prop:killed}, let
$S_B(t)=\PP_{\pi_N}(H_B>t)$. For every $t\ge0$,
\begin{equation}\label{chI:eq:stationary}
 \begin{aligned}
 S_B(t)&=\frac{e^{-\lambda_Bt}}{1+v_B}+E_B(t),\\
 0\le E_B(t)&\le
 \left(1-\frac{|B|}{N^2}-\frac1{1+v_B}\right)e^{-\gamma_Nt}
 \le v_Be^{-\gamma_Nt}.
 \end{aligned}
\end{equation}
If $|B^c|=1$, the remainder is zero.
\end{proposition}

\begin{proof}
The killed semigroup acts by
\[
 (e^{-tK_B}f)(x)=\sum_{y\in B^c}\PP_x(X_t=y,H_B>t)f(y),
 \qquad x\in B^c.
\]
Indeed, conditioning on the first jump gives the finite backward
equations with zero boundary on $B$, whose solution is $e^{-tK_B}$.
Take a real orthonormal eigenbasis
$(e_j)_{1\le j\le |B^c|}$ in $L^2(\pi_N)$, extended by zero on $B$.
We may choose
$e_1=\varphi_B/(1+v_B)^{1/2}$, with eigenvalue $\lambda_B$;
all other eigenvalues $\lambda_j$ are at least $\gamma_N$ by
\eqref{chI:eq:killed-gap}. It follows that
\[
 S_B(t)=\langle\ind_{B^c},e^{-tK_B}\ind_{B^c}\rangle_{\pi_N}
       =\sum_j a_je^{-\lambda_jt},\qquad
 a_j=\langle\ind_{B^c},e_j\rangle_{\pi_N}^{\,2}.
\]
Here $a_j\ge0$, $a_1=(1+v_B)^{-1}$, and Parseval's identity gives
\[
 \sum_{j\ge2}a_j
 =1-|B|/N^2-(1+v_B)^{-1}\le v_B/(1+v_B)\le v_B.
\]
The sum over $j\ge2$ is the required remainder; it is empty when
$|B^c|=1$.
\end{proof}

\subsection*{1.5. Hitting, local covering, and survival comparison}

We first control the extra time needed to cover a small set after
entering it. The starting point in the following covering estimate
belongs to the target set.

\begin{lemma}\label{chI:lem:hitting}
For all $x,y\in V_N$,
\begin{equation}\label{chI:eq:mean-hit}
 \EE_xH_y\le N^2\bigl(G_N(y,y)-G_N(x,y)\bigr).
\end{equation}
\end{lemma}

\begin{proof}
For fixed $y$, put $F(x)=N^2(G_N(y,y)-G_N(x,y))$.
Then $F\ge0$, $F(y)=0$ and $\cL_NF=-1$ off $y$.
Dynkin's formula at $t\wedge H_y$ gives
\[
 \EE_xF(X_{t\wedge H_y})+\EE_x(t\wedge H_y)=F(x).
\]
Nonnegativity and monotone convergence give \eqref{chI:eq:mean-hit}.
\end{proof}

\begin{proposition}\label{chI:prop:local-cover}
For $B\ne\varnothing$ and $x\in B$, put
$M_B=\max_{u,v\in B}\EE_uH_v$. Then
\begin{equation}\label{chI:eq:matthews}
 \EE_xT_B\le\left(\sum_{j=1}^{|B|-1}\frac1j\right)M_B,
\end{equation}
where an empty sum is zero. Consequently, if
$\operatorname{diam}_d B\le d_0<N/2$, with $d_0\ge0$, then
\begin{equation}\label{chI:eq:local-cover}
 \sup_{x\in B}\EE_xT_B\le CN^2(1+\log(d_0+1))^2.
\end{equation}
\end{proposition}

\begin{proof}
Equation \eqref{chI:eq:matthews} is Matthews' upper bound
\cite{Matthews88}, applied to the target set $B$ with $x\in B$.
If $\operatorname{diam}_dB\le d_0<N/2$, then
\eqref{chI:eq:green-bounds} and \eqref{chI:eq:mean-hit} give
\[
 M_B\le CN^2(1+\log(d_0+1)),\qquad |B|\le C(d_0+1)^2.
\]
Combining these with $\sum_{j=1}^{m}j^{-1}\le1+\log m$ for $m\ge1$
proves \eqref{chI:eq:local-cover}; for $|B|=1$, $T_B=0$.
\end{proof}

The following comparison transfers the stationary expansion to an
arbitrary starting point and retains the subsequent survival factor
in the loss from an early hit.

\begin{proposition}\label{chI:prop:survival-comparison}
Suppose $u\ge0$, $0\le\varepsilon<1$, and
\[
 \sup_{x,y}|N^2p_u(x,y)-1|\le\varepsilon.
\]
For every $B\subseteq V_N$, $y\in V_N$, and $T>2u$,
\begin{equation}\label{chI:eq:survival-comparison}
 \begin{aligned}
 \PP_y(H_B>T)&\le(1+\varepsilon)S_B(T-u),\\
 \PP_y(H_B>T)&\ge(1-\varepsilon)S_B(T-u)\\
 &\quad-(1+\varepsilon)\PP_y(H_B\le u)S_B(T-2u),
 \end{aligned}
\end{equation}
where $S_B(t)=\PP_{\pi_N}(H_B>t)$.
\end{proposition}

\begin{proof}
Put $f(w)=\PP_w(H_B>T-u)$. The Markov property at time $u$ gives
\[
 \PP_y(H_B>T)=\EE_yf(X_u)
             -\EE_y[\ind_{\{H_B\le u\}}f(X_u)].
\]
Since $f\ge0$,
\[
 (1-\varepsilon)S_B(T-u)
 \le\sum_w p_u(y,w)f(w)
 \le(1+\varepsilon)S_B(T-u).
\]
Dropping the nonnegative subtracted term proves the upper bound.
Apply that bound with initial point $w$ and horizon $T-u>u$ to obtain
\[
 \sup_w f(w)\le(1+\varepsilon)S_B(T-2u).
\]
Use this bound on the subtracted term to obtain the lower bound.
The argument also covers the empty and full target sets.
\end{proof}

In applications we take $u=h_N$ and $\varepsilon=L^{-8}$.
The remaining input is then an early-hit estimate for the particular
target geometry. Such estimates will be made where the targets
are introduced, using \eqref{chI:eq:occupation} and the strong Markov
property.

\clearpage
\endgroup

\begingroup
\mergedchaptersetup{II}{Gaussian fields and the critical measure}
\begin{center}
{\large\scshape Chapter II\par}\vspace{6pt}
{\LARGE Gaussian fields and the critical measure\par}
\end{center}
\vspace{6pt}
\noindent\emph{Standard pinned-cluster and high-point
results are cited explicitly; the disk-weight matching, square cap
comparison and exterior endpoint calculation are proved in \hyperlink{proof.appendix.E}{Appendix~E}.}
\par\medskip
\subsection*{2.1. Finite Gaussian fields and conditional decompositions}

We retain the walk normalization of Chapter~I. In particular,
$G_N$ is the mean-zero Green kernel for a walk of total jump rate
one. This section fixes the Gaussian fields associated with that
kernel and the decompositions used below. All fields in the section
are finite-dimensional; conditioning is Gaussian regression, including
conditioning on a prescribed value of a coordinate.
We use $\mathbb P$ and $\mathbb E$ for the Gaussian law under discussion;
walk expectations retain a starting-point subscript, as in $\mathbb E_x$.
For a random element $Y$, $\mathcal L(Y)$ denotes its law, and
$U\perp V$ means independence. Matrix products and vector norms use
the counting-measure conventions of Chapter~I.

For each nonzero $k\in\mathcal K_N$, let $\xi_k,\zeta_k$ be
independent standard normal random variables, independent also
between frequencies. Define
\begin{equation}\label{chII:eq:finite-field}
 \phi_N(x)=\sum_{k\in\mathcal K_N\setminus\{0\}}
 \frac{\xi_k\cos(2\pi k\cdot x/N)+\zeta_k\sin(2\pi k\cdot x/N)}
      {N\sqrt{\lambda_N(k)}}.
\end{equation}
The sum uses the full frequency set. The coefficient pairs for the
frequency classes of $k$ and $-k$ are independent whenever
$k\not\equiv-k\pmod N$; the normalization in
\eqref{chII:eq:finite-field} accounts for this convention.

\begin{lemma}\label{chII:lem:finite-field}
The field $\phi_N$ is centered Gaussian, has spatial mean zero,
and satisfies
\begin{equation}\label{chII:eq:finite-covariance}
 \mathbb E[\phi_N(x)\phi_N(y)]=G_N(x,y).
\end{equation}
Its law on $\{f:\pi_N(f)=0\}$ has density proportional to
$\exp\{-\mathcal E_N(f)/2\}$ with respect to Euclidean volume on
that subspace. It is invariant under translations of the torus and
under $f\mapsto-f$.
\end{lemma}

\begin{proof}
By \eqref{chII:eq:finite-field}, $\phi_N$ is centered Gaussian and
\[
 \mathbb E[\phi_N(x)\phi_N(y)]
 =\frac1{N^2}\sum_{k\ne0}
       \frac{\cos(2\pi k\cdot(x-y)/N)}{\lambda_N(k)}=G_N(x,y).
\]
Omission of the zero mode gives $\pi_N(\phi_N)=0$.
The precision operator on this subspace is $-\mathcal L_N$, giving
the asserted density; the covariance symmetries give the invariances.
\end{proof}

There are two useful ways to impose a root. Fix $r\in V_N$ and set
\begin{equation}\label{chII:eq:two-root-fields}
 \begin{aligned}
 \phi_N^{(r)}(x)
   &=\phi_N(x)-\frac{G_N(x,r)}{G_N(r,r)}\phi_N(r),\\
 h^r(x)&=\phi_N(x)-\phi_N(r).
 \end{aligned}
\end{equation}
The denominator is strictly positive by the Fourier formula.
The first field is a regression residual; the second is obtained
by choosing the representative that vanishes at $r$.

\begin{proposition}\label{chII:prop:root-fields}
The field $\phi_N^{(r)}$ is independent of $\phi_N(r)$ and has
the conditional law of $\phi_N$ given $\phi_N(r)=0$. Its covariance is
\begin{equation}\label{chII:eq:conditioned-covariance}
 G_N^{(r)}(x,y)
 =G_N(x,y)-\frac{G_N(x,r)G_N(r,y)}{G_N(r,r)}.
\end{equation}
Both its root value and spatial mean are zero. The covariance of
$h^r$ is instead
\begin{equation}\label{chII:eq:pinned-covariance}
 \widehat G_N^{\,r}(x,y)
 =G_N(x,y)-G_N(x,r)-G_N(r,y)+G_N(r,r).
\end{equation}
Restricted to $V_N\setminus\{r\}$, this matrix is the inverse of
the negative generator killed at $r$. Moreover,
\begin{equation}\label{chII:eq:root-reconstruction}
 \begin{aligned}
 h^r(x)
 &=\phi_N^{(r)}(x)
   +\left(\frac{G_N(x,r)}{G_N(r,r)}-1\right)\phi_N(r),\\
 \phi_N(x)&=h^r(x)
                       -\pi_N(h^r).
 \end{aligned}
\end{equation}
\end{proposition}

\begin{proof}
Gaussian regression gives independence of $\phi_N^{(r)}$ and $\phi_N(r)$,
\eqref{chII:eq:conditioned-covariance}, and the conditional-law version
\[
 \mathcal L(\phi_N\mid\phi_N(r)=s)
 =\mathcal L\left(\phi_N^{(r)}
                    +\frac{G_N(\cdot,r)}{G_N(r,r)}s\right).
\]
The zero column sums of $G_N$ give $\pi_N(\phi_N^{(r)})=0$.
For $h^r$, direct covariance calculation gives
\eqref{chII:eq:pinned-covariance}, and for $y\ne r$,
\[
 \widehat G_N^{\,r}(r,y)=0,\qquad
 -\mathcal L_N\widehat G_N^{\,r}(\cdot,y)
        =\mathbf1_{\{y\}}-\mathbf1_{\{r\}}.
\]
Restriction to $V_N\setminus\{r\}$ gives the killed inverse.
Substitution and $\pi_N(h^r)=-\phi_N(r)$ give
\eqref{chII:eq:root-reconstruction}.
\end{proof}

The distinction matters for domain decompositions. An ordinary
Dirichlet field arises naturally from $h^r$.
For the mean-zero field, conditioning on the exterior also fixes
the sum of the interior values, and removes one additional mode.

Let $\varnothing\ne D\subsetneq V_N$. Write
\[
 K_D=(-\mathcal L_N)_{D,D},\qquad G_D=K_D^{-1},\qquad
 b_D=G_D\mathbf1_D,\qquad v_D=\mathbf1_D^{\mathsf T}G_D\mathbf1_D.
\]
Here $K_D$ is the Dirichlet restriction of $-\mathcal L_N$ to $D$,
and $\mathbf1_D$ is the all-one vector indexed by $D$.
Vectors on $D$ will be extended by zero outside $D$ when used as
functions on $V_N$. The number $v_D$ is positive. An ordinary
Dirichlet field $\eta_D$ has covariance $G_D$ and density
proportional to $\exp\{-\mathcal E_N(f)/2\}$ on functions
vanishing on $D^c$. Put
\begin{equation}\label{chII:eq:domain-constraint}
 \begin{aligned}
 \eta_D^0&=\eta_D-\frac{b_D}{v_D}\sum_{x\in D}\eta_D(x),
 &q_D&=\frac{b_D}{\sqrt{v_D}},\\
 G_D^0&=G_D-\frac{b_Db_D^{\mathsf T}}{v_D},
 &\xi_D&=\frac1{\sqrt{v_D}}\sum_{x\in D}\eta_D(x).
 \end{aligned}
\end{equation}

\begin{proposition}\label{chII:prop:domain-fields}
The variable $\xi_D$ is standard normal and independent of
$\eta_D^0$, and
\begin{equation}\label{chII:eq:dirichlet-splitting}
 \eta_D=\eta_D^0+q_D\xi_D,\qquad
 \operatorname{Cov}(\eta_D^0)=G_D^0,\qquad
 \sum_{x\in D}\eta_D^0(x)=0.
\end{equation}
If $r\notin D$, the field $h^r$ admits the
decomposition
\begin{equation}\label{chII:eq:pinned-domain-splitting}
 h^r=\eta_D+H_D^r,
 \qquad \eta_D\perp H_D^r,
\end{equation}
where $H_D^r$ equals $h^r$ on $D^c$ and is its
harmonic extension to $D$.
For the mean-zero field there is instead a decomposition
\begin{equation}\label{chII:eq:zero-mean-splitting}
 \phi_N=\eta_D^0+H_D^{\mathrm{zm}},\qquad \eta_D^0\perp H_D^{\mathrm{zm}},
\end{equation}
with
\[
 H_D^{\mathrm{zm}}=\mathbb E[\phi_N\mid\sigma(\phi_N(x):x\in D^c)],\qquad
 \operatorname{Cov}(H_D^{\mathrm{zm}})|_{D\times D}
       =G_N|_{D\times D}-G_D+\frac{b_Db_D^{\mathsf T}}{v_D}.
\]
In these two distributional decompositions, $\eta_D$ and $\eta_D^0$
have the stated residual laws and may be realized separately.
\end{proposition}

\begin{proof}
The sum $\sum_D\eta_D$ has variance $v_D$ and covariance $b_D(x)$
with $\eta_D(x)$. Gaussian regression on that sum proves
\eqref{chII:eq:dirichlet-splitting} and independence. Its zero-sum
assertion also follows directly from $\sum_D b_D=v_D$.

By Proposition~\ref{chII:prop:root-fields}, the pinned field has precision
$(-\mathcal L_N)_{V_N\setminus\{r\},V_N\setminus\{r\}}$.
Given its values outside $D$, its
density as a function of the interior vector $f$ is proportional to
\[
 \exp\left\{-\frac12 f^{\mathsf T}K_Df
             -f^{\mathsf T}(-\mathcal L_N)_{D,D^c}
                                h^r|_{D^c}\right\}.
\]
Completing the square gives covariance $G_D$ and mean solving
$\mathcal L_NH_D^r=0$ on $D$, with the prescribed exterior
values. This mean is a linear function of the exterior. The
conditional covariance is deterministic, so the residual is
independent of the exterior and hence of $H_D^r$.

For $\phi_N$, Lemma~\ref{chII:lem:finite-field} gives the same quadratic
energy, but its conditional density is restricted to the affine
hyperplane
\[
 \sum_{x\in D}f(x)=-\sum_{x\in D^c}\phi_N(x).
\]
Starting with an ordinary Gaussian of covariance $G_D$ and
conditioning its sum to this value changes its mean but replaces
its covariance by $G_D^0$. Equivalently, complete the square on
the subspace of vectors of sum zero in $D$. Thus the conditional
residual has the law of $\eta_D^0$ and is independent of the
exterior. Subtracting its covariance from $G_N$ proves the
displayed covariance identity for $H_D^{\mathrm{zm}}$.
\end{proof}

The constrained Euler equation gives
$-\mathcal L_NH_D^{\mathrm{zm}}$ constant on $D$, by testing against vectors
supported in $D$ whose coordinates sum to zero. This constrained mode
produces the rank-one term in \eqref{chII:eq:domain-constraint}.

For the scalar argument later in the chapter, it is more convenient
to keep ordinary Dirichlet fields in many separate squares and
record a single correction for their common spatial mean.
Let $K\ge8$ be an integer, $N\ge2K$, and $m=\lfloor N/K\rfloor$. In the
coordinate representatives $\{0,\ldots,N-1\}^2$, take
\[
 D_i=im+\{1,\ldots,m-1\}^2,
 \qquad i\in\{0,\ldots,K-1\}^2,
 \qquad U=\bigcup_iD_i.
\]
These squares have no edges between distinct interiors, and
$0\notin U$. Let $\eta_i$ be the ordinary Dirichlet residual of
$h^0$ on $D_i$, and set
\[
 \overline\eta=\sum_i\pi_N(\eta_i),\qquad
 H=H_U^0-\pi_N(H_U^0),
\]
Here $H_U^0$ is the harmonic extension in
\eqref{chII:eq:pinned-domain-splitting} with root $0$.

\begin{proposition}\label{chII:prop:block-fields}
There is the exact identity
\begin{equation}\label{chII:eq:block-splitting}
 \phi_N=\sum_i\eta_i+H-\overline\eta.
\end{equation}
The fields $\eta_i$ are mutually independent and their entire
collection is independent of $H$. Furthermore,
\begin{equation}\label{chII:eq:fine-mean}
 \mathbb E\overline\eta=0,\qquad
 \operatorname{Var}(\overline\eta)\le\frac9{K^2},\qquad
 \mathbb P(|\overline\eta|\ge a)\le\frac9{K^2a^2}\quad(a>0).
\end{equation}
For bounded measurable functions $F_i$ with values in $[0,1]$,
\begin{equation}\label{chII:eq:block-product}
 \mathbb E\left[\prod_iF_i(H,\eta_i)\,\middle|\,H\right]
 =\prod_i\mathbb E[F_i(H,\eta_i)\mid H].
\end{equation}
\end{proposition}

\begin{proof}
Apply \eqref{chII:eq:pinned-domain-splitting} to $U$ with root $0$.
Since there are no edges between distinct $D_i$, its killed
generator is block diagonal. The residual consequently has the
law of independent Dirichlet fields on those squares, and is
independent of the whole harmonic extension. Centering this
decomposition, and then using \eqref{chII:eq:root-reconstruction},
gives \eqref{chII:eq:block-splitting}. Independence persists after
centering the harmonic extension. Conditional integration against
the product law of the $\eta_i$ proves \eqref{chII:eq:block-product}.

For the variance, independence and the Dirichlet covariances give
\begin{equation}\label{chII:eq:fine-mean-variance}
 \operatorname{Var}(\overline\eta)
 =N^{-4}\sum_i\sum_{x,y\in D_i}G_{D_i}(x,y).
\end{equation}
The row sum $b_{D_i}(x)=\sum_yG_{D_i}(x,y)$ is the mean exit
time from $D_i$. This identity follows either from the killed
semigroup integral or from the equations
$-\mathcal L_N b_{D_i}=1$ in $D_i$, with zero boundary values.
Each square is contained in a periodic Euclidean disk of radius
$2m$ about its corner. In that disk use the nonnegative function
\[
 z\longmapsto(2m+1)^2-|z|^2.
\]
The representatives here are taken relative to the disk center;
the disk and its neighboring vertices fit in the torus without
ambiguity because $2(2m+1)<N$ under our assumptions. A direct
calculation gives $\mathcal L_N|z|^2=1$. Stopping the one-step
identity for this function at the disk exit, as in the mean hitting-time
bound of Chapter~I, shows that the mean exit time is at most
$(2m+1)^2\le9m^2$. Exit from a contained square occurs no later.
Thus each row sum in \eqref{chII:eq:fine-mean-variance} is at most
$9m^2$. Since $\sum_i|D_i|\le N^2$ and $m\le N/K$, the
variance is at most $9/K^2$. The mean is zero, and Chebyshev's
inequality gives the last assertion of \eqref{chII:eq:fine-mean}.
\end{proof}

The mean correction in \eqref{chII:eq:block-splitting} is measurable
from the fine fields. We apply \eqref{chII:eq:block-product} before
height averaging, then remove the correction using
\eqref{chII:eq:fine-mean} and \eqref{chII:eq:dependent-translation}.

\subsection*{2.2. The periodic field and heat regularization}

Write $\mathbb T^2=\mathbb R^2/\mathbb Z^2$, with Lebesgue measure of
mass one, and put $g=2/\pi$ and $\alpha=2/\sqrt g$. The continuous
Laplacian in this chapter is $\Delta=\partial_1^2+\partial_2^2$.
For $s\in\mathbb R$, $H^s(\mathbb T^2)$ is the periodic Sobolev space
with squared Fourier norm
$\sum_{k\in\mathbb Z^2}(1+4\pi^2|k|^2)^s|\widehat u(k)|^2$,
where $\widehat u(k)=\langle u,e^{-2\pi i k\cdot x}\rangle$ is the
distributional Fourier coefficient.
Put $P_t^{\rm per}=e^{t\Delta}$ and write $p_t^{\rm per}(x,y)$ for its kernel with
respect to this Lebesgue measure.
For $k\in\mathbb Z^2\setminus\{0\}$ set $\lambda_k=4\pi^2|k|^2$
and choose the half-plane
$\mathcal K_+=\{k:k_1>0\}\cup\{(0,k_2):k_2>0\}$.
On a probability space carrying independent standard normal variables
$(\xi_k,\zeta_k)_{k\in\mathcal K_+}$, define
\begin{equation}\label{chII:eq:continuum-heat}
 h_\varepsilon(x)=\sum_{k\in\mathcal K_+}
 \sqrt{\frac8{\lambda_k}}e^{-\varepsilon^2\lambda_k}
 \bigl(\xi_k\cos(2\pi k\cdot x)+\zeta_k\sin(2\pi k\cdot x)\bigr).
\end{equation}
The factor $8$ belongs to the half-frequency convention. It is
consistent with the full-frequency discrete convention in
\eqref{chII:eq:finite-field}.

\begin{proposition}\label{chII:prop:periodic-heat}
There is a mean-zero Gaussian random distribution $h$, with covariance
$4(-\Delta)^{-1}$ on mean-zero functions, such that
$h_\varepsilon=e^{\varepsilon^2\Delta}h$. Almost surely,
$h\in H^{-\nu}(\mathbb T^2)$ for every $\nu>0$ and
$h_\varepsilon\to h$ in that space. For positive $\varepsilon,\delta$,
\begin{equation}\label{chII:eq:heat-covariance}
 \EE[h_\varepsilon(x)h_\delta(y)]
 =\sum_{k\in\mathcal K_+}\frac8{\lambda_k}
 e^{-(\varepsilon^2+\delta^2)\lambda_k}\cos(2\pi k\cdot(x-y)).
\end{equation}
The fields at positive scales are smooth. Their variance is independent
of $x$, and
\begin{equation}\label{chII:eq:heat-variance}
 V_\varepsilon:=\EE h_\varepsilon(x)^2
 =4\int_{2\varepsilon^2}^{\infty}
       \bigl(p_t^{\rm per}(0,0)-1\bigr)\,dt
 =g\log(1/\varepsilon)+O(1).
\end{equation}
\end{proposition}
\begin{proof}
The expected squared $H^{-\nu}$ norm of the unsmoothed series is a
constant times $\sum_{k\ne0}|k|^{-2}(1+|k|^2)^{-\nu}<\infty$.
Thus its nonnegative squared Fourier norm is finite almost surely;
the series defines $h$ and dominated convergence gives the asserted
heat convergence through all positive real scales. Take the countable
intersection of these probability-one events over $q\in\mathbb Q_{>0}$,
using the same Gaussian Fourier coefficients. For every $\nu>0$, choose
$0<q<\nu$ rational. The inequality
\[
 \|u\|_{H^{-\nu}(\mathbb T^2)}
 \le\|u\|_{H^{-q}(\mathbb T^2)}
\]
gives both membership and convergence for every $\nu>0$ on that event.
Exponential decay of the multipliers, followed by
Cauchy--Schwarz, gives absolute uniform convergence of every derivative
at positive scales. Independence gives \eqref{chII:eq:heat-covariance}.
Integrating the spectral expansion of the heat kernel gives the
first equality in \eqref{chII:eq:heat-variance}. For small $t$, the periodic
image formula gives
$p_t^{\rm per}(0,0)=(4\pi t)^{-1}+O(t^{-1}e^{-c/t})$;
for large $t$ its nonconstant modes decay exponentially.
The logarithmic term is therefore
$4(4\pi)^{-1}\log(1/(2\varepsilon^2))=g\log(1/\varepsilon)+O(1)$.
\end{proof}

All subsequent chaos measures use exact variance subtraction:
\begin{equation}\label{chII:eq:periodic-chaos-approximant}
 \mathcal Z_\varepsilon(dx)=\sqrt{\log(1/\varepsilon)}
 \exp\{\alpha h_\varepsilon(x)-\tfrac12\alpha^2V_\varepsilon\}\,dx,
 \qquad 0<\varepsilon<1.
\end{equation}
Replacing $V_\varepsilon$ by its logarithmic main term would change
the normalization of the limit. We keep the exact expression throughout.

\subsection*{2.3. The critical measure on the original field}

For a measure $\mu$ and an integrable test $f$, write
$\mu(f)=\int f\,d\mu$. Weak convergence of finite measures on
$\mathbb T^2$ is tested against continuous functions. On an open
domain, local weak convergence is tested against continuous functions
of compact support.

\begin{theorem}\label{chII:thm:critical-measure}
There is an $h$-measurable finite random measure $\mathcal Z$ such that
$\mathcal Z_\varepsilon\to\mathcal Z$ weakly in probability as
$\varepsilon\downarrow0$, through all positive real scales. This limit
is unique up to almost sure equality. Almost surely it is nonatomic
and charges every nonempty open set. In particular,
\begin{equation}\label{chII:eq:critical-mass}
 Z:=\mathcal Z(\mathbb T^2)\quad\hbox{satisfies}\quad 0<Z<\infty.
\end{equation}
\end{theorem}

We give the mechanism of the construction here; the heat-locality and
measurability arguments are proved in \hyperlink{proof.appendix.C}{Appendix~C}. The external
critical-chaos input is the independent-increment Dirichlet cutoff
theorem of Duplantier--Rhodes--Sheffield--Vargas \cite[Theorem~10]{chII:DRSV},
in the normalization stated there and rescaled to covariance
$4(-\Delta_D)^{-1}$. This input is used for existence, nontriviality
and nonatomicity.

For a centered smooth Gaussian regularization $U_\varepsilon$, write
\[
 \mathcal M_\varepsilon[U](dx)
 =\sqrt{\log(1/\varepsilon)}
  e^{\alpha U_\varepsilon(x)-\alpha^2\EE[U_\varepsilon(x)^2]/2}\,dx,
 \qquad 0<\varepsilon<1.
\]
Let $D$ be a bounded open Euclidean square, and let $\Delta_D$ be
the Dirichlet Laplacian.
Put $P_t^D=e^{t\Delta_D}$,
with Dirichlet heat kernel $p_t^D(x,y)$ relative to Lebesgue measure.
Let $W$ be space--time white noise on $(0,\infty)\times D$,
equivalently cylindrical Wiener noise on $L^2(D)$. Consider
\begin{equation}\label{chII:eq:dirichlet-three-fields}
 h_D=2\int_0^\infty P_{s/2}^D W(ds),\qquad
 X_\varepsilon=2\int_{2\varepsilon^2}^\infty P_{s/2}^D W(ds),\qquad
 Y_\varepsilon=P_{\varepsilon^2}^Dh_D.
\end{equation}
At each fixed scale $X_\varepsilon$ and $Y_\varepsilon$ have the same
covariance, namely $4\int_{2\varepsilon^2}^\infty P_s^D\,ds$.
The cutoff theorem gives a local weak probability limit $M_D$ for
the exact-variance chaos of $X_\varepsilon$.
Lemma~\ref{chII:lem:critical-limit-measurable} proves that
$M_D=F_D(h_D)$ for a Borel map $F_D$. The joint limits are
\[
 (Y_\varepsilon,\mathcal M_\varepsilon[Y])
 \overset d=(X_\varepsilon,\mathcal M_\varepsilon[X])
 \Longrightarrow(h_D,F_D(h_D)).
\]
Since $Y_\varepsilon-h_D\to0$ in probability in every fixed
negative Sobolev space under consideration, Slutsky's lemma gives
\[
 (h_D,\mathcal M_\varepsilon[Y])
 \Longrightarrow (h_D,F_D(h_D)).
\]
Lemma~\ref{chII:lem:graph-collapse} from Appendix~C upgrades this joint
limit to convergence in probability on the original field.

The following notation refers to continuum domains. Let $G_D(x,y)$
be the covariance kernel of $4(-\Delta_D)^{-1}$ and put
\[
 b_D(x)=\int_DG_D(x,y)\,dy,\qquad
 v_D=\int_Db_D(x)\,dx,\qquad q_D=b_D/\sqrt{v_D}.
\]
For the ordinary Dirichlet field define
$\xi_D=v_D^{-1/2}\int_Dh_D$ and $h_D^0=h_D-q_D\xi_D$.
The integral is the mean-square Gaussian pairing with $\mathbf1_D$.
For the periodic field, suppose additionally that $\overline D$ lies
in a torus coordinate chart. Let $\mathcal F_{D^c}$ be the completed
sigma-field generated by its pairings with smooth tests supported
in the exterior of $D$, and set
$H_D=\EE[h|_D\mid\mathcal F_{D^c}]$ in the sense of distributions.
The Gaussian domain projection has residual covariance
$G_D^0:=G_D-b_D\otimes b_D/v_D$, the continuum counterpart of
Proposition~\ref{chII:prop:domain-fields}. On the original periodic-field
space set $h_D^0:=h|_D-H_D$; this residual has the constrained
Dirichlet law above. On an extension adjoin a standard normal $\xi_D$
independent of the entire original field, and define
$h_D:=h_D^0+q_D\xi_D$. From this point on, $h_D^0$, $\xi_D$ and $h_D$
refer to this realization on the extended probability space. Its covariance is
\[
 \operatorname{Cov}(h_D)=G_D^0+q_D\otimes q_D=G_D,
\]
so $h_D$ has the ordinary Dirichlet-field law. Thus
\begin{equation}\label{chII:eq:continuum-chart}
 h_D=h_D^0+q_D\xi_D,\qquad h|_D=h_D^0+H_D,
 \qquad h_D^0\perp\xi_D,\quad h_D^0\perp H_D.
\end{equation}
Here $H_D$ is the binding field for the constrained residual.
Both $q_D$ and $H_D$ are smooth in the
interior. Write $q_{D,\varepsilon}=P_{\varepsilon^2}^Dq_D$ and
$H_{D,\varepsilon}=P_{\varepsilon^2}^DH_D$, and bind
\[
 Y_\varepsilon:=P_{\varepsilon^2}^Dh_D.
\]
Here $\mathcal M_\varepsilon[Y]$ denotes the exact-variance chaos
constructed from this $Y_\varepsilon$. Since $h_D$ has the ordinary
Dirichlet law, the Borel map $F_D$ constructed above applies to this
realization and transfers the same-field convergence
$\mathcal M_\varepsilon[Y]\to F_D(h_D)$ in probability.
Let $M^0_{D,\varepsilon}$ and
$M^{\rm per,D}_\varepsilon$ be the exact-variance chaos approximants
of $P_{\varepsilon^2}^Dh_D^0$ and $P_{\varepsilon^2}^D(h|_D)$, respectively.
Exact subtraction gives
\begin{equation}\label{chII:eq:chart-factors}
 \begin{aligned}
 M^0_{D,\varepsilon}
 &=e^{-\alpha q_{D,\varepsilon}\xi_D+\alpha^2q_{D,\varepsilon}^2/2}
       \mathcal M_\varepsilon[Y],\\
 M^{\rm per,D}_\varepsilon
 &=e^{\alpha H_{D,\varepsilon}-\alpha^2\EE H_{D,\varepsilon}^2/2}
       M^0_{D,\varepsilon}.
 \end{aligned}
\end{equation}
The smooth factors converge uniformly on compact interiors. Killed
and periodic heat smoothing have the same local limit, including
their variance counterterms. Limits on overlapping charts therefore
agree. A finite deep-chart cover and a partition of unity produce
the finite measure $\mathcal Z$ in the theorem. Positive smooth
factors preserve local positivity and nonatomicity; a countable
base gives full support. Uniqueness is uniqueness of a probability
limit in a metric space of finite measures. The remaining heat-locality and measurability arguments are proved
in \hyperlink{proof.appendix.C}{Appendix~C}.

\subsection*{2.4. High roots and contact probabilities}

Two Gaussian estimates will be needed, at quite different scales.
The first concerns a finite cluster with arbitrary internal geometry.
The second concerns the maximum of a reference square. They use
different common Gaussian coordinates and will be kept separate.

Put $L=\log N$, $v_N=G_N(0,0)$, and
\[
 R_N^{\rm loc}=e^{L^{1/6}},\qquad D_N^{\rm sep}=Ne^{-2L^{1/8}}.
\]
Fix a root $r\in V_N$. Distances in this paragraph use the periodic
Euclidean metric $d$ of Chapter~I: for nonempty sets,
$d(B,B')=\min_{x\in B,y\in B'}d(x,y)$ and $d(B,r)=d(B,\{r\})$.
An admissible cluster is a nonempty set $B\subset V_N$ with
$\operatorname{diam}_d B\le R_N^{\rm loc}$ and
$d(B,r)\ge D_N^{\rm sep}$. A pair is admissible if each cluster is
admissible and $d(B,B')\ge D_N^{\rm sep}$. For a deterministic nonnegative
profile $\ell=(\ell_x)_{x\in V_N}\in[0,\infty)^{V_N}$ and
$0<v\le1$, define, for $z=(z_x)_{x\in B}\in\mathbb R^B$,
\[
 \begin{gathered}
 F_B(z)=\operatorname{Leb}\bigcup_{x\in B}
 \left[-z_x-\sqrt{(v-2\ell_x)_+},\,-z_x+\sqrt{(v-2\ell_x)_+}\right],\\
 a_B^G(\ell;r,v)=(2\pi g)^{-1/2}\EE F_B(\phi_N^{(r)}).
 \end{gathered}
\]
With $\ell,r,v$ fixed, we write $a_B^G$ for this value.
Zero-length intervals have zero Lebesgue measure. Points with
$v<2\ell_x$ can be omitted from the union without changing its
length. Write
\[
 E_B=\{\exists x\in B:2\ell_x+(h^r(x))^2\le v\},
 \qquad S_r=\phi_N(r)/\sqrt{v_N}.
\]

\begin{proposition}\label{chII:prop:relative-contact}
There exist absolute constants $c,C>0$ and an integer $N_0\ge3$
such that, for every $N\ge N_0$, the following estimates hold uniformly
over all admissible clusters, nonnegative profiles and levels $0<v\le1$:
\begin{equation}\label{chII:eq:contact-bounds}
 \begin{aligned}
 c\,a_B^G/\sqrt L
 &\le\PP(E_B,\ |S_r|\le1)\le C\,a_B^G/\sqrt L,\\
 \PP(E_B\cap E_{B'},\ |S_r|\le1)
 &\le C\,a_B^G a_{B'}^G/L
 \end{aligned}
\end{equation}
for admissible pairs. These estimates remain valid for arbitrarily
small positive activities; zero activity gives zero contact probability.
\end{proposition}

Integrate \eqref{chII:eq:R15.22} and~\eqref{chII:eq:R15.24} over
$S_r\in[-1,1]$ with $M=1$. Their uniform relative error and uniform
two-cluster bound give the constants and threshold above.
For a finite family $\mathcal B$ of pairwise admissible clusters,
the truncated count $Y=\ind_{\{|S_r|\le1\}}\sum_B\ind_{E_B}$ satisfies
\[
 \begin{aligned}
 \EE Y&\ge\frac c{\sqrt L}\sum_{B\in\mathcal B}a_B^G,\\
 \EE Y^2&\le\frac C{\sqrt L}\sum_{B\in\mathcal B}a_B^G
            +\frac CL\left(\sum_{B\in\mathcal B}a_B^G\right)^2.
 \end{aligned}
\]
Hence
\begin{equation}\label{chII:eq:contact-amplification}
 \PP\Bigl(\bigcup_{B\in\mathcal B}E_B\Bigr)
 \ge c\frac{\sum_{B\in\mathcal B}a_B^G}
                 {\sqrt L+\sum_{B\in\mathcal B}a_B^G}.
\end{equation}
All these are Gaussian statements for fixed profiles and sets.

For the second estimate let $M\ge2$ be an integer and
$D_M=\{1,\ldots,M-1\}^2\subset\mathbb Z^2$. Let $\eta_M$
be the centered Gaussian field whose covariance $G_{D_M}$ is the Green
kernel of the planar total-rate-one walk killed on leaving $D_M$, and set
\begin{equation}\label{chII:eq:extremal-center}
 m_M=2\sqrt g\log M-\frac34\sqrt g\log\log M,
 \qquad Q_\delta=[\delta,1-\delta]^2.
\end{equation}
Conditional laws at a prescribed root value always mean the
Gaussian regression version.

\begin{proposition}\label{chII:prop:high-root}
There exists $b>0$ such that, for every $0<\delta<1/2$, some
$C_\delta<\infty$ has the following property. For every fixed $t\ge1$,
all sufficiently large even $M$, and every $y/M\in Q_\delta$,
\begin{equation}\label{chII:eq:high-root}
 \left|\log M\,\PP\left(\max_{x\in D_M}\eta_M(x)\le m_M+t
          \,\middle|\,\eta_M(y)=m_M+t\right)-bt\right|
 \le C_\delta\bigl(1+\log^2(t+2)\bigr).
\end{equation}
\end{proposition}

To define the limiting gap law, let $p_t^{\mathbb Z^2}$ be the
transition kernel of the planar nearest-neighbor walk with total
jump rate one. Its potential kernel is
\[
 a_{\rm lat}(x)=\int_0^\infty
  \bigl(p_t^{\mathbb Z^2}(0,0)-p_t^{\mathbb Z^2}(0,x)\bigr)\,dt.
\]
Let $\psi$ be the centered Gaussian field on $\mathbb Z^2$ with
$\psi(0)=0$ and covariance
\[
 \operatorname{Cov}(\psi(x),\psi(z))
 =a_{\rm lat}(x)+a_{\rm lat}(z)-a_{\rm lat}(x-z).
\]
Set $\omega_x:=\alpha a_{\rm lat}(x)-\psi(x)$. The probability law $\nu$ on $\mathbb R^{\mathbb Z^2}$ is defined by
the weak limit in the product topology
\[
 \mathcal L\bigl(\omega\mid \omega_x\ge0
              \text{ for every }|x|_2\le R\bigr)
 \Longrightarrow\nu\qquad(R\to\infty).
\]
The pinned-cluster theorem \cite[Theorem~2.3 and Proposition~5.8]{chII:BL}
gives this limit. Under $\nu$, $\omega_0=0$, $\omega_x\ge0$, and
every bounded sublevel set $\{x:\omega_x\le T\}$ is finite almost
surely. \hyperlink{proof.appendix.E}{Appendix~E} identifies this law with the limiting gaps below
the high maximum in the reference square. Below we use its finite-window contact consequence.

We must retain the deterministic finite parts. For an embedded
open square $D$ and $x\in D$, define
\[
 r_D(x)=4\int_0^1\left(p_t^D(x,x)-\frac1{4\pi t}\right)dt
            +4\int_1^\infty p_t^D(x,x)\,dt.
\]
For $Q=(0,1)^2$ this gives
\begin{equation}\label{chII:eq:heat-finite-parts}
 \begin{aligned}
 r_Q(z)&=4\int_0^1\left(p_t^Q(z,z)-\frac1{4\pi t}\right)dt
             +4\int_1^\infty p_t^Q(z,z)dt,\\
 r_{\rm per}&=4\int_0^1\left(p_t^{\rm per}(0,0)-1-\frac1{4\pi t}\right)dt
             +4\int_1^\infty\bigl(p_t^{\rm per}(0,0)-1\bigr)dt,\\
 w(z)&=e^{\alpha^2r_Q(z)/2},\qquad
 p_{\rm per}=e^{\alpha^2r_{\rm per}/2}.
 \end{aligned}
\end{equation}
The finite-part convention is equivalently
\[
 \begin{aligned}
 4\int_{2\varepsilon^2}^{\infty}p_t^D(x,x)\,dt
    &=g\log(1/\varepsilon)-\pi^{-1}\log2+r_D(x)+o(1),\\
 V_\varepsilon&=g\log(1/\varepsilon)-\pi^{-1}\log2+r_{\rm per}+o(1),
 \end{aligned}
\]
uniformly for $x$ in each compact subset of $D$. In particular the
common cutoff constant cancels when these variances are subtracted.
Set $w=0$ on the boundary. The function $w$ extends continuously
there; this follows either from the image kernels or from the
Dirichlet Green diagonal near the boundary. The remaining finite-part
constants are
\[
 c_{\rm lat}=\lim_{|x|\to\infty}(a_{\rm lat}(x)-g\log|x|),\qquad
 c_{\rm heat}=\frac1\pi\left\{\int_0^1\frac{e^{-1/(4t)}}t\,dt
                    +\int_1^\infty\frac{e^{-1/(4t)}-1}t\,dt\right\}.
\]
The first limit is uniform in direction; both integrals converge.
The lattice-to-heat
diagonal comparison reads, uniformly on compact interiors,
\begin{equation}\label{chII:eq:robin-comparison}
 G_{D_M}(y,y)=g\log M+c_{\rm lat}+c_{\rm heat}+r_Q(y/M)+o(1).
\end{equation}
Here $c_{\rm lat}$ is the constant term in the lattice potential
kernel and $c_{\rm heat}$ converts the planar Green finite part to
the heat convention in \eqref{chII:eq:heat-finite-parts}. Thus $c_{\rm lat}$ and $c_{\rm heat}$ are deterministic constants fixed by
the normalization in \eqref{chII:eq:heat-finite-parts}.
For $T\ge1$ and a nonempty finite set $W\subset\mathbb Z^2$, define
\begin{equation}\label{chII:eq:contact-leading-coefficient}
 c_* =\frac{b}{\sqrt{2\pi g}}
          e^{\alpha^2(c_{\rm lat}+c_{\rm heat})/2}>0,
 \qquad
 J_{W,T}=\int_{-1/2}^{T}e^{-\alpha u}
       \nu\{\exists w_0\in W:|u-\omega_{w_0}|\le1/2\}\,du.
\end{equation}

For $\rho>0$ use the increment seminorm
\[
 [d]_\rho=\sup_{\substack{z,z'\in\overline Q\\z\ne z'}}
           \frac{|d(z)-d(z')|}{|z-z'|^\rho}.
\]
This convention applies also when $\rho>1$; on the square, finite
seminorm then forces $d$ to be constant. For a continuous profile
$d$ on $\overline Q$, let $\eta_M^d(x)=\eta_M(x)+d(x/M)$ and
$D_{M,\delta}=\{x\in D_M:x/M\in Q_\delta\}$.
Let $y_M^d$ be the lexicographically first maximizer of $\eta_M^d$
on $D_M$, and define
\[
 \begin{split}
 p_M(\delta,t,T,W,d)=\PP\bigl\{&y_M^d\in D_{M,\delta},\quad
        \eta_M^d(y_M^d)-m_M-t\in[-1/2,T],\\
 &\exists w_0\in W:\ y_M^d+w_0\in D_{M,\delta},\\
 &\hspace{49mm}|\eta_M^d(y_M^d+w_0)-m_M-t|\le1/2\bigr\}.
 \end{split}
\]
Distinct site values have continuous differences, so ties have probability zero.

\begin{proposition}\label{chII:prop:reference-contact}
Fix $0<\delta<1/2$, $T\ge1$, a nonempty finite $W\subset\mathbb Z^2$,
and $B<\infty$, $\rho>0$. For every $e>0$ there is $t_0\ge1$ such
that, for every $t_1\ge t_0$, all sufficiently large even $M$
satisfy
\begin{equation}\label{chII:eq:reference-contact}
 \sup_{\substack{t_0\le t\le t_1\\
                  \|d\|_\infty+[d]_\rho\le B}}
 \left|\frac{p_M(\delta,t,T,W,d)}{t e^{-\alpha t}}
       -c_*J_{W,T}\int_{Q_\delta}e^{\alpha d(z)}w(z)\,dz\right|\le e.
\end{equation}
\end{proposition}
The order is large fixed height, then a compact height interval,
then large side length. \hyperlink{proof.appendix.E}{Appendix~E}, Section~\ref{chII:sec:root-density},
gives the root-density integration and the removal of the hard boundaries.

\subsection*{2.5. The critical Laplace asymptotic and scalar gluing}

Let $M_Q$ be the same-field exact-variance Dirichlet critical measure
constructed in Section~2.3, and write $\mathcal M(dz)=w(z)M_Q(dz)$.

\begin{proposition}\label{chII:prop:critical-laplace}
For $0<\delta<1/2$ and $B<\infty$,
\begin{equation}\label{chII:eq:critical-laplace}
 \lim_{s\downarrow0}\ \sup_{\substack{d\in C(\overline Q)\\\|d\|_\infty\le B}}
 \left|\frac{-\log\EE\exp\{-s\int_{Q_\delta}e^{\alpha d}\,d\mathcal M\}}
                   {s\log(1/s)}
      -\frac1\alpha\int_{Q_\delta}e^{\alpha d(z)}w(z)\,dz\right|=0.
\end{equation}
\end{proposition}
\begin{proof}
Fix $A>0$, a compact $K\Subset Q$ and a nonnegative continuous $f$
supported in $K$. Appendix~F constructs the barrier truncation
$M_b(f)$ and crossing event $E_b$, and proves
\eqref{chII:eq:barrier-package}, uniformly for bounded $\|f\|_\infty$.
Set $s=e^{-\alpha b}$ and fix $0<q<1$. These estimates give
\[
 sqb\int f-s^2C\|f\|_\infty^2(1+qb)^2e^{\alpha qb}
 \le1-\EE e^{-sM_Q(f)}
 \le s\left(b\int f+C\|f\|_\infty\right)+\PP(E_b).
\]
Divide by $s\log(1/s)=s\alpha b$, then let $b\to\infty$ and
$q\uparrow1$. Uniformly for bounded $\|f\|_\infty$,
\[
 \frac{-\log\EE e^{-sM_Q(f)}}{s\log(1/s)}
 \longrightarrow\frac1\alpha\int f.
\]
The uniform cutoff argument at the end of Appendix~F applies this
limit to $f=e^{\alpha d}w\ind_{Q_\delta}$ and proves
\eqref{chII:eq:critical-laplace}.
\end{proof}

Return to the square decomposition of Proposition~\ref{chII:prop:block-fields}.
Let $H_{N,K,i}(z)$ be the continuous piecewise bilinear (tensor-hat) interpolation of
the centered coarse field in unit coordinates of square $i$.
Thus $H_{N,K,i}(x/m)=H(im+x)$ for $x\in\{0,\ldots,m\}^2$,
with $m=\lfloor N/K\rfloor$ and torus coordinates understood modulo $N$.
On an
extension independent of this field let $(\mathcal M_i)_i$ be
independent copies of $\mathcal M$. Set
\begin{equation}\label{chII:eq:scalar-block-mass}
 B_{N,K,\delta}=\frac{K^{-4}}{p_{\rm per}}
       \sum_i\int_{Q_\delta}e^{\alpha H_{N,K,i}(z)}\,\mathcal M_i(dz).
\end{equation}
All deterministic variance corrections are the factors displayed in
this formula.

\begin{proposition}\label{chII:prop:scalar-gluing}
For probability laws $\mu,\nu$ on $\mathbb R$, let
\[
 d_{\rm BL}(\mu,\nu)
 =\sup_{\substack{\|f\|_\infty\le1\\\operatorname{Lip}(f)\le1}}
       |\mu(f)-\nu(f)|.
\]
The supremum is over real functions on $\mathbb R$, and
$\operatorname{Lip}(f)=\sup_{x\ne y}|f(x)-f(y)|/|x-y|$.
For every $e>0$ there are $\delta_0>0$ and $K_0$
such that, whenever $0<\delta\le\delta_0$, $\delta<1/2$ and
$K\ge K_0$, all sufficiently large $N$ satisfy
\begin{equation}\label{chII:eq:scalar-gluing}
 d_{\rm BL}(\mathcal L(B_{N,K,\delta}),\mathcal L(Z))<e.
\end{equation}
The threshold for $N$ may depend on $\delta,K$. The corresponding
family of masses is tight in this order of limits.
\end{proposition}
\begin{proof}
Put
\[
 D_{K,i}=(i+Q)/K,\qquad
 U_{K,\delta}=\bigcup_{i\in\{0,\ldots,K-1\}^2}(i+Q_\delta)/K.
\]
In the continuum decomposition of Appendix~G.1, $H_{K,i}$ is the
unit-coordinate coarse field and
$\overline\eta_K=\sum_i\int_{D_{K,i}}\eta_{K,i}$ is the fine mean.
For fixed $K,\delta$, that appendix proves
\[
 B_{N,K,\delta}\Rightarrow B_{K,\delta}^{\rm cont},\qquad
 B_{K,\delta}^{\rm cont}
       =e^{\alpha\overline\eta_K}\mathcal Z(U_{K,\delta})\quad\text{a.s.}
\]
In \eqref{chII:eq:scalar-block-mass}, $K^{-4}=K^{-2}K^{-2}$
is the product of the area and Dirichlet finite-part factors;
$p_{\rm per}^{-1}$ is the periodic normalization.
Appendix~G.1 verifies these factors.
In particular,
\[
 \bigl|B_{K,\delta}^{\rm cont}-Z\bigr|
 \le\bigl|e^{\alpha\overline\eta_K}-1\bigr|Z
        +\mathcal Z(U_{K,\delta}^c),\qquad
 \operatorname{Var}(\overline\eta_K)\le CK^{-2}.
\]
For $M,\eta>0$, the collar estimate \eqref{chII:eq:collar-mass} gives
\[
 \PP\{\mathcal Z(U_{K,\delta}^c)>\eta\}
 \le\PP(Z>M)+CM\delta/\eta,
\]
uniformly in $K$. Choose $M$, then a small collar width $\delta_0$,
and then $K_0$. The preceding bounds and $Z<\infty$ make
\[
 d_{\rm BL}\bigl(\mathcal L(B_{K,\delta}^{\rm cont}),\mathcal L(Z)\bigr)
 \le\EE\min\{2,|B_{K,\delta}^{\rm cont}-Z|\}<e/2
\]
uniformly for $0<\delta\le\delta_0$ and $K\ge K_0$.
The fixed-$K,\delta$ discrete limit then proves
\eqref{chII:eq:scalar-gluing} for all sufficiently large $N$.
Bounded continuous cutoffs give tightness in the same regime.
\end{proof}

\subsection*{2.6. The averaged Gaussian shell source}

We now give the scalar Gaussian limit used in Chapter~V. Its
relation to local times will be established by coarea in Chapter~III.
All objects in this section are defined directly from Gaussian fields.

Write $d_{\rm gr}=d_1$ for the graph distance of Chapter~I, the minimum
number of nearest-neighbor edges in a path on $V_N$.
For $U\subset V_N$, put $\rho_N=\tfrac12e^{L^{1/6}}$ and
$R_N=D_N^{\rm sep}$. Delete every $x\in U$ having some $y\in U$
with $\rho_N<d_{\rm gr}(x,y)<R_N$.
On the remaining vertices join pairs with $d_{\rm gr}(x,y)\le\rho_N$.
Let $\mathcal P_N(U)$ be the family of components of this graph.
Choose the lexicographically first vertex of each component in
$\{0,\ldots,N-1\}^2$; write $\operatorname{Rep}_N(U)$ for the set of these
representatives and $D_a(U)$ for the component represented by $a$.
Thus
\[
 \mathcal P_N(U)=\{D_a(U):a\in\operatorname{Rep}_N(U)\},\qquad
 n_N(U)=|\mathcal P_N(U)|=|\operatorname{Rep}_N(U)|,
\]
with $n_N(\varnothing)=0$.
Deletion precedes taking components. For $0\le\lambda\le1$ and a
nonnegative profile $W\in[0,\infty)^{V_N}$, define
\[
 \Psi_{N,\lambda}(W)
 =(1-\lambda)^{n_N(\{x:W(x)\le1/4\})},
\]
with the convention $0^0=1$. Squares of fields are taken pointwise.
Using the centering $m_N$ from \eqref{chII:eq:extremal-center} with
$M=N$, for $h\in\mathbb R$ and $\phi\in\mathbb R^{V_N}$ we have
\begin{equation}\label{chII:eq:canonical-band}
 \Psi_{N,\lambda}((\phi-m_N-h)^2)
 =(1-\lambda)^{n_N(\{x:|\phi(x)-m_N-h|\le1/2\})}.
\end{equation}

For nonempty finite sets $W_j\uparrow\mathbb Z^2$ and heights
$1\le T_j\uparrow\infty$,
monotone convergence in \eqref{chII:eq:contact-leading-coefficient} gives
$J_{W_j,T_j}\uparrow J_\infty$, independently of the exhausting
sequences, where
\begin{equation}\label{chII:eq:Gamma-definition}
 J_\infty=\int_{-1/2}^\infty e^{-\alpha u}
       \nu\{\exists x\in\mathbb Z^2:|u-\omega_x|\le1/2\}\,du,
 \qquad \Gamma=p_{\rm per}c_*J_\infty.
\end{equation}
Since $\omega_0=0$, the integrand equals $e^{-\alpha u}$ on
$[-1/2,1/2]$ and is bounded above by $e^{-\alpha u}$ throughout. Hence
\[
 0<\frac{e^{\alpha/2}-e^{-\alpha/2}}{\alpha}
 \le J_\infty\le\frac{e^{\alpha/2}}{\alpha}<\infty.
\]
Thus $\Gamma>0$. Later we define the deterministic local-time
comparison coefficient $\beta>0$ and set $\kappa=\beta\Gamma$.

\begin{proposition}\label{chII:prop:height-average}
For $0<\lambda<1$ and $\chi\in C_c(\mathbb R)$,
\begin{equation}\label{chII:eq:height-average}
 \int\chi(h)\EE \Psi_{N,\lambda}((\phi_N-m_N-h)^2)\,dh
 \longrightarrow\int\chi(h)\EE e^{-\lambda\Gamma e^{-\alpha h}Z}\,dh.
\end{equation}
\end{proposition}
\begin{proof}
Equations~\eqref{chII:eq:canonical-block-error} and
\eqref{chII:eq:block-mass-integral} of
\hyperlink{proof.appendix.G}{Appendix~G} give, for every $H_0>0$ and $e>0$
and with the parameter choices specified there,
\[
 \begin{split}
 \int_{-H_0}^{H_0}\Big|&\EE\Psi_{N,\lambda}((\phi_N+\overline\eta-m_N-h)^2)\\
   &-\EE e^{-\lambda\Gamma e^{-\alpha h}B_{N,K,2\delta}}\Big|\,dh
   <(2H_0+1)e.
 \end{split}
\]
The intermediate term is $\EE\prod_i(1-\lambda q_i(h,H))$.
The probabilities $q_i(h,H)$ use collar width $\delta$, whereas
the mass uses $2\delta$. Proposition~\ref{chII:prop:scalar-gluing},
applied with collar width $2\delta$, then gives the limit with $Z$.

To remove the dependent mean correction, use the following elementary
bound. If $0\le F(t,\omega)\le1$ and $R$ is any random variable on
the same space, then Fubini and translation in $t$ give
\begin{equation}\label{chII:eq:dependent-translation}
 \left|\int\chi(t)\EE[F(t-R,\omega)-F(t,\omega)]dt\right|
 \le\EE\|\chi(\cdot+R)-\chi\|_{L^1}.
\end{equation}
For $R=\overline\eta$ this tends to zero as $K\to\infty$, by
\eqref{chII:eq:fine-mean}, continuity of translations in $L^1$, and the
bound $2\|\chi\|_1$. The translation bound is pathwise in $R$ and therefore applies to an
arbitrary $R$ on the same probability space. It completes the averaged
limit.
\end{proof}

To specify the shell coordinate, recall $b_N=gN^2L$ and define the
central deterministic time
\begin{equation}\label{chII:eq:central-time}
 s_N=b_N\left(2L-\frac32\log L\right)
              +\frac{9g}{32}N^2(\log L)^2.
\end{equation}
Write $s_N(c)=s_N+cb_N$. In particular $2s_N/N^2=m_N^2$
exactly. Put $n=N^2$ and
\begin{equation}\label{chII:eq:shell-parameters}
 \begin{gathered}
 \mathcal{Q}_N=\frac1n\sum_x\phi_N(x)^2,\qquad
 B_N(\theta)=\frac{2(s_N+nv_N+\theta b_N)}n,\\
 \alpha_N(\theta,\phi)=\sqrt{\left(B_N(\theta)-\frac1n\sum_x\phi(x)^2\right)_+}.
 \end{gathered}
\end{equation}
For $a>0$ and $\sigma\in\{-1,1\}$ define the unnormalized branch
weight by
\begin{equation}\label{chII:eq:shell-weights}
 w_N^\sigma(\phi,a)=\frac{\sum_x(\phi(x)+\sigma a)_+}{n^2a}.
\end{equation}
The shell source, with $a=\alpha_N(\theta,\phi_N)$ inside the expectation, is
\begin{equation}\label{chII:eq:shell-source}
 S_N(\theta,\lambda)
 =n\EE\left[\ind_{\{\mathcal{Q}_N<B_N(\theta)\}}
        \sum_{\sigma=\pm1}w_N^\sigma(\phi_N,a)
                        \Psi_{N,\lambda}((\phi_N+\sigma a)^2)\right].
\end{equation}
The entire integrand is defined to be zero on
$\{\mathcal Q_N\ge B_N(\theta)\}$. The branch weights are evaluated
only on $\{\mathcal Q_N<B_N(\theta)\}$, where $a>0$.
In particular $S_N(\theta,\lambda)=0$ when
$s_N+nv_N+\theta b_N\le0$.

\begin{theorem}\label{chII:thm:shell-source}
For every $0<\lambda<1$ and $\chi\in C_c(\mathbb R)$,
\begin{equation}\label{chII:eq:shell-limit}
 \int\chi(\theta)S_N(\theta,\lambda)\,d\theta
 \longrightarrow
 \int\chi(\theta)\EE\exp\{-\lambda e^{-(\theta+1/2)}\Gamma Z\}\,d\theta.
\end{equation}
On each compact $\theta$ interval, $0\le S_N(\theta,\lambda)\le2$
for all sufficiently large $N$.
\end{theorem}
\begin{proof}
Equation~\eqref{chII:eq:shell-reference-error} proves that replacing the weighted two-branch source by
the unweighted positive-height band at amplitude $\alpha_N(\theta,\phi_N)$
has error tending to zero pointwise in $\theta$, with a uniform
compact bound. Define
\[
 a_N^0(\theta)=\sqrt{(B_N(\theta)-v_N)_+},\qquad
 R_N^E=\frac{\mathcal{Q}_N-v_N}{2gL}.
\]
There is the exact identity
$\alpha_N(\theta,\phi_N)=a_N^0(\theta-R_N^E)$.
Fourier summation gives $\EE|\mathcal{Q}_N-v_N|=O(1)$, hence
$\EE|R_N^E|=O(L^{-1})$. Apply \eqref{chII:eq:dependent-translation}
in the time variable to remove this random shift before changing
variables in the deterministic height. Uniformly on compact sets,
\[
 a_N^0(\theta)-m_N\longrightarrow\frac{\theta+1/2}{\alpha},
 \qquad \frac{d}{d\theta}a_N^0(\theta)\longrightarrow\frac1\alpha.
\]
Indeed $(a_N^0)^2=m_N^2+v_N+2\theta gL$, $v_N=gL+O(1)$,
and $m_N\sim2\sqrt gL$. The inverse change of variable sends
$\chi(\theta)d\theta$ to weights converging in $L^1$ to
$\alpha\chi(\alpha h-1/2)dh$. Apply
Proposition~\ref{chII:prop:height-average} with this limiting weight;
the $L^1$ error costs at most its norm because the band test lies
in $[0,1]$. Changing back to $\theta$ proves \eqref{chII:eq:shell-limit}.
The compact bound from \hyperlink{proof.appendix.H}{Appendix~H} permits dominated convergence
for the initial shell replacement.
\end{proof}

The output of this chapter is the same-field random measure
$\mathcal Z$, the relative contact estimates, and the scalar source
limit \eqref{chII:eq:shell-limit}, averaged against every
$\chi\in C_c(\mathbb R)$. The coarea identity \eqref{chIII:eq:coarea}
in Chapter~III relates this Gaussian source to local-time observables.

\clearpage
\endgroup

\begingroup
\mergedchaptersetup{III}{Local times and resampling}
\begin{center}
{\large\scshape Chapter III\par}\vspace{6pt}
{\LARGE Local times and resampling\par}
\end{center}
\renewcommand{\AppendixEquationPrefix}{A}
\vspace{6pt}

The finite-graph identity \eqref{chIII:eq:coarea} equates the image
measures of $(t+\|h\|_2^2/2,2\ell_t+h^2)$ and the Gaussian shell.
For resampling, the specified-endpoint law factors into selected
records and exterior coordinates; the observation window is imposed
by $J_I$ in \eqref{chIII:eq:window-functional}. Sections~3.4 and~3.6
derive probability and counting bounds from these two identities.

\subsection*{3.1. Canonical components and their activity}

Retain the rate-one walk on $V_N=(\ZZ/N\ZZ)^2$ and its local times
$\ell_t(x)=\int_0^t\ind_{\{X_u=x\}}\,du$. Write
$n=N^2$, $L=\log N$, $g=2/\pi$, and $b_N=gnL$.
The observation window is centered at
\begin{equation}\label{chIII:eq:center}
 s_N=b_N\left(2L-\frac32\log L\right)
             +\frac{9g}{32}n(\log L)^2,
 \qquad s_N(c)=s_N+cb_N.
\end{equation}
These are the shell coordinates fixed in Chapter~II. The root $r$ is
the starting vertex of the path and remains fixed as $t$ varies.

Use the graph distance $d_1$ of Chapter~I, and put
\[
 \rho_N=\frac12e^{L^{1/6}},\qquad R_N=Ne^{-2L^{1/8}}.
\]
For $U\subset V_N$ define its clean part by
\begin{equation}\label{chIII:eq:clean}
 U^\circ=\{x\in U:\text{there is no }y\in U
                  \text{ with }\rho_N<d_1(x,y)<R_N\}.
\end{equation}
Join distinct vertices of $U^\circ$ at distance at most $\rho_N$.
Choose the lexicographically first vertex of each component, using the
representatives $\{0,\ldots,N-1\}^2$ of $V_N$. Denote the resulting
set by $\operatorname{Rep}_N(U)$ and the component with representative $a$ by
$D_a(U)$. Thus deletion is performed first, and the graph components are then
taken on the remaining vertices.
In the notation of \eqref{chII:eq:canonical-band},
$n_N(U)=|\mathcal P_N(U)|=|\operatorname{Rep}_N(U)|$.

\begin{lemma}\label{chIII:lem:geometry}
Suppose $2\rho_N<R_N$. Every $D_a(U)$ has $d_1$-diameter at most
$\rho_N$. Distinct components are at $d_1$-distance at least $R_N$.
Every clean vertex belongs to exactly one of them.
\end{lemma}
\begin{proof}
Two clean vertices have distance at most $\rho_N$ or at least $R_N$.
Fix a vertex $x$ and follow an edge path in its component. If this path
first reaches a vertex at distance greater than $\rho_N$ from $x$,
that distance is at most $2\rho_N<R_N$, a contradiction. Starting
from either member of any pair gives the diameter assertion. Vertices
in different components cannot be at distance at most $\rho_N$,
since that would give an edge. The gap alternative therefore gives
their separation. The last assertion is the partition into graph
components. The scale hypothesis holds for all sufficiently large $N$.
\end{proof}

Let $\phi_N$ be the mean-zero Gaussian field of Chapter~II, with
covariance $G_N$, and let $v_N=G_N(r,r)$, independent of $r$.
Recall the two rooted fields from \eqref{chII:eq:two-root-fields}:
\begin{equation}\label{chIII:eq:rooted-fields}
 \begin{aligned}
 \phi_N^{(r)}(x)
     &=\phi_N(x)-\frac{G_N(x,r)}{v_N}\phi_N(r),\\
 h^r(x)&=\phi_N(x)-\phi_N(r).
 \end{aligned}
\end{equation}
The first has the conditional law of $\phi_N$ given $\phi_N(r)=0$;
it is independent of $\phi_N(r)$. The second is obtained by subtracting
the root value. In particular,
\begin{equation}\label{chIII:eq:field-connection}
 h^r(x)=\phi_N^{(r)}(x)
          +\left(\frac{G_N(x,r)}{v_N}-1\right)\phi_N(r).
\end{equation}
Both vanish at $r$, but their covariances differ. The activity uses
$\phi_N^{(r)}$; the coarea identity uses $h^r$.

For profiles $\ell:V_N\to[0,\infty)$ and $z:V_N\to\mathbb R$,
a set $D\subset V_N$ and a threshold $v>0$, define
\begin{equation}\label{chIII:eq:contact-union}
 \mathcal I_D(\ell,z;v)
 =\bigcup_{x\in D}
 \left[-z(x)-\sqrt{(v-2\ell(x))_+},
       -z(x)+\sqrt{(v-2\ell(x))_+}\right].
\end{equation}
Lebesgue length is denoted by $|\mathcal I_D|$. Set
\begin{equation}\label{chIII:eq:activity}
 a_D^G(\ell;r,v)
 =\frac1{\sqrt{2\pi g}}\,
       \EE\left|\mathcal I_D(\ell,\phi_N^{(r)};v)\right|.
\end{equation}
The expectation is over $\phi_N^{(r)}$, with $\ell,D,r,v$
fixed. Thus $a_D^G$ is a deterministic function of these arguments.
If $2\ell(x)>v$, the corresponding interval is a singleton and has
zero Lebesgue length. Consequently
\begin{equation}\label{chIII:eq:contact-integral}
 |\mathcal I_D(\ell,z;v)|
 =\int_{\mathbb R}\ind_{\{\exists x\in D:
                       2\ell(x)+(z(x)+q)^2\le v\}}\,dq.
\end{equation}
The two scalar sets can differ only at finitely many such singletons,
which do not affect this Lebesgue integral.

Write $U_u(t)=\{x:\ell_t(x)\le u\}$ and $U_0(t)=\{x:\ell_t(x)=0\}$.
For physical time $t\ge0$, define the two scalar path functionals
\[
 A_N(t)=\frac1{\sqrt L}
       \sum_{D\in\mathcal P_N(U_1(t))}a_D^G(\ell_t;r,1/4),
 \qquad C_N(t)=\frac{|\mathcal P_N(U_0(t))|}{\sqrt L},
\]
where $r=X_0$ is the original root. Their centered-time versions are,
for $s_N(c)\ge0$,
\begin{align}
 \widehat A_N(c)&=A_N(s_N(c)),\label{chIII:eq:AN}\\
 \widehat C_N(c)&=C_N(s_N(c)).\label{chIII:eq:CN}
\end{align}
For fixed $c$, $s_N(c)>0$ for all sufficiently large $N$. At every
observation time define the zero set by $U_0(t)=\{x:\ell_t(x)=0\}$; at
deterministic positive times it agrees almost surely with the unvisited
set.

\begin{lemma}\label{chIII:lem:regularity}
The activities in \eqref{chIII:eq:activity} are finite and nonnegative, and
are Borel functions of the occupation profile. For $\ell\ge0$,
\begin{equation}\label{chIII:eq:activity-bound}
 0\le a_D^G(\ell;r,v)
       \le\frac{2\sqrt v}{\sqrt{2\pi g}}\,|D|.
\end{equation}
The corresponding activity sum and component count are jointly
measurable in the path and observation time. In particular $\widehat A_N(c)$
and $\widehat C_N(c)$ are random variables.
\end{lemma}
\begin{proof}
The length of a finite union of intervals is at most the sum of their
lengths. This proves \eqref{chIII:eq:activity-bound} and integrability in
the Gaussian variable. Its length is Borel in the endpoints, either
by inclusion--exclusion for intervals or by integrating its jointly
Borel indicator. Gaussian integration preserves this measurability.
Each $\ell_t(x)$ is jointly measurable; it is continuous in $t$ on
every path. Membership in $U_u(t)$ is therefore measurable. There are
only finitely many subsets of $V_N$, and the component construction
is a deterministic map on this finite set. Summing over each possible
subset proves the assertion. The same argument, with the equality
$\ell_t(x)=0$, treats the zero set.
\end{proof}

\subsection*{3.2. Ordinary paths, bridges and integrated loops}

For $t>0$ write $\PP_t^{r,r}=\PP_r(\,\cdot\mid X_t=r)$, with tests
restricted to the path up to $t$. Since $p_t(r,r)\ge e^{-t}>0$, this
is a probability law. We keep the return probability outside the
bridge when time is integrated:
\begin{equation}\label{chIII:eq:loop-measure}
 \mathsf L_r(dt,d\omega)
 =\ind_{\{t>0,\,X_t(\omega)=r\}}\,dt\,\PP_r(d\omega)
 =\ind_{\{t>0\}}p_t(r,r)\,dt\,\PP_t^{r,r}(d\omega).
\end{equation}
This is a $\sigma$-finite measure, with finite mass on bounded time
intervals. In the first expression the infinite future is retained;
it can be integrated out whenever the test reads only the stopped path.

For $I=[l,u]\subset(0,\infty)$, $l<u$, and a nonnegative measurable
test $F(r,t,\omega)$, put
\begin{equation}\label{chIII:eq:window-law}
 \begin{split}
 \mathsf W_N(I;F)
   &=\frac1n\sum_{r\in V_N}\int_l^u
                  \EE_r[\ind_{\{X_t=r\}}F(r,t,\omega)]\,dt,\\
 \mathsf Q_N^I(F)&=\frac{\mathsf W_N(I;F)}{\mathsf W_N(I;1)}.
 \end{split}
\end{equation}
The same notation applies to bounded signed tests by linearity.
Translation invariance gives
\begin{equation}\label{chIII:eq:normalizer}
 \mathsf W_N(I;1)=\int_l^u p_t(0,0)\,dt.
\end{equation}
Under $\mathsf Q_N^I$ the root is uniform. Time has density
$p_t(0,0)/\mathsf W_N(I;1)$; given the root and time, the path is a bridge.

The heat-kernel estimate of Chapter~I implies the following useful
normalization rule. If
\[
 \eta_N(l)=\sup_{t\ge l}\sup_{x,y}|np_t(x,y)-1|,
\]
then
\begin{equation}\label{chIII:eq:window-normalization}
 (1-\eta_N(l))\frac{|I|}{n}
 \le \mathsf W_N(I;1)\le(1+\eta_N(l))\frac{|I|}{n}.
\end{equation}
A fixed central window is $I=[s_N(c_0),s_N(c_1)]$ with fixed
$c_0<c_1$. In particular $\mathsf W_N(I;1)=(1+o(1))|I|/n$ uniformly
when $c_0,c_1$ range over a bounded interval. For $0\le d<l$, put
$I^d=[l-d,u+d]$. If $\eta_N(l-d)=o(1)$ and $d=o(|I|)$,
then
\begin{equation}\label{chIII:eq:normalizer-ratio}
 \frac{\mathsf W_N(I^d;1)}{\mathsf W_N(I;1)}=1+o(1).
\end{equation}
Equation~\eqref{chIII:eq:window-normalization} controls the normalization
under window enlargement.

\begin{lemma}\label{chIII:lem:bridge-prefix}
Let $0<\delta<t$. On the path up to time $t-\delta$, the density of
$\PP_t^{r,r}$ with respect to $\PP_r$ is
\begin{equation}\label{chIII:eq:prefix-density}
 \frac{p_\delta(X_{t-\delta},r)}{p_t(r,r)}.
\end{equation}
If $\eta_N(\delta)<1$, this density lies between
$(1-\eta_N(\delta))/(1+\eta_N(\delta))$ and its reciprocal.
Hence expectations of tests bounded in absolute value by one differ
by at most $2\eta_N(\delta)/(1-\eta_N(\delta))$.
\end{lemma}
\begin{proof}
For a bounded test $F$ of the prefix, the Markov property gives
$\EE_r[F\ind_{\{X_t=r\}}]
=\EE_r[Fp_\delta(X_{t-\delta},r)]$. Division by $p_t(r,r)$
proves the formula. Apply the definition of $\eta_N(\delta)$ to
the numerator and denominator, and integrate the bound on the density
minus one.
\end{proof}

The comparison applies to the prefix up to $t-\delta$; the terminal
bridge constraint $X_t=r$ remains at time $t$.

There is also an energy-window measure. Write $\mathcal L(Y)$
for the law of a random element $Y$, and put $\gamma_r=\mathcal L(h^r)$.
For a field coordinate $h$ with law $\gamma_r$, independent of the
path under $\mathsf L_r$, define
\begin{equation}\label{chIII:eq:energy-time}
 \mathcal E_r(t,h)=t+\frac12\|h\|_2^2,
 \qquad \|h\|_2^2=\sum_xh(x)^2.
\end{equation}
The root-averaged product measure on $(r,t,\omega,h)$ is
\[
 \mathsf M_N=\frac1n\sum_{r\in V_N}
       \delta_r\otimes\mathsf L_r\otimes\gamma_r,
\]
where $\delta_r$ is unit mass at $r$. For a bounded Borel set
$J\subset(0,\infty)$ of positive $\mathsf M_N$-mass under the energy map,
define
\[
 \mathsf Q_N^{J,\mathrm{en}}(F)
 =\frac{\displaystyle\int
     \ind_{\{\mathcal E_r(t,h)\in J\}}F(r,t,\omega,h)\,d\mathsf M_N}
    {\displaystyle\int\ind_{\{\mathcal E_r(t,h)\in J\}}\,d\mathsf M_N}.
\]
The denominator is finite because the restriction forces
$0<t\le\sup J$. The displayed probability law is the normalized joint
restriction of $\mathsf M_N$ to $\{\mathcal E_r(t,h)\in J\}$;
path-field tests are integrated against this joint law.
Here the square norm in \eqref{chIII:eq:energy-time} is the
counting-measure norm.

\subsection*{3.3. The exact finite-graph coarea identity}

Fix $r$ for this subsection and abbreviate $h=h^r$ and $\phi=\phi_N$.
For a field $f$, the notation $f^2$ means the profile $x\mapsto f(x)^2$.
Write $\EE_\phi$ for expectation with respect to the law of $\phi_N$,
and $\EE_h$ for integration against $\gamma_r$. Let
\[
 a_s(\phi)=\sqrt{\frac{2s-\|\phi\|_2^2}{n}}
          \quad\text{when }\|\phi\|_2^2<2s.
\]
On this set define, for $\sigma\in\{-1,1\}$,
\begin{equation}\label{chIII:eq:shell-density}
 b_{r,s}^{\sigma}(\phi)
     =\frac{(\phi(r)+\sigma a_s(\phi))_+}{n a_s(\phi)},
 \qquad W_s^{\sigma}(\phi)(x)
     =(\phi(x)+\sigma a_s(\phi))^2.
\end{equation}
Outside this strict set, including its boundary, set
$a_s=b_{r,s}^\sigma=0$ and $W_s^\sigma\equiv0$.
Thus the amplitude of Chapter~II satisfies
\[
 \alpha_N(\theta,\phi)
 =a_{s_N+nv_N+\theta b_N}(\phi),
\]
whenever the shell level $s_N+nv_N+\theta b_N$ is positive.
Define a measure $\mathsf S_r$ on
$(0,\infty)\times[0,\infty)^{V_N}$ by
\begin{equation}\label{chIII:eq:shell-measure}
 \int F\,d\mathsf S_r
  =\int_0^\infty\EE_\phi\sum_{\sigma=\pm1}
       \ind_{\{\|\phi\|_2^2<2s\}}
       b_{r,s}^{\sigma}(\phi) F(s,W_s^{\sigma}(\phi))\,ds.
\end{equation}
The branch variable $\sigma\in\{-1,1\}$ carries counting measure.

\begin{theorem}\label{chIII:thm:coarea}
For every nonnegative Borel function $F$ of $(s,W)$,
\begin{equation}\label{chIII:eq:coarea}
 \int_0^\infty\EE_r\!\left[
   \ind_{\{X_t=r\}}\EE_h
    F\left(t+\tfrac12\|h\|_2^2,\,2\ell_t+h^2\right)\right]dt
       =\int F\,d\mathsf S_r.
\end{equation}
All operations are valid for extended nonnegative integrals;
equivalently, \eqref{chIII:eq:coarea} is an equality of image measures.
\end{theorem}
\begin{proof}
Let $\tau_r(\lambda)=\inf\{t:\ell_t(r)>\lambda\}$. Occupation of
the root gives the pathwise change of variable
\begin{equation}\label{chIII:eq:occupation-substitution}
 \int_0^\infty\ind_{\{X_t=r\}}H(t)\,dt
       =\int_0^\infty H(\tau_r(\lambda))\,d\lambda
\end{equation}
for every nonnegative Borel $H$, almost surely. Indeed, the local-time
clock has slope one on each holding interval at $r$ and slope zero
between such intervals. Its inverse translates the interiors of these
intervals, and their endpoints form a countable null set. Recurrence
on the finite graph makes the range of this clock all of $[0,\infty)$.

The finite-graph square-field identity, proved in \hyperlink{proof.appendix.I}{Appendix~I} by
Laplace transforms, is
\begin{equation}\label{chIII:eq:RK}
 2\ell_{\tau_r(\lambda)}+h^2
          \ \overset{\mathrm{law}}=\ (h+\sqrt{2\lambda})^2,
 \qquad \lambda\ge0,
\end{equation}
where the path and $h$ on the left are independent. Moreover
$\sum_x\ell_t(x)=t$. The first coordinate of the image in
\eqref{chIII:eq:coarea} is therefore half the sum of its second coordinate.
Applying \eqref{chIII:eq:occupation-substitution} and \eqref{chIII:eq:RK} changes
the left-hand side of \eqref{chIII:eq:coarea} into
\begin{equation}\label{chIII:eq:coarea-intermediate}
 \int_0^\infty\EE_\phi
 F\left(\tfrac12\|\phi-\phi(r)+\sqrt{2\lambda}\|_2^2,
                (\phi-\phi(r)+\sqrt{2\lambda})^2\right)d\lambda.
\end{equation}

For fixed $\phi$ put $q=\sqrt{2\lambda}-\phi(r)$. Then
$q> -\phi(r)$ and $d\lambda=(q+\phi(r))\,dq$.
Since $\sum_x\phi(x)=0$,
\[
 s=\tfrac12\|\phi+q\|_2^2
       =\tfrac12\|\phi\|_2^2+\tfrac n2 q^2.
\]
Split the $q$-integral at zero. On each nonempty branch,
$q=\sigma a_s(\phi)$ and $|dq/ds|=(n a_s(\phi))^{-1}$.
The restriction $q> -\phi(r)$ replaces its numerator by
$(\phi(r)+\sigma a_s(\phi))_+$. These substitutions give exactly
\eqref{chIII:eq:shell-measure}. The points $q=0$ and $q=-\phi(r)$ do not
contribute to the scalar integral. Tonelli justifies every interchange.
\end{proof}

For root-independent image tests, averaging \eqref{chIII:eq:coarea} over
$r$ replaces $b_{r,s}^{\sigma}$ by
\begin{equation}\label{chIII:eq:averaged-density}
 \overline b_s^\sigma(\phi)
       =\frac{\sum_r(\phi(r)+\sigma a_s(\phi))_+}
                    {n^2a_s(\phi)}.
\end{equation}
This is the branch weight of Chapter~II. The factor $n^{-2}$
contains both the Jacobian and the uniform root average.

On central shell intervals define the time densities and conditional
field--branch laws by
\[
 \begin{aligned}
 \varrho_{N,r}(s)&=\EE_\phi\sum_{\sigma=\pm1}b_{r,s}^\sigma(\phi),
 &\nu_{N,s}^r(d\phi,\{\sigma\})
   &=\frac{b_{r,s}^\sigma(\phi)}{\varrho_{N,r}(s)}
                          \,\mathcal L(\phi_N)(d\phi),\\
 \varrho_N(s)&=\EE_\phi\sum_{\sigma=\pm1}\overline b_s^\sigma(\phi),
 &\nu_{N,s}(d\phi,\{\sigma\})
   &=\frac{\overline b_s^\sigma(\phi)}{\varrho_N(s)}
                          \,\mathcal L(\phi_N)(d\phi).
 \end{aligned}
\]
Translation invariance and \eqref{chII:eq:shell-normalizer} give
\[
 \varrho_{N,r}(s)=\varrho_N(s),\qquad
 n\varrho_N(s)=1+o(1),\qquad
 \frac1n\sum_r\nu_{N,s}^r=\nu_{N,s},
\]
uniformly when $s=s_N+nv_N+\theta b_N$ and $\theta$ stays bounded.
In particular the normalizers are finite and positive there.
These probability laws live on $(\phi,\sigma)$; their square-profile
events are written $\nu_{N,s}^r(W_s^\sigma(\phi)\in E)$.
The corresponding joint root law is
$n^{-1}\sum_r\delta_r\otimes\nu_{N,s}^r$, which records the dependence of
the field and branch on the root.

The coarea identity applies to Borel tests of the image $(s,W)$.
Section~3.4 combines the probe lower bound
\eqref{chIII:eq:probe-probability} and image majorant
\eqref{chIII:eq:probe-majorant} to obtain
\eqref{chIII:eq:probe-transfer}.

\subsection*{3.4. Gaussian probes and counting on the image space}

Fix $r$ and write $h=h^r$, with law $\gamma_r$. For $a=(t,\omega)$ write
\[
 T_r(a,h)=\left(t+\tfrac12\|h\|_2^2,\,2\ell_t+h^2\right).
\]
The following form allows the time window, the selected objects and
their weights to depend on the source path.

\begin{proposition}\label{chIII:prop:counted-probes}
Let $\mathcal J$ be finite or countable. For $i\in\mathcal J$, let
$w_i(a)\ge0$ be measurable, let $p_i\ge0$ be finite constants, and
let $B_i(a,h)\ge0$ be measurable probes. Suppose that, whenever
$w_i(a)>0$,
\begin{equation}\label{chIII:eq:probe-probability}
 \EE_h B_i(a,h)\ge p_i.
\end{equation}
If a nonnegative Borel image test $F$ satisfies
\begin{equation}\label{chIII:eq:probe-majorant}
 \sum_{i\in\mathcal J}w_i(a)B_i(a,h)\le F(T_r(a,h))
\end{equation}
almost everywhere under $\mathsf L_r\otimes\gamma_r$, then
\begin{equation}\label{chIII:eq:probe-transfer}
 \int\sum_{i\in\mathcal J}p_iw_i(a)\,\mathsf L_r(da)
                  \le\int F\,d\mathsf S_r.
\end{equation}
Here $0\cdot\infty=0$ when extended nonnegative weights are used.
\end{proposition}
\begin{proof}
Multiply \eqref{chIII:eq:probe-probability} by $w_i(a)$, integrate over
$a$ and sum. Tonelli gives
\[
 \int\sum_i p_iw_i\,d\mathsf L_r
 \le\int\EE_h\sum_iw_i(a)B_i(a,h)\,\mathsf L_r(da)
 \le\int\EE_hF(T_r(a,h))\,\mathsf L_r(da).
\]
The last integral equals the right-hand side of
\eqref{chIII:eq:probe-transfer} by Theorem~\ref{chIII:thm:coarea}.
For countably many labels the same calculation follows by monotone
convergence from finite partial sums.
\end{proof}

For indicator weights and probes, $w_i(a)B_i(a,h)=1$ specifies a
source label and its Gaussian event. An energy restriction is included
in $B_i$ itself. The hypothesis is then
$\EE_h B_i(a,h)\ge p_i$ for every $a$ with $w_i(a)=1$;
the sum in \eqref{chIII:eq:probe-majorant} is bounded pointwise by
$F(T_r(a,h))$.

\begin{corollary}\label{chIII:cor:truncated-probes}
Let $I=[l,u]\subset(0,\infty)$ and $H\ge0$. Suppose the source
weights vanish unless $t\in I$, and the probes vanish unless
$\|h\|_2^2\le H$. If
\[
 \sum_iw_i(a)B_i(a,h)\le K(2\ell_t+h^2)
\]
for a nonnegative Borel $K$, then the conclusion of
Proposition~\ref{chIII:prop:counted-probes} holds with
\begin{equation}\label{chIII:eq:truncated-image}
 F(s,W)=\ind_{\{l\le s\le u+H/2\}}K(W).
\end{equation}
\end{corollary}
\begin{proof}
Whenever the left-hand side of the majorization is nonzero,
$t\in[l,u]$ and $0\le\|h\|_2^2/2\le H/2$. Hence the first
image coordinate belongs to $[l,u+H/2]$. Outside this event the
left-hand side is zero. Apply Proposition~\ref{chIII:prop:counted-probes}.
\end{proof}

For indicator weights and probes, an injective correspondence from
the labels with $w_i(a)B_i(a,h)=1$ to image components implies
\[
 \sum_i w_i(a)B_i(a,h)\le K(2\ell_t+h^2),
\]
where $K$ counts those components. A correspondence with bounded
multiplicity requires the corresponding factor in this inequality.
The required geometric bounds are proved in Chapter~IV; the
measure-theoretic argument uses exactly
\eqref{chIII:eq:probe-majorant}.

\subsection*{3.5. Specified-endpoint records and their durations}

Let $D$ be a nonempty proper subset of $V_N$. Write
$\tau_D=\inf\{t:X_t\notin D\}$ and
\[
 h_{D,b}(a)=\PP_a(X_{\tau_D}=b),\qquad a\in D,\quad b\notin D.
\]
Write $\mathcal L_a$ for law under $\PP_a$.
When $h_{D,b}(a)>0$, define the local kernel
\begin{equation}\label{chIII:eq:local-kernel}
 \mathsf B_D^{a,b}
     =\mathcal L_a\bigl((X_s)_{0\le s\le\tau_D},\tau_D
                            \mid X_{\tau_D}=b\bigr).
\end{equation}
Thus $\mathsf B_D^{a,b}$ is the law of the stopped path together with its duration,
conditioned on the specified exit vertex. Recall the natural filtration
$\mathcal F_s=\sigma(X_u:0\le u\le s)$. On $\{s<\tau_D\}$, the density relative
to the unconditioned stopped path is $h_{D,b}(X_s)/h_{D,b}(a)$.
Indeed, conditional expectation of the specified exit event given
$\mathcal F_s$ is $h_{D,b}(X_s)$ on this set.

For successive records, an \emph{endpoint word} is the ordered list
$w=((a_i,b_i))_{i=1}^m$ of their entrance--exit pairs. The domain
$D_i$ and entrance $a_i$ are determined by the initial vertex and
the preceding exits; for consecutive records $a_{i+1}=b_i$.
The word is feasible if $h_{D_i,b_i}(a_i)>0$ for every $i$. For a
fixed domain family $(D_i)_{i=1}^m$ and feasible word $w$, write
\begin{equation}\label{chIII:eq:word-product-kernel}
 P_{m,w}:=\bigotimes_{i=1}^m\mathsf B_{D_i}^{a_i,b_i}.
\end{equation}
When the domains and endpoint word are fixed by the surrounding
stratum, we abbreviate this law to $P_m$. For $m=0$ the product kernel
is the point mass at the empty record family, with total occupation and
total duration equal to zero. Its ordinary endpoint-word probability is
$p_w=\prod_{i=1}^m h_{D_i,b_i}(a_i)$, as proved in
\eqref{chIII:eq:word-factor}; the empty product has value one.

\begin{proposition}\label{chIII:prop:record-product}
Consider finitely many successive first-exit records. Their domains
may be deterministic functions of the preceding exit vertices. Fix
an endpoint word of positive probability. Conditional on this word, the records are independent with laws
\eqref{chIII:eq:local-kernel}. The future path from the last exit vertex
has the ordinary walk law started at that vertex and is independent of
the records. Concatenation in the original temporal order gives the
conditional path law and adds the record durations in that order.
\end{proposition}
\begin{proof}
At the first exit, the strong Markov property separates the stopped
record from the future path started at its endpoint. On fixing that
endpoint, the first factor is \eqref{chIII:eq:local-kernel}. The domain
for the next record is now fixed. Repeating the same argument gives the product law and an independent
future path started at the final exit vertex. The endpoint probabilities
multiply to the mass of the word and disappear upon conditioning.
The measurable extraction and reconstruction maps are inverse on their
full-measure record space. \hyperlink{proof.appendix.J}{Appendix~J} gives the coordinate
factorization \eqref{chIII:eq:word-factor} and the loop-time
decomposition \eqref{chIII:eq:actual-slice}, with its censored terminal segment.
\end{proof}

\begin{proposition}\label{chIII:prop:duration}
Let $v_N=G_N(0,0)$ and $k=|D|$. For $d\ge0$,
\begin{equation}\label{chIII:eq:unconditioned-duration}
 \sup_{a\in D}\PP_a(\tau_D>d)
           \le k\exp\{-d/(4v_Nk)\}.
\end{equation}
Let $\PP=\bigotimes_{i=1}^m\mathsf B_D^{a_i,b_i}$, where $m\ge1$,
and let $V_i$ be the duration of record $i$. Its local-time profile is
$\ell^{(i)}(x)=\int_0^{V_i}\ind_{\{X_s^{(i)}=x\}}\,ds$,
where $X^{(i)}$ is that record. Suppose
\[
 h_{D,b_i}(a_i)\ge h_*>0,\qquad
 \mathsf B_D^{a_i,b_i}(\text{avoids }K)\ge p_*>0
\]
for a nonempty $K\subset D$. Let $S$ be any measurable event
containing the event that all $m$ records avoid $K$. Then
\begin{equation}\label{chIII:eq:conditional-duration}
 \PP\bigl(\max_iV_i>d\mid S\bigr)
 \le\frac{mk}{h_*p_*^m}\exp\{-d/(4v_Nk)\}.
\end{equation}
In particular this applies to
$S=\{\min_{x\in K}\sum_i\ell^{(i)}(x)\le u_0\}$ for $u_0\ge0$.
\end{proposition}
\begin{proof}
The killed spectral bound \eqref{chIII:eq:domain-gap} in \hyperlink{proof.appendix.K}{Appendix~K}.1 proves
\eqref{chIII:eq:unconditioned-duration}, using $|G_N(x,y)|\le v_N$.
Conditioning on a specified endpoint multiplies this bound
by at most $h_*^{-1}$.
A union bound bounds the conditional numerator by
$mk h_*^{-1}\exp\{-d/(4v_Nk)\}$.
Independence and simultaneous avoidance give $\PP(S)\ge p_*^m$.
Division proves \eqref{chIII:eq:conditional-duration}. Avoidance makes
every local time in $K$ zero, which proves the last assertion.
\end{proof}

If, for every $i=1,\ldots,m$, there exists a simple path from $a_i$
to $b_i$ that stays in $D\setminus K$ until its last step and uses
at most $k+1$ jumps, then both lower bounds hold uniformly in $i$ with
$h_*=p_*=4^{-(k+1)}$. Consequently
\begin{equation}\label{chIII:eq:finite-duration}
 \PP\bigl(\max_iV_i>d\mid S\bigr)
 \le mk\exp\left\{-\frac d{4v_Nk}
                      +(m+1)(k+1)\log4\right\}.
\end{equation}
The factor $p_*^m$ in \eqref{chIII:eq:conditional-duration} accounts for
the conditioning on $S$. For the small domains used later, this gives
the following uniform bound.

\begin{corollary}\label{chIII:cor:small-duration}
Fix $C<\infty$. Under the feasible avoiding-path assumption just
stated, uniformly for
\[
 |D|\le e^{CL^{1/8}},\qquad m\le e^{L/8}+1,
 \qquad u_0\ge0,
\]
all sufficiently large $N$ satisfy
\begin{equation}\label{chIII:eq:small-duration}
 \PP\bigl(\max_iV_i>e^{L/2}\mid
          \min_{x\in K}\textstyle\sum_i\ell^{(i)}(x)\le u_0\bigr)
                    \le\exp\{-e^{0.49L}\}.
\end{equation}
\end{corollary}
\begin{proof}
Chapter~I gives $v_N=O(L)$. In \eqref{chIII:eq:finite-duration} the negative
term has magnitude $e^{L/2-o(L)}$, whereas
$(m+1)(k+1)\log4+\log(mk)=e^{L/8+o(L)}$ for $m\ge1$.
Their difference is at least $e^{0.49L}$ for all sufficiently large
$N$, uniformly in the displayed ranges.
\end{proof}

For larger disks it is more efficient to use their radius than their
cardinality. Choose a lift $\widetilde y\in\mathbb Z^2$ of $y\in V_N$,
and let $D$ be the image in $V_N$ of the Euclidean disk
\[
 B(\widetilde y,R)=\{z\in\mathbb Z^2:|z-\widetilde y|<R\}.
\]
If $R\ge1$ and $2(R+1)<N$, the walk up to exit has a planar
lift. In these lifted coordinates, write $y$ again for the center.
Optional stopping of
$|X_t-y|^2-t$ gives
$\sup_{a\in D}\EE_a\tau_D\le(R+1)^2$. Markov's inequality at
$2(R+1)^2$, followed by iteration at these deterministic times, gives
\begin{equation}\label{chIII:eq:disk-duration}
 \sup_{a\in D}\PP_a(\tau_D>d)
      \le 2\exp\left\{-\frac{(\log2)d}{2(R+1)^2}\right\}.
\end{equation}
The same conditioning argument replaces the right-hand side of
\eqref{chIII:eq:conditional-duration} by
\begin{equation}\label{chIII:eq:disk-conditional}
 \frac{2m}{h_*p_*^m}
       \exp\left\{-\frac{(\log2)d}{2(R+1)^2}\right\}.
\end{equation}
With the endpoint and avoidance constants appearing in
\eqref{chIII:eq:disk-conditional}, the same resampling argument applies
to both the cardinality-based and radius-based exit bounds. The
annular estimates used for the latter are supplied in Chapter~IV.

\subsection*{3.6. Exterior selection and one successful replacement}

\paragraph{Conditional-kernel convention.}
Whenever a conditional kernel appears in this section, it is used on a
stratum of positive conditional mass.

Fix a root $r$ and a positive endpoint word as in
Proposition~\ref{chIII:prop:record-product}. Split its records into a
selected family $x$ and the complementary information $y$, which
includes the unused future. The product record law factors as
$\mu(dx)\nu(dy)$. Write $V(x)\ge0$ for the total selected duration
and $E(y)\ge0$ for the total complementary completed duration.
For any bounded interval $I=[l,u]\subset\mathbb R$, $l<u$, and a measurable
$G:Y\times[0,\infty)\to[0,1]$, where $Y$ is the exterior record space, define
\begin{equation}\label{chIII:eq:window-functional}
 J_I(x,y)=\int_0^\infty
        \ind_{\{E(y)+V(x)+q\in I\}}G(y,q)\,dq.
\end{equation}
The function $G:Y\times[0,\infty)\to[0,1]$ can impose return to the root,
censor the terminal segment before its next entrance, and retain any
bounded test of the exterior record variable $Y$. Always
$0\le J_I\le |I|$.

For $d\ge0$, if $V(x),V(x')\le d$, then
\begin{equation}\label{chIII:eq:short-window}
 J_I(x,y)\le J_{I^d}(x',y),\qquad I^d=[l-d,u+d].
\end{equation}
Indeed the same residual time $q$ can be used on both sides and
$|V(x)-V(x')|\le d$. The bound is uniform in the exterior duration $E(y)\ge0$.

\begin{theorem}\label{chIII:thm:one-success}
Suppose $\mu,\nu$ are probability measures and $B$ is a measurable
bad event of the selected records. Let $d,\varrho,\varepsilon\ge0$. If
\begin{equation}\label{chIII:eq:success-hypothesis}
 \mu(B)\le\varrho,\qquad \mu(V>d)\le\varepsilon,
 \qquad \varrho+\varepsilon\le\tfrac12,
\end{equation}
then
\begin{equation}\label{chIII:eq:one-success}
 \int J_I\,d(\mu\otimes\nu)
 \le 2\int\ind_{B^c}(x)J_{I^d}(x,y)\,d(\mu\otimes\nu)
                                      +|I|\varepsilon.
\end{equation}
There is also a bad-event form. Retain the first two bounds in
\eqref{chIII:eq:success-hypothesis}, but replace
$\varrho+\varepsilon\le1/2$ by $\varepsilon\le1/2$. Then
\begin{equation}\label{chIII:eq:bad-transfer}
 \int\ind_B(x)J_I(x,y)\,d(\mu\otimes\nu)
 \le 2\varrho\int J_{I^d}\,d(\mu\otimes\nu)
                                      +|I|\varepsilon.
\end{equation}
Both bounds hold for a single record and for any finite selected
family, with $V$ its total duration.
\end{theorem}
\begin{proof}
Put $T=\{V\le d\}$. The first hypotheses give
$\mu(B^c\cap T)\ge1/2$. For fixed $y$ and $x\in T$, integrate
\eqref{chIII:eq:short-window} over $x'\in B^c\cap T$ to obtain
\[
 J_I(x,y)\le2\int_{B^c}J_{I^d}(x',y)\,\mu(dx').
\]
Integrate over $x\in T$ and then over $y$. On $T^c$ use
$J_I\le |I|$. This proves \eqref{chIII:eq:one-success}. For the second
bound integrate \eqref{chIII:eq:short-window} over $x'\in T$, whose
mass is at least $1/2$, and then integrate the resulting bound over
$x\in B\cap T$. Its mass is at most $\varrho$. The remaining
integral over $B\cap T^c$ is bounded by $|I|\varepsilon$.
\end{proof}

A \emph{certified random selection} with proposed selected index set
$C$ is a stratum whose selection event has the form
\begin{equation}\label{chIII:eq:selection-certificate}
 \Lambda_C(y)\cap\bigcap_{j\in C}\{x_j\in F_j(y)\},
\end{equation}
where $y$ contains the unselected records and the retained endpoint
data. The sets $F_j(y)$ are measurable, and each has positive
probability under its local kernel whenever the stratum has positive
conditional mass.

\begin{lemma}\label{chIII:lem:selection}
On the stratum \eqref{chIII:eq:selection-certificate}, conditionally on
$y$, the selected records are independent with laws
\begin{equation}\label{chIII:eq:selected-law}
 \bigotimes_{j\in C}\mathsf B_{D_j}^{a_j,b_j}(\,\cdot\mid F_j(y)).
\end{equation}
The exterior measure acquires the factor
$\ind_{\Lambda_C}(y)\prod_{j\in C}\mathsf B_{D_j}^{a_j,b_j}(F_j(y))$.
The conclusions of Theorem~\ref{chIII:thm:one-success} can therefore be
applied conditionally on this exterior whenever their bounds hold
uniformly, or integrated with their exterior-dependent errors.
\end{lemma}
\begin{proof}
Multiply the product measure by the indicator in
\eqref{chIII:eq:selection-certificate}. For fixed $y$, factor each selected
restriction into its mass and its normalized conditional kernel.
The product of the masses is exactly the factor stated in the lemma.
Tonelli proves the factorization against every nonnegative test, hence
the conditional-law assertion. The proof of
Theorem~\ref{chIII:thm:one-success} works separately for each $y$; integrate
it against this exterior measure.
\end{proof}

We also use the bounded-weight version of \eqref{chIII:eq:bad-transfer}.
If $0\le f\le1$, $\mu(f)\le\varrho$, and
$\mu(V>d)\le\varepsilon\le1/2$, the same short-record integration gives
\[
 \int f(x)J_I(x,y)\,d(\mu\otimes\nu)
 \le2\varrho\int J_{I^d}\,d(\mu\otimes\nu)+|I|\varepsilon.
\]
The bounded-weight inequality applies on each exterior fibre of
\eqref{chIII:eq:selection-certificate}, with a local probability as its
measurable weight $f$.

An exterior-measurable choice is the case in which each $F_j$ is the
whole record space. Any selected-record flag used by the selection rule
is included in the corresponding conditioning event $F_j(y)$, so the
selected law is conditioned on that flag.

The applications in Chapter~VI use several disjoint hosts at once.
We record the corresponding extraction in the same coordinates.

\begin{lemma}[Multi-host record extraction]\label{chIII:lem:multihost}
Fix the root and a finite ordered family of hosts
$(I_j,O_j)_{j\in\mathcal J}$ with $I_j\subset O_j\subsetneq V_N$.
Assume the exit domains together with the one-step boundaries required
by the construction are pairwise disjoint within the chosen color and
that the root lies outside every $O_j$. Starting after the preceding
host exit, let the walk first hit $\bigcup_j I_j$; the hit location
identifies a unique host $j$, and the associated record ends at the
first subsequent exit from $O_j$. Repeating this procedure extracts
successive host records on the infinite ordinary path in chronological
order. Their time intervals are disjoint. At an observation time, only
completed records are included; an unfinished terminal segment is
assigned to the exterior coordinate.

Fix a deterministic finite number of these records and condition on
their chronological host labels and entrance--exit endpoint word.
Before imposing an observation-time or closing restriction, their
local kernels form the product $\bigotimes_i\mathsf B_{O_{j_i}}^{a_i,b_i}$.
For a selected subfamily $\xi$, the exterior
coordinate $y$ retains the host labels, endpoint words, chronological
slots, the gap paths and their durations, all unselected records, and
the terminal censored future. These data determine the original order;
if $V(\xi)$ is the selected total duration and $E(y)$ the completed
exterior duration, reconstruction has total completed duration
$E(y)+V(\xi)$. Counts and any eligibility conditions depending only on
host labels, endpoints and unselected/exterior data are measurable from
$y$.

For a loop-time slice the resulting unnormalized mass has the form
\begin{equation}\label{chIII:eq:multihost-slice}
 \one_{\{E(y)+V(\xi)+q\in I\}}G(y,q)
 \left[\bigotimes_i\mathsf B_{O_{j_i}}^{a_i,b_i}(d\xi_i)\right]
       \nu(dy)\,dq,
\end{equation}
where $q$ is the residual terminal time and every terminal exterior
test is included in the measurable factor $G(y,q)$. The indicator
in \eqref{chIII:eq:multihost-slice} imposes the observation-time
restriction on the product kernel; the normalized restricted law
need not be a product. Certified selection is applied to the local
factors, and the window inequalities are applied with this weight.
\end{lemma}
\begin{proof}
The first-hit/first-exit recursion makes the record intervals disjoint.
Apply Proposition~\ref{chIII:prop:record-product} successively to their
chronological endpoint word, and then regroup the product factors by
host. The retained coordinate list gives measurable extraction and
reconstruction; the duration identity and \eqref{chIII:eq:multihost-slice}
then follow from \eqref{chIII:eq:window-functional}.
\end{proof}

For later applications, specialize the abstract integral to the loop
decomposition as follows. Take nonempty $I_0\subset O\subsetneq V_N$
with root $r\notin O$.
Starting at $r$, alternate first entrance into $I_0$ and first exit
from $O$. A completed pair consists of an inward record in
$V_N\setminus I_0$ and an outward record in $O$. Fix $m$ completed
pairs and an endpoint word $w$ of positive probability. Select the
$m$ outward records. The exterior contains the inward records and
the future $\eta$ after the last outward exit, started at that exit vertex. Put
\begin{equation}\label{chIII:eq:censored-kernel}
 G(y,q)=\ind_{\{q<H_{I_0}(\eta),\ X_q(\eta)=r\}}\,F(y,q),
 \qquad 0\le F\le1.
\end{equation}
Then \eqref{chIII:eq:window-functional}, integrated under the conditional
product record law, is precisely the loop-time mass on this word
with $m$ completed pairs and exterior test $F$. At a loop endpoint
the path cannot be in an unfinished outward record, since $r\notin O$.
The remaining terminal part is therefore the censored inward segment
in \eqref{chIII:eq:censored-kernel}. The case $m=0$ uses the empty product,
$E=V=0$ and the original future from $r$. \hyperlink{proof.appendix.J}{Appendix~J} proves the
identity in \eqref{chIII:eq:actual-slice}, including these endpoint conventions.

For a fixed count $m$ and word $w$, denote the conditional selected
and exterior laws just described by $\mu_{m,w}$ and $\nu_{m,w}$,
and their window functional by $J_I^{m,w}$. Thus
$\mu_{m,w}$ is the product of the $m$ outward kernels in $O$;
$\nu_{m,w}$ is the product of the inward kernels in $V_N\setminus I_0$
and the future-path law from the last outward exit. Write $p_{m,w}$ for the ordinary
probability of the first $2m$ endpoints, given by the product formula
\eqref{chIII:eq:word-factor}. Empty products have value one.

\begin{corollary}\label{chIII:cor:counted-window}
Let $M\ge0$ be an integer. For each $m=0,\ldots,M$ and each feasible word $w$, let
$B_{m,w}$ be a measurable bad event of its outward records. Assume
\eqref{chIII:eq:success-hypothesis}, with $B=B_{m,w}$, under that law, with a
common duration threshold $d$ and error $\varepsilon$. Let
the two unnormalized loop masses be
\[
 \begin{aligned}
 \mathcal W_{\le M}(I)
 &=\sum_{m=0}^M\sum_{w:p_{m,w}>0}p_{m,w}
       \int J_I^{m,w}\,d(\mu_{m,w}\otimes\nu_{m,w}),\\
 \mathcal W_{\le M}^{\mathrm{good}}(I^d)
 &=\sum_{m=0}^M\sum_{w:p_{m,w}>0}p_{m,w}
       \int\ind_{B_{m,w}^c}(x)J_{I^d}^{m,w}(x,y)
                        \,d(\mu_{m,w}\otimes\nu_{m,w}).
 \end{aligned}
\]
They retain the same exterior test $F$ from \eqref{chIII:eq:censored-kernel}.
Then
\begin{equation}\label{chIII:eq:counted-window}
 \mathcal W_{\le M}(I)
 \le2\mathcal W_{\le M}^{\mathrm{good}}(I^d)
                       +(M+1)|I|\varepsilon.
\end{equation}
More generally, with nonnegative deterministic weights $a_m$, the
error is $|I|\sum_{m=0}^M a_m\varepsilon_m$ if the bound on the
duration error for count $m$ is $\varepsilon_m$.
\end{corollary}
\begin{proof}
For fixed $m$, the word probabilities satisfy $\sum_wp_{m,w}=1$.
Multiply \eqref{chIII:eq:one-success} by $p_{m,w}$ and sum over $w$.
Since $\sum_w p_{m,w}=1$, summing over words preserves the error
$|I|\varepsilon$ for this count. Now sum over $m$. The loop slices for
different counts partition the loop-time integral, but the separate
word probability distributions each have total mass one. Their errors
therefore retain the factor $M+1$. Multiplication by $a_m$ before
the last sum gives the weighted assertion.
\end{proof}

For a root-dependent family of regions with $r\notin O(r)$ for
every root, suppose these estimates hold uniformly in $r$. They
can then be averaged over all roots. For a fixed region, decompose the root average into the two parts
$r\notin O$ and $r\in O$, and apply the displayed estimate to the first. In the former case,
write the root dependence of the preceding masses explicitly as
$\mathcal W_{\le M}(I;r)$ and
$\mathcal W_{\le M}^{\mathrm{good}}(I^d;r)$, with $F=1$.
Dividing their root-averaged comparison by
\eqref{chIII:eq:normalizer}, for $0\le d<l$, gives
\begin{equation}\label{chIII:eq:normalized-transfer}
 \begin{aligned}
 \frac{n^{-1}\sum_r\mathcal W_{\le M}(I;r)}{\mathsf W_N(I;1)}
 &\le2\frac{n^{-1}\sum_r\mathcal W_{\le M}^{\mathrm{good}}(I^d;r)}
                 {\mathsf W_N(I;1)}\\
 &\quad+\frac{(M+1)|I|\varepsilon}{\mathsf W_N(I;1)}.
 \end{aligned}
\end{equation}
Thus the probability-scale requirement on central windows is
$n(M+1)\varepsilon=o(1)$. For weighted counts, the corresponding
weighted sum of duration errors must be $o(n^{-1})$. The superexponential estimates above
permit the finite mixtures needed later. The main comparison uses one successful replacement and one window
enlargement.

\clearpage
\endgroup

\begingroup
\mergedchaptersetup{IV}{Geometry and tightness of low local times}
\begin{center}
{\large\scshape Chapter IV\par}\vspace{6pt}
{\LARGE Geometry and tightness of low local times\par}
\end{center}
\let\mathscr\mathcal
\allowdisplaybreaks[2]
\raggedbottom
We first prove that $|D_N^{\rm wide}(U_u(t))|/\sqrt L$ tends to zero.
The Gaussian contact bound \eqref{chIV:eq:detect-main}, combined with
counted coarea, then proves tightness of $M_{u,v}$.
The matching inequality \eqref{chIV:eq:matching} bounds its change
between two times by the early-parent weights whose neighbourhoods
are hit, together with $e_N(s)$. These estimates are established
before identifying the activity limit.

\subsection*{4.1. Traversals and harmonic occupation}
Retain the rate-one walk on $V_N=(\ZZ/N\ZZ)^2$, its local times
$\ell_t$, and the Green kernel of Chapter~I. Throughout,
\[
 n=N^2,\quad L=\log N,\quad g=2/\pi,\quad b_N=gnL,
 \qquad s_N(c)=s_N+cb_N,
\]
\[
 s_N=b_N(2L-\tfrac32\log L)+\tfrac{9g}{32}n(\log L)^2.
\]
The symbols $d$ and $d_1$ denote periodic Euclidean and graph distance,
respectively; $d\le d_1\le\sqrt2d$. In this chapter and Appendices~L--Q, write
\[
 B(y,R)=\{x\in V_N:d(x,y)<R\},\qquad
 B_{d_1}(y,R)=\{x\in V_N:d_1(x,y)\le R\}.
\]
For radii below $N/4$, $B(y,R)$ is identified with its Euclidean
lattice disk in a planar lift. Constants may depend on fixed thresholds
and fixed central windows, but are uniform in the root and the initial
vertex. Write $U_u(t)=\{x:\ell_t(x)\le u\}$.

For the local calculations only, fix
\[
 \gamma=1/8192,\qquad Q=L^{16\gamma},\qquad
 q\in[Q,L-Q],\qquad w=q^\gamma.
\]
Let $I_{q,y}=B(y,e^{q+w})$ and $O_{q,y}=B(y,e^{q+2w})$.
Starting with $\tau_0=0$, define alternating entrance and exit times
\[
 \sigma_i=\inf\{s\ge\tau_{i-1}:X_s\in I_{q,y}\},\qquad
 \tau_i=\inf\{s\ge\sigma_i:X_s\notin O_{q,y}\},\qquad i\ge1.
\]
A complete $q$-block is the path on $[\sigma_i,\tau_i]$; all returns
before $\tau_i$ belong to this block. At observation time $t$ put
$m_q(y)=\#\{i:\tau_i\le t\}$. We suppress the dependence on $t$
in this count and in the quantities below. Its endpoint word is the
ordered list $(X_{\sigma_i},X_{\tau_i})_{\tau_i\le t}$.
When the root lies outside $O_{q,y}$, each entrance is within one
lattice step of the boundary of $I_{q,y}$.

Let
\[
 \tau_{B(y,e^{q+1})}=\inf\{s\ge0:X_s\notin B(y,e^{q+1})\},\qquad
 \mu_q(z)=\PP_y(X_{\tau_{B(y,e^{q+1})}}=z).
\]
Thus $\mu_q$ is the exit harmonic measure from the smaller disk.
For the $i$th complete block, define
\[
 \ell^{(i)}(z)=\int_{\sigma_i}^{\tau_i}\one_{\{X_s=z\}}\,ds,
 \qquad A_i=\frac1g\sum_z\mu_q(z)\ell^{(i)}(z).
\]
Put
\[
 Y_q(y)=\sqrt{q^\gamma m_q(y)},\qquad
 \mathcal O_q(y)=\left(g^{-1}\sum_z\mu_q(z)\ell_t(z)\right)^{1/2}.
\]
For a loop with root outside $O_{q,y}$, all occupation in this
harmonic average comes from complete blocks: the terminal path is
outside $I_{q,y}$ before its next entrance. In particular, for all
sufficiently large $N$, $\supp\mu_q\subset I_{q,y}$ and
\[
 \mathcal O_q(y)^2=\sum_{i=1}^{m_q(y)}A_i
 \quad\text{on }\{X_t=r,\ r\notin O_{q,y}\}.
\]
On a fixed endpoint word,
the resampling kernel is the product of the kernels of Section~3.5;
the loop-time restriction is imposed separately as in Section~3.6.
We abbreviate each kernel $\mathsf B_{O_{q,y}}^{a,b}$ by $P^{a,b}$;
unsubscripted probabilities and expectations in local block estimates
refer to this endpoint-conditioned law or its specified product.

\begin{proposition}\label{chIV:prop:local}
Let $a$ be an entrance within one lattice step of the boundary of
$I_{q,y}$ and let $b$ be an exterior neighbour of $O_{q,y}$.
Uniformly over these endpoints, the unconditioned exit probability is
\[
 h_{O_{q,y},b}(a)=\PP_a(X_{\tau_{O_{q,y}}}=b)\ge cR^{-C},
 \qquad R=e^{q+2w},
\]
with $h_{D,b}$ as in Section~3.5. Under the conditioned block law
$P^{a,b}$, the duration $\tau=\tau_{O_{q,y}}$ satisfies
\[
 P^{a,b}(\tau>d)\le CR^Ce^{-cd/R^2}.
\]
The probability of avoiding $B(y,e^q)$ is bounded below by $c>0$;
the probability of visiting any specified point of this core is at
least $cw/q$. Under a product of $m$ specified block kernels, the
harmonic occupations $A_i$ defined above satisfy, for every integer $j\ge1$,
\begin{equation}\label{chIV:eq:local-moments}
 EA_i=w+O(we^{-cw}),\qquad E A_i^j\le j!(Cw)^j.
\end{equation}
For $Y=\sqrt{mw}$ and $Ce^{-cw}Y\le\Delta\le Y$,
\begin{equation}\label{chIV:eq:local-sqrt}
 P\left(\left|\sqrt{\sum_{i=1}^m A_i}-Y\right|>\Delta\right)
 \le2e^{-c\Delta^2/w}.
\end{equation}
The bounds remain uniform under the logarithmic radius perturbations
specified in \hyperlink{proof.appendix.L}{Appendix~L}. \end{proposition}

\hyperlink{proof.appendix.L}{Appendix~L}.1 proves the endpoint and duration assertions in
\eqref{chIV:eq:R6.1}--\eqref{chIV:eq:R6.4}. The Kac moment calculation
and Bernstein estimate are \eqref{chIV:eq:R7.1}--\eqref{chIV:eq:R7.4}
in \hyperlink{proof.appendix.M}{Appendix~M}.1.
In each Kac product, the conditioned Green kernel contributes the
endpoint factor at the final site, as in \eqref{chIV:eq:R7.1}--\eqref{chIV:eq:R7.2}.

Let $I_0$ be an interval contained in $[s_N-Wb_N,s_N+Wb_N]$,
of length at least $\delta b_N$, where $W,\delta>0$ are fixed.
For enlargements $I$ of this interval set
\begin{equation}\label{chIV:eq:loop-integral}
 \mathcal Q_I(F)=\frac1{|I_0|}\sum_r\int_I
             \EE_r[\one_{\{X_t=r\}}F(r,t,\omega)]\,dt.
\end{equation}
Thus $\mathcal Q_{I_0}(1)=1+o(1)$. This normalization is kept when the
window is enlarged. The exact normalized probability law is obtained
by division by $\mathcal Q_{I_0}(1)$. In the notation of
Section~3.2, the exact relation is
\[
 \mathcal Q_I(F)=\frac{n}{|I_0|}\mathsf W_N(I;F)
  =\frac{n\mathsf W_N(I;1)}{|I_0|}\mathsf Q_N^I(F).
\]
For a measure $\mu$ we write $\mu[F;E]=\int F\ind_E\,d\mu$.

For a centre $x$ and the scales considered here, write
\[
 \omega_{q,x}(\xi)=\PP_x(X_{\tau_{q,x}}=\xi),\qquad
 \tau_{q,x}=\inf\{t\ge0:X_t\notin B(x,e^{q+1})\}.
\]
These are the weights $\mu_q$ of Section~4.1 with centre $x$.
For a nonnegative vector $W$, put
\[
 \mathcal R_{q,x}^{\rm sq}(W)
    =\left(\frac1{2g}\sum_\xi\omega_{q,x}(\xi)W(\xi)\right)^{1/2}.
\]
For the pinned field $h=h^r$, write $P_h,E_h$ for its probability
and expectation with the root fixed.
\begin{lemma}[Coarse-scale probes]\label{chIV:lem:coarse-probes}
Fix $0<c_0\le C_0<\infty$ and a finite support-radius constant $C$.
Suppose $\Var h(x)\in[c_0L,C_0L]$. Let $\omega_q$ be probability weights supported
within periodic Euclidean distance $Ce^{q+1}$ of $x$, with $q\ge1$. Write
\[
 V_q(h)=\frac1{2g}\sum_\xi\omega_q(\xi)h(\xi)^2.
\]
There are constants $D,C_1,c>0$ such that
\begin{equation}\label{chIV:eq:coarse-probe-single}
 P_h\{|h(x)|\le1,\ V_q(h)\le C_1(q+1),\ \|h\|^2\le DnL\}
 \ge cL^{-1/2}.
\end{equation}
For fixed $\chi>0$ and a deterministic set of at most $L$ scales $q\ge a_{\rm min}\to\infty$
satisfying $L e^{-ca_{\rm min}^\chi}\to0$ for every fixed $c>0$, also
\begin{equation}\label{chIV:eq:coarse-probe-simultaneous}
 P_h\{|h(x)|\le1,\ V_q(h)\le C_1q^{1+\chi}\ \forall q,
                    \ \|h\|^2\le DnL\}\ge cL^{-1/2}.
\end{equation}
The constants may depend on $c_0,C_0,C$ and, in the simultaneous
assertion, $\chi$, but are uniform in $x,r$, the weights and the indicated scales.
\end{lemma}
The proof is in \hyperlink{proof.appendix.M}{Appendix~M}.2.

\begin{proposition}\label{chIV:prop:count}
Fix $C_1<\infty$ and put $m_L=\sqrt2L-3\log L/(4\sqrt2)$.
Uniformly for $Q\le q\le L-Q$, $q\le a\le C_1L$, and all $y$,
with constants allowed to depend on $C_1$,
\begin{multline}\label{chIV:eq:count-main}
 \mathcal Q_{I_0}\{d(r,y)\ge Ne^{-Q/2},\ Y_q(y)\le a\}\\
 \le C\sqrt L\exp\left\{-\frac{(m_L-a-L^{2\gamma}-C\log L)_+^2}
                                    {L-q+C}\right\}
       +e^{-cL\log L}.
\end{multline}
Outside an event of $\mathcal Q_{I_0}$-mass $o(1)$, simultaneously over the
integer scales $Q\le q\le L-Q$ and centres with
$d(r,y)\ge2Ne^{-Q/2}$,
\begin{equation}\label{chIV:eq:barrier-main}
 Y_q(y)\ge\sqrt2q-L^{3\gamma}.
\end{equation}
\end{proposition}

The coarea proof uses the event $B_{q,y}$ of
\hyperlink{proof.appendix.M}{Appendix~M}.3, with
$P_h(B_{q,y})\ge cL^{-1/2}$ and $\|h\|^2\le DnL$ on $B_{q,y}$.
The same $D$ is fixed before $N$. Take $\upsilon=a\in[q,C_1L]$.
Equation~\eqref{chIV:eq:R7.9} gives
\[
 \mathcal O_q(y)\le\upsilon+L^{2\gamma},\quad h\in B_{q,y}
 \quad\Longrightarrow\quad
 \mathcal R_{q,y}^{\rm sq}(2\ell_t+h^2)
       \le\upsilon+L^{2\gamma}+C_3=v_q.
\]
The Gaussian harmonic-average tail bounds this image event by the
exponent in \eqref{chIV:eq:count-main}.
To pass from occupation to the count, use the one-success window
inequality of Chapter~III, with success probability $1-o(1)$ supplied
by \eqref{chIV:eq:local-sqrt}. Lemmas~\ref{chIV:lem:coarse-probes} and
\ref{chIV:lem:low-harmonic-coarea} give the probe and its image bound;
\hyperlink{proof.appendix.M}{Appendix~M}.3 applies the window transfer
\eqref{chIII:eq:counted-window}. \hyperlink{proof.appendix.M}{Appendix~M}.4 proves the simultaneous
assertion by a spatial grid and chronological injection of crossings.

\subsection*{4.2. Weighted deviations and coarse repulsion}
The lower barrier alone does not control occupation at a soft vertex.
The inner-scale softness and outer-scale occupation deviation are
correlated. The following estimate keeps their common weight.

\begin{proposition}\label{chIV:prop:deviation}
Let
\[
 \Theta=\left\{\begin{gathered}
   Y_q(y)\ge\sqrt2q-L^{3\gamma}\ \text{ for every }q\in[Q,L-Q]\cap\mathbb Z\\
   \text{and every }y\in V_N\text{ with }d(r,y)\ge2Ne^{-Q/2}
 \end{gathered}\right\}.
\]
For each fixed $u>0$,
\begin{equation}\label{chIV:eq:deviation-main}
 \mathcal Q_{I_0}[\one_\Theta\,\#\mathcal E_N^{\rm rad}]=o(1),
\end{equation}
where
\[
 \mathcal E_N^{\rm rad}=\left\{(x,q)\in V_N\times\mathbb Z:\begin{gathered}
 \ell_t(x)\le u,\quad 2Q\le q\le L-2Q,\\
 |\mathcal O_q(x)-Y_q(x)|>L^{4\gamma},\quad d(r,x)\ge2Ne^{-Q/2}
 \end{gathered}\right\}.
\]
Moreover,
$\mathcal Q_{I_0}\{\exists x\in U_u(t):d(x,r)<2Ne^{-Q/2}\}=o(1)$.
\end{proposition}

To see the cancellation, fix $x\in V_N$ and an integer
$q\in[2Q,L-2Q]$, and set $l=\lfloor q-3q^\gamma\rfloor$.
In this calculation write $Y_k=Y_k(x)$ and $\mathcal O_k=\mathcal O_k(x)$ for $k=q,l$.
Fix a feasible $q$-word and let $P,E$ denote its product block law
and expectation. Write $Y_q=\sqrt2q+s$ and $X=s-(Y_l-\sqrt2l)$. Put
\[
 J_l=\{Y_l-\sqrt2l\ge-L^{3\gamma},\ Y_l\le10L\}.
\]
On the lower barrier, and after excluding counts above $10L$, the
nested-crossing estimate gives
\[
 P(X>a,\ J_l)\le2e^{-ca^2/q^\gamma},\qquad
 E(e^{\lambda X}\one_{J_l})\le C_\lambda e^{C_\lambda q^\gamma}.
\]
Cauchy--Schwarz and \eqref{chIV:eq:local-sqrt} therefore yield
\begin{equation}\label{chIV:eq:weighted-main}
 E\left[e^{2\sqrt2X}\one_{J_l}
       \one_{\{|\mathcal O_q-Y_q|>L^{4\gamma}\}}\right]
 \le C\exp\{Cq^\gamma-cL^{8\gamma}/q^\gamma\}.
\end{equation}
Resample the $l$-domain first and the $q$-domain second.
On $J_l$, \eqref{chIV:eq:R8.2} gives
\[
 p_{\rm soft}(Y_l)
 \le e^{-2l-2\sqrt2s+Cq^\gamma}e^{2\sqrt2X}.
\]
The bounded-weight window inequality and \eqref{chIV:eq:weighted-main}
give the bound $c_w$ of Appendix~N.3, with
\[
 c_w\le\min\{1,C\exp(-2l-2\sqrt2s+Cq^\gamma-cL^{8\gamma}/q^\gamma)\}.
\]
Multiplying by the main term of \eqref{chIV:eq:R8.4} gives
\[
 c_w\,C\sqrt L\,e^{-2(L-q)+2\sqrt2s+CL^{2\gamma}}
 \le C\sqrt L\,e^{-2L-c'L^{7\gamma}},
\]
by \eqref{chIV:eq:count-deviation-summable}. Summing over the outer
counts, vertices and scales, with the additive and duration errors
estimated in Appendix~N.3, proves the assertion.
\hyperlink{proof.appendix.N}{Appendix~N}.1--N.3 gives the endpointwise bounds, both high-count
exclusions, and the weighted estimate \eqref{chIV:eq:weighted-count-deviation}.

For each fixed softness threshold $u>0$, define the coarse exceptional set
\[
 \mathcal O_{N,u}(t)=\left\{x\in U_u(t):\begin{gathered}
 Y_q(x)\le\sqrt2q+L^{8\gamma}\text{ for some }q\in\mathbb Z,\\
 L^{1/10}\le q\le L-L^{1/10}
 \end{gathered}\right\}.
\]
\begin{proposition}\label{chIV:prop:repulsion}
For every fixed $u>0$ and $\eta>0$,
$\mathcal Q_{I_0}\{|\mathcal O_{N,u}(t)|>\eta\sqrt L\}=o(1)$.
\end{proposition}
\begin{proof}
Use the event
\[
 \mathcal G_N=\Theta\cap\{\mathcal E_N^{\rm rad}=\varnothing\}
   \cap\{\nexists x\in U_u(t):d(x,r)<2Ne^{-Q/2}\},
\]
where $\Theta$ and $\mathcal E_N^{\rm rad}$ are from
Proposition~\ref{chIV:prop:deviation} at the fixed threshold $u$.
Propositions~\ref{chIV:prop:count} and \ref{chIV:prop:deviation}
give $\mathcal Q_{I_0}(\mathcal G_N^c)=o(1)$.
With $\mathcal B_N(s,W)$ defined in
\hyperlink{proof.appendix.O}{Appendix~O}.4 at threshold $u$, the
probe majorant proved there and counted coarea give
\[
 \frac{\mathcal Q_{I_0}[|\mathcal O_{N,u}|;\mathcal G_N]}{\sqrt L}
 \le C\sup_{s\in I^+}E_{\nu_{N,s}}
        |\mathcal B_N(s,W_s^\sigma(\phi))|,
 \qquad I^+=I_0+[0,DnL/2].
\]
Here $D$ is the common energy cutoff from
\eqref{chIV:eq:coarse-probe-simultaneous}.
Equations~\eqref{chIV:eq:coarse-vertex-image-bound} and
\eqref{chIV:eq:R9.8} bound the supremum by
\[
 C(1+\log L)L^{C\varepsilon}
 \{L^{96\gamma-1/20}+Le^{-c(\log L)^2}\}+o(1)=o(1)
\]
for a sufficiently small fixed $\varepsilon>0$.
Indeed $1/20-96\gamma=49/1280>0$.
Markov's inequality, together with
$\mathcal Q_{I_0}(\mathcal G_N^c)=o(1)$, proves the assertion.
\end{proof}

\subsection*{4.3. Intermediate-distance defects}
For a set $U\subset V_N$ define its wide defects by
\begin{multline}\label{chIV:eq:wide-def}
 D_N^{\rm wide}(U)=\{x\in U:\exists z\in U,\\
       \tfrac12e^{L^{1/8}}<d(x,z)<2Ne^{-L^{1/8}}\}.
\end{multline}
The canonical construction of Chapter~III instead uses
\[
 \rho_N=\tfrac12e^{L^{1/6}},\qquad R_N=Ne^{-2L^{1/8}},
 \]
\[
 D_N^{\rm can}(U)=\{x\in U:\exists z\in U,
                             \rho_N<d_1(x,z)<R_N\}.
\]
For large $N$, $D_N^{\rm can}(U)\subset D_N^{\rm wide}(U)$, by the
two distance comparisons. The components of
$U\setminus D_N^{\rm can}(U)$ joined at graph distance at most $\rho_N$
have graph diameter at most $\rho_N$ and graph separation at least $R_N$.
Indeed a chain that first exceeded diameter $\rho_N$ would contain
a pair at distance in $(\rho_N,2\rho_N]$, which is forbidden.
We call these components the canonical parents.

\begin{theorem}\label{chIV:thm:defects}
For fixed $u>0$, $W\ge0$ and $\eta>0$,
\begin{equation}\label{chIV:eq:wide-window}
 \sup_r \PP_r\left(\sup_{|c|\le W}
       \frac{|D_N^{\rm wide}(U_u(s_N(c)))|}{\sqrt L}>\eta\right)
 \longrightarrow0.
\end{equation}
Moreover $|D_N^{\rm wide}(U_u(t))|/\sqrt L\to0$ in probability
under every normalized central loop window as above.
\end{theorem}
\begin{proof}
First work in the loop integral. Outside $\mathcal O_{N,u}$, a wide
witness pair $(x,z)$ at logarithmic distance $q$ has
$Y_q(x)=\sqrt2q+s$, $s\ge L^{8\gamma}$. With prescribed exterior
word the two-target softness estimate is
\[
 P(\ell_t(x)\le u,\ell_t(z)\le u\mid\text{word})
 \le\exp\{-4q-4\sqrt2s+Cq^\gamma\}.
\]
There are at most $CN^2e^{2q}$ ordered pairs at this scale.
Multiplication by \eqref{chIV:eq:count-main}, followed by summation over
the possible counts, leaves a factor $e^{-cL^{8\gamma}}$.
The exceptional count branches and duration errors are summable.
Together with Proposition~\ref{chIV:prop:repulsion}, this proves the
loop assertion. \hyperlink{proof.appendix.P}{Appendix~P} gives the two-target
calculation and the exponent summation.

For a stationary ordinary path, choose $\lceil\eta\sqrt L\rceil$
distinct defective vertices and one witness for each, if that many
exist. Append a bridge of length between $\delta b_N$ and
$2\delta b_N$. With probability $1-o(1)$ none of these vertices is
near either endpoint. On this endpoint-separation event, the bridge
hitting bound \eqref{chIV:eq:R10.6} applies to both vertices of every
selected pair. Linearity of conditional expectation bounds the expected
lost fraction by $C\delta+o(1)$. Markov's inequality then gives
survival of at least half the pairs with probability at least $1/2$
for small fixed $\delta$. The completed path is a loop; integrating
its terminal time proves the deterministic-time bound. The endpoint
exclusion and density normalization are proved in
\hyperlink{proof.appendix.P}{Appendix~P}.

Fix a sufficiently large constant $C_{\rm mix}\ge1$ for the
heat-kernel estimate \eqref{chI:eq:mixing}. Delete an initial segment
of length $C_{\rm mix}n\log L=o(b_N)$ to treat an arbitrary starting point. The remaining starting distribution
is uniformly close to stationary, and deletion can only enlarge
the soft set and its defects. Finally $U_u(t)$ decreases with $t$,
and $D_N^{\rm wide}(U)\subset D_N^{\rm wide}(V)$ when $U\subset V$.
The whole-window supremum in \eqref{chIV:eq:wide-window} is consequently
bounded by its earliest-time value.
\end{proof}

\subsection*{4.4. Detecting large activity}
Write $\phi=\phi_N$ and $\phi^{(r)}=\phi_N^{(r)}$ for the mean-zero
field and its regression residual from Chapter~II. Thus
\[
 \begin{aligned}
 \phi^{(r)}(x)&=\phi(x)-G_N(x,r)\phi(r)/v_N,
       &v_N&=G_N(r,r),\\
 h^r(x)&=\phi(x)-\phi(r).
 \end{aligned}
\]
Activity uses $\phi^{(r)}$, and coarea uses $h^r$, related by
\eqref{chIII:eq:field-connection}. Write $P_h,E_h$
for probability and expectation with respect to the pinned field $h^r$,
with the root and occupation vector held fixed. For a parent $D$
and a nonnegative occupation vector $\ell$, recall the activity of
Section~3.1:
\begin{multline}\label{chIV:eq:activity}
 a_D^G(\ell;r,v)=c_g E\Leb\bigcup_{x\in D}
 [-\phi^{(r)}(x)-\sqrt{(v-2\ell(x))_+},
  -\phi^{(r)}(x)+\sqrt{(v-2\ell(x))_+}],\\
 \qquad c_g=(2\pi g)^{-1/2}.
\end{multline}
Here $E$ integrates only the regressed Gaussian field
$\phi^{(r)}$; the occupation vector and root are fixed. For a
general occupation vector put
\[
 M_{u,v}(\ell;r)=\frac1{\sqrt L}
   \sum_{a\in\operatorname{Rep}_N(\{x:\ell(x)\le u\})}
       a_{D_a(\{x:\ell(x)\le u\})}^G(\ell;r,v).
\]
On a path write $M_{u,v}(t)=M_{u,v}(\ell_t;r)$, with $r=X_0$.
Thus $A_N(t)=M_{1,1/4}(t)$ and $\widehat A_N(c)=M_{1,1/4}(s_N(c))$. The same root is used when an
earlier prefix is evaluated.
For fixed $u,v$, write the normalized defect error as
\[
 e_N(t)=\frac{2c_g\sqrt v}{\sqrt L}
              |D_N^{\rm can}(U_u(t))|.
\]

\begin{lemma}\label{chIV:lem:envelope}
Uniformly in $\ell\ge0$, $r$ and $0<v\le1$,
$a_D^G(\ell;r,v)\le CL^{1/6}$. The normalized activity of parents
whose representatives are at graph distance at most $R_N/4$ from any one prescribed
vertex is at most $CL^{-1/3}$.
\end{lemma}
\begin{proof}
The union in \eqref{chIV:eq:activity} has length at most $2$ plus the
range of $\phi^{(r)}$ on $D$. Put $d_D=\operatorname{diam}_{d_1}D$.
Conditional increment variances are bounded by $C\log(2+d_D)$.
The Gaussian maximum bound and $|D|\le C(1+d_D)^2$ bound its
expected range by $C\log(2+d_D)\le CL^{1/6}$. At most one representative lies
in a graph ball of radius $R_N/4$, since representatives are
$R_N$-separated. Divide by $\sqrt L$.
\end{proof}

\begin{lemma}[Canonical-family colouring]\label{chIV:lem:parent-colouring}
Suppose $7R_N<N$ and $\rho_N\le R_N/4$. In any family of canonical
parents, let $x_D\in D$ be the canonical representative. At most $49$
parents have $d_1(x_D,r)<3R_N$. The remaining parents admit a partition
into at most $49$ colour classes, each consisting of pairwise admissible
clusters in the sense of Section~2.4.
\end{lemma}
\begin{proof}
By Lemma~\ref{chIII:lem:geometry}, the representatives are
$R_N$-separated in $d_1$. Embed them in the continuous flat torus of
side $N$, with its periodic $\ell^1$ metric. Their open balls of radius
$R_N/2$ are disjoint. Those with centres in a ball of radius $3R_N$
are contained in the concentric ball of radius $7R_N/2$. Since
$7R_N<N$, the planar area formula gives
\begin{equation}\label{chIV:eq:parent-packing}
 \#\{D:d_1(x_D,z)<3R_N\}
 \le\frac{2(7R_N/2)^2}{2(R_N/2)^2}=49,
 \qquad z\in V_N.
\end{equation}
Delete the indicated root-near parents. Join two remaining parents
when $d_1(x_D,x_{D'})<3R_N$. Each closed neighbourhood has at most
$49$ vertices, so greedy colouring uses at most $49$ colours.
Within any colour, using $R_N=D_N^{\rm sep}$ and
$2\rho_N=R_N^{\rm loc}$ from Section~2.4,
\begin{equation}\label{chIV:eq:parent-admissibility}
 \begin{aligned}
 \operatorname{diam}_d D&\le\rho_N=R_N^{\rm loc}/2,\\
 d(D,r)&\ge3R_N/\sqrt2-\rho_N\ge R_N,\\
 d(D,D')&\ge3R_N/\sqrt2-2\rho_N\ge R_N\qquad(D\ne D').
 \end{aligned}
\end{equation}
\end{proof}

\begin{proposition}\label{chIV:prop:detection}
There are fixed $a,p,D_0>0$ such that, for large $N$, every root,
every nonnegative occupation vector, every $u\ge0$, and $0<v\le1$, if
$M_{u,v}\ge B\ge1252$, then
\begin{equation}\label{chIV:eq:detect-main}
 P_h\left(\#\{x:2\ell(x)+h^r(x)^2\le v\}>aB,
                    \|h^r\|_2^2\le D_0nL\right)\ge p.
\end{equation}
\end{proposition}
\begin{proof}
Write $x_D$ for each canonical representative and apply
Lemma~\ref{chIV:lem:parent-colouring}. By the single-parent
bound in Lemma~\ref{chIV:lem:envelope}, the deleted normalized activity is
\[
 \frac1{\sqrt L}\sum_{D:\ d_1(x_D,r)<3R_N}a_D^G(\ell;r,v)
 \le49CL^{-1/3}<1
\]
for all sufficiently large $N$, uniformly in the profile and $0<v\le1$.
There are at most $49\le626$ colour classes. Choose the class
$\mathcal C$ with largest normalized activity, so that
\[
 V=\frac1{\sqrt L}\sum_{D\in\mathcal C}a_D^G(\ell;r,v)
       \ge(M_{u,v}-1)/626.
\]
For this family define the contact events and their truncated count by
\[
 \begin{gathered}
 E_D=\{\exists x\in D:2\ell(x)+h^r(x)^2\le v\},\\
 J=\one_{\{|S_r|\le1\}}\sum_{D\in\mathcal C}\one_{E_D},
 \qquad S_r=\frac{\phi(r)}{\sqrt{v_N}}.
 \end{gathered}
\]
Since $\phi(r)=-n^{-1}\sum_xh^r(x)$, this count is measurable under
the pinned law. Proposition~\ref{chII:prop:relative-contact}, applied to
the admissible family from \eqref{chIV:eq:parent-admissibility}, gives
fixed constants $a_0,C_0>0$, uniformly including zero activities, with
\begin{equation}\label{chIV:eq:contact-moments}
 a_0V\le E_hJ,\qquad E_hJ^2\le C_0(V+V^2).
\end{equation}
For $V\ge1$, Paley--Zygmund gives a fixed $p_0>0$ such that
$P_h(J>a_0V/2)\ge p_0$. The pinned covariance has trace
$2nG_N(r,r)\le CnL$. Fix $D_0$ so that
\[
 \sup_{r,N}P_h(\|h^r\|_2^2>D_0nL)\le p_0/2.
\]
For every $B\ge1252$ and $M_{u,v}\ge B$, the chosen class has
$V\ge B/1252\ge1$. Distinct parents give distinct low points.
Subtracting the energy tail proves \eqref{chIV:eq:detect-main} with
$a=a_0/2504$ and $p=p_0/2$.
\end{proof}

\begin{lemma}\label{chIV:lem:shell-count}
Fix $0\le W<\infty$ and let $s=s_N+nv_N+\theta b_N$, $|\theta|\le W$,
with $v_N=G_N(0,0)$. For every fixed $v\ge0$, the fixed-root
field--branch law $\nu_{N,s}^r$ of Section~3.3 satisfies
\[
 \lim_{K\to\infty}\limsup_{N\to\infty}
 \sup_{r,\,|\theta|\le W}
 \nu_{N,s}^r\!\left(\#\{x:W_s^\sigma(\phi)(x)\le v\}>K\right)=0.
\]
\end{lemma}
\begin{proof}
Write the shell image as $(\phi+\sigma a_s)^2$, with its exact
fixed-root density. The shell comparison \eqref{chII:eq:shell-reference-error}, with its
fixed-root extension in \hyperlink{proof.appendix.Q}{Appendix~Q}.2, bounds the
total-variation error from the reference field on the positive branch by
$C/\sqrt L$. The probability that $a_s$ differs by more than a
fixed sufficiently large $T$ from $\sqrt{2s/n-G_N(r,r)}$ is at
most $C e^{-cLT}$. On the complement, every image coordinate
below $v$ is a high coordinate of one of $\phi,-\phi$, at height
at least $m_N-C-T-\sqrt v$, where $m_N$ is the raw-field extremal
centre of Chapter~II. The high-point count bounds proved there
are uniform at any fixed height below $m_N$. First choose their
count cutoff and then take $N$ large to absorb the two shell
errors. This proves the assertion also for any other fixed upper
threshold, such as $2u+1$.
\end{proof}

We use the following deterministic comparison in the tightness and
time-comparison proofs. On a path rooted at $r$, take $0\le s\le t$.
Write $D^-$ for the early canonical
parents. Surviving clean points of
one early parent lie in one late parent; distinct early parents
cannot merge. Every additional late point was an early defect.
Its contribution to union length is at most $2\sqrt v$.
If an early parent $D$ has no point in the later soft set, the path
has visited $D$ between $s$ and $t$. Since
$D\subset B_{d_1}(z_D,2\rho_N)$ for its canonical representative
$z_D$, its disappearance implies
$\exists t'\in(s,t]:d_1(X_{t'},z_D)\le2\rho_N$.
Summing the corresponding early activities gives
\begin{equation}\label{chIV:eq:matching}
 \begin{aligned}
 |M_{u,v}(t)-M_{u,v}(s)|
 &\le L^{-1/2}\sum_{D^-}a_{D^-}^G(\ell_s;r,v)\\
 &\qquad{}\times\ind_{\{\exists t'\in(s,t]:\ d_1(X_{t'},z_{D^-})\le2\rho_N\}}
       +e_N(s),\qquad s\le t.
 \end{aligned}
\end{equation}
It also gives $M_{u,v}(t)\le M_{u,v}(s)+e_N(s)$.
For a spatial test, a change of representative contributes its
modulus of continuity over distances at most $\rho_N$, multiplied
by the early activity mass; Section~4.6 gives the resulting bound.

\subsection*{4.5. Tightness before limit identification}
\begin{theorem}\label{chIV:thm:tight}
For $u>0$, $0<v\le1$, $W\ge0$ and $\eta>0$, there is $B<\infty$
such that, for all sufficiently large $N$,
\begin{equation}\label{chIV:eq:mass-tight}
 \sup_r\sup_{|c|\le W}\PP_r(M_{u,v}(s_N(c))>B)\le\eta.
\end{equation}
For every root, uniformly over $l,r'\in[s_N-Wb_N,s_N+Wb_N]$,
\begin{equation}\label{chIV:eq:loop-tight}
 \frac n{b_N}\int_l^{r'}\EE_r[\one_{\{X_t=r\}}
                  \one_{\{M_{u,v}(t)>B\}}]\,dt\le\eta
 \qquad(l\le r').
\end{equation}
Both \eqref{chIV:eq:mass-tight} and \eqref{chIV:eq:loop-tight}
also hold with $M_{u,v}(t)$ replaced by $n_N(U_u(t))/\sqrt L$,
uniformly for $u$ in a fixed bounded nonnegative interval.
\end{theorem}
\begin{proof}
Fix the detection constants $a,p,D_0$ from
Proposition~\ref{chIV:prop:detection}. With this fixed energy cutoff,
choose $B$ using Lemma~\ref{chIV:lem:shell-count}. Applying
\eqref{chIV:eq:detect-main} on $M_{u,v}(t)>B$ sends the image time into
$[l,r'+D_0nL/2]$; on this interval the probability of more than $aB$
low coordinates is uniformly as small as required. More precisely, the density
in time, after integrating the field and branch, is
$\varrho_{N,r}(s)\le C/n$, and counted coarea gives
\[
 \begin{aligned}
 &p\int_l^{r'}\EE_r[\one_{\{X_t=r\}}\one_{\{M_{u,v}(t)>B\}}]dt\\
 &\quad\le\int_l^{r'+D_0nL/2}\varrho_{N,r}(s)
       \nu_{N,s}^r\!\left(\#\{x:W_s^\sigma(\phi)(x)\le v\}>aB\right)ds.
 \end{aligned}
\]
The enlarged interval has
length at most $(2W+D_0/(2g))b_N$. Division by the fixed detection
probability proves \eqref{chIV:eq:loop-tight}.

We give the deterministic-time transfer explicitly. At time $s$
let $E=\{M_{u,v}(s)>2B,\ e_N(s)\le1\}$, with the defect error
$e_N$ defined in Section~4.4.
For $T\in[\delta b_N,2\delta b_N]$, the endpoint-weighted hit
bound proved below is
\begin{equation}\label{chIV:eq:endpoint-hit}
 \PP_x(X_T=y,\ H_{B_{d_1}(z,2\rho_N)}\le T)
 \le\frac Cn\left(\delta+\frac{L^{1/6}+1}{L}\right)
\end{equation}
when $d_1(x,z),d_1(y,z)\ge R_N/4$.
For each early parent $D$, set
\[
 w_D=L^{-1/2}a_D^G(\ell_s;r,v),\qquad
 \widetilde w_D=
 \begin{cases}w_D/M_{u,v}(s),&\text{on }E,\\0,&\text{on }E^c.\end{cases}
 \qquad \sum_D\widetilde w_D=\one_E.
\]
These weights are measurable from the prefix up to $s$.
Lemma~\ref{chIV:lem:envelope} gives
\[
 \sum_{D:\ d_1(z_D,X_s)<R_N/4\ \text{or}\ d_1(z_D,r)<R_N/4}
       w_D\le CL^{-1/3}.
\]
Thus, for $B\ge4$ and large $N$, \eqref{chIV:eq:matching} implies
\[
 \begin{aligned}
 E\cap\{M_{u,v}(s+T)\le B\}&\subseteq\\
 &\left\{\sum_{\substack{D:\ d_1(z_D,X_s)\ge R_N/4\\
                              d_1(z_D,r)\ge R_N/4}}
 \widetilde w_D\,
 \one_{\{\exists t'\in(s,s+T]:d_1(X_{t'},z_D)\le2\rho_N\}}
       \ge\tfrac14\right\}.
 \end{aligned}
\]
Condition on the prefix, apply \eqref{chIV:eq:endpoint-hit} to these
far weights, and integrate over the prefix. Since
$\sum_D\widetilde w_D=\one_E$, the hit bound is multiplied by
$\PP_r(E)$. Uniform mixing for the closing endpoint then gives
\[
 \frac1{2n}\PP_r(E)\le
 \PP_r(X_{s+T}=r,M_{u,v}(s+T)>B)
 +\frac{4C}{n}\left(\delta+\frac{L^{1/6}+1}{L}\right)\PP_r(E).
\]
Here the first term on the left uses uniform mixing of the
closing endpoint. Choose a fixed $\delta>0$ so that the last
coefficient is at most $1/(4n)$. Integrate in $T$, and apply
\eqref{chIV:eq:loop-tight} with its error smaller than $\delta\eta/8$.
The remaining event $e_N(s)>1$ has probability $o(1)$ by
Theorem~\ref{chIV:thm:defects}. This proves \eqref{chIV:eq:mass-tight}.

For parent counts, take the canonical representatives of any
$\lceil B\sqrt L\rceil$ parents, $B\ge1$, and delete those with
$d_1(x_D,r)<3R_N$. By \eqref{chIV:eq:parent-packing}, at least
$\lceil B\sqrt L\rceil-49$ remain. Distinct remaining points $x,y$ satisfy
\[
 \min\{d(x,r),d(y,r),d(x,y)\}\ge R_N/\sqrt2.
\]
The Green bound \eqref{chI:eq:green-bounds} gives
\[
 |G_N(x,r)|+|G_N(y,r)|+|G_N(x,y)|\le CL^{1/8},
 \qquad v_N=gL+O(1).
\]
For the local covariance matrix $\Sigma_{x,y}$ of $(h^r(x),h^r(y))$,
the identity $h^r(z)=\phi(z)-\phi(r)$ yields
\begin{equation}\label{chIV:eq:parent-pinned-covariance}
 \begin{aligned}
 \operatorname{Var}h^r(x)&=2v_N-2G_N(x,r),\\
 \operatorname{Cov}(h^r(x),h^r(y))
 &=v_N+G_N(x,y)-G_N(x,r)-G_N(y,r),\\
 \Sigma_{x,y}&=v_N\begin{pmatrix}2&1\\1&2\end{pmatrix}+O(L^{1/8}),\\
 \det\Sigma_{x,y}&=3v_N^2+O(L^{9/8})\asymp L^2.
 \end{aligned}
\end{equation}
The matrix error is uniform entrywise. The one- and two-dimensional
Gaussian densities therefore give
\[
 \frac c{\sqrt L}\le P_h(|h^r(x)|\le1)\le\frac C{\sqrt L},\qquad
 P_h(|h^r(x)|\le1,|h^r(y)|\le1)\le\frac CL.
\]
The contact count over these points has mean at least $cB$ and
second moment at most $C(B+B^2)$, for all sufficiently large $N$.
Paley--Zygmund gives at least $cB$ contacts with a fixed positive
probability. Fix an energy cutoff by the same trace bound as above,
subtract its tail, and then choose the count cutoff $B$. The image
values at these contacts are at most $2u+1$.
Lemma~\ref{chIV:lem:shell-count} and coarea prove the loop count bound.
The same calculation applies to any soft subset that is $R_N$-separated in $d_1$,
so it also bounds its maximal packing number. For an ordinary
prefix choose $B\sqrt L$ separated soft points. The closing
bridge hits an expected fraction at most $C\delta+o(1)$, apart
from a bounded number near its endpoints. With a fixed positive
probability at least half survive. Integrate over the closing
time and use the loop packing bound. Every canonical family is
a separated packing, proving the remaining assertion.
\end{proof}

\subsection*{4.6. Time comparisons and fixed closed endpoints}
To prove \eqref{chIV:eq:endpoint-hit}, apply the strong Markov property
to the occupation time at $z$ upon hitting the filling ball.
Equation~\eqref{chI:eq:occupation}, the Green bounds
\eqref{chI:eq:green-bounds}, and the remainder estimate in the proof
of Lemma~\ref{chI:lem:mixing} give
\begin{equation}\label{chIV:eq:ordinary-hit}
 \PP_x(H_{B_{d_1}(z,2\rho_N)}\le T)
 \le C\left(T/b_N+(L^{1/6}+1)/L\right)
\end{equation}
for $d_1(x,z)\ge R_N/4$; enlarging the ball by its bounded boundary
layer changes only the constant. The denominator is at least
$gL-C\log(2+\rho_N)-C\ge cL$, and the numerator is bounded by
$(T+n)/n+C(1+L^{1/6})$. Split the hit at $T/2$. The heat kernel
on the other half is at most $C/n$; reverse the second half and
apply the same argument from $y$. This proves
\eqref{chIV:eq:endpoint-hit} with its endpoint factor retained.

\begin{theorem}\label{chIV:thm:oscillation}
For every $W\ge0$ and $\eta>0$,
\begin{equation}\label{chIV:eq:oscillation}
 \lim_{\delta\downarrow0}\limsup_{N\to\infty}
 \sup_r\sup_{|c|\le W}\PP_r\left(
   \sup_{|c'-c|\le\delta}|\widehat A_N(c')-\widehat A_N(c)|>\eta\right)=0.
\end{equation}
For each fixed $c_0$ and $\eta>0$,
\begin{equation}\label{chIV:eq:fixed-loop}
 \lim_{\delta\downarrow0}\limsup_{N\to\infty}
 \sup_{r,\,|c-c_0|\le\delta}\PP_{s_N(c)}^{r,r}
       (|\widehat A_N(c)-\widehat A_N(c_0-3\delta)|>\eta)=0.
\end{equation}
The earlier value in \eqref{chIV:eq:fixed-loop} is read from the same
closed path, with its original root. Both assertions extend to continuous $f:\mathbb T^2\to[0,1]$,
when $\widehat A_N(c)$ is replaced by $\widehat A_N(c;f)$, where
\[
 \begin{gathered}
 A_N(t;f)=\frac1{\sqrt L}
 \sum_{a\in\operatorname{Rep}_N(U_1(t))}
 f(a/N)\,a_{D_a(U_1(t))}^G(\ell_t;r,1/4),\\
 \widehat A_N(c;f)=A_N(s_N(c);f).
 \end{gathered}
\]
Here $a/N$ denotes the corresponding point of the unit torus.
\end{theorem}
\begin{proof}
For \eqref{chIV:eq:oscillation}, put $s=s_N(c-\delta)$ and
$t=s_N(c+\delta)$, and restrict to $A_N(s)\le B$.
Conditional on $\mathcal F_s$, the subsequent walk starts at $X_s$.
An early parent $D$ is far when $d_1(a_D,X_s)\ge R_N/4$;
the activity $a_D^G(\ell_s;r,1/4)$ retains the fixed root $r=X_0$.
Equation~\eqref{chIV:eq:ordinary-hit} and
Lemma~\ref{chIV:lem:envelope} give, on $\{A_N(s)\le B\}$,
\[
 \begin{aligned}
 &\EE_r\left[
 \frac1{\sqrt L}\sum_{D\in\mathcal P_N(U_1(s))}
 a_D^G(\ell_s;r,1/4)
 \one_{\{\exists u\in[s,t]:X_u\in B_{d_1}(a_D,2\rho_N)\}}
 \,\middle|\,\mathcal F_s\right]
 \\ &\qquad\le CB\left(2\delta+(L^{1/6}+1)/L\right)+CL^{-1/3}.
 \end{aligned}
\]
The far-parent bound is applied from $X_s$.
The other representatives satisfy $d_1(a_D,X_s)<R_N/4$.
Their normalized activity is at most $CL^{-1/3}$,
uniformly in both $r$ and $X_s$.
A hit on any shorter subinterval is included in this one hit event. Apply \eqref{chIV:eq:matching} to compare every time
with the earliest time, then use Markov's inequality, the mass
cap from Theorem~\ref{chIV:thm:tight}, and the whole-window defect
bound. Choose $B$ before $\delta$, then let $N$ grow.

For \eqref{chIV:eq:fixed-loop}, put $s=s_N(c_0-3\delta)$ and
$t=s_N(c)$. Then $2\delta b_N\le t-s\le4\delta b_N$. Relative
to the ordinary path rooted at $r$, the closed-path prefix has
density
\begin{equation}\label{chIV:eq:prefix-density}
 p_{t-s}(X_s,r)/p_t(r,r)=1+o(1),
\end{equation}
uniformly in the endpoints for every fixed $\delta>0$.
Thus ordinary mass caps and defect bounds transfer to this
prefix. Conditional on the prefix, divide the endpoint-weighted
bound \eqref{chIV:eq:endpoint-hit} by $p_{t-s}(X_s,r)\sim1/n$.
For representatives satisfying
\[
 \min\{d_1(a_D,X_s),d_1(a_D,r)\}\ge R_N/4,
\]
the far hit budget is
at most $CB(2\delta+(L^{1/6}+1)/L)$. Lemma~\ref{chIV:lem:envelope}
bounds the two remaining contributions, centered at $X_s$ and $r$,
by $O(L^{-1/3})$. Matching and Markov prove the assertion, for
every sufficiently small fixed $\delta$. For a spatial test the
additional term is at most
$B\sup_{d_1(x,y)\le\rho_N}|f(x/N)-f(y/N)|=o(1)$.
\end{proof}

There is also a vanishing-length comparison useful for scalar
statistics. With the fixed $C_{\rm mix}$ chosen in the proof of
Theorem~\ref{chIV:thm:defects}, couple the ordinary and closed-path prefixes
of length $t-C_{\rm mix}nL^{1/8}$, where $C_{\rm mix}nL^{1/8}=o(b_N)$, using
\eqref{chIV:eq:prefix-density}, and complete them by an ordinary segment
and a bridge respectively. On a parent-count cap $B\sqrt L$,
\eqref{chIV:eq:ordinary-hit} and Lemma~\ref{chIV:lem:envelope} bound each
expected normalized activity loss by
\[
 CBL^{-2/3}+CL^{-1/3}=o(1).
\]
The corresponding normalized zero-parent loss is
$CBL^{-5/6}+CL^{-1/2}$. Theorem~\ref{chIV:thm:defects} gives
$e_N(t)=o_P(1)$ uniformly on compact central time intervals.
For the ordinary path $\omega$ and closed path $\widetilde\omega$
in the coupling above, for every fixed $W<\infty$ and
$\varepsilon>0$,
\[
 \begin{aligned}
 \lim_{N\to\infty}\sup_{r,\,|c|\le W}\PP\bigl(
 &|A_N(t)(\omega)-A_N(t)(\widetilde\omega)|\\
 &+|C_N(t)(\omega)-C_N(t)(\widetilde\omega)|>\varepsilon\bigr)=0,
 \qquad t=s_N(c).
 \end{aligned}
\]
Here $\PP$ is the coupling law for each $(r,c)$.

We record the remaining finite-measure form needed for marking.
Set
\[
 \mathfrak m_{N,W}(F)=b_N^{-1}\sum_r
   \int_{s_N-Wb_N}^{s_N+Wb_N}\EE_r[\one_{\{X_t=r\}}F]\,dt.
\]
Its mass is at most $4W$ eventually. Theorem~\ref{chIV:thm:tight} gives
simultaneous parent-count and activity caps under this measure.
The loop defect assertion and bounded truncation give, for fixed $C$,
\begin{equation}\label{chIV:eq:truncated-defect}
 \mathfrak m_{N,W}\left[1\wedge
       C|D_N^{\rm wide}(U_u(t))|/\sqrt L\right]\longrightarrow0.
\end{equation}
Indeed split at $|D_N^{\rm wide}|/\sqrt L=\varepsilon$, obtaining
$4WC\varepsilon$ plus the exceptional measure, and then let
$\varepsilon\downarrow0$.

For the mean-zero field $\phi=\phi_N$, the standardized root coordinate is
\[
 z=\phi_N(r)/\sqrt{G_N(r,r)}.
\]
With the source path held fixed, let
$q_N(r,t,z)$ be the conditional Gaussian probability, given this coordinate,
of a pinned square contact $2\ell_t(x)+h^r(x)^2\le v$ at a vertex with
$d_1(x,r)\le5R_N$; its dependence on the path is suppressed.
For fixed $v$, \hyperlink{proof.appendix.Q}{Appendix~Q}.2 proves
\begin{equation}\label{chIV:eq:near-contact}
 \mathfrak m_{N,W}\left[\int q_N(r,t,z)\,\gamma_1(dz)\right]
 \longrightarrow0,
\end{equation}
where $\gamma_1$ is standard Gaussian measure. The same bound
holds when its integral is restricted to any root-coordinate
window. Equation~\eqref{chIV:eq:near-contact} is the contact estimate
with the root coordinate integrated against $\gamma_1$.

\clearpage
\endgroup

\begingroup
\mergedchaptersetup{V}{The limiting law of local-time activity}
\begin{center}
{\large\scshape Chapter V\par}\vspace{6pt}
{\LARGE The limiting law of local-time activity\par}
\end{center}
\vspace{6pt}

The proof proceeds through
\[
 \eqref{chV:eq:source-limit}\ \Longrightarrow\
 \eqref{chV:eq:master}\ \Longrightarrow\
 \eqref{chV:eq:half-line}\ \Longrightarrow\
 \eqref{chV:eq:average}.
\]
The energy shift is $R(G)=(G^2-1)/2$. The time comparison
\eqref{chV:eq:uniform-time-removal} then gives
\eqref{chV:eq:small-laplace}; tightness and uniqueness yield
\eqref{chV:eq:main}.

The proof uses the estimates of Chapters~I--IV. Chapter~IV supplies
ordinary and loop-window activity tightness and the time comparisons
recalled in Section~5.5.

\subsection*{5.1. Scalar tests, root coordinates and truncation}

Retain $V_N=(\ZZ/N\ZZ)^2$, the rate-one walk, and its local times.
Write $n=N^2$, $L=\log N$, $g=2/\pi$, $b_N=gnL$ and
$v_N=G_N(0,0)=gL+O(1)$. Let $\phi=\phi_N$ be the mean-zero
Gaussian field of Chapter~II. The time coordinates are
\begin{equation}\label{chV:eq:times}
 s_N=b_N\left(2L-\frac32\log L\right)
       +\frac{9g}{32}n(\log L)^2,
 \qquad s_N(c)=s_N+cb_N.
\end{equation}
Use $U^\circ$, $\mathcal P_N(U)$ and $D_a(U)$ from Section~3.1,
with $\rho_N=\tfrac12e^{L^{1/6}}$ and $R_N=Ne^{-2L^{1/8}}$.
Lemma~\ref{chIII:lem:geometry} gives, for large $N$,
\[
 \begin{gathered}
 \diam_{d_1}D_a(U)\le\rho_N,\qquad
 \dist_{d_1}(D_a(U),D_{a'}(U))\ge R_N\quad(a\ne a'),\\
 n_N(U)=|\mathcal P_N(U)|.
 \end{gathered}
\]

The activity is the same root-dependent functional as before:
\begin{equation}\label{chV:eq:activity}
 \widehat A_N(c)=\frac1{\sqrt L}\sum_{a\in\operatorname{Rep}_N(U_1(s_N(c)))}
 a^G_{D_a(U_1(s_N(c)))}(\ell_{s_N(c)};r,1/4),
 \qquad U_1(t)=\{x:\ell_t(x)\le1\}.
\end{equation}
Here $r=X_0$ remains the root at all observation times. Write
$\phi^{(r)}=\phi_N^{(r)}$ for the regression residual of Chapter~II:
$\phi^{(r)}(x)=\phi(x)-G_N(x,r)\phi(r)/v_N$. Explicitly,
\begin{equation}\label{chV:eq:cluster-activity}
 a_D^G(\ell;r,v)=\frac1{\sqrt{2\pi g}}
       \EE|\mathcal I_D(\ell,\phi^{(r)};v)|.
\end{equation}
The expectation in this definition is only over the Gaussian field.
The resulting $\widehat A_N(c)$ is a measurable nonnegative path functional.

Let $Z=\mathcal Z(\mathbb T^2)$ be the critical mass in
\eqref{chII:eq:critical-mass}, and let $\Gamma>0$ be the shell coefficient
in \eqref{chII:eq:Gamma-definition}. Thus $0<Z<\infty$
almost surely, and we put
\begin{equation}\label{chV:eq:target}
 \Lambda(a)=\EE e^{-a\Gamma Z},\qquad a>0.
\end{equation}
Since the Laplace integrand takes values in $[0,1]$, $\Lambda$ is
continuous by bounded convergence.

For $s_N(c)>0$, define the fixed-root Laplace tests
\[
 F_N^r(c,a)=\EE_r[e^{-ae^c\widehat A_N(c)}\mid X_{s_N(c)}=r],\qquad
 J_N^r(c,a)=\EE_r e^{-ae^c\widehat A_N(c)}.
\]
Their uniform root averages are
\begin{align}
 F_N(c,a)&=\frac1n\sum_{r\in V_N}F_N^r(c,a),\label{chV:eq:F}\\
 J_N(c,a)&=\frac1n\sum_{r\in V_N}J_N^r(c,a).
 \label{chV:eq:J}
\end{align}
Set $F_N^r(c,a)=F_N(c,a)=1$ when $s_N(c)\le0$. For every fixed
central window, $s_N(c)>0$ for all sufficiently large $N$. The initial distribution in \eqref{chV:eq:J}
is the uniform stationary law $\pi_N$.

\begin{theorem}\label{chV:thm:main}
Under $\PP_{\pi_N}$,
\begin{equation}\label{chV:eq:main}
 \widehat A_N(0)\ \Longrightarrow\ \Gamma Z.
\end{equation}
Consequently, for every $a>0$,
$\EE_{\pi_N}e^{-a\widehat A_N(0)}\longrightarrow\Lambda(a)$.
The constant $\Gamma$ is exactly the one fixed in Chapter~II.
\end{theorem}

The following boundedness and Lipschitz estimates yield the common
subsequence used in Section~5.3. For $a,b\ge d>0$,
\begin{equation}\label{chV:eq:amplitude-bound}
 0\le F_N(c,a)\le1,\qquad
 |F_N(c,a)-F_N(c,b)|\le\frac{|a-b|}{ed}.
\end{equation}
Indeed $\sup_{x\ge0}xe^{-dx}=1/(ed)$, applied to
$x=e^c\widehat A_N(c)$. The same bounds hold for $\Lambda$ and for $J_N$.
We shall also use
\begin{equation}\label{chV:eq:log-bound}
 |e^{-ae^cx}-e^{-ae^{c'}x}|\le e^{-1}|c-c'|,
 \qquad x\ge0,
\end{equation}
obtained by differentiating with respect to the logarithm of the
amplitude.

For clarity we record the precise probabilistic inputs used below.
For $W\ge0$, define the finite, unnormalized loop-window measure by
\begin{equation}\label{chV:eq:window}
 \mathfrak m_{N,W}(H)=\frac1{b_N}\sum_r
 \int_{s_N-Wb_N}^{s_N+Wb_N}
 \EE_r[\ind_{\{X_t=r\}}H(r,t,\omega)]\,dt.
\end{equation}
Here and below $N$ is large enough that the lower endpoint is
positive. The heat-kernel estimates give
$\mathfrak m_{N,W}(1)\le4W$ eventually. At physical time $t>0$,
the activity is $A_N(t)$ and the number of soft parents is
$n_N(U_1(t))$. Thus $A_N(s_N(c))=\widehat A_N(c)$ on a path
rooted at $r$.
Theorem~\ref{chIV:thm:tight} supplies, for every fixed
$W$ and $\varepsilon>0$, finite $B,H>0$ such that
\begin{equation}\label{chV:eq:caps}
 \limsup_N\mathfrak m_{N,W}
  \{n_N(U_1(t))>B\sqrt L\text{ or }A_N(t)>H\}\le\varepsilon.
\end{equation}
In Appendix~R, on $m\le B\sqrt L$, the covariance error is bounded by
\eqref{chV:eq:hs-decay}. On $\sum_i\alpha_i\le H$, the contact and
product errors are bounded by \eqref{chV:eq:sum-error} and
\eqref{chV:eq:bernoulli-error}.

Retain the wide-defect set $D_N^{\rm wide}(U)$ of
\eqref{chIV:eq:wide-def}. Its Euclidean thresholds are
$\tfrac12 e^{L^{1/8}}$ and $2Ne^{-L^{1/8}}$; in particular
$D_N^{\rm can}(U)\subset D_N^{\rm wide}(U)$ for large $N$.
For each fixed $C,W<\infty$, its integrated consequence is
\begin{equation}\label{chV:eq:defect-input}
 \mathfrak m_{N,W}\left[
   \min\left\{1,\frac{C|D_N^{\rm wide}(U_1(t))|}{\sqrt L}\right\}\right]
 \longrightarrow0.
\end{equation}

To state the remaining contact input correctly, distinguish the
regressed field in \eqref{chV:eq:cluster-activity} from
$h^r=\phi-\phi(r)$. For $S_r=\phi(r)/\sqrt{v_N}$, conditional
Gaussian regression gives
\begin{equation}\label{chV:eq:conditional-field}
 h^{r,z}(x)=\phi^{(r)}(x)
   +\left(\frac{G_N(x,r)}{\sqrt{v_N}}-\sqrt{v_N}\right)z,
 \qquad h^r\mid\{S_r=z\}\ \overset{\rm law}=\ h^{r,z}.
\end{equation}
The residual field is independent of the standard normal $S_r$.
With $r,z$ and the source path fixed, let $q_N(r,t,z)$ be the
Gaussian probability, taken only over the residual field
$\phi^{(r)}$, of
\begin{equation}\label{chV:eq:near-event}
 \exists x:\ d_1(x,r)\le5R_N,\qquad
             2\ell_t(x)+h^{r,z}(x)^2\le1/4.
\end{equation}
The dependence of $q_N(r,t,z)$ on the source path is suppressed.
With $\gamma_1=\mathcal N(0,1)$ denoting the standard Gaussian
probability measure on $\mathbb R$, the near-root input is
\begin{equation}\label{chV:eq:near-input}
 \mathfrak m_{N,W}\left[\int q_N(r,t,z)\,\gamma_1(dz)\right]
 \longrightarrow0.
\end{equation}
Accordingly, \eqref{chV:eq:near-input} is the $\gamma_1$-average
of the conditional Gaussian contact kernel $q_N(r,t,z)$ defined from
\eqref{chV:eq:conditional-field}.

\subsection*{5.2. From the Gaussian shell to the master identity}

Recall the square-profile test of Section~2.6:
\begin{equation}\label{chV:eq:square-test}
 \Psi_{N,\lambda}(W)=(1-\lambda)^{n_N(\{x:W(x)\le1/4\})},
 \qquad 0<\lambda<1,
\end{equation}
This is the canonical component-count Laplace test with constant
test function $-\log(1-\lambda)$, whose marking weight is $\lambda$.

We recall the source, including its normalization. Set
\begin{align}
 \mathcal{Q}_N&=n^{-1}\|\phi\|_2^2,
 &B_N(\theta)&=2(s_N+nv_N+\theta b_N)/n,\notag\\
 \alpha_N(\theta,\phi)&=\sqrt{(B_N(\theta)-\mathcal{Q}_N)_+},
 &w_N^\sigma(\phi,a)&=\frac{\sum_r(\phi(r)+\sigma a)_+}{n^2a}.
 \label{chV:eq:shell-data}
\end{align}
Here the square norm uses counting measure, and $\EE$ in the shell
formulas integrates $\phi_N$. The quotient defining $w_N^\sigma$
is used for $a>0$; set $w_N^\sigma(\phi,a)=0$ for $a\le0$.
With $\alpha_N=\alpha_N(\theta,\phi)$ and both signs summed,
\begin{equation}\label{chV:eq:source}
 S_N(\theta,\lambda)=n\EE\left[
 \ind_{\{\mathcal{Q}_N<B_N(\theta)\}}\sum_{\sigma=\pm1}
 w_N^\sigma(\phi,\alpha_N)
 \Psi_{N,\lambda}((\phi+\sigma\alpha_N)^2)\right].
\end{equation}
The entire integrand is zero on $\{\mathcal Q_N\ge B_N(\theta)\}$.
On fixed central $\theta$-windows, the field--branch notation of
Section~3.3 gives, for all sufficiently large $N$,
\[
 S_N(\theta,\lambda)
 =n\varrho_N(s)\,E_{\nu_{N,s}}\Psi_{N,\lambda}(W_s^\sigma(\phi)),
 \qquad s=s_N+nv_N+\theta b_N.
\]
The averaged shell theorem, Theorem~\ref{chII:thm:shell-source}, states
\begin{equation}\label{chV:eq:source-limit}
 \int\chi(\theta)S_N(\theta,\lambda)\,d\theta
 \longrightarrow
 \int\chi(\theta)\Lambda(\lambda e^{-(\theta+1/2)})\,d\theta,
 \qquad \chi\in C_c(\mathbb R).
\end{equation}
We use this averaged statement throughout.

The energy coordinate produced by coarea is
\begin{equation}\label{chV:eq:energy}
 U_N^r=\frac{\|h^r\|_2^2/2-nv_N}{b_N}.
\end{equation}
Since $\sum_x\phi(x)=0$, $\|h^r\|_2^2=\|\phi\|_2^2+n\phi(r)^2$.
Consequently, with $R(z)=(z^2-1)/2$,
\begin{align}
 U_N^r&=R(S_r)+\epsilon_N^r,\label{chV:eq:energy-decomposition}\\
 \epsilon_N^r&=\frac{\mathcal{Q}_N-v_N}{2gL}
       +\left(\frac{v_N}{gL}-1\right)R(S_r),
 &\sup_r\EE|\epsilon_N^r|&=O(L^{-1}).\label{chV:eq:energy-error}
\end{align}
The estimate follows from $\EE|\mathcal{Q}_N-v_N|=O(1)$ and the Green
asymptotic. In the joint source integral, this error is controlled
by the $L^1$ translation bound \eqref{chV:eq:translation-error}.

\begin{proposition}\label{chV:prop:master}
Let $\mu_{\rm en}$ be the law of $(G^2-1)/2$ for a standard
normal $G$. For $0<\lambda<1$ and $\chi\in C_c(\mathbb R)$,
\begin{multline}\label{chV:eq:master}
 \int\chi(\theta)\int
 F_N(\theta-r,\lambda e^{-(\theta+1/2)})\,
              \mu_{\rm en}(dr)\,d\theta\\
 \longrightarrow\int\chi(\theta)
                 \Lambda(\lambda e^{-(\theta+1/2)})\,d\theta.
\end{multline}
\end{proposition}
\begin{proof}
Initially take $\chi\ge0$. Apply the exact coarea identity of
Chapter~III to the square test and the shell-time test
$\chi((s-s_N-nv_N)/b_N)$. Sum over roots and divide by $b_N$.
The factor $n$ in
\eqref{chV:eq:source} then gives the exact equality
\begin{multline}\label{chV:eq:exact-source}
 \int\chi(\theta)S_N(\theta,\lambda)\,d\theta\\
 =\frac1{b_N}\sum_r\int_0^\infty
 \EE_r\left[\ind_{\{X_t=r\}}\EE_\phi\left[
 \chi((t-s_N)/b_N+U_N^r)
 \Psi_{N,\lambda}(2\ell_t+(h^r)^2)\right]\right]dt.
\end{multline}
For each fixed root $r$, $\EE_\phi$ integrates the Gaussian law
of $\phi_N$, independently of the path under $\PP_r$.
The time argument is $c+U_N^r$, with $c=(t-s_N)/b_N$.

First discard $0<t<b_N$. The short-time estimate proved in
\hyperlink{proof.appendix.S}{Appendix~S}.1, specifically \eqref{chV:eq:pinned-tail} and
\eqref{chV:eq:initial-mass}, makes its contribution $o(1)$. On $t\ge b_N$, replace
$U_N^r$ by $R(S_r)$ only in $\chi$. The absolute error is bounded by
\begin{equation}\label{chV:eq:translation-error}
 2\sup_r\EE\|\chi(\,\cdot+U_N^r)
                   -\chi(\,\cdot+R(S_r))\|_1=o(1).
\end{equation}
Here $np_t(r,r)\le2$ for $t\ge b_N$; translation continuity in
$L^1$ and \eqref{chV:eq:energy-error} prove the last assertion.

Disintegrate with respect to $S_r=z$. In the next estimate,
the expectation is only over the Gaussian residual $\phi^{(r)}$,
with $r,z$ and the path held fixed.
For every fixed $B,H,M>0$ and $\varepsilon>0$, the
deterministic marking estimate in \hyperlink{proof.appendix.R}{Appendix~R} gives, whenever
$n_N(U_1(t))\le B\sqrt L$, $A_N(t)\le H$, and $|z|\le M$,
\begin{equation}\label{chV:eq:marking}
 \left|\EE\Psi_{N,\lambda}(2\ell_t+(h^{r,z})^2)-e^{-\lambda e^{-z^2/2}A_N(t)}\right|
 \le \varepsilon+2q_N(r,t,z)
       +\min\{1,C_M|D_N^{\rm wide}(U_1(t))|/\sqrt L\}
\end{equation}
for $N\ge N_0(B,H,M,\varepsilon)$, uniformly in $0<\lambda<1$,
with a constant $C_M=C(M)$.

Fix $M$ and suppose $\supp\chi\subset[-H_0,H_0]$. Put
$W=H_0+(M^2+1)/2$. Then
\[
 |z|\le M,\quad\chi(c+R(z))\ne0\quad\Longrightarrow\quad |c|\le W.
\]
After integration under \eqref{chV:eq:window} and $\gamma_1(dz)$,
the last two terms in \eqref{chV:eq:marking} are bounded by
\[
 \begin{aligned}
 &2\|\chi\|_\infty\mathfrak m_{N,W}
       \left[\int q_N(r,t,z)\,\gamma_1(dz)\right]\longrightarrow0,\\
 &\|\chi\|_\infty\mathfrak m_{N,W}
       \left[\min\{1,C_M|D_N^{\rm wide}(U_1(t))|/\sqrt L\}\right]
       \longrightarrow0,
 \end{aligned}
\]
by \eqref{chV:eq:near-input} and \eqref{chV:eq:defect-input}.
The complement of the caps contributes at most
\[
 \|\chi\|_\infty\mathfrak m_{N,W}
   \{n_N(U_1(t))>B\sqrt L\text{ or }A_N(t)>H\}.
\]
For this fixed window choose the caps using \eqref{chV:eq:caps} and fix
the marking tolerance. Let $N\to\infty$, then decrease the cap error
and marking tolerance. Outside $|z|\le M$ integrate first in $c$;
the contribution is at most a constant times
$\|\chi\|_1\gamma_1\{|z|>M\}$. Thus the marking error tends to zero,
first with $M$ fixed and then as $M\to\infty$.

The expression that remains is
\begin{equation}\label{chV:eq:marked-loop}
 \frac1n\sum_r\int_{s_N(c)\ge b_N}
 np_{s_N(c)}(r,r)\int\chi(c+R(z))
 \EE_r[e^{-\lambda e^{-z^2/2}\widehat A_N(c)}\mid X_{s_N(c)}=r]
 \,\gamma_1(dz)\,dc.
\end{equation}
Uniformly on all $t\ge b_N$, $np_t(r,r)=1+o(1)$. Since
the integral of $|\chi(c+R(z))|$ in $c$ is $\|\chi\|_1$,
removing this endpoint factor costs $o(1)$. Extend the $c$-integral
to the whole line. The support condition gives
\[
 c\le1-s_N/b_N,\quad\chi(c+R(z))\ne0
 \quad\Longrightarrow\quad R(z)\ge s_N/b_N-H_0-1.
\]
Thus the absolute contribution from the omitted range is at most
\[
 \|\chi\|_1\gamma_1\{z:R(z)\ge s_N/b_N-H_0-1\}\longrightarrow0.
\]

Finally put $\theta=c+R(z)$. The amplitude identity is
\begin{equation}\label{chV:eq:amplitude-identity}
 (\lambda e^{-(\theta+1/2)})e^c=\lambda e^{-z^2/2}.
\end{equation}
Thus the extended integral is exactly the left side of
\eqref{chV:eq:master}. Equation~\eqref{chV:eq:source-limit} proves the
claim for nonnegative $\chi$. For signed $\chi$, apply the same
argument to $\chi_+$ and $\chi_-$, both continuous with compact
support, and subtract. Every final kernel is bounded by one, so
these signed integrals are finite.
\end{proof}

\subsection*{5.3. A common subsequence of bounded tests}

We apply weak-star compactness to the bounded deterministic functions
$F_N$.

\begin{lemma}\label{chV:lem:subsequence}
Every sequence $N_j\to\infty$ has a subsequence, still denoted by
$N_j$, and a jointly measurable $f:\mathbb R\times(0,\infty)
\to[0,1]$ such that, for every $a>0$ and $\psi\in L^1(\mathbb R)$,
\begin{equation}\label{chV:eq:weak-star}
 \int\psi(c)F_{N_j}(c,a)\,dc\longrightarrow
 \int\psi(c)f(c,a)\,dc.
\end{equation}
The same jointly measurable version $f$ satisfies
\[
 |f(c,a)-f(c,b)|\le\frac{|a-b|}{ed}
 \qquad(c\in\mathbb R,\ d>0,\ a,b\ge d).
\]
\end{lemma}
\begin{proof}
For each positive rational $a$ and integer $m\ge1$, the functions
$F_{N_j}(\cdot,a)$ lie in the unit ball of $L^\infty([-m,m])$.
This ball is weak-star sequentially compact because $L^1([-m,m])$
is separable. Diagonal extraction over these two countable indices
gives consistent measurable limits for all positive rational
amplitudes. Consistency on overlapping intervals follows by testing
against functions supported on the overlap.

Positivity, the upper bound one, and each rational instance of
\eqref{chV:eq:amplitude-bound} pass to the limits: if a claimed
inequality failed on a set of positive finite measure, its
indicator would contradict the weak-star inequality. Remove the
union of the exceptional sets for all rational pairs and rational
$d>0$. Outside this single null set, the rational-amplitude
function extends uniquely and locally Lipschitz continuously to
all positive amplitudes. On the null set define $f(c,a)=0$ for
every $a$. Rational approximation gives joint measurability.

For fixed real $a>0$, approximate it by rationals in a compact
subinterval of $(0,\infty)$. The error in either integral is at
most $\|\psi\|_1|a-b|/(ed)$, uniformly in $j$. This proves
\eqref{chV:eq:weak-star} for compactly supported $\psi$. For general
$\psi\in L^1$, first restrict to $[-m,m]$; the two tails together
cost at most $2\int_{|c|>m}|\psi(c)|\,dc$. Let $m\to\infty$.
\end{proof}

\begin{lemma}\label{chV:lem:variable}
The limit in Lemma~\ref{chV:lem:subsequence} satisfies, for every
$0<\lambda<1$ and $\chi\in C_c(\mathbb R)$,
\begin{multline}\label{chV:eq:limit-master}
 \int\chi(\theta)\int f(\theta-r,\lambda e^{-(\theta+1/2)})
                    \,\mu_{\rm en}(dr)\,d\theta\\
 =\int\chi(\theta)\Lambda(\lambda e^{-(\theta+1/2)})\,d\theta.
\end{multline}
\end{lemma}
\begin{proof}
On the compact support of $\chi$, the continuous amplitude
$a(\theta)=\lambda e^{-(\theta+1/2)}$ stays in some $[d,D]$
with $d>0$. On $\supp\chi$, approximate it uniformly by a simple
function $a_{\rm simple}$ taking finitely many positive rational
values $a_i\ge d/2$ on measurable sets $E_i$. For a sufficiently
fine mesh, \eqref{chV:eq:amplitude-bound} bounds the error in each
integral of the passage to the limit by
\[
 \frac{2\|\chi\|_1}{ed}
 \|a-a_{\rm simple}\|_{L^\infty(\supp\chi)}.
\]
For a single simple-function term, Fubini gives the test
\[
 \psi_i(c)=\int\chi(c+r)\ind_{E_i}(c+r)\,\mu_{\rm en}(dr),
 \qquad \|\psi_i\|_1\le\|\chi\ind_{E_i}\|_1.
\]
Apply Lemma~\ref{chV:lem:subsequence} on its common subsequence to
these finitely many $L^1$ tests and rational amplitudes $a_i$. Sum
over $i$, then refine the amplitude mesh. Proposition~\ref{chV:prop:master}
gives \eqref{chV:eq:limit-master}.
\end{proof}

\subsection*{5.4. One-sided Gamma inversion}

The shift by $1/2$ in the shell source matches the support of the
energy law. Indeed $R(G)+1/2=G^2/2$ has density
\begin{equation}\label{chV:eq:gamma-density}
 k(u)=\frac{e^{-u}}{\sqrt{\pi u}}\ind_{\{u>0\}}.
\end{equation}
The change of variables from the two Gaussian half-lines gives
this density with total mass one. Hence, for bounded measurable $v$,
\begin{equation}\label{chV:eq:energy-convolution}
 \int v(\theta-r)\,\mu_{\rm en}(dr)
       =\int_0^\infty k(u)v(\theta+1/2-u)\,du.
\end{equation}
With $x=\theta+1/2$, the arguments on the right are $x-u<x$
for $u>0$; thus the convolution uses only $v|_{(-\infty,x)}$.

\begin{proposition}\label{chV:prop:identification}
For the common subsequential limit $f$ and every $a>0$,
\begin{equation}\label{chV:eq:half-line}
 f(c,a)=\Lambda(a)\quad\text{for almost every }c<-\log a.
\end{equation}
\end{proposition}
\begin{proof}
Put
\[
 H(x,a)=\int_0^\infty k(u)f(x-u,a)\,du.
\]
Extend $k$ by zero on $(-\infty,0]$. Since $|f|\le1$ and
$\|k\|_1=1$, for $a,b\ge d>0$,
\[
 \begin{aligned}
 |H(x+h,a)-H(x,a)|&\le\|k(\cdot+h)-k\|_1\longrightarrow0
       \quad(h\to0),\\
 |H(x,a)-H(x,b)|&\le |a-b|/(ed).
 \end{aligned}
\]
Hence $H$ is jointly continuous on $\mathbb R\times(0,\infty)$.

With $x=\theta+1/2$, \eqref{chV:eq:limit-master} says that for each
$\lambda\in(0,1)$ the continuous function
\[
 H(x,\lambda e^{-x})-\Lambda(\lambda e^{-x})
\]
has integral zero against every compactly supported continuous
test. It is therefore zero everywhere. Fix $a>0$ and
$x<-\log a$. Choosing $\lambda=ae^x\in(0,1)$ gives
\begin{equation}\label{chV:eq:gamma-constant}
 H(x,a)=\Lambda(a),\qquad x<-\log a.
\end{equation}
This step explains the half-line restriction.

Let $v(x)=f(x,a)-\Lambda(a)$. Then $v$ is bounded and
$k*v=0$ on $(-\infty,-\log a)$. Direct integration yields
\begin{equation}\label{chV:eq:double-gamma}
 (k*k)(u)=\frac{e^{-u}}\pi\int_0^u\frac{ds}{\sqrt{s(u-s)}}
             =e^{-u},\qquad u>0.
\end{equation}
For example, $s=u\sin^2 t$ turns the integral into $\pi$.
For $x<-\log a$, all arguments $x-u$ in a second convolution remain
below $-\log a$. Fubini is justified by boundedness of $v$ and
$\|k\|_1=1$, so
\[
 0=(k*(k*v))(x)=\int_0^\infty e^{-u}v(x-u)\,du
              =e^{-x}\int_{-\infty}^x e^yv(y)\,dy.
\]
The last integral is locally absolutely continuous and its
derivative is $e^xv(x)$ almost everywhere. Since it vanishes
throughout $(-\infty,-\log a)$, $v=0$ almost everywhere there.
\end{proof}

\begin{corollary}\label{chV:cor:average}
Along the common subsequence, for $a>0$, $\delta>0$ and
$c+\delta<-\log a$,
\begin{equation}\label{chV:eq:average}
 \frac1{2\delta}\int_{c-\delta}^{c+\delta}F_{N_j}(x,a)\,dx
                 \longrightarrow\Lambda(a).
\end{equation}
\end{corollary}
\begin{proof}
Use $(2\delta)^{-1}\ind_{[c-\delta,c+\delta]}$ as the $L^1$
test in \eqref{chV:eq:weak-star}, then apply
Proposition~\ref{chV:prop:identification}.
\end{proof}

The conclusion gives convergence of interval averages on one common
subsequence; Section~5.5 supplies the time comparison needed at $c=0$.

\subsection*{5.5. Recovering the deterministic observation time}

We recall the time-comparison inputs in a form suited to this
step. Fix a central coordinate $c_0$ and $\eta>0$, and for
$\delta>0$ put $b=c_0-3\delta$. Theorem~\ref{chIV:thm:oscillation} proves
\begin{align}
 &\lim_{\delta\downarrow0}\limsup_N\sup_r
 \PP_r\{|\widehat A_N(c_0)-\widehat A_N(b)|>\eta\}=0,\label{chV:eq:ordinary-time}\\
 &\lim_{\delta\downarrow0}\limsup_N
 \sup_{\substack{r\in V_N\\|c-c_0|\le\delta}}
 \PP_{s_N(c)}^{r,r}\{|\widehat A_N(c)-\widehat A_N(b)|>\eta\}=0.
 \label{chV:eq:bridge-time}
\end{align}
In the second line $\widehat A_N(b)$ is evaluated on the prefix of the
same closed path, with the same root $r$. The available estimates
allow the errors to be made small for all sufficiently small
fixed $\delta$, with the subsequent threshold in $N$ depending
on $\delta$.

The bridge comparison is the scalar case of
Theorem~\ref{chIV:thm:oscillation}; its proof uses
\eqref{chIV:eq:prefix-density} and the endpoint-weighted hit estimate
\eqref{chIV:eq:endpoint-hit}, with the return factor retained before
normalizing the bridge law.

\begin{lemma}\label{chV:lem:prefix}
For fixed $c_0$ and $\delta>0$, the law of the prefix up to
$s_N(b)$ under $\PP_{s_N(c)}^{r,r}$ has density $1+o(1)$ relative
to its ordinary law under $\PP_r$, uniformly over
$r\in V_N$, $|c-c_0|\le\delta$ and the prefix endpoint.
\end{lemma}
\begin{proof}
The Markov property gives the exact density on this prefix:
\begin{equation}\label{chV:eq:prefix-density}
 \frac{p_{s_N(c)-s_N(b)}(X_{s_N(b)},r)}{p_{s_N(c)}(r,r)}.
\end{equation}
The numerator time is between $2\delta b_N$ and
$4\delta b_N$. Since $b_N/n=gL\to\infty$, the heat-kernel
mixing estimate of Chapter~I gives $np_t(x,r)=1+o(1)$ uniformly
in $x,r$ and in these times. It also applies to the denominator.
Taking the ratio gives
\[
 M_N:=\sup_{\substack{r,x\in V_N\\|c-c_0|\le\delta}}
 \left|\frac{p_{s_N(c)-s_N(b)}(x,r)}{p_{s_N(c)}(r,r)}-1\right|
 \longrightarrow0.
\]
Every measurable $[0,1]$-valued prefix test therefore has its
expectation changed by at most $M_N$.
\end{proof}

\begin{proposition}\label{chV:prop:time-removal}
For every fixed $a>0$ and $c_0\in\mathbb R$,
\begin{equation}\label{chV:eq:uniform-time-removal}
 \lim_{\delta\downarrow0}\limsup_N
 \sup_{|c-c_0|\le\delta}|F_N(c,a)-J_N(c_0,a)|=0.
\end{equation}
In particular the same conclusion holds after replacing $F_N(c,a)$
by its symmetric average over $[c_0-\delta,c_0+\delta]$.
\end{proposition}
\begin{proof}
At a fixed root $r$, the successive intermediate expectations are
\[
 \begin{aligned}
 &F_N^r(c,a),\\
 &\EE_r[e^{-ae^c\widehat A_N(b)}\mid X_{s_N(c)}=r],\\
 &\EE_r e^{-ae^c\widehat A_N(b)},\\
 &\EE_r e^{-ae^{c_0}\widehat A_N(b)},\\
 &J_N^r(c_0,a).
 \end{aligned}
\]
The bridge expectations use the same rooted path and its prefix
up to $s_N(b)$.
For $q\ge0$ and nonnegative $x,y$, $|x-y|\le\eta$ implies
$|e^{-qx}-e^{-qy}|\le q\eta$; on the complementary event the
difference is at most one. If the failure probabilities in
\eqref{chV:eq:ordinary-time} and \eqref{chV:eq:bridge-time} are at most
$E$, and the prefix density differs from one by at most $M_N$,
the four comparisons give
\begin{equation}\label{chV:eq:four-errors}
 |F_N^r(c,a)-J_N^r(c_0,a)|
 \le a(e^c+e^{c_0})\eta+2E+M_N+e^{-1}|c-c_0|.
\end{equation}
The fixed-root tests are those defined in Section~5.1. The last
term is \eqref{chV:eq:log-bound}.

Average over roots and take the supremum in $c$. With $\delta$
fixed, Lemma~\ref{chV:lem:prefix} gives $M_N\to0$. For any fixed
$\eta$, the two time inputs make $E$ vanish as $\delta\downarrow0$
after the upper limit in $N$. The remaining upper bound is
$2ae^{c_0}\eta$. Let $\eta\downarrow0$. Integrating the same
uniform bound proves the assertion about averages.
\end{proof}

\begin{corollary}\label{chV:cor:small-laplace}
For every $0<a<1$,
\begin{equation}\label{chV:eq:small-laplace}
 J_N(0,a)=\EE_{\pi_N}e^{-a\widehat A_N(0)}\longrightarrow\Lambda(a).
\end{equation}
\end{corollary}
\begin{proof}
Start with an arbitrary subsequence $N_j\to\infty$ and extract
the common subsequence of Lemma~\ref{chV:lem:subsequence}. For any
$0<\delta<-\log a$, Corollary~\ref{chV:cor:average} gives
\[
 \frac1{2\delta}\int_{-\delta}^{\delta}F_{N_j}(c,a)\,dc
                      \longrightarrow\Lambda(a).
\]
The triangle inequality and the preceding limit give
\[
 \begin{aligned}
 \limsup_j|J_{N_j}(0,a)-\Lambda(a)|
 &\le\limsup_j\sup_{|c|\le\delta}|J_{N_j}(0,a)-F_{N_j}(c,a)|\\
 &\longrightarrow0\qquad(\delta\downarrow0),
 \end{aligned}
\]
by Proposition~\ref{chV:prop:time-removal}. Thus the extracted
subsequence converges to $\Lambda(a)$.
If the full sequence failed to converge, a subsequence staying
a positive distance from $\Lambda(a)$ would contradict this
conclusion. This proves \eqref{chV:eq:small-laplace}.
\end{proof}

\subsection*{5.6. Uniqueness of the scalar probability law}

Let $\mu_N^A$ be the law of $\widehat A_N(0)$ under $\PP_{\pi_N}$.
The ordinary activity tightness from Chapter~IV gives, for every
$\varepsilon>0$, an $H<\infty$ with
\begin{equation}\label{chV:eq:law-tight}
 \limsup_N\mu_N^A((H,\infty))\le\varepsilon.
\end{equation}
Since all variables are nonnegative, this is tightness on
$\mathbb R$. The finitely many omitted initial indices can be
absorbed by increasing $H$.

\begin{lemma}\label{chV:lem:laplace-unique}
Two probability measures on $[0,\infty)$ having the same Laplace
transform on $(0,1)$ are equal.
\end{lemma}
\begin{proof}
For any finite measure $\nu$ on $[0,\infty)$, its Laplace
transform is real analytic at every $a_0>0$. In fact for
$|h|<a_0$ the exponential series may be integrated against
$e^{-a_0x}\nu(dx)$, since its series of absolute values is
$e^{-(a_0-|h|)x}\le1$. On smaller neighborhoods this yields
the convergent Taylor expansion. Thus equality on $(0,1)$
extends to all of the connected interval $(0,\infty)$ by
the identity theorem for real analytic functions.

Push the two measures forward by $x\mapsto e^{-x}$, regarding
their images as measures on $[0,1]$. Equality of the Laplace
transforms at every positive integer gives equality of all
positive integer moments; equality of total masses supplies
the zeroth moment. Polynomial approximation of continuous
functions on $[0,1]$ shows that the image measures coincide.
The map is a Borel isomorphism from $[0,\infty)$ to $(0,1]$,
so the original measures coincide as well.
\end{proof}

\begin{proof}[Proof of Theorem~\ref{chV:thm:main}]
Take any subsequence. By \eqref{chV:eq:law-tight} it has a further
subsequence converging weakly to a probability measure $\mu^A$.
The closed set $[0,\infty)$ has full measure under every
$\mu_N^A$, and hence also under $\mu^A$. For $a>0$ use the bounded
continuous function $x\mapsto e^{-a\max\{x,0\}}$ on $\mathbb R$.
Corollary~\ref{chV:cor:small-laplace} implies
\[
 \int e^{-ax}\mu^A(dx)=\Lambda(a),\qquad 0<a<1.
\]
Lemma~\ref{chV:lem:laplace-unique} identifies $\mu^A$ with the law of
$\Gamma Z$. Every subsequential limit is therefore the same,
and tightness gives \eqref{chV:eq:main}. Testing this weak convergence
with the same bounded continuous function proves Laplace
convergence for every $a>0$.
\end{proof}

This completes the first path-side distribution limit. The next
chapters compare the activity with the number of unvisited
components; the coefficient in that comparison is separate
from $\Gamma$.

\clearpage
\endgroup

\begingroup
\mergedchaptersetup{VI}{Fine structure and local conditional laws}
\begin{center}
{\large\scshape Chapter VI\par}\vspace{6pt}
{\LARGE Fine structure and local conditional laws\par}
\end{center}
\let\mathscr\mathcal
The activity of a soft component and the indicator that it contains an
unvisited point are two functions of the same local occupation profile.
We compare their means under a common local law, then show that their
centered sum is small. The comparison is made on actual path components:
the local regions must cover the retained components, false regions must
be removed, and the activity of the omitted components must be controlled.
The resulting residual bound supplies Chapter~VII with the local
comparison needed to identify the limiting coefficient and establish
its positivity.

\subsection*{6.1. Exact local kernels and rare soft events}
\label{chVI:sec:local}
Retain the rate-one walk on $V_N=(\ZZ/N\ZZ)^2$, its occupation
$\ell_t$, and
\[
 n=N^2,\quad L=\log N,\quad g=2/\pi,\quad b_N=gnL,\quad
 s_N=b_N(2L-\tfrac32\log L)+\tfrac{9g}{32}n(\log L)^2.
\]
Write $r$ for the fixed initial vertex, with the path convention of
Chapter~III. Fix $W\ge0$ and $\delta>0$. Throughout, a
physical window $I=[l_N,r_N]$ satisfies, eventually,
\begin{equation}\label{chVI:eq:window}
 s_N-Wb_N\le l_N\le r_N\le s_N+Wb_N,\qquad
 \delta b_N\le l_N,\qquad \delta b_N\le r_N-l_N.
\end{equation}
Write $\Qw$ for the root-averaged loop law in
\eqref{chIII:eq:window-law}, and $\mathsf Q_{N,r}^I$ for the same law
conditioned on its initial vertex $r$. Both laws include the random
observation time $t\in I$. When an independent uniform grid origin
$o\in V_N$ is used, the law is
$\Qaug=\Qw\otimes\operatorname{Unif}(V_N)$.
All limits in this chapter hold for every sequence of admissible windows.

For this chapter and Appendices~T--Y, let $d_2=d$ be the periodic
Euclidean distance and retain $d_1$ for graph distance. Put
\[
 B(y,R):=\{x\in V_N:d_2(x,y)<R\},\qquad
 \mathbb B(y,q):=\{x\in V_N:d_2(x,y)\le\lfloor e^q\rfloor\}.
\]
Thus $B$ is the open lattice disk used for auxiliary local kernels,
whereas $\mathbb B$ is the closed disk with integer radius used for
the original count profile. In a planar chart the same definitions
apply in $\ZZ^2$. For a lattice domain $D$, its one-step
boundaries are
\[
 \partial_iD:=\{x\in D:\exists y\notin D,\ x\sim y\},\qquad
 \partial_oD:=\{y\notin D:\exists x\in D,\ x\sim y\},
\]
where $x\sim y$ means nearest neighbors. Write
$\diam_j D=\max_{x,y\in D}d_j(x,y)$, $j=1,2$, for nonempty $D$.

Here are the upstream inputs, in the precise forms used below.
The Green estimates of Chapter~I, together with
Proposition~\ref{chIV:prop:local} and
\hyperlink{proof.appendix.L}{Appendix~L}.1, supply the endpoint-uniform
estimates for the local record kernels.
Chapter~II supplies the Dirichlet maximum window, the one- and two-point
Gaussian contact estimates, and the relative contact comparison. The additional radial
repulsion estimates and energy-truncated probes needed for fine
clustering are established in \hyperlink{proof.appendix.Y}{Appendix~Y}.
Chapter~III supplies exact product disintegration of complete local paths,
the certified-selection rule and its two-copy physical-window transfer,
and the energy-truncated coarea identity.
Chapter~IV supplies, uniformly in these windows, tightness of
$A_N$ and of the normalized soft packing number, and the vanishing
normalized wide-defect count. The soft one-point estimate in
\hyperlink{proof.appendix.Q}{Appendix~Q}.1 gives
$\mathsf Q_{N,r}^I(\ell_t(x)\le u)\le C_uN^{-2}L^{3/2}$.
For the first-layer estimate, Appendix~X also uses
\cite[Proposition 5.8]{chVI:LS26} with a central second-layer count
$O(k^2/\mathfrak h)<e^{c\mathfrak h}$.

An endpoint word is an ordered list $w=((a_j,b_j))_{j=1}^m$ of
entrance and exit vertices. For fixed domains $(D_j)$ it is feasible
when $h_{D_j,b_j}(a_j)>0$ for every $j$, and its rate-one product kernel is
\[
 P_{m,w}^{(1)}:=\bigotimes_{j=1}^m\mathsf B_{D_j}^{a_j,b_j}.
\]
When the domains and word are fixed by the current stratum, write
$P_m^{(1)}$ and $E_m^{(1)}$ for this law and its expectation. Holding
times and durations remain part of each local record.
For the same domains, endpoints and number of paths, define
\[
 (\Theta_g\omega)(s):=\omega(gs),\qquad
 P_m^{(g)}=(\Theta_g^{\otimes m})_*P_m^{(1)}.
\]
For each record $j$, its duration satisfies
\[
 V_j^{(g)}=\frac{V_j^{(1)}}g.
\]
Thus a physical duration cutoff $V_j^{(1)}\le d_N$ becomes
$V_j^{(g)}\le d_N/g$. Unless a superscript is displayed, $P_m,E_m$
in this chapter denote the physical rate-one kernels; Appendix~T uses
the rate-$g$ kernels and is transferred through the deterministic
pushforward above. For $m=0$ we use the empty-product convention of Section~3.5.

Let $L^{\rm loc}(x)$ be the sum of the occupation times at $x$ over
these local paths, in rate-one units. Local kernel estimates use
$\lambda(x)=L^{\rm loc}(x)/g$. Under $\Theta_g$ this is the
occupation field of a planar walk of total rate $g$, with Green
diagonal $q+O(1)$ in a disk of radius $e^q$ and Ray--Knight Gaussian
covariance one half of that Green kernel. Marks and activities use
the original rate-one units.

\begin{proposition}[An all-count soft denominator]\label{chVI:prop:denominator}
Fix $u>0$ and $0<\gamma_-\le\gamma_+<1/2$. There are
$q_0,C<\infty$, depending only on these parameters, such that
the following holds for every $q\ge q_0$ and
$q^{\gamma_-}\le h\le q^{\gamma_+}$. Consider $m$ exact-exit
paths from the one-step inner boundary of $B(y,e^{q+h})$ to the
one-step outer boundary of $B(y,e^{q+2h})$. On the torus, assume
that this outer disk and its one-step boundary lie in one planar
chart. Put
\[
 S=\{\min_{B(y,e^q)}\lambda\le u\},\qquad
 w=\sqrt{hm},\qquad s=w-\sqrt2q.
\]
Uniformly in every feasible endpoint word and every $m\ge0$,
\begin{equation}\label{chVI:eq:denominator}
 P_m(S)\ge
 \exp\{-2\sqrt2s_+-s_+^2/q-Ch(1+w^2/q^2)\}.
\end{equation}
Let $C_*>0$ be fixed and define the ordered soft-pair count
\[
 H_2:=\#\{(x,x')\in B(y,e^q)^2:
 \lambda(x),\lambda(x')\le u,\ d_2(x,x')\ge e^{q-C_*h}\}.
\]
Then, for $\gamma_+<\nu\le1$,
\begin{equation}\label{chVI:eq:pairs}
 \sup_{w\ge\sqrt2q+q^\nu}E_m(H_2\mid S)\le e^{-cq^\nu}.
\end{equation}
For the second estimate the constants may depend on
$(u,\gamma_-,\gamma_+,C_*,\nu)$, and the bound is uniform over $m$,
feasible endpoint words, centers and admissible ambient tori.
For fixed $C_0,c_0>0$, the estimates hold uniformly for the core,
entrance and exit logarithmic radii
\[
 q+\delta_K,\qquad q+h+d+\delta_A,\qquad q+2h+d+\delta_D,
 \qquad d\in[0,1],\quad
 \max_{\bullet\in\{K,A,D\}}|\delta_\bullet|\le C_0e^{-c_0h}.
\]
The constants may also depend on $C_0,c_0$; the annular width is
$h+\delta_D-\delta_A=h+O(e^{-c_0h})$.
\end{proposition}
\begin{proof}
\hyperlink{proof.appendix.T}{Appendix~T} proves a critical-count seed by a one-dimensional coarse
Gaussian shift and transfers it to each prescribed word by a relative
likelihood estimate. Appending paths that avoid a selected soft point
gives \eqref{chVI:eq:denominator} above the critical count; monotonicity gives
the lower-count case. The two-point numerator retains its negative
quadratic term before division by this denominator. This yields
\eqref{chVI:eq:pairs} also at arbitrarily large counts.
\end{proof}

For the same exact-exit domains and endpoint word, let $\omega_{q,y}$
be the exit distribution from $\{x:d_2(x,y)<e^{q+1}\}$ for a walk
started at $y$, and let $\Sigma_{q,y}$ be its support.
For record $j$, write $L^{(j)}$
for its rate-one occupation profile and define
\[
 Y_j:=g^{-1}\sum_{\xi\in\Sigma_{q,y}}
       \omega_{q,y}(\xi)L^{(j)}(\xi),\qquad
 Y:=\sum_{j=1}^mY_j,\qquad \mu_m:=E_mY.
\]
Thus $Y=\sum_\xi\omega_{q,y}(\xi)\lambda(\xi)$. Each $Y_j$
depends only on the $j$th coordinate of the product law $P_m$.
Kac's formula and the endpoint Green bound give
\begin{equation}\label{chVI:eq:bernstein}
 P_m(|Y-\mu_m|>z)
 \le 2\exp\{-c\min(z^2/(mh^2),z/h)\},\qquad z\ge0.
\end{equation}
For $m\ge1$, the single-record kernels give, uniformly in
$j$ and its prescribed endpoints,
\[
 E_m\exp\{\theta(Y_j-E_mY_j)\}\le\exp(Ch^2\theta^2),
 \qquad |\theta|\le c/h.
\]
Independence under $P_m$ permits multiplication of these bounds. Optimize
\[
 \theta=\min\{cz/(mh^2),c/h\}.
\]
The harmonic mean is
$mh(1+O(e^{-ch}))$ in the local units. Thus the corresponding
square-root deviations have exponent $-c\Delta^2/h$ in their
admissible range. In particular, for every $m\ge1$,
\[
 v\ge\sqrt{mh},\quad \Delta>0,\quad
 \Delta\ge C\sqrt{mh}\,e^{-ch}
 \quad\Longrightarrow\quad
 P_m(\sqrt Y>v+\Delta)\le2e^{-c\Delta^2/h}.
\]
The one-sided Bernstein calculation in Appendix~T after
\eqref{chVI:eq:R17.16} proves this throughout the displayed range,
including the small-count branch.

\subsection*{6.2. Retained parents and omitted activity}
\label{chVI:sec:retention}
Write $S_t=U_1(t)$ and recall its canonical component family
$\mathcal P_N(U_1(t))$ from Chapter~III. We call these
components \emph{parents}; $a_D$ is the canonical representative of
$D\in\mathcal P_N(U_1(t))$. The connection and separation scales are
\[
 \rho_N=\tfrac12e^{L^{1/6}},\qquad R_N=Ne^{-2L^{1/8}}.
\]
Thus $\diam_1D\le\rho_N$ and distinct parents are $R_N$-separated.
Set $c_g=(2\pi g)^{-1/2}$ and $v=1/4$. Use the Gaussian interval-union
activity $a_D^G$ of Chapter~III, evaluated under the root-conditioned
residual law $\phi_N^{(r)}$. For each parent write
\begin{equation}\label{chVI:eq:parentmarks}
 J_D=\one_{\{D\cap\{\ell_t=0\}\ne\varnothing\}},\quad
 a_D^G(t)=a_D^G(\ell_t;r,v),\quad
 A_N(t)=L^{-1/2}\sum_{D\in\mathcal P_N(U_1(t))}a_D^G(t).
\end{equation}
In particular, with $B_*=(2c_g\sqrt v)^{-1}$,
\begin{equation}\label{chVI:eq:loweractivity}
 0\le J_D\le B_*a_D^G(t).
\end{equation}
An exact-zero point supplies an interval of length $2\sqrt v$ in
every realization of the Gaussian residual, which proves this bound.

The following deterministic atlas makes retention measurable from
the soft-set geometry and the fixed root.
Put $D_N^{\rm mac}=Ne^{-L^{1/8}}$ and
$s=\min\{N,\max(1,\lfloor D_N^{\rm mac}/1000\rfloor)\}$.
The centers are the coordinate products of
$\{js:j\in\ZZ,\ 0\le j<\lfloor N/s\rfloor\}$ modulo $N$.
Write $d_{\rm circ}(a,b)=\min_{j\in\ZZ}|a-b+jN|$ for the
distance between coordinate residues. The core and its two buffers are
\[
 \begin{aligned}
 P_z&:=\{x:\max_{j=1,2}d_{\rm circ}(x_j,z_j)<2s\},\\
 P_z^{\rm buf}&:=\{x:\max_{j=1,2}d_{\rm circ}(x_j,z_j)<3s\},\\
 P_z^{\rm ext}&:=\{x:\max_{j=1,2}d_{\rm circ}(x_j,z_j)<4s\}.
 \end{aligned}
\]
These cores cover the torus. Fix a finite coloring of the centers
such that distinct centers of the same color have $d_2$-distance at
least $20s$; the number of colors is bounded independently of $N$.
In particular, their extended buffers are disjoint. Fix also an
ordering of the centers and set
\[
 z_*(x):=\min\{z:x\in P_z\}
\]
in that ordering. We call $z_*(x)$ the \emph{stable owner} of $x$:
it is a deterministic choice of one covering core.
For a center $z$, its \emph{soft cloud} is $S_t\cap P_z^{\rm buf}$;
the buffer is active when this set is nonempty.

For $K\in\mathbb N$, define
\[
 D_K=\lceil10+1024e^K\rceil,\qquad
 Q_K=[-D_K,D_K]^2\cap\ZZ^2.
\]
Let $\mathcal B_{N,K}(r)$ be the union of the sets
\[
 \begin{gathered}
 \{x:|d_{\rm circ}(x_j,z_j)-2s|\le2D_K\}
       \quad(z\text{ a core center},\ j=1,2),\\
 \{x:d_{\rm circ}(x_j,0)\le2D_K\}\quad(j=1,2),\qquad
 \{x:|d_1(x,r)-D_N^{\rm mac}|\le2D_K\}.
 \end{gathered}
\]
These are the core-boundary, seam, and root-distance bands.
A parent $D$ is retained when
\begin{equation}\label{chVI:eq:retained}
\begin{gathered}
 a_D\notin\mathcal B_{N,K}(r),\qquad
 d_1(a_D,r)>D_N^{\rm mac}+2D_K,\qquad D\subset a_D+Q_K,\\
 D\subset P_{z_*(a_D)},\qquad
 S_t\cap P_{z_*(a_D)}^{\rm buf}=D.
\end{gathered}
\end{equation}
Write $\mathcal P_{N,K}^{\rm ret}=\mathcal P_{N,K}^{\rm ret}(t;r)$
for this family. First define this family from the path and root;
independently, sample the grid origin used in the covering construction
below. Define its
omitted activity by
\begin{equation}\label{chVI:eq:omitted}
 O_{N,K}=L^{-1/2}\sum_{D\in\mathcal P_N(U_1(t))\setminus\mathcal P_{N,K}^{\rm ret}}a_D^G(t).
\end{equation}

\phantomsection\label{chVI:def:intrinsic}
For the covering disks used below, fix finitely many square charts
$U_\alpha\subset V_N$, fixed planar lifts of these charts to $\ZZ^2$,
and $\delta_0>0$ such that
\[
 V_N=\bigcup_\alpha
 \{x\in U_\alpha:d_2(x,U_\alpha^c)\ge\delta_0N\}.
\]
The charts, their lifts and $\delta_0$ are deterministic. Each lift is
injective on every local domain and one-step boundary used below and
preserves nearest-neighbor adjacency there. For $q\ge0$ and an origin
$a\in V_N$, put
\[
 \begin{split}
 s_q&:=\lfloor e^q\rfloor,\qquad
 \mathcal G_{N,q}(a):=
 \{a+s_q(i,j)\pmod N:i,j\in\ZZ,\ 0\le i,j\le\lfloor N/s_q\rfloor\},\\
 \gamma_{\rm loc}&:=1/8192,\qquad
 D_{q,y}:=B(y,e^{q+2q^{\gamma_{\rm loc}}-4e^{-q^{\gamma_{\rm loc}}}}),\\
 \mathscr A_{N,q}(a)&:=\{(\alpha,y):y\in\mathcal G_{N,q}(a),\
 d_2(y,U_\alpha^c)\ge\delta_0N/2,\
 D_{q,y}\cup\partial_oD_{q,y}\subset U_\alpha\}.
 \end{split}
\]
Members of $\mathscr A_{N,q}(a)$ are the available locators; their
core disks are $B(y,e^q)$. Repeated grid points modulo $N$ occur
only once. Choose a fixed sufficiently large integer $K_{\rm loc}$
as in Lemma~\ref{chVI:lem:intrinsic-geometry}. For a cloud
$C_z:=S_t\cap P_z^{\rm buf}$ with $S_t\cap P_z\ne\varnothing$,
define
\[
 q_z:=\inf\{q\in\ZZ:q\ge K_{\rm loc},\
 \exists(\alpha,y)\in\mathscr A_{N,q}(0),\ C_z\subset B(y,e^q)\},
 \qquad\inf\varnothing:=+\infty.
\]
This is its intrinsic scale. With the fixed deterministic order
$\prec$, for finite $q_z$ bind
\[
 (\alpha_z,y_z):=\min_{\prec}
 \{(\alpha,y)\in\mathscr A_{N,q_z}(0):C_z\subset B(y,e^{q_z})\}.
\]
The same construction applies with any prescribed origin $a$.

\begin{proposition}[Omitted activity]\label{chVI:prop:omitted}
For every $o,\rho>0$ there is $K_0\ge1$ such that for each
$K\ge K_0$, for all sufficiently large $N$,
\begin{equation}\label{chVI:eq:omittedprob}
 \Qw(O_{N,K}>o)\le\rho.
\end{equation}
The same order of limits makes
$E_{\Qw}[\min(1,O_{N,K})]$ tend to zero. The normalized number
of omitted zero-bearing parents is at most $B_*O_{N,K}$.
\end{proposition}
\begin{proof}
In this proof define the local exception
\[
 \mathcal E_{N,K}:=
 \{S_t\cap\mathcal B_{N,K}(r)\ne\varnothing\}
 \cup\{\exists x\in S_t:d_1(x,r)\le D_N^{\rm mac}+2D_K\}.
\]
The two regions have volumes at most $C_Kn/s+C_KN$ and
$C(D_N^{\rm mac}+2D_K)^2$. The one-point soft estimate gives
\[
 \Qw(\mathcal E_{N,K})
 \le C_KL^{3/2}\left(s^{-1}+N^{-1}
             +\frac{(D_N^{\rm mac}+2D_K)^2}{n}\right)
 =o_N(1)\qquad(K\text{ fixed}).
\]
Sample a mean-zero Gaussian field $\phi$ with covariance $G_N$,
independently of the root and path, and put $h=\phi-\phi(r)$.
Write $P_h$ for its law given the root, with the path held fixed.
For the scales $q_z$ just defined, \hyperlink{proof.appendix.U}{Appendix~U} proves that
$\{2\ell_t+h^2\le v\}$ meets $\bigcup_{z:q_z\ge K}(S_t\cap P_z)$
with probability tending to zero under this path--field law,
first as $N\to\infty$, then as $K\to\infty$;
see \eqref{chVI:eq:R18.3}.

On $\mathcal E_{N,K}^c$, for fixed $K$ and sufficiently large $N$,
there is a fixed $c_0>0$, depending only on the deterministic cover, such that
\begin{equation}\label{chVI:eq:omitted-core-inclusion}
 \bigcup_{D\in\mathcal P_N(U_1(t))\setminus\mathcal P_{N,K}^{\rm ret}}D
 \subset\bigcup_{z:q_z\ge K-c_0}(S_t\cap P_z).
\end{equation}
If an omitted parent is too large for $a_D+Q_K$, use its diameter
and a core containing its $3R_N$-neighborhood.
Otherwise use its stable owner. If that owner's whole cloud fit in
a disk of radius $e^{K-c_0}$, its diameter would be less than
$\rho_N/2$ for large $N$. A defective point would then have a witness
in the same buffer, contradicting this diameter bound. All cloud
points would be clean and connected to $D$, contradicting nonretention.
This argument includes extra buffer points that are themselves defects.

Still on $\mathcal E_{N,K}^c$, for fixed $K$ and sufficiently large $N$,
\[
 d_2(D,r)\ge\frac{D_N^{\rm mac}+2D_K}{\sqrt2}-\rho_N\ge R_N,
 \qquad
 \frac{D_N^{\rm mac}}{R_N}=e^{L^{1/8}}\longrightarrow\infty.
\]
All omitted representatives have $d_1(a_D,r)>3R_N$. Thus
Lemma~\ref{chIV:lem:parent-colouring} partitions them into at most
$49$ classes $\mathcal D_j$, without removing a parent, with
\[
 d_2(D,D')\ge3R_N/\sqrt2-2\rho_N\ge R_N
 \quad(D\ne D'\in\mathcal D_j),\qquad
 \diam_2D\le\rho_N<R_N^{\rm loc}.
\]
Set $u_j=L^{-1/2}\sum_{D\in\mathcal D_j}a_D^G(t)$ and
$u_*=\max_{1\le j\le49}u_j$, inserting empty classes if necessary.
Then $O_{N,K}\le49u_*$. Apply
\eqref{chII:eq:contact-amplification} to a maximizing class:
\begin{equation}\label{chVI:eq:amplify}
 \begin{split}
 \frac{O_{N,K}}{1+O_{N,K}}
 &\le49\frac{u_*}{1+u_*}\\
 &\le C\,P_h\bigl(\{2\ell_t+h^2\le v\}
       \cap\!\bigcup_{D\in\mathcal P_N(U_1(t))\setminus\mathcal P_{N,K}^{\rm ret}}
       D\ne\varnothing\mid\text{path}\bigr),
 \qquad\text{on }\mathcal E_{N,K}^c.
 \end{split}
\end{equation}
The activity $a_D^G(t)$ uses the root-conditioned residual
$\phi_N^{(r)}$ of Chapter~III; the relative contact comparison
is precisely what relates these two Gaussian objects.
Average \eqref{chVI:eq:amplify} on $\mathcal E_{N,K}^c$ and use
$\min(1,O_{N,K})\le1$ on $\mathcal E_{N,K}$.
Equations~\eqref{chVI:eq:omitted-core-inclusion} and
\eqref{chVI:eq:R18.3}, with cutoff $K-c_0$, give
\[
 \lim_{K\to\infty}\limsup_{N\to\infty}
 E_{\Qw}[\min(1,O_{N,K})]=0,
\]
with exceptional contribution at most $\Qw(\mathcal E_{N,K})=o_N(1)$.
Here $\min(1,x)\le2x/(1+x)$; the probability assertion follows from
$\one_{\{x>o\}}\le(1+o)x/[o(1+x)]$.
The last assertion follows from \eqref{chVI:eq:loweractivity}.
\end{proof}

\subsection*{6.3. A common local law and its two marks}
\label{chVI:sec:common}
Fix $K$ and $0<\epsilon\le1$. Choose $k$ so that the quantities below
satisfy $\mathfrak h\ge16$ and $T\ge2$, and subsequently let
$N\to\infty$. Use
\begin{equation}\label{chVI:eq:exponents}
 \gamma=1/20,\qquad \eta=1/12,\qquad \eta'=1/24,\qquad
 \mathfrak h=k^\gamma,\quad T=\lfloor(k^2-k)/\mathfrak h\rfloor.
\end{equation}
Layers are indexed by $i=0,\ldots,T-1$. The five additional
layers $i=T,\ldots,T+4$ are called \emph{outer guards}; their
crossing-count conditions are retained as exterior data when paths
inside layer $T-1$ are resampled. Put $q_i=k+i\mathfrak h$ and set
\[
 B_i^-(x)=\mathbb B(x,q_i),\qquad
 B_i^+(x)=\mathbb B(x,q_i+\mathfrak h-4e^{-\mathfrak h}).
\]
For every layer $i$ set $\tau_{i,0}=0$ and recursively
\[
 \sigma_{i,j}=\inf\{s\ge\tau_{i,j-1}:X_s\in B_i^-(x)\},\qquad
 \tau_{i,j}=\inf\{s\ge\sigma_{i,j}:X_s\notin B_i^+(x)\},
 \qquad j\ge1,
\]
with $\inf\varnothing=\infty$. If $X_0\in B_i^-(x)$ then
$\sigma_{i,1}=0$. A complete traversal, or \emph{record}, is the
stopped segment on $[\sigma_{i,j},\tau_{i,j}]$ with its duration.
Only records with $\tau_{i,j}\le t$ are counted; an unfinished terminal
segment belongs to the exterior/censored future in the record
decomposition. Thus
\[
 m_i(x,t)=\#\{j:\tau_{i,j}\le t\},\qquad
 w_i(x)=\sqrt{\mathfrak h m_i(x,t)}.
\]
The \emph{count corridor} with parameter $\zeta$ is the closed condition
\begin{equation}\label{chVI:eq:corridor}
 (q_i+\mathfrak h)^{1/2-\zeta}\le w_i(x)-\sqrt2q_i
                  \le(q_i+\mathfrak h)^{1/2+\zeta}.
\end{equation}
Use $\eta'$ at $i=0,T-1$ and $\eta$ at the other layers, including
the five outer guards. Both endpoints of every integer count interval
in \eqref{chVI:eq:corridor} are included. Inner counts enter the
local corridor flag; terminal and outer-guard counts enter the exterior data.

For every real $q$, define
\[
 Q(q):=\{z\in\ZZ^2:-\lfloor e^q\rfloor\le2z_j<\lfloor e^q\rfloor,\
 j=1,2\}.
\]
A field on $Q(k)$ is regarded as a function on $\ZZ^2$ by extension
by zero. For a nonnegative such field $F$, define
$E(F):=\{z\in Q(k):F(z)\le1\}$ and, when $E(F)\ne\varnothing$,
$b(F):=\min_{\rm lex}E(F)$. The \emph{anchor} is thus the first
soft vertex in lexicographic order. The local event is
\begin{equation}\label{chVI:eq:localflag}
 E_0=\{E(F)\ne\varnothing,\ b(F)\in Q(k-\epsilon),\
 E(F)\subset b(F)+Q_K\subset Q(k-\epsilon/2)\}.
\end{equation}
On $E_0$, the centered mark is
$e=(D,f):=(E(F)-b(F),F(b(F)+\cdot))$.
In particular $0\in D\subset Q_K$, $f\ge0$, and
$D=\{z\in Q(k)-b(F):f(z)\le1\}$ on this event. Thus the
soft-set definition has domain $Q(k)-b(F)$; zero extension specifies
the ambient field coordinate of the mark.
Use the mark space
$\{D:D\subset Q_K\}\times[0,\infty)^{\ZZ^2}$ with the discrete
and product Borel sigma-algebras. The empty mark is $(\varnothing,0)$.
For any such mark define
\begin{equation}\label{chVI:eq:localmarks}
\begin{gathered} J(e)=\one_{\{\exists z\in D:f(z)=0\}},\\
 H(e)=c_g E_\psi\left|\bigcup_{z\in D}
 [-\psi(z)-\sqrt{(v-2f(z))_+},-\psi(z)+\sqrt{(v-2f(z))_+}]\right|.\end{gathered}
\end{equation}
Here $|\cdot|$ is Lebesgue length, and $E_\psi$ is expectation
over the centered Gaussian field $\psi$, with
$\psi(0)=0$ and
$\operatorname{Cov}(\psi(x),\psi(y))=
a_{\rm lat}(x)+a_{\rm lat}(y)-a_{\rm lat}(x-y)$, where
$a_{\rm lat}$ is the rate-one potential kernel
defined in Chapter~II, Section~2.4. The mark $e$ is
held fixed in this expectation; in particular $J(\varnothing,0)=
H(\varnothing,0)=0$.
The interval-union definition gives
\begin{equation}\label{chVI:eq:markbounds}
 0\le J\le1,\qquad B_*^{-1}J\le H\le B_*^{-1}|Q_K|.
\end{equation}

\phantomsection\label{chVI:def:common-measure}
The reference law is constructed on the planar cell, as follows.
In this paragraph $B_i^\pm=B_i^\pm(0)\subset\ZZ^2$, and
$\mathbb P_x^{\rm pl}$ denotes rate-one planar random walk. With
$x_k=(\lfloor e^{k+\mathfrak h}\rfloor,0)$, define the two endpoint laws
\[
 \begin{aligned}
 \pi_k^-(a)&:=\mathbb P_{x_k}^{\rm pl}
  (X_{H_{B_0^-}}=a\mid H_{B_0^-}<\tau_{B_1^+}),\\
 \pi_k^+(b)&:=\mathbb P_0^{\rm pl}(X_{\tau_{B_0^+}}=b).
 \end{aligned}
\]
Here $H_D$ and $\tau_D$ are first entrance into and first exit from
$D$. Let $\mathcal W_k$ be the feasible words
$w=((a_j,b_j))_{j=1}^m\in(\partial_iB_0^-\times\partial_oB_0^+)^m$
whose counts satisfy the $i=0$ narrow corridor
\eqref{chVI:eq:corridor}. For each word let $\mathcal K_w$ be the
law, under its product of exact-exit kernels in $B_0^+$, of
$F=\sum_{j=1}^m L^{(j)}|_{Q(k)}$. On the Borel space
$\mathscr X_k:=[0,\infty)^{Q(k)}$, set
\[
 \begin{gathered}
 \alpha(w):=\sqrt{\mathfrak h m}-\sqrt2k,\qquad
 p_k(w):=\prod_{j=1}^m\pi_k^-(a_j)\pi_k^+(b_j),\\
 \mathcal R_k:=\sum_{w\in\mathcal W_k}
   p_k(w)\alpha(w)e^{2\sqrt2(\alpha(w)-1)}\mathcal K_w.
 \end{gathered}
\]
The count corridor makes $\mathcal W_k$ finite and $\alpha(w)>0$.
Thus $\mathcal R_k$ is a finite measure. For the centered-mark map
$\Theta(F)=(E(F)-b(F),F(b(F)+\cdot))$ on $E_0$, define
\[
 \mathcal O_{K,k,\epsilon}:=
 \begin{cases}
 \Theta_\#\bigl(\mathcal R_k(\,\cdot\,\cap E_0)/\mathcal R_k(E_0)\bigr),
       &\mathcal R_k(E_0)>0,\\
 \delta_{(\varnothing,0)},&\mathcal R_k(E_0)=0.
 \end{cases}
\]
As above, cell fields are extended by zero when $\Theta$ is applied.
\hyperlink{proof.appendix.V}{Appendix~V} analyzes this finite mixture; \eqref{chVI:eq:reference}
records the same measure.

\phantomsection\label{chVI:def:terminal-measure}
For a feasible terminal word $z$ of length $m(z)$ at layer $T-1$,
let $P_z$ be the product of its exact-exit kernels in $B_{T-1}^+$.
Let $\mathscr C_{\rm int}$ be the intersection of the count
corridors at layers $0,\ldots,T-2$, with exponent $\eta'$ at
layer $0$ and $\eta$ at the other layers. If $F$ is the total
record occupation on $Q(k)$, define, for Borel $B\subset\mathscr X_k$,
\[
 \begin{gathered}
 \mathcal M_z(B):=P_z(F\in B,\mathscr C_{\rm int}),\qquad
 b_z:=\sqrt{\mathfrak h m(z)}-\sqrt2q_{T-1},\\
 g_z:=\sqrt{2/\pi}\frac{b_z}{Tk^{3/2}}
 \exp\{-2(T-1)\mathfrak h+2\sqrt2-2\sqrt2b_z-b_z^2/(T\mathfrak h)\}.
 \end{gathered}
\]
This factor is $\mathcal L_k(1,b_z)$ in
\eqref{chVI:eq:leading}, and is positive on the narrow terminal corridor.

\phantomsection\label{chVI:def:total-variation}
For probability measures $\mu,\nu$ on the same measurable space, use
\[
 d_{\rm TV}(\mu,\nu):=\sup_{B\ {\rm measurable}}|\mu(B)-\nu(B)|.
\]
Thus $|\mu(f)-\nu(f)|\le2\|f\|_\infty d_{\rm TV}(\mu,\nu)$
for every bounded measurable $f$.

\begin{proposition}[Common conditional mark law]\label{chVI:prop:common}
The law $\mathcal O_{K,k,\epsilon}$ is the planar reference law defined
above. Define its coefficient
\begin{equation}\label{chVI:eq:coefficient}
 \beta_{K,k,\epsilon}=
 \frac{\mathcal O_{K,k,\epsilon}(J)}{\mathcal O_{K,k,\epsilon}(H)}
 \quad\text{if }\mathcal O_{K,k,\epsilon}(H)>0,
\end{equation}
On the branch $\mathcal O_{K,k,\epsilon}(H)=0$, set
$\beta_{K,k,\epsilon}=0$. Then
$0\le \beta_{K,k,\epsilon}\le B_*$ and
$\mathcal O_{K,k,\epsilon}(J-\beta_{K,k,\epsilon}H)=0$.
For all sufficiently large $k$, there are errors
$0\le e_k=e(k)<1$ and
$d_{N,K,k,\epsilon}=d(N,K,k,\epsilon)\ge0$,
uniform in terminal words, centers and roots, with the limits stated
below.
The terminal data are the exact count and ordered entrance/exit
word at layer $T-1$, satisfying its narrow corridor. For each such
word, consider the product of first-exit paths conditioned on the
corridors at layers $0,\ldots,T-2$, $E_0$, and individual duration
at most $e^{L/2}$. This conditional law is used only when the
conditioning event has positive product-law probability.
Uniformly over these positive-mass strata, its mark law differs
in total variation from $\mathcal O_{K,k,\epsilon}$ by
\begin{equation}\label{chVI:eq:lawerror}
 \frac{2e_k}{1-e_k}+d_{N,K,k,\epsilon},\qquad
 e_k\longrightarrow0,\quad
 \sqrt L\,d_{N,K,k,\epsilon}\longrightarrow0.
\end{equation}
The first limit is $k\to\infty$; the second is $N\to\infty$
with $K,k,\epsilon$ fixed. 
\end{proposition}
\begin{proof}
\hyperlink{proof.appendix.V}{Appendix~V}, equations \eqref{chVI:eq:wordfactor} and
\eqref{chVI:eq:leading}, proves for the measures just defined
the uniform, unnormalized comparison
\begin{equation}\label{chVI:eq:relative}
 (1-e_k)g_z\mathcal R_k(B)\le\mathcal M_z(B)
          \le(1+e_k)g_z\mathcal R_k(B),\qquad g_z>0,
\end{equation}
for every Borel field event $B\subset\mathscr X_k$. Condition on
$E_0$ and apply $\Theta$. If $\mathcal R_k(E_0)>0$, the two
conditional measures have density ratio between
$(1-e_k)/(1+e_k)$ and its reciprocal. Strata with $\mathcal R_k(E_0)=0$ have zero physical mass by
\eqref{chVI:eq:relative}.

At fixed $k$ there are finitely many feasible terminal words,
because the terminal corridor bounds the number of paths.
On the record space of word $z$, put
\[
 E_z:=\mathscr C_{\rm int}\cap\{F\in E_0\},\qquad
 D_{N,z}:=\bigcap_{j=1}^{m(z)}\{V_j\le e^{L/2}\},\qquad
 \mathcal Z_+:=\{z:P_z(E_z)>0\}.
\]
If $\mathcal Z_+=\varnothing$, set $d_{N,K,k,\epsilon}=0$;
there are no positive-mass conditional laws to compare.
Suppose henceforth that $\mathcal Z_+\ne\varnothing$, and define
\[
 p_*:=\min_{z\in\mathcal Z_+}P_z(E_z)>0.
\]
Exact planar identification at fixed $k$ makes $p_*$ independent
of $N$, centers and roots. The endpoint-conditioned exit tail
and the finite terminal count cap give
\[
 \begin{aligned}
 d_{N,K,k,\epsilon}
 &:=\max_{z\in\mathcal Z_+}P_z(D_{N,z}^c\mid E_z)\\
 &\le p_*^{-1}\max_{z\in\mathcal Z_+}
       \sum_{j=1}^{m(z)}P_z(V_j>e^{L/2})
 \le C_{K,k,\epsilon}e^{-c_ke^{L/2}}.
 \end{aligned}
\]
For all sufficiently large $N$,
$P_z(E_z\cap D_{N,z})>0$ for every $z\in\mathcal Z_+$.
On each of these record spaces,
\[
 d_{\rm TV}\bigl(P_z(\,\cdot\mid E_z),
                 P_z(\,\cdot\mid E_z\cap D_{N,z})\bigr)
 =P_z(D_{N,z}^c\mid E_z).
\]
The same bound holds after applying the centered-mark map.
Together with the empty-index definition, this proves
$\sqrt L\,d_{N,K,k,\epsilon}\to0$ at each fixed $K,k,\epsilon$.
Finally \eqref{chVI:eq:markbounds} proves every coefficient assertion,
including the branch of zero activity mean.
\end{proof}

We now define the centers at which this law is used. On the grid
$\mathcal G_{N,k}(o)$ from Section~6.2, set
\[
 R_k^{\rm ext}:=\lfloor\exp(q_{T+4}+\mathfrak h-4e^{-\mathfrak h})\rfloor,
 \qquad w_k:=2(R_k^{\rm ext}+1)+2.
\]
The integer $R_k^{\rm ext}$ is the radius of the outermost original
disk. The larger protective disk $B(x,2(R_k^{\rm ext}+1))$
contains that disk, its one-step boundary, and the shifted hosts
used in \hyperlink{proof.appendix.X}{Appendix~X}. For $w\ge0$ let the boundary-strip set be
\[
 \begin{split}
 \mathcal B_N(o,w):={}&
 \bigcup_{\substack{z\text{ a macro center}\\j=1,2;\ b\in\{2s,3s\}}}
       \{x:|d_{\rm circ}(x_j,z_j)-b|\le w\}\\
 &\quad\cup\bigcup_{j=1,2}\{x:d_{\rm circ}(x_j,o_j)\le w\}.
 \end{split}
\]
The admissible centers are exactly
\[
 \mathcal X_{N,k}^{\rm adm}(r,o):=
 \{x\in\mathcal G_{N,k}(o):x\notin\mathcal B_N(o,w_k),\
 d_1(x,r)>D_N^{\rm mac}+2w_k,\ 4(w_k+1)<N\}.
\]
The seam condition removes the incomplete grid cells. Membership
in every macro core and buffer is constant on
$B(x,2(R_k^{\rm ext}+1))$ for an admissible $x$, so $z_*$ is
constant there. This protective disk is root-free and lies in
$P_{z_*(x)}^{\rm buf}$.

At a mesh center $x$, apply $E_0$ to the local field
$F_x(z)=\ell_t(x+z)$, $z\in Q(k)$, and denote its mark by $e_x$.
A \emph{physical raw center} satisfies admissibility, all the
corridors \eqref{chVI:eq:corridor}, the event $F_x\in E_0$, and
\begin{equation}\label{chVI:eq:vacancy}
 \ell_t(y)>1\quad
 (y\in P_{z_*(x)}^{\rm buf}\setminus B_{T-1}^+(x)).
\end{equation}
Write $\mathcal R^{\rm phys}$ for these centers and $\mathcal R$
for the subfamily whose every completed layer-$(T-1)$ record has
duration at most $e^{L/2}$. The two families agree outside an event
whose window probability tends to zero at fixed $k$.
In particular, $\mathcal R\subset\mathcal R^{\rm phys}$ always.
The inner-corridor restriction is
$\mathcal M_z(B)=P_z(F\in B,\mathscr C_{\rm int})$.
For $\mathcal R_k(E_0)>0$, the common mark law is
\[
 \mathcal O_{K,k,\epsilon}
 =\Theta_\#\bigl(\mathcal R_k(\,\cdot\,\cap E_0)/\mathcal R_k(E_0)\bigr),
\]
with the empty-mark convention defined above when $\mathcal R_k(E_0)=0$.
Selected factors in Section~6.4 additionally restrict record durations.
Admissibility, terminal/outer data and \eqref{chVI:eq:vacancy}
enter the exterior certificate $\Lambda_C$ in \eqref{chVI:eq:certificate}.
The original-path containment condition is specified below.

A physical raw center is \emph{false} if its owner's soft cloud
is not contained in its cell. Thus
\[
 \mathcal F_N^{\rm false}:=\{x\in\mathcal R^{\rm phys}:
 S_t\cap P_{z_*(x)}^{\rm buf}\not\subset x+Q(k)\}.
\]
A center in $\mathcal R\setminus\mathcal F_N^{\rm false}$ is called \emph{true}.
For a true center the entire owner cloud equals the soft set in
its cell.

\subsection*{6.4. Concentration of the full raw sum}
\label{chVI:sec:concentration}
For the marks $e_x$ of the raw centers, put
\begin{equation}\label{chVI:eq:centeredsum}
 \Delta_{N,K,k,\epsilon}=
 \sum_{x\in\mathcal R}(J(e_x)-\beta_{K,k,\epsilon}H(e_x)).
\end{equation}
The index set in \eqref{chVI:eq:centeredsum} is the full raw family;
colors enter its selection and concentration estimates.

\begin{proposition}[False centers]\label{chVI:prop:false}
For fixed $K,\epsilon>0$ and $\tau,\rho>0$, for all sufficiently
large $k$, then all sufficiently large $N$,
\begin{equation}\label{chVI:eq:false}
 \Qw(|\mathcal F_N^{\rm false}|\ge\tau\sqrt L)\le\rho,
\end{equation}
uniformly in the origin.
\end{proposition}
\begin{proof}
Choose $B$ by \eqref{chVI:eq:R18.7}. For the active-owner count
and path guard of Appendix~W, \eqref{chVI:eq:R13.20} gives, in
the body's indexing,
\[
 \begin{aligned}
 &\EE_{\Qw}[|\mathcal F_N^{\rm false}|;
       \ \#\mathbb W_N^{\rm mac}(t)\le B\sqrt L,\ \mathsf G_{N,k}]\\
 &\qquad\le CB\sqrt L\sum_{i=2}^{T+4}e^{C\mathfrak h-cq_i^{1/3}}
       +o_N(\sqrt L).
 \end{aligned}
\]
The cap failure can be made arbitrarily small by choosing $B$;
the guard failure is $o_N(1)$ at fixed $k$. Since
$\gamma=1/20<1/3$, first choose $k$ to make the sum small,
then $N$ for the remaining errors. Markov's inequality proves
\eqref{chVI:eq:false}.
\end{proof}

\begin{proposition}[Full centered sum]\label{chVI:prop:concentration}
Fix $K$, $\epsilon>0$ and $\tau,\rho>0$. For all sufficiently
large $k$, then all sufficiently large $N$, uniformly in the
deterministic origin $o$,
\begin{equation}\label{chVI:eq:concentration}
 \Qw(|\Delta_{N,K,k,\epsilon}|\ge\tau\sqrt L)\le\rho.
\end{equation}
The same assertion holds under $\Qaug$.
\end{proposition}
\begin{proof}
Write
\[
 \mathcal O^{\rm active}=\mathcal O^{\rm active}(t)
 :=\{z:S_t\cap P_z^{\rm buf}\ne\varnothing\}
\]
for the buffer-active macro owners. The false-center estimate in
Proposition~\ref{chVI:prop:false} and the deterministic bound
\begin{equation}\label{chVI:eq:rawcap}
 |\mathcal R^{\rm phys}|\le256\,|\mathcal O^{\rm active}|+|\mathcal F_N^{\rm false}|
\end{equation}
give tightness of $|\mathcal R|/\sqrt L$ from the soft packing
input and Proposition~\ref{chVI:prop:false}. The factor $256$ is a
mesh-overlap bound. Choose $B$ so that
$\Qw(|\mathcal R|>B\sqrt L)$ is arbitrarily small.

Use the macro-owner colors of Section~6.2 and, independently, a
deterministic coloring of the protective disks
$B(x,2(R_k^{\rm ext}+1))$, disjoint within each fine color.
There are at most $J_k<\infty$ pairs of colors, uniformly in $N$
at fixed $k$. Within one pair, different owners have disjoint
extended buffers; candidates with the same owner have disjoint
protective disks and hence disjoint terminal disks.

Fix such a pair, the root, origin, candidate order and a feasible
discrete endpoint-word stratum. In this proof write $\mathcal I$
for its deterministic admissible candidate set,
$M=\lceil B\sqrt L\rceil$, $D_x=B_{T-1}^+(x)$, and $\xi_x$ for
the complete terminal-record coordinate. Let $\mathscr C_{\rm int}(x)$ denote the centered-at-$x$ version
of the event $\mathscr C_{\rm int}$ from Section~6.3: it is the
intersection of the count corridors at layers $0,\ldots,T-2$, with
exponent $\eta'$ at layer $0$ and $\eta$ at the remaining layers.
These crossing counts are read from the complete terminal-record
coordinate $\xi_x$. Define
\[
 \begin{aligned}
 A_x:={}&\mathscr C_{\rm int}(x)\cap\{F_x\in E_0\}\\
       &{}\cap\{\text{each terminal record at }x
                    \text{ has duration }\le e^{L/2}\}.
 \end{aligned}
\]
Complete records reconstruct the occupation in $x+Q(k)$, so
$A_x$ is a local event and supplies a soft point in that cell.
The terminal count, its narrow corridor and the outer guards remain
exterior data.
For an ordered set $S$, let $\operatorname{First}_M(S)$ consist
of its first $\min(M,|S|)$ elements. The selection rule is
\[
 \operatorname{Sel}:=\operatorname{First}_M(\mathcal R\cap\mathcal I).
\]

For a proposed $C\subset\mathcal I$, let $y$ be the multi-host
collapsed exterior of Lemma~\ref{chIII:lem:multihost}, including
all unselected records, and put $z_*(C)=\{z_*(x):x\in C\}$.
If two points of $C$ have the same owner, set $\Lambda_C=\varnothing$:
the local flags supply soft points in their respective cells,
and their terminal disks are disjoint. Each soft point is therefore
outside the other terminal disk, contradicting \eqref{chVI:eq:vacancy}.

Suppose the owners in $C$ are distinct. Their vacancy and
terminal/outer tests are measurable from $y$: other selected
disks lie in disjoint buffers, and disjoint protective disks
preserve the terminal/outer tests of every candidate.
On $\bigcap_{x\in C}A_x$, every unselected candidate with owner
in $z_*(C)$ fails vacancy. Every candidate with owner outside
$z_*(C)$ has its entire raw test determined by $y$.
For $\zeta\in\{\eta',\eta\}$ let
$\mathscr C_i^{\zeta}(x)$ denote the count-corridor event
\eqref{chVI:eq:corridor} at layer $i$ with exponent $\zeta$.
Admissibility is already encoded by $C\subset\mathcal I$. Define the
selected exterior event explicitly by
\[
 \mathcal E_C(y):=
 \bigcap_{x\in C}\left(
   \mathscr C_{T-1}^{\eta'}(x)
   \cap\bigcap_{i=T}^{T+4}\mathscr C_i^{\eta}(x)
   \cap\{\ell_t(v)>1:\ v\in
       P_{z_*(x)}^{\rm buf}\setminus B_{T-1}^+(x)\}
 \right).
\]
The multi-host exterior of Lemma~\ref{chIII:lem:multihost} retains the
terminal count, the outer-guard counts and the occupation outside the
selected terminal disks, so $\mathcal E_C$ is measurable from $y$.
Conditions involving the censored residual time are read through the
factor $G(y,q)$ in \eqref{chIII:eq:multihost-slice}, not through a
selected record. Define
\[
 \begin{aligned}
 \mathcal R_C^{\rm out}(y)
 &:=\{x\in\mathcal I:z_*(x)\notin z_*(C),\
                   x\text{ satisfies its raw tests in }y\},\\
 \Lambda_C(y)
 &:=\mathcal E_C(y)\cap
       \{\operatorname{First}_M(C\cup\mathcal R_C^{\rm out}(y))=C\}.
 \end{aligned}
\]
Then
\begin{equation}\label{chVI:eq:certificate}
 \{\operatorname{Sel}=C\}
 =\Lambda_C(y)\cap\bigcap_{x\in C}\{\xi_x\in A_x\}.
\end{equation}
For the forward direction, selection supplies every local flag
$\xi_x\in A_x$, while the retained exterior records the terminal,
outer-guard and vacancy conditions and hence lies in $\Lambda_C$.
Conversely, on the right every selected candidate passes its local
flag, every same-owner unselected candidate fails vacancy, and every
other candidate is exactly the member of $\mathcal R_C^{\rm out}(y)$
read from the exterior. The $\operatorname{First}_M$ condition in
$\Lambda_C$ therefore reconstructs $\operatorname{Sel}=C$.
If $P_x$ is the exact kernel with the fixed terminal word at $x$,
Lemma~\ref{chIII:lem:selection} gives, on positive-mass strata,
\[
 \bigotimes_{x\in C}P_x(\,\cdot\mid A_x),
 \qquad
 \text{exterior multiplier }\
 \one_{\Lambda_C(y)}\prod_{x\in C}P_x(A_x).
\]
The local conditional product law is formed before the loop-time
weight; the exterior multiplier is then integrated with the same
$J_I$ factor from \eqref{chIII:eq:window-functional}.

Let $P_C$ be the conditional product law on the selected record
coordinates in \eqref{chVI:eq:certificate}. Write
$f_x:=J(e_x)-\beta_{K,k,\epsilon}H(e_x)$ and $\mu_x:=E_{P_C}f_x$. Then $|f_x|\le (1+|Q_K|)$ and
\[
 |\mu_x|\le2(1+|Q_K|)\left(\frac{2e_k}{1-e_k}+d_{N,K,k,\epsilon}\right).
\]
For every $a>0$, Hoeffding's inequality on a selected fiber gives
\begin{equation}\label{chVI:eq:hoeffding}
 P_C\left(\left|\sum_{x\in C}(f_x-\mu_x)\right|\ge a\sqrt L\right)
 \le2\exp\{-a^2L/(8(1+|Q_K|)^2\lceil B\sqrt L\rceil)\}.
\end{equation}
We transfer this conditional-product bound with the $J_I$ weight.
If $m_k$ is the largest count allowed by the terminal corridor,
the selected duration lies in
\[
 [0,d_N^{\rm sel}],\qquad
 d_N^{\rm sel}:=\lceil B\sqrt L\rceil m_ke^{L/2}=o(b_N).
\]
Apply \eqref{chIII:eq:short-window} to two independent selected
copies, as in \hyperref[chVI:def:weighted-window]{Appendix~U.3},
with the indicator of the Hoeffding event as the nonnegative test.
Its normalized physical-window probability is bounded by its
conditional product probability times
$\mathsf W_N(I^{d_N^{\rm sel}};1)/\mathsf W_N(I;1)=1+o_N(1)$.

Apply \eqref{chVI:eq:hoeffding} with $a=\tau/(2J_k)$ and take a
union over colors. At fixed $k$ this tends to zero with $N$.
On the original raw-count cap, the sum of the means is bounded by
\[
 2(1+|Q_K|) B\sqrt L
 \left(\frac{2e_k}{1-e_k}+d_{N,K,k,\epsilon}\right).
\]
On this original-sample cap, each color selects all its raw
centers, so the selected sets partition $\mathcal R$ and count
each raw center exactly once. After fixing the raw-count cap $B$,
choose $k$ so that the local-law error is within the prescribed budget,
and then take $N$ sufficiently large.
Write $G:=\{|\mathcal R|\le B\sqrt L\}$ and let $C_j$ be the
selected set for color pair $j$, $1\le j\le J_k$. Define
\[
 H_j:=\left\{
   \left|\sum_{x\in C_j}(f_x-\mu_x)\right|
       \ge\frac{\tau\sqrt L}{2J_k}
       \right\}.
\]
Each $H_j$ has the physical-window bound obtained above from its
selection certificate. On $G$, the sets $C_j$ partition
$\mathcal R$. Once the displayed mean bound is at most
$\tau\sqrt L/2$, the triangle inequality gives
\[
 \{|\Delta_{N,K,k,\epsilon}|\ge\tau\sqrt L\}\cap G
       \subseteq\bigcup_{j=1}^{J_k}H_j.
\]
Consequently, the complete probability budget is
 \[
 \begin{aligned}
  \Qw(|\Delta_{N,K,k,\epsilon}|\ge\tau\sqrt L)
  &\le\Qw(|\mathcal R|>B\sqrt L)\\
  &\quad+2J_k(1+o_N(1))
  \exp\!\left\{-\frac{\tau^2L}
  {32J_k^2(1+|Q_K|)^2\lceil B\sqrt L\rceil}\right\}.
 \end{aligned}
 \]
\end{proof}

\subsection*{6.5. Coverage, false centers and relative activity}
\label{chVI:sec:coverage}
The false-center bound is Proposition~\ref{chVI:prop:false},
stated before the concentration estimate. We now compare true
centers with retained parents.

A true center $x$ is assigned to a parent $D$ when
$S_t\cap P_{z_*(x)}^{\rm buf}=D$. Define the uncovered family by
\[
 \mathcal U:=\{D\in\mathcal P_{N,K}^{\rm ret}:
 \text{no }x\in\mathcal R\setminus\mathcal F_N^{\rm false}
 \text{ is assigned to }D\}.
\]

\begin{proposition}[Coverage of original retained parents]\label{chVI:prop:coverage}
For fixed $K$, $\tau,\rho,\epsilon_0>0$, there is
$0<\epsilon\le\min(1,\epsilon_0)$ such that for every sufficiently
large $k$, then every sufficiently large $N$,
\begin{equation}\label{chVI:eq:coverage}
 \Qaug(|\mathcal U|\ge\tau\sqrt L)\le\rho.
\end{equation}
\end{proposition}
\begin{proof}
In \hyperlink{proof.appendix.X}{Appendix~X}, the two directional
bounds \eqref{chVI:eq:R21.11}--\eqref{chVI:eq:R21.12} compare
the crossing-count and harmonic-occupation failures at layers above
the first. Equation~\eqref{chVI:eq:R21.14} bounds the number of
parents for which these local implications fail by $o_P(\sqrt L)$.
For the remaining harmonic failures, the energy-truncated Gaussian
probes imply the image event $\mathcal A_{N,k}$; the shell repulsion
estimate of \hyperlink{proof.appendix.Y}{Appendix~Y} gives
\eqref{chVI:eq:R21.16}. Finally, \eqref{chVI:eq:R21.18} gives the
first-layer corridor from the central second-layer corridor.

Choose a separate constant $B_{\rm par}$ using the parent-count
part of Theorem~\ref{chIV:thm:tight}, so that
\[
 \Qw(|\mathcal P_N(U_1(t))|>B_{\rm par}\sqrt L)\le\rho/3
\]
for all sufficiently large $N$. For every origin-free event $E$,
$(\Qw\otimes\operatorname{Unif}(V_N))(E)=\Qw(E)$, so the same
parent-count cap holds under $\Qaug$. On the complementary
cap, uniform random translation loses in expectation at most
$C(\epsilon+D_Ke^{-k})B_{\rm par}\sqrt L$ parents to trimming.
Macro-boundary, seam and root exclusions have probability $o_N(1)$
for fixed $k$. Choose $B_{\rm par}$ from parent-count tightness, then $\epsilon$, next $k$, and finally $N$ from the corresponding thresholds.
The resulting physical raw center contains the entire owner cloud
and hence is not false. Under $\Qw$, the loop-endpoint and
positive-time conditions hold almost surely. The duration-tail
estimate gives, at fixed local parameters and uniformly in the root and origin,
\[
 \mathsf Q_{N,r}^I(\mathcal R^{\rm phys}\ne\mathcal R)=o_N(1).
\]
Averaging over the origin gives the same bound under $\Qaug$.
\end{proof}

\begin{lemma}[Relative local activity]\label{chVI:lem:relativeactivity}
For fixed $Q_K$ and every $\varepsilon_{\rm rel}>0$, for all sufficiently large $N$,
every true raw center $x$ assigned to parent $D$ satisfies
\begin{equation}\label{chVI:eq:relativeactivity}
 |H(e_x)-a_D^G(t)|\le \varepsilon_{\rm rel} a_D^G(t).
\end{equation}
\end{lemma}
\begin{proof}
Union length is invariant under a common translation of its
interval centers. Subtract the anchor value of the finite
root-conditioned residual. Chapter~I's Green difference estimates
give convergence, in total variation, of its increment vector on
$Q_K\setminus\{0\}$ to the potential-kernel field in
\eqref{chVI:eq:localmarks}, uniformly at these root-far anchors.
The covariance is nondegenerate on the distinct nonzero
coordinates, so entrywise covariance convergence gives this total
variation convergence.

Write $e_x=(E,f)$, with $E\subset Q_K$, let
$r_z:=\sqrt{(v-2f(z))_+}$, and define
$U(\xi):=|\bigcup_{z\in E}[-\xi_z-r_z,-\xi_z+r_z]|$.
Write $\mu_\infty$ for the potential-kernel increment law and $\mu_N$ for
the finite root-conditioned increment law on $E$. Under either law,
$2r_z\le U\le\sum_{y\in E}2r_y$. Consequently
\[
 \left|E_{\mu_\infty} U-E_{\mu_N} U\right|
 \le2d_{\rm TV}(\mu_\infty,\mu_N)\sum_{z\in E}2r_z
 \le2|Q_K|d_{\rm TV}(\mu_\infty,\mu_N)E_{\mu_N} U.
\]
Multiplication by $c_g$ proves \eqref{chVI:eq:relativeactivity};
when every radius is zero, both sides vanish.
\end{proof}

\subsection*{6.6. The deterministic budget for Chapter VII}
\label{chVI:sec:budget}
Let $U_0(t):=\{x:\ell_t(x)=0\}$ and set
$C_N(t):=L^{-1/2}|\operatorname{Rep}_N(U_0(t))|$, where
$\operatorname{Rep}_N$ is the canonical representative set of Chapter~III.
For each fixed observation time $t$, all path-dependent quantities
in this section are evaluated at that same time. We suppress fixed
arguments when convenient, in particular
$A_N=A_N(t)$, $C_N=C_N(t)$ and
$D_N^{\rm wide}=D_N^{\rm wide}(U_1(t))$, the widened soft-defect
set of Chapter~IV.
The latter set contains every canonical soft defect and
$|D_N^{\rm wide}|/\sqrt L\to0$ in the present windows.
The following statement is deterministic once the local parameters
have been fixed.

\begin{theorem}[Scalar residual budget]\label{chVI:thm:budget}
Fix $K,k\in\mathbb N$ with $k^{1/20}\ge16$,
$0<\epsilon\le1$, and $\varepsilon_{\rm rel}>0$.
For all sufficiently large $N$, all roots, origins, paths and
times $t\ge0$, and every $0\le a\le B_*$,
\begin{equation}\label{chVI:eq:budget}
\begin{split}
 |C_N(t)-aA_N(t)|\le{}&
  \frac{|D_N^{\rm wide}|}{\sqrt L}
  +\frac{|\sum_{x\in\mathcal R}(J(e_x)-aH(e_x))|}{\sqrt L}\\
  &+(1+|Q_K|)\frac{|\mathcal F_N^{\rm false}|}{\sqrt L}+a\varepsilon_{\rm rel}A_N(t)\\
 &+(B_*+a)\left(O_{N,K}+B_*^{-1}|Q_K|\frac{|\mathcal U|}{\sqrt L}\right).
\end{split}
\end{equation}
For fixed local parameters, this is a deterministic pathwise
inequality; $\mathcal U$ is the corresponding uncovered family.
\end{theorem}
\begin{proof}
Fix a true center $x$ and, within this proof, write
\[
 z=z_*(x),\qquad D=S_t\cap P_z^{\rm buf}=S_t\cap(x+Q(k)).
\]
By $E_0$, $D$ lies in a translate of $Q_K$, so
$\diam_1D\le4D_K<\rho_N$ for sufficiently large $N$ at fixed $K,k$.
Since $x\in P_z$, the core/buffer margin gives
\[
 R_N=o(s),\qquad \tfrac12\lfloor e^k\rfloor+R_N<s.
\]
If $y\in D$ had a canonical-defect witness $w\in S_t$ with
$\rho_N<d_1(y,w)<R_N$, then $w\in P_z^{\rm buf}$, hence $w\in D$,
contradicting $\diam_1D<\rho_N$. Thus every vertex of $D$ is clean.
The set $D$ is connected at radius $\rho_N$, and the same buffer
margin prevents a $\rho_N$-edge to a soft point outside $D$.
Therefore $D$ is one canonical parent.
Distinct true centers have disjoint half-open cells and hence
distinct nonempty parent clouds. This gives an injection into
the full parent family.

The false-center terms satisfy
\[
 \left|\sum_{x\in\mathcal R\cap\mathcal F_N^{\rm false}}
       (J(e_x)-aH(e_x))\right|
 \le(1+|Q_K|)|\mathcal F_N^{\rm false}|.
\]
In the next sum, $(x,D)$ ranges over true centers and their assigned
parents $D=S_t\cap P_{z_*(x)}^{\rm buf}$. Since $J(e_x)=J_D$,
\eqref{chVI:eq:relativeactivity} and the injection above give
\[
 \left|\sum_{(x,D)}(J(e_x)-aH(e_x)-J_D+aa_D^G(t))\right|
 \le a\varepsilon_{\rm rel}\sum_{(x,D)}a_D^G(t)
 \le a\varepsilon_{\rm rel}\sqrt L\,A_N(t).
\]
Unassigned retained parents are precisely $\mathcal U$ and each has
activity at most $B_*^{-1}|Q_K|$. By \eqref{chVI:eq:loweractivity},
\[
 \begin{aligned}
 \sum_{\substack{D\in\mathcal P_N(U_1(t))\\D\ \mathrm{unassigned}}}
     |J_D-aa_D^G(t)|
 &\le(B_*+a)\sum_{\substack{D\in\mathcal P_N(U_1(t))\\D\ \mathrm{unassigned}}}
     a_D^G(t)\\
 &\le(B_*+a)\bigl(\sqrt L\,O_{N,K}+B_*^{-1}|Q_K||\mathcal U|\bigr).
 \end{aligned}
\]
Divide these bounds by $\sqrt L$.

Finally a zero in a clean soft parent is clean in the zero set.
Its parent has diameter at most $\rho_N$, so all such zeros form
one canonical zero component. Distinct parents cannot merge at
this connection radius. Each remaining zero component contains
a canonical soft defect; select one such vertex to obtain an
injection into the defect set. Hence
\[
 0\le C_N(t)-L^{-1/2}\sum_{D\in\mathcal P_N(U_1(t))}J_D
 \le |D_N^{\rm wide}|/\sqrt L.
\]
Combining the two inequalities proves \eqref{chVI:eq:budget}.
\end{proof}

For positive thresholds $d,o,A,\tau$ and relative error $\varepsilon_{\rm rel}>0$,
the exact six-event consequence is the following. If
\begin{equation}\label{chVI:eq:thresholdbudget}
 3(1+|Q_K|)d+B_*\varepsilon_{\rm rel}A+2B_*o\le\tau,
\end{equation}
then, at $a=\beta_{K,k,\epsilon}$, the event
$\{|C_N-aA_N|>\tau\}$ is contained in the union of
\begin{equation}\label{chVI:eq:sixevents}
\begin{gathered}
 \{|D_N^{\rm wide}|>d\sqrt L\},\quad
 \{|\Delta_{N,K,k,\epsilon}|>d\sqrt L\},\quad
 \{|\mathcal F_N^{\rm false}|>d\sqrt L\},\\
 \{A_N>A\},\qquad\{|\mathcal U|>d\sqrt L\},\qquad
 \{O_{N,K}>o\}.
\end{gathered}
\end{equation}

For fixed $(\tau,\rho)$, choose successively
$(o,K)\longrightarrow d\longrightarrow\epsilon\longrightarrow A
\longrightarrow\varepsilon_{\rm rel}\longrightarrow k$
so that the deterministic budget holds and each of the six exceptional
probabilities is at most $\rho/6$. With these parameters fixed, take
$N$ sufficiently large. The residual event is origin-free, so the
product-measure identity above transfers the augmented-law bound to
$\Qw$. This provides the finite-accuracy comparison
needed in Chapter~VII. Showing that the resulting coefficients
converge to a single positive constant, and recovering the
deterministic observation time, are the tasks of that chapter.

\clearpage
\endgroup

\begingroup
\mergedchaptersetup{VII}{The number of unvisited components}
\begin{center}
{\large\scshape Chapter VII\par}\vspace{6pt}
{\LARGE The number of unvisited components\par}
\end{center}
\let\mathscr\mathcal
We identify the number of canonical unvisited components with a
deterministic multiple of the local-time activity. Chapter~VI gives
this comparison to any fixed accuracy, but its coefficient depends
on the local parameters. We first place all these comparisons on
one common window and show that their coefficients have a unique
limit. The hit-and-defect increment estimate transfers this comparison
to the deterministic observation time. Packing and bounded-count
local clearing give nondegeneracy of the zero count and hence
positivity of the identified coefficient.

\subsection*{7.1. Canonical zero components and nondegeneracy}
Retain the rate-one walk on $V_N=(\ZZ/N\ZZ)^2$, and write
\[
 n=N^2,\qquad L=\log N,\qquad g=2/\pi,\qquad b_N=gnL.
\]
\[
 s_N=b_N(2L-\tfrac32\log L)+\tfrac{9g}{32}n(\log L)^2.
\]
Use the fixed-root path convention of Chapter~III, with $r=X_0$. Unless
another law is specified, probabilities refer to the stationary
ordinary path law $\PP_{\pi_N}$, where $\pi_N$ is uniform on $V_N$.
For $u\ge0$ put $U_u(t)=\{x:\ell_t(x)\le u\}$ and let
$U_0(t)=\{x:\ell_t(x)=0\}$.

Here and below, a canonical component has the meaning fixed in
Chapters~III--IV. More explicitly, set
\begin{equation}\label{chVII:eq:scales}
 \rho_N=\tfrac12e^{L^{1/6}},\qquad R_N=Ne^{-2L^{1/8}}.
\end{equation}
Write $d_1$ for graph distance and $d$ for periodic Euclidean distance.
From a set $U$ remove
\[
 D_N^{\rm can}(U):=\{x\in U:\exists y\in U,\
                       \rho_N<d_1(x,y)<R_N\},
\]
put $U^\circ=U\setminus D_N^{\rm can}(U)$, and let
$\mathcal P_N(U)$ be the connected components of the graph on
$U^\circ$ with an edge between distinct vertices at $d_1$-distance
at most $\rho_N$. Equivalently,
$\mathcal P_N(U)=\{D_a(U):a\in\operatorname{Rep}_N(U)\}$ in the notation
of Section~3.1. For all sufficiently large $N$, each component
has $d_1$-diameter at most $\rho_N$ and distinct components have
$d_1$-distance at least $R_N$. Retain $D_N^{\rm wide}(U)$ for the
wide-defect set of Chapter~IV. For all sufficiently large $N$,
\[
 D_N^{\rm can}(U)\subset D_N^{\rm wide}(U).
\]
Indeed $d\le d_1\le\sqrt2d$, while
$\rho_N/\sqrt2>\tfrac12e^{L^{1/8}}$ and
$R_N<2Ne^{-L^{1/8}}$ eventually.
Define
\begin{equation}\label{chVII:eq:objects}
 C_N(t)=\frac{|\mathcal P_N(U_0(t))|}{\sqrt L},\qquad
 A_N(t)=\frac1{\sqrt L}\sum_{D\in\mathcal P_N(U_1(t))}
                  a_D^G(\ell_t;r,1/4).
\end{equation}
As in Chapter~III, these arguments are physical times:
$\widehat A_N(c)=A_N(s_N+cb_N)$ and
$\widehat C_N(c)=C_N(s_N+cb_N)$.
Thus $C_N(t)$ is the normalized number of canonical components of
$U_0(t)$, with ``canonical'' defined by \eqref{chVII:eq:scales}.
The Gaussian activity in \eqref{chVII:eq:objects} is the one from
\eqref{chIII:eq:activity}, using the mean-zero Gaussian field $\phi=\phi_N$
and its root-conditioned residual
\[
 \phi_N^{(r)}(x)
 =\phi(x)-\frac{G_N(x,r)}{G_N(r,r)}\phi(r).
\]
Set $A_N=A_N(s_N)$, $C_N=C_N(s_N)$ and $U_i=U_i(s_N)$ for
$i=0,1$. The input from Chapter~V is
\begin{equation}\label{chVII:eq:activitylimit}
 A_N\ \Longrightarrow\ \Gamma Z,\qquad
 Z=\mathcal Z(\mathbb T^2),\qquad 0<Z<\infty\quad\hbox{a.s.},
 \qquad \Gamma>0.
\end{equation}
Also retain the coefficient bound from Chapter~VI:
\begin{equation}\label{chVII:eq:cap}
 c_g=(2\pi g)^{-1/2},\qquad
 B_*=(2c_g\sqrt{1/4})^{-1}.
\end{equation}

For a finite set $U\subset V_N$, define its packing number by
\[
 \Pack_N(U):=\max\{|A|:A\subset U,\
 d_1(x,y)\ge R_N\text{ for all distinct }x,y\in A\}.
\]
In particular $\Pack_N(\varnothing)=0$.

\begin{proposition}[Nondegeneracy]\label{chVII:prop:nondeg}
There is $\zeta>0$ such that, for all sufficiently large $N$,
\begin{equation}\label{chVII:eq:nondeg}
 \PP_{\pi_N}(C_N\ge\zeta)\ge\tfrac12.
\end{equation}
\end{proposition}
\begin{proof}
\hyperlink{proof.appendix.Z}{Appendix~Z} proves that for every $\eta>0$ there is $b>0$ with
\begin{equation}\label{chVII:eq:packingtail}
 \PP_{\pi_N}\bigl(\Pack_N(U_0(s_N))<b\sqrt L\bigr)\le\eta
\end{equation}
eventually. Appendix~Z first transfers \eqref{chVII:eq:activitylimit}
using the activity increment \eqref{chVII:eq:incA}, then applies
bounded-count selection and local clearing
(Lemmas~\ref{chVII:lem:selection} and \ref{chVII:lem:clearing}).
Appendix~AA derives the increment input from tightness, defect
control and the endpoint-weighted bridge-hit bound.

Delete the points of a maximizing packing that belong to
$D_N^{\rm can}(U_0(s_N))$. Each surviving canonical component
contains at most one remaining packing point, since its diameter
is less than $R_N$. Consequently
\begin{equation}\label{chVII:eq:packingcount}
 C_N\ge
 \frac{\Pack_N(U_0(s_N))-|D_N^{\rm can}(U_0(s_N))|}{\sqrt L}.
\end{equation}
Since $U_0\subset U_1$, the defect set in this bound is contained
in $D_N^{\rm can}(U_1)$, and hence in $D_N^{\rm wide}(U_1)$,
the widened soft-defect set controlled in Chapter~IV. Take $\eta=1/4$ in
\eqref{chVII:eq:packingtail}. The probability that its normalized
defect count exceeds $b/2$ is eventually at most $1/4$.
Equation~\eqref{chVII:eq:nondeg} follows with $\zeta=b/2$.
\end{proof}

The packing step uses the set monotonicity
\[
 U\subset V\quad\Longrightarrow\quad\Pack_N(U)\le\Pack_N(V).
\]

\subsection*{7.2. A finite-accuracy comparison on physical windows}
For a bounded interval $I=[l_N,r_N]\subset(0,\infty)$ of positive
length, recall the normalized loop measure of Section~3.2:
for a measurable event $F$ of the root, path and terminal time,
\begin{equation}\label{chVII:eq:windowlaw}
 \begin{aligned}
 \mathsf Q_N^I(F)
 &=\frac1{n\,\mathsf W_N(I;1)}\sum_{r\in V_N}
       \int_I \EE_r[\one_{\{X_t=r\}}\one_F(r,\omega,t)]\,dt,\\
 \mathsf W_N(I;1)&=\frac1n\sum_r\int_I p_t(r,r)\,dt.
 \end{aligned}
\end{equation}
Here $\omega$ denotes the path and $r=X_0$. Translation invariance
makes the root uniform, while the time marginal has density
$p_t(0,0)/\mathsf W_N(I;1)$ on $I$.
For fixed $W\ge0$ and $\delta>0$, the admissible windows satisfy
\begin{equation}\label{chVII:eq:admissible}
 s_N-Wb_N\le l_N\le r_N\le s_N+Wb_N,\qquad
 \delta b_N\le l_N,\qquad \delta b_N\le r_N-l_N.
\end{equation}
All assertions concerning these windows are eventual assertions
as $N\to\infty$, with $W,\delta$ fixed.

\begin{proposition}[Window residual]\label{chVII:prop:residual}
Fix an admissible sequence $I_N$, and let $\tau,\varepsilon>0$.
There is a deterministic $a\in[0,B_*]$, independent of $N$,
such that
\begin{equation}\label{chVII:eq:residual}
 \mathsf Q_N^{I_N}(|C_N(t)-aA_N(t)|>\tau)\le\varepsilon
\end{equation}
for all sufficiently large $N$.
\end{proposition}
\begin{proof}
We use Theorem~\ref{chVI:thm:budget}, with its anchor set $Q_K$,
omitted activity $O_{N,K}(t)$, uncovered family $\mathcal U$,
false-center family $\mathcal F_N^{\rm false}$, and centered sum
$\Delta_{N,K,k,\epsilon}$. All these objects are evaluated at the sampled terminal time $t$
and at the same grid origin and local parameters.
That theorem gives the following deterministic implication:
if each of the six quantities
\begin{equation}\label{chVII:eq:six}
 \frac{|D_N^{\rm wide}(U_1(t))|}{\sqrt L},\quad
 \frac{|\Delta_{N,K,k,\epsilon}|}{\sqrt L},\quad
 \frac{|\mathcal F_N^{\rm false}|}{\sqrt L},\quad
 A_N(t),\quad \frac{|\mathcal U|}{\sqrt L},\quad O_{N,K}(t)
\end{equation}
is bounded respectively by $d,d,d,H,d,o$, then
\begin{equation}\label{chVII:eq:budget}
 |C_N(t)-aA_N(t)|
 \le3(1+|Q_K|)d+B_*\varepsilon_{\rm rel}H+2B_*o.
\end{equation}
Here $\varepsilon_{\rm rel}$ is the relative finite-field comparison error and
$a=\beta_{K,k,\epsilon}$ is the common local coefficient.

Choose $o>0$ with $2B_*o<\tau/4$, and then choose $K$ so
that the omitted-activity exception has probability at most
$\varepsilon/6$. Choose $d>0$ with
$3(1+|Q_K|)d<\tau/4$, and then the coverage parameter
$\epsilon>0$. Choose $H>0$ by activity tightness, and then
$\varepsilon_{\rm rel}>0$ with $B_*\varepsilon_{\rm rel}H<\tau/4$. Finally choose one sufficiently
large local scale $k$, simultaneously satisfying the coverage,
false-center and centered-sum estimates. With these local parameters
fixed, set $a:=\beta_{K,k,\epsilon}\in[0,B_*]$ and then take $N$
sufficiently large. For such $N$, each exception in \eqref{chVII:eq:six} has probability at most
$\varepsilon/6$.

The estimates involving the grid are applied under
$\mathsf Q_N^{I_N}\otimes\operatorname{Unif}(V_N)$, with one
independent common origin. Their union therefore has probability
at most $\varepsilon$. Marginalizing the auxiliary origin proves
\eqref{chVII:eq:residual}.
\end{proof}

The coefficient in Proposition~\ref{chVII:prop:residual} may depend on
the requested accuracy. We therefore place all finite-accuracy
comparisons on one common window and identify a unique coefficient there.

\subsection*{7.3. A common window and a single coefficient}
Fix, once and for all,
\begin{equation}\label{chVII:eq:common}
 I_N^*=[s_N,s_N+2b_N],\qquad \mathsf Q_N^*=\mathsf Q_N^{I_N^*}.
\end{equation}
This interval satisfies \eqref{chVII:eq:admissible} with $W=2$ and
$\delta=1$ for all sufficiently large $N$.

\begin{lemma}[Activity on the common window]\label{chVII:lem:positive}
The family $A_N(t)$ is tight under $\mathsf Q_N^*$, and there
are $h,p>0$ such that
\begin{equation}\label{chVII:eq:positive}
 \mathsf Q_N^*(A_N(t)\ge h)\ge p
\end{equation}
eventually.
\end{lemma}
\begin{proof}
Tightness is the loop-window activity bound from Chapter~IV.
For the lower bound, \hyperlink{proof.appendix.Z}{Appendix~Z}, Lemma~\ref{chVII:lem:terminallower},
gives $h>0$ and $0<\delta\le1$ such that
\[
 \mathsf Q_N^{[s_N+\delta b_N,s_N+2\delta b_N]}
       (A_N(t)\ge h)\ge\tfrac12.
\]
By the restriction estimate of Lemma~\ref{chVII:lem:restriction},
\[
 \tfrac12\le (6/\delta)\mathsf Q_N^*(A_N(t)\ge h).
\]
Thus \eqref{chVII:eq:positive} holds with $p=\delta/12$.
\end{proof}

\begin{lemma}[Separation of coefficients]\label{chVII:lem:separation}
Let $X,Y$ be real random variables on a probability space with
probability measure $\mu$. For $a,b\in\mathbb R$, $h,p>0$ and
$e,f,r,s\ge0$, suppose
\[
 \mu(X\ge h)\ge p,\quad
 \mu(|Y-aX|>e)\le r,\quad \mu(|Y-bX|>f)\le s,
 \qquad r+s<p.
\]
Then
\begin{equation}\label{chVII:eq:separation}
 |a-b|h\le e+f.
\end{equation}
\end{lemma}
\begin{proof}
The event on which $X\ge h$ and both residual bounds hold
has probability at least $p-r-s>0$. On this event,
$|a-b|h\le |a-b|X\le |Y-aX|+|Y-bX|\le e+f$.
The claimed inequality is between deterministic numbers.
\end{proof}

\begin{proposition}[Identification on the common window]\label{chVII:prop:common}
There is a unique $\beta\in[0,B_*]$ such that
\begin{equation}\label{chVII:eq:commonlimit}
 C_N(t)-\beta A_N(t)\longrightarrow0
 \quad\hbox{in probability under }\mathsf Q_N^*.
\end{equation}
\end{proposition}
\begin{proof}
Let $e_m=(m+1)^{-1}$. For each $m$, apply
Proposition~\ref{chVII:prop:residual} on this same $I_N^*$ with
threshold and failure probability both $e_m$. Choose one resulting
parameter tuple $(K_m,k_m,\epsilon_m)$ and set
\[
 \beta_m:=\beta_{K_m,k_m,\epsilon_m}\in[0,B_*].
\]
For the fixed tuple $(K_m,k_m,\epsilon_m)$,
Proposition~\ref{chVII:prop:residual} yields a deterministic
$N_m<\infty$ such that, for every $N\ge N_m$,
\begin{equation}\label{chVII:eq:approx}
 \mathsf Q_N^*(|C_N(t)-\beta_mA_N(t)|>e_m)\le e_m
\end{equation}
Fix $m,j$ with $e_m+e_j<p$. Choose one $N$ large enough for
both instances of \eqref{chVII:eq:approx} and \eqref{chVII:eq:positive}.
Lemma~\ref{chVII:lem:separation} gives
\[
 |\beta_m-\beta_j|h\le e_m+e_j.
\]
It follows that $(\beta_m)$ is Cauchy. Its limit $\beta$ belongs
to the closed interval $[0,B_*]$. Letting $j\to\infty$ in the
deterministic coefficient bound gives
\[
 |\beta_m-\beta|\le e_m/h\qquad(e_m<p).
\]
This bound holds for every member $\beta_m$ of the constructed
approximating sequence with $e_m<p$.

To pass from $\beta_m$ to $\beta$, fix $\tau,\varepsilon>0$.
Choose $H>0$ so that
$\mathsf Q_N^*(A_N(t)>H)\le\varepsilon/2$ eventually.
Choose a fixed $m$ with
\[
 e_m<\min(\tau/2,\varepsilon/2),\qquad
 |\beta_m-\beta|H<\tau/2.
\]
Outside the activity-tail event and the residual event in
\eqref{chVII:eq:approx}, the residual with coefficient $\beta$ is
less than $\tau$. Its failure probability is at most
$\varepsilon$. This proves \eqref{chVII:eq:commonlimit}.
If another coefficient has the same property, apply
Lemma~\ref{chVII:lem:separation} with arbitrarily small residual
thresholds and failure probabilities. Their difference is zero.
\end{proof}

\subsection*{7.4. Returning to the deterministic time}
For $0<\delta\le1$ put
\begin{equation}\label{chVII:eq:short}
 I_N^\delta=[s_N+\delta b_N,s_N+2\delta b_N],\qquad
 \mathsf Q_N^\delta=\mathsf Q_N^{I_N^\delta}.
\end{equation}
The positive lag between $s_N$ and the left endpoint is used
both for endpoint mixing and for the bridge hitting estimate.

\begin{lemma}[Restriction of windows]\label{chVII:lem:restriction}
For each fixed $0<\delta\le1$ and every measurable event $E$,
\begin{equation}\label{chVII:eq:restriction}
 \mathsf Q_N^\delta(E)\le\frac6\delta\mathsf Q_N^*(E)
\end{equation}
for all sufficiently large $N$.
\end{lemma}
\begin{proof}
The unnormalized measure in \eqref{chVII:eq:windowlaw}, restricted
to $I_N^\delta$, is at most its restriction to $I_N^*$.
The heat-kernel estimates of Chapter~I give, uniformly on
these intervals,
$(2n)^{-1}\le p_t(r,r)\le3/(2n)$ eventually. Hence
\[
 \mathsf W_N(I_N^*;1)\le\frac{3b_N}{n},\qquad
 \mathsf W_N(I_N^\delta;1)\ge\frac{\delta b_N}{2n}.
\]
Dividing the two normalizers proves \eqref{chVII:eq:restriction}.
\end{proof}

\begin{lemma}[Short-window increments]\label{chVII:lem:increments}
For every $\tau,\varepsilon>0$ there is $\delta_0\in(0,1]$
such that, for each fixed $0<\delta\le\delta_0$,
\begin{align}
 \mathsf Q_N^\delta(|A_N(t)-A_N(s_N)|>\tau)&\le\varepsilon,
                                                    \label{chVII:eq:incA}\\
 \mathsf Q_N^\delta(|C_N(t)-C_N(s_N)|>\tau)&\le\varepsilon
                                                    \label{chVII:eq:incC}
\end{align}
for all sufficiently large $N$.
\end{lemma}
\hyperlink{proof.appendix.AA}{Appendix~AA} proves both
\eqref{chVII:eq:incA} and \eqref{chVII:eq:incC} by the hit-and-defect
argument, using tightness and the endpoint-weighted bridge-hit bound.

\begin{lemma}[Prefix comparison]\label{chVII:lem:prefix}
Fix $0<\delta\le1$. For every $\mathcal F_{s_N}$-measurable event
$E_N$ and every $0<\theta\le1/2$,
\begin{equation}\label{chVII:eq:prefix}
 \big|\mathsf Q_N^\delta(E_N)-\PP_{\pi_N}(E_N)\big|\le4\theta
\end{equation}
for all sufficiently large $N$, uniformly in such events.
\end{lemma}
\begin{proof}
By the Markov property, the unnormalized mass of a prefix
event at a terminal time $t$ equals
\[
 \EE_{\pi_N}\bigl[\one_{E_N}
                  p_{t-s_N}(X_{s_N},X_0)\bigr].
\]
Uniformly for $t\in I_N^\delta$, the heat kernel at lag
$t-s_N\ge\delta b_N$ lies between $(1-\theta)/n$ and
$(1+\theta)/n$ eventually. The same bounds apply to
$p_t(r,r)$. Therefore the numerator in \eqref{chVII:eq:windowlaw}
lies between
$(|I_N^\delta|/n)(1\pm\theta)\PP_{\pi_N}(E_N)$,
and its denominator lies between
$(|I_N^\delta|/n)(1\pm\theta)$.
Division gives \eqref{chVII:eq:prefix}; for example,
$2\theta/(1-\theta)\le4\theta$.

\end{proof}

\begin{proposition}[Central proportionality]\label{chVII:prop:central}
For the coefficient of Proposition~\ref{chVII:prop:common},
\begin{equation}\label{chVII:eq:central}
 C_N-\beta A_N\longrightarrow0
 \quad\hbox{in }\PP_{\pi_N}\hbox{-probability}.
\end{equation}
\end{proposition}
\begin{proof}
Fix $\epsilon,\varrho>0$. Apply Lemma~\ref{chVII:lem:increments}
with threshold $\epsilon/[4(1+B_*)]$ and choose $\delta\in(0,1]$
so that each increment exception has probability at most
$\varrho/8$. Choose $0<\theta\le\min(1/2,\varrho/32)$.
After $\delta$ has been fixed, use \eqref{chVII:eq:commonlimit}
with threshold $\epsilon/2$ and failure probability
$\varrho\delta/48$. Lemma~\ref{chVII:lem:restriction} then gives
\[
 \mathsf Q_N^\delta
 (|C_N(t)-\beta A_N(t)|>\epsilon/2)\le\varrho/8.
\]
On the complement of these three exceptions,
\begin{align*}
 |C_N(s_N)-\beta A_N(s_N)|
 &\le |C_N(t)-\beta A_N(t)|+|C_N(t)-C_N(s_N)|\\
 &\hspace{28mm}+\beta|A_N(t)-A_N(s_N)|\\
 &\le \epsilon/2+\epsilon/4=3\epsilon/4.
\end{align*}
The event that the central residual exceeds $\epsilon$ is
$\mathcal F_{s_N}$-measurable. By Lemma~\ref{chVII:lem:prefix}, its ordinary-path
probability is at most
$3\varrho/8+4\theta\le\varrho/2\le\varrho$ eventually.
This proves \eqref{chVII:eq:central}.
\end{proof}

\subsection*{7.5. Positivity and the scalar limit}
\begin{theorem}[The limiting component count]\label{chVII:thm:limit}
There is a deterministic $\beta\in(0,B_*]$ such that
\begin{equation}\label{chVII:eq:finalratio}
 C_N-\beta A_N\longrightarrow0\quad\hbox{in probability}.
\end{equation}
With $\kappa=\beta\Gamma>0$,
\begin{equation}\label{chVII:eq:distribution}
 C_N\Longrightarrow\kappa Z,
 \qquad
 \EE_{\pi_N}e^{-qC_N}\longrightarrow\EE e^{-q\kappa Z}
 \quad(q>0).
\end{equation}
The expectation on the right is with respect to the law of the
critical Gaussian measure $\mathcal Z$.
\end{theorem}
\begin{proof}
Proposition~\ref{chVII:prop:central} proves
\eqref{chVII:eq:finalratio} with $\beta\ge0$.
If $\beta=0$, it gives $C_N\to0$ in probability. In particular,
$\PP_{\pi_N}(C_N>\zeta/2)\to0$, contradicting
Proposition~\ref{chVII:prop:nondeg}. Thus $\beta>0$.
By \eqref{chVII:eq:activitylimit},
$\beta A_N\Longrightarrow\beta\Gamma Z$; Slutsky's theorem
and \eqref{chVII:eq:finalratio} yield the distributional assertion.

For the last assertion apply weak convergence to the bounded
continuous function
\[
 x\longmapsto\exp\{-q\max(x,0)\},\qquad x\in\mathbb R.
\]
Both $C_N$ and $\kappa Z$ are nonnegative, so this function agrees
with the desired Laplace test on their supports.
\end{proof}

The limiting scalar law depends on the critical measure through
its total mass $Z$.

\subsection*{7.6. The set passed to the covering argument}
Let
\[
 \mathcal V_N(t)=V_N\setminus\{X_u:0\le u\le t\}
\]
be the unvisited set. At a deterministic positive time $t$,
\begin{equation}\label{chVII:eq:zerosvisited}
 \mathcal V_N(t)=U_0(t)\quad\hbox{almost surely}.
\end{equation}
Indeed, the rate-one walk has finitely many jumps on $[0,t]$
almost surely, every holding time is strictly positive, and
the probability of a jump at the prescribed time $t$ is zero.
Every vertex visited by time $t$ is consequently occupied for
a positive interval before $t$. The converse is immediate
from the definition of occupation. The same equivalence holds for almost every time on each nonexplosive
path and hence under the time-integrated loop laws by Fubini.

At time $s_N$ write
\begin{equation}\label{chVII:eq:decomposition}
 \mathcal V_N(s_N)=\mathcal D_N^{\rm rem}\ \dot\cup\
                 D_{N,1}\ \dot\cup\cdots\dot\cup\ D_{N,K_N},
 \qquad K_N=\sqrt L\,C_N,
\end{equation}
Here $\dot\cup$ denotes disjoint union,
$\mathcal D_N^{\rm rem}=D_N^{\rm can}(U_0(s_N))$, and the $D_{N,j}$
enumerate $\mathcal P_N(U_0(s_N))$ in the fixed order of their
representatives. Thus $K_N$ is an integer, possibly zero. Chapter~IV and
$U_0\subset U_1$ give
\begin{equation}\label{chVII:eq:deleted}
 \frac{|\mathcal D_N^{\rm rem}|}{\sqrt L}\longrightarrow0
 \quad\hbox{in probability}.
\end{equation}
The components in \eqref{chVII:eq:decomposition} have diameter
at most $\rho_N$ and pairwise separation at least $R_N$.
Together with $K_N/\sqrt L\Longrightarrow\kappa Z$, these
are the geometric and scalar inputs to Chapter~VIII.
The vertices of $\mathcal D_N^{\rm rem}$ remain part of the set to be
covered. Chapter~VIII bounds the probability that any of them is
still unvisited at the future observation time; see the remainder
term in \eqref{chVIII:eq:sandwich}.

By translation invariance, the laws of $A_N$, $C_N$ and
$C_N-\beta A_N$ are independent of the prescribed initial vertex. Translation of the torus
sends the law started at one vertex to the law started at
another and preserves occupation, graph distance, and the
number of canonical components. It also preserves the
activity sum: the Green kernel and the root-conditioned
Gaussian law translate with the root. The choice of a
representative for a component does not affect either total
mass. Thus the laws of $A_N$, $C_N$ and the residual are
identical for all initial vertices and for their uniform
mixture. This translation argument gives uniformity in the initial vertex for
the assertions of this chapter.

\clearpage
\endgroup

\begingroup
\mergedchaptersetup{VIII}{The limiting cover-time distribution}
\begin{center}
{\large\scshape Chapter VIII\par}\vspace{6pt}
{\LARGE The limiting cover-time distribution\par}
\end{center}
\let\mathscr\mathcal
At the central observation time the unvisited set consists, up to
a small remainder, of well-separated components of small diameter.
Chapter~VII identifies their normalized number. We now cover one
representative of each component, estimate the additional time
needed to cover the component itself, and restart the walk at the
observation time. The argument starts from the scalar component
limit \eqref{chVII:eq:distribution} and the geometric remainder estimate
\eqref{chVII:eq:deleted}.

\subsection*{8.1. Target geometry and the remaining time}
Retain the rate-one walk on $V_N=(\ZZ/N\ZZ)^2$ and set
\[
 n=N^2,\qquad L=\log N,\qquad g=2/\pi,\qquad b_N=gnL,
\]
\[
 a_N=b_N(2L-\log L),\qquad t_N(z)=a_N+zb_N.
\]
The central time is
\[
 s_N=b_N(2L-\tfrac32\log L)+\tfrac{9g}{32}n(\log L)^2.
\]
For $U\subset V_N$ let $H_U$ be its entrance time and
$T_U$ its cover time, with $H_\varnothing=\infty$ and
$T_\varnothing=0$. In particular $T_N=T_{V_N}$.
The finite hitting means in \eqref{chI:eq:mean-hit} imply
$T_N=\max_{x\in V_N}H_{\{x\}}<\infty$ almost surely.
Write $d$ for periodic Euclidean distance and $d_1$ for graph
distance, as in Section~1.2; in particular $d\le d_1\le2d$.

We recall the canonical construction, now for an arbitrary
deterministic target $U$. Put
\begin{equation}\label{chVIII:eq:scales}
 \rho_N=\tfrac12e^{L^{1/6}},\qquad
 R_N=Ne^{-2L^{1/8}},
\end{equation}
and
\[
 U^\circ=\{x\in U:\text{there is no }w\in U
                   \text{ with }\rho_N<d_1(x,w)<R_N\}.
\]
Join vertices of $U^\circ$ at graph distance at most $\rho_N$,
choose the lexicographically first vertex of each component in
$\{0,\ldots,N-1\}^2$, as in Section~3.1, and write
$\operatorname{Rep}_N(U)$ for the representatives and $D_a=D_a(U)$
for the component represented by $a$. Thus
$U\setminus U^\circ=D_N^{\rm can}(U)$ from Chapter~VII. Set
\begin{equation}\label{chVIII:eq:masses}
 M_N^{\rm comp}(U)=\frac{|\operatorname{Rep}_N(U)|}{\sqrt L},\qquad
 r_N(U)=\frac{|U\setminus U^\circ|}{\sqrt L},\qquad
 \Phi_m(z)=e^{-me^{-z}}\quad(m\ge0).
\end{equation}
For $U=\varnothing$ both ratios are zero and $\Phi_0=1$.

For large $N$, every $D_a$ has graph diameter at most $\rho_N$,
and different components have graph separation at least $R_N$.
Indeed, the first vertex of a component path more than $\rho_N$
from its initial vertex would be at distance at most
$2\rho_N<R_N$, contradicting the definition of $U^\circ$.
Vertices in different components are more than $\rho_N$ apart,
so the same definition makes their distance at least $R_N$.
Consequently
\begin{equation}\label{chVIII:eq:geometry}
 d(a,b)\ge R_N/2\quad(a\ne b\in\operatorname{Rep}_N(U)),\qquad
 |\{a\in\operatorname{Rep}_N(U):d(y,a)<R_N/8\}|\le1.
\end{equation}
The second assertion follows from the triangle inequality.

It is useful first to work with the unshifted remaining clock
\[
 \bar u_N(z)=b_N(\tfrac12\log L+z).
\]
The actual remaining time is
\begin{equation}\label{chVIII:eq:clock}
 u_N(z):=t_N(z)-s_N=\bar u_N(z-\delta_N^{\rm shift}),\qquad
 \delta_N^{\rm shift}=\frac9{32}\frac{(\log L)^2}{L}\longrightarrow0.
\end{equation}
For each fixed $W$, both $\bar u_N(z)$ and $u_N(z)$ are positive
uniformly on $|z|\le W$ for sufficiently large $N$.

\subsection*{8.2. Survival of a fixed number of separated points}
Put $\omega_N=L^{-7/8}\log L$. The covering argument uses the
following fixed-target survival estimate from Chapter~I.

\begin{proposition}[Fixed-target survival]\label{chVIII:prop:survival}
Fix an integer $k\ge1$ and $W\ge0$. Uniformly over $k$-point sets
$A\subset V_N$ with pairwise Euclidean separation at least
$R_N/2$, and over $|z|\le W$,
\begin{align}
 \PP_y(H_A>\bar u_N(z))&\le C_{k,W}L^{-k/2}
                              &&(y\in V_N),\label{chVIII:eq:upper}\\
 \PP_y(H_A>\bar u_N(z))&=L^{-k/2}e^{-kz}
                       (1+O_{k,W}(\omega_N))
                              &&(d(y,A)\ge R_N/8).
                                                        \label{chVIII:eq:relative}
\end{align}
\end{proposition}
The proof is in \hyperlink{proof.appendix.AB}{Appendix~AB}. Briefly, the Green matrix
on $A$ equals $gL$ times the identity plus entries of order
$L^{1/8}$. Capacity and the killed spectral expansion give
the stationary survival probability to relative error
$O_{k,W}(\omega_N)$. A mixing interval transfers this
estimate to a distant starting point. The early-hit term retains
the subsequent survival factor, giving the relative error
$O_{k,W}(\omega_N)$ in \eqref{chVIII:eq:relative}.

\subsection*{8.3. Factorial moments and the probability of no survivors}
For fixed $y,U,z$, remove the representatives too close to
the initial point and write
\[
 \operatorname{Rep}_N^y(U)=\{a\in\operatorname{Rep}_N(U):d(y,a)\ge R_N/8\},
 \qquad M=|\operatorname{Rep}_N(U)|,\quad s=|\operatorname{Rep}_N^y(U)|.
\]
By \eqref{chVIII:eq:geometry}, $M-1\le s\le M$. Let
\[
 q=L^{-1/2}e^{-z},\qquad \lambda=Mq=M_N^{\rm comp}(U)e^{-z},\qquad
 Y=\sum_{a\in\operatorname{Rep}_N^y(U)}\one_{\{H_a>\bar u_N(z)\}}.
\]
For integers $j,k\ge0$, use the falling factorial
$(j)_k=j(j-1)\cdots(j-k+1)$, with $(j)_0=1$.

\begin{lemma}\label{chVIII:lem:moments}
For every fixed $k\ge1$ and $B,W\ge0$,
\begin{equation}\label{chVIII:eq:moments}
 \left|\EE_y(Y)_k-\lambda^k\right|
 \le C_{k,B,W}(\omega_N+L^{-1/2})
\end{equation}
eventually, uniformly in $y,U,z$ with $M_N^{\rm comp}(U)\le B$ and
$|z|\le W$.
\end{lemma}
\begin{proof}
Expanding over distinct ordered tuples and using
\eqref{chVIII:eq:relative} gives
\[
 \EE_y(Y)_k
 =k!\sum_{\substack{A\subset\operatorname{Rep}_N^y(U)\\|A|=k}}
                  \PP_y(H_A>\bar u_N(z))
 =(s)_kq^k(1+O_{k,W}(\omega_N)).
\]
For $M\ge1$, a tuple among the $M^k$ ordered tuples
that is not counted by $(s)_k$ either uses the one removed
point or has a repeated coordinate. Thus
\[
 0\le M^k-(s)_k
 \le\left(k+\binom{k}{2}\right)M^{k-1}\le k^2M^{k-1}.
\]
Since $Mq\le Be^W$ and $q\le e^WL^{-1/2}$, multiplication
by $q^k$ and the preceding expansion prove
\eqref{chVIII:eq:moments}.
\end{proof}

\begin{proposition}[Covering the representatives]\label{chVIII:prop:void}
For every $B,W\ge0$,
\begin{equation}\label{chVIII:eq:void}
 \sup_{\substack{y\in V_N,\ U\subset V_N\\M_N^{\rm comp}(U)\le B,\ |z|\le W}}
 \left|\PP_y(T_{\operatorname{Rep}_N(U)}\le \bar u_N(z))
                          -\Phi_{M_N^{\rm comp}(U)}(z)\right|
 \longrightarrow0.
\end{equation}
\end{proposition}
\begin{proof}
For an integer $J\ge0$, Bonferroni's inequalities give
\begin{equation}\label{chVIII:eq:bonferroni}
 \sum_{k=0}^{2J+1}\frac{(-1)^k}{k!}\EE_y(Y)_k
 \le\PP_y(Y=0)\le
 \sum_{k=0}^{2J}\frac{(-1)^k}{k!}\EE_y(Y)_k.
\end{equation}
These follow pointwise by evaluating the sums at an
integer value of $Y$: for integers $j\ge1$ and $0\le r<j$,
$\sum_{k=0}^r(-1)^k\binom jk=(-1)^r\binom{j-1}{r}$;
the sum is zero for $r\ge j$, and is one for $j=0$.

For every integer $r\ge0$, Taylor's formula,
using that all derivatives of $e^{-x}$ have absolute value
at most one for $x\ge0$, gives
\[
 \sup_{0\le\lambda\le Be^W}
 \left|e^{-\lambda}-\sum_{k=0}^{r}\frac{(-\lambda)^k}{k!}\right|
 \le\frac{(Be^W)^{r+1}}{(r+1)!}.
\]
Use Lemma~\ref{chVIII:lem:moments} in the finitely many terms of
\eqref{chVIII:eq:bonferroni}. Uniformly in the stated parameters,
\begin{align}
 |\PP_y(Y=0)-e^{-\lambda}|
 &\le C_{J,B,W}(\omega_N+L^{-1/2})\notag\\
 &\quad+\frac{(Be^W)^{2J+1}}{(2J+1)!}
           +\frac{(Be^W)^{2J+2}}{(2J+2)!}.\label{chVIII:eq:truncation}
\end{align}
Fix $J$, let $N\to\infty$, and then let $J\to\infty$.
The event $Y=0$ says that all far representatives have
been visited. At most one representative was removed,
and its survival probability is at most $C_WL^{-1/2}$
by the uniform-in-$y$ bound \eqref{chVIII:eq:upper} with $k=1$. Therefore
\[
 0\le\PP_y(Y=0)-\PP_y(T_{\operatorname{Rep}_N(U)}\le \bar u_N(z))
 \le C_WL^{-1/2},
\]
which proves \eqref{chVIII:eq:void}.
\end{proof}

\subsection*{8.4. Covering entire components and the remainder}
\begin{theorem}[A deterministic configuration bound]\label{chVIII:thm:configuration}
Fix $B,W\ge0$. For every $\epsilon>0$ there is $C>0$
such that, for all sufficiently large $N$, all $y\in V_N$,
all $U\subset V_N$ with $M_N^{\rm comp}(U)\le B$, and all $|z|\le W$,
\begin{equation}\label{chVIII:eq:configuration}
 \left|\PP_y(T_U\le u_N(z))-\Phi_{M_N^{\rm comp}(U)}(z)\right|
 \le\epsilon+C r_N(U).
\end{equation}
\end{theorem}
\begin{proof}
On one future walk, for each $a\in\operatorname{Rep}_N(U)$ put
\[
 \tau_a:=\inf\{u\ge0:D_a\subset\{X_{H_a+s}:0\le s\le u\}\}.
\]
Thus $\tau_a$ is the time needed to cover $D_a$ after its
representative is first hit. Proposition~\ref{chI:prop:local-cover},
the strong Markov property at $H_a$, and
$\operatorname{diam}_d(D_a)\le\rho_N<N/2$ imply
\begin{equation}\label{chVIII:eq:clearmean}
 \EE_y\tau_a=\EE_aT_{D_a}
 \le Cn(1+\log(\rho_N+1))^2\le CnL^{1/3}.
\end{equation}
Define the clearing allowance
\[
 \eta_N=L^{-1/12}.
\]
Since $|\operatorname{Rep}_N(U)|\le B\sqrt L$, Markov's inequality and a union
bound give
\begin{equation}\label{chVIII:eq:cleartail}
 \PP_y(\exists a\in\operatorname{Rep}_N(U):\tau_a>\eta_N b_N)
 \le\frac{\sum_{a\in\operatorname{Rep}_N(U)}\EE_y\tau_a}{\eta_Nb_N}
 \le CB\eta_N.
\end{equation}
For $w\in\mathbb R$ with $\bar u_N(w-\eta_N)\ge0$, on the same path,
\[
 \{T_{\operatorname{Rep}_N(U)}\le \bar u_N(w-\eta_N)\}
 \cap\bigcap_a\{\tau_a\le\eta_N b_N\}
 \subset\{T_{U^\circ}\le \bar u_N(w)\}.
\]
For bounded $w$ and sufficiently large $N$, the single-point
bound \eqref{chVIII:eq:upper} gives
\[
 \PP_y(T_{U\setminus U^\circ}>\bar u_N(w))
 \le\sum_{x\in U\setminus U^\circ}\PP_y(H_x>\bar u_N(w))
 \le C_Wr_N(U).
\]
Writing $F_N^{\rm rep}(y,U,w)=\PP_y(T_{\operatorname{Rep}_N(U)}\le \bar u_N(w))$,
we obtain
\begin{equation}\label{chVIII:eq:sandwich}
 \begin{split}
 F_N^{\rm rep}(y,U,w-\eta_N)-CB\eta_N-C_Wr_N(U)
 &\le\PP_y(T_U\le \bar u_N(w))\\
 &\le F_N^{\rm rep}(y,U,w).
 \end{split}
\end{equation}
The upper bound follows from $\operatorname{Rep}_N(U)\subset U$.

Set
\[
 e_N(B,W+2):=
 \sup_{\substack{y\in V_N,\ U\subset V_N\\
 M_N^{\rm comp}(U)\le B,\ |w|\le W+2}}
 \left|F_N^{\rm rep}(y,U,w)-\Phi_{M_N^{\rm comp}(U)}(w)\right|.
\]
This deterministic number lies in $[0,1]$ and tends to zero by
\eqref{chVIII:eq:void}.
For every $m\ge0$,
\begin{equation}\label{chVIII:eq:modulus}
 0\le\Phi_m'(z)=me^{-z}e^{-me^{-z}}\le e^{-1},\qquad
 |\Phi_m(z)-\Phi_m(w)|\le e^{-1}|z-w|.
\end{equation}
Insert $w=z-\delta_N^{\rm shift}$ in \eqref{chVIII:eq:sandwich}. For large
$N$, both $w$ and $w-\eta_N$ belong to $[-W-2,W+2]$, and
$\bar u_N(w-\eta_N)>0$, $\bar u_N(w)>0$.
Equations
\eqref{chVIII:eq:clock}, \eqref{chVIII:eq:void} and \eqref{chVIII:eq:modulus}
give
\begin{equation}\label{chVIII:eq:totalbudget}
 \left|\PP_y(T_U\le u_N(z))-\Phi_{M_N^{\rm comp}(U)}(z)\right|
 \le e_N(B,W+2)+CB\eta_N
          +\frac{\delta_N^{\rm shift}+\eta_N}{e}+C_{W+1}r_N(U).
\end{equation}
First choose the Bonferroni order $J$ in \eqref{chVIII:eq:truncation};
with $J$ fixed, take $N$ large in its finitely many moment estimates.
This makes $e_N(B,W+2)$ arbitrarily small. Since
$\eta_N,\delta_N^{\rm shift}\to0$, the first three terms tend to zero, proving
\eqref{chVIII:eq:configuration} with $C=C_{W+1}$.
\end{proof}

\subsection*{8.5. Restarting at the central time}
In the rest of the chapter, the original walk has law
$\PP_{\pi_N}$ unless an initial law is stated explicitly. Recall
$C_N=C_N(s_N)$ from Chapter~VII and define
\[
 U_N=\{x:\ell_{s_N}(x)=0\},\qquad Y_N=X_{s_N},\qquad
 R_N^{\rm rem}=r_N(U_N).
\]
These are measurable with respect to the history
$\mathcal F_{s_N}$, and $M_N^{\rm comp}(U_N)=C_N$, with $C_N$ exactly
as in Chapter~VII. By \eqref{chVII:eq:zerosvisited}, $U_N$ is almost surely the
unvisited set at time $s_N$.

Define, for all $z\in\mathbb R$,
\begin{equation}\label{chVIII:eq:kernel}
 F_N^{\rm cond}(z)=k_N(Y_N,U_N,z),\qquad
 k_N(y,U,z)=\PP_y(T_U\le [u_N(z)]_+),
 \qquad [t]_+=\max(t,0).
\end{equation}
For $u_N(z)>0$ this is the restart kernel at $s_N$. For every fixed
compact $z$-interval, $[u_N(z)]_+=u_N(z)$ for all sufficiently large
$N$.

\begin{proposition}[The conditional kernel]\label{chVIII:prop:restart}
For each $z$ with $u_N(z)>0$, under any initial law,
\begin{equation}\label{chVIII:eq:restart}
 F_N^{\rm cond}(z)=\EE\bigl[\one_{\{T_N\le t_N(z)\}}
                         \mid\mathcal F_{s_N}\bigr]
 \quad\hbox{a.s.}
\end{equation}
\end{proposition}
\begin{proof}
The map $(Y_N,U_N)$ takes values in a finite set, so
$F_N^{\rm cond}(z)$ is past measurable and lies in $[0,1]$.
Given the prefix, append an independent walk $X'$ from
$Y_N$, for duration $u_N(z)$, and write
\[
 \ell'_t(x):=\int_0^t\one_{\{X'_s=x\}}\,ds,
 \qquad t\ge0,\quad x\in V_N.
\]
The deterministic-time Markov property, including the exponential residual
holding time at $s_N$, gives the law of the original
future. Occupation adds under this concatenation:
\[
 \ell_{s_N+u_N(z)}(x)=\ell_{s_N}(x)+\ell'_{u_N(z)}(x).
\]
Both terms are nonnegative. With $T_{U_N}(X')$ denoting the
cover time evaluated on the appended path, almost surely
\[
 \begin{split}
 \{\forall x\in V_N:\ell_{s_N+u_N(z)}(x)>0\}
 &=\{\forall x\in U_N:\ell'_{u_N(z)}(x)>0\}\\
 &=\{T_{U_N}(X')\le u_N(z)\}.
 \end{split}
\]
The last equality uses $u_N(z)>0$.

To retain the conditioning information, let $w$ be any
bounded nonnegative $\mathcal F_{s_N}$-measurable function.
Concatenation leaves $w$ unchanged. Integrating first
over the independent suffix therefore gives
\begin{equation}\label{chVIII:eq:weightedrestart}
 \EE[w\one_{\{T_N\le t_N(z)\}}]=\EE[wF_N^{\rm cond}(z)].
\end{equation}
This is the defining conditional-expectation identity
and proves \eqref{chVIII:eq:restart}.
\end{proof}

Let $Z=\mathcal Z(\mathbb T^2)$ be the total mass of the critical
measure constructed in Chapter~II, and let $\kappa=\beta\Gamma>0$
be the deterministic constant of Chapter~VII. The two scalar
inputs, under $\PP_{\pi_N}$, are
\begin{equation}\label{chVIII:eq:inputs}
 C_N\Longrightarrow\kappa Z,
 \qquad R_N^{\rm rem}\longrightarrow0\quad\hbox{in probability},
 \qquad 0<Z<\infty\ \hbox{a.s.}
\end{equation}
By Theorem~\ref{chVII:thm:limit},
$C_N\Rightarrow\kappa Z$ with $\kappa=\beta\Gamma>0$ and
$Z=\mathcal Z(\Ttwo)$ as defined in Chapter~II. For the second, the deterministic
inclusion
\[
 U_N\setminus U_N^\circ
 \subset D_N^{\rm can}(U_1(s_N))
 \subset D_N^{\rm wide}(U_1(s_N))
\]
holds for large $N$. The first inclusion uses $U_N\subset
U_1(s_N)$; the second compares the canonical graph-distance
thresholds with the wider Euclidean thresholds of Chapter~IV.
Its vanishing normalized cardinality gives the assertion
about $R_N^{\rm rem}$. The bounded-error argument below uses the
vanishing of its fixed-threshold tail probabilities.

For $W\ge0$, define the random uniform error
\[
 E_{N,W}:=\sup_{|z|\le W}|F_N^{\rm cond}(z)-e^{-C_Ne^{-z}}|.
\]

\begin{proposition}[Uniform conditional error]\label{chVIII:prop:window}
For every fixed $W\ge0$,
\begin{equation}\label{chVIII:eq:windowerror}
 \EE_{\pi_N}\left[
   \sup_{|z|\le W}|F_N^{\rm cond}(z)-e^{-C_Ne^{-z}}|\right]
 \longrightarrow0.
\end{equation}
\end{proposition}
\begin{proof}
The variable $E_{N,W}$ is a measurable function of
the finite-valued pair $(Y_N,U_N)$: for each deterministic
pair the supremum is just a real number. Since both
functions being compared belong to $[0,1]$,
$0\le E_{N,W}\le1$, so $E_{N,W}$ is integrable.

The sequence $(C_N)$ is tight by \eqref{chVIII:eq:inputs}. It
also follows directly from the Laplace-transform convergence in
\eqref{chVII:eq:distribution}:
\begin{equation}\label{chVIII:eq:laplacetight}
 \limsup_N\PP_{\pi_N}(C_N>B)
 \le\frac{1-\EE e^{-\kappa Z/B}}{1-e^{-1}}
 \longrightarrow0\quad(B\to\infty).
\end{equation}
Indeed $(1-e^{-1})\one_{\{C_N>B\}}\le1-e^{-C_N/B}$;
first take $N\to\infty$, then use bounded convergence
and $Z<\infty$.

Fix $\epsilon>0$. Choose $B$ so that
$\PP_{\pi_N}(C_N>B)\le\epsilon/8$ eventually.
Apply Theorem~\ref{chVIII:thm:configuration} at this $B,W$
with deterministic error $\epsilon/4$, obtaining $C>0$.
On $\{C_N\le B,R_N^{\rm rem}\le\epsilon/(4C)\}$,
the deterministic estimate holds simultaneously for every
$|z|\le W$, and hence $E_{N,W}\le\epsilon/2$.
By \eqref{chVIII:eq:inputs}, $\PP_{\pi_N}(R_N^{\rm rem}>\epsilon/(4C))\le\epsilon/8$
eventually. The unit bound on the complementary event gives
\[
 \EE_{\pi_N}E_{N,W}
 \le\epsilon/2+\PP_{\pi_N}(C_N>B)+\PP_{\pi_N}(R_N^{\rm rem}>\epsilon/(4C))
 \le3\epsilon/4.
\]
This proves \eqref{chVIII:eq:windowerror}.
\end{proof}

\subsection*{8.6. The cover-time limit}
Define
\begin{equation}\label{chVIII:eq:limitcdf}
 F(z)=\EE e^{-\kappa Ze^{-z}},\qquad z\in\mathbb R.
\end{equation}
Since $0\le e^{-\kappa Ze^{-z}}\le1$ and $0<Z<\infty$ almost surely,
this expectation is well defined on the probability space of the
critical measure.

\begin{theorem}[Cover-time distribution]\label{chVIII:thm:main}
For the rate-one walk, the normalized variable is exactly
\[
 \frac{T_N-a_N}{b_N}
 =\frac{T_N}{gN^2\log N}-2\log N+\log\log N.
\]
With the constant
$\kappa=\beta\Gamma>0$ of Chapter~VII,
\begin{equation}\label{chVIII:eq:main}
 \sup_{x\in V_N}\sup_{z\in\mathbb R}
 \left|\PP_x\left(\frac{T_N-a_N}{b_N}\le z\right)
                          -\EE e^{-\kappa Ze^{-z}}\right|
 \longrightarrow0.
\end{equation}
Equivalently, the limit has the law of
$G+\log(\kappa Z)$, where $G$ is a standard Gumbel random
variable independent of $Z$.
\end{theorem}
\begin{proof}
For a fixed $z$, $u_N(z)>0$ eventually. By
Proposition~\ref{chVIII:prop:restart} and
\eqref{chVIII:eq:windowerror},
\[
 \left|\PP_{\pi_N}(T_N\le t_N(z))
               -\EE_{\pi_N}e^{-C_Ne^{-z}}\right|
 \le\EE_{\pi_N}E_{N,|z|}\longrightarrow0.
\]
The Laplace limit of Chapter~VII, at the fixed positive
parameter $q=e^{-z}$, identifies the limit as $F(z)$.

The distribution of $T_N$ is the same for every initial
vertex. Indeed, construct walks from $x$ and $y$ using
the same holding times and jump directions. Their
positions differ at all times by $y-x$. Translation
permutes the torus, so their cover times agree pathwise.
Their common distribution is also that of the stationary
mixture. Thus the pointwise convergence just proved is
already uniform in the initial vertex.

The function $F$ is nondecreasing and continuous by
bounded convergence; in fact \eqref{chVIII:eq:modulus} gives
$|F(z)-F(w)|\le |z-w|/e$. Since $0<\kappa Z<\infty$
almost surely, bounded convergence also gives
\begin{equation}\label{chVIII:eq:tails}
 \lim_{z\to-\infty}F(z)=0,\qquad
 \lim_{z\to\infty}F(z)=1.
\end{equation}
To upgrade to uniform convergence on $\mathbb R$, fix
$\epsilon>0$ and choose $a<b$ such that $F(a)<\epsilon$
and $1-F(b)<\epsilon$. Choose a finite grid
$a=z_0<\cdots<z_m=b$ with
$F(z_{j+1})-F(z_j)<\epsilon$. At all these finitely
many points the common distribution function
$F_N^{\rm cov}(z):=\PP_{\pi_N}((T_N-a_N)/b_N\le z)$
differs from $F$ by at most $\epsilon$ for large
$N$. For $z_j\le z\le z_{j+1}$, monotonicity places
$F_N^{\rm cov}(z)$ between its endpoint values and shows
$|F_N^{\rm cov}(z)-F(z)|\le2\epsilon$. For $z<a$ and $z>b$,
the same bound follows from the endpoint estimates,
the range $[0,1]$, and the chosen tails. The same endpoint bounds
control the range $u_N(z)\le0$, completing the uniform convergence in
$z$ and proving \eqref{chVIII:eq:main}.

Finally take an independent $G$ with
$\PP(G\le w)=e^{-e^{-w}}$. Conditional on $Z$,
\[
 \PP(G+\log(\kappa Z)\le z\mid Z)
 =\exp\{-e^{-(z-\log(\kappa Z))}\}
 =e^{-\kappa Ze^{-z}}.
\]
Taking expectations shows that the limiting law is
$G+\log(\kappa Z)$.
\end{proof}

\clearpage
\endgroup
\MergedMainfalse
\section*{Appendices}
\mergedcontentsline{toc}{section}{Appendices}
\CompactAppendixHeadings
\begingroup
\mergedchaptersetup{I}{Random walk and potential theory}
\appendix
\renewcommand{\theequation}{A.\arabic{equation}}
\renewcommand{\theHequation}{A.\arabic{equation}}
\setcounter{equation}{0}
\section*{Appendix A. Proof of Lemma~\ref{chI:lem:fourier-sums}}
\hypertarget{proof.appendix.A}{}
\mergedcontentsline{toc}{subsection}{Appendix A. Proof of Lemma~\ref{chI:lem:fourier-sums}}
\label{chI:app:fourier}

The notation for rows and frequency cutoffs is local to this proof. Write
$I_N=\ZZ\cap[-N/2,N/2)$, so $\cK_N=I_N^2$.

\begin{proof}
We begin with the diagonal sum. For $j\ge1$, monotonicity of
$(j^2+t^2)^{-1}$ on $t\ge0$ gives
\[
 \frac{\pi}{4j}=\int_0^j\frac{dt}{j^2+t^2}
 \le\sum_{a=0}^{j-1}\frac1{j^2+a^2}
 \le\frac{\pi}{4j}+\frac1{j^2}.
\]
The square shell $\{k\in\ZZ^2:|k|_\infty=j\}$ has reciprocal
square sum
\[
 \sum_{|k|_\infty=j}|k|^{-2}
 =8\sum_{a=0}^{j-1}\frac1{j^2+a^2}-\frac2{j^2}
 =\frac{2\pi}{j}+O(j^{-2}).
\]
The counting identity used here is
\[
 \sum_{|k|_\infty=j}|k|^{-2}
 =\frac4{j^2}+\frac4{2j^2}
       +8\sum_{a=1}^{j-1}\frac1{j^2+a^2}.
\]
Together with $\sum_{j=1}^m j^{-1}=\log m+O(1)$ and
$\sum_{j\ge1}j^{-2}<\infty$, it gives
\[
 \sum_{0<|k|_\infty\le m}|k|^{-2}=2\pi\log m+O(1).
\]
Since
\[
 \begin{gathered}
 \{k:|k|_\infty\le\lfloor(N-1)/2\rfloor\}\subseteq\mathcal K_N
 \subseteq\{k:|k|_\infty\le\lfloor N/2\rfloor\},\\
 N/4\le\lfloor(N-1)/2\rfloor\le\lfloor N/2\rfloor\le N/2,
 \end{gathered}
\]
both bounding sums equal $2\pi\log N+O(1)$, proving the first assertion.

For the off-diagonal estimate, choose a shortest representative
$(z_1,z_2)$ of $z\ne0$ and interchange coordinates if necessary
so that $|z_1|\ge d(z,0)/\sqrt2$. Put
$\theta=2\pi z_1/N$. For $b\ge1$, consider the row
\[
 W_b=\sum_{a\in I_N}\frac{e^{ia\theta}}{a^2+b^2}.
\]
Its absolute sum satisfies
\begin{equation}\label{chI:eq:row-absolute}
 |W_b|\le b^{-2}+2\int_0^\infty\frac{dt}{t^2+b^2}
          \le C/b.
\end{equation}
To use oscillation, note that any consecutive geometric sum has
absolute value at most $2/|1-e^{i\theta}|$. If $w_j$ is nonnegative
and decreasing, summation by parts therefore gives
\[
 \left|\sum_{j=1}^m w_je^{ij\theta}\right|
 \le\frac{2w_1}{|1-e^{i\theta}|}.
\]
Writing $U_j=\sum_{a=1}^je^{ia\theta}$, the identity is
\[
 \sum_{j=1}^m w_je^{ij\theta}
 =w_mU_m+\sum_{j=1}^{m-1}(w_j-w_{j+1})U_j,\qquad
 w_m+\sum_{j=1}^{m-1}(w_j-w_{j+1})=w_1.
\]
Apply it to each of the positive and negative parts of $I_N$,
with $w_j=(j^2+b^2)^{-1}$.
The two parts may have different lengths when $N$ is even; the
same estimate applies to each. Since
\[
 |1-e^{i\theta}|=2\sin(\pi|z_1|/N)\ge4|z_1|/N,
\]
we obtain the complementary bound
\begin{equation}\label{chI:eq:row-oscillatory}
 |W_b|\le b^{-2}\left(1+\frac4{|1-e^{i\theta}|}\right)
          \le C\frac N{|z_1|b^2}.
\end{equation}
The row with $b=0$, with its zero term omitted, has absolute sum
at most $2\sum_{j\ge1}j^{-2}\le C$.

Now group $\Sigma_N(z)$ by its second coordinate. There are at most
two rows with each positive $|b|$. Since $N/|z_1|\ge2$, the row bounds give
\[
 \begin{aligned}
 |\Sigma_N(z)|
 &\le C+C\sum_{b=1}^{\lfloor N/2\rfloor}
                  \min\{b^{-1},N/(|z_1|b^2)\}\\
 &\le C+C\sum_{1\le b\le N/|z_1|}b^{-1}
       +\frac{CN}{|z_1|}\sum_{b>N/|z_1|}b^{-2}\\
 &\le C\bigl(1+\log(N/|z_1|)\bigr).
 \end{aligned}
\]
The last line uses, for $u\ge2$,
\[
 \sum_{1\le b\le u}\frac1b\le1+\log u,\qquad
 \sum_{b>u}\frac1{b^2}
 \le\int_{\lfloor u\rfloor}^{\infty}\frac{dt}{t^2}\le\frac2u.
\]
Since $|z_1|\ge d(z,0)/\sqrt2$,
$\log(N/|z_1|)\le\log(N/d(z,0))+\tfrac12\log2$.
This proves the second assertion of \eqref{chI:eq:quadratic-sums}.

Finally, on the square shell $|k|_\infty=j$ there are at most
$8j$ points and $|k|\ge j$. Thus, for $u\ge1$,
\[
 \begin{aligned}
 \sum_{k\ne0}\min\{u^{-2},|k|^{-2}\}
 &\le8\sum_{j=1}^{\lfloor N/2\rfloor}
                      j\min\{u^{-2},j^{-2}\}\\
 &\le8u^{-2}\sum_{1\le j\le u}j
       +8\sum_{u<j\le N/2}j^{-1}\\
 &\le C\bigl(1+\log\max\{N/u,1\}\bigr).
 \end{aligned}
\]
The first sum is $O(1)$, and the second is bounded by a harmonic
integral. This proves the last assertion.
\end{proof}

\renewcommand{\theequation}{B.\arabic{equation}}
\renewcommand{\theHequation}{B.\arabic{equation}}
\setcounter{equation}{0}
\renewcommand{\thelemma}{B.\arabic{lemma}}
\renewcommand{\theHlemma}{B.\arabic{lemma}}
\setcounter{lemma}{0}
\section*{Appendix B. Proof of Lemma~\ref{chI:lem:unit-differences}}
\hypertarget{proof.appendix.B}{}
\mergedcontentsline{toc}{subsection}{Appendix B. Proof of Lemma~\ref{chI:lem:unit-differences}}
\label{chI:app:differences}

We prove Lemma~\ref{chI:lem:unit-differences} using periodic
summation by parts on the frequency torus; all frequency differences
are interpreted modulo $N$.

Choose an even $C^\infty$ function $\eta:\mathbb R\to[0,1]$,
equal to one on $[-1,1]$, zero outside $(-2,2)$, and nonincreasing
on $[0,\infty)$. Write
\[
 \rho(v_1,v_2)=\eta(v_1)\eta(v_2),\qquad
 \alpha(v)=\rho(v)-\rho(2v).
\]
Then $\alpha$ is smooth, vanishes when $|v|_\infty\le1/2$ or
$|v|_\infty\ge2$, and is even in each coordinate. In this
appendix $R$ denotes a frequency radius. For a real frequency
$\xi=(\xi_1,\xi_2)$, put
\[
 \sigma_R(\xi)_i=\frac N{\pi R}\sin\frac{\pi\xi_i}{N},\qquad
 F_R(\xi)=\rho(\sigma_R(\xi)),\qquad
 A_R(\xi)=F_R(\xi)-F_{R/2}(\xi)=\alpha(\sigma_R(\xi)).
\]
These functions are $N$-periodic in each coordinate, since
translation by $N$ changes the sign of the sine and $\rho$ is
even in each coordinate. This choice of cutoff is useful even
when a block reaches the edge of the frequency square.

Fix $m\in\{1,2\}$ and unit directions $e_1^*,\ldots,e_m^*$,
where the stars distinguish these directions from the coordinate
basis. Define, for $k\ne0$,
\[
 a_R(k)=\frac{A_R(k)}{N^2\lambda_N(k)}
                  \prod_{\nu=1}^m(\chi_k(e_\nu^*)-1),
 \qquad a_R(0)=0.
\]
For an array $a:\cK_N\to\mathbb C$, write
\[
 \widehat a(z)=\sum_{k\in\cK_N}a(k)\chi_k(z),\qquad
 \Delta_i a(k)=a(k+e_i)-a(k),
\]
where addition of frequencies is modulo $N$.
The transform here is the unnormalized finite sum, and
$\Delta_i^p$ denotes the $p$-fold iterate, with $\Delta_i^0a=a$.

\begin{lemma}\label{chI:lem:block}
Uniformly in $R\ge1/2$, $N\ge3$, and $z\ne0$,
\begin{equation}\label{chI:eq:block-envelope}
 |\widehat a_R(z)|\le C(R/N)^m
                \min\{1,(Rd(z,0)/N)^{-3}\}.
\end{equation}
\end{lemma}

\begin{proof}
We first verify the frequency difference estimates
\begin{equation}\label{chI:eq:block-l1}
 \sum_{k\in\cK_N}|\Delta_i^p a_R(k)|
       \le C(R/N)^m R^{-p},\qquad p=0,3.
\end{equation}
If $R\ge N$, then $F_R=F_{R/2}=1$, since both arguments of
$\rho$ have absolute value at most $2/\pi<1$; thus $a_R=0$.
It suffices to consider $1/2\le R<N$.

Here are the derivative and support details for
\eqref{chI:eq:block-l1}. Put $\delta=R/N$ and
\[
 S_\delta(v)_i=\begin{cases}
       \sin(\pi\delta v_i)/(\pi \delta),&\delta\ne0,\\
       v_i,&\delta=0.
       \end{cases}
\]
Define the smooth function
\[
 H(w)=\frac{\alpha(w)}{\pi^2|w|^2},
\]
extended by zero near $w=0$. For a unit direction $e$, also set
\[
 J_{\delta,e}(v)=\begin{cases}
 (e^{2\pi i \delta e\cdot v}-1)/\delta,&\delta\ne0,\\
 2\pi i e\cdot v,&\delta=0.
 \end{cases}
\]
For all real $\delta$, including $\delta=0$,
\[
 J_{\delta,e}(v)=2\pi i(e\cdot v)
       \int_0^1e^{2\pi it\delta e\cdot v}\,dt,\qquad
 S_\delta(v)_i=v_i\int_0^1\cos(\pi t\delta v_i)\,dt.
\]
These formulas give joint smoothness in $(\delta,v)$.
Extend $\lambda_N$ to real frequencies by the same sine formula,
and for each lattice unit direction $e^*\in\mathbb Z^2$ put
$\chi_\xi(e^*)=\exp(2\pi i\xi\cdot e^*/N)$. The smooth periodic
extension of the multiplier can then be written as
\begin{equation}\label{chI:eq:scaled-block}
 a_R(\xi)=\frac{(R/N)^m}{R^2}\,
       H(S_\delta(\xi/R))\prod_{\nu=1}^m J_{\delta,e_\nu^*}(\xi/R).
\end{equation}
Indeed, $\lambda_N(\xi)=\pi^2\delta^2|S_\delta(\xi/R)|^2$;
the cutoff removes every zero of this symbol.

For centered $\xi$, we have $|\delta(\xi_i/R)|\le1/2$ and
\[
 |S_\delta(\xi/R)_i|\ge(2/\pi)|\xi_i/R|.
\]
Where $H(S_\delta(v))$ or any of its derivatives is nonzero, its
argument is in the closed square $|S_\delta(v)|_\infty\le2$.
It follows that $|v_i|\le\pi$. On the compact set
$0\le \delta\le1$, $|v_i|\le\pi$, all derivatives through order
three of the smooth product in \eqref{chI:eq:scaled-block} are
bounded independently of $\delta$. Outside this set the product and
these derivatives vanish on the centered fundamental square.
Hence, throughout the frequency torus,
\begin{equation}\label{chI:eq:block-derivative}
 |\partial_i^p a_R(\xi)|\le
 C(R/N)^m R^{-2-p},\qquad 0\le p\le3.
\end{equation}
The periodic extension is smooth at the seams: $H$ is even in
each sine coordinate, while the exponential increments are
$N$-periodic. Thus the derivative estimate holds also along
segments crossing a seam.

Iterating the fundamental theorem of calculus gives
\[
 \Delta_i^p a_R(k)=\int_{[0,1]^p}
   \partial_i^p a_R(k+(t_1+\cdots+t_p)e_i)\,dt_1\cdots dt_p
 \quad(p\ge1).
\]
For $0\le p\le3$ and $R\ge\tfrac12$, with all translations
taken modulo $N$,
\[
 \begin{aligned}
 \supp(\Delta_i^pa_R)&\subseteq
       \bigcup_{j=0}^p(\supp a_R-je_i),\\
 \#\supp(\Delta_i^pa_R)&\le C(R+p+1)^2\le CR^2.
 \end{aligned}
\]
Indeed $\supp a_R$ has a centered representative in
$[-\pi R,\pi R]^2$; if this square covers the torus, then
$N^2\le CR^2$. The preceding integral identity and
\eqref{chI:eq:block-derivative} give
\[
 \sum_k|\Delta_i^pa_R(k)|
 \le CR^2\,C(R/N)^mR^{-2-p}
 \le C(R/N)^mR^{-p}.
\]
This proves \eqref{chI:eq:block-l1}.

Choose a coordinate $i$ for which a shortest representative of
$z$ satisfies $|z_i|\ge d(z,0)/\sqrt2$. A change of variable in
the finite periodic sum gives
\[
 \widehat{\Delta_i a_R}(z)
       =(e^{-2\pi i z_i/N}-1)\widehat a_R(z).
\]
The factor has absolute value at least $4|z_i|/N$, hence at
least $C^{-1}d(z,0)/N$. Applying the identity three times and using
\eqref{chI:eq:block-l1} yields
\[
 |\widehat a_R(z)|\le
 C(R/N)^m R^{-3}(N/d(z,0))^3.
\]
The case $p=0$ of \eqref{chI:eq:block-l1} supplies the bound
$C(R/N)^m$. Taking the smaller of the two proves
\eqref{chI:eq:block-envelope}.
\end{proof}

\begin{proof}[Proof of Lemma~\ref{chI:lem:unit-differences}]
Take dyadic radii $R_j=2^{j-1}$, $j=0,\ldots,J$, where
$R_J\ge N$. For every $k\in\cK_N\setminus\{0\}$, at least
one coordinate has absolute value at least one, and
\[
 \left|\frac N{\pi/4}\sin\frac{\pi k_i}{N}\right|
 \ge\frac8\pi>2.
\]
Thus $F_{1/4}(k)=0$, while $F_{R_J}(k)=1$. Telescoping the
cutoffs gives $\sum_{j=0}^J A_{R_j}(k)=1$. Therefore
\[
 \frac1{N^2\lambda_N(k)}\prod_{\nu=1}^m(\chi_k(e_\nu^*)-1)
       =\sum_{j=0}^J a_{R_j}(k)\qquad(k\ne0).
\]
By \eqref{chI:eq:green}, the transform of the left side, extended by
zero at $k=0$, is precisely the first or second spatial difference
in \eqref{chI:eq:unit-differences}. For $d(z,0)>0$,
Lemma~\ref{chI:lem:block} bounds its absolute value by
\[
 C\sum_{j=0}^J(R_j/N)^m
       \min\{1,(R_jd(z,0)/N)^{-3}\}.
\]
Since $R_{j+1}=2R_j$ and $m\in\{1,2\}$,
\[
 \begin{aligned}
 \sum_{R_jd(z,0)\le N}(R_j/N)^m
       &\le\frac{d(z,0)^{-m}}{1-2^{-m}},\\
 d(z,0)^{-3}\sum_{R_jd(z,0)>N}(R_j/N)^{m-3}
       &\le\frac{d(z,0)^{-m}}{1-2^{m-3}}.
 \end{aligned}
\]
Each empty sum is zero. Substitution proves
\eqref{chI:eq:unit-differences}, uniformly over the unit directions.
\end{proof}


\par\endgroup

\begingroup
\mergedchaptersetup{II}{Gaussian fields and the critical measure}
\appendixchapter{C}{Proof of Theorem~\ref{chII:thm:critical-measure}}

We use \cite[Theorem~10]{chII:DRSV} for existence, nontriviality and
nonatomicity of the independent-increment Dirichlet cutoff. The heat
replacement and field-measurability steps are proved below from the
exact covariance identity \eqref{chII:eq:R26.2} and the Gaussian tail
sigma-fields.
In Appendices~C--H, unsubscripted $\mathbb P,\mathbb E$ (also
written $P,E$) denote the Gaussian law specified in the local
construction; walk laws retain their starting-point subscripts. We write
$\mathbb E[Y;A]=\mathbb E[Y\mathbf1_A]$.

\subsection{Hilbert-space construction and fixed-scale smoothing}
For $D=(0,\ell)^2$, up to translation, let
\[
 \lambda_{m,n}=\pi^2(m^2+n^2)/\ell^2,\qquad
 e_{m,n}(x)=\frac2\ell\sin(\pi mx_1/\ell)\sin(\pi nx_2/\ell),
 \quad m,n\ge1.
\]
These functions form an orthonormal basis of $L^2(D,dx)$.
Use a single index $j$ when convenient and define the spectral Sobolev
space $\mathcal H_D^s$ by completion in the norm
\[
 \|u\|_{\mathcal H_D^s}^2=\sum_j(1+\lambda_j)^s|\langle u,e_j\rangle|^2,
 \qquad P_t^D=e^{t\Delta_D}.
\]
For $s<0$, the bracket is the distributional pairing extending the
$L^2(D,dx)$ inner product. The Gaussian integrals below are taken on
the probability space $(\Omega,\mathcal F,P)$ carrying the white noise.
With space--time white noise $W$ on $(0,\infty)\times D$, define
\begin{equation}\label{chII:eq:R26.1}
 h_D=2\int_0^\infty P_{s/2}^D W(ds),\qquad
 X_\varepsilon=2\int_{2\varepsilon^2}^\infty P_{s/2}^D W(ds),\qquad
 Y_\varepsilon=P_{\varepsilon^2}^Dh_D.
\end{equation}
The covariance of $h_D$ is $4(-\Delta_D)^{-1}$ and has logarithmic coefficient
$g=2/\pi$ in the normalization used throughout the manuscript.

\begin{lemma}[Negative-Sobolev convergence and smoothing]\label{chII:lem:dirichlet-smoothing}
Fix $0<\delta<1/2$. The variables in \eqref{chII:eq:R26.1} are well defined in
$\mathcal H_D^{-\delta}$ and
\[
 E\|h_D\|_{\mathcal H_D^{-\delta}}^2
 =4\sum_j\frac{(1+\lambda_j)^{-\delta}}{\lambda_j}<\infty.
\]
For every decreasing sequence $\varepsilon_k\downarrow0$,
$X_{\varepsilon_k}\to h_D$ almost surely and in
$L^2(\Omega;\mathcal H_D^{-\delta})$.  Moreover
$Y_\varepsilon\to h_D$ almost surely in $\mathcal H_D^{-\delta}$ as
$\varepsilon\downarrow0$.  For fixed $t>0$ and integer $m\ge0$,
$P_t^D:\mathcal H_D^{-\delta}\to C^m(\overline D)$ is bounded.
\end{lemma}
\begin{proof}
The two-dimensional eigenvalue growth gives the displayed summability.  If
$\mathcal F_k$ is the sigma-field of the white noise on
$[2\varepsilon_k^2,\infty)\times D$, then
$X_{\varepsilon_k}=E(h_D\mid\mathcal F_k)$ as a square-integrable Hilbert-space
random variable.  The increasing Hilbert martingale converges almost surely and in
$L^2$ to $h_D$ because the union of these sigma-fields reveals all white noise entering
$h_D$.  The heat convergence follows directly from the spectral multiplier
$e^{-\varepsilon^2\lambda_j}\to1$ and dominated convergence in
$\mathcal H_D^{-\delta}$.

For $|\beta|\le m$, the sine basis satisfies
$\|\partial^\beta e_j\|_\infty\le C_m\lambda_j^{m/2}$.  Hence
\[
 \|P_t^Du\|_{C^m}
 \le C_m\left(\sum_j e^{-2t\lambda_j}(1+\lambda_j)^{m+\delta}\right)^{1/2}
 \|u\|_{\mathcal H_D^{-\delta}},
\]
which proves the smoothing assertion.
\end{proof}

The exact covariance computation is
\begin{equation}\label{chII:eq:R26.2}
 \Cov Y_\varepsilon
 =P_{\varepsilon^2}^D4(-\Delta_D)^{-1}P_{\varepsilon^2}^D
 =4\int_{2\varepsilon^2}^\infty P_s^D\,ds
 =\Cov X_\varepsilon.
\end{equation}
Thus $X_\varepsilon\overset d=Y_\varepsilon$ as smooth Gaussian fields for every
fixed $\varepsilon$.

\subsection{The critical limit is measurable from the total field}
For $U_\varepsilon$ a centered smooth Gaussian field, write
\[
 \mathcal M_\varepsilon[U](dx)=\sqrt{\log(1/\varepsilon)}
 \exp\{\alpha U_\varepsilon(x)-\alpha^2\EE[U_\varepsilon(x)^2]/2\}\,dx.
\]
By \cite{chII:DRSV}, Theorem 10, after the deterministic normalization convention fixed here,
\begin{equation}\label{chII:eq:drsv-cutoff-limit}
 \mathcal M_\varepsilon[X]\longrightarrow M_D
 \quad\text{locally weakly in probability}.
\end{equation}
The limit is nonzero, locally finite and nonatomic.

Expand the space--time white noise in the Dirichlet eigenbasis.  Let
$W_j$ be independent real white noises on $(0,\infty)$ and set
\[
 k_j(s)=2e^{-\lambda_js/2},\qquad
 \mathfrak H=\bigoplus_{j\ge1}L^2(0,\infty),\qquad
 \mathcal W(v)=\sum_j\int_0^\infty v_j(s)\,dW_j(s).
\]
Here $\mathbf e_j$ denotes the $j$th coordinate vector of the Hilbert direct sum;
$\mathcal W$ is its isonormal Gaussian process. Then
$\langle h_D,e_j\rangle=\mathcal W(k_j\mathbf e_j)$. For $a>0$ define the
closed subspace
\begin{equation}\label{chII:eq:tail-subspaces}
 K_a=\bigoplus_{j\ge1}\Bigl(L^2(0,a)\oplus
      \operatorname{span}\{k_j\one_{[a,\infty)}\}\Bigr),\qquad
 K_0=\bigoplus_{j\ge1}\operatorname{span}\{k_j\}.
\end{equation}
The first summand in each coordinate is understood as functions supported in $(0,a)$.

\begin{lemma}[Intersection of the Gaussian tail sigma-fields]\label{chII:lem:gaussian-tail-intersection}
For any deterministic $a_n\downarrow0$, the spaces $K_{a_n}$ decrease and
\[
 \bigcap_{n\ge1}K_{a_n}=K_0.
\]
If $\mathcal G_a$ denotes the completion of
$\sigma(\mathcal W(v):v\in K_a)$, then
\begin{equation}\label{chII:eq:gaussian-sigma-intersection}
 \bigcap_{n\ge1}\mathcal G_{a_n}=\mathcal G_0
 \quad\text{modulo null sets}.
\end{equation}
Moreover $\mathcal G_0=\sigma(h_D)$ after completion.
\end{lemma}
\begin{proof}
If $v\in\cap_nK_{a_n}$, then for each coordinate $j$ and each $n$ the restriction
$v_j|_{[a_n,\infty)}$ is a scalar multiple of $k_j$.  The tail intervals overlap and
$k_j$ never vanishes, so the scalar is independent of $n$.  Sending $n\to\infty$
gives $v_j=c_jk_j$ a.e.; square summability gives $v\in K_0$.  The reverse inclusion
is immediate.

For decreasing closed subspaces, the orthogonal projections
$P_{K_{a_n}}v$ converge strongly to $P_{K_0}v$.  For every $v\in\mathfrak H$,
Gaussian regression gives
\[
 E\bigl(e^{i\mathcal W(v)}\mid\mathcal G_{a_n}\bigr)
 =\exp\left\{i\mathcal W(P_{K_{a_n}}v)
 -\tfrac12\|v-P_{K_{a_n}}v\|^2\right\}.
\]
The right side converges in $L^2$ to the same expression with $K_0$.  Reverse
martingale convergence and density of Gaussian exponentials in the Gaussian
$L^2$ space prove \eqref{chII:eq:gaussian-sigma-intersection}.  Finally the variables
$\mathcal W(k_j\mathbf e_j)=\langle h_D,e_j\rangle$ generate $\mathcal G_0$.
\end{proof}

\begin{lemma}[Field measurability of the DRSV limit]\label{chII:lem:critical-limit-measurable}
The random measure $M_D$ in \eqref{chII:eq:drsv-cutoff-limit} is measurable with respect
to $\sigma(h_D)$.  Equivalently, for the Polish space of locally finite measures on
$D$ with the local weak topology, there is a Borel map $F_D$ such that
\begin{equation}\label{chII:eq:critical-graph}
 M_D=F_D(h_D)\quad\text{almost surely}.
\end{equation}
\end{lemma}
\begin{proof}
Fix $a>0$.  If $2\varepsilon^2<a$, then in coordinate $j$ the kernel
$k_j\one_{[2\varepsilon^2,\infty)}$ is the sum of a function supported in $(0,a)$
and a scalar multiple of $k_j\one_{[a,\infty)}$.  Hence
$X_\varepsilon$, and therefore $\mathcal M_\varepsilon[X]$, is $\mathcal G_a$-measurable.
A probability limit of eventually $\mathcal G_a$-measurable random variables is
measurable with respect to the completed $\mathcal G_a$.  Apply this to a countable
convergence-determining family of compactly supported continuous test functions.
Thus $M_D$ is measurable with respect to every $\mathcal G_{a_n}$ and hence, by
Lemma \ref{chII:lem:gaussian-tail-intersection}, with respect to $\mathcal G_0=\sigma(h_D)$.
The Doob--Dynkin lemma on standard Borel spaces yields the map $F_D$.
\end{proof}

\subsection{A graph-supported joint limit forces convergence on the same field}
\begin{lemma}[Collapse onto a measurable graph]\label{chII:lem:graph-collapse}
Let $E,F$ be Polish spaces, with a compatible metric $d_F$ on $F$.
Let $U$ be $E$-valued, let $V_n$ be $F$-valued on the
same probability space and let $F_0:E\to F$ be Borel.  If
\[
 (U,V_n)\Longrightarrow(U,F_0(U)),
\]
then $d_F(V_n,F_0(U))\to0$ in probability.
\end{lemma}
\begin{proof}
Fix $\eta,\rho>0$.  By Lusin's theorem choose a compact $K\subset E$ with
$P(U\in K)>1-\rho$ such that $F_0|_K$ is continuous.  The set
\[
 \{(u,v):u\in K,\ d_F(v,F_0(u))\ge\eta\}
\]
is closed.  Portmanteau gives
$\limsup_nP(U\in K,\ d_F(V_n,F_0(U))\ge\eta)\le0$, because the limit pair is on the graph.
Thus $\limsup_nP(d_F(V_n,F_0(U))\ge\eta)\le\rho$.  Send $\rho\downarrow0$.
\end{proof}

\begin{proposition}[Same-field Dirichlet heat chaos]\label{chII:prop:dirichlet-heat-chaos}
On the probability space of \eqref{chII:eq:R26.1},
\begin{equation}\label{chII:eq:R26.3}
 \mathcal M_\varepsilon[Y]\longrightarrow M_D
 \quad\text{locally weakly in probability as }\varepsilon\downarrow0.
\end{equation}
In particular the heat-chaos limit is the same $h_D$-measurable random measure as in
\eqref{chII:eq:critical-graph}.
\end{proposition}
\begin{proof}
By Lemma \ref{chII:lem:dirichlet-smoothing} and \eqref{chII:eq:drsv-cutoff-limit},
\[
 (X_\varepsilon,\mathcal M_\varepsilon[X])
 \longrightarrow(h_D,M_D)
\]
in probability in $\mathcal H_D^{-\delta}$ times the local-weak measure space.
For each fixed $\varepsilon$, \eqref{chII:eq:R26.2} and the fact that
$U\mapsto\mathcal M_\varepsilon[U]$ is a deterministic measurable functional of the
smooth field imply
\[
 (Y_\varepsilon,\mathcal M_\varepsilon[Y])
 \overset d=(X_\varepsilon,\mathcal M_\varepsilon[X]).
\]
Hence the left pair converges in law to $(h_D,F_D(h_D))$.  Since
$Y_\varepsilon\to h_D$ in probability in $\mathcal H_D^{-\delta}$, Slutsky's lemma
gives
\[
 (h_D,\mathcal M_\varepsilon[Y])\Longrightarrow(h_D,F_D(h_D)).
\]
Lemma \ref{chII:lem:graph-collapse} now gives \eqref{chII:eq:R26.3} along every null sequence.
If the full limit failed, one could select a null sequence along which one fixed metric
error probability stays bounded away from zero, a contradiction.
\end{proof}

\subsection{A uniform heat-locality estimate on squares}
View $D$ as an embedded square in $\mathbb T^2$, and extend $p_t^D(x,\cdot)$ by zero
outside $D$. This zero-trace extension belongs to $H^1(\mathbb T^2)$.
The notation $K\Subset D$ means that $K$ is a compact subset of the
interior of $D$; distances to $\partial D$ are Euclidean in the chart.

\begin{lemma}[Uniform Sobolev heat locality]\label{chII:lem:sobolev-heat-locality}
For $K\Subset D$ and $0<\delta<1$, there are $c,C,A>0$ such that, for $0<t\le1$,
\begin{equation}\label{chII:eq:sobolev-heat-locality}
 \sup_{x\in K}\|p_t^{\rm per}(x,\cdot)-\one_Dp_t^D(x,\cdot)\|_{H^\delta(\mathbb T^2)}
 \le Ct^{-A}e^{-c/t}.
\end{equation}
The same estimate holds after one derivative in $x$.
\end{lemma}
\begin{proof}
Translate $D$ to $(0,\ell)^2$ and put
$g_t(u)=(4\pi t)^{-1/2}e^{-u^2/(4t)}$.  The one-dimensional image formula is
\[
 p_t^{(0,\ell)}(x,y)=\sum_{m\in\mathbb Z}
 [g_t(x-y+2m\ell)-g_t(x+y+2m\ell)].
\]
The two-dimensional killed kernel is the product of these formulas.  If
$\dist(K,\partial D)>0$, the unreflected principal term cancels in the periodic-minus-killed
difference and every remaining image travels a distance bounded below by a fixed multiple
of $\dist(K,\partial D)$.  Gaussian derivative estimates and summation of image tails give
$Ct^{-A}e^{-c/t}$ bounds in $L^2$ for the difference and its first weak derivative in the
integration variable, uniformly in $x$, also after one $x$ derivative.  Because the killed kernel has zero trace on $\partial D$, its zero
extension belongs to $H^1(\Ttwo)$ with weak gradient equal to the zero
extension of its interior gradient. Interpolation between $L^2$ and
$H^1$ proves \eqref{chII:eq:sobolev-heat-locality}.
\end{proof}

\begin{corollary}[Field and counterterm locality]\label{chII:cor:heat-field-locality}
For the mean-zero periodic field ${h}$ of covariance $4(-\Delta)^{-1}$, put
${h}_t^{\rm per}=P_t^{\rm per}{h}$ and
${h}_t^D(x)=\langle{h},\one_Dp_t^D(x,\cdot)\rangle$.  Then, for $K\Subset D$,
\begin{equation}\label{chII:eq:heat-field-locality}
 \sup_{x\in K}|{h}_t^{\rm per}(x)-{h}_t^D(x)|\to0\quad\text{a.s.},\qquad
 \sup_{x\in K}|\Var{h}_t^{\rm per}(x)-\Var{h}_t^D(x)|\to0.
\end{equation}
\end{corollary}
\begin{proof}
For $0<\delta<1/2$, Fourier summation gives ${h}\in H^{-\delta}(\mathbb T^2)$
almost surely.  Pairing with Lemma \ref{chII:lem:sobolev-heat-locality} controls the field
difference uniformly by $C\|{h}\|_{H^{-\delta}}t^{-A}e^{-c/t}$.
The Gaussian $L^2$ estimate and the separate variance bounds
$O(\log(1/t))$ give the second limit in
\eqref{chII:eq:heat-field-locality}.
\end{proof}

\subsection{Mean-zero conditioning and the periodic heat chaos}
Use $G_D,b_D,v_D,q_D$ from Section~2.3. The Cameron--Martin
inner product of the periodic field is
$\frac14\int_{\mathbb T^2}\nabla f\cdot\nabla g$ on mean-zero functions.
The directions invisible to exterior observations form the closed space
\[
 \left\{f\in H_0^1(D):\int_D f=0\right\},
\]
with zero extension to the torus. Projecting onto this space gives
the conditional residual covariance
\[
 G_D^0(x,y)=G_D(x,y)-\frac{b_D(x)b_D(y)}{v_D},
 \qquad \int_DG_D^0(x,y)\,dx=0.
\]
Indeed $G_D$ inverts $-\Delta_D/4$, and conditioning its Gaussian
pairing with $\mathbf1_D$ removes exactly the displayed rank-one term.
On the original periodic-field space take
$H_D=\EE[h|_D\mid\mathcal F_{D^c}]$ and $h_D^0=h|_D-H_D$.
Adjoin a standard normal $\xi_D$ independent of that entire space
and set $h_D=h_D^0+q_D\xi_D$. The symbols $h_D^0$, $\xi_D$ and $h_D$
will henceforth denote this realization on the extended space. Since
\[
 \operatorname{Cov}(h_D)=G_D^0+q_D\otimes q_D=G_D,
\]
this $h_D$ is an ordinary Dirichlet field. On this extension,
\begin{equation}\label{chII:eq:R26.4}
 h_D=h_D^0+q_D\xi_D,\quad \xi_D\perp h_D^0,\quad \xi_D\sim N(0,1),\qquad
 {h}|_D=h_D^0+H_D,\quad H_D\perp h_D^0.
\end{equation}
In distributional kernel notation,
\[
 \operatorname{Cov}(H_D)=G_{\rm per}|_{D\times D}-G_D
                      +q_D\otimes q_D.
\]
The logarithmic singularities cancel in this difference. Interior
elliptic regularity, or the heat-kernel derivative bounds, gives a
smooth version of $H_D$ on compact interiors; $q_D$ is smooth there
because $-\Delta_Dq_D=4/\sqrt{v_D}$.

\begin{lemma}[Heat approximation of the binding field]\label{chII:lem:binding-heat}
For $K\Subset D$, $P_t^DH_D\to H_D$ uniformly on $K$ almost surely and the corresponding
variances converge uniformly.  Also $P_t^Dq_D\to q_D$ uniformly on $K$.
\end{lemma}
\begin{proof}
The diagonal regular-part identity and domain monotonicity give
$\Var H_D(x)\le C(1+|\log d(x,\partial D)|)$. Consequently
\[
 \EE\|H_D\|_{L^1(D)}^2
 \le\left(\int_D\sqrt{\Var H_D(x)}\,dx\right)^2<\infty.
\]
Choose $\vartheta\in C_c^\infty(D)$ equal to one on a neighborhood
of $K$. The heat approximate identity applied to the smooth local
field $\vartheta H_D$ gives uniform convergence on $K$.
Separated-support heat tails bound the remaining part by
\[
 \sup_K|P_t^D((1-\vartheta)H_D)|
       \le Ct^{-1}e^{-c/t}\|H_D\|_{L^1(D)}.
\]
The same estimates in $L^2(\Omega)$ give
$\sup_{x\in K}\EE|P_t^DH_D(x)-H_D(x)|^2\to0$; bounded
interior variances then give
\[
 \sup_{x\in K}|\Var(P_t^DH_D(x))-\Var(H_D(x))|\longrightarrow0.
\]
Apply the deterministic version of this argument to $q_D$.
\end{proof}

Let $q_{D,\varepsilon}=P_{\varepsilon^2}^Dq_D$ and similarly for $H_D$, and set
\[
 Y_\varepsilon:=P_{\varepsilon^2}^Dh_D.
\]
Thus $\mathcal M_\varepsilon[Y]$ is the exact-variance chaos of this
realized smoothing. The Borel map $F_D$ from Lemma~\ref{chII:lem:critical-limit-measurable}
applies to the present $h_D$; equality in law with the ordinary
Dirichlet construction transfers Proposition~\ref{chII:prop:dirichlet-heat-chaos}
to this realization. Independence and exact variance subtraction then give, samplewise,
\begin{equation}\label{chII:eq:R26.5}
 \begin{aligned}
 M^0_{D,\varepsilon}
 &=e^{-\alpha q_{D,\varepsilon}\xi_D+\alpha^2q_{D,\varepsilon}^2/2}
                 \mathcal M_\varepsilon[Y],\\
 M^{\mathrm{per},D}_\varepsilon
 &=e^{\alpha H_{D,\varepsilon}-\alpha^2\Var(H_{D,\varepsilon})/2}
                 M^0_{D,\varepsilon}.
 \end{aligned}
\end{equation}
Proposition \ref{chII:prop:dirichlet-heat-chaos} and uniform convergence of the smooth factors
therefore give local weak probability limits $M_D^0$ and
\[
 \mathcal Z|_D
 =e^{\alpha H_D-\alpha^2\Var(H_D)/2}M_D^0.
\]
These measurability statements follow from the approximants:
$M^0_{D,\varepsilon}$ is a function of $h_D^0$, and
$M^{\mathrm{per},D}_\varepsilon$ is a function of $h|_D$.
Taking almost surely convergent subsequences for a countable family
of local tests shows that the respective limits have the same
measurability. Hence the periodic limit is $\sigma(h|_D)$-measurable.
By Corollary \ref{chII:cor:heat-field-locality}, killed and periodic heat smoothing have the same
local limit.  A finite deep-chart cover and a partition of unity thus yield
\begin{equation}\label{chII:eq:R26.6}
 \sqrt{\log(1/\varepsilon)}
 e^{\alpha{h}_\varepsilon(x)-\alpha^2V_\varepsilon/2}\,dx
 \longrightarrow \mathcal Z(dx)\quad\text{weakly in probability}.
\end{equation}
Limits agree on overlaps because they are limits of the same periodic regularization.

\begin{lemma}[Positivity and full support of the periodic critical chaos]\label{chII:lem:periodic-chaos-support}
Almost surely every nonempty open $V\subset\mathbb T^2$ satisfies $\mathcal Z(V)>0$;
in particular $0<\mathcal Z(\mathbb T^2)<\infty$ almost surely.
\end{lemma}
\begin{proof}
Choose a Euclidean square $U$ with $\overline U\subset V$ inside a deep chart.  Use the constrained chart decomposition \eqref{chII:eq:R26.4} on $U$.
Multiply the ordinary Dirichlet chaos first by
$e^{-\alpha q_U\xi_U+\alpha^2q_U^2/2}$, and then by the positive binding-field
factor in \eqref{chII:eq:R26.5}.  Theorem 10 of \cite{chII:DRSV} gives a nontrivial positive critical measure for the ordinary Dirichlet field;
Proposition \ref{chII:prop:dirichlet-heat-chaos} places it on the heat regularization of that same field.
The smooth-shift rule \eqref{chII:eq:R26.5} multiplies it by a strictly positive continuous random
density.  Hence the chart measure, and therefore the periodic limit, charges $U$ almost surely.
A countable base gives full support.  Finiteness follows from the finite chart cover.
\end{proof}

\appendixchapter{D}{Proof of Proposition~\ref{chII:prop:relative-contact}}

The root $r\in V_N$ is fixed throughout this appendix, and
$L=\log N$. All vectors and matrix norms use counting measure;
$\|M\|_{\rm F}=(\sum_{i,j}|M_{ij}|^2)^{1/2}$ is the Frobenius norm.
Inequalities between symmetric matrices are in the quadratic-form order:
$M\ge cI$ means $u^{\mathsf T}Mu\ge c\|u\|_2^2$ for every $u$.
Put $n=N^2$. The Green operator is
$G_N=(-\mathcal L_N)^\dagger$, where the dagger denotes the
Moore--Penrose pseudoinverse, acting by counting-measure summation.
Write $v_N=G_N(x,x)=gL+O(1)$.
For nonempty sets $U,V\subseteq V_N$, write
\[
 \begin{gathered}
 d(U,V)=\min_{x\in U,y\in V}d(x,y),\qquad
 d(U,r)=d(U,\{r\}),\\
 \diam U=\operatorname{diam}_d U=\max_{x,y\in U}d(x,y).
 \end{gathered}
\]
Throughout this section, $\mathcal C,\mathcal C'$ are nonempty deterministic sets satisfying
\begin{equation*}
\diam\mathcal C,\diam\mathcal C'\le R_N^{\rm loc},\qquad
d(\mathcal C,r),d(\mathcal C',r),d(\mathcal C,\mathcal C')\ge D_N^{\rm sep}.
\end{equation*}
\paragraph{Standing setup.} Throughout Appendix~D the clusters and
occupation profiles are deterministic and satisfy the stated diameter
and separation conditions; all estimates are uniform over this class.
We use
\begin{equation}
\begin{gathered}
0\le {v_N}-G_N(x,y)\le C\log(2+d(x,y)),\\
|G_N(x,y)|\le C(1+\log(N/d(x,y)))\quad(d(x,y)\ge{D_N^{\rm sep}}/2),\\
|G_N(x,y)-G_N(x',y)|\le C {R_N^{\rm loc}}/{D_N^{\rm sep}},\\
|G_N(x,y)-G_N(x',y)-G_N(x,y')+G_N(x',y')|
 \le C {R_N^{\rm loc}}^2/{D_N^{\rm sep}}^2.
\end{gathered}
\label{chII:eq:R15.1}
\end{equation}
Here the difference bounds hold for $x,x'\in\mathcal C$ and $y,y'\in\mathcal C'$.
The first-difference bound also holds with $y=r$.
The first line follows also from Fourier expansion and
$1-\cos(2\pi k\cdot(x-y)/N)\le\min(C|k|^2d^2/N^2,2)$.
The remaining lines are the Green difference and off-diagonal estimates of Chapter~I.

\subsection{Exact positivity of unconditioned equilibrium weights}
Write $G_{\mathcal C}=G_N[\mathcal C,\mathcal C]$ for the principal submatrix indexed by $\mathcal C$, and
$\mathbf1$ for the all-one vector on the index set of the matrix multiplying it.
Since nonzero eigenvalues of $-\mathcal L_N$ are at most $2$, every vector supported
on the proper subset $\mathcal C$ satisfies
\begin{equation}
u^{\mathsf T}G_Nu
\ge\frac12\left(\|u\|_2^2-\frac{(\sum u)^2}{n}\right).
\label{chII:eq:R15.2}
\end{equation}
Thus $G_{\mathcal C}$ is positive definite. Define
\begin{equation}
e_{\mathcal C}^0=\frac{G_{\mathcal C}^{-1}\mathbf1}
{\mathbf1^{\mathsf T}G_{\mathcal C}^{-1}\mathbf1},\qquad
V_{\mathcal C}^0=(\mathbf1^{\mathsf T}G_{\mathcal C}^{-1}\mathbf1)^{-1}.
\label{chII:eq:R15.3}
\end{equation}
Extend $e_{\mathcal C}^0$ by zero and set $f=G_Ne_{\mathcal C}^0$. Then
\begin{equation}
f|_{\mathcal C}=V_{\mathcal C}^0,\qquad
-\mathcal L_Nf=e_{\mathcal C}^0-\frac1n\mathbf1.
\label{chII:eq:R15.4}
\end{equation}
The positive-charge identity of Chapter~I, with target $\mathcal C$, gives
\begin{equation}
e_{\mathcal C}^0(x)=\frac1n+\frac14\sum_{y\sim x}(V_{\mathcal C}^0-f(y))\ge\frac1n.
\label{chII:eq:R15.5}
\end{equation}
Consequently $\sum e_{\mathcal C}^0=\|e_{\mathcal C}^0\|_1=1$, and \eqref{chII:eq:R15.1} implies
\begin{equation}
{v_N}-C\log(2+\operatorname{diam}\mathcal C)
\le V_{\mathcal C}^0=(e_{\mathcal C}^0)^{\mathsf T}G_{\mathcal C}e_{\mathcal C}^0\le {v_N}.
\label{chII:eq:R15.6}
\end{equation}
For $\phi\sim N(0,G_N)$, write $\phi_{\mathcal C}=(\phi(x))_{x\in\mathcal C}$ and
$S_r=\phi(r)/\sqrt{v_N}$. Define $Z_{\mathcal C}^0=(e_{\mathcal C}^0)^{\mathsf T}\phi_{\mathcal C}$ and
$\eta_{\mathcal C}^0=\phi_{\mathcal C}-Z_{\mathcal C}^0\mathbf1$. Gaussian orthogonality gives
\begin{equation*}
Z_{\mathcal C}^0\perp\eta_{\mathcal C}^0,\qquad \Var Z_{\mathcal C}^0=V_{\mathcal C}^0,\qquad
\Var\eta_{\mathcal C}^0(x)=v_N-V_{\mathcal C}^0.
\end{equation*}

\subsection{Increment spectrum and root conditioning}
Fix $x_0\in\mathcal C$ and let $B_{\mathcal C}\phi_{\mathcal C}=(\phi_x-\phi_{x_0})_{x\in\mathcal C\setminus\{x_0\}}$. Set
\begin{equation}
\Sigma_{\mathcal C}=B_{\mathcal C}G_{\mathcal C}B_{\mathcal C}^{\mathsf T},\quad
q_{\mathcal C}=(G_N(x,r):x\in\mathcal C),\quad
d_{\mathcal C}=B_{\mathcal C}q_{\mathcal C},\quad
c_{\mathcal C}=(e_{\mathcal C}^0)^{\mathsf T}q_{\mathcal C}.
\label{chII:eq:R15.7}
\end{equation}
For every increment coefficient vector $b$,
$\mathbf1^{\mathsf T}B_{\mathcal C}^{\mathsf T}b=0$ and
$\|B_{\mathcal C}^{\mathsf T}b\|_2^2=\|b\|_2^2+(\sum b)^2$. Thus
\begin{equation}
\Sigma_{\mathcal C}\ge\frac12 I.
\label{chII:eq:R15.8}
\end{equation}
Put $\delta_N=R_N^{\rm loc}/D_N^{\rm sep}=e^{-L+L^{1/6}+2L^{1/8}}$,
and $d=|\mathcal C|\le C(R_N^{\rm loc})^2$. By \eqref{chII:eq:R15.1} and \eqref{chII:eq:R15.5},
\begin{equation}
|c_{\mathcal C}|\le CL^{1/8},\qquad
\|d_{\mathcal C}\|_2\le C\sqrt d\,\delta_N.
\label{chII:eq:R15.9}
\end{equation}
Conditional on $\phi(r)$, the centered block covariance is
\begin{equation}
K_{\mathcal C}=G_{\mathcal C}-\frac{q_{\mathcal C}q_{\mathcal C}^{\mathsf T}}{v_N},\qquad
\Sigma_{\mathcal C}^r=\Sigma_{\mathcal C}-\frac{d_{\mathcal C}d_{\mathcal C}^{\mathsf T}}{v_N}.
\label{chII:eq:R15.10}
\end{equation}
The matrix $G_N[\mathcal C\cup\{r\}]$ is positive definite because $|\mathcal C|+1<n$.
Its Schur complement is $K_{\mathcal C}\succ0$. Moreover $\Sigma_{\mathcal C}^r\ge I/4$ for sufficiently large $N$.
All formulas below use the standard empty-matrix convention when the increment space has dimension zero.

Let $\psi_{\mathcal C}=\phi_{\mathcal C}-q_{\mathcal C}\phi(r)/v_N$.
The covariance of $(e_{\mathcal C}^0)^{\mathsf T}\psi_{\mathcal C}$ with $B_{\mathcal C}\psi_{\mathcal C}$ is
$-c_{\mathcal C}d_{\mathcal C}^{\mathsf T}/v_N$. Removing this regression gives
\begin{equation}
e_{\mathcal C}^r=e_{\mathcal C}^0+
\frac{c_{\mathcal C}}{{v_N}}B_{\mathcal C}^{\mathsf T}(\Sigma_{\mathcal C}^r)^{-1}d_{\mathcal C}.
\label{chII:eq:R15.11}
\end{equation}
Its total mass is $1$, and $Z_{\mathcal C}=(e_{\mathcal C}^r)^{\mathsf T}\psi_{\mathcal C}$ is independent of all internal increments.
Gaussian regression, equivalently Sherman--Morrison, yields
\begin{equation}
V_{\mathcal C}:=\operatorname{Var}Z_{\mathcal C}
=V_{\mathcal C}^0-\frac{c_{\mathcal C}^2}
{{v_N}-d_{\mathcal C}^{\mathsf T}\Sigma_{\mathcal C}^{-1}d_{\mathcal C}}.
\label{chII:eq:R15.12}
\end{equation}
The denominator is $v_N-o(1)$ and $L^{-3/4}\to0$; therefore
\begin{equation}
V_{\mathcal C}=gL+O(\log(2+\operatorname{diam}\mathcal C)).
\label{chII:eq:R15.13}
\end{equation}
Uniqueness of the mass-one coordinate orthogonal to increments identifies
\[
 \begin{aligned}
 e_{\mathcal C}^r&=K_{\mathcal C}^{-1}\mathbf1/(\mathbf1^{\mathsf T}K_{\mathcal C}^{-1}\mathbf1),\\
 V_{\mathcal C}&=(\mathbf1^{\mathsf T}K_{\mathcal C}^{-1}\mathbf1)^{-1}.
 \end{aligned}
\]
Since $\|B_{\mathcal C}^{\mathsf T}b\|_1\le2\sqrt d\,\|b\|_2$,
\begin{equation}
\|e_{\mathcal C}^r-e_{\mathcal C}^0\|_1
\le C\frac{L^{1/8}}{L}\,d\,\delta_N
=e^{-L+o(L)},\qquad
\|e_{\mathcal C}^r\|_1\le1+e^{-L+o(L)}.
\label{chII:eq:R15.14}
\end{equation}
We use the $\ell^1$ control in \eqref{chII:eq:R15.14} for the regression weights $e_{\mathcal C}^r$.
For $\eta_{\mathcal C}=\psi_{\mathcal C}-Z_{\mathcal C}\mathbf1$, independence and $K_{\mathcal C}(x,x)\le v_N$ give
\begin{equation}
\max_x\operatorname{Var}\eta_{\mathcal C}(x)
\le C\log(2+\operatorname{diam}\mathcal C).
\label{chII:eq:R15.15}
\end{equation}

\subsection{Conditional two-cluster density}
Write $d'=|\mathcal C'|$, in parallel with $d=|\mathcal C|$.
For two separated clusters, \eqref{chII:eq:R15.1} and \eqref{chII:eq:R15.14}, including the root correction, give
\begin{equation}
\begin{aligned}
|\operatorname{Cov}(Z_{\mathcal C},Z_{\mathcal C'})|&\le CL^{1/8},\\
\|\operatorname{Cov}(B_{\mathcal C}\psi_{\mathcal C},B_{\mathcal C'}\psi_{\mathcal C'})\|_{\rm F}
&\le C\sqrt{dd'}\,\delta_N^2=e^{-2L+o(L)},\\
\|\operatorname{Cov}(Z_{\mathcal C},B_{\mathcal C'}\psi_{\mathcal C'})\|_2
&\le C\sqrt{d'}\,\delta_N=e^{-L+o(L)}.
\end{aligned}
\label{chII:eq:R15.16}
\end{equation}
The covariance with increments of the same cluster is exactly zero.
The joint increment covariance is at least $I/8$, since its off-diagonal block is exponentially small.
Write the joint covariance, with the two common coordinates first,
as $\left(\begin{smallmatrix}V_*&C_*\\C_*^{\mathsf T}&\Sigma_*\end{smallmatrix}\right)$.
Then
\[
 \Sigma_*\ge I/8,\qquad
 \|C_*\|_{\rm F}=e^{-L+o(L)},\qquad
 V_*=gL I_2+O(L^{1/6}+L^{1/8}),
\]
and the conditional covariance is the deterministic Schur complement
\[
 V_*-C_*\Sigma_*^{-1}C_*^{\mathsf T},\qquad
 \|C_*\Sigma_*^{-1}C_*^{\mathsf T}\|\le8\|C_*\|_{\rm F}^2
       =e^{-2L+o(L)}.
\]
Its smallest eigenvalue is therefore at least $cL$.
Thus, uniformly over all increment values and common-coordinate values,
\begin{equation}
\sup_{z,z'} f_{Z_{\mathcal C},Z_{\mathcal C'}\mid\text{all increments}}(z,z')
\le \frac CL.
\label{chII:eq:R15.17}
\end{equation}
Here $f$ denotes conditional density with respect to Lebesgue measure
on $\mathbb R^2$, with both clusters' increment vectors fixed.

After whitening, the off-diagonal block of the centered joint increment covariance has
Frobenius norm $e^{-2L+o(L)}$. Gaussian relative entropy and Pinsker give the same order in total variation.
Fix $M>0$. Given $S_r=s$, $|s|\le M$, the nonconstant increment mean is at most
$C_M\sqrt{d+d'}\,\delta_N/\sqrt L$ in norm. Gaussian translation yields
a common total-variation error bound $\tau_{N,M}$, both for the shifted
one-cluster increment law and for the shifted two-cluster law relative
to the product of its centered marginals. With
$\|\mu-\nu\|_{\rm TV}=\sup_A|\mu(A)-\nu(A)|$, we may take
\begin{equation}
\tau_{N,M}=e^{-L+o(L)}.
\label{chII:eq:R15.18}
\end{equation}
uniformly in all admissible shapes.

\subsection{Relative contact estimates at arbitrarily small activity}
Fix a nonnegative occupation vector $\ell$ and $0<v_{\rm lev}\le1$.
In this argument write $a_{\mathcal C}^G$ for
$a_{\mathcal C}^G(\ell;r,v_{\rm lev})$. Put
$r_x=\sqrt{(v_{\rm lev}-2\ell_x)_+}$, $r_* =\max_{x\in\mathcal C}r_x$, and
\begin{equation}
\begin{gathered}
F_{\mathcal C}(z)=\operatorname{Leb}
\bigcup_{x\in\mathcal C}[-z_x-r_x,-z_x+r_x],\\
a_{\mathcal C}^G=(2\pi g)^{-1/2}\mathbb E F_{\mathcal C}(\eta_{\mathcal C}).
\end{gathered}
\label{chII:eq:R15.19}
\end{equation}
The length of a union is between its longest member and the sum of its member lengths, so $2r_*\le F_{\mathcal C}\le2dr_*$ and $F_{\mathcal C}(z+a\mathbf1)=F_{\mathcal C}(z)$.
Let $\iota_{\mathcal C}:\mathbb R^{\mathcal C\setminus\{x_0\}}\to\mathbb R^{\mathcal C}$ insert zero at $x_0$, and define
\[
D_{\mathcal C}=B_{\mathcal C}\psi_{\mathcal C},\qquad
\mu_{\mathcal C}(s)=B_{\mathcal C}(q_{\mathcal C}-v_N\mathbf1)\,s/\sqrt{v_N},\qquad
\widetilde F_{\mathcal C}(u)=F_{\mathcal C}(\iota_{\mathcal C}u).
\]
Then $F_{\mathcal C}(\eta_{\mathcal C})=\widetilde F_{\mathcal C}(D_{\mathcal C})$ exactly. Define the contact event by
\[
A_{\mathcal C}(v_{\rm lev})=\{\exists x\in\mathcal C:2\ell_x+h^r(x)^2\le v_{\rm lev}\}.
\]
The notation $\ell_x$ denotes a fixed nonnegative occupation profile.
Translation invariance of $F_{\mathcal C}$ also gives
$a_{\mathcal C}^G=(2\pi g)^{-1/2}\mathbb EF_{\mathcal C}(\psi_{\mathcal C})$,
which agrees with the activity of Section~2.4.
The increment spectral bound and Gaussian translation imply
\begin{equation}
\frac{|\mathbb E\widetilde F_{\mathcal C}(D_{\mathcal C}+\mu_{\mathcal C}(s))-\mathbb E\widetilde F_{\mathcal C}(D_{\mathcal C})|}
{\mathbb E\widetilde F_{\mathcal C}(D_{\mathcal C})}
\le C d\,\tau_{N,M}=o(1)
\quad(r_*>0).
\label{chII:eq:R15.20}
\end{equation}
For the conditional integration, put
\[
 \begin{aligned}
 c_{\mathcal C}^r&=(e_{\mathcal C}^r)^{\mathsf T}q_{\mathcal C},\\
 u_{\mathcal C}(s)&=(c_{\mathcal C}^r-v_N)s/\sqrt{v_N},\\
 m_{\mathcal C}(x;s)&=(q_{\mathcal C}(x)-c_{\mathcal C}^r)s/\sqrt{v_N}.
 \end{aligned}
\]
Given $S_r=s$, the same-field decomposition is
\[
 h^r|_{\mathcal C}=Z_{\mathcal C}\mathbf1
       +u_{\mathcal C}(s)\mathbf1+\eta_{\mathcal C}+m_{\mathcal C}(s),
 \qquad Z_{\mathcal C}\perp\eta_{\mathcal C}.
\]
Equations~\eqref{chII:eq:R15.9} and~\eqref{chII:eq:R15.14} imply, uniformly for $|s|\le M$,
\[
 u_{\mathcal C}(s)=-\sqrt{v_N}s+O_M(L^{1/8}/\sqrt L),\qquad
 \|m_{\mathcal C}(s)\|_\infty\le C_M\delta_N/\sqrt L.
\]
If $I_{\mathcal C}(z)=\bigcup_{x\in\mathcal C}[-z_x-r_x,-z_x+r_x]$
and $\varphi_{u,V}$ is the $N(u,V)$ density, Fubini gives the exact formula
\[
 \PP(A_{\mathcal C}(v_{\rm lev})\mid S_r=s)
 =\EE\int_{I_{\mathcal C}(\eta_{\mathcal C}+m_{\mathcal C}(s))}
                   \varphi_{u_{\mathcal C}(s),V_{\mathcal C}}(z)\,dz.
\]
Put $\mathcal C_v:=\{x\in\mathcal C:2\ell_x\le v_{\rm lev}\}$ and
take the interval union over $\mathcal C_v$. On $\{\max|\eta_{\mathcal C}|\le L^{1/3}\}$,
the integrand, uniformly over the entire interval union and $|s|\le M$, is
\begin{equation}
\frac{e^{-s^2/2}}{\sqrt{2\pi gL}}(1+o_M(1)).
\label{chII:eq:R15.21}
\end{equation}
Indeed $V_{\mathcal C}=gL+O(L^{1/6})$ and $|z|\le L^{1/3}+O_M(1)$
give $(z-u_{\mathcal C}(s))^2/(2V_{\mathcal C})=s^2/2+o_M(1)$.
By \eqref{chII:eq:R15.15} and $d\le e^{O(L^{1/6})}$, the residual exceptional probability is
at most $e^{-cL^{1/2}}$. Its relative length loss is at most $de^{-cL^{1/2}}=o(1)$.
The global common-coordinate density bound $C/\sqrt L$ controls its contact contribution as well.
Combining these bounds with \eqref{chII:eq:R15.20} gives
\begin{equation}
\PP(A_{\mathcal C}(v_{\rm lev})\mid S_r=s)
=\frac{e^{-s^2/2}}{\sqrt L}\,a_{\mathcal C}^G(1+o_M(1)),
\qquad |s|\le M.
\label{chII:eq:R15.22}
\end{equation}
uniformly whenever $a_{\mathcal C}^G>0$. If $a_{\mathcal C}^G=0$,
the contact event has probability zero.

For two clusters, condition on all increments and apply \eqref{chII:eq:R15.17} to the two interval unions.
Their length product satisfies
\begin{equation}
F_{\mathcal C}F_{\mathcal C'}\le4dd'r_*r_*',\qquad
\mathbb EF_{\mathcal C}\,\mathbb EF_{\mathcal C'}\ge4r_*r_*'.
\label{chII:eq:R15.23}
\end{equation}
The total-variation contribution satisfies a denominator-free bound
by $Cdd'\tau_{N,M}$ times the product of the two expected interval
lengths. Using \eqref{chII:eq:R15.19} gives the required two-cluster
estimate.
Consequently,
\begin{equation}
\PP(A_{\mathcal C}(v_{\rm lev})\cap A_{\mathcal C'}(v_{\rm lev})
\mid S_r=s)
\le\frac{C_M}{L}a_{\mathcal C}^G a_{\mathcal C'}^G,
\qquad |s|\le M,
\label{chII:eq:R15.24}
\end{equation}
uniformly for large $N$.

\subsection{Amplification with a truncated contact count}
For a fixed $M>0$ and a finite deterministic family $\mathcal J$ of
pairwise separated clusters satisfying the stated geometry, abbreviate
$A_{\mathcal C}=A_{\mathcal C}(v_{\rm lev})$ and define
\begin{equation}
Y_M=\mathbf1_{\{|S_r|\le M\}}
\sum_{\mathcal C\in\mathcal J}\mathbf1_{A_{\mathcal C}}.
\label{chII:eq:R15.25}
\end{equation}
Integrating \eqref{chII:eq:R15.22}--\eqref{chII:eq:R15.24} over a bounded root-mode interval gives
\begin{equation*}
\begin{aligned}
\EE Y_M&\ge\frac{c_M}{\sqrt L}\sum_{\mathcal C\in\mathcal J}a_{\mathcal C}^G,\\
\EE Y_M^2&\le\frac{C_M}{\sqrt L}\sum_{\mathcal C\in\mathcal J}a_{\mathcal C}^G
 +\frac{C_M}{L}\left(\sum_{\mathcal C\in\mathcal J}a_{\mathcal C}^G\right)^2.
\end{aligned}
\end{equation*}
Cauchy--Schwarz yields
\begin{equation}
\PP\left(\bigcup_{\mathcal C\in\mathcal J}A_{\mathcal C}\right)
\ge \PP(Y_M>0)
\ge c'_M\frac{\sum_{\mathcal C\in\mathcal J}a_{\mathcal C}^G}{\sqrt L+\sum_{\mathcal C\in\mathcal J}a_{\mathcal C}^G}.
\label{chII:eq:R15.26}
\end{equation}
This is a Gaussian statement for fixed sets and profiles. It is precisely the amplification used in \eqref{chII:eq:contact-amplification}.

\appendixchapter{E}{Proofs of Propositions~\ref{chII:prop:high-root} and~\ref{chII:prop:reference-contact}}

This appendix records the finite-root calculation behind
Propositions~\ref{chII:prop:high-root} and~\ref{chII:prop:reference-contact}.
The auxiliary disk radius is fixed throughout, and the resulting
coefficient is uniform over the test profiles of
Proposition~\ref{chII:prop:reference-contact}.

\subsection{Radial decomposition and its normalization}\label{chII:sec:radial-construction}

For an integer $R\ge2$, let
$D=D_R=\{-R+1,\ldots,R-1\}^2$, up to a lattice translation.
Its boundary side length is $2R$. Let $y$ be at
distance at least $\delta R$ from its boundary, where $0<\delta<1$ is fixed.
In this appendix $G_A$ is the Green kernel of the planar
total-rate-one walk killed on leaving the finite set $A$; fields on
nested sets are coupled by the domain Markov decomposition. Write
$B(y,s)=\{x\in\mathbb Z^2:|x-y|<s\}$ and $B_s=B(0,s)$;
the radius parameter of $B=B(y,s)$ is denoted by $\operatorname{rad}(B)=s$.
Fix an absolute $C_G\ge1$ such that
\[
 |G_{B(y,s)}(y,y)-g\log s|\le C_G,\qquad s\ge1.
\]
Such a constant exists by the killed Green representation
$G_{B(y,s)}(y,y)=\EE_y a_{\rm lat}(X_{\tau_{B(y,s)}}-y)$,
where $\tau_{B(y,s)}$ is the first exit time of the planar
total-rate-one walk from the disk.
The exit distance lies in $[s,s+1]$, and
$a_{\rm lat}(z)=g\log|z|+O(1)$ for $|z|\ge1$.
Choose once and for all $r\ge4000$, large enough for the
interior disk estimates used below, and use concentric disks
with radii $r\rho^j$, where $\rho=\exp((2C_G+1)/g)>1$.
Then each disk-clock
increment lies in $[1,4C_G+1]$, and its mean logarithmic step is
$g\log\rho=2C_G+1$. Choose $n$ by
$r\rho^n\le\delta R<r\rho^{n+1}$.
The domain Markov decomposition at successive disk boundaries
gives a Gaussian radial coordinate with independent increments;
its clock is the diagonal killed Green function. Given the value
at $y$, its middle part is a Gaussian bridge in this clock.
To specify the actual cap probability, let $\eta_D$ be the ordinary
Dirichlet field and, for a prescribed root height $u$, define
\[
 \eta_D^{y,u}(x)=\eta_D(x)
   -\frac{G_D(x,y)}{G_D(y,y)}\eta_D(y)
   +\frac{G_D(x,y)}{G_D(y,y)}u.
\]
Set
\[
 \mathrm{cap}(D,y,u)=\{\eta_D^{y,u}(x)\le u\text{ for all }x\in D\}.
\]
This field has the Gaussian regression law conditional on
$\eta_D(y)=u$. The shorthand $\mathrm{cap}$ below uses the
prescribed root height of the radial calculation; for the
reference-square statement it is $u=m_{2R}+t$.

For a buffer depth $k$ with $n\ge2k+64$, keep the inner disk of
radius $r\rho^{k+2}$ and the outer disk of radius $r\rho^{n-k}$.
Denote these disks by $A\subset B$ and put
$V=G_B(y,y)-G_A(y,y)>0$.
The Green estimate gives
\begin{equation}\label{chII:eq:radial-clock}
 V=g\log R+O_{r,\delta,k}(1).
\end{equation}
\subsection{The retained fields and the two weights}\label{chII:sec:endpoint-weights}

Define the endpoint factors and the retained data for the middle
bridge comparison as follows. Couple the ordinary Dirichlet fields on nested
domains by their domain Markov decomposition. Write $\eta_A$
for the field on a disk $A$ and $\eta_D$ for the field on the
containing square. For the fixed root $y$, set
\[
 v_D=G_D(y,y),\qquad
 X_A(x)=\eta_A(x)-\frac{G_A(x,y)}{v_D}\eta_D(y).
\]
For $A\subset B\subset D$, the data retained on either side
of the middle bridge are
\begin{equation}\label{chII:eq:retained-endpoint-data}
 \left(X_A,\ X_B(y),\ (\Delta_D^{y,s}(x))_{x\in D\setminus B}\right),
 \qquad
 \Delta_D^{y,s}(x)=s-\eta_D(x)
                       +\frac{G_D(x,y)}{v_D}(\eta_D(y)-s).
\end{equation}
Thus $\Delta_D^{y,s}$ is the gap below the prescribed root
height in the regression version of the conditioned field.
It includes every point of the square outside $B$.

Write $r_j=r\rho^j$, so $A=B(y,r_{k+2})$ and $B=B(y,r_{n-k})$.
Set $l_k=(k+3)^{1/6}$ and $h_k=(k+3)^2$. For a field $v$ on $A$
write
\[
 \omega_v(x)=\alpha a_{\rm lat}(x-y)-(v(x)-v(y)),\qquad
 \mathcal C_k(v)=\{\omega_v(x)\ge0\text{ for }|x-y|<r_{k+1}\}.
\]
The use of $r_{k+1}$ inside a field on $r_{k+2}$ leaves one
full buffer layer. For a finite displacement set $W\subset\mathbb Z^2$
with $y+W\subset A$, let $f$ be a measurable $[0,1]$-valued function of
$(\omega_v(y+z))_{z\in W}$. This is a local, or cylinder, test;
we abbreviate its value as $f(\omega_v)$. Define
\begin{align}
 j_{k,f}(v)&=v(y)f(\omega_v)
           \ind_{\{l_k\le v(y)\le h_k\}}\ind_{\mathcal C_k(v)},
 \label{chII:eq:inner-weight-defined}\\
 o_{R,k,t}&=X_B(y)\ind_{\{l_k\le X_B(y)\le h_k\}}
       \ind_{\{\Delta_D^{y,m_{2R}+t}(x)\ge0\ \forall x\in D\setminus B\}},
 \label{chII:eq:outer-weight-defined}\\
 J_k(f)&=\EE j_{k,f}(\eta_A),\qquad
 O_{R,k}(t)=\EE o_{R,k,t}.\label{chII:eq:weights-defined}
\end{align}
Both integrands are nonnegative and at most $h_k$. The inner
integral uses the unconditioned disk field. Lemma~\ref{chII:lem:endpoint-relative-factor} compares the joint law of $X_A$
and the retained exterior data with the product of the unconditioned
inner field law and that exterior law on bounded endpoint windows.

Define the flat comparison event by
\[
 \begin{aligned}
 \mathrm{flat}={}&\mathcal C_k(X_A)
 \cap\{X_A(y),X_B(y)\in[l_k,h_k]\}\\
 &\cap\{\Delta_D^{y,m_{2R}+t}(x)\ge0:\ x\in D\setminus B\}\\
 &\cap\{X_{B(y,r_j)}(y)\ge0:\ k+2\le j\le n-k\}.
\end{aligned}
\]
The last condition is imposed at the discrete domain clocks.
When the buffer depth needs to be displayed, write this event as
$\mathrm{flat}_k(D,y,t)$.
In the square expectation below, the local test is $f(\omega_{X_A})$.
The centered Gaussian bridge between the two retained
endpoints is independent of all the data
\eqref{chII:eq:retained-endpoint-data}. Indeed its covariance with
each inner field value, each endpoint and each exterior
conditioned gap is zero by the nested Green-kernel covariance
identity; joint Gaussianity then gives independence. Endpoint
integration below uses the joint Gaussian law of the two endpoints.

Endpoint integration supplies nonnegative inner and outer weights
$J_k(f)$ and $O_{R,k}(t)$. Here $f$ is a bounded local test of
the drops from $y$, supported at finitely many lattice sites;
$J_k(1)$ denotes the test-free inner weight. The outer weight
includes all cap constraints in the square outside the retained disk.

The quantitative estimates used in this calculation are
\begin{equation}\label{chII:eq:radial-estimates}
 \begin{gathered}
 J_k(1)\longrightarrow I\in(0,\infty),\\
 |O_{R,k}(t)-t|\le C_\delta(1+\log^2(t+2))+e,
 \qquad t\ge1,\\
 \left|V\,\EE[\ind_{\mathrm{flat}}f]-2J_k(f)O_{R,k}(t)\right|
       \le e\,2J_k(f)O_{R,k}(t),\qquad 0\le f\le1.
 \end{gathered}
\end{equation}
For fixed $t,e$, choose $k$ and then the domain size for the second
estimate. For the third, choose $k$ from $e$ and the finite support of
$f$, and then choose the total depth. The event $\mathrm{flat}$ uses endpoint
cutoffs, the buffered inner cap, the full exterior cap, and
nonnegativity at the sampled middle root coordinates. For the test-free cap,
its symmetric difference from the actual conditioned square cap
has probability at most $e/(n+1)$ for any prescribed $e>0$
after enlarging $k$ and then $n$. The two directions and the
seven joint error estimates are specified below.

The three lines of \eqref{chII:eq:radial-estimates} follow from the
inner limit \eqref{chII:eq:inner-two-scale-input}, the outer estimate
\eqref{chII:eq:completed-outer-mass}, and the relative bridge comparison
\eqref{chII:eq:sampled-bridge-relative}, respectively. The following
subsections prove these estimates and the two directions of the cap
replacement on the original sample.

The numerical factor $2$ originates in the continuous Brownian
bridge calculation. For a continuous bridge with variance clock
length $V$ and positive endpoint gaps $u,v$, reflection gives
\begin{equation}\label{chII:eq:bridge-reflection}
 \PP\{\hbox{the bridge stays below its affine boundary}\}
       =1-e^{-2uv/V}.
\end{equation}
On fixed bounded endpoint windows this is
$2uv/V+O((uv/V)^2)$. For the sampled middle bridge, Section~E.9 proves the matching
upper bound by filling the clock intervals with Brownian bridges. Endpoint integration requires a relative estimate uniform in the
outer weight. The filling argument below supplies the matching upper
bound and specifies which angular constraints can be flattened. 

\begin{lemma}\label{chII:lem:endpoint-relative-factor}
Keep the inner disk $A$ and a root window $[-M,M]$, $M>0$, fixed.
Let $B$ expand inside $D$. For every $\varepsilon>0$, if $F\ge0$ is a bounded measurable
test of $X_A$ supported on $|X_A(y)|\le M$, and $Q\ge0$ is a
bounded measurable test of the exterior data supported on
$|X_B(y)|\le M$, then
\begin{equation}\label{chII:eq:endpoint-relative-factor}
 |\EE[F(X_A)Q]-\EE[F(\eta_A)]\EE Q|
 \le\varepsilon\EE[F(\eta_A)]\EE Q
\end{equation}
whenever $G_B(y,y)\ge b_0(A,M,\varepsilon)$, uniformly over the
terminal square and the admissible exterior constraint.
\end{lemma}
\begin{proof}
Put $a=G_A(y,y)$, $b=G_B(y,y)$ and $k(x)=G_A(x,y)$.
The field
\[
 Y(x)=\eta_A(x)-k(x)\eta_B(y)/b
\]
is independent of the exterior data in
\eqref{chII:eq:retained-endpoint-data}; its covariance is
$G_A-k\otimes k/b$. Direct subtraction gives
$X_A=Y+kX_B(y)/b$. Conditional on an exterior value with
$X_B(y)=z$, this is a Gaussian field of mean $kz/b$.
Conditioning further on its root value $u$, its remaining
coordinates have exactly the same conditional law as
$\eta_A$ given $\eta_A(y)=u$. Only the root density changes.
Its density relative to the latter root law is
\begin{equation}\label{chII:eq:scalar-endpoint-density}
 \left(1-\frac ab\right)^{-1/2}
 \exp\left\{-\frac{(u-az/b)^2}{2a(1-a/b)}+\frac{u^2}{2a}\right\}.
\end{equation}
For fixed $a,M$, this tends uniformly to one on
$|u|,|z|\le M$ as $b\to\infty$. Hence the conditional
integral of $F$ lies between $(1-\varepsilon)\EE F(\eta_A)$
and $(1+\varepsilon)\EE F(\eta_A)$ throughout the support of
$Q$. Multiply by $Q$, integrate, and use its nonnegativity. The resulting
the resulting bounds are multiplicative inequalities between the
corresponding finite measures. Finally
$b=g\log\operatorname{rad}(B)+O(1)$, uniformly in the
ambient square, by the disk Green estimate.
\end{proof}

Taking $F=j_{k,f}$ and $Q=o_{R,k,t}$ proves the relative
factorization of the retained leading moment. In particular,
it holds for every nonnegative bounded measurable local test and
retains all square exterior constraints.

\subsection{The pinned cap and cluster law from the literature}\label{chII:sec:pinned-input}

The pinned-field inputs are Biskup and Louidor
\cite[Theorems~2.3--2.4 and Proposition~5.8]{chII:BL}.
Their Green function is that of the
discrete-time simple random walk. The total-rate-one convention
of this paper gives the same killed Green function: each visit
has mean holding time one. Thus both conventions have
$g=2/\pi$, the same potential kernel $a_{\rm lat}$, and the
same parameter $\alpha=2/\sqrt g$.
Use the pinned Gaussian field $\psi$ defined in Section~2.4.
Its law is denoted by $\nu^0$ in \cite{chII:BL}. Centered Gaussian
symmetry allows us to use the drop convention
$\omega(x)=\alpha a_{\rm lat}(x)-\psi(x)$ instead of their plus sign.

Theorem~2.4 of \cite{chII:BL} gives a constant $c_{\rm pin}\in(0,\infty)$
such that
\begin{equation}\label{chII:eq:imported-pinned-cap}
 \PP\{\omega(x)\ge0:\ |x|\le R\}
 \sim\frac{c_{\rm pin}}{\sqrt{\log R}}.
\end{equation}
Theorem~2.3 and Proposition~5.8 of the same paper give the weak
limit of these conditional laws and its finite-sublevel property:
\begin{equation}\label{chII:eq:conditional-local-limit}
 \begin{gathered}
 \mathcal L(\omega\mid\omega\ge0\text{ on }B_R)
 \Longrightarrow\nu,\qquad
 \nu\{\omega(0)=0,\ \omega\ge0\}=1,\\
 \#\{x:\omega(x)\le T\}<\infty\quad\nu\text{-a.s.}
 \end{gathered}
\end{equation}
Proposition~5.8 supplies the last assertion for each finite $T$.
Together these conclusions give the cluster law used in Section~2.4.

The same limits hold for any fixed equivalent lattice norm: the
corresponding balls are nested between radii $cR$ and $CR$. By \eqref{chII:eq:imported-pinned-cap}, the probabilities of
the two nested cap events have ratio tending to one. The total
variation distance between their conditional laws is therefore
at most the relative probability of their difference, which tends
to zero. The intermediate cap is squeezed in the same way.
In particular, write $\mathcal C_n=\{\omega(x)\ge0:
|x|\le r\rho^{n+1}\}$ for the infinite-field disk cap. Then
\begin{equation}\label{chII:eq:imported-geometric-cap}
 \sqrt{n+1}\,\EE[f(\omega);\mathcal C_n]
 \longrightarrow\frac{c_{\rm pin}}{\sqrt{\log\rho}}\nu(f)
\end{equation}
for every bounded continuous cylinder test $f$, where $\nu(f)=\int f\,d\nu$.

Set $a_0=2/\sqrt{2\pi g\log\rho}$. Then
\eqref{chII:eq:imported-geometric-cap} gives the retained inner-weight
normalization
\begin{equation}\label{chII:eq:inner-two-scale-input}
 J_k(f)\longrightarrow I\nu(f),\qquad
 I=\tfrac12\sqrt{2\pi g}\,c_{\rm pin}.
\end{equation}
The identification of our buffered $J_k$ with this limit is
proved below directly on the infinite pinned sample. We use
the cited cap limit to identify the value of these disk weights.

\subsection{Standard spatial estimates and the radial interface}

Recall $D_M=\{1,\ldots,M-1\}^2$ and its ordinary Dirichlet
field $\eta_M$ from Section~2.4. We use
\cite[Lemmas B.12--B.13]{chII:BL} in the following form:
\begin{align}
 \PP\{\max_{x\in A}\eta_D(x)\ge m_M+z-y\}
 &\le C\sqrt{|A|/M^2}\,z e^{-\alpha(z-y)},
 \quad A\subset D\subset D_M,\ z\ge1,\ y\ge0,
 \label{chII:eq:dirichlet-subset-cited}\\
 \PP\{|\max_{D_M}\eta_M-m_M|>u\}
 &\le C e^{-cu},\qquad u\ge0.
 \label{chII:eq:dirichlet-maximum-cited}
\end{align}
The sources cited there are Bramson--Ding--Zeitouni for the
subset bound and Ding for the two-sided maximum tails.
After replacing the cited centering by $m_M+O(1)$, the bounds take
the displayed form under the present side-length convention.

The formula \eqref{chII:eq:extremal-center} defines $m_s$ also for real
$s>1$. For $s>0$, put $Q_s=\{x\in\mathbb Z^2:|x|_\infty<s\}$,
and use the same notation for a common lattice translate of these
concentric squares. Write $\eta_{Q_s}$ for the ordinary Dirichlet
field on $Q_s$.

Choose $K$ so large that
\eqref{chII:eq:dirichlet-maximum-cited} gives
$\PP\{\max_{Q_s}\eta_{Q_s}>m_s-K\}\ge1/2$.
In $Q_{2s}$ write $\eta_{Q_{2s}}=\eta_{Q_s}+H$ on $Q_s$,
with $H$ independent of $\eta_{Q_s}$. Select a maximizing site
using only the latter field. At this site the independent
centered Gaussian $H$ is nonnegative with probability at least
one half, also in the degenerate case. Hence
\begin{equation}\label{chII:eq:core-highpoint-seed-input}
 \PP\{\max_{Q_s}\eta_{Q_{2s}}>m_s-K\}\ge1/4.
\end{equation}
Converting between radius and side-length conventions changes the
centering by $O(1)$, which is absorbed into $K$. This gives the
ordinary deep-core seed by citation and domain Markov comparison.

For radial caps, \cite[Propositions~5.1--5.2]{chII:BL} gives the
pinned inner asymptotic and two-ended high-root factorization, the
latter with additive error $\varepsilon/n$ for bounded height offsets
and continuous cylinder tests in the stated range. Lemma~4.23 there
applies to outer domains between two dyadic squares. For logarithmic
curves and sampled clocks we use
\cite[Propositions~4.7--4.14 and Lemmas~4.15--4.19]{chII:BL}
under their stated hypotheses. The disk matching and buffered
square-cap comparison for the present domains are proved below.

The construction keeps the complete square exterior and gives
constants uniform over interior root positions. The literature
supplies the pinned law and standard maximum bounds; the annular
controls, rare-event budgets and endpoint calculations supply
the disk and finite-square identifications.

\subsection{A surviving-root bound}

We retain a short one-dimensional bound needed for the local
boundary calculation. All barriers in this subsection are deterministic.
Fix $1\le U<\infty$. Let $S_0=0$ be a centered Gaussian walk with
independent increments of variances in $[1,U]$.
For an increasing sequence $b_j$, $b_0\ge0$,
write
\[
 E_j=\{S_i\le b_i:0\le i\le j\},\quad p_j=\PP(E_j),\quad
 K_j=b_0+2\sqrt U+\sum_{i<j}(b_{i+1}-b_i)p_i.
\]
The budget $K_j$ is increasing. The elementary Gaussian overshoot
bound implies
\begin{equation}\label{chII:eq:moving-gap-budget}
 \EE[(b_j-S_j+2\sqrt U)\ind_{E_j}]\le K_j.
\end{equation}
Indeed, if $G$ is centered normal of variance $\sigma^2\le U$
and $d\ge0$, then
\[
 \EE[(d-G+2\sqrt U)\ind_{\{G\le d\}}]\le d+2\sqrt U.
\]
To verify it, subtract the right side. The result is
$\EE[(G-d-2\sqrt U)\ind_{\{G>d\}}]\le0$:
the conditional overshoot has tail at most
$e^{-x^2/(2\sigma^2)}$ and mean at most
$\sigma\sqrt{\pi/2}<2\sqrt U$.
Apply this conditionally with $d=b_{j+1}-S_j$ on $E_j$,
then telescope the nonnegative barrier increments. This proves
\eqref{chII:eq:moving-gap-budget}.

Split at $h=\lfloor j/2\rfloor$, $j\ge2$.
The independent increment from $h$ to $j$ has density bounded
by $(2\pi(j-h))^{-1/2}$. On $E_j$, split the final gap at a
threshold $a>0$. The part with $b_j-S_j\le a$ has probability
at most $ap_h/\sqrt{2\pi(j-h)}$; the other part has probability
at most $K_j/a$, by \eqref{chII:eq:moving-gap-budget}.
With $a=\sqrt{2\pi(j-h)}/4$,
\[
 p_j\le\tfrac14p_{\lfloor j/2\rfloor}+4K_j/\sqrt j.
\]
Induction, using monotonicity of $K_j$ and $p_0,p_1\le1$,
gives $p_j\le16(K_j+1)/\sqrt{j+1}$.
The recurrence for $K_j$ now yields
\[
 K_{j+1}+1\le(K_j+1)
 \left(1+\frac{16(b_{j+1}-b_j)}{\sqrt{j+1}}\right).
\]
Multiplying these inequalities and using $1+x\le e^x$ proves
\begin{equation}\label{chII:eq:moving-survival-proved}
 \PP(E_m)\le
 \frac{16(b_0+2\sqrt U+1)}{\sqrt{m+1}}
 \exp\left\{16\sum_{j<m}\frac{b_{j+1}-b_j}{\sqrt{j+1}}\right\}.
\end{equation}

Set $\ell_{\log}(j)=1+\log(j+1)$. Then
\[
 16B\sum_{j\ge0}\frac{\ell_{\log}(j+1)^2-\ell_{\log}(j)^2}{\sqrt{j+1}}
 \le CB<\infty.
\]
Finiteness follows from the mean-value bound
$\ell_{\log}(j+1)^2-\ell_{\log}(j)^2\le C(1+\log(j+2))/(j+1)$.
For $b_j=H+B\ell_{\log}(j)^2+dj$, $H,B,d\ge0$, the exponent in
\eqref{chII:eq:moving-survival-proved} is at most
$CB+32d\sqrt m$. In particular a logarithmic-square barrier
has survival probability at most $C_{B,U}(H+B+1)/\sqrt{m+1}$.
After a cut at index $o$, use
\[
 \ell_{\log}(o+j)^2\le2\ell_{\log}(o)^2+2\ell_{\log}(j)^2.
\]
The initial gap increases by $2B\ell_{\log}(o)^2$; there is a
constant $C=C(B,U)$ valid for every cut index $o$.

\subsection{Annular controls on the original Gaussian sample}\label{chII:sec:annular-controls}

We use the following construction for both the infinite pinned field
and the square field conditioned at its root. Put $r_j=r\rho^j$.
In a finite clock we prepend the zero field, so its time-zero
coordinate is zero; disk subscripts below retain their geometric
indices.
For the square, the chain is the empty domain, the disks of
radii $r_0,\ldots,r_{n-1}$, and finally the square itself.
Leaving the last disk one scale inside $r_n$ ensures the
positive lower bound for the terminal clock step. Write
$D_{-1}=\varnothing$, $D_j=B(y,r_j)$ for the disk domains,
and $D_n=D$ for the terminal square in the finite chain.
Set $\eta_{D_{-1}}=0$, $S_{-1}=v_{-1}=0$ and, with fields
extended by zero, put for $j\ge0$
\[
 Z_j=\eta_{D_j}-\eta_{D_{j-1}},\qquad R_j=Z_j(y),\qquad
 S_j=\sum_{i=0}^jR_i=\eta_{D_j}(y),\qquad v_j=G_{D_j}(y,y).
\]
In the square, $S_n=\eta_D(y)$ and $v_n=G_D(y,y)$.
The domain increments $Z_j$ are independent. Define
\[
 p_j(x)=\frac{\Cov(Z_j(x),R_j)}{\Var R_j},\qquad
 \zeta_j(x)=Z_j(x)-p_j(x)R_j.
\]
The whole residual family is independent of the whole root vector:
different increments are orthogonal, and the displayed regression
removes the covariance within each increment. Applying the finite
Gaussian criterion and then a monotone-class argument proves the
same assertion for the infinite family. In the square, regression
on the final root turns the partial root sums into the bridge
$W_j=S_j-(v_j/v_n)S_n$, with clock $v_j=G_{D_j}(y,y)$.
It leaves the residual family unchanged and independent of $W$.
The disk clock steps lie in $[1,U]$; the initial and terminal
steps obey these bounds after increasing the fixed $U=U(r,\delta)$.

For the infinite chain take $y=0$. The pinned field is realized on
this same space as
$\psi(x)=\lim_{j\to\infty}(\eta_{D_j}(x)-\eta_{D_j}(0))$ in $L^2$
at each fixed site; the Green difference estimates give this limit
and identify its covariance with the pinned field of
Section~\ref{chII:sec:pinned-input}. Put
\[
 \omega_\ell(x)=\alpha a_{\rm lat}(x)
   -\bigl(\eta_{D_{\ell+2}}(x)-\eta_{D_{\ell+2}}(0)\bigr),
       \qquad x\in D_{\ell+2},
\]
and, for integers $m\ge\ell+2$, define the buffered event
\[
 F_{\ell,m}
 =\{\omega_\ell(x)\ge0:\ |x|<r_{\ell+1}\}
  \cap\{l_\ell\le S_{\ell+2}\le h_\ell\}
  \cap\{S_j\ge0:\ \ell+2\le j\le m\}.
\]
The corresponding buffered event for the conditioned square is
$\mathrm{flat}_\ell(D,y,t)$ of Section~\ref{chII:sec:endpoint-weights}.

With $\ell_{\log}(j)=1+\log(j+1)$ as above, the spatial estimates
needed in this construction are
\begin{align}
 \PP\{\sup_{|x-y|<r_i}|\zeta_{i+d+2}(x)|
                      >C e^{-cd}u\}&\le C e^{-cu^2},
                  &&u\ge1,\label{chII:eq:completed-far-tail}\\
 \PP\{|\max_{r_j\le|x-y|<r_{j+1}}
            (\zeta_{j+1}+\zeta_{j+2})-m_{r_{j+1}}|>u\}
                  &\le C e^{-cu},&&u\ge0.
                  \label{chII:eq:completed-local-tail}
\end{align}
The annular estimate is uniform in the retained depth. For the
finite chain it is applied through level $n-2$; the terminal two layers
use their separate root bounds.

We give the spatial justification, including the lower tail.
An increment generated outside a disk is harmonic inside it.
The Green representation bounds its variance on a fixed interior
fraction by a constant. Subtracting its regression at the center
and applying the oscillation contraction on successive concentric
disks gives variance at most $Ce^{-2cd}$ on the disk in
\eqref{chII:eq:completed-far-tail}. The same argument for differences
gives a canonical metric bounded by
$Ce^{-cd}(|x-x'|/r_i)^{1/2}$. Dyadic nets in the rescaled disk
have $O(4^q)$ edges at depth $q$, of variance at most
$Ce^{-2cd}2^{-q}$. Gaussian tails, with edge thresholds
$C e^{-cd}(u+\sqrt{q+1})2^{-q/4}$, and summation over $q$
prove \eqref{chII:eq:completed-far-tail} with constants uniform in
the lattice size. The terminal square increment
has the same estimate: its inner disk is a fixed distance from
the square boundary, and the same harmonic argument applies.

For the local upper tail add back the two root components.
Conditional on the residuals, at a residual maximizer their sum
is a centered Gaussian, nonnegative with probability at least
one half. On the annulus the still smaller disk field vanishes;
the restored field is the ordinary field on a disk of radius
comparable to $r_{j+1}$. Embed that disk in a comparable square
and apply \eqref{chII:eq:dirichlet-subset-cited}. The bounded change
in the center gives the upper half of
\eqref{chII:eq:completed-local-tail}.

For the lower half, first obtain a core estimate from
\eqref{chII:eq:core-highpoint-seed-input}. Pack $M^2$ disjoint squares
of radius $2s$ in $Q_R$, with $40sM\le R$. Their Dirichlet
fields are independent. With probability at least
$1-e^{-cM^2}$ one small core has a point above $m_s-K$.
Let $x_*$ be the lexicographically first maximizer of the selected
Dirichlet core, hence measurable with respect to the Dirichlet-field
sigma-field. Conditional on that sigma-field, the independent binding
field $H(x_*)$ is centered Gaussian with variance at most
$g\log(R/s)+C$. Consequently
\[
 \PP\{\max_{Q_R}\eta_{Q_{2R}}\le m_s-K-a\}
 \le e^{-cM^2}+\exp\{-a^2/[2(g\log(R/s)+C)]\}.
\]
Take $M=2^k$, $s$ comparable to $R2^{-k-6}$ and $a=k+1$.
Since $m_R-m_s\le C(k+6)$, this proves an exponential lower
tail at offsets $C_0(k+1)+C_1$ until $s$ reaches a fixed size.
Choose $C_0$ larger than the centering coefficient. The last
such offset then exceeds $m_R+c\log R$; beyond it the single
central coordinate, of variance $O(\log R)$, gives the same
exponential bound. Monotonicity fills the intervening offsets.
Comparable dyadic squares and one further binding comparison
give, for all $R$,
\begin{equation}\label{chII:eq:completed-core-lower}
 \PP\{\max_{Q_R}\eta_{Q_{2R}}<m_R-u\}\le Ce^{-cu}.
\end{equation}
Choose a square core of radius comparable to $r_{j+1}$ inside
the annulus, with its double still inside the next disk.
At its fine-field maximizer the binding value has uniformly
bounded variance. The first removed root component is bounded
by $C|R_{j+1}|$ there; the next residual is independent of
that fine field and has uniformly bounded variance there.
Split the offset into three parts and use
\eqref{chII:eq:completed-core-lower} and the two Gaussian tails.
This proves the lower half of \eqref{chII:eq:completed-local-tail}.

Here is an explicit choice of the control events. The index set is
$\mathcal J=\mathbb N_0$ for the infinite chain and
$\mathcal J=\{0,\ldots,m\}$ for the finite chain, relabeling its terminal
index as $m$. Set $d_j=j$ in the first case and $d_j=\min(j,m-j)$
in the second. Write $Y_j=R_j$ in the first case and
$Y_j=W_j-W_{j-1}$, with $W_{-1}=0$, in the second. Define
\[
 \begin{gathered}
 \mathcal A=\{j\ge0:j+2\in\mathcal J\},\qquad
 \mathcal F=\{(i,d)\in\mathbb N_0^2:i+d+2\in\mathcal J\},\\
 M_j=\left|\max_{r_j\le|x-y|<r_{j+1}}
        (\zeta_{j+1}(x)+\zeta_{j+2}(x))-m_{r_{j+1}}\right|,\\
 V_{i,d}=\sup_{|x-y|<r_i}|\zeta_{i+d+2}(x)|.
 \end{gathered}
\]
For a finite chain, omit $j=m-2$ from $\mathcal A$; the terminal
two annuli use the root bounds as stated above.
Choose fixed constants $C_{\rm ctl},c_{\rm ctl}>0$ from
\eqref{chII:eq:completed-far-tail}--\eqref{chII:eq:completed-local-tail},
enlarging the former if necessary. For $K\ge1$, put
\[
 \begin{split}
 G_K={}&\bigcap_{j\in\mathcal J}
   \{|Y_j|\le C_{\rm ctl}\ell_{\log}(\max(K,d_j))\}\\
 &\cap\bigcap_{j\in\mathcal A}
   \{M_j\le C_{\rm ctl}\ell_{\log}(\max(K,d_j))^2\}\\
 &\cap\bigcap_{(i,d)\in\mathcal F}
   \{V_{i,d}\le C_{\rm ctl}e^{-c_{\rm ctl}d}
                         \ell_{\log}(\max(K,d_i,d))\}.
 \end{split}
\]
These are increasing measurable events, and, for every $p>0$,
\begin{equation}\label{chII:eq:completed-control-tail}
 \PP(G_K^c)\le C_p(K+1)^{-p},\qquad
 \PP(\bigcup_KG_K)=1.
\end{equation}
Indeed all failure probabilities are bounded by
$C\exp\{-c\ell_{\log}(\max(K,i,d))^2\}$, with at most two indices
at each finite depth. Multiplication by any prescribed power
of $(K+1)(i+1)(d+1)$ remains bounded; the union bounds are
therefore summable uniformly in $m$.

Let $B$ be a deterministic constant, depending on the disk geometry
and on the prescribed bounded root-height range in the square case,
large enough to dominate the following logarithmic allowances.
Telescoping the domain increments now gives the spatial
comparison used below. On the $j$th annulus, the maximum of
the centered field plus its root coordinate differs from
$m_{r_{j+1}}$ by at most $CB\ell_{\log}(\max(K,d_j))^2$ on $G_K$.
For the infinite field the root coordinate is $S_j$; for the
conditioned square it is $W_{j+1}$. To check this, keep the
two adjacent residuals in \eqref{chII:eq:completed-local-tail};
the other terms are bounded by \eqref{chII:eq:completed-far-tail}
and the root regression estimate
$|p_i(x)-1|\le Ce^{-c(i-j)}$. The series
$\sum_{a\ge0}e^{-ca}\ell_{\log}(a)$ is finite. The deterministic drift
differs from $m_{r_{j+1}}$ by $O(\ell_{\log}(d_j))$ on the annulus.
In the square this follows from
$G_D(y,y)-G_D(x,y)=a_{\rm lat}(x-y)+O_\delta(|x-y|/R)$
and the regression factor $(m_{2R}+t)/G_D(y,y)$; the error
is uniform for bounded $t$. Thus the actual cap forces
\begin{equation}\label{chII:eq:completed-root-envelope}
 S_j\ge-B\ell_{\log}(\max(K,j))^2\quad\hbox{or}\quad
 W_j\ge-B\ell_{\log}(\max(K,d_j))^2,
\end{equation}
respectively. On $G_K$, each buffered event is also contained in
the corresponding root-envelope event: use its retained inner and
exterior inequalities near the ends, and its nonnegative roots in the
middle. Conversely, on $G_K$ the corresponding positive barrier
implies every spatial cap inequality outside a fixed inner disk.
At the finitely many indices adjacent to a cutoff, use the
controlled increment to move to the nearest constrained root;
its logarithmic cost is absorbed in the same square-logarithm
allowance.

\subsection{Joint rare-event bounds and terminal conditioning}

We record the one-dimensional estimates used with these controls.
They are stated for Gaussian clocks with steps in $[1,U]$.
If a walk is constrained below $H+B\ell_{\log}(j)^2$, its survival is
at most $C_B(1+H)/\sqrt{m+1}$ by
\eqref{chII:eq:moving-survival-proved}. For $a,b\ge0$, let $q_m(a,b)$
be the density at $b$ of $a+S_m$, restricted to
\[
 a+S_j\ge-H-B\ell_{\log}(\min(j,m-j))^2,
 \qquad 0\le j\le m.
\]
Splitting this killed transition into three pieces and applying the
survival bound at its two ends gives
\begin{equation}\label{chII:eq:completed-killed-kernel}
 q_m(a,b)\le C_B(a+H+1)(b+H+1)(m+1)^{-3/2},
 \qquad a,b\ge0.
\end{equation}
A fixed multiple of $B$ is absorbed in $C_B$.
The middle free density is at most $C/\sqrt m$;
the two killed pieces have masses bounded by the survival bound.
For a bridge from zero to zero divide by its free endpoint
density, which is comparable to $m^{-1/2}$. In particular
\begin{equation}\label{chII:eq:completed-curved-bridge}
 \PP\{W_j\ge-H-B\ell_{\log}(d_j)^2\ \forall j\}
          \le C_B(1+H)^2/(m+1).
\end{equation}
Splitting instead at a contact index, integrating its endpoint
height, and using \eqref{chII:eq:completed-killed-kernel} on both
sides gives the following contact bounds. With $H,B,\Delta\ge0$ fixed,
the mass of contact with $\Delta\ell_{\log}(j)^2$ at some $q\le j\le m$
in the one-ended survival event is at most
\[
 \frac{C_{H,B,\Delta}}{\sqrt{m+1}}
 \left(\sum_{j\ge q}\frac{\ell_{\log}(j)^6}{(j+1)^{3/2}}
                     +\frac{\ell_{\log}(m)^6}{\sqrt{m+1}}\right).
\]
For a bridge the analogous bound, at indices $q\le d_j$, is
\begin{equation}\label{chII:eq:completed-contact-budget}
 \frac{C_{H,B,\Delta}}{m+1}
       \sum_{j\ge q}\frac{\ell_{\log}(j)^6}{(j+1)^{3/2}}.
\end{equation}
The endpoint integration has length $O(\ell_{\log}(j)^2)$ and two
linear endpoint factors, giving $\ell_{\log}(j)^6$. In the one-ended
case split at $m/2$: the first half yields the displayed
series, and the second half is bounded by
$C\ell_{\log}(m)^6/(m+1)$. Both normalized errors tend to zero, first
as $m\to\infty$, then $q\to\infty$.
At a specified depth $j$, replacing $\Delta\ell_{\log}(j)^2$ by
$j^{1/6}$ gives normalized error at most
$C(j^{1/6}+\ell_{\log}(j)^2+1)^3/(j+1)^{3/2}=o(1)$.
Conditionally on fixed initial data, the future increments have the
same centered Gaussian-walk law started from the conditioned height;
the preceding estimates therefore apply with that height in the linear
factor.

In the infinite chain let $E_m$ be either $\mathcal C_m$ or
$F_{\ell,m}$; in the square let it be either
$\mathrm{cap}(D,y,m_{2R}+t)$ or $\mathrm{flat}_\ell(D,y,t)$,
with total radial depth $m=n$ and an admissible buffer depth $\ell$.
For the infinite chain put $a_m=\sqrt{m+1}$ and for the finite
conditioned square put $a_m=m+1$. We claim
\begin{equation}\label{chII:eq:completed-joint-tail}
 \PP(E_m\cap G_K^c)\le\frac{C_p}{a_m(K+1)^p},
 \qquad \PP(E_m)\le C/a_m.
\end{equation}
Here and in the finite version all constants are uniform in
the buffer depth, and in bounded ranges of the prescribed
root offset. 

Decompose by the first successful control index
$G_k\setminus G_{k-1}$. On that event use
\eqref{chII:eq:completed-root-envelope}. Spatial failure at index
$k-1$ is independent of the whole root path, so it multiplies
the survival or bridge estimate. For a root failure in the
infinite chain, expose its prefix if the failed increment is
in the first half; its controlled prefix sum is at most
$C(j+1)\ell_{\log}(\max(k,j))$. Conditional future survival contributes
at most $C(j+1)^2\ell_{\log}(k)^2/a_m$. If the failure is in the second
half, retain survival in the first half, which is independent
of that increment. Gaussian tails then make the sum over
failed indices smaller than every power of $k+1$.

For the square bridge, a failure between levels $k-1$ and $k$
can occur only among the first or last $k+3$ increments:
elsewhere the thresholds coincide. Expose these outside
increments and the two middle endpoint values. Their absolute
values are at most $C(k+3)\ell_{\log}(k+3)$ on $G_k$.
The centered middle bridge is independent of the exposed data.
If $k\le m/8$, its length is at least $m/2$; the inherited
allowance is at most
$H+2B\ell_{\log}(k+3)^2+C(k+3)\ell_{\log}(k+3)$.
Applying \eqref{chII:eq:completed-curved-bridge} conditionally gives
\[
 \PP(E_m\cap(G_k\setminus G_{k-1}))
 \le\frac{C(k+4)^4}{m+1}\PP(G_{k-1}^c).
\]
For $k>m/8$, use \eqref{chII:eq:completed-control-tail} directly.
Summing with its exponent increased by $p+6$ proves the first
part of \eqref{chII:eq:completed-joint-tail}. Taking $K=2$ and adding
the controlled curved-survival bound proves the second part.

It follows from \eqref{chII:eq:completed-contact-budget} and
\eqref{chII:eq:completed-joint-tail} that, on either cap event,
roots stay above any fixed logarithmic-square curve away from
the ends, with error $o_q(1)/a_m$. Given the error tolerance, choose $K$ so that the $G_K^c$
contribution is below half the tolerance and then choose $q$ so that the
contact-series tail is below the other half. The endpoint
window $[l_\ell,h_\ell]$ can also be imposed at cost
$o_\ell(1)/a_m$: the lower endpoint follows from the specified
depth estimate above, while the root increment bounds on $G_K$
give an upper height $C\ell \ell_{\log}(\ell)<h_\ell$ for large $\ell$.
For each fixed $\ell$, let $m\to\infty$ in the infinite chain;
for the square the two symmetric depths are treated
together. This establishes both the cap-weighted contact
budgets and the cap-weighted endpoint-window deletion.

\subsection{Identification of the buffered disk weights}

We match the disk weights directly to the cited pinned law.
Use $D_j=B(0,r_j)$, $\omega_\ell$ and $F_{\ell,m}$ from
Section~\ref{chII:sec:annular-controls}. Their gaps use the field on
$D_{\ell+2}$, leaving one layer beyond the inner cap.
The root process has variance $g\log r_j+O(1)$.
Conditional on the retained annular data its future increments
are independent. Reflection and the Gaussian overshoot estimate
used in \eqref{chII:eq:sampled-survival-moment} imply
\begin{equation}\label{chII:eq:completed-inner-relative}
 \left|\sqrt{m+1}\,\EE[f(\omega_\ell);F_{\ell,m}]
                    -a_0J_\ell(f)\right|
 \le\varepsilon a_0J_\ell(f),\qquad
 a_0=\frac{2}{\sqrt{2\pi g\log\rho}},
\end{equation}
for all sufficiently large $\ell$, then sufficiently large
$m$, uniformly over measurable $0\le f\le1$ of the retained
data. To see the relative assertion, for a future starting at
$a\in[l_\ell,h_\ell]$ continuous survival gives the lower
coefficient $\sqrt{2/\pi}\,a$. For the upper bound split off
an initial fraction of the future, fill the remaining Gaussian
intervals, and use the surviving first moment at the split.
The latter is at most $a+C\sqrt U$ by the overshoot identity;
the filling allowance is $O(\sqrt{\log m})$ times a survival
probability of order $m^{-1/2}$ and tends to zero. Let the split
fraction tend to zero. The upper coefficient is at most
$\sqrt{2/\pi}(a+C\sqrt U)$. Since $l_\ell\to\infty$, the
relative excess vanishes uniformly in the endpoint window.
The future clock divided by $m$ tends to $g\log\rho$.
Multiplication by the nonnegative retained weight and integration
proves \eqref{chII:eq:completed-inner-relative}.

The projection and boundary errors can be controlled on the
same rare-event scale. For a fixed window inside $r_{q+1}$,
the far residual sum and the root regression differences give,
on $G_K$,
\begin{equation}\label{chII:eq:completed-projection}
 \max_{x\in B_{r_{q+1}}}|\omega(x)-\omega_\ell(x)|
 \le C \ell_{\log}(\max(K,q))(1+\ell-q)e^{-c(\ell-q)}.
\end{equation}
In a finite square add
$\max_{|x|\le r_{q+1}}|
(m_{2R}+t)(1-G_D(y+x,y)/G_D(y,y))
-\alpha a_{\rm lat}(x)|$, which tends to zero as $R\to\infty$
after $q,\ell$ are fixed. The tail sum includes the terminal square layer.

For the infinite actual cap, the fixed finite-dimensional gap laws
are absolutely continuous, so their boundary slabs vanish with their
width. Indeed the pinned covariance is the Green kernel
killed at zero, hence is nonnegative. Regress all its other
coordinates on the chosen nonzero coordinate. Their joint cap
probability is nonincreasing in that coordinate. Integrating
the resulting one-dimensional normal density gives, for small
$\eta>0$,
\[
 \PP(\mathcal C_m,\ 0\le\omega(x)\le\eta)
       \le C_x\eta\PP(\mathcal C_m).
\]
The same argument works on each finite cap; passage to the
cited conditional limit proves the assertion for $\nu$.
For the projected field, its gap has uniformly bounded density
at each fixed nonzero site. Future survival, conditional on
the retained data and its endpoint window, is at most
$C(\ell+3)^2/\sqrt{m+1}$. Thus
\begin{equation}\label{chII:eq:completed-inner-slab}
 \PP(F_{\ell,m},\ 0\le\omega_\ell(x)\le\eta)
       \le C_x\eta(\ell+3)^2/\sqrt{m+1}.
\end{equation}

Discard bad control and late root contacts by
\eqref{chII:eq:completed-joint-tail}, choose the local window, and
then take $\ell$ large. Outside that window the positive root
barrier enforces all the spatial constraints in either
direction. Inside it use \eqref{chII:eq:completed-projection} with
margin $(\ell+3)^{-4}$ and the two boundary bounds just proved.
The window deletion has already been proved above. The union
of these errors gives
\begin{equation}\label{chII:eq:completed-infinite-comparison}
 \lim_{\ell\to\infty}\limsup_{m\to\infty}
 \sqrt{m+1}\,\PP(\mathcal C_m\mathbin\triangle F_{\ell,m})=0.
\end{equation}
For a bounded uniformly continuous cylinder test, its change
under projection is bounded by its modulus of continuity on
$G_K$; on the complement use the joint bound. It follows from
\eqref{chII:eq:completed-inner-relative}--\eqref{chII:eq:completed-infinite-comparison}
and \eqref{chII:eq:imported-geometric-cap} that
\[
 a_0J_\ell(f)\longrightarrow
          \frac{c_{\rm pin}}{\sqrt{\log\rho}}\nu(f).
\]
For $f=1$ this also proves positivity and boundedness of the
weights. Bounded uniformly continuous cylinder tests determine
weak convergence of every finite-dimensional law. Tightness
then extends the assertion to bounded continuous cylinder tests,
by uniform approximation on compact sets. This proves
\eqref{chII:eq:inner-two-scale-input} with its stated constant $I$.
The measurable-test estimate needed for boundary deletion is
\eqref{chII:eq:completed-inner-relative}.

\subsection{From the sampled middle bridge to the coefficient two}

We prove the matching discrete-clock upper bound and hence the
two-sided relative bridge estimate. Let $0=v_0<\cdots<v_m=T$ satisfy
$1\le v_{j+1}-v_j\le U$, and let $S_j$ be a Gaussian walk
with these increment variances. Such a walk can be filled by
Brownian bridges on its clock intervals: the fillings are
independent of all sampled values, and of each other. For
one interval of length at most $U$, reflection gives
\[
 \PP\{\min\text{(centered filling)}<-q\}\le e^{-2q^2/U}.
\]
Thus, by a union bound, a filling of $j$ intervals stays
above $-q$ with probability at least $1-\varepsilon$ if
$q^2\ge(U/2)\log(j/\varepsilon)$. Independence from the
samples is essential here.

In particular the ordinary sampled survival probability
$p_j(a)=\PP(a+S_i\ge0,\ 0\le i\le j)$ satisfies
\begin{equation}\label{chII:eq:sampled-flat-survival}
 p_j(a)\le C\frac{a+\sqrt{U\log(2j+2)}}{\sqrt j},
 \qquad a\ge0,\quad j\ge1.
\end{equation}
Write $\Phi$ for the standard normal distribution function in this
subsection. Choose a filling success probability at least $1/2$.
On sampled survival and successful filling, the Brownian
path stays above $-a-q$. Its probability is
$2\Phi((a+q)/\sqrt{v_j})-1\le
\sqrt{2/\pi}(a+q)/\sqrt{v_j}$, and $v_j\ge j$.

A sharper first-moment bound is also elementary:
\begin{equation}\label{chII:eq:sampled-survival-moment}
 \EE[(a+S_j)\ind_{\{a+S_i\ge0, i\le j\}}]
                  \le a+C\sqrt U.
\end{equation}
To prove it, stop at the first negative sampled value or at
$j$. The stopped expectation is $a$. On a crossing from a
nonnegative height $d$ by a Gaussian increment of variance
$\sigma^2\le U$, the mean overshoot, conditional on crossing,
is at most $C\sigma$, uniformly over $d\ge0$. One can see
this directly from its tail:
\[
 \PP(-d-\sigma G>x\mid\sigma G<-d)
 =\frac{\Phi(-(d+x)/\sigma)}{\Phi(-d/\sigma)}
 \le e^{-x^2/(2\sigma^2)},\quad x\ge0.
\]
For the inequality, write both Gaussian tails as integrals
from $d/\sigma$ onward and translate the numerator by
$x/\sigma$; the density ratio is at most
$e^{-x^2/(2\sigma^2)}$. Integrating this tail and summing
over the disjoint first-crossing events bounds the total
expected negative overshoot by $C\sqrt U$. The finite
stopping identity proves \eqref{chII:eq:sampled-survival-moment}.

Now condition the walk from height $a>0$ to height $b>0$
over clock length $T$, and let $P_{m}(a,b)$ be the probability
that its sampled heights are nonnegative. Split off
$j=\lfloor m/k\rfloor$ sampled intervals from each end,
where $k>2U$ is fixed. The middle clock length $V$ satisfies
$V\ge T-2Uj$, hence $V/T\ge1-2U/k$.
Fill only the middle intervals, with failure probability
$\varepsilon=1/k$ and $q=O_{U,k}(\sqrt{\log(m+1)})$.
For middle endpoint heights $x,y\ge0$, the sampled killed
transition density is bounded by $(1-1/k)^{-1}$ times
the continuous killed density above $-q$. Write
$p_V(z)=(2\pi V)^{-1/2}e^{-z^2/(2V)}$. Reflection gives
the latter as
\[
 p_V(y-x)-p_V(y+x+2q)
 \le\frac{2(x+q)(y+q)}{\sqrt{2\pi V}\,V}.
\]
The inequality uses $1-e^{-u}\le u$ and $p_V\le(2\pi V)^{-1/2}$.

Integrate against the two end-segment killed densities.
By \eqref{chII:eq:sampled-flat-survival} and
\eqref{chII:eq:sampled-survival-moment}, each endpoint integral
is at most
\[
 a+C\sqrt U+
 Cq\frac{a+\sqrt{U\log(2j+2)}}{\sqrt j}
\]
and the analogous expression with $b$. Divide by the free
endpoint density $p_T(b-a)$. With $k$ fixed, the last term
tends to zero as $m\to\infty$, uniformly for $a,b$ in any
fixed compact interval of positive numbers. We obtain
\begin{equation}\label{chII:eq:sampled-bridge-upper}
 \limsup_{m\to\infty}\sup_{a,b\in[l_k,h_k]}
 \frac{T P_m(a,b)}{2ab}
 \le\frac{(1+C\sqrt U/l_k)^2}
              {(1-1/k)(1-2U/k)^{3/2}}.
\end{equation}
The estimates are uniform in the admissible clock.
The lower bound is continuous survival,
$P_m(a,b)\ge1-e^{-2ab/T}$, which gives normalized lower
limit one at fixed $k$. Since $l_k=(k+3)^{1/6}\to\infty$,
the right side of \eqref{chII:eq:sampled-bridge-upper} tends to
one as $k\to\infty$. Therefore, given $\varepsilon>0$,
choose $k$ sufficiently large, and then $m$ sufficiently
large, to have
\begin{equation}\label{chII:eq:sampled-bridge-relative}
 |T P_m(a,b)-2ab|\le\varepsilon\,2ab,
 \qquad a,b\in[l_k,h_k].
\end{equation}
Rescaling clocks with a different fixed positive lower step
bound gives the same conclusion.

The middle bridge is independent of the retained data.
Multiply \eqref{chII:eq:sampled-bridge-relative} by their
nonnegative cap and test weight, and integrate. This compares
$V\EE[\ind_{\mathrm{flat}}f]$ with
$2\EE[j_{k,f}(X_A)o_{R,k,t}]$ with relative error.
Lemma~\ref{chII:lem:endpoint-relative-factor} compares the latter
expectation with $2J_k(f)O_{R,k}(t)$, again relatively.
Choosing each tolerance smaller than one quarter of the
desired tolerance proves the third line of
\eqref{chII:eq:radial-estimates}.  It also shows exactly where
the increasing lower endpoint cutoff is used.

\subsection{Stopping mass with a dependent exterior barrier}

The stopping argument behind the outer factor is useful in an
abstract finite form. Let $S_0=0,S_1,\ldots,S_m$ be a centered
Gaussian walk with independent increments and variance clock
$j-i\le\Var(S_j-S_i)\le U(j-i)$. The lower bound just fixes
the units of the clock. Let $\beta_j$ be the observed exterior
barrier and define
\[
 E_j=\{\beta_i\le t+S_i:\ 1\le i\le j\},\qquad E_0=\Omega.
\]
For the joint law of $\beta$ and $S$, assume the innovation condition
\begin{equation}\label{chII:eq:stopping-innovation}
 \EE[\ind_{E_j}(S_{j+1}-S_j)]=0,\qquad 0\le j<m.
\end{equation}
Assume $t\ge0$ and $B>0$. Suppose there is a nonnegative envelope $\mathcal E$ with
$\EE \mathcal E^2\le M$ and
\begin{equation}\label{chII:eq:barrier-envelope}
 |\beta_j|\le \mathcal E+B(1+\log(j+1))^2.
\end{equation}
For fixed $U,B,M$, the following stopping bound applies to every
pair satisfying the displayed clock and barrier hypotheses. The next
subsection verifies these hypotheses for the reverse square
exploration.

For a Gaussian walk with bounded clock steps, the filling
argument in \eqref{chII:eq:sampled-flat-survival} gives
$\PP(\min_{i\le j}S_i\ge-a)\le
C(a+\sqrt{U\log(2j+2)})/\sqrt j$ for $a\ge0$.
On $\mathcal E\le(j+1)^{1/4}$, \eqref{chII:eq:barrier-envelope} embeds $E_j$
in such a flat event. Markov's inequality and
$(1+\log x)^2\le81x^{1/4}$, $x\ge1$, give
\begin{equation}\label{chII:eq:quarter-survival}
 \PP(E_j)\le Q(t+2)(j+1)^{-1/4},
\end{equation}
for a finite $Q=Q(U,B,M)$.

Let $J=\max_{0\le j<m}(|S_{j+1}-S_j|
-B(1+\log(j+2))^2)_+$. Gaussian tails and the union bound
show that $\EE J\le C(U,B)$ uniformly in $m$: integrate
\[
 \PP(J>x)\le\sum_{j\ge0}
 2\exp\{-[x+B(1+\log(j+2))^2]^2/(2U)\}.
\]
The sum of the integrals is finite, since its terms decrease
faster than any power of $j$. On first failure at $j+1$,
the preceding successful inequality and the failed inequality
give
\begin{equation}\label{chII:eq:first-failure-overshoot}
 |t+S_{j+1}|\le \mathcal E+J+2B(1+\log(j+2))^2.
\end{equation}
The first-failure sets $E_j\setminus E_{j+1}$ are disjoint.
To sum their logarithmic costs, use
\eqref{chII:eq:quarter-survival} and summation by parts. More
explicitly, let $\tau$ be the first-failure index, set to $m+1$
if no failure occurs by $m$. The tail
integral of $(1+\log(\tau+1))^2$ is bounded using
\[
 \PP(\tau>j)\le\min\{1,Q(t+2)(j+1)^{-1/4}\}.
\]
Split the integral at $\log(j+1)=4\log(Q(t+2)+1)$.
Below it the integrand is bounded by one, while above it
the remaining tail is exponentially decreasing in
$\log(j+1)$. This gives a bound
$C(1+\log^2(t+2))$, also for the truncated failure sum at
any finite $m$. Summing \eqref{chII:eq:first-failure-overshoot}
therefore gives
\[
 \sum_{j<m}\EE[\ind_{E_j\setminus E_{j+1}}|t+S_{j+1}|]
       \le C(1+\log^2(t+2)).
\]
Here $\EE \mathcal E\le1+M$.

Telescoping \eqref{chII:eq:stopping-innovation} yields the exact
finite stopping identity
\[
 \EE[\ind_{E_m}(t+S_m)]
   =t-\sum_{j<m}\EE[\ind_{E_j\setminus E_{j+1}}(t+S_{j+1})].
\]
Thus
\begin{equation}\label{chII:eq:stopping-mass-derived}
 |\EE[\ind_{E_m}S_m]-t|
 \le C(1+\log^2(t+2))+t\PP(E_m).
\end{equation}

We also need to remove the endpoint tails in the original
window. For $l\ge0$, $h,u>0$, on $E_m$,
$S_m\ge-t-\mathcal E-B(1+\log(m+1))^2$. Split the part $S_m<l$
at $\mathcal E\le u$, and use $\EE[\mathcal E\ind_{\mathcal E>u}]\le M/u$.
For $S_m>h$, use
$\EE[S_m\ind_{S_m>h}]\le\EE S_m^2/h\le Um/h$.
Together with \eqref{chII:eq:quarter-survival} and
\eqref{chII:eq:stopping-mass-derived}, this gives
\begin{multline}\label{chII:eq:stopping-window-derived}
 |\EE[S_m\ind_{E_m}\ind_{\{l\le S_m\le h\}}]-t|
 \le C(1+\log^2(t+2))\\
 +[l+2t+B(1+\log(m+1))^2+u]
                  Q(t+2)(m+1)^{-1/4}+M/u+Um/h.
\end{multline}
Choose $l$ of order $m^{1/6}$, $h$ of order $m^2$, and
$u$ of order $m^{1/8}$. At fixed $t$, all terms on the
second line tend to zero.

In the square application the endpoint in the original
definition of $O_{R,k}$ is $X=q_R S_m$, where
$q_R=G_B(y,y)/G_D(y,y)\in(0,1]$. At fixed retained
depth, the Green difference $G_D(y,y)-G_B(y,y)$ stays
bounded while $G_D(y,y)\to\infty$. Hence $q_R\to1$.
The window $l\le X\le h$ is exactly
$l/q_R\le S_m\le h/q_R$. Multiply
\eqref{chII:eq:stopping-window-derived} for this window by
$q_R$; the additional difference between $q_R t$
and $t$ is $t(1-q_R)=o(1)$ at fixed depth and $t$.
For large square size, $q_R\ge1/2$, so the lower
window term is bounded by $2l$. For fixed retained depth first take the square-size limit; then let
$k\to\infty$ with $l=(k+3)^{1/6}$ and $h=(k+3)^2$.

\subsection{The reverse exploration keeps the entire square exterior}

Here we verify the innovation condition used above. Let
$D_0=D\supset D_1\supset\cdots\supset D_L$ be the reverse
disk exploration, with $y\in D_L$, and put
$v_j=G_{D_j}(y,y)$, $V=v_0$. Define
\[
 S_j=\frac{V}{v_j}\eta_{D_j}(y)-\eta_D(y).
\]
Then $S_0=0$ and, for $i\le j$, nested-domain covariance gives
\[
 \operatorname{Cov}(S_i,S_j)=V^2/v_i-V.
\]
Thus $S_j$ is a centered Gaussian walk in the reverse clock
$q_j=V^2/v_j-V$. Its increments are independent. For one
step, with $A=D_{j+1}$ and $B=D_j$,
\begin{equation}\label{chII:eq:reverse-innovation-explicit}
 S_{j+1}-S_j=\frac{V}{G_A(y,y)}
 \left(\eta_A(y)-\frac{G_A(y,y)}{G_B(y,y)}\eta_B(y)\right).
\end{equation}
The expression in parentheses has zero covariance with
$\eta_D(y)$ and with every $\eta_D(x)$ for $x\in D\setminus B$.
For the latter assertion, the nested-domain covariance is
$G_A(y,x)=G_B(y,x)=0$, because $x\notin B$.
It is therefore independent of all conditioned exterior gaps
$\Delta_D^{y,s}(x)$, $x\in D\setminus B$.

Define the actual barrier by
\[
 \beta_j=\max_{x\in D_{j-1}\setminus D_j}
                \{S_j+t-\Delta_D^{y,s}(x)\}.
\]
For an empty annulus set $\beta_j=S_j+t$, which imposes no
constraint. The annular inequality $\beta_j\le t+S_j$ is
exactly nonnegativity of every gap in that annulus. It follows
inductively that
\[
 E_j=\{\Delta_D^{y,s}(x)\ge0:\ x\in D\setminus D_j\}.
\]
By \eqref{chII:eq:reverse-innovation-explicit}, the next increment
is centered and independent of the vector determining $E_j$.
The current survival event is measurable before the next centered
increment, so its product with that increment has expectation zero.
This proves \eqref{chII:eq:stopping-innovation}.

The first reverse step groups the exterior down to the disk at
index $n-3$; subsequent steps remove one geometric annulus.
Forward clock steps are between $1$ and a fixed $U_0$.
For fixed $L$, choose the square size so that
$V\ge2U_0(L+3)$. Then $v_j\ge V/2$, and
\[
 q_{j+1}-q_j=\frac{V^2(v_j-v_{j+1})}{v_jv_{j+1}}
\]
is bounded below by the forward step and above by four times
that step. Including the grouped initial step gives a common
upper bound $12U_0$. Hence the reverse-clock increments lie in the required bounded
range and the stopping estimate applies.

We now verify the random envelope on the actual square sample.
Put $V=G_D(y,y)$ and let $r_{n-3}\le\delta R/8$ be the
first retained reverse disk radius. Denote by $M_1$ the maximum
of $\eta_D$ on $D\setminus B(y,r_{n-3})$. Uniformly over
interior roots,
\begin{equation}\label{chII:eq:completed-exterior-maximum}
 \PP\{|M_1-m_{2R}|>u\}\le C_\delta e^{-c_\delta u},
 \qquad \EE(M_1-m_{2R})^2\le C_\delta.
\end{equation}
For the upper tail use the whole-square maximum estimate
\eqref{chII:eq:dirichlet-subset-cited}. For the lower tail place
a square core of radius $c\delta R$ between distances
$\delta R/8$ and $\delta R/3$ from $y$, with its double
inside $D$. Apply \eqref{chII:eq:completed-core-lower} in that
core. At its Dirichlet maximizer the independent binding field
has bounded variance, since the two domain sizes are comparable
at fixed $\delta$. Splitting $u$ into two equal parts gives
the lower tail. Finally integrate $2u\PP(|M_1-m_{2R}|>u)$.

The root profile on the first exterior layer satisfies
\[
 0\le G_D(x,y)/V\le K_\delta/V,
                 \qquad |x-y|\ge r_{n-3}.
\]
This is the off-diagonal Green bound, because $r_{n-3}/R$
is bounded below at fixed $r,\delta$. For $s=m_{2R}+t$,
the first actual barrier is consequently bounded in absolute
value by
\[
 F=|S_1|+|M_1-m_{2R}|+\frac{K_\delta}{V}(|m_{2R}|+t+|\eta_D(y)|).
\]
To verify this, substitute
$\Delta_D^{y,s}(x)=s-\eta_D(x)
+[G_D(x,y)/V](\eta_D(y)-s)$ into its defining maximum.
Both the upper and lower comparisons use the same maximum
$M_1$; the profile term is bounded pointwise by
$K_\delta(|\eta_D(y)|+|s|)/V$. The reverse clock gives
$\EE S_1^2\le12U_0$, and
$\EE\eta_D(y)^2=V$. Choose $R$ sufficiently large after
fixed $t$ that $V\ge K_\delta^2$ and $V\ge K_\delta t$.
Since $|m_{2R}|/V\le\alpha+1$ eventually, we obtain
\begin{equation}\label{chII:eq:completed-first-envelope}
 \EE F^2\le4\{12U_0+C_\delta+
              [K_\delta(\alpha+1)+1]^2+1\}\le C_\delta.
\end{equation}
The right side is independent of $t$ and of the reverse horizon.

For the remaining annuli use the controls of
\eqref{chII:eq:completed-control-tail} with exponent six. On the
annulus $r_j\le|x-y|<r_{j+1}$, let $\beta^{(j)}$ be its
actual barrier and let $S^{(j)}$ be the reverse root at its
inner disk. The telescoping comparison above, applied before
discarding the root-profile term, gives on $G_K$
\begin{equation}\label{chII:eq:completed-reverse-annular}
 |\beta^{(j)}|\le
 |S^{(j)}-W_{j+1}|+
 B \ell_{\log}(\max(K,\min(j+1,n-j)))^2+a_jt,
\end{equation}
where $a_j=\sup_{r_j\le|x-y|<r_{j+1}}G_D(x,y)/V$ on the annulus in $D$.
Indeed the centered conditional field has annular maximum
$-W_{j+1}+m_{r_{j+1}}+O(B\ell_{\log}(\max(K,d_j))^2)$.
Subtract its deterministic high-root drift and use the
Green estimate for that drift at $t=0$; the remaining height
term is exactly bounded by $a_jt$. Both sides of the local
maximum estimate are needed for the absolute value in
\eqref{chII:eq:completed-reverse-annular}.

For any nonempty finite set $\mathcal I$ of such annular indices define
\[
 \mathcal E_{\mathcal I}=\max_{j\in\mathcal I}\bigl(
 |\beta^{(j)}|-|S^{(j)}-W_{j+1}|-a_jt
                -B \ell_{\log}(\min(j+1,n-j))^2\bigr)_+.
\]
Set $\mathcal E_\varnothing=0$. This is a measurable finite maximum. Since
$\ell_{\log}(\max(K,d))^2\le \ell_{\log}(K)^2+\ell_{\log}(d)^2$, on $G_K$ it is at most
$B\ell_{\log}(K)^2$. Decompose according to the first successful
control index. Equation~\eqref{chII:eq:completed-control-tail}
gives
\begin{equation}\label{chII:eq:completed-interior-envelope}
 \EE \mathcal E_{\mathcal I}^2\le B^2\ell_{\log}(2)^4+
       C B^2\sum_{K\ge2}\frac{\ell_{\log}(K+1)^4}{(K+1)^6}\le CB^2.
\end{equation}
The displayed bound is uniform in $t$ and in finite annular index
sets, while the profile contribution is the explicit term $a_jt$ in
\eqref{chII:eq:completed-reverse-annular}.

For the first $L$ reverse layers, take $\mathcal I$ to be their annular
indices after the first layer. The Green bound gives
$a_j\le K_{\delta,L}/V$. The correction from the forward
bridge to the reverse root is
$(V-v_j)S_j/V$. Thus its maximum $J$ satisfies
\[
 \EE J^2\le
 \frac{(L+1)[U_0(L+3)]^2(12U_0L)}{V^2}.
\]
Here we bounded a maximum of squares by their sum and used
$\EE S_j^2\le12U_0L$. Choose $R\ge R_0(L,t,\delta)$ so that this bound and
$K_{\delta,L}t/V$ are at most one. Forward
annular depth is at most $L+2$; more precisely at reverse
index $i$ it is at most $i+2$, so $\ell_{\log}(i+2)^2\le9\ell_{\log}(i)^2$.
Equations~\eqref{chII:eq:completed-first-envelope} and
\eqref{chII:eq:completed-interior-envelope} therefore give
\begin{equation}\label{chII:eq:completed-exterior-envelope}
 |\beta_i|\le\mathcal E+9B\ell_{\log}(i)^2\quad(1\le i\le L),\qquad
 \mathcal E=F+\mathcal E_{\mathcal I}+J+1,\qquad
 \EE\mathcal E^2\le C_\delta(1+B^2).
\end{equation}
All square exterior constraints are included. The envelope constants
are uniform in $L,t$; the required starting square size may depend
on $L,t$. Thus \eqref{chII:eq:stopping-window-derived} applies.

With the uniform envelope \eqref{chII:eq:completed-exterior-envelope},
take $L=k-2$ so that the last reverse disk is exactly
$B(y,r_{n-k})$ in the definition of $O_{R,k}$.
Use $l=l_k$, $h=h_k$, and $u=(L+1)^{1/8}$ there, including
the endpoint rescaling already proved. The error terms other
than $C_\delta(1+\log^2(t+2))$ vanish as $k\to\infty$,
after taking $R\to\infty$ at each fixed $k,t$. Hence
\begin{equation}\label{chII:eq:completed-outer-mass}
 \limsup_{k\to\infty}\limsup_{R\to\infty}
 |O_{R,k}(t)-t|\le C_\delta(1+\log^2(t+2)).
\end{equation}
Equation~\eqref{chII:eq:completed-outer-mass} is the second estimate
in \eqref{chII:eq:radial-estimates}.

\subsection{Local boundary layers with the cap weight retained}

We first prove a relative boundary estimate for Gaussian rectangles.
Let $(X_i)_{i\in I}$ be a finite centered Gaussian family, and fix
$o\in I$ with variance $v>0$ and
$\operatorname{Cov}(X_o,X_i)\ge0$ for all $i$. For thresholds $a_i$,
put $E=\{X_i\le a_i\text{ for all }i\}$. If $0\le\eta\le a_o$,
then
\begin{equation}\label{chII:eq:relative-rectangle-boundary-proved}
 \PP(E,\ X_o\ge a_o-\eta)
 \le \frac{2\eta}{\sqrt{2\pi v}}\PP(E).
\end{equation}
Indeed, write $X_i=c_iX_o+Y_i$, where $c_i\ge0$ and the family
$Y$ is independent of $X_o$. The function
$q(u)=\PP\{Y_i\le a_i-c_i u\text{ for all }i\ne o\}$ is
nonincreasing. Set $b=a_o-\eta\ge0$, and let $\varphi_v$ denote
the centered Gaussian density of variance $v$. The numerator in
\eqref{chII:eq:relative-rectangle-boundary-proved} is at most
$q(b)\eta/\sqrt{2\pi v}$, whereas
\[
 \PP(E)=\int_{-\infty}^{a_o}\varphi_v(u)q(u)\,du
 \ge q(b)\int_{-\infty}^b\varphi_v(u)\,du\ge q(b)/2.
\]

Apply this to the ordinary field in a finite domain $D$,
conditioned to have value $s$ at $y$. Gaussian regression gives
\[
 \eta_D(x)\mid\{\eta_D(y)=s\}
 \ \stackrel{\mathrm{law}}=
 h(x)+s\frac{G_D(x,y)}{G_D(y,y)},\qquad
 \operatorname{Cov}(h(x),h(z))=G_{D\setminus\{y\}}(x,z).
\]
The covariance identity follows either by the Schur complement
or by stopping the killed walk on first reaching $y$.
The latter representation also shows that every covariance is
nonnegative. In the ordinary normalization the diagonal at any
nonroot interior vertex is at least one, because it includes the
initial visit. With
\[
 d_s(x)=s\left(1-\frac{G_D(x,y)}{G_D(y,y)}\right),\qquad
 \mathcal A=\{h(x)\le d_s(x),\ x\in D\},
\]
write $\mathrm{gap}(x)=d_s(x)-h(x)$ for the actual conditioned gap.
Consequently, for $x\ne y$ and
$0\le\eta\le d_s(x)$,
\begin{equation}\label{chII:eq:actual-relative-boundary-proved}
 \PP(\mathcal A,\ d_s(x)-h(x)\le\eta)
 \le2\eta\PP(\mathcal A).
\end{equation}
Let $\mathcal A$ denote the full square-cap event; every radial
estimate below is applied under this event.

The positivity requirement holds uniformly on any fixed punctured
decoration window. To see this for a square of radius $R$, with
$y$ at distance at least $\delta R$ from its boundary, write
$V=G_D(y,y)$. The Green estimates give
\[
 V=g\log R+O_\delta(1),\qquad
 V-G_D(y+z,y)=a_{\rm lat}(z)+O_\delta(|z|/R),\qquad
 \frac{m_R+t}{V}\longrightarrow\alpha,
\]
uniformly for bounded $t$ and fixed $z$. Here $a_{\rm lat}$ is
the lattice potential kernel in the ordinary normalization.
The second estimate follows by the killed Green representation:
the potential difference at exit is $O(|z|/(\delta R))$, since
all exit points are at distance at least $\delta R$ from $y$.
Thus
$d_{m_R+t}(y+z)\to\alpha a_{\rm lat}(z)>0$ for $z\ne0$.
Choose a positive margin below half the minimum of these finitely
many limits. Uniform convergence then gives the corresponding lower
bound for all allowed roots and heights once $R$ is large. A union bound gives
\begin{equation}\label{chII:eq:actual-finite-boundary-proved}
 \PP(\mathcal A,\ \min_{z\in W}\mathrm{gap}(y+z)\le\eta)
 \le 2|W|\eta\PP(\mathcal A),\qquad 0\notin W.
\end{equation}
The finite-square Gaussian regression estimate together with
\eqref{chII:eq:completed-joint-tail} gives the required
$C_W\eta/(n+1)$ bound.

The projected boundary requires a different estimate because
the retained inner functional carries a surviving-root weight.
Translate the root to zero. In this paragraph let
$D_j=B(0,r\rho^j)$ and couple the fields $\eta_{D_j}$ by
their domain Markov decomposition. Write $S_j=\eta_{D_j}(0)$
for the radial root process. Fix $x\ne0$, let $k=\ell+2$, and put
\[
 Y_k=\eta_{D_k}(x)-\eta_{D_k}(0),\qquad
 \omega_\ell(x)=\alpha a_{\rm lat}(x)-Y_k.
\]
The potential-kernel representation gives
$\operatorname{Var}Y_k\to2a_{\rm lat}(x)>0$ as
$k\to\infty$.
Thus $Y_k$ has a Gaussian density bounded by a fixed constant
$C_x$, uniformly for all sufficiently large $k$.
In particular,
\begin{equation}\label{chII:eq:projection-interval-proved}
 \PP\{0\le \omega_\ell(x)\le\eta\}\le C_x\eta.
\end{equation}
We insert this bound directly into the joint survival probability,
keeping its unconditioned slab factor.

Let $F_{\ell,n}$ be the buffered one-ended flat event. Its
retained conditions include $\omega_\ell(x)\ge0$ and
$0\le S_k\le(k+1)^2$; its future roots remain nonnegative.
Let $\mathcal G_k$ be generated by the roots through $k$ and all
residuals. The projected slab is $\mathcal G_k$-measurable, and the
future root increments are independent of $\mathcal G_k$. The flat instance of
\eqref{chII:eq:moving-survival-proved} bounds their conditional
survival by
$C\{(k+1)^2+1\}/\sqrt{n-k+1}$.
For $n\ge2k$, integrating and using
\eqref{chII:eq:projection-interval-proved} proves
\begin{equation}\label{chII:eq:buffered-slab-survival-proved}
 \PP(F_{\ell,n},\ \omega_\ell(x)\le\eta)
 \le \frac{C_x\eta\{(\ell+3)^2+1\}}{\sqrt{n+1}}.
\end{equation}

Denote by $J_\ell(f)$ the corresponding buffered leading
functional. Apply \eqref{chII:eq:completed-inner-relative} with
relative error $1/2$ and the measurable slab test
$f=\mathbf1_{\{\omega_\ell(x)\le\eta\}}$.
There is $\ell_0$ such that for $\ell\ge\ell_0$ there is $n_0(\ell)$ with, for $n\ge n_0(\ell)$, its lower inequality and
\eqref{chII:eq:buffered-slab-survival-proved} give
\begin{equation}\label{chII:eq:buffered-leading-slab-proved}
 J_\ell(\mathbf1_{\{\omega_\ell(x)\le\eta\}})
 \le C_x\eta\{(\ell+3)^2+1\}.
\end{equation}
The tail and total depths have ratio tending to one, so replacing
one by the other changes the bound by a factor $1+o(1)$; fixed
geometric constants are absorbed into $C_x$.

Put $\mathcal B=\mathrm{flat}_\ell(D,y,t)$,
$A=B(y,r_{\ell+2})$, and abbreviate $O_{R,\ell}=O_{R,\ell}(t)$,
with $t$ fixed. The corresponding projected gap in this square event
is $\omega_{X_A}(y+x)$.
The relative factorization is used as an inequality of nonnegative measures with its exterior weight retained.
The relative factorization gives, for an arbitrarily small fixed
$e\in(0,1/2)$ and sufficiently large $n$ after $\ell$,
\[
 \EE[\mathbf1_{\mathcal B}f]\le (1+e)\frac{2}{V}J_\ell(f)O_{R,\ell},
 \qquad
 \PP(\mathcal B)\ge (1-e)\frac{2}{V}J_\ell(1)O_{R,\ell}.
\]
Combining the two finite-measure inequalities and the positive
inner normalization gives the stated conditional bound without a
separate division by the exterior weight.
The positive inner limit gives $J_\ell(1)\ge c>0$ for all large
$\ell$. Hence the fraction of the flat mass in the slab is at
most $C_x\eta\{(\ell+3)^2+1\}$. With
$\eta_\ell=(\ell+3)^{-4}$ this tends to zero. Combining with
the flat mass upper bound $C/(n+1)$ and summing over the fixed
punctured window yields
\begin{equation}\label{chII:eq:comparison-finite-boundary-proved}
 \PP(\mathcal B,\ \min_{x\in W}\omega_{X_A}(y+x)\le\eta_\ell)
 \le \frac{o_\ell(1)}{n+1}.
\end{equation}
For every fixed finite window there is $\ell_0$ such that, for
$\ell\ge\ell_0$ and $n\ge n_0(\ell)$, the displayed estimate holds. The actual-cap version follows from
\eqref{chII:eq:actual-finite-boundary-proved} with the same margins.
Both mass bounds are supplied by \eqref{chII:eq:completed-joint-tail}.

\subsection{The two directions of the cap replacement}

We specify the error decomposition used when replacing the
actual square cap $\mathcal A=\mathrm{cap}(D,y,m_{2R}+t)$ by the
buffered flat event $\mathcal B=\mathrm{flat}_\ell(D,y,t)$.
These events are defined in Sections~\ref{chII:sec:radial-construction}
and~\ref{chII:sec:endpoint-weights}. We prove that $\PP(\mathcal A\triangle\mathcal B)=o((n+1)^{-1})$
by bounding each difference-event contribution on $\mathcal A$ or
$\mathcal B$.
Fix a bounded height range $|t|\le T$ and a tolerance $e>0$.
Return to the forward bridge $W_j$ of
Section~\ref{chII:sec:annular-controls}, with total index $m=n$;
thus $W_j=X_{B(y,r_j)}(y)$ for $j<n$, $W_n=0$, and
$d_j=\min(j,n-j)$. Fix an integer $K\ge1$.
Let $b_*=b_*(T)$ be the annular-comparison threshold from
Section~\ref{chII:sec:annular-controls}, chosen uniformly for
$|t|\le T$.
For integers $q\ge\max(K,3)$, $\ell\ge q+2$ and $n\ge2\ell+64$, set
\[
 \mathcal L=\bigcup_{q\le j\le n-q}
       \{W_j\le b_*\ell_{\log}(d_j)^2\},\qquad
 \mathcal W=\{W_{\ell+2},W_{n-\ell}\in[l_\ell,h_\ell]\}.
\]
Here the contact inequality is with the positive comparison barrier.
Let $W_q=B(0,r_{q+1})\setminus\{0\}$ and
$\epsilon_\ell=(\ell+3)^{-4}$. With $A=B(y,r_{\ell+2})$, define
\[
 \begin{aligned}
 B_A&=\{\min_{z\in W_q}\Delta_D^{y,m_{2R}+t}(y+z)\le\epsilon_\ell\},\\
 B_B&=\{\min_{z\in W_q}\omega_{X_A}(y+z)\le\epsilon_\ell\}.
 \end{aligned}
\]
These gaps are nonnegative on $\mathcal A$ and $\mathcal B$,
respectively, so the displayed events select their boundary layers.

On $G_K$, for $\ell\ge\ell_0$ and $n\ge n_0(\ell)$ the
buffered local projection differs from the actual gap field by at most
$\epsilon_\ell$. Away from $\mathcal L$, the positive root barrier
gives both the remote projected inequalities and the remote actual cap
inequalities.
The relevant annuli have $q+1\le j<n-\ell$; their roots $W_{j+1}$
and all retained middle roots lie in the range defining $\mathcal L$.
Combining the local projection bound with the common exterior
inequalities gives the difference-event inclusions
\begin{align*}
 \mathcal A\setminus\mathcal B
 &\subset(\mathcal A\cap G_K^c)\cup(\mathcal A\cap\mathcal L)
       \cup(\mathcal A\cap\mathcal W^c)\cup(\mathcal A\cap B_A),\\
 \mathcal B\setminus\mathcal A
 &\subset(\mathcal B\cap G_K^c)\cup(\mathcal B\cap\mathcal L)
       \cup(\mathcal B\cap B_B).
\end{align*}
In the definitions, $\mathcal B$ includes $\mathcal W$ and the
finite boundary events are indexed away from the root, whose gap is
identically zero.

The estimates established above apply to these inclusions:
\eqref{chII:eq:completed-joint-tail} gives
joint control tails $C/((K+1)(n+1))$ on either event;
\eqref{chII:eq:completed-contact-budget} and the endpoint argument
give contact and window errors arbitrarily small relative to
$(n+1)^{-1}$. Equation~\eqref{chII:eq:completed-projection}, with
its finite-square drift term, gives the local projection bound
on $G_K$. The preceding subsection proves the two finite-boundary
budgets on their respective events.

Given $e,T$, choose successively the control depth, local disk,
retained depth and total depth so that each term in the displayed
symmetric-difference inclusions is at most $e/[8(n+1)]$, uniformly over
the prescribed height range. The seven terms in the inclusions sum to at most
$7e/[8(n+1)]$. Therefore
\begin{equation}\label{chII:eq:flat-symmetric-budget}
 \PP(\mathcal A\mathbin\triangle\mathcal B)\le\frac{e}{n+1}.
\end{equation}
 Equation~\eqref{chII:eq:flat-symmetric-budget} also
bounds the error of every common measurable test $0\le f\le1$.
For a continuous local gap test use
\eqref{chII:eq:completed-projection} on $G_K$ and the joint bound
on $G_K^c$ to replace its actual gaps by projected gaps.
Tightness of the normalized inner weights, already proved
above, permits restriction to a compact set of gap vectors;
uniform continuity there removes this last test error.

\subsection{Assembly of the high-root estimate}

Fix $t\ge1$ and choose $k$ so that $|J_k(1)-I|\le1/(t+1)$, the
outer estimate holds with $e=1$, and the relative bridge error is at
most $\eta=(t+C_\delta(1+\log^2(t+2))+2)^{-1}$. With this $k$ fixed,
let $R\to\infty$. Since $J_k(1)\le I+1$ and
$O_{R,k}(t)\le t+C_\delta(1+\log^2(t+2))+1$, multiplication by $V$ bounds the resulting bridge error by an
absolute $C_\delta$ uniformly for $t\ge1$.

Furthermore
\[
 |J_k(1)O_{R,k}(t)-It|
 \le |J_k(1)-I|t+J_k(1)|O_{R,k}(t)-t|
 \le 1+(I+1)(C_\delta(1+\log^2(t+2))+1).
\]
Equation~\eqref{chII:eq:radial-clock} replaces $\log R/V$ by
$1/g$ with the stated fixed-$(k,t)$ error, while the cap replacement
contributes at most $2\log R/(n+1)$. It follows that
\[
 \left|\log R\,\PP(\mathrm{cap})-\frac{2I}{g}t\right|
 \le C_\delta(1+\log^2(t+2)).
\]
Hence $b=2I/g>0$. With $M=2R$, the square and its prescribed root
height $m_{2R}+t=m_M+t$ are unchanged. Since $\log M-\log R=\log2$ and the cap probability is
$O(1/\log R)$, the prefactor change contributes $o_R(1)$. Hence
\eqref{chII:eq:high-root} has coefficient $b=2I/g$. The ordinary
Dirichlet field on the reference square and its domain residual in
a larger torus have the same covariance, the inverse of the killed
generator on that square. This identifies the reference-square law.

\subsection{Decoration and the density at the root}\label{chII:sec:root-density}

The local law follows from the retained weights with their cap
factor kept. For a bounded continuous cylinder test $f$, the relative
factorization, cap replacement and $J_k(f)-\nu(f)J_k(1)\to0$ give
\[
 \log M\left|\EE[f(\Delta^{y,m_M+a});\mathrm{cap}(D_M,y,m_M+a)]
       -\nu(f)\PP(\mathrm{cap}(D_M,y,m_M+a))\right|\longrightarrow0.
\]
Here $\Delta^{y,s}$ is the gap in the regression version at height
$s$; the test uses fixed displacements from $y$. The convergence is uniform over roots in a fixed core and over
$a$ in a fixed compact interval, by the compact-uniform control depths
in \eqref{chII:eq:flat-symmetric-budget} and the exterior bound
\eqref{chII:eq:completed-outer-mass}.
The scalar cap estimate is also uniform on compact height intervals:
$\mathrm{cap}(D_M,y,s)$ increases with $s$, because
$0\le G_{D_M}(x,y)/G_{D_M}(y,y)\le1$, so finite height meshes extend
\eqref{chII:eq:high-root} from fixed heights to each compact interval.

A profile changes the cap on the whole square. Put
\[
 a_y(x)=1-\frac{G_{D_M}(x,y)}{G_{D_M}(y,y)},\qquad
 \Delta^{y,s,d}(x)=\Delta^{y,s}(x)+d(y/M)-d(x/M),
\]
and let $\mathcal A_y^d(s)=\{\Delta^{y,s,d}(x)\ge0:\ x\in D_M\}$.
For $R=M/2$, the Green estimates imply, uniformly for $y/M\in Q_\delta$
and all sufficiently large $M$,
\[
 a_y(x)\ge
 \begin{cases}
 (8g\log R)^{-1},&0<|x-y|\le\sqrt R,\\
 1/8,&|x-y|>\sqrt R.
 \end{cases}
\]
For the first bound use
$G_{D_M}(y,y)-G_{D_M}(x,y)=a_{\rm lat}(x-y)+O_\delta(|x-y|/R)$,
$a_{\rm lat}(z)\ge1/2$ for $z\ne0$, and $G_{D_M}(y,y)\le2g\log R$.
Between $\sqrt R$ and a fixed fraction of $R$, the same formula gives
a deficit at least $(g/4)\log R$; farther out the off-diagonal Green
bound is $O_\delta(1)$. These yield the second line.
Since $(\sqrt R/(2R))^\rho\log R\to0$, every profile with
$\|d\|_\infty+[d]_\rho\le B$ consequently satisfies
\[
 |d(x/M)-d(y/M)|\le16B\,a_y(x)\qquad(x\in D_M).
\]
For $M\ge M_0(\delta,\rho)$ the displayed bound holds uniformly
over the allowed $d,y,s$, with coefficient $16B$. The identity
$\Delta^{y,s+v}=\Delta^{y,s}+v a_y$ now gives, on the same sample,
\[
 \begin{gathered}
 \mathcal A_y^0(s-16B)\subset\mathcal A_y^d(s)
                    \subset\mathcal A_y^0(s+16B),\\
 \mathcal A_y^d(s)\mathbin\triangle\mathcal A_y^0(s)
 \subset\mathcal A_y^0(s+16B)\setminus\mathcal A_y^0(s-16B).
 \end{gathered}
\]
For $s=m_M+t+a$, $|a|\le H$, and $t\ge H+16B+1$, the compact-height
cap bounds therefore imply
\[
 \begin{aligned}
 |\log M\,\PP(\mathcal A_y^d(s))-bt|
     &\le C_{\delta,B,H}(1+\log^2(t+2)),\\
 \log M\,\PP(\mathcal A_y^d(s)\mathbin\triangle\mathcal A_y^0(s))
     &\le C_{\delta,B,H}(1+\log^2(t+2)).
 \end{aligned}
\]
On a fixed displacement window the profile changes the gap by at most
\[
 B\left(\frac{\max_{w\in W}|w|}{M}\right)^\rho.
\]
Combining its vanishing modulus
with the preceding cap-weighted local limit proves the corresponding
profiled local estimate after division by $t/\log M$: choose the
lower height first, then take the domain limit uniformly on each fixed
compact height interval. This is the
whole-square comparison required for $p_M$.

The density part of the calculation is explicit. If
$z=y/M\in Q_\delta$ and $u$ stays in a fixed compact interval,
then \eqref{chII:eq:robin-comparison} and Gaussian density expansion give
\begin{equation}\label{chII:eq:root-density}
 \begin{split}
 &\frac{M^2}{\log M}\,
 \frac{\exp\{-(m_M+t+u-d(z))^2/(2G_{D_M}(y,y))\}}
      {\sqrt{2\pi G_{D_M}(y,y)}}\\
 &\qquad\longrightarrow
 \frac{e^{\alpha^2(c_{\rm lat}+c_{\rm heat})/2}}{\sqrt{2\pi g}}
       w(z)e^{-\alpha(t+u-d(z))}.
 \end{split}
\end{equation}
The convergence is uniform on a fixed core and compact intervals
of $t,u,d(z)$. To check the logarithmic prefactor, the leading
exponent is $-2\log M+\tfrac32\log\log M$; the normal density
denominator removes one half-power of $\log M$. Thus the root
density has size $M^{-2}\log M$, which cancels the $1/\log M$
in the cap probability.

Let $W_{M,y}=\{w\in W:y+w\in D_{M,\delta}\}$, let
$\varphi_{G_{D_M}(y,y)}$ be the centered normal density of variance
$G_{D_M}(y,y)$, and set $s_y(u)=m_M+t+u-d(y/M)$.
Uniqueness of the maximizer gives
\[
 \begin{aligned}
 p_M(\delta,t,T,W,d)
 &=\sum_{y\in D_{M,\delta}}\int_{-1/2}^T
   \varphi_{G_{D_M}(y,y)}(s_y(u))\\
 &\qquad\cdot\PP\bigl(\mathcal A_y^d(s_y(u)),\\
 &\hspace{24mm}\exists w\in W_{M,y}:
      |u-\Delta^{y,s_y(u),d}(y+w)|\le1/2\bigr)\,du.
 \end{aligned}
\]
The profiled cap and gap estimates replace the probability by
$(bt/\log M)\nu\{\exists w\in W:|u-\omega_w|\le1/2\}$,
after the boundary replacements below. Their relative error vanishes
as the lower height tends to infinity, since $u,d$ are bounded and
$(1+\log^2(t+2))/t\to0$. Multiplication by
\eqref{chII:eq:root-density}, followed by summation over $y$, leaves
\[
 c_*t e^{-\alpha t}
       J_{W,T}\,M^{-2}\sum_{y/M\in Q_\delta}
                              e^{\alpha d(y/M)}w(y/M).
\]
For profiles bounded in $\|d\|_\infty+[d]_\rho$ as in Section~2.4,
the last Riemann sum converges
uniformly to its integral. Indeed the oscillation of $d$ on a
mesh square is bounded by $B(\sqrt2/M)^\rho$; on a fixed core
$w$ is uniformly continuous.

There are two boundary replacements. First restrict roots to
$Q_{\delta+\tau}$, so $(y+W)/M\subset Q_\delta$
once $M$ is large. The upper bound uses $Q_\delta$ itself. The
deterministic difference of position integrals is at most
$e^{\alpha B}\|w\|_\infty$ times the intervening strip area.
Second use continuous upper and lower ramps with transition width
$\eta$ around $|u-\omega_w|=1/2$. Their difference is supported in
$\bigcup_{w\in W}\{u:||u-\omega_w|-1/2|\le\eta\}$, whose Lebesgue
measure is at most $4|W|\eta$. Consequently its limiting weighted
integral is at most $4|W|e^{\alpha/2}\eta$, uniformly in the gap vector.
For fixed $\tau,\eta>0$, the compact family of ramp tests permits
the domain limit uniformly in $u$. The limiting error is the sum of an
$O(\eta)$ term and the collar contribution; sending $\eta\downarrow0$
and then $\tau\downarrow0$ proves \eqref{chII:eq:reference-contact}.

\appendixchapter{F}{Proof of Proposition~\ref{chII:prop:critical-laplace}}

We first prove the barrier estimates in \eqref{chII:eq:barrier-package}
and then insert them into the Laplace-transform argument.

\subsection{The cutoff, its first moment, and its crossing event}\label{chII:sec:barrier-construction}

Use the independent-increment Dirichlet realization in
\eqref{chII:eq:dirichlet-three-fields}. Fix an upper heat time $A>0$
and a compact core $K\Subset Q$. All spatial tests in this appendix
are nonnegative continuous functions supported in $K$.
For $0<a\le A$, write $X_a(x)=2\int_a^\infty P_{s/2}^QW(ds)(x)$ and use the
variance clock for the increments between $A$ and $a$.
Here $a$ is heat time, corresponding to $a=2\varepsilon^2$ in
the earlier notation $X_\varepsilon$.
The barrier at height $b$ is imposed on the compensated field
$X_a(x)-\alpha\EE X_a(x)^2$ at every finer scale; the initial
coarse field is included. Write $v_a(x)=\EE X_a(x)^2$ and
$R_{b,a}(x)=b+\alpha v_a(x)-X_a(x)$. The derivative barrier density is
\[
 D_{b,a}(x)=e^{\alpha X_a(x)-\alpha^2v_a(x)/2}R_{b,a}(x)
       \ind_{\{R_{b,u}(x)\ge0\text{ for every }a\le u\le A\}}.
\]
Let $\mathcal F_a=\sigma(W|_{[a,\infty)\times Q})$, completed.
For fixed $x$, $X_a(x)-X_A(x)$ is a Brownian motion, independent
of $\mathcal F_A$, at clock $v_a(x)-v_A(x)$ as $a$ decreases.
Exponential tilting removes the drift of the compensated gap;
the killed Brownian identity $\EE_r[R_t;\tau_0>t]=r$ gives
\[
 \EE[D_{b,a}(x)\mid\mathcal F_{a'}]=D_{b,a'}(x),
 \qquad 0<a\le a'\le A.
\]
Thus the gap-weighted density is a nonnegative martingale.
Its spatial integrals are
bounded in $L^2$ on compact interiors by
\eqref{chII:eq:barrier-pair}, so they have a local mean-square limit
$M_b$, for each fixed $b>0$; explicitly,
$M_b(f)=\lim_{a\downarrow0}\int f(x)D_{b,a}(x)\,dx$ in $L^2$.
Write
\[
 E_b=\left\{\sup_{\substack{x\in K\\0<a\le A}}
       (X_a(x)-\alpha v_a(x))>b\right\}
\]
for the all-time crossing event on $K$.
To identify its normalization, distinguish the derivative density
from the truncated Wick density. For $a=2\varepsilon^2\le A$, put
\[
 W_{b,\varepsilon}(f)=\sqrt{\log(1/\varepsilon)}
 \int_Q f(x)e^{\alpha X_a(x)-\alpha^2v_a(x)/2}
     \ind_{\{R_{b,u}(x)\ge0\ \forall u\in[a,A]\}}\,dx.
\]
The fixed-barrier Brownian survival calculation gives
$W_{b,\varepsilon}(f)\to M_b(f)$ in probability. Its constant uses
$\PP^{\rm BM}_x(\tau_0>v)\sim\sqrt{2/\pi}\,x/\sqrt v$
as $v\to\infty$ for fixed $x>0$, where $\PP^{\rm BM}_x$ is
standard Brownian law started at $x$ and $\tau_0$ is its first
hitting time of zero.
Since $v\sim g\log(1/\varepsilon)$, the conversion factor
from derivative mass to the $\sqrt{\log(1/\varepsilon)}$ Wick
normalization is $\sqrt{2/(\pi g)}=1$. We prove the convergence in $L^2$.
At a fixed intermediate cutoff, conditional expectation gives the
survival approximation. Explicitly, for
$a=2\varepsilon^2<a'\le A$, Gaussian tilting and reflection give
\[
 \begin{split}
 \EE[W_{b,\varepsilon}(f)\mid\mathcal F_{a'}]
 =\sqrt{\log(1/\varepsilon)}\int_Q
 &f(x)e^{\alpha X_{a'}(x)-\alpha^2v_{a'}(x)/2}
       \ind_{\{R_{b,u}(x)\ge0\ \forall u\in[a',A]\}}\\
 &\cdot\left[2\Phi\left(\frac{R_{b,a'}(x)}
                     {\sqrt{v_a(x)-v_{a'}(x)}}\right)-1\right]dx.
 \end{split}
\]
Here $\Phi$ is the standard normal distribution function. Put
$D_{b,a'}(f)=\int_Q f(x)D_{b,a'}(x)\,dx$.
For fixed $a'$, the survival asymptotic and
$2\Phi(r)-1\le\sqrt{2/\pi}\,r$ for $r\ge0$ give
\[
 \EE[W_{b,\varepsilon}(f)\mid\mathcal F_{a'}]
       \longrightarrow D_{b,a'}(f)\quad\hbox{in }L^2.
\]
Indeed, on $K$ the ratio of the logarithmic normalization to
$\sqrt{v_a-v_{a'}}$ is uniformly bounded for small $\varepsilon$.
The conditional expectation is therefore bounded by a constant
multiple of $D_{b,a'}(f)$, which is in $L^2$ by
\eqref{chII:eq:barrier-pair}. Dominated convergence proves the assertion.
Since $D_{b,a'}(f)$ is $\mathcal F_{a'}$-measurable,
\begin{equation}\label{chII:eq:barrier-coarse-pairing}
 \EE[W_{b,\varepsilon}(f)D_{b,a'}(f)]
       \longrightarrow\EE[D_{b,a'}(f)^2].
\end{equation}
The spatial pair estimates in the next subsection give
\begin{equation}\label{chII:eq:barrier-second-moment-match}
 \EE[W_{b,\varepsilon}(f)^2]-\EE[D_{b,a}(f)^2]
       \longrightarrow0,\qquad a=2\varepsilon^2.
\end{equation}
Thus $\EE[W_{b,\varepsilon}(f)^2]\to\EE[M_b(f)^2]$, and the
$L^2$ norms of $W_{b,\varepsilon}(f)$ are eventually bounded by some
$C_f<\infty$. By Cauchy--Schwarz,
\[
 \bigl|\EE[W_{b,\varepsilon}(f)(M_b(f)-D_{b,a'}(f))]\bigr|
       \le C_f\|M_b(f)-D_{b,a'}(f)\|_2.
\]
First let $\varepsilon\downarrow0$ in
\eqref{chII:eq:barrier-coarse-pairing}, then let $a'\downarrow0$.
The mean-square convergence of $D_{b,a'}(f)$ yields
$\EE[W_{b,\varepsilon}(f)M_b(f)]\to\EE[M_b(f)^2]$. Consequently,
\[
 \EE|W_{b,\varepsilon}(f)-M_b(f)|^2
 =\EE[W_{b,\varepsilon}(f)^2]+\EE[M_b(f)^2]
      -2\EE[W_{b,\varepsilon}(f)M_b(f)]\longrightarrow0.
\]
At every cutoff,
\[
 0\le W_{b,\varepsilon}(f)\le\mathcal M_\varepsilon[X](f),\qquad
 W_{b,\varepsilon}(f)=\mathcal M_\varepsilon[X](f)\quad\hbox{on }E_b^c.
\]
Passing to their probability limits proves
$0\le M_b(f)\le M_Q(f)$ and $M_b(f)=M_Q(f)$ on $E_b^c$.

Under the one-point exponential tilt, the compensation removes
the drift of the variance-clock Brownian motion. Reflection at
the barrier then identifies the limiting mean with the initial
gap normalizer
\[
 \begin{aligned}
 q_b(x)&=\EE\left[e^{\alpha X_A(x)-\alpha^2v_A(x)/2}
              (b+\alpha v_A(x)-X_A(x))_+\right]\\
       &=\EE(b-X_A(x))_+.
 \end{aligned}
\]
The second equality is the Gaussian exponential-tilt identity.
Denote the standard normal density and distribution function by
$\varphi,\Phi$, respectively. Then
\[
 \begin{aligned}
 q_b(x)&=b\Phi(b/\sqrt{v_A(x)})
             +\sqrt{v_A(x)}\varphi(b/\sqrt{v_A(x)}),\\
 0\le q_b(x)-b&\le\sqrt{v_A(x)}/\sqrt{2\pi}.
 \end{aligned}
\]
In the present normalization this gives
\begin{equation}\label{chII:eq:barrier-first}
 \EE M_b(f)=\int f(x)q_b(x)\,dx,
 \qquad b\le q_b(x)\le b+C
\end{equation}
for nonnegative $f$ supported in a fixed core. The upper error
is the mean absolute value of the initial coarse Gaussian and
is uniformly bounded on that core. For fixed $b$, local mean-square convergence gives the cutoff
limit; the resulting estimate is then used as $b\to\infty$.

\subsection{The quantitative two-point and crossing inputs}

For the derivative barrier density above, the quantitative
two-point estimate used for the limit is,
with $u=\max(a,|x-y|^2)$, where $u_0=u_0(K,A)\in(0,A)$ is fixed, and
\begin{equation}\label{chII:eq:barrier-pair}
 \EE[D_{b,a}(x)D_{b,a}(y)]
 \le C(1+b)^2e^{\alpha b}
 \left[1+\frac{\ind_{\{u\le u_0\}}}
 {u\{1+\tfrac12\log(A/u)\}^{9/8}}\right],
 \qquad b\ge1.
\end{equation}
Fix a deterministic logarithmic time $T_0\ge1$ after which the
cell estimates below apply. Put
\[
 Y_0=\sup_{\substack{x\in K\\0\le t\le T_0}}
 \bigl(X_{Ae^{-2t}}(x)-\alpha v_{Ae^{-2t}}(x)\bigr).
\]
Thus $Y_0$ controls
the initial interval, and $E_b$ is the all-time crossing event
used for tests supported in $K$.
The associated all-time crossing bound on a fixed core is
\begin{equation}\label{chII:eq:barrier-crossing}
 \PP(E_b)\le\PP\{Y_0>b\}+C(b+C)^{3/4}e^{-\alpha b},
 \qquad \EE e^{\theta Y_0}<\infty\quad(\theta>0).
\end{equation}
The late first-crossing estimate is summable in time, with
weights $(1+j)^{-9/8}$; its sum gives the second term above.
The near-diagonal bound in \eqref{chII:eq:barrier-pair} is obtained
by splitting at time $|x-y|^2$ and applying the two-point
barrier estimate to the common coarse part. Bounded covariance
variation of the fine part costs a fixed change in $b$, absorbed
by $(1+b)^2e^{\alpha b}$. At distances below the cutoff the same
argument uses $u=a$. Away from the diagonal the covariance is
bounded and gives the constant term.

Here is the derivation of these exponents and the estimates.
At a fixed root $x\in K$ and cutoff $0<a\le A$, define probability laws by
\[
 \frac{d\PP^{\rm Wick}_{x,a}}{d\PP}
   =e^{\alpha X_a(x)-\alpha^2v_a(x)/2},\qquad
 \frac{d\PP^{\rm root}_{b,x,a}}{d\PP}
   =\frac{D_{b,a}(x)}{\EE D_{b,a}(x)}.
\]
Since the displayed root density has positive normalizer, it
defines a probability law; the martingale identity makes these laws
consistent on earlier cutoff sigma-fields. Below,
$\EE^{\rm root}$ and $\PP^{\rm Wick}$ suppress the fixed root and
the cutoff used in the displayed expression.
Under the one-point derivative tilt, the gap below the barrier
is a three-dimensional Bessel process in the variance clock,
with the initial gap distributed according to its size-biased
Gaussian law. This follows directly from the killed Brownian
transition density
\[
 p_v^+(a,r)=\frac1{\sqrt{2\pi v}}
 \left(e^{-(r-a)^2/(2v)}-e^{-(r+a)^2/(2v)}\right),\qquad a,r>0:
\]
its Doob transform is $(r/a)p_v^+(a,r)$. Since
$1-e^{-2ar/v}\le2ar/v$, this density is at most
$Cr^2v^{-3/2}$, uniformly in the initial $a$. There exists $u_0=u_0(K,A)>0$ such
that, for $0<u\le u_0$, the clock at cutoff
$u$ is bounded below by a constant times
$1+\tfrac12\log(A/u)$. Therefore, for every nonnegative $F$
with $\int_0^\infty r^2F(r)dr<\infty$,
\begin{equation}\label{chII:eq:bessel-damped}
 \EE^{\rm root}F(R_{b,u})
 \le C\{1+\tfrac12\log(A/u)\}^{-3/2}
                   \int_0^\infty r^2F(r)dr.
\end{equation}
Integrating against the size-biased initial-gap law gives the same
bound uniformly over later rooted cutoffs.

For $|x-y|^2\le u$, heat-kernel differences bound the variance of
$E_u=X_u(y)-X_u(x)$ and its covariance with the root clock by
fixed constants. Gaussian exponential tilting, including the
initial-gap size bias, consequently gives, for every $q\in\mathbb R$,
\[
 \EE^{\rm root}e^{qE_u}
 \le e^{q^2C/2+\alpha|q|C}(1+|q|C/b),\qquad b\ge1.
\]
Indeed the tilt shifts the initial Gaussian gap by a quantity of
absolute value at most $|q|C$; dividing its shifted positive-part
mean by the original normalizer costs at most $1+|q|C/b$.
Apply H\"older with exponents $4/3$ and $4$ to the exponentially
damped gap factors $(1+R_{b,u})e^{-\alpha R_{b,u}}$ and
$R_{b,u}^{-1}e^{-\alpha R_{b,u}}$ times $e^{qE_u}$.
Their $4/3$ powers satisfy the integrability condition in
\eqref{chII:eq:bessel-damped}; near zero the more singular weighted
integrand is $r^{2-4/3}$, which is integrable. The decay becomes
$\{1+\tfrac12\log(A/u)\}^{-9/8}$.

Writing the pair moment as the root normalizer times an expectation
under this rooted law, the second Wick factor is bounded by
$Ce^{\alpha b+\alpha^2\EE X_u(x)^2/2}$ times those damped gap
factors and the exponential increment. The normalizer is at most
$b+C$, and at criticality
$e^{\alpha^2\EE X_u(x)^2/2}\le C/u$.
This proves the near-diagonal estimate at cutoff $u$.
If $a<u=|x-y|^2\le u_0$, condition on $\mathcal F_u$. The
fine cross-covariance has total variation
\[
 4\int_a^u p_s^Q(x,y)\,ds
 \le C\int_0^u s^{-1}e^{-|x-y|^2/(8s)}\,ds\le L_0,
\]
where $L_0=L_0(K,A)$. In the conditional
two-point Gaussian tilt, each gap is shifted by at most $\alpha L_0$.
The killed endpoint-payoff bound consequently gives
\[
 \EE[D_{b,a}(x)D_{b,a}(y)\mid\mathcal F_u]
 \le e^{\alpha^2L_0}(1+L_0)
       D_{b+1+\alpha L_0,u}(x)D_{b+1+\alpha L_0,u}(y).
\]
Apply the preceding coarse estimate at the raised barrier
$b+1+\alpha L_0$ and absorb this fixed shift in its constant.
This gives \eqref{chII:eq:barrier-pair} with
$u=\max(a,|x-y|^2)$. For $u>u_0$ the same conditioning uses a
fixed time and yields the constant term in that bound.

For \eqref{chII:eq:barrier-second-moment-match}, put
\[
 \begin{gathered}
 V_{b,a}(x)=e^{\alpha X_a(x)-\alpha^2v_a(x)/2}
       \ind_{\{R_{b,u}(x)\ge0\ \forall u\in[a,A]\}},\qquad
 L_\varepsilon=\log(1/\varepsilon),\quad a=2\varepsilon^2,\\
 K_\varepsilon^W(x,y)=L_\varepsilon\EE[V_{b,a}(x)V_{b,a}(y)],\qquad
 K_a^D(x,y)=\EE[D_{b,a}(x)D_{b,a}(y)].
 \end{gathered}
\]
At distinct points of $K$, condition at a fixed coarse heat time $a'$,
let $a\downarrow0$, and then let $a'\downarrow0$.
The two-point survival calculation and
$4\int_0^{a'}p_s^Q(x,y)\,ds\to0$ for $x\ne y$ give
\[
 K_\varepsilon^W(x,y)-K_a^D(x,y)\longrightarrow0\qquad(x\ne y).
\]
The killed transition and inverse-gap bounds following
\eqref{chII:eq:bessel-damped}, together with \eqref{chII:eq:barrier-pair}, give
\[
 \begin{gathered}
 \sup_{\substack{(x,y)\in K^2\\|x-y|>r}}
       \{K_\varepsilon^W(x,y)+K_a^D(x,y)\}
       \le C_{b,K,A,r}\qquad(0<\varepsilon<\varepsilon_r),\\
 \lim_{r\downarrow0}\limsup_{\varepsilon\downarrow0}
 \int_{\substack{(x,y)\in K^2\\|x-y|\le r}}
       (K_\varepsilon^W+K_a^D)(x,y)\,dx\,dy=0.
 \end{gathered}
\]
The derivative term uses the integrable envelope in
\eqref{chII:eq:barrier-pair}; the Wick term uses the inverse-gap and
killed survival bounds. Dominated convergence on $|x-y|>r$ gives
\[
 \begin{aligned}
 \limsup_{\varepsilon\downarrow0}\int_{K^2}|K_\varepsilon^W-K_a^D|
 &\le\limsup_{\varepsilon\downarrow0}
   \int_{\substack{(x,y)\in K^2\\|x-y|\le r}}(K_\varepsilon^W+K_a^D)
 \longrightarrow0\quad(r\downarrow0).
 \end{aligned}
\]
Multiplication by $f(x)f(y)$ and Fubini give
\eqref{chII:eq:barrier-second-moment-match}.

For the crossing estimate, use logarithmic time
$a=Ae^{-2t}$ and divide a fixed core into $O(e^{2j})$ cells of
diameter $O(e^{-j})$ on each interval $[j,j+1]$. Choose a root
at the initial time of each cell. A first crossing in the cell
has root survival up to time $j$. Under the ordinary Wick tilt
at that root, reflection gives the past-strip bound
\[
 \PP^{\rm Wick}\{\hbox{survival to }j,
                 \hbox{terminal gap}\le u\}
       \le C(b+C)u^2j^{-3/2}.
\]
To check it, integrate $p_v^+(a,r)$ for $0<r\le u$ and then
average the initial $a$; the preceding kernel bound before the
Doob transform is $Carv^{-3/2}$. Its initial mean is at most
$b+C$. The spatial and temporal oscillation on a natural cell,
under this same tilt, has Gaussian tail $Ce^{-cj_0^2}$ above
$C+j_0$. After rescaling the cell to unit size, the centered
increments satisfy
\[
 \Var(\Delta X)\le C\bigl(|\Delta x_1|^2+|\Delta x_2|^2+|\Delta t|\bigr).
\]
Use a parabolic dyadic mesh with spatial step $2^{-l}$ and time
step $2^{-2l}$. Its adjacent and parent edges number $O(2^{4l})$,
and their standard deviations are $O(2^{-l})$. Gaussian edge
bounds at thresholds $C(1+j_0)\sqrt{l+1}\,2^{-l}$, summed over
edges and levels, give the asserted tail since
$\sum_l\sqrt{l+1}\,2^{-l}<\infty$. The Wick tilt leaves the
covariance unchanged, and the mean of the compensated field
is uniformly bounded on the cell. All bounds are uniform in
the cell and time.

In oscillation bin $j_0$, inverse Wick density is bounded by
$Ce^{-\alpha b+\alpha(C+j_0+1)}e^{-2j}$. Keep this factor
outside H\"older's inequality. Apply H\"older to the past-strip
and oscillation probabilities with powers $3/4$ and $1/4$.
The cell-bin probability is at most
\[
 Ce^{-\alpha b+\alpha(C+j_0+1)}e^{-2j}
 \{(b+C)(C+j_0+1)^2j^{-3/2}\}^{3/4}
                              e^{-cj_0^2/4}.
\]
Summation in $j_0$ is finite. Summation over the $O(e^{2j})$
cells cancels $e^{-2j}$ exactly and leaves
$C(b+C)^{3/4}e^{-\alpha b}j^{-9/8}$.
Summing in $j$ proves the late part of
\eqref{chII:eq:barrier-crossing}. On the fixed initial time interval
the smooth Gaussian field has the same Gaussian mesh bound
for its supremum; its compensated supremum is $Y_0$ and has
all positive exponential moments. This completes that estimate.

\subsection{Integration and the limit at fixed barrier}

The spatial kernel in \eqref{chII:eq:barrier-pair} is uniformly
integrable in $y$. On the shells
$Ae^{-2(j+1)}<|x-y|^2\le Ae^{-2j}$ its integral is bounded by
$C(1+j)^{-9/8}$; the innermost disk $|x-y|^2\le a$ has area
$O(a)$ and is bounded by the corresponding innermost shell
term. Summing, then integrating against $|f(x)f(y)|$, gives
\[
 \EE\left(\int f(x)D_{b,a}(x)dx\right)^2
       \le C\|f\|_\infty^2(1+b)^2e^{\alpha b}.
\]
At fixed $b$, local mean-square convergence passes this bound
to $M_b(f)$. Both the derivative densities and the truncated Wick
densities increase with $b$, so the limiting barrier measures do
as well. The truncated Wick identities above prove the domination
and the good-event identity in \eqref{chII:eq:barrier-package}.

For the crossing term, fix $\theta>\alpha$ in the exponential moment
bound in \eqref{chII:eq:barrier-crossing}. Markov's inequality gives
\[
 \frac{\PP(Y_0>b)}{be^{-\alpha b}}
 \le\frac{\EE e^{\theta Y_0}}b e^{-(\theta-\alpha)b}\longrightarrow0.
\]
The second term divided by $be^{-\alpha b}$ is $O(b^{-1/4})$.
This proves $\PP(E_b)=o(be^{-\alpha b})$.

Combining \eqref{chII:eq:barrier-first} with the preceding bounds gives
\begin{equation}\label{chII:eq:barrier-package}
 \begin{gathered}
 0\le M_b(f)\le M_Q(f),\qquad
 b\int f\le\EE M_b(f)\le b\int f+C\|f\|_\infty,\\
 \EE M_b(f)^2\le C\|f\|_\infty^2(1+b)^2e^{\alpha b},\qquad
 \PP(E_b)=o(be^{-\alpha b}),\\
 M_b(f)=M_Q(f)\quad\hbox{on }E_b^c.
 \end{gathered}
\end{equation}
Here $C=C(K,A)$; the same event $E_b$ serves every such $f$, and the
estimates are uniform over $f$ with bounded $\|f\|_\infty$.

Set $s=e^{-\alpha b}$ and fix $0<q<1$. From
$u-u^2\le1-e^{-u}\le u$ for $u\ge0$, we obtain
\[
 sqb\int f-s^2C\|f\|_\infty^2(1+qb)^2e^{\alpha qb}
 \le 1-\EE e^{-sM_Q(f)}
 \le s\left(b\int f+C\|f\|_\infty\right)+\PP(E_b).
\]
Divide by $s\log(1/s)=s\alpha b$. The second-moment remainder
vanishes exponentially since $q<1$, and the upper remainders vanish
by the crossing bound. Let $b\to\infty$, then $q\uparrow1$.
The estimates are uniform for $\|f\|_\infty$ bounded. The defect
itself tends uniformly to zero, and
$u\le-\log(1-u)\le u+2u^2$ for $0\le u\le1/2$
replaces the defect by the logarithm.

Choose continuous inner and outer cutoffs for $Q_\delta$, with
transition in a strip of width $\tau$, contained in a fixed larger
compact interior. Apply the result to each cutoff times
$e^{\alpha d}w$. The error between the two deterministic integrals
is bounded by $e^{\alpha B}\|w\|_\infty$ times the strip area,
uniformly in $d$. Monotonicity sandwiches the hard-core Laplace transform. Let
$\tau\downarrow0$; the passage uses the uniform $L^\infty$ bound on the
profiles together with bounded truncation.

\appendixchapter{G}{Proofs of Propositions~\ref{chII:prop:scalar-gluing} and~\ref{chII:prop:height-average}}

We use the squares $D_i$, fine fields $\eta_i$, coarse field $H$ and
mean correction $\overline\eta$ of Section~2.1. Put
$m=\lfloor N/K\rfloor$ and $H_i=H_{N,K,i}$. For $0<\delta<1/2$ the
retained cores are
\[
 D_{i,\delta}=\{im+x:x\in\{1,\ldots,m-1\}^2,\ x/m\in Q_\delta\}.
\]
Their complement $V_N\setminus\bigcup_iD_{i,\delta}$ is the deleted collar;
it includes the grid boundaries and any remainder from
rounding $N/K$. Conditional probabilities in the block comparison
integrate over the independent fine fields with the entire coarse
field fixed; the conditioning sigma-field is $\mathcal H=\sigma(H)$.

\subsection{Bounded tests in the coarse-field limit}\label{chII:sec:coarse-limit}

For the continuum comparison put
\[
 D_{K,i}=(i+Q)/K,\qquad
 U_{K,\delta}=\bigcup_{i\in\{0,\ldots,K-1\}^2}(i+Q_\delta)/K.
\]
The continuum orthogonal square decomposition on the original
periodic-field space is
\[
 h=\sum_i\eta_{K,i}+H_K-\overline\eta_K,\qquad
 \overline\eta_K=\sum_i\int_{D_{K,i}}\eta_{K,i},\qquad
 H_K\perp(\eta_{K,i})_i.
\]
Here the $\eta_{K,i}$ are independent ordinary Dirichlet fields on
$D_{K,i}$, with covariance $G_{D_{K,i}}$ as in Section~2.3 and zero
extension outside their squares. The integrals are Gaussian pairings.
The field $H_K$ is smooth inside each square. Write
$H_{K,i}(z)=H_K((i+z)/K)$.

For a finite positive measure $\mu$ and continuous $u,v$, set
$L_\mu(u)=\exp\{-a\int e^{\alpha u}d\mu\}$, $a\ge0$.
Along the segment $u+t(v-u)$, differentiation gives
\[
 \left|\frac d{dt}L_\mu(u+t(v-u))\right|
 \le\alpha\|u-v\|_\infty\sup_{s\ge0}se^{-s}
 =\frac\alpha e\|u-v\|_\infty.
\]
Thus the estimate is uniform over all finite measures $\mu$ and all
$a$ in the stated range, and it survives integration over any law of
finite measures. This
is the bounded-test estimate used to approximate coarse fields
on a fixed mesh, then remove the mesh approximation by their
uniform modulus of continuity. Convergence of the covariance
matrices handles the finite mesh.

At lattice points $x\in D_i$, $y\in D_j$, independence in
\eqref{chII:eq:block-splitting} gives the covariance identity
\[
 \operatorname{Cov}(H(x),H(y))
 =G_N(x,y)-\ind_{\{i=j\}}G_{D_i}(x,y)
   +N^{-2}\{b_{D_i}(x)+b_{D_j}(y)\}
   -\operatorname{Var}(\overline\eta).
\]
The Green limits, including their diagonal finite parts, and the
Riemann sums for $b_{D_i}$ give convergence to the analogous continuum
kernel at fixed $K$ on any product of closed cores. On one core the
logarithmic singularities in $G_N-G_{D_i}$ cancel. The same Green
increment bounds, extended by the tensor-hat weights, give
\[
 \sup_{N\ge N_0(K,\delta)}\EE|H_i(z)-H_i(z')|^{2p}
 \le C_{K,\delta,p}|z-z'|^p,
 \qquad z,z'\in Q_\delta,\quad p>2,
\]
together with a bounded variance at one point. Gaussian moment
identities and the two-dimensional Kolmogorov criterion therefore
give tightness in $C(Q_\delta)$. The limiting covariance identifies
every finite-dimensional limit, proving the joint convergence of all
$K^2$ cores used in Proposition~\ref{chII:prop:scalar-gluing}.

The decomposition above is the continuum orthogonal projection
counterpart of \eqref{chII:eq:block-splitting}. The fine mean is
generally dependent on the fine chaos. With $b_D,v_D$ as in Section~2.3,
\[
 \operatorname{Var}(\overline\eta_K)=\sum_i v_{D_{K,i}}
   =O(K^{-2}),\qquad
 \sup_{x\in D_{K,i}}b_{D_{K,i}}(x)=O(K^{-2}).
\]
The first bound follows by scaling the integrated Dirichlet Green
kernel; the second is the continuum exit-time bound. For
$x\in D_{K,i}$, the covariance of the regularized fine field with
$\overline\eta_K$ tends to $b_{D_{K,i}}(x)$. Comparing variance
finite parts therefore gives
\[
 \operatorname{Var}(H_K(x))-2b_{D_{K,i}}(x)
        +\operatorname{Var}(\overline\eta_K)
   =r_{\rm per}-r_{D_{K,i}}(x).
\]
Let $M_{D_{K,i}}$ be the exact-variance chaos measurable from
$\eta_{K,i}$, and define on the unit square
\[
 \mathcal M_i^K(dz):=K^2w(z)\,
 \bigl[(x\mapsto Kx-i)_*M_{D_{K,i}}\bigr](dz).
\]
The map $x\mapsto Kx-i$ sends $D_{K,i}$ to $Q$; Dirichlet scaling
therefore gives $\mathcal L(\mathcal M_i^K)=\mathcal L(\mathcal M)$.
The smooth change rule and heat locality from
Appendix~C give the measure identity on compact interiors
\[
 \mathcal Z(dx)
 =e^{\alpha(H_K(x)-\overline\eta_K)
       -\alpha^2(r_{\rm per}-r_{D_{K,i}}(x))/2}
       M_{D_{K,i}}(dx).
\]
 Under $x=(i+z)/K$, the physical
heat scale $\varepsilon$ becomes $K\varepsilon$ and
\[
 \sqrt{\frac{\log(1/\varepsilon)}{\log(1/(K\varepsilon))}}
       \longrightarrow1\qquad(\varepsilon\downarrow0,\ K\text{ fixed}).
\]
Thus the exact-variance fine measure contributes the physical
area factor $K^{-2}$, while
$r_{D_{K,i}}(x)=r_Q(z)-g\log K$ contributes the finite-part factor
$K^{-2}$. Their product is $K^{-4}$, the normalization displayed above. Consequently
\[
 \mathcal Z\bigl((i+Q_\delta)/K\bigr)
 =e^{-\alpha\overline\eta_K}\frac{K^{-4}}{p_{\rm per}}
     \int_{Q_\delta}e^{\alpha H_{K,i}(z)}\,\mathcal M_i^K(dz).
\]
Field measurability, Dirichlet scaling and independence give
\[
 \mathcal L\bigl((H_{K,i})_i,(\mathcal M_i^K)_i\bigr)
 =\mathcal L\bigl((H_{K,i})_i\bigr)
       \otimes\mathcal L(\mathcal M)^{\otimes K^2}.
\]
Together with the fixed-$K,\delta$ coarse-field convergence proved
above, bounded-test integration gives
\[
 B_{N,K,\delta}\Rightarrow B_{K,\delta}^{\rm cont}
 :=\frac{K^{-4}}{p_{\rm per}}\sum_i
       \int_{Q_\delta}e^{\alpha H_{K,i}(z)}\,\mathcal M_i^K(dz).
\]
Summing the disjoint-core measure identities yields
$B_{K,\delta}^{\rm cont}=e^{\alpha\overline\eta_K}\mathcal Z(U_{K,\delta})$.
Since $\overline\eta_K\to0$ in probability and $Z<\infty$ almost
surely,
\[
 \begin{gathered}
 |e^{\alpha\overline\eta_K}-1|Z\longrightarrow0
       \quad\hbox{in probability},\\
 \EE\min\{1,|e^{\alpha\overline\eta_K}-1|Z\}\longrightarrow0.
 \end{gathered}
\]

For collar removal, fix $M>0$ and define the finite deterministic measure
$\mu_M(E)=\EE[\ind_{\{Z\le M\}}\mathcal Z(E)]$.
Translation invariance of the construction implies that $\mu_M$
is a constant multiple of Lebesgue measure, of total mass at
most $M$. Hence for any deterministic union $E$ of collars,
\begin{equation}\label{chII:eq:collar-mass}
 \PP\{\mathcal Z(E)>e\}
 \le\PP(Z>M)+e^{-1}M|E|.
\end{equation}
The total area of the collars is $O(\delta)$, uniformly in $K$.
Given an error tolerance, choose a mass cutoff $M$, then the collar
width $\delta$, and then $K$ for the mean correction. With $K,\delta$
fixed, the discrete limit yields the uniform comparison in
Proposition~\ref{chII:prop:scalar-gluing}.

\subsection{The exact conditional product and its error bounds}

For $h\in\mathbb R$, define
\[
 q_i(h,H)=\PP\{\exists x\in D_{i,\delta}:
       |\eta_i(x)+H(x)-m_N-h|\le1/2\mid\mathcal H\}.
\]
Conditional independence gives
\begin{equation}\label{chII:eq:contact-product}
 \EE\!\left[(1-\lambda)^{\sum_i\ind_{\{
       \exists x\in D_{i,\delta}:|\eta_i(x)+H(x)-m_N-h|\le1/2\}}}
       \,\middle|\,\mathcal H\right]
 =\prod_i(1-\lambda q_i(h,H)).
\end{equation}
For numbers $u_i,v_i\in[0,1]$, telescoping the product gives
\begin{equation}\label{chII:eq:product-error}
 \left|\prod_i(1-u_i)-\prod_i(1-v_i)\right|
       \le\min\left(1,\sum_i|u_i-v_i|\right).
\end{equation}
Deleting indices is the case $v_i=0$ on those indices in
\eqref{chII:eq:product-error}.

If $B$ is the number of conditionally independent contacts and
$R=\sum_iq_i$, then for $0<c\le1$ and $L,D>0$,
\begin{equation}\label{chII:eq:conditional-rate-tail}
 (1-e^{-cL})\PP(R>L)\le cD+\PP(B>D).
\end{equation}
Indeed, $\prod_i(1-cq_i)\le e^{-cR}$ and the conditional product identity give
\[
 (1-e^{-cL})\PP(R>L)
 \le\EE[1-(1-c)^B]\le cD+\PP(B>D).
\]
Given an error tolerance, choose the mass cutoff $D$, then a small
$c$, and finally $L$ so that the displayed rate tail is below that
tolerance.

\subsection{Weighted Gaussian deletion on the original coarse field}

For the weighted deletion argument, first normalize the conditional
profile weights and then apply the Gaussian tail estimate. Let
$p_i(h,H)=\PP\{\max_{x\in D_{i,\delta}}(\eta_i(x)+H(x))>m_N+h\mid\mathcal H\}$
be the conditional exceedance probability in the retained core, and put
\[
 R(h,H)=\sum_i p_i(h,H),\qquad
 w_i(h,c,H)=\frac{p_i(h,H)}{1+R(h-2c,H)},\qquad c\ge0.
\]
These are exceedance probabilities, as opposed to the band
probabilities in \eqref{chII:eq:contact-product}. For deterministic
perturbations $v_i:V_N\to\mathbb R$ with $\|v_i\|_\infty\le c$,
monotonicity of the exceedance event gives
\[
 p_i(h,H+v_i)\le p_i(h-c,H),\qquad
 R(h-c,H)\le R(h-2c,H+v_i).
\]
For block-dependent perturbations $v_i$ with the stated uniform
bound, monotonicity gives the common comparison
\begin{equation}\label{chII:eq:normalized-shift-budget}
 \sum_i w_i(h,c,H+v_i)
 \le\frac{R(h-c,H)}{1+R(h-c,H)}\le1.
\end{equation}

\begin{lemma}\label{chII:lem:weighted-gaussian-tail}
Let $\mathcal J$ be a finite index set, let $\xi$ be a standard
Gaussian vector, let $0\le f_i\le1$, $i\in\mathcal J$, be measurable,
and let $a_i$ be deterministic vectors of the same dimension. For $u\in\mathbb R$, suppose
$t\ge0$, $t^2\|a_i\|^2\le V$, and
$\sum_i f_i(\xi+ta_i)\le C$ for every $\xi$. Then
\begin{equation}\label{chII:eq:weighted-gaussian-tail}
 \sum_i\EE[\ind_{\{\langle a_i,\xi\rangle\ge u\}}f_i(\xi)]
 \le C e^{-tu+V/2}.
\end{equation}
\end{lemma}
\begin{proof}
On $\langle a_i,\xi\rangle\ge u$,
\[
 1\le e^{-tu+V/2}
       e^{t\langle a_i,\xi\rangle-t^2\|a_i\|^2/2}.
\]
Multiplication by $f_i\ge0$ and Gaussian translation give
\[
 \EE[\ind_{\{\langle a_i,\xi\rangle\ge u\}}f_i(\xi)]
 \le e^{-tu+V/2}\EE f_i(\xi+ta_i).
\]
Sum and use the translated budget. By the Cameron--Martin
translation identity, the displayed mean-one exponential sends every
bounded measurable test to its translated Gaussian expectation.
\end{proof}

Apply this to the Fourier Gaussian coordinates of the actual
coarse field. If $Y_i$ is a centered linear Gaussian functional,
translation in its coefficient vector by $t$ changes $H(x)$ by
$t\Cov(H(x),Y_i)$. Thus, whenever
\[
 t^2\Var(Y_i)\le V,\qquad
 \sup_x|t\Cov(H(x),Y_i)|\le c,
\]
\eqref{chII:eq:normalized-shift-budget} verifies the hypothesis with
$f_i(\xi)=w_i(h,c,H(\xi))$ and $C=1$. In particular
\begin{equation}\label{chII:eq:actual-weighted-tail}
 \sum_i\EE[\ind_{\{Y_i\ge u\}}w_i(h,c,H)]\le e^{-tu+V/2}.
\end{equation}

Here is the countable-grid application. Let
$\Pi_{\delta/2}:\overline Q\to Q_{\delta/2}$ clamp each coordinate to
$[\delta/2,1-\delta/2]$. Take the deterministic grid center
$z_i=(\lfloor m/2\rfloor/m,\lfloor m/2\rfloor/m)$, which lies in
$Q_{\delta/2}$ for all sufficiently large $m$. Define the extended profile
\[
 d_i(z)=H_i(\Pi_{\delta/2} z)-H_i(z_i),\qquad z\in\overline Q.
\]
It agrees with $H_i-H_i(z_i)$ on $Q_\delta$. Because the clamping map is $1$-Lipschitz and maps the square into
the retained core, the interior covariance and increment bounds extend
to the whole square. For sufficiently large $m$ with
$\delta m\ge4$, the Green difference estimates, followed by bilinear
interpolation, give constants $C_1,C_2(\delta)$ such that
\begin{equation}\label{chII:eq:profile-covariance-input}
 \Var(d_i(z)-d_i(z'))\le C_1\delta^{-1}|z-z'|,
 \qquad
 |\Cov(H(x),d_i(z)-d_i(z'))|^2\le C_2(\delta)|z-z'|.
\end{equation}
The second estimate holds for the full coarse field at every
$x\in V_N$. Both estimates follow
by applying the covariance difference bounds to the corner
values in each interpolation cell; the interpolation
coefficients are nonnegative, have sum one, and only adjacent
corner differences remain after subtraction.

Take $\mathcal N_j=\{0,2^{-j},\ldots,1\}^2$. For $j\ge1$, map a
point of $\mathcal N_j$ to $\mathcal N_{j-1}$ by rounding both
coordinates down. Let $\mathcal E_j$ consist of these parent edges
and the horizontal and vertical nearest-neighbor edges in $\mathcal N_j$;
$\mathcal E_0$ consists just of the latter edges.
Their lengths are at most $24\,2^{-j}$ and their number at most
$2048\,256^j$. Put
\[
 \mathcal B_i(B)=\bigcup_{j\ge0}\ \bigcup_{(z,z')\in\mathcal E_j}
       \{|d_i(z)-d_i(z')|\ge B(j+1)(3/4)^j\}.
\]
Apply the weighted-tail estimate with $t=(4/3)^j$. The identities
$(4/3)^jB(j+1)(3/4)^j=B(j+1)$ and
$(4/3)^{2j}2^{-j}=(8/9)^j\le1$ show that
\eqref{chII:eq:profile-covariance-input} verifies the hypotheses of
\eqref{chII:eq:actual-weighted-tail} with constants $V_\delta=V(\delta)$ and $c_\delta=c(\delta)$,
uniformly in $j,B,K,N$. Apply the estimate to both signs of
every edge, then sum over its labels. For $B\ge16$,
\begin{align}
 \sum_i\EE[\ind_{\mathcal B_i(B)}w_i(h,c_\delta,H)]
 &\le\sum_{j\ge0}4096\,256^j e^{-B(j+1)+V_\delta/2}\notag\\
 &\le8192e^{V_\delta/2-B}.
 \label{chII:eq:bad-profile-weighted}
\end{align}
The last inequality uses $256e^{-16}\le1/2$.
Summing first over finitely many levels and then sending the level
cutoff to infinity, monotone convergence gives the displayed bound.

On the complement of $\mathcal B_i(B)$, chains from a point
to its nearest net points have successive increments bounded
by $B(j+1)(3/4)^j$. Since $\sum_j(j+1)(3/4)^j<\infty$, these chains
converge. If $|z-z'|$ is of order $2^{-j}$, join their two
chains through a bounded number of scale-$j$ edges. The
resulting difference is bounded by
$CB(j+1)(3/4)^j$. Continuity extends the estimate from the
nets to the whole square. For each
$0<\rho<\log(4/3)/\log2$ and the centered profile $d_i$,
\[
 \mathcal B_i(B)^c\quad\Longrightarrow\quad
 \|d_i\|_\infty+[d_i]_\rho\le C_\rho B.
\]
Thus these profiles satisfy the compact-class hypotheses of
Proposition~\ref{chII:prop:reference-contact}.

\begin{lemma}\label{chII:lem:normalized-deletion}
For any nonnegative random $R$, measurable $0\le p_i\le1$,
and events $E_i$, set
$D=\sum_i p_i\ind_{E_i}$. For $Q\ge0$,
\begin{equation}\label{chII:eq:normalized-deletion}
 \EE\min(1,D)\le\PP(R>Q)
       +(1+Q)\sum_i\EE\frac{p_i\ind_{E_i}}{1+R}.
\end{equation}
\end{lemma}

Take $R=R(h-2c_\delta,H)$, $p_i=p_i(h,H)$ and
$E_i=\mathcal B_i(B)$. The weighted Gaussian bound and
\eqref{chII:eq:normalized-deletion} give
\begin{equation}\label{chII:eq:bad-deletion-final}
 \EE\min\left(1,\sum_i p_i(h,H)\ind_{\mathcal B_i(B)}\right)
 \le\PP\{R(h-2c_\delta,H)>Q\}
       +(1+Q)8192e^{V_\delta/2-B}.
\end{equation}
Apply \eqref{chII:eq:conditional-rate-tail} to the independent
exceedance indicators at this fixed lower height. Their count
is tight by the high-point count estimate. Fix $\delta$. Given $\varepsilon$, choose $Q$ so that the first
term is below $\varepsilon/2$, then choose $B=B(Q,\delta,\varepsilon)$
so that the second term is below $\varepsilon/2$. A band contact is contained in the exceedance event at its lower
band endpoint, so the band-product contribution is bounded by the same
deletion estimate at that height.

The copied chaos masses require their own deletion budget. On the
product space of $H$ and the independent $\mathcal M_i$, put
\[
 b_i(H,\mathcal M)=\frac{K^{-4}}{p_{\rm per}}
       \int_{Q_\delta}e^{\alpha H_i}\,d\mathcal M_i,
 \qquad B_*:=\sum_i b_i=B_{N,K,\delta}.
\]
These are finite almost surely. If $\|v_i\|_\infty\le c$ on every
core, where the full-field perturbation $v_i$ may depend on the
selected index, then samplewise
\[
 \frac{b_i(H+v_i,\mathcal M)}{1+B_*(H+v_i,\mathcal M)}
 \le e^{2\alpha c}\frac{b_i(H,\mathcal M)}{1+B_*(H,\mathcal M)},
 \qquad
 \sum_i\frac{b_i(H+v_i,\mathcal M)}{1+B_*(H+v_i,\mathcal M)}
 \le e^{2\alpha c}.
\]
Condition on the fine measures and apply
Lemma~\ref{chII:lem:weighted-gaussian-tail} to these bounded weights.
The same two-sign grid sum as in \eqref{chII:eq:bad-profile-weighted} gives
\[
 \sum_i\EE\frac{b_i\ind_{\mathcal B_i(B)}}{1+B_*}
 \le8192e^{V_\delta/2+2\alpha c_\delta-B}.
\]
Consequently, for $B_*^{\rm bad}:=\sum_i b_i\ind_{\mathcal B_i(B)}$,
\[
 \EE\min(1,B_*^{\rm bad})
 \le\PP(B_*>Q)+(1+Q)8192e^{V_\delta/2+2\alpha c_\delta-B}.
\]
Choose $Q$ from scalar-mass tightness and then choose $B$ from the weighted profile bound. For every $a\ge0$,
\[
 0\le\EE e^{-a(B_*-B_*^{\rm bad})}-\EE e^{-aB_*}
 \le\max(1,a)\EE\min(1,B_*^{\rm bad}).
\]
These estimates also restore mass on a deeper core, since each
deleted deeper-core mass is bounded by its mass on $Q_\delta$.
They use finite sampled masses and bounded normalized weights only.

\subsection{A deterministic repair of the whole-square maximum}

The following profile penalty compares core band contacts with the
whole-square maximizing-root event in the reference estimate. Let $0<\rho<1$, $\|d\|_\infty+[d]_\rho\le B$, and
fix $0<\zeta<\delta$. Set
\[
 \vartheta(z)=\min\{1,\operatorname{dist}_\infty(z,Q_\delta)/\zeta\},
 \qquad 0\le\vartheta\le1,\qquad [\vartheta]_\rho\le\zeta^{-1}.
\]
The last bound follows from the $1$-Lipschitz property of the
distance function and $|z-z'|_\infty\le1$ on $\overline Q$.
Thus $\vartheta=0$ on $Q_\delta$ and $\vartheta=1$ outside
$Q_{\delta-\zeta}$. Replace $d$ by
\[
 d^-=d-(B+1+D)\vartheta,\qquad D\ge0.
\]
In particular,
\[
 \|d^-\|_\infty+[d^-]_\rho
 \le B+(B+1+D)(1+\zeta^{-1})=:B_{\rm pen}.
\]
Using $D_{M,\delta}$ from Section~2.4, define
\[
 \begin{aligned}
 S_M&=\{x\in D_M:x/M\in Q_{\delta-\zeta}\setminus Q_\delta\},\\
 O_M&=\{\exists x\in D_M\setminus D_{M,\delta}:\
                \eta_M^{d^-}(x)\ge m_M+t-1/2\}.
 \end{aligned}
\]
Because the penalty vanishes on $Q_\delta$, the penalized and
original profiles agree there. A point outside this core whose penalized value is
at least $m_M+t-1/2$ must either lie in the strip
$Q_{\delta-\zeta}\setminus Q_\delta$ and have unprofiled value
greater than $m_M+t-B-1$, or have unprofiled value at least
$m_M+t+D$. The reference subset and upper-tail bounds therefore
give, for $t\ge\max(1,D)$,
\begin{equation}\label{chII:eq:outer-penalty-budget}
 \PP(O_M)
 \le C\left(\sqrt{|S_M|/M^2}\,e^{\alpha(B+1)}
                 +2e^{-\alpha D}\right)te^{-\alpha t},
\end{equation}
Here $S_M$ is the strip at lattice scale.
For the second term we used
$(t+D)e^{-\alpha(t+D)}\le2e^{-\alpha D}te^{-\alpha t}$ when $t\ge D$.
The strip satisfies $|S_M|/M^2\le C\zeta+C/M$ by counting
its four rows of rectangles. At fixed $B$, choose $D$ and $\zeta$ so that the first two error
terms are small, and then take $M$ large for the lattice-strip error. The coefficient in
\eqref{chII:eq:outer-penalty-budget} becomes arbitrarily small.

For a finite nonempty $W\subset\mathbb Z^2$ and $T\ge1$, put
\[
 \begin{aligned}
 C_M&=\{\exists x\in D_{M,\delta}:|\eta_M^d(x)-m_M-t|\le1/2\},\\
 H_M&=\{\max_{x\in D_{M,\delta}}\eta_M^d(x)>m_M+t+T\},\\
 P_M&=\{\exists x,y\in D_{M,\delta}:x-y\notin W,\
           |\eta_M^d(x)-m_M-t|\le1/2,\
           \eta_M^d(y)\ge m_M+t-1/2\}.
 \end{aligned}
\]
Let $A_M^-$ be the event defining $p_M(\delta,t,T,W,d^-)$ in Section~2.4.
Since $d^-=d$ on $Q_\delta$, comparison with the maximizer $y_M^{d^-}$ gives
\[
 A_M^-\subseteq C_M\subseteq A_M^-\cup O_M\cup H_M\cup P_M.
\]
Choose the penalty parameters so that $\PP(O_M)\le\varepsilon te^{-\alpha t}$.
Taking probabilities gives
\begin{equation}\label{chII:eq:contact-repair-explicit}
 \begin{split}
 0&\le\PP(C_M)-p_M(\delta,t,T,W,d^-)\\
  &\le\varepsilon te^{-\alpha t}+\PP(H_M)+\PP(P_M).
 \end{split}
\end{equation}
Here $\eta_M^{d^-}|_{D_{M,\delta}}=\eta_M^d|_{D_{M,\delta}}$,
so $H_M$ and $P_M$ are events for the original core field.

\subsection{Height averaging under a small coupled perturbation}

We isolate the analytic part of the adjacent-square argument.
Fix $H_0\ge0$. For a finite family of finite sets $S_i$ and real heights $a_x$, put
\[
 P_a(h)=\prod_i
 \left(1-\lambda\ind_{\{\exists x\in S_i:|a_x-h|\le1/2\}}\right),
 \qquad 0\le\lambda\le1.
\]
Suppose $|a_x-b_x|\le e\le1/2$ on $\bigcup_iS_i$.
If the two band indicators at a site differ at height $h$,
then $h$ lies in one of the two intervals of radius $e$
centered at $a_x-1/2$ and $a_x+1/2$. Their total length is
at most $4e$. Moreover, if $|h|\le H_0$, that site has
$a_x\ge-H_0-1$. The products can differ only at such a site,
and both lie in $[0,1]$. Hence
\begin{equation}\label{chII:eq:height-perturbation-budget}
 \int_{-H_0}^{H_0}|P_a(h)-P_b(h)|\,dh
 \le4e\,\#\{x\in\bigcup_iS_i:a_x\ge-H_0-1\}.
\end{equation}
For random coupled fields define
$N_a=\#\{x\in\bigcup_i S_i:a_x\ge-H_0-1\}$ and
$E_e=\{\max_{x\in\bigcup_iS_i}|a_x-b_x|>e\}$, using the convention $\max\varnothing=0$.
For $Q\ge0$, truncating $N_a$ at $Q$ and applying Fubini gives
\begin{multline}\label{chII:eq:height-perturbation-random}
 \int_{-H_0}^{H_0}\EE|P_a(h)-P_b(h)|\,dh\\
 \le4eQ+2H_0\{\PP(E_e)+\PP(N_a>Q)\}.
\end{multline}

For adjacent sides $m$ and $m+1$, realize the two ordinary
Dirichlet fields on the same ambient torus using their nested
domains. Their difference on the smaller square is the
independent harmonic increment, with covariance
$G_{D_{m+1}}-G_{D_m}$. The interior diagonal Green expansion
implies that its variance tends uniformly to zero on each
fixed core. The interior difference estimates give tightness
of these harmonic functions in the uniform topology there:
use a fixed mesh, bound mesh values by their vanishing
Gaussian variances, and use the common increment modulus to
pass from the mesh to the core. Thus the supremum of the
increment tends to zero in probability.

To apply \eqref{chII:eq:height-perturbation-random}, denote the coupled
fields in block $i$ by $\eta_{i,m}$ and $\eta_{i,m+1}$, and fix
$\delta<\delta'<1/2$. Use the labeled common core
\[
 S_i=\{(i,x):x\in D_m,\ x/m\in Q_{\delta'},\
                              x/(m+1)\in Q_{\delta'}\}.
\]
For the original extended coarse profile $d_i$ and its center value
$H_i(z_i)$, set
\[
 \begin{aligned}
 a_{i,x}&=\eta_{i,m}(x)+d_i(x/m)+H_i(z_i)-m_N,\\
 b_{i,x}&=\eta_{i,m+1}(x)+d_i(x/(m+1))+H_i(z_i)-m_N.
 \end{aligned}
\]
Thus $|a_{i,x}-h|\le1/2$ is the original band event on $S_i$;
$|b_{i,x}-h|\le1/2$ is the adjacent-square band event, with effective
height $m_N-m_{m+1}+h-H_i(z_i)$. On a fixed good profile class,
\[
 \max_{(i,x)\in\bigcup_iS_i}|a_{i,x}-b_{i,x}|
 \le\max_{(i,x)\in\bigcup_iS_i}|\eta_{i,m+1}(x)-\eta_{i,m}(x)|
                 +B(\sqrt2/m)^\rho\longrightarrow0
 \quad\hbox{in probability}
\]
at fixed $K,\delta,\delta'$. Take $N_a$ and $E_e$ from these labeled heights in
\eqref{chII:eq:height-perturbation-random}. The high-point estimate
bounds $\PP(N_a>Q)$, and the displayed convergence bounds $\PP(E_e)$.
Given a tolerance, first choose $Q$, then $e$ with $4eQ$ below that
tolerance, and then let $m\to\infty$ at fixed $K$. Sites outside the
common core lie in the collars between $Q_\delta$ and $Q_{\delta'}$,
apart from a vanishing lattice displacement; their contribution is
bounded by the collar high-point estimate. Let $\delta'\downarrow\delta$
after the side-length limit. This proves the height-integral comparison
for the original odd side from that for the adjacent even side.
Finally $m_{m+1}-m_m\to0$ and the deterministic scalar-mass centering
factor tends to one at fixed $K$; the latter replacement follows by
restricting the total mass to a fixed cutoff and then removing it.

\subsection{Standard high-point bounds and their periodic transfer}

We use the Dirichlet subset estimate
\eqref{chII:eq:dirichlet-subset-cited} and the following square case of the
separation theorem in \cite[Lemma B.11]{chII:BL}. For the reference square
$D_M$ and its Dirichlet field $\eta_M$ of Section~2.4, and each fixed $H$,
\begin{equation}\label{chII:eq:dirichlet-separation-cited}
 \lim_{r\to\infty}\limsup_{M\to\infty}
 \PP\{\exists x,y\in D_M:\ \eta_M(x),\eta_M(y)\ge m_M-H,
                      \ r<|x-y|<M/r\}=0.
\end{equation}
For squares this is the Ding--Zeitouni high-point separation
theorem, as recorded in the cited lemma. We transfer the cited
subset and separation bounds to the periodic mean-zero field by
the constrained decomposition below.

Cover the torus by a fixed number of cores $D'_j$ of embedded
squares $D_j$ of side comparable to $N$. Choose the cover with
a fixed positive margin: each point lies in some $D'_j$ at
distance at least $cN$ from its boundary. A square grid of
centers of spacing $N/J$, for a fixed sufficiently large integer $J$,
provides such a cover. Choose the nearest lattice centers; for
large $N$ their $O(1)$ displacement is smaller than half the fixed
margin. Thus every pair at distance $o(N)$ lies in a common core.
All these domains are proper embedded planar squares.

Use the exact constrained decomposition of
Proposition~\ref{chII:prop:domain-fields} on one such square $D$:
\[
 \phi_N=\eta_D^0+H_D^{\mathrm{zm}},\qquad
 \eta_D^0\perp H_D^{\mathrm{zm}},\qquad
 \eta_D=\eta_D^0+q_D\xi_D,
\]
where $\xi_D$ can be taken independent of both displayed fields.
The latter field $\eta_D$ has the ordinary Dirichlet law.
The profile $q_D=b_D/\sqrt{v_D}$ is nonnegative and uniformly
bounded: exit-time estimates give $b_D\le CN^2$, and the
same estimates on a fixed inner square give $v_D\ge cN^4$.
On any fixed interior core, the covariance identity in that
proposition and the Green difference estimates give
\begin{equation}\label{chII:eq:fixed-core-coarse-bounds}
 \sup_{x\in D'}\Var(H_D^{\mathrm{zm}}(x))\le C,
 \qquad
 \EE\bigl(H_D^{\mathrm{zm}}(x)-H_D^{\mathrm{zm}}(y)\bigr)^2
 \le C|x-y|/N.
\end{equation}
Put $F_D=G_N-G_D$. The stopped Green identity gives
\[
 F_D(x,z)=\mathbb E_xG_N(X_{\tau_D},z)-b_D(x)/N^2.
\]
When $z$ is in a slightly larger interior core, the exit point
is at distance at least $cN$ from $z$, so the boundary values
in this formula are uniformly bounded. Interior harmonic
gradient estimates and the Poisson estimate for $b_D$ give
$|F_D(x,z)|\le C$ and
$|F_D(x,z)-F_D(y,z)|\le C|x-y|/N$ on the smaller core.
By symmetry the same holds in the second variable. Apply these
bounds to the four terms of the increment variance. The
rank-one correction contributes $(q_D(x)-q_D(y))^2$, which
obeys the same bound since $|q_D(x)-q_D(y)|\le C|x-y|/N$.
The rank-one term therefore obeys the same increment bound, proving
\eqref{chII:eq:fixed-core-coarse-bounds}.

Interpolate to the unit core and use dyadic
nets with $O(4^k)$ edges at level $k$. Each edge variance is at
most $C2^{-k}$. The probability that some level-$k$ edge exceeds
$A(t+\sqrt{k+1})2^{-k/4}$ is at most
$C4^k\exp\{-cA^2(t+\sqrt{k+1})^22^{k/2}\}$.
Choose $A$ fixed and large and sum over $k$. Chaining to a
fixed base point and using its bounded variance proves
\begin{equation}\label{chII:eq:fixed-core-sup-tail}
 \PP\{\sup_{D'}|H_D^{\mathrm{zm}}|>C(1+t)\}
 \le Ce^{-ct^2},\qquad t\ge0.
\end{equation}
For every fixed exponential parameter, the corresponding moment of
this supremum is bounded uniformly over $N$ and the finite interpolated
patch family.

For a prescribed subset $A\subset D'$, condition on
$H_D^{\mathrm{zm}}$. Since $q_D\ge0$, adding $q_D\xi_D$ on
the independent event $\xi_D\ge0$ can only increase every
coordinate. Therefore, for every deterministic $u$,
\[
 \PP\{\max_A\eta_D^0\ge u\}\le2\PP\{\max_A\eta_D\ge u\}.
\]
Set $h=\max(0,\sup_{D'}H_D^{\mathrm{zm}})+C_0$, where $C_0$
absorbs the bounded centering difference between the patch side
and $N$. Apply \eqref{chII:eq:dirichlet-subset-cited} with upward
parameter $z\ge1$ and downward parameter $y+h\ge0$. Conditional
integration and \eqref{chII:eq:fixed-core-sup-tail} give
\[
 \PP\{\max_A\phi_N\ge m_N+z-y\}
 \le C\sqrt{|A|/N^2}\,z e^{-\alpha(z-y)}\EE e^{\alpha h}
 \le C'\sqrt{|A|/N^2}\,z e^{-\alpha(z-y)}.
\]
Partition an arbitrary torus subset among the finitely many
cores and use Cauchy--Schwarz on the sum of their square-root
cardinalities. This proves \eqref{chII:eq:periodic-subset-input}
with its full prescribed-subset factor and height dependence.
The ordinary reference-square input is precisely
\eqref{chII:eq:dirichlet-subset-cited}.

For separation use the same coupling, now on
$\{\sup_{D'}|H_D^{\mathrm{zm}}|+\|q_D\|_\infty|\xi_D|\le A\}$.
Its complementary probability tends to zero as $A\to\infty$,
uniformly in $N$. Two periodic heights at least $m_N-H$ in
the same core then give two ordinary Dirichlet heights at
least $m_M-H-A-C_0$, where $M$ is the patch side. On the
embedded patch Euclidean and periodic distances agree up to
fixed constants. Pairs at distances between $R$ and $N/R$
therefore fall within the separation range in
\eqref{chII:eq:dirichlet-separation-cited} after changing $R$ by a
fixed factor. For $R$ sufficiently large every such pair is
in a common core by construction. For fixed $A$, apply the cited estimate uniformly over the finite
patch family; after the spatial limits, let $A\to\infty$. We obtain
\begin{equation}\label{chII:eq:periodic-separation-proved}
 \lim_{R\to\infty}\limsup_{N\to\infty}
 \PP\{\exists x,y:\phi_N(x),\phi_N(y)\ge m_N-H,
                         R<d(x,y)<N/R\}=0.
\end{equation}
The argument uses the exact rank-one correction from mean-zero
conditioning together with the finite-patch Dirichlet separation
estimate.

\subsection{From actual events to conditional error sums}

Recall $\mathcal H=\sigma(H)$. If $E_i$ is an
event determined by the $i$th fine field and $\mathcal H$, put
$q_i=\PP(E_i\mid\mathcal H)$ and $R=\sum_iq_i$.
Conditional independence of the ordinary fine fields gives
\[
 \PP\Bigl(\bigcup_iE_i\,\Bigm|\,\mathcal H\Bigr)
 =1-\prod_i(1-q_i)\ge1-e^{-R}
 \ge(1-e^{-1})\min\{1,R\}.
\]
The last inequality follows by concavity on $[0,1]$ and
monotonicity on $[1,\infty)$. Thus
\begin{equation}\label{chII:eq:actual-union-budget}
 \EE\min\{1,\textstyle\sum_iq_i\}
 \le (1-e^{-1})^{-1}\PP(\textstyle\bigcup_iE_i).
\end{equation}
This estimates precisely the truncated error in the product
formula.

Write $\bar\eta$ for the fine-field mean correction in the exact
decomposition of \S2.1. On $|\bar\eta|<1$, a fine-plus-coarse
point above $m_N+h-1/2$, with $|h|\le H$, is an actual periodic
point above $m_N-(H+2)$. Consequently the collar-error union is
contained in
\[
 \{\exists x\text{ in the collar}:\phi_N(x)\ge m_N-(H+2)\}
 \ \cup\ \{|\bar\eta|\ge1\}.
\]
The fine-mean variance estimate and Chebyshev's inequality control
the second term. For the first, the periodic subset estimate
\begin{equation}\label{chII:eq:periodic-subset-input}
 \PP\{\max_{x\in A}\phi_N(x)\ge m_N+z-y\}
 \le C\sqrt{|A|/N^2}\,z e^{-\alpha(z-y)},
 \quad z\ge1,\quad y\ge0,
\end{equation}
with $z=1$, $y=H+3$, gives
$C_H\sqrt{|A|/N^2}$. A collar of relative width $\delta$
has $|A|/N^2\le C\delta+o_N(1)$ at fixed block number.
Equation~\eqref{chII:eq:actual-union-budget} therefore proves the
conditional collar budget after $N\to\infty$ and then
$\delta\downarrow0$, with the fine-mean error removed in the
prescribed block limit.

Similarly, the union of core peaks above an upper cutoff $T$
is contained in the actual maximum event at height
$m_N+T-H-1$, together with $\{|\bar\eta|\ge1\}$.
Taking $A=V_N$ in \eqref{chII:eq:periodic-subset-input}, with
$z=T-H$ and $y=1$, bounds its probability by
$C(T-H)e^{-\alpha(T-H-1)}$ when $T\ge H+1$.
The displayed right-hand side tends to zero as $T\to\infty$.
Applying \eqref{chII:eq:actual-union-budget} gives the
upper-peak conditional error budget.

For the window error, suppose the retained window contains
$\{z\in\mathbb Z^2:|z|\le r\}$, the closed Euclidean lattice ball,
and $K\ge\max(2,2r)$. A true core contact
outside that window and its core maximum produce two actual
periodic points above $m_N-(H+2)$, unless $|\bar\eta|\ge1$.
Their periodic Euclidean distance is greater than $r$: within one
block the lifted displacement realizes this distance, and a displacement
of norm at most $r$ belongs to the retained window. Their distance is
less than $2N/K\le N/r$ because both belong to the same block.
Hence the entire error union is
contained in
\begin{equation}\label{chII:eq:window-pair-inclusion}
 \{\exists x,y:\ \phi_N(x),\phi_N(y)\ge m_N-(H+2),\
             r<d(x,y)<N/r\}\ \cup\ \{|\bar\eta|\ge1\}.
\end{equation}
The periodic two-scale separation estimate, followed by
\eqref{chII:eq:actual-union-budget}, proves the required small
conditional error. For each fixed $r$, take $K\ge2r$ and then the $N$-limit; the
resulting error then tends to zero as $r\to\infty$.

To deduce count tightness from \eqref{chII:eq:periodic-subset-input},
let $X$ be the
number of points above $m_N-H$ in a finite set of $n=N^2$
vertices. For an integer $1\le B\le n$, choose an independent
uniform subset $A$ of size $\lceil n/B\rceil$. Conditional
on $X\ge B$, the probability that $A$ misses all high points
is at most $(1-B/n)^{\lceil n/B\rceil}\le e^{-1}$. Averaging the deterministic
subset bound over $A$ gives
\[
 (1-e^{-1})\PP(X\ge B)
 \le C_H\sqrt{\lceil n/B\rceil/n}\le C_H\sqrt{2/B}.
\]
For $B>n$ the assertion is automatic. This gives the count cutoff used in the height-perturbation
estimate.

The subset estimate follows from the cited Dirichlet bound and
the constrained periodic transfer above. The pair input is
\eqref{chII:eq:periodic-separation-proved}. The subset term uses the periodic transfer above, and the pair term
uses \eqref{chII:eq:periodic-separation-proved}.

\subsection{The block estimates needed before multiplication}

For reference we state the complete block comparison exit. For
$H_0<\infty$, $0<\lambda\le1$ and $e>0$, there is $0<\delta_0<1/4$
such that, for $0<\delta\le\delta_0$, all sufficiently large
$K$, followed by sufficiently large $N$, satisfy
\begin{equation}\label{chII:eq:block-mass-integral}
 \int_{-H_0}^{H_0}\left|
 \EE\prod_i(1-\lambda q_i(h,H))
 -\EE e^{-\lambda\Gamma e^{-\alpha h}B_{N,K,2\delta}}
                        \right|dh<e.
\end{equation}
The probabilities $q_i(h,H)$ use cores with collar width $\delta$;
the mass uses $2\delta$. In addition,
\begin{equation}\label{chII:eq:canonical-block-error}
 \sup_{|h|\le H_0,\,0\le\lambda\le1}
 \left|\EE \Psi_{N,\lambda}((\phi_N+\overline\eta-m_N-h)^2)
       -\EE\prod_i(1-\lambda q_i(h,H))\right|<e.
\end{equation}

The block comparison combines the collar, upper-peak, nonlocal-gap,
weighted-profile and parity estimates from the preceding subsections
with the reference-cap asymptotic of Appendix~E.

Write $m=\lfloor N/K\rfloor$, $H_i=H_{N,K,i}$, and
\[
 z_i=(\lfloor m/2\rfloor/m,\lfloor m/2\rfloor/m).
\]
Thus $z_i$ is the deterministic lattice center in unit square
coordinates, including when $m$ is odd. The symbol $m_m$ means
the centering in \eqref{chII:eq:extremal-center} evaluated at side $m$.
On the good coarse-profile event write the effective fine height
as $t_i=m_N-m_m+h-H_i(z_i)$ and the
remaining profile as the clamped extension $d_i$ defined above;
on the retained core it is exactly $H_i(z)-H_i(z_i)$. Choose the conditional-rate cutoff from its tightness bound and
then the profile cutoff from \eqref{chII:eq:bad-deletion-final}. The reference estimate is then
applied only to a bounded H\"older class.

The lower effective-height restriction has a separate Gaussian
budget. At the block centers the covariance formula above gives,
for a fixed absolute $C$ and sufficiently large $m$,
\[
 \operatorname{Var}(H_i(z_i))\le v_K:=g\log K+C,\qquad
 m_N-m_m\longrightarrow2\sqrt g\log K\quad(N\to\infty,\ K\text{ fixed}).
\]
For a fixed threshold $t_0$, set
$d_K=2\sqrt g\log K-H_0-t_0-1$. Once $d_K>0$, Mills' bound and a union
bound give
\[
 \limsup_{N\to\infty}\PP\{\exists i,\ \exists h\in[-H_0,H_0]:t_i(h)<t_0\}
 \le\frac{K^2\sqrt{v_K}}{d_K\sqrt{2\pi}}
                       e^{-d_K^2/(2v_K)}
 =O_{H_0,t_0}((\log K)^{-1/2}).
\]
Indeed $2\log K-d_K^2/(2v_K)$ stays bounded. After $K$ is fixed,
choose $D_K$ so that $K^2v_K/D_K^2$ is smaller than the remaining
error budget. Chebyshev's inequality then bounds the probability of
$\max_i|H_i(z_i)|\ge D_K$, and
\[
 T_K=\max\{t_0,2\sqrt g\log K+1+H_0+D_K\}
\]
is a deterministic upper bound for all effective heights on the
remaining event. Thus, for fixed $t_0$ and $K$, $T_K$ is a deterministic upper
height bound; the compact-height estimate then determines the required
side-length threshold.

Apply \eqref{chII:eq:contact-repair-explicit} to each retained block,
then \eqref{chII:eq:product-error} to the conditional products.
The summed penalty errors are bounded by
\eqref{chII:eq:unpenalized-rate-budget} and its adjacent-side version below;
the bad-profile errors satisfy \eqref{chII:eq:bad-deletion-final}.
For the window and height cutoffs, monotone convergence gives
$J_{W,T}\to J_\infty$ as $W\uparrow\mathbb Z^2$ and $T\uparrow\infty$,
with the integrable envelope $e^{-\alpha u}$.

For the strip between $Q_\delta$ and $Q_{2\delta}$, apply
\eqref{chII:eq:actual-union-budget} and
\eqref{chII:eq:periodic-subset-input} to its contact events.
For a deterministic retained profile $d$, put
\[
 \begin{gathered}
 c_\lambda=\lambda c_*J_\infty,\qquad
 I(d)=\int_{Q_{2\delta}}e^{\alpha d(z)}w(z)\,dz,\\
 \ell(t,d)=-\log\EE\exp\left\{-c_\lambda e^{-\alpha t}
                  \int_{Q_{2\delta}}e^{\alpha d}\,d\mathcal M\right\}.
 \end{gathered}
\]
For $\|d\|_\infty\le B$, $I(d)\ge e^{-\alpha B}\int_{Q_{2\delta}}w>0$.
Equation~\eqref{chII:eq:critical-laplace} implies, uniformly in these profiles,
\[
 \frac{\ell(t,d)}{te^{-\alpha t}}\longrightarrow c_\lambda I(d),
 \qquad
 \log\frac1{c_\lambda e^{-\alpha t}}=\alpha t-\log c_\lambda.
\]
Given $\varepsilon,\tau>0$, choose $W,T$ so that $J_{W,T}$
approximates $J_\infty$, then choose a fixed $t_0$ and the side-length
threshold for the compact interval $[t_0,T_K]$. Together with
\eqref{chII:eq:reference-contact}, this gives
\[
 |\lambda p_m(2\delta,t,T,W,d)-\ell(t,d)|\le\varepsilon\ell(t,d),
 \qquad \lambda p_m(2\delta,t,T,W,d)\le\tau.
\]
Before choosing the penalty, fix the bound $B$ on the original good
profiles and choose a buffered core $Q_{\delta_*}$ with
$2\delta<\delta_*<1/2$. The reference estimate with $T=1$, $W=\{0\}$
gives a constant $c=c(\delta_*,B)>0$ such that
\[
 cte^{-\alpha t}\le p_m(\delta_*,t,1,\{0\},d)
 \qquad(\|d\|_\infty+[d]_\rho\le B),
\]
for sufficiently large fixed $t$, uniformly on subsequent compact
height intervals and sufficiently large even sides. For example one
may take $c=\tfrac12c_*J_{\{0\},1}e^{-\alpha B}\int_{Q_{\delta_*}}w$.
This $c$ depends on the original profile class and is independent of
$D,\zeta$ and $B_{\rm pen}$.
For $|h|\le H_0$ put
\[
 R_*:=\sum_i\PP\{\max_{x\in D_i}(\eta_i(x)+H(x))
                    >m_N-H_0-1\mid\mathcal H\}.
\]
The periodic high-point bound and \eqref{chII:eq:conditional-rate-tail}
give a deterministic cutoff $R_0$ for $R_*$ with arbitrarily small
exceptional probability. Every reference window contains a band point,
so the preceding lower bound, summed over the retained indices, gives
\begin{equation}\label{chII:eq:unpenalized-rate-budget}
 \sum_{i\in\mathcal I}t_i e^{-\alpha t_i}\le R_*/c.
\end{equation}
For odd $m$, use the adjacent even side and write
$\widehat t_i=t_i+m_m-m_{m+1}$. On the buffered core, the coupling
from the preceding subsection and the profile displacement imply
the same event inclusion except on a coupling-failure event.
Define
\[
 A_{N,K}:=\sum_i\PP\left\{
  \max_{\substack{x\in D_m\\x/(m+1)\in Q_{\delta_*}}}
       |\eta_{i,m+1}(x)-\eta_{i,m}(x)|>1/4\right\}.
\]
For large $m$ the profile discrepancy is at most $1/4$.
Then $A_{N,K}\to0$ at fixed $K$, and
\[
 \sum_{i\in\mathcal I}\widehat t_i e^{-\alpha\widehat t_i}
       \le(R_*+A_{N,K})/c.
\]
On $R_*\le R_0$ and $A_{N,K}\le1$, choose a local tolerance
$0<\eta\le\varepsilon c/(R_0+1)$, and only then choose $D,\zeta$
in \eqref{chII:eq:outer-penalty-budget} to make its coefficient at most
$\eta$. The total penalty error is at most $\eta(R_0+1)/c\le\varepsilon$
for either parity. Let $B_{\rm pen}$ be the resulting profile bound,
and increase the height and side-length thresholds for that class.
For each $i$, define $d_i^-=d_i-(B+1+D)\vartheta$, with $\vartheta$
from the penalty construction at core width $2\delta$. Since
$\vartheta=0$ on $Q_{2\delta}$,
\[
 d_i^-|_{Q_{2\delta}}=d_i|_{Q_{2\delta}},\qquad
 I(d_i^-)=I(d_i),\qquad \ell(t_i,d_i^-)=\ell(t_i,d_i).
\]
Apply the reference and Laplace estimates to this bounded penalized
class, with the window and height cutoffs chosen for the same tolerance.
For $0\le u\le\tau<1$,
$0\le-\log(1-u)-u\le\tau u/(1-\tau)$.
On the retained block set $\mathcal I$, write
$u_i=\lambda p_m(2\delta,t_i,T,W,d_i^-)$.
The same formulas use $m+1$ and $\widehat t_i$ when $m$ is odd.
Then
\[
 \begin{aligned}
 &\left|\sum_{i\in\mathcal I}-\log(1-u_i)-\sum_{i\in\mathcal I}\ell(t_i,d_i)\right|\\
 &\qquad\le\left(\varepsilon+\frac{\tau(1+\varepsilon)}{1-\tau}\right)
       \sum_{i\in\mathcal I}\ell(t_i,d_i).
 \end{aligned}
\]
If $\varepsilon+\tau(1+\varepsilon)/(1-\tau)\le1/2$, the two
exponential products differ by at most
$2\{\varepsilon+\tau(1+\varepsilon)/(1-\tau)\}$,
independently of the number of blocks. The penalty sum is at most
$\varepsilon$ by \eqref{chII:eq:unpenalized-rate-budget} and its
adjacent-side version. The bad-profile error is bounded by
\eqref{chII:eq:bad-deletion-final}; the upper-peak and collar errors
by \eqref{chII:eq:actual-union-budget} and
\eqref{chII:eq:periodic-subset-input}. For the nonlocal-gap error,
use \eqref{chII:eq:window-pair-inclusion},
\eqref{chII:eq:periodic-separation-proved} and
\eqref{chII:eq:actual-union-budget}.

The deterministic centering factor is
\[
 \gamma_{N,K}=K^4e^{-\alpha(m_N-m_m)}
 =\left(\frac{Km}{N}\right)^4
       \left(\frac{\log N}{\log m}\right)^{3/2}
 \longrightarrow1\quad(N\to\infty,\ K\text{ fixed}).
\]
Thus $c_\lambda e^{-\alpha t_i}
=c_\lambda\gamma_{N,K}K^{-4}e^{-\alpha h}e^{\alpha H_i(z_i)}$.
Conditional independence of the copied fine measures gives exactly
\[
 \prod_{i\in\mathcal I}e^{-\ell(t_i,d_i)}
 =\EE\left[\exp\left\{-\lambda\Gamma\gamma_{N,K}e^{-\alpha h}
       \frac{K^{-4}}{p_{\rm per}}
       \sum_{i\in\mathcal I}\int_{Q_{2\delta}}e^{\alpha H_i}\,d\mathcal M_i
                  \right\}\,\middle|\,H\right].
\]
The bound
$\sup_{x\ge0}|e^{-\gamma x}-e^{-x}|\le |\log\gamma|/\mathrm e$
for $\gamma>0$ removes $\gamma_{N,K}$; the mass-deletion estimate
then restores all blocks. The adjacent even-square comparison gives the same integrated estimate
for odd $m$. The preceding estimates determine the thresholds successively from
the compact height window and tolerance through the collar, rate, mass,
profile and penalty cutoffs, and finally through $K$, its upper-height
bound and $N$.

To see the geometry behind \eqref{chII:eq:canonical-block-error},
exclude the following events: $|\overline\eta|\ge1$; a point
above $m_N-H_0-2$ in the deleted collar; and two such high points
at distance between $a$ and $N/a$. The high-point collar and
two-scale separation estimates make the last two probabilities
small by taking the collar narrow and $a$ large, then $K,N$
large. On the complement each relevant cluster lies within one
core, has Euclidean diameter at most $a$ and hence graph diameter
at most $2a$, while distinct clusters have graph distance at least
$N/a$. Take $K>2a$. For fixed $a,K$ and large $N$,
$2a<\rho_N<R_N<N/a$ and each core has graph diameter less
than $N/a$. Thus each contacting core contains exactly one
canonical component, and none is deleted. The counts agree
simultaneously for all $|h|\le H_0$. Their bounded transforms
therefore differ in mean by at most the excluded probability.

With \eqref{chII:eq:block-mass-integral} and
\eqref{chII:eq:canonical-block-error} available, the proof of
Proposition~\ref{chII:prop:height-average} is complete: multiply by
a bounded compactly supported test, use scalar gluing, and
remove the mean by \eqref{chII:eq:dependent-translation}.

\appendixchapter{H}{Proof of Theorem~\ref{chII:thm:shell-source}}

The shell replacement follows from Gaussian energy concentration,
a density bound near the shell edge, and the algebraic estimates
below.

\subsection{Energy concentration and a density bound}

Let $(\lambda_{N,j})_{j=1}^{n-1}$ be the nonzero eigenvalues of
$-\mathcal L_N$, counted with multiplicity; these are the values
$\lambda_N(k)$ for $k\in\mathcal K_N\setminus\{0\}$.
Diagonalize the finite covariance on the real mean-zero Fourier
space. With independent standard normals $\zeta_j$,
\begin{equation}\label{chII:eq:energy-diagonal}
 \mathcal{Q}_N=\sum_{j=1}^{n-1}a_{N,j}\zeta_j^2,
 \qquad a_{N,j}=\frac1{n\lambda_{N,j}},\qquad
 \sum_ja_{N,j}=v_N.
\end{equation}
The eigenvalue bounds of Chapter~I imply
$\sup_{N,j}a_{N,j}<\infty$ and
$\sup_N\sum_ja_{N,j}^2<\infty$. It follows that
\begin{equation}\label{chII:eq:energy-moments}
 \operatorname{Var}(\mathcal{Q}_N)=2\sum_ja_{N,j}^2\le C,
 \qquad \EE e^{c(\mathcal{Q}_N-v_N)}\le C
\end{equation}
for sufficiently small fixed $c>0$. For the second bound use
\[
 \log\EE e^{c(\mathcal{Q}_N-v_N)}
 =\sum_j\left[-ca_{N,j}-\tfrac12\log(1-2ca_{N,j})\right]
 \le Cc^2\sum_ja_{N,j}^2.
\]
In particular $\EE|\mathcal{Q}_N-v_N|\le C$.

The energy has a density $f_N$ satisfying
\begin{equation}\label{chII:eq:energy-density}
 f_N(q)\le C e^{-c(q-v_N)},\qquad q\ge0.
\end{equation}
One direct verification uses exponential tilting in
\eqref{chII:eq:energy-diagonal}. Under the tilted law, each coefficient
becomes $a_{N,j}/(1-2ca_{N,j})$. Select two real Fourier
coordinates of the same lowest nonzero frequency. Their
coefficients are bounded above and below by positive constants,
uniformly for large $N$; their sum of squares therefore has a
bounded exponential density. Convolution with the remaining
nonnegative energy preserves that density bound. Undoing the
tilt and using \eqref{chII:eq:energy-moments} proves
\eqref{chII:eq:energy-density}. After increasing $C$ if necessary, the bound holds for all
$N\ge3$.

\subsection{Normalization and a bounded-test replacement}

For the deterministic identities in this subsection, let $\phi$ be
a mean-zero vector, $a>0$, and abbreviate
$\mathcal Q_N=n^{-1}\sum_x\phi(x)^2$ and
$w_N^\pm=w_N^\pm(\phi,a)$. Then
\begin{equation}\label{chII:eq:two-branch-algebra}
 w_N^+(\phi,a)+w_N^-(\phi,a)
 =\frac1n+\frac{\sum_x(|\phi(x)|-a)_+}{n^2a}
 \le\frac1n\left(1+\frac{\sqrt{\mathcal{Q}_N}}a\right).
\end{equation}
The equality follows by summing
$(x+a)_++(x-a)_+=x+a+(|x|-a)_+$; the inequality is
Cauchy--Schwarz. Also
\begin{equation}\label{chII:eq:one-branch-algebra}
 nw_N^+-1=\frac{\sum_x(-\phi(x)-a)_+}{na},\qquad
 0\le\max\{nw_N^-,nw_N^+-1\}\le\frac{\sqrt{\mathcal{Q}_N}}a.
\end{equation}
Here the first identity uses
$(x+a)_+=x+a+(-x-a)_+$ and $\sum_x\phi(x)=0$;
each of the two tail sums is bounded by $\sum_x|\phi(x)|$,
which gives the displayed inequalities.

Now set $\phi=\phi_N$, so $\mathcal Q_N$ is again the random energy
of Section~2.6. Fix $|\theta|\le H_0$ and abbreviate $B=B_N(\theta)$.
Uniformly on this interval $B\asymp L^2$ and $v_N=O(L)$.
On $\mathcal{Q}_N\le B/2$, $a=\sqrt{B-\mathcal{Q}_N}\ge\sqrt{B/2}$, and
\[
 \EE\left[\ind_{\{\mathcal{Q}_N\le B/2\}}\frac{\sqrt{\mathcal{Q}_N}}a\right]
 \le\sqrt{2v_N/B}=O(L^{-1/2}).
\]
Near the edge, \eqref{chII:eq:energy-density} gives
\begin{align*}
 \EE\left[\ind_{\{B/2<\mathcal{Q}_N<B\}}
            \left(1+\frac{\sqrt{\mathcal{Q}_N}}{\sqrt{B-\mathcal{Q}_N}}\right)\right]
 &\le C e^{-c(B/2-v_N)}
       \int_{B/2}^B\left(1+\frac{\sqrt B}{\sqrt{B-q}}\right)dq\\
 &\le CB e^{-c(B/2-v_N)}=o(1).
\end{align*}
Moreover $\PP(\mathcal{Q}_N\ge B)=o(1)$. Consequently
\begin{equation}\label{chII:eq:shell-normalizer}
 n\EE\left[\ind_{\{\mathcal{Q}_N<B\}}(w_N^++w_N^-)\right]\longrightarrow1
\end{equation}
uniformly on compact $\theta$ intervals, and these normalizers
are eventually at most $2$.

The same estimates work for any measurable $[0,1]$-valued
tests $T^+,T^-$ on the field space, with arbitrary extensions to
$\{\mathcal Q_N\ge B\}$. More precisely, \eqref{chII:eq:one-branch-algebra} gives
\[
 \begin{split}
 &\left|n\EE[\ind_{\{\mathcal Q_N<B\}}(w_N^+T^++w_N^-T^-)]
                        -\EE T^+\right|\\
 &\quad\le\PP(\mathcal Q_N\ge B)
       +2\EE\left[\ind_{\{\mathcal Q_N<B\}}
                    \frac{\sqrt{\mathcal Q_N}}{\sqrt{B-\mathcal Q_N}}\right]
       =O_{H_0}(L^{-1/2}).
 \end{split}
\]
The bound is uniform over the tests and $|\theta|\le H_0$; their
dependence on the random amplitude causes no change in the estimate.
For the square-band tests, sign symmetry of $\phi_N$ and the
fact that $a$ depends only on its energy give
$\EE T^+=\EE T^-$. The negative-shift test is exactly
$\Psi_{N,\lambda}((\phi_N-a)^2)$. Thus
\begin{equation}\label{chII:eq:shell-reference-error}
 S_N(\theta,\lambda)
 -\EE \Psi_{N,\lambda}((\phi_N-\alpha_N(\theta,\phi_N))^2)\longrightarrow0.
\end{equation}
The compact normalizer bound also proves $0\le S_N\le2$.

\subsection{Removal of the energy shift}

Since $B_N(\theta)$ is affine with slope $2gL$,
\[
 B_N(\theta)-\mathcal{Q}_N=B_N\left(\theta-\frac{\mathcal{Q}_N-v_N}{2gL}\right)-v_N.
\]
Taking the common positive-part square-root convention gives the
corresponding exact amplitude identity. For $\chi\in C_c(\mathbb R)$,
\eqref{chII:eq:dependent-translation} bounds the integrated
replacement error by
$\EE\|\chi(\cdot+R_N^E)-\chi\|_1$. For every $a>0$ this is at most
\[
 \sup_{|r|\le a}\|\chi(\cdot+r)-\chi\|_1
       +2\|\chi\|_1\frac{\EE|R_N^E|}{a}.
\]
Consequently
$\lim_{a\downarrow0}\limsup_{N\to\infty}
\EE\|\chi(\cdot+R_N^E)-\chi\|_1=0$ by
$\EE|R_N^E|=O(L^{-1})$ and translation continuity in $L^1$.

For the deterministic change of variable, if
$h=a_N^0(\theta)-m_N$, then on every relevant compact interval
\[
 \theta=\frac{2m_Nh+h^2-v_N}{2gL},\qquad
 \frac{d\theta}{dh}=\frac{m_N+h}{gL}.
\]
These converge uniformly to $\alpha h-1/2$ and $\alpha$,
respectively. The transformed compact supports remain inside
one fixed compact interval. Uniform continuity of $\chi$ and
the bounded Jacobians give convergence of the transformed
weights in $L^1$. This proves the final change of variables in
Theorem~\ref{chII:thm:shell-source}.


\par\endgroup

\begingroup
\mergedchaptersetup{III}{Local times and resampling}
\renewcommand{\AppendixEquationPrefix}{A}
\appendixchapter{I}{Proof of the square-field identity \eqref{chIII:eq:RK}}

This appendix proves \eqref{chIII:eq:RK}. All matrices and vectors in the
appendix are indexed by $r^c:=V_N\setminus\{r\}$ unless stated otherwise.
Fix $r\in V_N$, put $v_N=G_N(r,r)$ and $h=h^r$, and take $h$
independent of the walk under $\PP_r$ when both appear in an expectation.
The symbols $B,C,U$ below are local matrix notation. The global
quantities $A_N,C_N$ retain their meaning.

\subsection{The pinned covariance}

Let $P=I+\cL_N$ be the transition matrix on $V_N$ of the embedded
discrete-time walk. Put $B=(-\cL_N)_{r^c,r^c}$. It is positive definite: extend a vector
by zero at $r$ and use the Dirichlet form. Zero energy forces a
constant vector on the connected torus, hence the zero vector.
Its inverse is the Green matrix of the walk killed at $r$.

For $x,y\ne r$ define
\begin{equation}\label{chIII:eq:pinned-covariance}
 C(x,y)=G_N(x,y)-G_N(x,r)-G_N(r,y)+v_N.
\end{equation}
Extend the right-hand side by zero at $x=r$. The identity
$-\cL_NG_N(\cdot,y)=\ind_y-n^{-1}$ gives
$-\cL_NC(\cdot,y)=\ind_y-\ind_r$. Restriction to $r^c$ proves
$BC=I$. Thus $C=B^{-1}$ and $h^r|_{r^c}$ is a centered Gaussian
vector with covariance $C$.

Let $U$ be diagonal with entries $u_x\ge0$. The vector
\begin{equation}\label{chIII:eq:response}
 f=(B+U)^{-1}B\mathbf1=(I+CU)^{-1}\mathbf1
\end{equation}
solves $(B+U)f=B\mathbf1$. Probabilistically,
\begin{equation}\label{chIII:eq:response-probability}
 f(x)=\EE_x\exp\left\{-\int_0^{H_r}u_{X_t}\,dt\right\},
                  \qquad x\ne r.
\end{equation}
To check this, condition on the first holding time and jump. The
holding time is exponential of rate one, so
$(1+u_x)f(x)=P(x,r)+\sum_{y\ne r}P(x,y)f(y)$.
This is precisely the displayed linear system; uniqueness follows
from positive definiteness of $B+U$. In particular $0\le f\le1$.

\subsection{Root excursions and their Laplace transform}

The successive holding times at $r$ are independent exponentials
of rate one. Measured in root local time, the departure points
therefore form a Poisson process of rate one. At each departure the
walk jumps according to $P(r,\cdot)$ and then runs until its next
return to $r$. The strong Markov property makes these excursions
independent marks, independent also of the holding-time clock.
Every excursion returns almost surely on the finite irreducible
graph. These assertions also follow directly by multiplying the
densities of the holding times and applying the strong Markov
property at successive returns.

Let $M_\lambda$ count departures whose root-local-time coordinate
lies in $[0,\lambda)$. Then $M_\lambda$ is Poisson with mean $\lambda$. The
Laplace transform of the occupation in one excursion is
$\sum_{x\ne r}P(r,x)f(x)$. Consequently
\begin{equation}\label{chIII:eq:excursion-laplace}
 \EE_r\exp\left\{-\sum_{x\ne r}u_x
                      \ell_{\tau_r(\lambda)}(x)\right\}
 =\exp\left\{-\lambda
            \left(1-\sum_{x\ne r}P(r,x)f(x)\right)\right\}.
\end{equation}
This follows from the Poisson identity
$\EE a^{M_\lambda}=e^{\lambda(a-1)}$, with $a\in[0,1]$.

Since $B\mathbf1=(P(x,r))_{x\ne r}$ and $B$ is symmetric,
\[
 1-\sum_{x\ne r}P(r,x)f(x)
 =\mathbf1^{\mathsf T}B(\mathbf1-f)
 =\mathbf1^{\mathsf T}Uf.
\]
Here $\mathbf1^{\mathsf T}B\mathbf1=\sum_{x\ne r}P(r,x)=1$ and
$B(\mathbf1-f)=Uf$ by \eqref{chIII:eq:response}. Define
\begin{equation}\label{chIII:eq:laplace-coefficient}
 q(U)=\mathbf1^{\mathsf T}U(I+CU)^{-1}\mathbf1.
\end{equation}
Since $I+CU=B^{-1}(B+U)$ and both $B$ and $B+U$ are invertible,
$I+CU$ is invertible. Moreover $q(U)=\sum_{x\ne r}u_x f(x)\ge0$.
Restoring a potential $u_r\ge0$ at the root multiplies
\eqref{chIII:eq:excursion-laplace} by $e^{-\lambda u_r}$, because
$\ell_{\tau_r(\lambda)}(r)=\lambda$. We obtain
\begin{equation}\label{chIII:eq:full-occupation-laplace}
 \EE_r e^{-\sum_xu_x\ell_{\tau_r(\lambda)}(x)}
                      =e^{-\lambda(u_r+q(U))}.
\end{equation}

\subsection{Gaussian integration and uniqueness}

If $\eta\sim N(0,C)$ and $a\in\mathbb R$, completing the square
in its ordinary finite-dimensional density yields
\begin{equation}\label{chIII:eq:gaussian-square}
 \EE e^{-\frac12(\eta+a\mathbf1)^{\mathsf T}U(\eta+a\mathbf1)}
 =\det(I+CU)^{-1/2}e^{-a^2q(U)/2}.
\end{equation}
Indeed, the new quadratic matrix is $B+U$ and the linear term is
$-aU\mathbf1$. Its integration gives exponent
\[
 -\frac{a^2}2\mathbf1^{\mathsf T}
     [U-U(B+U)^{-1}U]\mathbf1
       =-\frac{a^2}2q(U),
\]
because $(B+U)^{-1}U=I-(I+CU)^{-1}$. The determinant ratio is
$\det(B)^{1/2}/\det(B+U)^{1/2}$. For every positive semidefinite diagonal $U$, $B+U$ is positive
definite.

Apply \eqref{chIII:eq:gaussian-square} with $a=0$ to the independent
field on the occupation side. By \eqref{chIII:eq:full-occupation-laplace},
the joint Laplace transform of $2\ell_{\tau_r(\lambda)}+h^2$ is
\begin{equation}\label{chIII:eq:joint-laplace}
 \EE e^{-\frac12\sum_xu_x(2\ell_{\tau_r(\lambda)}(x)+h(x)^2)}
 =e^{-\lambda(u_r+q(U))}\det(I+CU)^{-1/2}.
\end{equation}
The root coordinate of $(h+\sqrt{2\lambda})^2$ is $2\lambda$.
Applying \eqref{chIII:eq:gaussian-square} with $a=\sqrt{2\lambda}$ to
the other coordinates gives precisely \eqref{chIII:eq:joint-laplace}.
Thus, for every $(u_x)_{x\in V_N}\in[0,\infty)^{V_N}$, the
transform on the left uses the product law $\PP_r\otimes\gamma_r$,
whereas that on the right uses $\gamma_r$ alone. This is an
equality of the occupation-profile laws.

These Laplace transforms determine the laws on $[0,\infty)^{V_N}$.
One direct verification maps $z_x$ to $e^{-z_x}$ in $[0,1]$.
Equality of Laplace transforms at nonnegative integer vectors gives
equality of every polynomial moment of the image laws. Polynomials
are dense in the continuous functions on the compact cube, so the
image laws agree. Both give mass one to $(0,1]^{V_N}$, on which
the map has the Borel inverse $z_x=-\log y_x$. Their original laws
therefore agree. Thus the two random vectors have the same law on the full
coordinate space. Equation~\eqref{chIII:eq:RK} follows as an equality of
image measures.

\appendixchapter{J}{Proof of Proposition~\ref{chIII:prop:record-product}}

We prove the conditional product law in
Proposition~\ref{chIII:prop:record-product} and its loop-time
decomposition used in Corollary~\ref{chIII:cor:counted-window}.
The latter retains the censored terminal segment in
\eqref{chIII:eq:censored-kernel}.

\subsection{Stopped driving sequences}

Let $\mathsf D$ be the product law of independent pairs
$(e_j,\delta_j)_{j\ge0}$,
where $e_j$ has exponential law of rate one and $\delta_j$ is uniform
on the four nearest-neighbor directions, independently of $e_j$.
Write $\EE_{\rm drv}$ for integration under $\mathsf D$.
Given the initial vertex, construct the walk from these pairs.
The walk holds for $e_j$ and then makes jump $\delta_j$. If its
first exit from $D$ occurs on jump $k$, a record consists of $k$
and these first $k$ pairs. Its duration is
$V=\sum_{j=0}^{k-1}e_j$ and its endpoint is the vertex after jump
$k$. The unused sequence starts with pair $k$. Denote the stopped
record by $\xi$ and the unused sequence by $\eta$.

For every proper domain $D$, the first-exit index is finite almost
surely: from each vertex there is a path to $D^c$ of at most $n$ jumps,
so the probability of no exit in $jn$ jumps is at most
$(1-4^{-n})^j$. Nonexplosion and positive holding times then make the
record map well defined almost surely. On each fixed exit-index event,
the record coordinates use finitely many jump variables and a finite
sum of holding times, hence are measurable.

For a bounded nonnegative test $f$ of the stopped record, a bounded
nonnegative test $g$ of the unused sequence, and a feasible endpoint
$b$, independence of the driving pairs gives
\begin{equation}\label{chIII:eq:one-record-factor}
 \EE_a[\ind_{\{X_{\tau_D}=b\}}f(\xi)g(\eta)]
   =h_{D,b}(a)\,\mathsf B_D^{a,b}(f)\,\EE_{\rm drv}[g(\eta)].
\end{equation}
Partition by the first-exit jump index $k$.
The event of first exit on jump $k$ through $b$, and the stopped
record, use only the first $k$ driving pairs. The unused sequence
has its original independent product law. Sum the resulting
identities over $k$. A monotone-class argument extends from cylinder
tests to all nonnegative measurable tests.

\subsection{Endpoint words and selected coordinates}

Iterate \eqref{chIII:eq:one-record-factor}. On a fixed feasible word of
$k\ge1$ successive endpoints, the next domain and starting vertex are
deterministic. For tests $f_1,\ldots,f_k$ of the records and $g$
of the final unused sequence, one obtains
\begin{equation}\label{chIII:eq:word-factor}
 \EE_{a_1}\left[\ind_{\{\text{word}=w\}}
                   \prod_{j=1}^kf_j(\xi_j)\,g(\eta)\right]
 =p_w\prod_{j=1}^k\mathsf B_{D_j}^{a_j,b_j}(f_j)\,\EE_{\rm drv}g,
 \quad p_w=\prod_{j=1}^kh_{D_j,b_j}(a_j).
\end{equation}
Taking every test equal to one proves the formula for the word
mass. Dividing by $p_w>0$ proves the product law. Tests that do not
factor follow by a monotone-class argument. Concatenation is the inverse of record extraction on the
full-measure record space. Therefore the factorization is an identity
for the reconstructed path law.

Under this product law, a deterministic selected index set $C$ gives
the measure decomposition
\[
 \left(\bigotimes_{j\in C}\mathsf B_{D_j}^{a_j,b_j}\right)
 \otimes\left[
       \left(\bigotimes_{j\notin C}\mathsf B_{D_j}^{a_j,b_j}\right)
                  \otimes\mathsf D\right].
\]
Given selected coordinates $x$ and exterior coordinates $y$, let
$\mathcal C(x,y)$ reinsert the records in their original index order.
Its completed duration is $T(\mathcal C(x,y))=V(x)+E(y)$. Restriction by a
certificate \eqref{chIII:eq:selection-certificate} yields
\eqref{chIII:eq:selected-law} by factoring its indicators. The fixed endpoint word and collapsed exterior determine the
product kernels; their duration coordinates are integrated under those
kernels. The physical-time restriction is inserted
through the weight $J_I$ below.

\subsection{Completed traversals and the censored future}

Fix nonempty $I_0\subset O$ and $r\notin O$. Let $T_0=0$ and define
successively $T_{2j+1}$ as the first entrance into $I_0$ after
$T_{2j}$, and $T_{2j+2}$ as the first exit from $O$ after
$T_{2j+1}$. Write $\xi_1,\ldots,\xi_{2m}$ for the completed
records up to $T_{2m}$, and $w$ for their endpoints. The domains
alternate between $V_N\setminus I_0$ and $O$, both proper.
All these stopping times are finite almost surely.

Except at the countable set of traversal endpoints, every $t$ for
which $X_t=r$ belongs to a unique interval
$(T_{2m},T_{2m+1})$, $m\ge0$. It cannot lie strictly between
$T_{2m+1}$ and $T_{2m+2}$ because the walk is then in $O$.
On the interval for count $m$, put $q=t-T_{2m}$. Concatenation gives
\[
 t=\sum_{j=1}^{2m}V(\xi_j)+q,
 \qquad 0<q<H_{I_0}(\eta),\qquad X_q(\eta)=r.
\]
Let $\mu_w$ be the product of the $m$ outward-record kernels and
let $\nu_w$ be the product of the $m$ inward-record kernels and the
independent unused driving-sequence law, all for this fixed word.
Translation of Lebesgue measure gives $dt=dq$. Applying
\eqref{chIII:eq:word-factor} to a nonnegative exterior test $F$ therefore
gives the exact conditional slice formula
\begin{equation}\label{chIII:eq:actual-slice}
 \begin{split}
 &\EE_r\left[\left.
  \int_{T_{2m}}^{T_{2m+1}}\ind_{\{t\in I,\ X_t=r\}}
                  F(y,t-T_{2m})\,dt\ \right|\ \text{word}=w\right]
 \\
 &=\int\mu_w(dx)\nu_w(dy)\int_0^\infty
    \ind_{\{E(y)+V(x)+q\in I\}}
    \ind_{\{q<H_{I_0}(\eta),\ X_q(\eta)=r\}}F(y,q)\,dq.
 \end{split}
\end{equation}
Insertion of any measurable test of the outward records multiplies
the integrand on the right by that test of $x$. This proves the
claimed identification with \eqref{chIII:eq:window-functional}.
In particular, with $F=1$ and positive finite slice mass, the
normalized law of $(x,y)$ in this loop slice is
\[
 \frac{J_I^{m,w}(x,y)}
      {\int J_I^{m,w}\,d(\mu_w\otimes\nu_w)}
                  \,\mu_w(dx)\nu_w(dy).
\]
Equivalently, the normalized loop slice has Radon--Nikodym density
proportional to $J_I^{m,w}$ with respect to the product law
$\mu_w\otimes\nu_w$. Apply the window comparison pointwise to $J_I^{m,w}(x,y)$ under
$\mu_w\otimes\nu_w$ and then integrate the resulting inequality.

For each fixed $m$, the possible endpoint words of the first $2m$
records partition the infinite path space; consequently
$\sum_wp_{m,w}=1$. Summing their slice formulas, and then summing
over $m$, recovers the rooted-loop measure
\eqref{chIII:eq:loop-measure}. For $m=0$ use the unique empty word with mass one, empty record
products and the original future from $r$; the sum over $m\ge0$ then
contains the zero-count term.

\appendixchapter{K}{Proof of Proposition~\ref{chIII:prop:duration}}

We prove the duration bounds in Proposition~\ref{chIII:prop:duration}
and the disk version \eqref{chIII:eq:disk-duration}--\eqref{chIII:eq:disk-conditional}.
Let $\varnothing\ne D\subsetneq V_N$,
put $k=|D|$ and $v_N=G_N(0,0)$, and retain $\tau_D$ and $h_{D,b}$
from Section~3.5. The matrices below are local to this appendix.

\subsection{A killed spectral bound for an arbitrary domain}

Choose $z\notin D$, and let $G_D=(-\cL_N)_{D,D}^{-1}$ be the
killed Green matrix. It is positive definite. The occupation
representation and $\tau_D\le H_z$ give
\begin{equation}\label{chIII:eq:killed-diagonal}
 \begin{split}
 G_D(x,x)&=\EE_x\int_0^{\tau_D}\ind_{\{X_t=x\}}\,dt\\
 &\le G_N(x,x)-G_N(x,z)-G_N(z,x)+G_N(z,z)\le4v_N.
 \end{split}
\end{equation}
The middle expression is the Green matrix killed at $z$, as proved
in \eqref{chIII:eq:pinned-covariance}. The last bound follows from
$|G_N(x,z)|\le v_N$, the covariance Cauchy--Schwarz inequality.

If $\lambda_D$ is the smallest eigenvalue of $(-\cL_N)_{D,D}$,
then
\begin{equation}\label{chIII:eq:domain-gap}
 \lambda_D^{-1}=\|G_D\|_{2\to2}
                 \le\operatorname{tr}G_D\le4v_Nk.
\end{equation}
The killed semigroup is the matrix exponential of
$-(-\cL_N)_{D,D}$: conditioning on the first holding time gives
its backward equation and zero boundary condition. Hence, in
counting-measure Euclidean norm,
\[
 \PP_a(\tau_D>d)
  =\langle\ind_a,e^{-d(-\cL_N)_{D,D}}\mathbf1\rangle
  \le\sqrt{k}\,e^{-d\lambda_D}
  \le k e^{-d/(4v_Nk)}.
\]
This proves \eqref{chIII:eq:unconditioned-duration} for every starting
vertex in $D$ and every $d\ge0$.

\subsection{The cost of the two conditionings}

For a feasible endpoint, the duration $\tau_D$ under
$\mathsf B_D^{a,b}$ satisfies
\[
 \mathsf B_D^{a,b}(\tau_D>d)
  =\frac{\PP_a(\tau_D>d,X_{\tau_D}=b)}{h_{D,b}(a)}
  \le\frac{k}{h_{D,b}(a)}e^{-d/(4v_Nk)}.
\]
Under the product law $\PP=\bigotimes_{i=1}^m\mathsf B_D^{a_i,b_i}$,
let $V_i$ be the duration of record $i$. Fix a nonempty $K\subset D$
and a measurable event $S$ containing simultaneous avoidance of $K$.
If record $i$ avoids $K$ with probability at least $p_i>0$, a union
bound followed by division by $\PP(S)$ gives the more precise
version of \eqref{chIII:eq:conditional-duration}:
\begin{equation}\label{chIII:eq:heterogeneous-duration}
 \PP(\max_iV_i>d\mid S)
 \le\frac{k\sum_{i=1}^m h_{D,b_i}(a_i)^{-1}}
                {\prod_{i=1}^m p_i}\,e^{-d/(4v_Nk)}.
\end{equation}
The probability of $S$ is positive under these assumptions.
Suppose that, for every $i=1,\ldots,m$, there is a simple path from
$a_i$ to $b_i$ staying in $D\setminus K$ until its last step and using
at most $k+1$ jumps. Following its prescribed directions has probability
at least $4^{-(k+1)}$. Integrating the unrestricted holding times
contributes a factor one.
For every $i$ this gives
\[
 h_{D,b_i}(a_i)\ge4^{-(k+1)},\qquad
 \mathsf B_D^{a_i,b_i}(\text{avoids }K)
 \ge\frac{4^{-(k+1)}}{h_{D,b_i}(a_i)}\ge4^{-(k+1)}.
\]
Substitution proves \eqref{chIII:eq:finite-duration}.
For $m=0$, sums of occupation and duration are zero and the
long-duration event is empty, so no division is necessary.

Sometimes the duration cutoff is imposed after soft conditioning.
Put $V=\sum_{i=1}^mV_i$, the total duration of the records under consideration.
For $\PP(S)>0$, use the conditional joint law
$\mu(A)=\PP(A\cap S)/\PP(S)$ of the records. If
$\mu(V>d)\le\varepsilon<1$, its short-duration restriction obeys
\begin{equation}\label{chIII:eq:duration-truncation}
 \mu(\,\cdot\mid V\le d)\le\frac{1}{1-\varepsilon}\mu
\end{equation}
as measures. The proof of Theorem~\ref{chIII:thm:one-success}
integrates over the successful short records directly and pays the
original long-duration error once.

\subsection{Disk exit times and the duration used in a window}

For the disk in \eqref{chIII:eq:disk-duration}, lift the walk to
$\mathbb Z^2$ until exit. The generator applied to $|x-y|^2$ is one.
Stopping its martingale at $t\wedge\tau_D$ yields
\[
 \EE_a(t\wedge\tau_D)
   =\EE_a|X_{t\wedge\tau_D}-y|^2-|a-y|^2
   \le(R+1)^2.
\]
Letting $t\uparrow\infty$ by monotone convergence proves the mean
bound, including finiteness of the exit time. Markov's inequality gives
$\sup_a\PP_a(\tau_D>2(R+1)^2)\le1/2$. At successive times $2j(R+1)^2$ the
Markov property gives $\sup_a\PP_a(\tau_D>2j(R+1)^2)\le2^{-j}$.
Taking $j=\lfloor d/[2(R+1)^2]\rfloor$ proves
\eqref{chIII:eq:disk-duration}. Endpoint restriction and simultaneous
avoidance then give \eqref{chIII:eq:disk-conditional}.

The window lemma uses the total duration of the selected records. If there are at most $M$ selected records and
each is at most $d_0$, their sum is at most $Md_0$. Thus
\begin{equation}\label{chIII:eq:duration-sum}
 \PP\left(\sum_{i=1}^mV_i>Md_0\right)
        \le\PP(\max_{i\le m}V_i>d_0),\qquad m\le M.
\end{equation}
The same inclusion holds under any of the conditional laws above.
This is why the enlargement is by $Md_0$ when a per-record tail
estimate at $d_0$ is supplied. With an empty selected family there
is no enlargement and no duration error.


\par\endgroup

\begingroup
\mergedchaptersetup{IV}{Geometry and tightness of low local times}
\let\mathscr\mathcal
\allowdisplaybreaks[2]
\raggedbottom
\appendixchapter{L}{Proof of Proposition~\ref{chIV:prop:local}: endpoint estimates}\label{chIV:sec:R6}

We prove the endpoint, duration and hitting estimates in
Proposition~\ref{chIV:prop:local}. Appendix~M.1 proves its
occupation estimates. The remaining subsections derive the
conditional-duration and loop-window bounds used in the proof of
Proposition~\ref{chIV:prop:count}.

Use the lattice potential kernel $a_{\rm lat}$ of Chapter~II and its
expansion $a_{\rm lat}(z)=g\log|z|+c_{\rm lat}+O(|z|^{-2})$, together with scale-uniform
interior Harnack and oscillation estimates for lattice disks. Here $g=2/\pi$ and the total jump rate is $1$.
\subsection{Uniform disk-kernel estimates}
In the planar coordinates below, $B(y,R)=\{z\in\mathbb R^2:|z-y|<R\}$;
the walk domains are its lattice points.
Let $q\to\infty$, $0<\gamma<1$, $w=q^\gamma$, and
\begin{equation*}
{r_{\rm in}}=e^q,\qquad \rho=e^{q+w},\qquad R=e^{q+2w},\qquad
{\mathsf K}=B(y,{r_{\rm in}})\cap\mathbb Z^2,\quad D=B(y,R)\cap\mathbb Z^2.
\end{equation*}
Take $a$ with $\bigl||a-y|-\rho\bigr|\le2$, and an exterior neighbor $b$ of $D$.
Use $\mathsf B_D^{a,b}$ for the
stopped-path kernel of \eqref{chIII:eq:local-kernel}, abbreviated
$P^{a,b}$ in the main text. The duration is sampled as part of the stopped record.
Here and below $G_D$ is the counting-measure Green kernel of the walk
killed at $\tau_D$.
There exist $c,C>0$, independent of $q,y,a,b$, such that, for all sufficiently large $q$,
\begin{equation}
\mathbb P_a(X_{\tau_D}=b)\ge cR^{-C},
\label{chIV:eq:R6.1}
\end{equation}
\begin{equation}
\mathsf B_D^{a,b}(\tau_D>t)
\le CR^C e^{-ct/R^2},\qquad t\ge0,
\label{chIV:eq:R6.2}
\end{equation}
\begin{equation}
\mathsf B_D^{a,b}(H_{\mathsf K}>\tau_D)\ge c,
\label{chIV:eq:R6.3}
\end{equation}
\begin{equation}
\inf_{x\in {\mathsf K}}\mathsf B_D^{a,b}(H_x<\tau_D)
\ge c\frac{w}{q}.
\label{chIV:eq:R6.4}
\end{equation}
For the first estimate, write $h_b(z)=h_{D,b}(z)=\PP_z(X_{\tau_D}=b)$.
An interior neighbor of $b$ has $h_b\ge1/4$. Following a bounded
number of lattice edges inward, and then retracing them to $b$,
gives a point at lattice distance at least $4$ from $D^c$ with
$h_b\ge c$.
A chain of $O(\log R)$ Harnack balls of geometrically increasing radii gives
$h_b(y)\ge cR^{-C}$; interior Harnack then applies to $|a-y|\le R/4$.

Optional stopping of $|X_t-y|^2-t$, first at bounded stopping times, gives
\begin{equation*}
\sup_{z\in D}\EE_z\tau_D\le(R+1)^2
\quad\Longrightarrow\quad
\sup_{z\in D}\PP_z(\tau_D>t)\le2e^{-ct/R^2}.
\end{equation*}
The implication follows by Markov's inequality and strong Markov iteration.
Division by \eqref{chIV:eq:R6.1} proves \eqref{chIV:eq:R6.2}.
Harnack oscillation decay also gives, for some fixed $\xi>0$,
\begin{equation}
\sup_{z,z'\in B(y,2\rho)}
\left|\frac{h_b(z)}{h_b(z')}-1\right|
\le C(\rho/R)^\xi\le Ce^{-\xi w}.
\label{chIV:eq:R6.5}
\end{equation}
Stopping the potential kernel at $H_{\mathsf K}\wedge\tau_D$, with radii
$r_{\rm in}+O(1)$ and $R+O(1)$, yields
\begin{equation*}
\mathbb P_a(H_{\mathsf K}<\tau_D)
=\frac{\log R-\log\rho+O({r_{\rm in}}^{-1})}
{\log R-\log {r_{\rm in}}+O({r_{\rm in}}^{-1})}
=\frac12+o(1).
\end{equation*}
Exit conditioning multiplies this probability by $1+O(e^{-\xi w})$, by the strong Markov property
and \eqref{chIV:eq:R6.5}; this proves \eqref{chIV:eq:R6.3}.
For $x\in\mathsf K$, the potential representation gives
\begin{equation*}
G_D(a,x)=gw+O(1),\qquad
G_D(x,x)=g(q+2w)+O(1).
\end{equation*}
uniformly, since $r_{\rm in}/\rho,r_{\rm in}/R\to0$.
Use $\PP_a(H_x<\tau_D)=G_D(a,x)/G_D(x,x)$ and \eqref{chIV:eq:R6.5} for \eqref{chIV:eq:R6.4}.
All estimates remain uniform under $o(1)$ perturbations of the three logarithmic radii.
Torus disks of radius below $N/4$ have the same planar kernels.

\subsection{Durations conditional on a soft event}
For $m$ independent blocks with arbitrary feasible endpoint pairs and durations $\delta_i$,
let $\ell_i=\ell^{(i)}$ be their local-time profiles as in Section~4.1,
and put $\ell=\sum_{i=1}^m\ell_i$. Write
$\widehat{\mathbb P}=\bigotimes_{i=1}^m\mathsf B_D^{a_i,b_i}$
and let $\widehat{\mathbb E}$ denote its expectation. For $u_0\ge0$, define
\begin{equation*}
S=\{\exists x\in {\mathsf K}:\ell(x)\le u_0\}.
\end{equation*}
Avoidance of the whole core by all blocks implies $S$. Thus \eqref{chIV:eq:R6.2}--\eqref{chIV:eq:R6.3} give
\begin{equation}
\widehat{\mathbb P}(S)\ge c^m,\qquad
\widehat{\mathbb P}(\max_{i\le m}\delta_i>d\mid S)
\le CmR^C c^{-m}e^{-cd/R^2}.
\label{chIV:eq:R6.6}
\end{equation}
For $m=0$ the exceptional probability is zero.
When this kernel is used inside a global path, all core occupation must be reconstructed by these blocks.
Take $\gamma=1/8192$, $Q=L^{16\gamma}$, and
\begin{equation*}
Q\le q\le L-Q,\qquad M=\lceil L^4\rceil,\qquad d=ne^{-Q}.
\end{equation*}
Then $R^2\le ne^{-2Q+4L^\gamma}$ and $d/R^2\ge e^{Q-4L^\gamma}$.
For $m\le M$, the logarithm of the bound in \eqref{chIV:eq:R6.6} is at most
\begin{equation*}
O(L+L^4)-ce^{Q-4L^\gamma}\le-e^{Q/2}
\end{equation*}
for sufficiently large $N$. Consequently, with or without conditioning on $S$,
\begin{equation*}
\begin{gathered}
\max\{\widehat{\mathbb P}(\max_{i\le m}\delta_i>d),\widehat{\mathbb P}(\max_{i\le m}\delta_i>d\mid S)\}
\\
 \le\varepsilon_N:=\exp\{-e^{Q/2}\},\qquad
N^aL^b\varepsilon_N\to0\quad(a,b<\infty).
\end{gathered}
\end{equation*}

\subsection{Excessive block counts at a soft vertex}
Exit conditioning depends only on the jump chain. After a block first hits $x\in\mathsf K$,
the first holding time at $x$ is therefore still an independent $\operatorname{Exp}(1)$ variable.
By \eqref{chIV:eq:R6.4},
\begin{equation*}
\widehat{\mathbb E}e^{-\ell_i(x)}
\le1-\frac12\widehat{\mathbb P}(H_x<\tau_D)
\le1-cw/q.
\end{equation*}
Exponential Markov applied to the first $M$ independent blocks gives
\begin{equation}
\widehat{\mathbb P}\left(\sum_{i=1}^{M}\ell_i(x)\le u_0\right)
\le e^{u_0-cMw/q}\le e^{u_0-cL^3}.
\label{chIV:eq:R6.7}
\end{equation}
The last estimate uses only $w\ge1$ and $q\le L$.
Subsequent occupation increments are nonnegative, so the same bound applies when the total count is at least $M$.

\subsection{Transfer to integrated loop windows}
Let the root $r$ lie outside $D$. For a loop with at most $M$
completed blocks, use the slices $m=0,\ldots,M$ of
\eqref{chIII:eq:actual-slice}. On each fixed count and endpoint word,
the selected law is the untruncated product of its outward kernels.
The terminal censoring in $J_I$ encodes the time-dependent count
and weights this untruncated product law.
If only the first $M$ blocks are removed, the exterior instead
retains the entire future after their last exit, and its terminal
test allows all residual times $q\ge0$ with $X_q=r$. This describes
the event that at least $M$ blocks have completed.

Use the entrance--exit recursion $T_j$ of Appendix~J on the
current inner and outer domains. The ordinary duration bound gives
\begin{equation*}
 \PP_r\{\exists\,1\le i\le M:T_{2i}\le t,
                     \ T_{2i}-T_{2i-1}>d\}\le2Me^{-cd/R^2}.
\end{equation*}
\begin{equation*}
 \PP_t^{r,r}\{\exists\,1\le i\le M:T_{2i}\le t,
             \ T_{2i}-T_{2i-1}>d\}\le CnMe^{-cd/R^2}\le\varepsilon_N.
\end{equation*}
Here $p_t(r,r)\asymp n^{-1}$ uniformly at central times.
On the complementary event, the removed duration is at most $Md$.
For a replacement conditioned on $S$, use \eqref{chIV:eq:R6.6};
otherwise use the unconditional bound. If a normalized short-duration
kernel is used, its factor is at most $(1-\varepsilon_N)^{-1}$,
as in \eqref{chIII:eq:duration-truncation}.
Alternatively, the untruncated one-success or bad-event inequality
of Chapter~III pays $|I|\varepsilon_N$ on each count. Its endpoint
words have total probability one. Thus after summing $m\le M$,
averaging roots, and using the normalization $1/|I_0|$ of
\eqref{chIV:eq:loop-integral}, the additive error is at most
\[
 \frac{|I|}{|I_0|}\,n(M+1)\varepsilon_N
 \le Cn(M+1)e^{-e^{Q/2}}.
\]
The enlarged window is $I+[-Md,Md]$ in either formulation.
Since
\begin{equation*}
Md=n\lceil L^4\rceil e^{-Q}=o(b_N),
\end{equation*}
the enlargement is negligible. For finitely many nested domains, resample successively and enlarge
the window after each step; only disjoint fibers are treated simultaneously.
Even after this count and root normalization, a union over
$N^4L^C$ deterministic vertex-domain-scale choices is summable:
\[
 N^4L^C\{e^{-cL^3}+n(M+1)e^{-e^{Q/2}}\}\longrightarrow0.
\]
Uniform normalization follows from $p_t(r,r)\sim N^{-2}$.

\appendixchapter{M}{Proofs of Propositions~\ref{chIV:prop:local} and~\ref{chIV:prop:count}}\label{chIV:sec:R7}

The first subsection proves the occupation estimates
\eqref{chIV:eq:local-moments}--\eqref{chIV:eq:local-sqrt}, completing
Proposition~\ref{chIV:prop:local}. We then prove the low-count
bound and simultaneous lower barrier in Proposition~\ref{chIV:prop:count}.

\subsection{A direct moment-generating-function estimate}
Use the Appendix~L geometry with radii $e^q,e^{q+w},e^{q+2w}$, $w=q^\gamma$.
Let $\mu_q$ be the exit harmonic measure from $B(y,e^{q+1})$ started at $y$, and define
the one-block average, normalized by $g$, by
\begin{equation*}
A_i=\frac1g\sum_z\mu_q(z)\ell_i(z).
\end{equation*}
For $m$ independent blocks with a fixed feasible entrance--exit word,
probabilities and expectations in this subsection are under the product
of the endpoint kernels in Appendix~L; $\ell_i=\ell^{(i)}$ is the
local-time profile defined in Section~4.1. Then
\begin{equation}
\left|\mathbb EA_i-w\right|\le Cwe^{-cw},\qquad
\mathbb E A_i^j\le j!(Cw)^j\quad(j\ge1).
\label{chIV:eq:R7.1}
\end{equation}
Indeed, with $h_b(z)=\PP_z(X_{\tau_D}=b)$, the potential representation and \eqref{chIV:eq:R6.5} give
\begin{equation*}
G_D(a,z)/g=w+O(e^{-cw}),\qquad
h_b(z)/h_b(a)=1+O(e^{-cw}).
\end{equation*}
on $\supp\mu_q$, proving the first-moment identity. Moreover,
\begin{equation}
\sup_{z\in\operatorname{supp}\mu_q}
\sum_{z'}\mu_q(z')G_D(z,z')/g\le Cw.
\label{chIV:eq:R7.2}
\end{equation}
To verify its scale, combine
\begin{equation*}
G_D(z,z')\le g\log\frac R{|z-z'|\vee1}+C,\qquad
\sum_{z'}\mu_q(z')\log(|z-z'|\vee1)\ge q-C.
\end{equation*}
The harmonic-measure potential lower bound follows by optional stopping: outside the inner ball
its harmonic extension equals the center potential, and inside it adds the killed Green function.
Thus the average Green scale is $\log R-q+O(1)=2w+O(1)$.
For a block with endpoints $a,b$, Kac's formula reads
\[
 EA_i^j=\frac{j!}{g^j}\sum_{z_1,\ldots,z_j}
   \mu_q(z_1)\cdots\mu_q(z_j)
   G_D(a,z_1)\prod_{k=1}^{j-1}G_D(z_k,z_{k+1})
   \frac{h_b(z_j)}{h_b(a)}.
\]
The endpoint factors telescope to $h_b(z_j)/h_b(a)$.
Its uniform bound on $\supp\mu_q$,
followed by \eqref{chIV:eq:R7.2} at each summation, proves
\eqref{chIV:eq:R7.1}.
The centered moment series therefore gives
\begin{equation*}
\log\mathbb E e^{\theta(A_i-\mathbb EA_i)}
\le C\theta^2w^2,\qquad |\theta|\le c/w.
\end{equation*}
for fixed $c,C>0$. Multiply over independent blocks and optimize the exponential parameter:
\begin{equation}
\mathbb P\left(\left|\sum_{i=1}^mA_i-\sum_{i=1}^m\mathbb EA_i\right|>z\right)
\le2\exp\left\{-c\min\left(\frac{z^2}{mw^2},\frac zw\right)\right\}.
\label{chIV:eq:R7.3}
\end{equation}
Put $Y=\sqrt{wm}$, $H=\sqrt{\sum_iA_i}$.
For $Ce^{-cw}Y\le\Delta\le Y$, the corresponding occupation deviation is at least a constant times $Y\Delta$; hence
\begin{equation}
\mathbb P(|H-Y|>\Delta)\le2e^{-c\Delta^2/w}.
\label{chIV:eq:R7.4}
\end{equation}
The same calculation, uniformly for $wm\le\upsilon^2$ and
$Ce^{-cw}\upsilon\le\Delta\le\upsilon$, gives
\begin{equation}
\mathbb P(H>\upsilon+\Delta)\le2e^{-c\Delta^2/w}.
\label{chIV:eq:R7.5}
\end{equation}
For $m=0$ the probability is zero. These estimates hold for each exact exterior word.
Under $o(1)$ radius perturbations replace $w$ by the actual logarithmic width.

Put $Q=L^{16\gamma}$, $d=N^2e^{-Q}$ and $M=\lceil L^4\rceil$;
thus $Md=o(b_N)$.

\subsection{Energy-truncated probes}
Write $h=h^r$ for the pinned field, with the root fixed.
We prove Lemma~\ref{chIV:lem:coarse-probes} from Section~4.1.
\begin{proof}
For $|u|\le1$, regression, the Green difference estimate and conditional-variance
minimization give
\[
 \Var(h(\xi)\mid h(x))\le\Var(h(\xi)-h(x))\le C(q+1),
 \qquad |E(h(\xi)\mid h(x)=u)|\le C.
\]
The last inequality uses $\Var h(\xi)\le CL$ and $\Var h(x)\ge c_0L$.
For the pinned covariance $C_r$, $\operatorname{tr}C_r=2nv_N$ and
$\|C_r(\cdot,x)\|_2^2\le\Var h(x)\operatorname{tr}C_r$. Consequently
\begin{equation}\label{chIV:eq:probe-regression}
 E(V_q\mid h(x)=u)\le C(q+1),\qquad
 E(\|h\|^2\mid h(x)=u)\le CnL.
\end{equation}
Choose $D$ and $C_1$ so the two conditional Markov errors total at most $1/2$.
The density of $h(x)$ on $[-1,1]$ is at least $cL^{-1/2}$, proving
\eqref{chIV:eq:coarse-probe-single}.

For the simultaneous assertion, Jensen's inequality for the probability weights and
the one-dimensional Gaussian-square formula imply, for sufficiently small $c>0$,
\[
 E\left(e^{cV_q/(q+1)}\mid h(x)=u\right)
 \le\sum_\xi\omega_q(\xi)
       E\left(e^{c h(\xi)^2/[2g(q+1)]}\mid h(x)=u\right)\le C.
\]
Thus the conditional probability of any radial failure is at most
$CL e^{-c'a_{\rm min}^\chi}=o(1)$. The common energy restriction loses at most $1/4$ after
increasing $D$. Integrate over $|u|\le1$; the same global cutoff is retained in every
probe. These are the one-point detection bounds used in counted coarea.
\end{proof}

\subsection{Integrated low-count marginals}\label{chIV:sec:shell-envelope}
Let $I_0$ be a central physical window with $|I_0|\asymp b_N$. Define
\begin{equation*}
\mathcal Q_I(F)=\frac{N^2}{|I_0|}\frac1{N^2}
\sum_{r\in V_N}\int_I p_t(r,r)\,
\mathbb E_t^{r,r}F\,dt.
\end{equation*}
This positive integral has total mass $1+o(1)$ at $I=I_0$.

We use the following notation for the root-averaged shell of Chapter~III.
Set
\[
 \mathcal{Q}_N(\phi)=n^{-1}\|\phi\|_2^2,\qquad
 a_s(\phi)=\sqrt{(2s/n-\mathcal{Q}_N(\phi))_+},
\]
and, on $\mathcal{Q}_N(\phi)<2s/n$, put
\[
 \overline b_s^\sigma(\phi)
   =\frac{\sum_r(\phi(r)+\sigma a_s(\phi))_+}{n^2a_s(\phi)},
 \qquad \sigma\in\{-1,1\}.
\]
Set these densities to zero outside that strict set. Define
\[
 \begin{split}
 \varrho_N(s)&=\int\sum_{\sigma=\pm1}
                  \overline b_s^\sigma(\phi)\,\mathcal L(\phi_N)(d\phi),\\
 \nu_{N,s}(d\phi,\{\sigma\})
   &=\frac{\overline b_s^\sigma(\phi)}{\varrho_N(s)}\,\mathcal L(\phi_N)(d\phi).
 \end{split}
\]
Here $s$ ranges over central shells, where $\varrho_N(s)>0$.
Thus $\nu_{N,s}$ is a probability law on the field and branch.
On the strict shell support,
$W_s^\sigma(\phi)=(\phi+\sigma a_s(\phi))^2$; the square image is zero
off that support, as in Chapter~III. The root-averaged coarea
image is the pushforward of
$\varrho_N(s)\,ds\,\nu_{N,s}(d\phi,d\sigma)$ by
$(s,\phi,\sigma)\mapsto(s,W_s^\sigma(\phi))$. The two branches are summed, not
averaged. Chapter~II gives $n\varrho_N(s)=1+o(1)$ uniformly on central
shell intervals, in particular $\varrho_N(s)\le C/n$.

For rare image events and counts, use the following quantitative
shell domination. Fix a
compact central shell interval and $\varepsilon>0$. In the next
calculation put
\[
 B=2s/n,\qquad a_0(s)=\sqrt{B-v_N},\qquad
 k_N=v_N+c_*\varepsilon L\log L,
\]
where $c_*>0$ is a sufficiently small fixed constant.
Uniformly on this interval, for all sufficiently large $N$,
\[
 k_N\le B/2,\qquad
 \mathcal Q_N\le k_N\ \Longrightarrow\
 |a_s-a_0(s)|
 =\frac{|\mathcal Q_N-v_N|}{\sqrt{B-\mathcal Q_N}+\sqrt{B-v_N}}
 \le\varepsilon\log L.
\]
Let $0\le T_s(\phi,a,\sigma)\le K_0$ be a measurable test, with
deterministic $K_0<\infty$, and let $H_s^\sigma(\phi)\ge0$ satisfy
\[
 T_s(\phi,a,\sigma)\le H_s^\sigma(\phi)
 \quad\text{whenever }|a-a_0(s)|\le\varepsilon\log L.
\]
The two-branch identity \eqref{chII:eq:two-branch-algebra} and the
energy density \eqref{chII:eq:energy-density} imply
\[
 n\sum_\sigma\overline b_s^\sigma\le2
 \qquad(\mathcal Q_N\le B/2).
\]
\[
 \begin{aligned}
 &n\EE_\phi\!\left[\one_{\{k_N<\mathcal Q_N<B\}}
                     \sum_\sigma\overline b_s^\sigma\right]\\
 &\quad\le Ce^{-c(k_N-v_N)}
       \int_{k_N}^{B}\left(1+\frac{\sqrt B}{\sqrt{B-q}}\right)dq\\
 &\quad\le CB e^{-c(k_N-v_N)}.
 \end{aligned}
\]
Dividing by $n\varrho_N(s)\ge1/2$ gives the uniform envelope bound
\[
 E_{\nu_{N,s}}T_s(\phi,a_s,\sigma)
 \le C\EE_\phi[H_s^+(\phi)+H_s^-(\phi)]
       +CK_0L^2e^{-c\varepsilon L\log L}.
\]
We use $K_0=1$ below and $K_0=n$ for the image count in
Appendix~O. Both signs remain in the Gaussian envelope; the quotient
in each branch is evaluated only on $\mathcal Q_N<B$.

Use $\omega_{q,x}$ and $\mathcal R_{q,x}^{\rm sq}$ as defined in
Section~4.1.
For $Q\le q\le L-Q$, let $m_q(y)$ be the complete block count and
$Y_q(y)=\sqrt{q^\gamma m_q(y)}$. Impose
\begin{equation}
d(r,y)\ge Ne^{-Q/2}.
\label{chIV:eq:R7.7}
\end{equation}
so the domain excludes the root. For each fixed $C_1<\infty$, uniformly over $q\le\upsilon\le C_1L$,
\begin{equation}
\mathcal Q_{I_0}\bigl(\mathbf1_{\{\text{\eqref{chIV:eq:R7.7}},\,Y_q(y)\le \upsilon\}}\bigr)
\le C\sqrt L\exp\left\{
-\frac{(m_L-\upsilon-L^{2\gamma}-C\log L)_+^2}{L-q+C}
\right\}+e^{-cL\log L},
\label{chIV:eq:R7.8}
\end{equation}
where $m_L=\sqrt2L-\frac{3}{4\sqrt2}\log L$.
Constants may depend on the window and $C_1$, but not on $q,y$.

\begin{lemma}[Low-harmonic mass on the coarea source]\label{chIV:lem:low-harmonic-coarea}
Let $I^*$ be $I_0$ or a central enlargement by $o(b_N)$, and put
\[
 E_{q,y,\upsilon}(t,\omega)
 =\{d(r,y)\ge Ne^{-Q/2},\ \mathcal O_q(y)\le\upsilon+L^{2\gamma}\}.
\]
For $Q\le q\le L-Q$ and $q\le\upsilon\le C_1L$,
\begin{equation}\label{chIV:eq:low-harmonic-coarea}
 \mathcal Q_{I^*}(E_{q,y,\upsilon})
 \le C\sqrt L\exp\left\{
 -\frac{(m_L-\upsilon-L^{2\gamma}-C\log L)_+^2}{L-q+C}\right\}
 +e^{-cL\log L}.
\end{equation}
The constants are uniform over these enlargements, $q,y$ and the indicated thresholds.
\end{lemma}
\begin{proof}
For each fixed eligible root $r$, the Green estimates give
$cL\le\Var h(y)\le CL$. Apply Lemma \ref{chIV:lem:coarse-probes} to obtain a Borel probe
\[
 \begin{gathered}
 B_{q,y}=\left\{|h(y)|\le1,\quad
 \frac1{2g}\sum_z\mu_q(z)h(z)^2\le C_2q,\quad\|h\|^2\le DnL\right\},
 \\ P_h(B_{q,y})\ge cL^{-1/2}.
\end{gathered}
\]
The constants $D,C_2$ are fixed before $N$. On $E_{q,y,\upsilon}\cap B_{q,y}$,
with $W=2\ell_t+h^2$,
\begin{equation}\label{chIV:eq:R7.9}
 \mathcal R_{q,y}^{\rm sq}(W)
 \le\sqrt{(\upsilon+L^{2\gamma})^2+C_2q}
 \le\upsilon+L^{2\gamma}+C_3=:v_q,
\end{equation}
because $\upsilon\ge q$. The image event
$\mathcal B_{q,y}(s,W)=\{\mathcal R_{q,y}^{\rm sq}(W)\le v_q\}$ is Borel
measurable in the image coordinates $(s,W)$. On the source put
$s:=t+\|h\|_2^2/2$; then $s\in I^{*+}:=I^*+[0,DnL/2]$.
Use the counted-probe identity of Chapter~III for one label, sum over eligible roots,
and then enlarge that root sum to all roots. The unnormalized identity gives
\[
 \frac{c}{\sqrt L}\,\mathcal Q_{I^*}(E_{q,y,\upsilon})
 \le\frac{n}{|I_0|}\int_{I^{*+}}\varrho_N(s)
        E_{\nu_{N,s}}\one_{\mathcal B_{q,y}(s,W_s^\sigma(\phi))}\,ds.
\]
Here $|I^{*+}|/|I_0|\le C_D$ uniformly in the enlargement, and $\varrho_N(s)\le C/n$.
The eligible-root source has first been bounded by the full
root-averaged image. The normalizing factor is consequently
$n/|I_0|$, as displayed.

On the strict shell support, Cauchy--Schwarz for the probability
weights $\mu_q$ gives the image inclusion
\[
 \mathcal B_{q,y}(s,W_s^\sigma(\phi))\subseteq
 \left\{\left|\frac1{\sqrt{2g}}\sum_z\mu_q(z)\phi_N(z)
               +\frac{\sigma a_s}{\sqrt{2g}}\right|\le v_q\right\}.
\]
For the centered Gaussian average $Z_{q,y}=(2g)^{-1/2}\sum_z\mu_q(z)\phi_N(z)$,
the periodic Green estimate and the harmonic-measure potential give
\[
 \Var Z_{q,y}
 =\frac1{2g}\sum_{z,z'}\mu_q(z)\mu_q(z')G_N(z,z')
 \le\frac{L-q}{2}+C.
\]
Indeed $G_N(z,z')\le g[L-\log(|z-z'|\vee1)]+C$ on these
planar circles, and the double harmonic average of the logarithm
is at least $q-C$. Also $a_0(s)/\sqrt{2g}=m_L+O(1)$ on the
central shell interval. The envelope bound just proved, with
$K_0=1$, reduces either branch to
\[
 \{|Z_{q,y}|\ge(m_L-v_q-C\log L)_+\}.
\]
The Gaussian tail is at most
$2\exp\{-(m_L-v_q-C\log L)_+^2/(L-q+C)\}$, with a larger
constant covering the zero positive part. The shell remainder is
$CL^2e^{-c\varepsilon L\log L}$; multiplication by $\sqrt L$
and the bounded window ratio leaves $e^{-c'L\log L}$.
This proves
\eqref{chIV:eq:low-harmonic-coarea}.
\end{proof}

To deduce \eqref{chIV:eq:R7.8}, remove all complete $q$-blocks and take
$A=\{d(r,y)\ge Ne^{-Q/2},\allowbreak Y_q(y)\le\upsilon\}$ and
$G=\{\mathcal O_q(y)\le\upsilon+L^{2\gamma}\}$.
The count and root condition are exterior-measurable; on $A$,
$m_q\le C_1^2L^2/q^\gamma<M$.
For each exact exterior word, \eqref{chIV:eq:R7.5} gives
$P(G^c)\le2e^{-cL^{3\gamma}}$.
Its threshold conditions follow from $\upsilon\ge q\ge Q$;
at $m_q=0$ success is automatic.
Appendix~L gives, under every such untruncated product kernel,
$P(\sum_i\tau_i>Md)\le\varepsilon_N=e^{-e^{Q/2}}$.
Since $2e^{-cL^{3\gamma}}+\varepsilon_N\le1/2$ eventually, apply
\eqref{chIII:eq:counted-window} on each fixed count, then sum its
endpoint probabilities and the counts $0\le m\le M$. It gives
\[
 \mathcal Q_{I_0}(A)
 \le2\mathcal Q_{I_0+[-Md,Md]}
       (d(r,y)\ge Ne^{-Q/2},\ \mathcal O_q(y)\le\upsilon+L^{2\gamma})
       +n(M+1)\varepsilon_N.
\]
Apply Lemma \ref{chIV:lem:low-harmonic-coarea} to the enlarged-window
term and absorb the factor $2$ in $C$. The normalized duration error
is smaller than $e^{-cL\log L}$, proving \eqref{chIV:eq:R7.8}
with its stated precision.

\subsection{A global lower count barrier}
For each integer $q\in[Q,L-Q]$, choose a deterministic grid with spacing comparable to $e^q$
and at most $Ce^{2(L-q)}$ points.
Moving a center $x$ to its nearest grid point $y$, perturb the inner radius inward and the outer
radius outward by logarithmic amounts $O(e^{-q^\gamma})$, so that
\begin{equation*}
B(y,\rho')\subset B(x,\rho),\qquad
B(x,R)\subset B(y,R').
\end{equation*}
Each crossing from the smaller inner ball to outside the larger outer ball contains a distinct
original crossing; chronological matching is injective. Thus the grid count is at most the original count.
The new width is $w+O(e^{-w})$, and for square-root counts $O(L)$ the normalization error is
$O(Le^{-w}/w)=o(1)$.
Apply \eqref{chIV:eq:R7.8} with $\upsilon=\sqrt2q-L^{3\gamma}+o(1)\ge q$.
Then
\begin{equation*}
\frac{(m_L-\upsilon-L^{2\gamma}-C\log L)^2}{L-q+C}
\ge2(L-q)+cL^{3\gamma}
\end{equation*}
uniformly for large $N$, since its numerator is
$(\sqrt2(L-q)+L^{3\gamma}-L^{2\gamma}-O(\log L))^2$ and $L-q\ge Q$.
The grid factor $e^{2(L-q)}$ cancels, leaving
$O(L^{3/2})e^{-cL^{3\gamma}}\to0$; the remainder is summable over all grid points.
Therefore
\begin{equation*}
 \mathcal Q_{I_0}\!\left(\left\{\begin{gathered}
 \exists x\in V_N,\ q\in[Q,L-Q]\cap\mathbb Z:\\
 d(x,r)\ge2Ne^{-Q/2},\qquad Y_q(x)<\sqrt2q-L^{3\gamma}
 \end{gathered}\right\}\right)=o(1).
\end{equation*}
The stronger root-distance restriction absorbs the grid displacement.
The excluded root neighborhood is handled by the one-point bridge soft bound when this result is used below.

\appendixchapter{N}{Proof of Proposition~\ref{chIV:prop:deviation}}\label{chIV:sec:R8}

We prove \eqref{chIV:eq:deviation-main}, with the root-near error
estimated at the end. Fix a softness threshold $u_0\ge0$.
Use $\gamma=1/8192$, $Q=L^{16\gamma}$ and entrance/exit logarithmic radii
$q+q^\gamma,q+2q^\gamma$. The integral $\mathcal Q_I$ is
\eqref{chIV:eq:loop-integral}, with $I=I_0$ or a central enlargement
by $o(b_N)$.
Let $Y_q(x)=\sqrt{q^\gamma m_q(x)}$ and recall
$\mathcal O_q(x)=(g^{-1}\sum_z\mu_q(z)\ell_t(z))^{1/2}$ from Section~4.1.
The weights are the exit harmonic measure of $B(x,e^{q+1})$.
Let $\Theta$ be the event that $Y_k(x)\ge\sqrt2k-L^{3\gamma}$ for every integer $k\in[Q,L-Q]$
and every $x\in V_N$ with $d(x,r)\ge2Ne^{-Q/2}$.
Appendix~M gives $\mathcal Q_I(\Theta^c)\to0$.
We prove
\begin{equation}
\mathcal Q_I\!\left[
\mathbf1_\Theta\sum_{\substack{x:\ d(x,r)\ge2Ne^{-Q/2}}}
\sum_{q\in[2Q,L-2Q]\cap\mathbb Z}
\mathbf1_{\{\ell_t(x)\le u_0,\ |\mathcal O_q(x)-Y_q(x)|>L^{4\gamma}\}}
\right]\longrightarrow0.
\label{chIV:eq:R8.1}
\end{equation}
\subsection{Endpointwise local inputs}
Use the endpointwise binomial--geometric and binomial--exponential
comparisons in \cite[Lemmas A.6--A.7, equations (A.56)--(A.58),
(A.86)--(A.87)]{chIV:LS24}, followed by their Chernoff estimates.
Their local statements allow arbitrary fixed $0<c_1<c_2$, hence $c_1=1,c_2=2$,
uniformly in each complete entrance-exit word.
Fix the centre $x$ and abbreviate $Y_k=Y_k(x)$ and $\mathcal O_k=\mathcal O_k(x)$
in the local calculations below. Fix a feasible
$l$-word $\mathbf w=((a_i,b_i))_{i=1}^m$ and retain
the complementary path. Let
$K_l^{\mathbf w}=\bigotimes_{i=1}^m P^{a_i,b_i}$ be the untruncated
product law of its complete $l$-blocks. If $\ell^{(i)}(x)$ is the
occupation at $x$ during block $i$, write, for $w'=\sqrt{ml^\gamma}$,
\[
 p_{\rm soft}(w')
   =K_l^{\mathbf w}\left(\sum_{i=1}^m\ell^{(i)}(x)\le u_0\right),
\]
suppressing the fixed endpoints from the notation. The bound is uniform
over them. When the root is outside the outer disk, the complementary
path contributes no occupation at $x$, so this is the softness event
on the reconstructed loop. We first evaluate the endpoint-conditioned local kernel
$K_l^{\mathbf w}$; the physical-time restriction is then introduced by
the window-transfer theorem.
For $w'=\sqrt2l+s'$, $s'\ge-l$,
\begin{equation}
p_{\rm soft}(w')
\le\exp\left\{-2l-2\sqrt2s'
-\frac{(s')^2}{l}(1-Cl^{\gamma-1})+Cl^\gamma\right\}.
\label{chIV:eq:R8.2}
\end{equation}
The occupation convention in that reference is $\ell/g$, so its
threshold is $u_0/g$; for $u_0=0$, bound by any fixed positive threshold.
The one-site input can also be read from the exact transform
\[
 E^{a,b}\exp\{-\theta\ell^{(i)}(x)/g\}
 =1-p_{a,b}\frac{\theta}{\lambda_x+\theta},\qquad
 p_{a,b}=\frac{G_D(a,x)h_b(x)}{G_D(x,x)h_b(a)},\quad
 \lambda_x=\frac g{G_D(x,x)}.
\]
It retains every return before exit. The Green and Harnack estimates
give $p_{a,b}\ge(l^\gamma/l)(1-Cl^{\gamma-1})$ and
$\lambda_x\le2/l$. Multiplication of these transforms over the
blocks and Chernoff optimization give \eqref{chIV:eq:R8.2}, retaining
the negative quadratic term until the linear errors have been absorbed.

For the nested scale $l=\lfloor q-3q^\gamma\rfloor$, so that
$B(x,e^{l+2l^\gamma})\subset B(x,e^q)$, let
$K_q^{\mathbf v}$ be the product block law of a feasible $q$-word
$\mathbf v$; its block count fixes $Y_q=w$.
Let $N_i$ count complete $l$-traversals within outer block $i$.
For every retained exterior compatible with a completed loop rooted
outside the outer domain, chronological inclusion gives
\[
 m_l\ge\sum_{i=1}^m N_i,\qquad
 \{Y_l\le w'\}\subseteq
 \left\{l^\gamma\sum_{i=1}^m N_i\le(w')^2\right\}.
\]
The outer blocks start outside the inner exit domain, so their
complete inner traversals occupy disjoint chronological slots.
We bound the event on the right under $K_q^{\mathbf v}$; the same
bound therefore holds for the reconstructed $Y_l$ on every such
exterior fibre. If $w-w'\ge1$ and
$0\le w'<w\le e^{\delta l^\gamma}$ for a fixed sufficiently small $\delta>0$, then
\begin{equation}
K_q^{\mathbf v}(Y_l\le w')
\le2\exp\left\{-\frac{(w-w')^2}{4[q+2q^\gamma-(l+l^\gamma)]}\right\}.
\label{chIV:eq:R8.3}
\end{equation}
To make the endpoint and rounding errors explicit, put
\[
 G=q+2q^\gamma-(l+l^\gamma)\asymp q^\gamma,\quad
 p=\frac{q^\gamma-\epsilon_q}{G},\quad
 e=\frac{l^\gamma+\epsilon_q}{G},\quad
 0\le\epsilon_q\le Cq^\gamma e^{-cl^\gamma}.
\]
Annular stopping and \eqref{chIV:eq:R6.5} give a first inner
entrance probability at least $p$ and an escape probability at most
$e$ after every inner traversal. Strong Markov iteration gives
$P(N_i\ge j+1)\ge p(1-e)^j$. Thus
\[
 E e^{-\theta l^\gamma N_i}
 \le1-p+p\frac{e}{e+\theta l^\gamma},\qquad
 K_q^{\mathbf v}(Y_l\le w')
 \le\exp\!\left\{-\left(\sqrt{mp}-w'\sqrt{e/l^\gamma}\right)_+^2\right\}.
\]
Here $m=w^2/q^\gamma$. For $\delta<c/2$ the error in the
square-root gap, multiplied by $\sqrt G$, is at most
$Cw\epsilon_q/l^\gamma=o(1)$. A gap $w-w'\ge1$ therefore
loses at most a factor $2$, proving \eqref{chIV:eq:R8.3}.
These are killed planar kernels. Appendix~L supplies their root-free torus and time-window versions.

\subsection{Nested scales and high-count branches}
For $q$ in the target range of \eqref{chIV:eq:R8.1}, set $l=\lfloor q-3q^\gamma\rfloor$.
Then $l\in[Q,L-Q]$ and $l+2l^\gamma<q$.
After removal of the $l$-domain blocks, both $Y_q$ and $\mathcal O_q$ are exterior-measurable;
the averaging circle for $\mathcal O_q$ is outside that domain.

First exclude soft branches with $Y_q>10L$ or $Y_l>10L$.
For $k=q,l$, remove only the first $\lceil100L^2/k^\gamma\rceil<L^4$ blocks.
Softness for the original path implies softness for these blocks alone, and \eqref{chIV:eq:R8.2} gives
\begin{equation*}
p_{\rm soft}\le\exp\left\{-\frac{100L^2}{k}(1-o(1))+O(k^\gamma)\right\}
\le e^{-40L}
\end{equation*}
On an ordinary infinite path, condition on the realized word of
these first blocks and integrate against its probability law. The
word-uniform bound $e^{-40L}$ therefore also holds after this averaging. The original high-count soft event implies this
prefix softness event. At a central bridge endpoint its probability
is at most $Cn e^{-40L}$, by division by $p_t(r,r)\asymp n^{-1}$.
Hence the expected number of high-count soft vertex--scale pairs
under $\mathcal Q_I$ is at most $Cn^2L e^{-40L}=o(1)$.
Henceforth $Y_q=w\le10L$, $Y_l=w'\le10L$, and
\begin{equation*}
s=w-\sqrt2q,\qquad s'=w'-\sqrt2l.
\end{equation*}
On $\Theta$, $s,s'\ge-L^{3\gamma}$. Expanding the positive-part square in \eqref{chIV:eq:R7.8} gives
\begin{equation}
\mathcal Q_{I'}(Y_q=w,\ \text{root-far})
\le CL^{1/2}e^{-2(L-q)+2\sqrt2s+CL^{2\gamma}}
+e^{-cL\log L}.
\label{chIV:eq:R8.4}
\end{equation}
uniformly over the finitely enlarged central windows $I'$.
When its positive part vanishes, the resulting upper-bound exponent is already positive;
otherwise use $s\ge-L^{3\gamma}$ and $L-q\ge Q$.

\subsection{A weighted deviation estimate}
Fix an exact $q$-word with $\sqrt2q-L^{3\gamma}\le w\le10L$, and denote its
untruncated local kernel by $K$. Put
$X=s-(Y_l-\sqrt2l)$ and
$J_l=\{Y_l-\sqrt2l\ge-L^{3\gamma},\ Y_l\le10L\}$.
For $u\ge0$, the event $X>u$ requires
$Y_l<w-\sqrt2(q-l)-u$. If this threshold is below $\sqrt2l-L^{3\gamma}$,
the event is disjoint from $J_l$; otherwise \eqref{chIV:eq:R8.3} applies:
the count gap is at least $1$, its denominator is at most $Cq^\gamma$,
and $w\le10L<e^{\delta l^\gamma}$.
Rounding a threshold down to an attainable count only strengthens the bound.
Thus, for every fixed $\lambda>0$,
\begin{equation}\label{chIV:eq:excess-exponential-moment}
 \begin{split}
 K(X>u,J_l)&\le2e^{-cu^2/q^\gamma},\\
 E_K[e^{\lambda X}\one_{J_l}]
 &\le1+2\lambda\int_0^\infty e^{\lambda u-cu^2/q^\gamma}\,du
 \le C_\lambda e^{C_\lambda q^\gamma}.
 \end{split}
\end{equation}
For $E_q=\{|\mathcal O_q-w|>L^{4\gamma}\}$, the threshold conditions of
\eqref{chIV:eq:R7.4} hold since $w\ge q$, $L^{4\gamma}=o(q)$ and
$we^{-cq^\gamma}=o(L^{4\gamma})$. Cauchy--Schwarz, with
$\lambda=4\sqrt2$ above, yields
\begin{equation}\label{chIV:eq:weighted-count-deviation}
 E_K[e^{2\sqrt2X}\one_{J_l}\one_{E_q}]
 \le C\exp\{Cq^\gamma-cL^{8\gamma}/q^\gamma\}.
\end{equation}

First resample the $l$-domain. The two counts and $E_q$ are exterior
at this step. Retain only $s,s'\ge-L^{3\gamma}$ from $\Theta$ and discard its
other conditions. On $J_l$, \eqref{chIV:eq:R8.2} gives
$p_{\rm soft}(Y_l)\le e^{-2l-2\sqrt2s+Cq^\gamma}e^{2\sqrt2X}$.
On the $q$-product space retain the bounded weight
\[
 f=\min\{1,e^{-2l-2\sqrt2s+Cq^\gamma}e^{2\sqrt2X}\}
                   \one_{J_l}\one_{E_q}.
\]
Then resample the $q$-domain using the bounded-weight form of
\eqref{chIII:eq:bad-transfer} and apply
\eqref{chIV:eq:weighted-count-deviation} to bound $E_Kf$.
The bounded local weight therefore has expectation at most
\[
c_w=\min\{1,\ C\exp(-2l-2\sqrt2s+Cq^\gamma-cL^{8\gamma}/q^\gamma)\}.
\]
The two window transfers multiply this bound by at most a fixed constant.
These are successive resamplings, with total window enlargement
$o(b_N)$. The duration-failure error at each step, for a local
test bounded by $1$, is at most $nL^C e^{-e^{Q/2}}$ after count/root
normalization. By \eqref{chIII:eq:duration-truncation}, the two
conditionings multiply the nonnegative expectation bound by at most
$(1-\varepsilon_N)^{-2}$, where $\varepsilon_N=e^{-e^{Q/2}}$;
this changes only $C$.
The joint term is controlled directly by the weighted expectation and
Cauchy--Schwarz estimate above.

Multiply $c_w$ by \eqref{chIV:eq:R8.4}. In its main term the factors
$e^{2\sqrt2s}$ and $e^{-2\sqrt2s}$ cancel, and the exponent is at most
\begin{equation}\label{chIV:eq:count-deviation-summable}
 -2L+2(q-l)+CL^{2\gamma}+Cq^\gamma-cL^{8\gamma}/q^\gamma
 \le-2L-c'L^{7\gamma}.
\end{equation}
The additive remainder is multiplied by $c_w\le1$.
Only the at most $CL^2$ outer counts must now be summed, along with
at most $L$ scales and $e^{2L}$ vertices. Both displayed errors and
the duration failures are summable. Adding the high-count branches
proves \eqref{chIV:eq:R8.1}.

Together with Appendix~M, this gives
\begin{equation*}
 \mathcal Q_{I_0}\!\left(\left\{\begin{gathered}
 \exists x\in U_{u_0}(t),\ q\in[2Q,L-2Q]\cap\mathbb Z:\\
 d(x,r)\ge2Ne^{-Q/2},\qquad
 |\mathcal O_q(x)-Y_q(x)|>L^{4\gamma}
 \end{gathered}\right\}\right)=o(1).
\end{equation*}
The expected number of soft vertices in the excluded root neighborhood is at most
$CL^{3/2}e^{-Q}=o(1)$.
The Appendix~O tolerance gives the available interval
$[L^{32\gamma},L-L^{32\gamma}]$, with $32\gamma<1/10$.

\appendixchapter{O}{Proof of Proposition~\ref{chIV:prop:repulsion}}\label{chIV:sec:R9}

The first three subsections establish the Gaussian estimate used in
Proposition~\ref{chIV:prop:repulsion}. The last transfers it to the
local-time exceptional set $\mathcal O_{N,u}(t)$ by counted coarea.

Use the field normalization of \cite{chIV:LS24},
$\Phi_N^{\rm LS}=\phi_N/\sqrt{2g}$, with point variance $L/2+O(1)$.
All Gaussian probabilities and expectations in this appendix refer to
the field currently specified; conditioning is displayed explicitly.
Set
\begin{equation*}
m_L=\sqrt2L-\frac3{4\sqrt2}\log L,\qquad
A=L^{32\gamma},\quad r_*=L^{1/10},\quad \gamma=1/8192.
\end{equation*}
For each integer $q$ and real field $f$, let
\[
 \mathcal R_q(f;x)=\left(\sum_\xi\omega_{q,x}(\xi)f(\xi)^2\right)^{1/2},
\]
with the exit weights of Appendix~M. In particular
$\mathcal R_q(f;x)=\mathcal R_{q,x}^{\rm sq}(2g f^2)$.
Fix $C_0,C_1<\infty$.
A Gaussian vertex is bad if some $|a-m_L|\le\varepsilon\log L+C_0$ satisfies
\begin{equation*}
|{\Phi_N^{\rm LS}}(x)-a|\le C_0,
\end{equation*}
\begin{equation*}
\mathcal R_q({\Phi_N^{\rm LS}}-a;x)\ge\sqrt2q-C_1L^{8\gamma}
\quad\text{for every integer }q\in[A,L-A],
\end{equation*}
\begin{equation}
\mathcal R_q({\Phi_N^{\rm LS}}-a;x)\le\sqrt2q+C_1L^{8\gamma}
\quad\text{for some integer }q\in[r_*,L-r_*].
\label{chIV:eq:R9.1}
\end{equation}
For sufficiently small fixed $\varepsilon>0$, we show that the expected number of bad vertices tends to zero.

\subsection{A Brownian bridge bound}
Let $B$ be a standard Brownian bridge from $0$ to $0$ of length $T$;
thus $B$ is centered Gaussian with covariance
$\operatorname{Cov}(B_s,B_t)=\min(s,t)-st/T$.
For $0<r<T/3$ and $a,b>0$,
\begin{equation}
\mathbb P\left(\sup_{[0,T]}B\le a,\
\sup_{[r,T-r]}B\ge a-b\right)
\le C\frac{a^2b}{T\sqrt r}.
\label{chIV:eq:R9.2}
\end{equation}
To prove it, replace $B$ by $a-B$, a bridge from $a$ to $a$ required to stay nonnegative
and enter $[0,b]$ at an intermediate time. In this subsection
$p_t(u,v)=(2\pi t)^{-1/2}e^{-(u-v)^2/(2t)}$ is the Brownian transition
density. The kernel killed at $0$ is
\begin{equation*}
p_t^+(u,v)=\frac1{\sqrt{2\pi t}}
\left(e^{-(u-v)^2/(2t)}-e^{-(u+v)^2/(2t)}\right)
\le Cuvt^{-3/2}.
\end{equation*}
By time reversal it suffices to consider $[r,T/2]$. Let
$\tau=\inf\{t\ge r:a-B_t\le b\}$, including $\tau=r$ when
$a-B_r\le b$. On survival up to $\tau\le T/2$, the remaining killed density
to $a$ is at most $CabT^{-3/2}$.
Survival to $r$ has probability at most $Ca/\sqrt r$.
Multiply, divide by $p_T(a,a)=(2\pi T)^{-1/2}$, and multiply by $2$ for reversal.
In particular, $a,b\le Ch$ gives $Ch^3/(T\sqrt r)$.

\subsection{Dirichlet harmonic averages as a Gaussian bridge}
Let $D$ be a lattice disk or square of diameter comparable to $N$,
and restrict centers to a region at distance at least $cN$ from its
boundary, for a fixed $c>0$. Work in the normalization
$\Phi_N^{\rm LS}=\phi_N/\sqrt{2g}$.
In this appendix $h_D$ denotes the discrete Dirichlet field with covariance
$G_D/(2g)$, equivalently the law of $\eta_D/\sqrt{2g}$ in the discrete
notation of Chapter~II; $G_D$ is the counting-measure killed Green kernel.
For fixed $x$ and integer scales $q$ with $B(x,e^{q+1}+2)\subset D$, let
$M_q=\sum_z\omega_{q,x}(z)h_D(z)$ be the
exit harmonic average on $B(x,e^{q+1})$, and set
\begin{equation*}
Z_q=h_D(x)-M_q,\quad
v_q=\operatorname{Var}(Z_q)=\tfrac12q+O(1),\quad
\upsilon_D=\operatorname{Var}(h_D(x))=\tfrac12L+O(1).
\end{equation*}
The process $Z_q$ has independent increments and covariance
\[\Cov(Z_q,Z_p)=v_{q\wedge p}.\]
Its terminal value is $h_D(x)$ at variance time $\upsilon_D$.
Conditional on $h_D(x)=m_L+u$, a standard bridge $\mathbb B$ of length $\upsilon_D$ therefore satisfies
\begin{equation}
M_q-m_L=u-\frac{v_q}{\upsilon_D}(m_L+u)-\mathbb B_{v_q}
=-\sqrt2q+e_q-\mathbb B_{v_q},
\quad |e_q|\le C(1+\log L+|u|).
\label{chIV:eq:R9.3}
\end{equation}
The potential expansion gives $v_{q+1}-v_q=1/2+o(1)$ uniformly, hence bounded gaps in variance time.
The exact square-average decomposition is
\begin{equation}
\mathcal R_q(h_D-m_L;x)^2
=(M_q-m_L)^2+V_q,\qquad
V_q=\sum_z\mu_q(z)(h_D(z)-M_q)^2.
\label{chIV:eq:R9.4}
\end{equation}
Here $\mu_q=\omega_{q,x}$ for the fixed center $x$.
For $|u|\le A/10$, each boundary residual $h_D(z)-M_q$ has conditional variance at most $Cq$
and uniformly bounded conditional mean. Indeed,
$\Cov(h_D(z)-M_q,h_D(x))=O(1)$ by regular-Green variation and exit rounding,
while $(m_L+u)/\upsilon_D=O(1)$.
Jensen and the one-dimensional Gaussian-square moment formula give
\begin{equation*}
\mathbb E(e^{cV_q/q}\mid h_D(x)=m_L+u)\le C.
\end{equation*}
Consequently, except with conditional probability $Ce^{-c(\log L)^2}$,
$V_q\le Cq(\log L)^2$ at all at most $L$ scales.

Consider the lower square-root barrier on $[A,L-A]$ and an upper barrier at some scale
in $[r_*,L-r_*]$, allowing errors $o(A)$.
By \eqref{chIV:eq:R9.4}, the lower barrier implies
$|M_q-m_L|\ge\sqrt2q-o(A)-C(\log L)^2$.
Together with \eqref{chIV:eq:R9.3}, this gives $\mathbb B_{v_q}\ge-CA$,
apart from the opposite-sign branch $\mathbb B_{v_q}\le-cq$, whose union probability is $Ce^{-c'A}$.
The upper barrier gives $\mathbb B_{v_q}\le CA$ at an intermediate scale.
After reversing sign, the bridge stays below $CA$ at the grid times and reaches at least $-CA$ in the middle.
Endpoint intervals of variance length $O(A)$ contribute $Ce^{-cA}$;
between-grid fluctuations above $A$ contribute $CL e^{-cA^2}$.
Apply \eqref{chIV:eq:R9.2} with $a,b=O(A)$, $T\asymp L$, and $r\asymp r_*$ to obtain
\begin{equation}
C\frac{A^3}{L\sqrt{r_*}}+Ce^{-c(\log L)^2}.
\label{chIV:eq:R9.5}
\end{equation}
uniformly over deep centers, $|u|=o(A)$, and barrier errors $o(A)$.

\subsection{Adding the periodic exterior field}
In each fixed deep chart,
\begin{equation*}
{\Phi_N^{\rm LS}}=h_D^0+H_N,
\qquad h_D=h_D^0+q_N Z,
\end{equation*}
Here $h_D^0$ has the law of $h_D$ conditioned to have discrete average
zero, that is, the law of $\eta_D^0/\sqrt{2g}$ from Chapter~II.
Also $H_N\perp h_D^0$, and the local variable $Z\sim N(0,1)$ is
adjoined independently of the original pair $(h_D^0,H_N)$.
Thus both displayed equalities hold on this enlarged probability
space, and $h_D$ is independent of $H_N$. In this subsection the normalized torsion is
\[
 q_N(x)=\frac{\sum_{y\in D}G_D(x,y)}
              {\sqrt{2g\sum_{u,v\in D}G_D(u,v)}},\qquad x\in D.
\]
It is uniformly bounded inside.
The covariance of $H_N$ is periodic Green minus Dirichlet Green plus the torsion rank-one term.
Let $\|H_N\|_\infty$ denote the supremum on a fixed interior region
containing the centers and all averaging circles. Cancellation of
singularities, interior increment estimates, and Gaussian maximum bounds give
\begin{equation}
\mathbb P(\|H_N\|_\infty>t)\le Ce^{-ct^2},\qquad t\ge0.
\label{chIV:eq:R9.6}
\end{equation}
on a fixed deep compact set containing all circles, whose largest radius is
$Ne^{-A+O(1)}=o(N)$. Finitely many such charts cover the torus.

Since $\PP(|Z|\le1)>0$, for a periodic bad-vertex event $E_x$
at the fixed center $x$,
\[
 \PP(E_x\mid H_N)
 =\frac{\PP(E_x,\ |Z|\le1\mid H_N)}{\PP(|Z|\le1)}.
\]
Radial norms and endpoints differ between periodic and Dirichlet fields by at most
$\|H_N\|_\infty+C$.
On $\|H_N\|_\infty\le L^{12\gamma}$, these are $o(A)$ errors, so the periodic event implies
the Dirichlet event of \eqref{chIV:eq:R9.5}.
If $E_{D,x}(H_N)$ denotes the resulting Dirichlet event with these
enlarged endpoint and radial tolerances, the inclusion gives
\[
 \PP(E_x\mid H_N)\le\frac{\PP(E_{D,x}(H_N)\mid H_N)}{\PP(|Z|\le1)}.
\]
The last probability uses the ordinary Dirichlet law of $h_D$,
independent of $H_N$, with $h_D=h_D^0+q_NZ$ on the stated joint space.
For every real $u$, the Dirichlet endpoint density at $m_L+u$ satisfies
\begin{equation}
f_{D,x}(m_L+u)\le CN^{-2}L e^{C|u|}.
\label{chIV:eq:R9.7}
\end{equation}
by expansion of $-(m_L+u)^2/(L+O(1))$ and removal of the negative quadratic term.
The $\varepsilon\log L$ amplitude interval costs a length $O(1+\log L)$ and a factor $L^{C\varepsilon}$.
For fixed $H_N$ the extra cost is at most $e^{C\|H_N\|_\infty}$, integrable uniformly by \eqref{chIV:eq:R9.6}.
After summation over vertices, the truncated-field contribution is at most
\begin{equation}
C(1+\log L)L^{C\varepsilon}
\left(L^{3\cdot32\gamma-1/20}+L e^{-c(\log L)^2}\right).
\label{chIV:eq:R9.8}
\end{equation}
On $\|H_N\|_\infty>L^{12\gamma}$, retain only the endpoint bound to obtain
\begin{equation*}
CL^{1+C\varepsilon}(1+\log L)
\mathbb E\left[e^{C\|H_N\|_\infty};\ \|H_N\|_\infty>L^{12\gamma}\right]
\le CL^{1+C\varepsilon}(1+\log L)e^{-cL^{24\gamma}}\to0.
\end{equation*}
The factor $N^{-2}$ is kept before summing vertices.
Since
\begin{equation*}
\frac1{20}-3\cdot32\gamma
=\frac{49}{1280}>0,
\end{equation*}
one can choose fixed sufficiently small $\varepsilon>0$ so that \eqref{chIV:eq:R9.8} tends to zero.
Summing finitely many charts proves the claimed bad-vertex expectation bound.
Gaussian symmetry treats the opposite sign.

\subsection{Coarse-repulsion mass under integrated loops}
Fix $u_0\ge0$ and use the set $\mathcal O_{N,u_0}(t)$ defined in
Section~4.2, with $u=u_0$. At $u_0=0$ the same definition applies.
With $\Theta$ from Proposition~\ref{chIV:prop:deviation}, define
\[
 \begin{split}
 \mathcal G_N={}&\Theta\cap\left\{|\mathcal O_q(x)-Y_q(x)|\le L^{4\gamma}
   \text{ for all }\substack{x\in U_{u_0}(t),\ d(x,r)\ge2Ne^{-Q/2}\\
                    q\in[2Q,L-2Q]\cap\mathbb Z}\right\}\\
 &\cap\{\nexists x\in U_{u_0}(t):d(x,r)<2Ne^{-Q/2}\}.
 \end{split}
\]
The dependence on $(t,\omega,r)$ is suppressed. Appendices~M.4 and
N.3 and the one-point bound in Appendix~Q.1 give
$\mathcal Q_I(\mathcal G_N^c)=o(1)$.

For a fixed sufficiently large constant $C_*$, define, for nonnegative $W:V_N\to\mathbb R$,
\[
 \mathcal B_N(s,W)=\left\{x\in V_N:\begin{gathered}
 W(x)\le2u_0+1,\\
 \mathcal R_{q,x}^{\rm sq}(W)\ge\sqrt2q-C_*L^{8\gamma}
       \quad\forall q\in[A,L-A]\cap\mathbb Z,\\
 \mathcal R_{q,x}^{\rm sq}(W)\le\sqrt2q+C_*L^{8\gamma}
       \quad\exists q\in[r_*,L-r_*]\cap\mathbb Z
 \end{gathered}\right\}.
\]
This set reads only the image profile $W$; $s$ records its shell time.

\begin{lemma}[Vertex-count coarea transfer]\label{chIV:lem:coarse-vertex-coarea}
For $C_*$ sufficiently large there is a fixed $D<\infty$ such that
\begin{equation}\label{chIV:eq:coarse-vertex-transfer}
 \frac1{\sqrt L}\mathcal Q_I[|\mathcal O_{N,u_0}(t)|;\mathcal G_N]
 \le C\sup_{s\in I^+}E_{\nu_{N,s}}|\mathcal B_N(s,W_s^\sigma(\phi))|,
 \qquad I^+=I+[0,DnL/2],
\end{equation}
where $D$ is fixed. Moreover,
\begin{equation}\label{chIV:eq:coarse-vertex-image-bound}
 \sup_{s\in I^+}E_{\nu_{N,s}}|\mathcal B_N(s,W_s^\sigma(\phi))|
 \le C E|\mathcal B_N^{\rm G}|+CnL^2e^{-c\varepsilon L\log L}=o(1).
\end{equation}
Here $\mathcal B_N^{\rm G}$ is the ordinary Gaussian bad-vertex set of
\eqref{chIV:eq:R9.1}, with fixed enlarged constants and sufficiently small fixed
$\varepsilon>0$.
\end{lemma}
\begin{proof}
Fix the deterministic label $x$ and the source path. On
$E_x:=\mathcal G_N\cap\{x\in\mathcal O_{N,u_0}(t)\}$ the root-distance restriction implies
$\Var h(x)\asymp L$. Apply \eqref{chIV:eq:coarse-probe-simultaneous} with
minimum scale $A=L^{32\gamma}$ and $\chi=\gamma$ to the same pinned Gaussian field at all coarse scales. It gives
\[
 B_x=\left\{\begin{gathered}|h(x)|\le1,\quad
 \frac1{2g}\sum_\xi\omega_{q,x}(\xi)h(\xi)^2\le Cq^{1+\gamma}
     \\ (q\in[A,L-A]\cap\mathbb Z),\quad
 \|h\|^2\le DnL\end{gathered}\right\},
\]
with $P_h(B_x)\ge cL^{-1/2}$ uniformly in the path and $x$ on $E_x$.
Here $Le^{-cA^\gamma}\to0$; the probability bound is
\eqref{chIV:eq:coarse-probe-simultaneous}.

On $E_x$, Appendices~M.4 and N.3 give the required lower bounds on $\mathcal O_q$. At the chosen low-count
scale they give $\mathcal O_q\le\sqrt2q+L^{8\gamma}+L^{4\gamma}$. On $B_x$,
\[
 \mathcal R_{q,x}^{\rm sq}(2\ell_t+h^2)^2
 =\mathcal O_q(x)^2+\frac1{2g}\sum_\xi\omega_{q,x}(\xi)h(\xi)^2.
\]
Adding the last term changes a norm of order $q$ by $O(q^\gamma)$, and
$q^\gamma\le L^\gamma=o(L^{8\gamma})$. The lower bounds only improve.
With $s=t+\|h\|_2^2/2$, choosing $C_*$ to dominate these fixed errors gives, pointwise,
\[
 \one_{\{t\in I\}}\sum_x\one_{E_x}\one_{B_x}
 \le\one_{\{s\in I^+\}}|\mathcal B_N(s,2\ell_t+h^2)|.
\]
Apply the counted-probe identity of Chapter~III, then average the root and normalize by
$|I_0|/n$. This proves \eqref{chIV:eq:coarse-vertex-transfer}, using
$\varrho_N(s)\le C/n$ and the bounded interval-length ratio.

For each shell time $s$, the image count is a Borel function of the
profile $W$ with values in $[0,n]$. Apply
the quantitative envelope bound in Appendix~M.3 with $K_0=n$
and amplitude tolerance $\varepsilon\sqrt{2g}\log L$.
On this amplitude window the square profile is
$(\phi_N+\sigma a_s)^2$, with
$|a_s/\sqrt{2g}-m_L|\le\varepsilon\log L+C$.
For each sign, the image bad-vertex count is bounded pointwise by
the count in \eqref{chIV:eq:R9.1} for $-\sigma\Phi_N^{\rm LS}$,
with enlarged fixed constants. Gaussian symmetry equates the two
envelope expectations. The resulting remainder is
$CnL^2e^{-c\varepsilon L\log L}$, which proves
\eqref{chIV:eq:coarse-vertex-image-bound}. The bound \eqref{chIV:eq:R9.8} tends to zero for
small enough fixed $\varepsilon$; the exponentially small count-transfer remainder
is summable even after its factor $n$.
\end{proof}

For $\eta>0$, Markov's inequality on $\mathcal G_N$ now yields
\[
 \mathcal Q_I\{|\mathcal O_{N,u_0}(t)|>\eta\sqrt L\}
 \le\mathcal Q_I(\mathcal G_N^c)
   +\frac{\mathcal Q_I[|\mathcal O_{N,u_0}(t)|;\mathcal G_N]}{\eta\sqrt L}=o(1).
\]
Equivalently, under the normalized root-averaged physical-loop window law,
\begin{equation}\label{chIV:eq:R9.9}
 |\mathcal O_{N,u_0}(t)|/\sqrt L\longrightarrow0\quad\text{in probability}.
\end{equation}
The coarea argument is applied label-by-label to deterministic vertex
labels, using image-measurable probe majorants.

\appendixchapter{P}{Proof of Theorem~\ref{chIV:thm:defects}}\label{chIV:sec:R10}

Let $b_N=gN^2L$, with $s_N$ as fixed in the main text. For fixed
$0\le u_0,W<\infty$, define
\begin{equation*}
U_{u_0}(t)=\{x:\ell_t(x)\le u_0\},
\end{equation*}
\begin{equation*}
\begin{split}
D_N^{\rm wide}(U_{u_0}(t))
&=\bigl\{x\in U_{u_0}(t):\exists z\in U_{u_0}(t),\\
&\qquad\tfrac12e^{L^{1/8}}<d(x,z)<2Ne^{-L^{1/8}}\bigr\}.
\end{split}
\end{equation*}
This is the wide-defect set at softness level $u_0$.
The distance $d$ and the wide-defect set are those of
\eqref{chIV:eq:wide-def}. We prove
\begin{equation}
\lim_{N\to\infty}\sup_{a\in V_N}
\mathbb P_a\left(\sup_{|c|\le W}
\frac{|D_N^{\rm wide}(U_{u_0}(s_N+cb_N))|}{\sqrt L}>\varepsilon\right)=0
\quad(\varepsilon>0).
\label{chIV:eq:R10.1}
\end{equation}

\subsection{Summation of repelled witness pairs}
Work first in the Appendix~M root-averaged loop integral.
By \eqref{chIV:eq:R9.9}, $|\mathcal O_{N,u_0}(t)|/\sqrt L\to0$ in probability.
Appendix~N excludes a soft vertex with count above $10L$ at any coarse scale with probability $o(1)$;
the one-point bridge bound excludes root-near soft vertices.
For a wide-distance pair $(x,z)$, choose $q=\lceil\log|x-z|\rceil$, shifted by a bounded amount if needed
to have $z\in B(x,e^q)$. Then
\begin{equation*}
L^{1/8}-O(1)\le q\le L-L^{1/8}+O(1),
\qquad \log|x-z|\ge q-O(1).
\end{equation*}
which lies in $[L^{1/10},L-L^{1/10}]$.
If $x\notin\mathcal O_{N,u_0}(t)$, write
\begin{equation*}
Y_q(x)=w=\sqrt2q+s,\qquad L^{8\gamma}\le s\le10L.
\end{equation*}
The two-target estimate follows from the endpointwise calculation of
\cite[Lemma A.8 and Proposition 3.7]{chIV:LS24}, with annulus coefficients
$1,2$. Its separation hypothesis holds because
$\log|x-z|\ge q-\delta q^\gamma$. For a fixed feasible $q$-word with $Y_q=w$, let $p_{x,z}(w)$ be
the probability, under its product endpoint kernel, that the summed
block local times at both $x$ and $z$ are at most $u_0$. The endpoints
are suppressed from the notation. Uniformly over these words,
\begin{equation}
p_{x,z}(w)\le
\exp\{-4q-4\sqrt2s+Cq^\gamma\},
\label{chIV:eq:R10.2}
\end{equation}
with LS threshold $u_0/g$ and the favorable negative square term dropped.

The coefficient follows directly from the two-site Green matrix.
For one block with endpoints $a,b$, write $z_1=x,z_2=z$,
$D=B(x,e^{q+2q^\gamma})$, and $h_b(y)=\PP_y(X_{\tau_D}=b)$.
In this calculation put $h=q^\gamma$ and
\[
 K_{ij}=\frac{G_D(z_i,z_j)h_b(z_j)}{g h_b(z_i)},\qquad
 c_j=\frac{G_D(a,z_j)h_b(z_j)}{g h_b(a)},\qquad i,j\in\{1,2\}.
\]
The point separation and Appendix~L give, uniformly in the endpoints,
\[
 K_{ii}=q+O(h),\quad 0\le K_{12},K_{21}\le Ch,\quad
 c_1+c_2\ge2h(1-Ch/q).
\]
For a sufficiently large fixed $C$, $\det K>0$ and
\[
 K^{-1}\mathbf1\ge e\mathbf1,\qquad e=(q+Ch)^{-1}.
\]
For example its first row sum is $(K_{22}-K_{12})/\det K$;
the displayed Green bounds give the claimed lower bound, and
the second row is identical with the indices interchanged.
The off-diagonal entries of $K^{-1}$ are nonpositive. Resolvent
comparison, or direct inversion of its two-by-two matrix, gives
$(I+\theta K)^{-1}\mathbf1\ge e\mathbf1/(e+\theta)$ for $\theta\ge0$.
The exact occupation transform is
\[
 E^{a,b}\exp\{-\theta[\ell^{(i)}(x)+\ell^{(i)}(z)]/g\}
 =1-\theta(c_1,c_2)(I+\theta K)^{-1}\mathbf1
 \le1-p\frac{\theta}{\lambda+\theta},
\]
where $p=(2h/q)(1-C'h/q)\ge0$ and $\lambda=2/q$:
$p\le e(c_1+c_2)$ and $e\le\lambda$. Multiplying this joint two-site
Laplace transform over the independent blocks gives, for every $U\ge0$,
\[
 P\!\left(\frac1g\sum_i[\ell^{(i)}(x)+\ell^{(i)}(z)]\le U\right)
 \le\inf_{\theta>0}\exp\{\theta U-mp\theta/(\lambda+\theta)\}
 =\exp\{-(\sqrt{mp}-\sqrt{U\lambda})_+^2\}.
\]
Substitute $m=w^2/h$, $U=2u_0/g$ and $w=\sqrt2q+s$.
Keeping the negative quadratic term while absorbing the errors yields
\[
 p_{x,z}(w)\le
 \exp\{-4q-4\sqrt2s-(2-C_{u_0}h/q)s^2/q+C_{u_0}h\}.
\]
Since $h/q\to0$, dropping that term proves \eqref{chIV:eq:R10.2},
also for $u_0=0$.

Apply Appendix~L window transfer and the marginal count bound \eqref{chIV:eq:R7.8}.
There are at most $CN^2e^{2q}$ ordered pairs at scale $q$.
Set $d=L^{2\gamma}+C\log L$. Up to polynomial factors in $L$, the exponent is
\begin{equation}
2(L-q)-\frac{(\sqrt2(L-q)-s-d)_+^2}{L-q+C}-4\sqrt2s+Cq^\gamma.
\label{chIV:eq:R10.3}
\end{equation}
If $s+d\le\sqrt2(L-q)$, then
\begin{equation*}
\frac{(\sqrt2(L-q)-s-d)^2}{L-q+C}
\ge\frac{(\sqrt2(L-q)-s-d)^2}{L-q}-C',
\end{equation*}
and the exponent is at most
$-2\sqrt2s+2\sqrt2d+Cq^\gamma+C'\le-cs$.
If $s+d>\sqrt2(L-q)$, use $2(L-q)<\sqrt2(s+d)$ to obtain
$-3\sqrt2s+\sqrt2d+Cq^\gamma\le-cs$.
Here $s\ge L^{8\gamma}$ while $d+q^\gamma=o(L^{8\gamma})$.
The count and scale choices are polynomial in $L$, so the repelled-pair expectation is at most
$L^{C'}e^{-cL^{8\gamma}}+o(1)=o(1)$.
The remainder $e^{-cL\log L}$ remains summable even over $N^4$ pairs; Appendix~L handles duration errors.
Outside an $o(1)$ event every defect vertex belongs to $\mathcal O_{N,u_0}(t)$, hence
\begin{equation}
|D_N^{\rm wide}(U_{u_0}(t))|/\sqrt L\to0
\quad\text{in probability in each fixed central loop window}.
\label{chIV:eq:R10.4}
\end{equation}

\subsection{A short closing-bridge hitting bound}
Let $Q=L^{16\gamma}$ and $R_0=Ne^{-Q/2}$.
For $d(y,z)\ge R_0$, Green and heat-kernel bounds give
\begin{equation*}
\mathbb E_y\ell_{T+N^2}(z)\le(T+N^2)/N^2+C(Q+1),\qquad
\mathbb E_z\ell_{N^2}(z)\ge gL-C.
\end{equation*}
On $H_z\le T$, the interval of length $N^2$ after $H_z$ lies in $[0,T+N^2]$.
The strong Markov property implies
\begin{equation}
\mathbb P_y(H_z\le T)
\le C\left(T/b_N+(Q+1)/L\right).
\label{chIV:eq:R10.5}
\end{equation}
For $T\in[\delta b_N,2\delta b_N]$, fixed $\delta>0$, and both endpoints $y,y'$ at distance at least $R_0$ from $z$,
\begin{equation}
\mathbb P_y(H_z\le T\mid X_T=y')
\le C\delta+o(1).
\label{chIV:eq:R10.6}
\end{equation}
Split a hit into the first or second half. In the first half, the remaining heat kernel is
uniformly comparable to $N^{-2}$; apply \eqref{chIV:eq:R10.5}.
Reverse time for the second half and apply the same argument from $y'$.

\subsection{From integrated loops to stationary ordinary paths}
Write $\mathbb P_{\rm loop,t}=n^{-1}\sum_r\mathbb P_t^{r,r}$ for
the fixed-time bridge law with a uniform root. Fix $t_0=s_N+cb_N$.
The probability of a soft vertex within $R_0$ of $X_0$ is at most
$CL^{3/2}e^{-Q}=o(1)$ by the one-point bound and a vertex sum.
Stationary time reversal gives the same estimate at $X_{t_0}$.
For sufficiently small fixed $\delta>0$ and $t\in[t_0+\delta b_N,t_0+2\delta b_N]$,
the loop prefix density relative to a stationary ordinary prefix is
\begin{equation*}
\frac{p_{t-t_0}(X_{t_0},X_0)}{p_t(X_0,X_0)}=1+o(1)
\end{equation*}
uniformly. Its remainder, conditional on the prefix, is a closing bridge from $X_{t_0}$ to $X_0$.
If the prefix has at least $\lceil\varepsilon\sqrt L\rceil$ defect vertices,
select exactly this many distinct vertices and one soft witness for each.
By \eqref{chIV:eq:R10.6}, each pair is hit with conditional probability at most $C\delta+o(1)$.
Linearity of conditional expectation bounds the expected lost fraction by $C\delta+o(1)$.
Choose $\delta$ small; conditional Markov gives probability at least $1/2$ that at least half the pairs survive.
Their occupation times and distances are unchanged, so
\begin{equation*}
\mathbb P_{\rm loop,t}\left(|D_N^{\rm wide}(U_{u_0}(t))|\ge\varepsilon\sqrt L/3\right)
\ge\tfrac12\mathbb P_{\pi_N}\left(|D_N^{\rm wide}(U_{u_0}(t_0))|\ge\varepsilon\sqrt L\right)-o(1).
\end{equation*}
Integrate over the stated central window and use \eqref{chIV:eq:R10.4}.
This proves the deterministic-$t_0$ stationary result.

\subsection{Arbitrary starting points and the whole central interval}
With $C_{\rm mix}$ as in Theorem~\ref{chIV:thm:defects}, set
$\Delta_N=C_{\rm mix}N^2\log L=o(b_N)$.
The distribution of $X_{\Delta_N}$ approaches $\pi_N$ uniformly in the starting point.
Deleting the initial segment only enlarges the soft set; the defect set is increasing in the soft set.
Hence the original defect count at $t$ is at most that of the suffix at $t-\Delta_N$.
For $|c|\le W$ and large $N$, $t-\Delta_N\ge s_N-(W+1)b_N$.
Apply the stationary result and the mixing error.
Finally, both soft and defect sets decrease with time, so the supremum over $|c|\le W$
is bounded by the earliest-time count. This proves \eqref{chIV:eq:R10.1}.

\appendixchapter{Q}{Proofs of Lemma~\ref{chIV:lem:shell-count} and the contact estimate \eqref{chIV:eq:near-contact}}

The first subsection gives the one-point bound used to exclude
root-near soft vertices in Proposition~\ref{chIV:prop:deviation}.
The second proves the fixed-root shell comparison needed for
Lemma~\ref{chIV:lem:shell-count} and the averaged contact estimate
\eqref{chIV:eq:near-contact}. The last records the small-radius
perturbation of Proposition~\ref{chIV:prop:local} used for the host
records in Appendix~Z.

\subsection{The uniform soft one-point bound}
For every fixed $0\le u,W<\infty$, uniformly in $a,x$ and $|c|\le W$,
\[
 \PP_a(\ell_{s_N(c)}(x)\le u)\le Cn^{-1}L^{3/2},\qquad
 \PP_{s_N(c)}^{a,a}(\ell_{s_N(c)}(x)\le u)\le Cn^{-1}L^{3/2}.
\]
To prove this, fix $\theta>0$ and let
$H_{x,\theta}=-\mathcal L_N+\theta\one_{\{x\}}$, where
$\one_{\{x\}}$ here denotes multiplication by the indicator of $\{x\}$.
Write $\lambda_k=\lambda_N(k)$ for the Fourier eigenvalues of
$-\mathcal L_N$ from Chapter~I. The rank-one resolvent identity gives
the lowest eigenvalue $\lambda$ of $H_{x,\theta}$ through
\[
 \frac1{n\lambda}=\theta^{-1}
          +\frac1n\sum_{k\ne0}(\lambda_k-\lambda)^{-1}.
\]
The Fourier estimates of Chapter~I give
$n^{-1}\sum_{k\ne0}\lambda_k^{-1}=gL+O(1)$ and
$n^{-1}\sum_{k\ne0}\lambda_k^{-2}\le Cn$.
The secular equation first gives $\lambda\le C/(nL)$, below
half the spectral gap for large $N$, and then
$\lambda=[n(gL+O(1))]^{-1}$.
With the same fixed mixing constant, set $\Delta=C_{\rm mix}n\log L$. Positivity and domination by the ordinary
semigroup imply
\[
 e^{-tH_{x,\theta}}(a,b)\le
 \|e^{-\Delta H_{x,\theta}}\delta_a\|_2
 \|e^{-\Delta H_{x,\theta}}\delta_b\|_2 e^{-\lambda(t-2\Delta)}
 \le \frac Cn e^{-\lambda(t-2\Delta)}.
\]
The last inequality follows from
$\|e^{-\Delta H_{x,\theta}}\delta_a\|_2^2\le p_{2\Delta}(a,a)\le C/n$.
At $t=s_N(c)$ this is at most $Cn^{-2}L^{3/2}$.
Feynman--Kac and exponential Markov, followed either by summing
over $b$ or by dividing the diagonal kernel by $p_t(a,a)\asymp1/n$,
give the two assertions. Integrating gives their loop-window
versions. In particular a disk of radius $Ne^{-Q/2}$ around the
root contains a soft point with probability at most
$CL^{3/2}e^{-Q}=o(1)$. Stationary reversal supplies the same
bound at the other endpoint of an ordinary prefix.

\subsection{Near-root contact removal with the original normalization}
At each fixed root $r$, the source field and path are integrated
under the product measure
\[
 \mathsf L_r(dt,d\omega)\,\gamma_r(dh),\qquad
 \gamma_r=\mathcal L(h^r).
\]
After coarea, the field--branch law is $\nu_{N,s}^r$; the
root-summed image measure is displayed below.
The trace bound gives
$\sup_rP(\|h^r\|_2^2>D_0nL)\le C/D_0$. Under the source measure this
tail costs at most $4WC/D_0$. On its complement, the coarea image
time is in the enlarged interval
$[s_N-Wb_N,s_N+Wb_N+D_0nL/2]$. The contact event maps into
\[
 \{\exists x:\ d_1(x,r)\le5R_N,\ (\phi(x)+\sigma a_s)^2\le v\}.
\]
The root-summed image measure on this interval is
\[
 n\varrho_N(s)\,ds\left(\frac1n\sum_r
                  \delta_r\otimes\nu_{N,s}^r\right).
\]
Here the laws on the right are on $(r,\phi,\sigma)$, and the
square event is evaluated at $W_s^\sigma(\phi)$. The root weight
is retained inside this joint law. On every compact shell interval
its comparison with a uniform independent root and the positive
reference branch has error $o(1)$.
To verify this also for a fixed root, use the exact densities from
Chapter~III,
\[
 n b_{r,s}^\sigma=\frac{(\phi(r)+\sigma a_s)_+}{a_s},
 \qquad a_s=\sqrt{B-\mathcal{Q}_N},\quad B=2s/n,\quad \mathcal{Q}_N=\|\phi\|_2^2/n.
\]
On $\mathcal{Q}_N\le B/2$, the sum of the absolute errors from the density
one on the positive branch and zero on the negative branch is
at most $2|\phi(r)|/a_s$. Its expectation is $O(L^{-1/2})$.
For the edge $B/2<\mathcal{Q}_N<B$, translation invariance and conditional
Jensen give $E(|\phi(r)|\mid \mathcal{Q}_N)\le\sqrt{\mathcal{Q}_N}$.
The energy-density estimate of Chapter~II therefore bounds the
edge error by
\[
 C e^{-c(B/2-G_N(r,r))}
 \int_{B/2}^B\left(1+\sqrt{B/(B-q)}\right)dq=o(L^{-1/2}).
\]
The omitted tail $\mathcal{Q}_N\ge B$ has the same bound.
Normalization by $n\varrho_{N,r}(s)=1+O(L^{-1/2})$ preserves
this estimate. With total variation defined by $\|\mu-\nu\|_{\rm TV}=\sup_E|\mu(E)-\nu(E)|$,
the resulting statement is
\[
 \sup_{r,\,|\theta|\le W'}
 \|\nu_{N,s}^r-\mathcal L(\phi_N)\otimes\delta_{+1}\|_{\rm TV}
 \le C_{W'}L^{-1/2},\qquad s=s_N+nv_N+\theta b_N,
\]
for every fixed $W'<\infty$. It proves the fixed-root comparison
used in Lemma~\ref{chIV:lem:shell-count}, and, after averaging the
root, the joint comparison required here. By sign symmetry,
either deterministic reference sign gives the same square-test law.

For a uniform independent root, condition first on the field.
On the event that the number of image low coordinates is at most
$K$, the union of their $5R_N$ neighbourhoods occupies at most
$CKR_N^2/n=CKe^{-4L^{1/8}}$ of the torus. The complementary
event has arbitrarily small probability by
Lemma~\ref{chIV:lem:shell-count}. Choose $K$, then take $N$ large.
Multiplication by the bounded image density and the enlarged
window length proves that the source contact integral is
$o(1)+4WC/D_0$. Finally let $D_0\to\infty$.
Disintegration over the original Gaussian root coordinate proves
\eqref{chIV:eq:near-contact}; restricting the root integral only
decreases it.

\subsection{Complete host records with the original radius shrinkage}
Some later resampling steps use the host radii
\begin{equation}\label{chIV:eq:host-radii}
 r_{\rm in}=e^{q+h},\qquad
 r_{\rm out}=e^{q+2h-4e^{-h}},\qquad h=q^{1/8192}.
\end{equation}
The harmonic average is still taken at radius $e^{q+1}$. Keep
every return until the specified outer exit. For every count $m$
and every feasible entrance--exit word, with
$2(r_{\rm out}+1)<N$, let $\ell^{(i)}$ be the local-time profile
of record $i$ and put
\[
 S_m=\frac1g\sum_{i=1}^m\sum_z\mu_q(z)\ell^{(i)}(z).
\]
The probabilities and expectations below refer to the product of the
specified endpoint kernels. The Green calculation of Appendix~M
applied to these exact radii gives
\[
 |ES_m-mh|\le C e^{-ch}mh,\qquad
 P(|S_m-ES_m|>z)\le
 2\exp\{-c\min(z^2/(mh^2),z/h)\}.
\]
Indeed the logarithmic width is $h-4e^{-h}$, and its difference
from $h$ is absorbed by the displayed relative error. The
conditioned Green factor is unchanged. If $Y=\sqrt{mh}$ and
$2Ce^{-ch}Y\le\Delta\le Y$, then
\begin{equation}\label{chIV:eq:host-lower}
 P(\sqrt{S_m}<Y-\Delta)\le2e^{-c\Delta^2/h}.
\end{equation}
To check the centering, $ES_m\ge Y^2-Y\Delta/2$, while
$(Y-\Delta)^2\le Y^2-Y\Delta$. Thus a square-root deficit
forces a centered deficit of at least $Y\Delta/2$; substitute
this into Bernstein. The case $m=0$ is empty. In particular the
centre in \eqref{chIV:eq:host-lower} is $\sqrt{mh}$, with the width
perturbation absorbed in the displayed relative mean error.

The complete-count and duration guards used with these hosts
are taken on the original sample. On the ordinary infinite path,
conditional on their realized endpoint word, the first $M$ complete
records have the endpoint product law of
Chapter~III. The point-hit estimate of Appendix~L and
exponential Markov bound excessive counts at a soft point by
$\exp\{u-cMh/q\}$. The duration bound is obtained before
conditioning on the soft event, using the original block law;
when a replacement is conditioned on softness, its avoidance
denominator is retained as in \eqref{chIV:eq:R6.6}. Use duration cutoff
$e^{L/2}$ and put
\[
 M=\lceil e^{L/8}\rceil,\qquad J=\lceil e^{4L}\rceil.
\]
For $q\le L^{1/8}+10$, the count bound at a soft core point and
the duration bound for at most $M$ records are each at most
$\exp\{-e^{L/16}\}$ for large $N$. The logarithms of the endpoint
denominators and the outer radius are $O(q+h)$, whereas
$e^{L/2}/r_{\rm out}^2\ge e^{L/4}$ eventually. The total central time
is at most $e^{3L}$, so the rate-one Poisson clock satisfies
$P(\text{at least }J\text{ jumps})\le e^{-e^{4L}/2}$.
On its complement, at most $2J+1$ chronological word lengths are
possible. Summing over these lengths, the hosts and core points
costs at most $(2J+1)n^2$. Thus either failure has ordinary-path
budget at most
\[
 e^{-e^{4L}/2}+(2J+1)n^2e^{-e^{L/16}}\le e^{-20L}.
\]
The original root-window normalization costs at most $2n$, by
the lower heat-kernel bound at times at least $\delta b_N$.
It follows that each failure has normalized probability at most
$2n e^{-20L}\le N^{-17}$ for large $N$.
For the fixed-host guard the resulting bound for each of its two
failure types and each colour is $N^{-17}$, uniformly over the
single chosen grid origin and admissible root windows. The
parameters are chosen at fixed $q$ before sending $N$ to infinity.

The exact union at this step has $15377$ geometric
colours and $k(q)$ locator colours, where $k(q)<\infty$ for each
fixed host scale. Its probability is bounded by
$2\cdot15377\,k(q)N^{-17}$, hence tends to zero for every such
$q$. One then takes the host scale large. Each endpointwise
estimate is integrated against the law of the observed endpoint
word. An independent random grid origin can be averaged afterwards
by the same uniform bound. Window resampling retains its original
duration displacement and normalization, as in Chapter~III.


\par\endgroup

\begingroup
\mergedchaptersetup{V}{The limiting law of local-time activity}
\appendixchapter{R}{Proof of the marking estimate \eqref{chV:eq:marking}}

Fix the occupation profile, root and conditional root value.
The corresponding canonical parents are deterministic; contact
events and the square set remain random under the Gaussian law.
We prove \eqref{chV:eq:marking} using the relative finite-cluster
contact estimate and Green covariance bounds from Chapters~I--II,
retaining the full supports of the parents. Unless another law is
specified, probabilities and expectations are over $\phi^{(r)}$
and its indicated transforms.
In particular $h^{r,z}$ has the shifted Gaussian law in
\eqref{chV:eq:conditional-field}.

\subsection{The Gaussian block comparison}

Let $D_1,\ldots,D_m$ be canonical parents whose representatives
are at graph distance greater than $4R_N$ from the root, with
$m\le B\sqrt L$. The conditional field is \eqref{chV:eq:conditional-field}.
For each block take its canonical representative $a_i$ as anchor,
and start from the vector
\[
 \left(h^{r,z}(a_i),\ (h^{r,z}(x)-h^{r,z}(a_i))_{x\in D_i\setminus\{a_i\}}\right).
\]
Here is an explicit invertible change of coordinates for this
regression. Put $\psi=\phi^{(r)}$, let $B_i$ take the differences
from $a_i$ on $D_i$, and use the mass-one regression weights
$e_i=e_{D_i}^r$ from \eqref{chII:eq:R15.11}. Define
\[
 \begin{gathered}
 Z_i=e_i^{\mathsf T}\psi|_{D_i},\qquad
 V_i=\Var Z_i,\qquad \Sigma_i=\Cov(B_i\psi|_{D_i}),\\
 \sum_{x\in D_i}e_i(x)=1,\quad \|e_i\|_1\le2,\quad
 \Cov(Z_i,B_i\psi|_{D_i})=0,\quad
 V_i\ge(g/2)L,\quad\Sigma_i\ge I/4.
 \end{gathered}
\]
These bounds follow from \eqref{chII:eq:R15.10}--\eqref{chII:eq:R15.14};
the distance hypotheses there allow the fixed-factor changes between
$d$ and $d_1$. The centered coordinates
$(Z_i/\sqrt{V_i},\Sigma_i^{-1/2}B_i\psi|_{D_i})$ have identity
covariance. Their affine counterparts for $h^{r,z}$ retain its
transformed mean. The inverse recovers all coordinates on $D_i$
because the common row has total mass one; for a singleton the
increment space is empty.

Here is the quantitative consequence of the Green difference
estimates. Put $d_N=e^{L^{1/6}}$ and let $\Sigma_{ij}$ be a
cross-covariance matrix in these normalized coordinates. Write
$\|A\|_{\rm HS}^2=\sum_{i,j}|A_{ij}|^2$ for the Hilbert--Schmidt norm
and $\|A\|_{\rm op}=\sup_{\|v\|_2=1}\|Av\|_2$ for the Euclidean
operator norm. For $i\ne j$,
\begin{align}\label{chV:eq:cross-rate}
 \|\Sigma_{ij}\|_{\rm HS}^2\le C\bigg\{
 &\frac{(1+2L^{1/8})^2}{L^2}
 +\frac{d_N^2}{L}\left(\frac{d_N}{R_N}\right)^2
 +d_N^4\left(\frac{d_N}{R_N}\right)^4\bigg\}.
\end{align}
To see how the dimension enters, separate a cross-covariance
matrix into its common--common, common--increment, and
increment--increment entries. There is one of the first kind,
at most $9d_N^2$ in each mixed row, and at most $(9d_N^2)^2$
of the last kind. Before increment whitening, with each common
variance already normalized, the corresponding entry bounds are
$C(1+2L^{1/8})/L$,
$C(d_N/R_N)/\sqrt L$, and $C(d_N/R_N)^2$.
The common rows use $\|e_i\|_1\le2$; the root correction is the
rank-one term $-G_N(x,r)G_N(r,y)/v_N$ in every covariance.
Before increment whitening, squaring and summing gives the right
side of \eqref{chV:eq:cross-rate}. Since
$\|\Sigma_i^{-1/2}\|_{\rm op}\le2$, the inequality
$\|ACB\|_{\rm HS}\le\|A\|_{\rm op}\|C\|_{\rm HS}\|B\|_{\rm op}$
retains that Hilbert--Schmidt bound after whitening, uniformly in the
block dimensions.
For fixed supports $D_i$, the covariance depends on $r$ but
not on the occupation values or $z$; the latter change tests and
means. If the occupation profile changes the supports, the
covariance can change, but \eqref{chV:eq:cross-rate} remains uniform
over every family satisfying the stated geometric and count bounds.

After whitening each diagonal block exactly, the joint covariance
is $I+E$, with zero diagonal blocks in $E$. Summing
\eqref{chV:eq:cross-rate} over at most $B^2L$ ordered pairs gives
\begin{equation}\label{chV:eq:hs-decay}
 \|E\|_{\rm HS}^2\le C B^2\left\{
 L^{-3/4}+e^{-2L+4L^{1/6}+4L^{1/8}}
       +Le^{-4L+8L^{1/6}+8L^{1/8}}\right\}=o(1).
\end{equation}
This is why the full growing supports can be used.

Compare this joint law with the product of its exact block marginals. They have the same mean. Once
$\|E\|_{\rm op}\le1/2$, the relative entropy is
\[
 \frac12\{\operatorname{tr}E-\log\det(I+E)\}
 \le C\|E\|_{\rm HS}^2.
\]
The inequality follows by diagonalizing $E$ and using
$0\le t-\log(1+t)\le Ct^2$ for $|t|\le1/2$.
Pinsker's inequality bounds the total variation distance by
$C\|E\|_{\rm HS}=o(1)$. It follows that for any block tests
$0\le T_i\le1$ measurable with respect to block $i$,
\begin{equation}\label{chV:eq:product-comparison}
 \left|\EE\prod_{i=1}^mT_i-\prod_{i=1}^m\EE T_i\right|=o(1),
\end{equation}
uniformly in all these tests and means. Empty families satisfy
the same statement with both products equal to one.

\subsection{The relative contact estimate and the mass cap}

For a nonnegative occupation profile $\ell$ and fixed $|z|\le M$,
let $E_i$ be the event that
$2\ell(x)+h^{r,z}(x)^2\le1/4$ somewhere in $D_i$. Write
\[
 p_i=\PP(E_i),\qquad
 \alpha_i=\frac{a^G_{D_i}(\ell;r,1/4)}{\sqrt L}.
\]
The relative contact estimate of Chapter~II supplies, for every
$\tau>0$, eventually and uniformly over these parents,
\begin{equation}\label{chV:eq:relative-contact}
 |p_i-e^{-z^2/2}\alpha_i|\le\tau\alpha_i,
 \qquad \max_i p_i\le C_M L^{-1/3}.
\end{equation}
In particular the estimate also covers $\alpha_i=0$. Use
$T_i=1-\lambda\ind_{E_i}$ in \eqref{chV:eq:product-comparison}.
On the mass cap $\sum_i\alpha_i\le H$, and with $\max_i p_i\le t<1$,
\begin{align}
 \left|\lambda\sum_i p_i-\lambda e^{-z^2/2}\sum_i\alpha_i\right|
     &\le\tau H,\label{chV:eq:sum-error}\\
 0\le\sum_i[-\log(1-\lambda p_i)-\lambda p_i]
     &\le\frac{t}{2(1-t)}\sum_i p_i
       \le\frac{t(1+\tau)H}{2(1-t)}.\label{chV:eq:bernoulli-error}
\end{align}
The scalar inequality used here follows by summing
$\sum_{k\ge2}x^k/k\le x^2/(2(1-x))$.
Since $e^{-x}$ is one-Lipschitz on $[0,\infty)$, these bounds
and \eqref{chV:eq:product-comparison} imply
\begin{equation}\label{chV:eq:far-marking}
 \EE\prod_i(1-\lambda\ind_{E_i})
       =\exp\{-\lambda e^{-z^2/2}\sum_i\alpha_i\}+o(1),
\end{equation}
uniformly in $0<\lambda<1$, with $B,H,M$ fixed. Choose
$\tau,t$ small after $H$ is fixed, and then take $N$ large enough
for both the product comparison at count cap $B$ and the
relative-contact bounds; these choices are independent of $\lambda$.

\subsection{Restoring near-root parents}

For a fixed path at time $t$, take $\ell=\ell_t$ and, for
$a\in\operatorname{Rep}_N(U_1(t))$, put
\[
 E_a=\{\exists x\in D_a(U_1(t)):
             2\ell_t(x)+h^{r,z}(x)^2\le1/4\}.
\]
These are the preceding contact events indexed by the canonical representatives.
Representatives at distance at most $4R_N$ from $r$ form an
$R_N$-separated set. A grid packing of this graph ball gives
at most $625$ such representatives, for all sufficiently large
$N$. Each of their parents lies in the ball of radius $5R_N$
about $r$.

For these near-root parents we prove
\[
 \sum_{\substack{a\in\operatorname{Rep}_N(U_1(t))\\d_1(a,r)\le4R_N}}
 \frac{a^G_{D_a(U_1(t))}(\ell_t;r,1/4)}{\sqrt L}
 \le625\,C L^{-1/3}=o(1),
\]
uniformly in the occupation profile.
Indeed a union of the intervals in \eqref{chV:eq:cluster-activity}
has length at most $1$ plus the range of the field on the parent.
For a parent $D=\{x_1,\ldots,x_m\}$ with anchor $a\in D$, put
$Y_i=\phi^{(r)}(x_i)-\phi^{(r)}(a)$. These variables are centered
Gaussian, and subtracting the anchor does not change the range.
The Green bound and the fact that Gaussian conditioning decreases
increment variance give $\operatorname{Var}Y_i\le\sigma^2:=C\log(2+\rho_N)$
for a fixed constant $C$. There are $m\le Cd_N^2$ sites. The union bound
$\PP(\max_i|Y_i|>t)\le2m e^{-t^2/(2\sigma^2)}$, integrated
above $\sigma\sqrt{2\log(2m)}$, gives
\[
 \EE\max_i|Y_i|\le C\sigma\sqrt{\log(2m)}=O(L^{1/6}).
\]
Thus each parent has normalized activity at most $CL^{-1/3}$.
Summing over at most $625$ parents proves the displayed bound.

For the near-root indices,
\[
 \begin{aligned}
 &\left|\prod_{a\in\operatorname{Rep}_N(U_1(t))}(1-\lambda\ind_{E_a})
   -\prod_{\substack{a\in\operatorname{Rep}_N(U_1(t))\\d_1(a,r)>4R_N}}
                     (1-\lambda\ind_{E_a})\right|\\
 &\qquad\le\ind_{\bigcup_{a:\,d_1(a,r)\le4R_N}E_a},\qquad
 \PP\left(\bigcup_{a:\,d_1(a,r)\le4R_N}E_a\right)\le q_N(r,t,z).
 \end{aligned}
\]
Restoring their activity in the exponential changes it by at most
their total normalized activity, hence by $o(1)$. It follows
from \eqref{chV:eq:far-marking} that, on the count and mass caps,
\begin{equation}\label{chV:eq:full-marking}
 \left|\EE\prod_{a\in\operatorname{Rep}_N(U_1(t))}
          (1-\lambda\ind_{E_a})
        -e^{-\lambda e^{-z^2/2}A_N(t)}\right|
 \le o(1)+q_N(r,t,z).
\end{equation}

\subsection{Matching the square set to its parents}

Let $T=\{x:2\ell_t(x)+h^{r,z}(x)^2\le1/4\}$ and
$U=U_1(t)$. Then $T\subset U$. Since
$D_N^{\rm can}(U)\subset D_N^{\rm wide}(U)$, we have
\[
 T\cap D_N^{\rm wide}(U)=\varnothing
 \quad\Longrightarrow\quad
 \bigl(T\subset U^\circ\ \text{and}\ T=T^\circ\bigr).
\]
For the second assertion, any pair in $T$ at graph distance strictly
between $\rho_N$ and $R_N$ would have both endpoints in
$D_N^{\rm wide}(U)$. Indeed, periodic Euclidean distance lies between
$d_1/\sqrt2$ and $d_1$, and the canonical annulus is contained in
the wider Euclidean annulus for large $N$.
Within each $D_a(U)$, the vertices of $T$ are at distance at most
$\rho_N$ and form one component; vertices in different parents are
separated by at least $R_N$. Therefore
\[
 n_N(T)=\sum_{a\in\operatorname{Rep}_N(U)}\ind_{E_a}
 \quad\text{whenever }T\cap D_N^{\rm wide}(U)=\varnothing.
\]
Since both tests lie in $[0,1]$,
\[
 \left|\Psi_{N,\lambda}(2\ell_t+(h^{r,z})^2)
       -\prod_{a\in\operatorname{Rep}_N(U)}(1-\lambda\ind_{E_a})\right|
 \le\ind_{\{T\cap D_N^{\rm wide}(U)\ne\varnothing\}}.
\]

For a site outside the $5R_N$ root neighborhood, the conditional
variance in \eqref{chV:eq:conditional-field} is at least $cL$ by the
Green bounds. Its density is therefore at most $C/\sqrt L$,
uniformly in its mean. Its contact interval has length at most
one, so a union bound gives
\begin{equation}\label{chV:eq:defect-contact}
 \PP(T\cap D_N^{\rm wide}(U)\ne\varnothing)
 \le q_N(r,t,z)+\min\{1,C_M|D_N^{\rm wide}(U)|/\sqrt L\}.
\end{equation}
For the near part, \eqref{chV:eq:near-event} gives
\[
 \PP\{\exists x\in T\cap D_N^{\rm wide}(U):d_1(x,r)\le5R_N\}
 \le q_N(r,t,z).
\]
For the far part, the union bound and the probability bound $1$
give the minimum in \eqref{chV:eq:defect-contact}.
Combine \eqref{chV:eq:defect-contact} with
\eqref{chV:eq:full-marking} to obtain \eqref{chV:eq:marking}.

The argument also explains why the defect count is truncated
before integration and why a near-root contact probability
appears twice.

\appendixchapter{S}{Proof of Proposition~\ref{chV:prop:master}: coarea error estimates}

We supply the error estimates in the proof of
Proposition~\ref{chV:prop:master}: the passage from
\eqref{chV:eq:exact-source} to \eqref{chV:eq:marked-loop},
removal of the endpoint factor and extension of the time integral.
Fix $\chi\in C_c(\mathbb R)$ and choose $H_0<\infty$ with
$\operatorname{supp}\chi\subset[-H_0,H_0]$. Constants may depend on this support.

\subsection{Field energy and the initial time interval}

Diagonalize $G_N$ on the mean-zero subspace. In this subsection
$\lambda_j$ denote its $n-1$ positive eigenvalues, namely the reciprocals
of the nonzero eigenvalues $\lambda_N(k)$ of $-\mathcal L_N$. Then
$\|\phi\|_2^2=\sum_j\lambda_j\xi_j^2$ for independent standard
normal coordinates $\xi_j\sim\mathcal N(0,1)$. All probabilities and
expectations in this energy calculation are over these Gaussian coordinates. The estimates from Chapter~I imply
\begin{equation}\label{chV:eq:spectral-energy}
 \max_j\lambda_j\le Cn,\qquad
 \sum_j\lambda_j=nv_N,\qquad
 \sum_j\lambda_j^2\le Cn^2.
\end{equation}
Thus $\EE \mathcal{Q}_N=v_N$ and $\Var(\mathcal{Q}_N)=2n^{-2}\sum_j\lambda_j^2\le C$.
Cauchy--Schwarz gives the energy estimate used in
\eqref{chV:eq:energy-error}.

We also need a tail much further from the mean. Choose a fixed
$a>0$ small enough that $2a\lambda_j/n\le1/2$. By independence
and $-\log(1-x)\le2x$ on $[0,1/2]$,
\[
 \log\EE\exp\{a\|\phi\|_2^2/n\}
       =-\frac12\sum_j\log(1-2a\lambda_j/n)\le CL.
\]
Consequently $\PP(\|\phi\|_2^2\ge cnL^2)
\le e^{-c'L^2+CL}$. Since $\phi(r)$ has variance $v_N=O(L)$,
the identity $\|h^r\|_2^2=\|\phi\|_2^2+n\phi(r)^2$ gives,
for every fixed $c>0$ and suitable $c'>0$,
\begin{equation}\label{chV:eq:pinned-tail}
 \sup_r\PP(\|h^r\|_2^2\ge cnL^2)\le C e^{-c'L}.
\end{equation}

If $0<t\le b_N$ and the test in \eqref{chV:eq:exact-source} is
nonzero, then
\[
 \frac12\|h^r\|_2^2=s_N+nv_N+\theta b_N-t,
 \qquad |\theta|\le H_0.
\]
Because $s_N\sim2gnL^2$, the right side is at least $c nL^2$
for all sufficiently large $N$. The field and path are
independent, and the square test is at most one. Also
\begin{equation}\label{chV:eq:initial-mass}
 \frac1{b_N}\sum_r\int_0^{b_N}p_t(r,r)\,dt
 \le 1+\frac{nv_N}{b_N}=O(1).
\end{equation}
Indeed $\int_0^T(p_t(r,r)-1/n)dt\le v_N$ by the nonnegative
Fourier expansion on the diagonal. Multiplying
\eqref{chV:eq:pinned-tail} by \eqref{chV:eq:initial-mass} and
$\|\chi\|_\infty$ proves that the initial coarea contribution
vanishes. On the remaining range $t\ge b_N$, the heat-kernel bound
gives $np_t(r,r)\le2$ for all sufficiently large $N$.

\subsection{Energy replacement without decorrelation}

Translation invariance yields
\[
 \|\chi(\cdot+U_N^r)-\chi(\cdot+R(S_r))\|_1
       =\|\chi(\cdot+\epsilon_N^r)-\chi\|_1.
\]
For $\delta>0$, split according to $|\epsilon_N^r|\le\delta$
to bound its expectation by
\[
 \sup_{|h|\le\delta}\|\chi(\cdot+h)-\chi\|_1
     +2\|\chi\|_1\PP(|\epsilon_N^r|>\delta).
\]
By \eqref{chV:eq:energy-error} this is at most
\[
 \sup_{|h|\le\delta}\|\chi(\cdot+h)-\chi\|_1
   +\frac{C\|\chi\|_1}{\delta L},
\]
uniformly in the root. The error is controlled by the unconditional
$L^1$ translation bound together with the uniform energy estimate.
First let $N\to\infty$, using \eqref{chV:eq:energy-error}
uniformly in $r$, and then $\delta\downarrow0$, using continuity
of translations in $L^1$.

For each fixed field realization, bound the square test by
one and integrate the resulting absolute difference over
paths. On $t\ge b_N$, $p_t(r,r)\le2/n$, so after the change
$dt=b_Ndc$ the total error is at most the quantity in
\eqref{chV:eq:translation-error}. This proves the replacement even
though the energy error is correlated with the square test.

\subsection{Integrating the marking error}

On $|z|\le M$, the support condition gives a fixed window
$W=H_0+(M^2+1)/2$. Under \eqref{chV:eq:window}, the capped
error \eqref{chV:eq:marking}, integrated in $\gamma_1(dz)$, is at most
\begin{multline}\label{chV:eq:integrated-marking}
 \|\chi\|_\infty\bigg[
 4W\varepsilon+
 \mathfrak m_{N,W}\{n_N(U_1(t))>B\sqrt L\text{ or }A_N(t)>H\}\\
 +2\mathfrak m_{N,W}\!\left[\int q_N(r,t,z)\,\gamma_1(dz)\right]
 +\mathfrak m_{N,W}\!\left[\min\{1,C_M|D_N^{\rm wide}(U_1(t))|/\sqrt L\}\right]
 \bigg].
\end{multline}
The near-root and defect terms vanish by
\eqref{chV:eq:near-input} and \eqref{chV:eq:defect-input}, respectively.
For fixed $M$ and hence fixed $W$, choose $B,H$ by
\eqref{chV:eq:caps}, then $\varepsilon$, and finally the required
threshold in $N$ for \eqref{chV:eq:marking}.
For $|z|>M$, integrate over all $c$ after bounding $np_t\le2$;
the difference of the two tests is at most one. This part is
bounded by $2\|\chi\|_1\gamma_1\{|z|>M\}$.

One may replace the last integrated input by the corresponding
defect tail bound. For any $\eta>0$,
\[
 \begin{aligned}
 &\mathfrak m_{N,W}[\min\{1,C|D_N^{\rm wide}(U_1(t))|/\sqrt L\}]\\
 &\quad\le4WC\eta
       +\mathfrak m_{N,W}\{|D_N^{\rm wide}(U_1(t))|>\eta\sqrt L\}.
 \end{aligned}
\]
Let $N\to\infty$ at fixed $\eta$, then $\eta\downarrow0$.
Gaussian disintegration gives, for each fixed root and path,
\[
 \int q_N(r,t,z)\,\gamma_1(dz)
 =\gamma_r\{h:\exists x:d_1(x,r)\le5R_N,
                        \ 2\ell_t(x)+h(x)^2\le1/4\},
\]
where $\gamma_r=\mathcal L(h^r)$ is the pinned-field law of
Section~3.2.

\subsection{Endpoint normalization and the early-time envelope}

Let
$\eta_N=\sup_{t\ge b_N,\,r}|np_t(r,r)-1|$.
The heat-kernel bound gives $\eta_N\to0$ since $b_N/n=gL\to\infty$.
Every conditional exponential test in \eqref{chV:eq:marked-loop}
lies in $[0,1]$, so removing the endpoint factor costs at most
$\eta_N\|\chi\|_1$.

The left cutoff in the central time coordinate is $1-s_N/b_N\to-\infty$.
To extend the remaining integral from $c\ge1-s_N/b_N$ to all $c$,
bound its absolute value on the omitted part by
\begin{align*}
 &\int\gamma_1(dz)\int_{c<1-s_N/b_N}|\chi(c+R(z))|\,dc\\
 &\quad\le\|\chi\|_1\PP\{R(G)\ge-H_0-1+s_N/b_N\}\longrightarrow0.
\end{align*}
The implication used here is simply
$c+R(z)\ge-H_0$ wherever the integrand is nonzero. This
also controls the portion with nonpositive physical time,
where $F_N$ was assigned the value one.


\par\endgroup

\begingroup
\mergedchaptersetup{VI}{Fine structure and local conditional laws}
\let\mathscr\mathcal
\appendixchapter{T}{Proof of Proposition~\ref{chVI:prop:denominator}}
All kernels condition on exact entrances, exits and block counts, with variable durations.
A coarse Gaussian shift gives the critical seed; annular likelihood comparison
transfers it to every prescribed word. Added avoidance paths then give a
denominator for every count.

\subsection{Kernel units and the coarse-field critical seed}
Let $L^{\rm loc}$ be the accumulated rate-one occupation of the
local paths and put $\ell=L^{\rm loc}/g$, $g=2/\pi$. For a rate-one
record $\omega$ define
\[
 (\Theta_g\omega)(s):=\omega(gs).
\]
With the same spatial domains, entrances, exits and record count,
\[
 P_m^{(g)}=(\Theta_g^{\otimes m})_*P_m^{(1)},\qquad
 \ell=\frac{L^{\rm loc}}g,\qquad
 V_j^{(g)}=\frac{V_j^{(1)}}g.
\]
Accordingly a physical cutoff $V_j^{(1)}\le d_N$ is represented here
as $V_j^{(g)}\le d_N/g$. Write
$\mathcal L_{\rm LS}=g\mathcal L_{\rm pl}$ for the generator
of the time-changed planar walk, where $\mathcal L_{\rm pl}$ is the
rate-one nearest-neighbor generator on $\ZZ^2$. Write $P_x,E_x$
for the law and expectation of the rate-$g$ planar walk started
at $x$. The superscript
$\mathrm{LS}$ below denotes this local normalization: $G_D^{\rm LS}$
is its counting-measure Green kernel killed on leaving $D$.
Its diagonal is $\log R+O(1)$ in a disk of radius $R$, and the
Ray--Knight field has covariance $G_D^{\rm LS}/2$. Fix
\begin{equation}\label{chVI:eq:R17.1}
 \mathsf K=B(y,e^q),\qquad A=B(y,e^{q+\mathfrak h}),\qquad
 D=B(y,e^{q+2\mathfrak h}).
\end{equation}
Use the open lattice disks $B$ of Section~6.1. The estimates below
also allow the closed integer-radius disks $\mathbb B$: rounding
changes a physical radius by at most one, and hence a logarithmic
radius by $O(e^{-q})$.
Let $P_m=P_m^{(g)}$ be the product of $m$ first-exit path laws for
the rate-$g$ walk, with feasible prescribed entrances in $\partial_iA$
and exits in $\partial_oD$; write $E_m=E_m^{(g)}$ for its expectation.
All probabilities and expectations in the Gaussian seed argument
refer instead to the Gaussian fields explicitly introduced there.
For fixed $u>0$, put
\[
 S=\{\min_{\mathsf K}\ell\le u\},\quad W=\sqrt{\mathfrak h m},\quad
 m_c=\lceil2q^2/\mathfrak h\rceil.
\]
The required conclusion is
\begin{equation}\label{chVI:eq:R17.2}
 \inf_{\text{feasible words of length }m_c}P_{m_c}(S)\ge e^{-C\mathfrak h},
\end{equation}
uniformly when
\begin{equation}\label{chVI:eq:R17.3}
 q^{\gamma_-}\le\mathfrak h\le q^{\gamma_+},\qquad
 0<\gamma_-\le\gamma_+<1/2.
\end{equation}
In particular $\log q=o(\mathfrak h)$ and $\mathfrak h=o(q)$ uniformly.

Let $f_D$ be the centered Gaussian
field with covariance $G_D^{\rm LS}/2$, and take the Ray--Knight shift
\[
 a=\sqrt2(q+2\mathfrak h).
\]
It suffices to produce one $x\in\mathsf K$ for which $|f_D(x)+a|\le\sqrt u$.

Translate $y$ to the origin and choose once and for all a square
$Q_q=y+[-4e^q,4e^q]^2\cap\ZZ^2$.  For large $q$, $\mathsf K$ is a fixed-distance
fraction inside $Q_q$ and $Q_q\subset D$.  The Gaussian Markov property gives
\begin{equation}\label{chVI:eq:coarse-decomposition}
 f_D=f_Q+H\quad\text{on }Q_q,\qquad f_Q\perp H,
\end{equation}
where $f_Q$ is the zero-boundary field in $Q_q$ with the same covariance normalization.
Put $Y=H(y)$.

\begin{lemma}[Coarse mode and bounded residual]\label{chVI:lem:coarse-mode}
Uniformly in \eqref{chVI:eq:R17.3},
\begin{equation}\label{chVI:eq:coarse-mode-estimates}
 \Var Y=\mathfrak h+O(1),\qquad
 b_x:=\frac{\Cov(H(x),Y)}{\Var Y}=1+O(\mathfrak h^{-1})
 \quad(x\in\mathsf K).
\end{equation}
Writing
\[
 H(x)=b_xY+\widehat H(x),
\]
we have $Y\perp(f_Q,\widehat H)$ and
\begin{equation}\label{chVI:eq:coarse-residual-sup}
 \sup_{q,\mathfrak h}E\max_{x\in\mathsf K}|\widehat H(x)|<\infty.
\end{equation}
\end{lemma}
\begin{proof}
The binding covariance is
\[
 C_H(x,x')=\tfrac12\bigl(G_D^{\rm LS}(x,x')-G_{Q_q}^{\rm LS}(x,x')\bigr).
\]
On $\mathsf K\times\mathsf K$ the logarithmic singularities cancel.  The potential-kernel
expansion and interior Green differences, after rescaling $Q_q$ by $e^{-q}$, give
\begin{equation}\label{chVI:eq:binding-cov-regularity}
 C_H(y,y)=\mathfrak h+O(1),\qquad
 \sup_{x\in\mathsf K}|C_H(x,y)-C_H(y,y)|\le C,
\end{equation}
and, for $x,x'\in\mathsf K$,
\[
 E\bigl[(H(x)-H(x'))^2\bigr]
 \le C\min\{1,|x-x'|_2^2e^{-2q}\}.
\]
Indeed the outer-minus-inner Green kernel is discrete harmonic in either variable throughout
a fixed neighborhood of the core and its oscillation on the rescaled compact set is bounded;
the first and mixed discrete-difference bounds scale as $e^{-q}$ and $e^{-2q}$.
The first two assertions of \eqref{chVI:eq:coarse-mode-estimates} follow from
\eqref{chVI:eq:binding-cov-regularity}.

Regression makes $\widehat H$ independent of $Y$; independence from $f_Q$ follows from
\eqref{chVI:eq:coarse-decomposition}.  Moreover $\widehat H(y)=0$ and its canonical metric obeys
the same increment bound, up to a constant.  Cover $\mathsf K$ at dyadic mesh
$2^{-j}e^q$ and telescope from a nearest parent point.  There are $O(2^{2j})$ increments
at level $j$, each with standard deviation $O(2^{-j})$; the expected maximum increment is
$O(2^{-j}\sqrt{j+1})$.  The sum over $j$ is finite, including the final lattice scale.
This proves \eqref{chVI:eq:coarse-residual-sup}.
\end{proof}

\begin{lemma}[Interior fine-field minimum window]\label{chVI:lem:fine-minimum-window}
There is $C_0<\infty$ such that, for all large $q$,
\begin{equation}\label{chVI:eq:fine-minimum-window}
 P\left(\left|\min_{x\in\mathsf K}f_Q(x)+\sqrt2q\right|\le C_0\log q\right)
 \ge\frac78.
\end{equation}
\end{lemma}
\begin{proof}
By symmetry it suffices to control the maximum of $f_Q$ on the fixed-fraction interior ball
$\mathsf K$.  The Dirichlet extremal-process theorem of \cite{chII:BL}, applied to the square
$Q_q$ and then restricted to any fixed smaller square contained in $\mathsf K$, gives tightness
after the standard centering
$\sqrt2q-(3/(4\sqrt2))\log q+O(1)$ in the present LS units.  The maximum on
$\mathsf K$ lies between that on the smaller square and that on $Q_q$.  Thus its difference
from $\sqrt2q$ is $O_P(\log q)$, which yields the displayed probability after enlarging $C_0$.
\end{proof}

Enlarge $C_0$ so that Lemma \ref{chVI:lem:coarse-mode} gives
$P(\max_{\mathsf K}|\widehat H|\le C_0)\ge7/8$.  On the intersection with
\eqref{chVI:eq:fine-minimum-window}, which has probability at least $3/4$, put
$F=f_Q+\widehat H$.  Then
\begin{equation}\label{chVI:eq:F-min-window}
 \min_{\mathsf K}F=-\sqrt2q+O(\log q),
\end{equation}
and $F$ is independent of $Y$.

\begin{lemma}[One-dimensional coarse shift seed]\label{chVI:lem:critical-gaussian-seed}
There is a fixed $C_{\rm seed}<\infty$ such that, uniformly in \eqref{chVI:eq:R17.3},
\begin{equation}\label{chVI:eq:R17.4}
 P\left(\min_{x\in\mathsf K}|f_D(x)+a|\le\sqrt u\right)
 \ge p_q:=q^{-C_{\rm seed}}e^{-4\mathfrak h}.
\end{equation}
Consequently
\begin{equation}\label{chVI:eq:R17.5}
 p_q=\exp\{-4\mathfrak h+O(\log q)\}\ge e^{-5\mathfrak h}
\end{equation}
for all sufficiently large $q$.
\end{lemma}
\begin{proof}
Condition on a realization of $F$ satisfying \eqref{chVI:eq:F-min-window}.  Define
\[
 T_F(z)=\min_{x\in\mathsf K}\{F(x)+b_xz\}.
\]
For large $q$, all $b_x$ lie in
$[1-C/\mathfrak h,1+C/\mathfrak h]$.  Thus, for $z'>z$,
\[
 (1-C/\mathfrak h)(z'-z)
 \le T_F(z')-T_F(z)
 \le(1+C/\mathfrak h)(z'-z).
\]
Hence $T_F$ is a strictly increasing onto map.  Put
$z_0=-a-\min F$.  Since
$|(b_x-1)z_0|=O(1)$ uniformly in $x$ by \eqref{chVI:eq:F-min-window},
$T_F(z_0)=-a+O(1)$.  Let $z_F$ be the unique point with $T_F(z_F)=-a$.
Then
\[
 z_F=z_0+O(1)=-2\sqrt2\,\mathfrak h+O(\log q).
\]
The preimage of $[-a-\sqrt u,-a+\sqrt u]$ under $T_F$ is an interval of length at least
$2\sqrt u/(1+C/\mathfrak h)\ge c_u$.  On this interval, using
$Y\sim N(0,\mathfrak h+O(1))$,
\[
 \frac{z^2}{2\Var Y}
 =4\mathfrak h+O(\log q)+O((\log q)^2/\mathfrak h).
\]
Because $\mathfrak h\ge q^{\gamma_-}$, the last two terms and the Gaussian density prefactor
cost only a fixed power of $q$.  Thus the conditional probability is at least
$q^{-C}e^{-4\mathfrak h}$.  Multiply by the fixed positive probability of
\eqref{chVI:eq:F-min-window} and the residual-supremum event, and enlarge $C_{\rm seed}$.
\end{proof}

\subsubsection{Wired excursions, their intensity, and the exact condenser normalization}
Wire $D^c$ to $\partial$, with conductances $g/4$ on planar edges and boundary conductances
added at the wired vertex. Use actual occupation at $\partial$ as the excursion clock.
For a set $A$, write $H_A:=\inf\{t\ge0:X_t\in A\}$ and
$\tau_D:=\inf\{t\ge0:X_t\notin D\}$. An excursion is parametrized
from its initial interior vertex; its lifetime is $\zeta:=H_\partial$.
Write $\operatorname{cap}_{\rm LS}(A,D)$ for the condenser capacity
in this conductance network: it is the Dirichlet energy
$\frac12\sum_{x,z}c(x,z)(v(x)-v(z))^2$ of the function equal to
$1$ on $A$, $0$ at $\partial$, and harmonic elsewhere.

\begin{lemma}[Excursion thinning and the condenser]\label{chVI:lem:wired-thinning}
Let $M'$ count the excursions before wired local time $t'$ that hit $A$.
Each such excursion contributes its suffix from first entrance in $A$ to exit from $D$.
Then
\[
 M'\sim\operatorname{Poi}(t'\operatorname{cap}_{\rm LS}(A,D)),
\]
and, conditional on $M'$ and all entrance/exit labels, these suffixes have the product
of the corresponding first-exit kernels. Moreover,
\begin{equation}\label{chVI:eq:R17.6}
 \operatorname{cap}_{\rm LS}(A,D)^{-1}
 =\mathfrak h+O(e^{-cq}).
\end{equation}
\end{lemma}
\begin{proof}
If $c_\partial=\sum_zc(\partial,z)$, successive holding times at $\partial$ are independent
$\operatorname{Exp}(c_\partial)$ variables. Their partial sums in local-time coordinates form
a rate-$c_\partial$ Poisson process. Strong Markov gives independent excursion marks with
characteristic measure
\[
 \mathbf n_\partial(F)=\sum_zc(\partial,z)E_z[F(X_{[0,H_\partial]})].
\]
The event $\{H_A<\zeta\}$ is a measurable event of a single excursion. Its thinning intensity is
\[
 \mathbf n_\partial(H_A<\zeta)
 =\sum_zc(\partial,z)P_z(H_A<H_\partial)
 =\operatorname{cap}_{\rm LS}(A,D).
\]
The final equality is the equilibrium flux at $\partial$, equivalently the Dirichlet energy
of the function equal to $1$ on $A$ and $0$ on $D^c$. Marked-Poisson conditioning and strong
Markov at the first entrance into $A$ give the asserted conditional product. All occupation
in $\mathsf K\subset A$ is contained in these suffixes.

Here is the accuracy needed for the capacity. Translate $y$ to $0$ and let
$\mathfrak a_{\rm LS}:=a_{\rm lat}/g$, with $c_{\rm lat}$ as in Chapter~II.
Then $\mathcal L_{\rm LS}\mathfrak a_{\rm LS}=\delta_0$ and
$\mathfrak a_{\rm LS}(z)=\log|z|+c_{\rm lat}/g+O(|z|^{-2})$.
Let $v=1$ on $A$, $v=0$ on $D^c$ and $\mathcal L_{\rm LS}v=0$ on $D\setminus A$.
The inward and outward fluxes of $v$ are positive measures of equal mass
$\operatorname{cap}_{\rm LS}(A,D)$. Summation by parts gives
\[
 1=\sum_zv(z)\mathcal L_{\rm LS}\mathfrak a_{\rm LS}(z)
   =\sum_z\mathfrak a_{\rm LS}(z)\mathcal L_{\rm LS}v(z)
   =\operatorname{cap}_{\rm LS}(A,D)
     (\overline{\mathfrak a}_{\rm out}-\overline{\mathfrak a}_{\rm in}).
\]
The bars denote averages with respect to the two normalized flux
measures. Their supports have
radii $e^{q+2\mathfrak h}+O(1)$ and
$e^{q+\mathfrak h}+O(1)$, respectively. Thus their difference is
$\mathfrak h+O(e^{-q-\mathfrak h})$; the potential-kernel constant cancels.
This proves \eqref{chVI:eq:R17.6}.
\end{proof}

Take $t'=a^2=2(q+2\mathfrak h)^2$, and denote the corresponding wired occupation by $\ell'$.
Let $f,f'$ each have the law of $f_D$, with $f$ independent of the
wired walk. Ray--Knight in these units is
$\ell'+f^2\overset d=(f'+a)^2$. Therefore
\[
 P(\min_{\mathsf K}\ell'\le u)
 \ge P\left(\min_{x\in\mathsf K}|f'(x)+a|\le\sqrt u\right)
 \ge p_q
\]
by Lemma \ref{chVI:lem:critical-gaussian-seed}.  Writing
$\lambda'=t'\operatorname{cap}_{\rm LS}(A,D)$ and using \eqref{chVI:eq:R17.6},
\[
 \sqrt{\lambda'}-\sqrt{m_c}
 =2\sqrt{2\mathfrak h}\,(1+o(1)).
\]
Hence the Poisson lower-tail bound
$P(\operatorname{Poi}(\lambda)\le b)\le
\exp\{-(\sqrt\lambda-\sqrt b)^2\}$ gives
\begin{equation}\label{chVI:eq:R17.7}
 P(M'<m_c)\le\exp\{-(8+o(1))\mathfrak h\}=o(p_q).
\end{equation}

\subsubsection{Annular transition kernels and local records}
Put $B_-=B(y,e^{q+\mathfrak h/2})$. Write $H_D(a,b)=P_a(X_{\tau_D}=b)$.
After translation let $r_-,r_A,r_D$ be the largest Euclidean norms
of lattice points in $B_-,A,D$, respectively. Each differs by
$O(1)$ from its specified radius, as follows by taking a lattice
point on a coordinate axis. For each of these disks, if $r$ is
this largest norm, the disk is exactly $\{z:|z|_2\le r\}$ and
$r^2$ is an integer.
Their logarithms differ from $q+\mathfrak h/2,q+\mathfrak h,q+2\mathfrak h$ by $O(e^{-q})$.

\begin{lemma}[Annular entrance and exit comparisons]\label{chVI:lem:annular-harmonic-kernels}
For the first-exit law from $D$ conditioned on its exit label $b$, let
$\Pi^\downarrow_{a,b}$ be the law of first entrance into $B_-$ or exit from $D$,
identifying every final exit with $\dagger$; $a$ is on either one-step boundary of $A$.
Let $\Pi^\uparrow_{z,b}$ be the law of first exit from $A$, starting from
$z\in\partial_i B_-$. Uniformly in feasible endpoints,
\begin{equation}\label{chVI:eq:annular-transition-ratios}
 \begin{gathered}
 \sup_{a,a',b,b',z}
 \left|\frac{\Pi^\downarrow_{a,b}(z)}{\Pi^\downarrow_{a',b'}(z)}-1\right|
 \le\varepsilon_q,\qquad
 \sup_{z,b,b',w}
 \left|\frac{\Pi^\uparrow_{z,b}(w)}{\Pi^\uparrow_{z,b'}(w)}-1\right|
 \le\varepsilon_q,\\
 \varepsilon_q\le C e^{-c\mathfrak h},\qquad
 \Pi^\downarrow_{a,b}(\dagger)\ge1/4.
 \end{gathered}
\end{equation}
Ratios are taken on their common positive support, including $\dagger$ in the downward kernel.
\end{lemma}
\begin{proof}
For $z\in\partial_i B_-$, let
$K_z(v)=P_v(X_{H_{B_-}}=z)$ for the recurrent planar walk, and let $K_z(\infty)$ be its
harmonic measure from infinity. The precise external harmonic-measure input is
\cite{chII:GL}, Lemma A.8, equation (A.35):
\[
 K_z(v)=K_z(\infty)\bigl(1+O(r_-/|v|)\bigr),\qquad |v|\ge2r_-.
\]
It applies to the present lattice disk exactly: choose its logarithmic scale as
$\kappa=\log(r_--1/2+\zeta_{r_-})$, where
$0<\zeta_{r_-}<\sqrt{r_-^2+1}-r_-$.
Then the complement of $B_\kappa^-$ in the discretization of that lemma is precisely
$\{|v|\le r_-\}$. In its parameter $l$, take $e^l=|v|/\sqrt2$; for large $q$ its required
fixed lower bound on $l-\kappa$ is satisfied. Thus the constant is independent of
$\mathfrak h/q^{\gamma_-}$ and $\mathfrak h/q^{\gamma_+}$.

Set $p_a=P_a(H_{B_-}<\tau_D)$. Optional stopping of the potential kernel in
$D\setminus B_-$ gives
\[
 p_a=\frac{\log r_D-\log r_A+O(r_-^{-1})}{\log r_D-\log r_-}
     =\frac23+O(e^{-q}/\mathfrak h).
\]
Strong Markov at $H_{B_-}\wedge\tau_D$ yields the exact identity
\[
 P_a(H_{B_-}<\tau_D,X_{H_{B_-}}=z)
 =K_z(a)-E_a[\one_{\{\tau_D<H_{B_-}\}}K_z(X_{\tau_D})].
\]
Consequently this probability equals
$K_z(\infty)[p_a+O(r_-/r_A+r_-/r_D)]$.
Since $p_a$ is bounded away from zero, this is a relative error $O(e^{-\mathfrak h/2})$.

For $u,v\in B(0,2r_A+2)$ and $b\in\partial_oD$, the positive discrete harmonic
function $H_D(\cdot,b)$ satisfies
\[
 H_D(u,b)/H_D(v,b)=1+O(r_A/r_D)=1+O(e^{-\mathfrak h}).
\]
For each vertex $v$ on the comparison paths inside $B(0,2r_A+3)$,
the ball of radius $r_D/4$ about $v$ lies in $D$
for large $q$. Apply \cite{chVI:LL}, Theorem 6.3.8(b), equation (6.19), to the nonnegative
harmonic function $H_D(\cdot,b)$ on that ball:
\[
 |H_D(v+e,b)-H_D(v,b)|\le Cr_D^{-1}H_D(v,b),\qquad |e|_1=1.
\]
The constant is uniform over all $b\in\partial_oD$. Multiplying
the adjacent ratios along $O(r_A)$ steps
gives $\exp\{O(r_A/r_D)\}=1+O(r_A/r_D)$.
The Doob-transform formulas are
\[
 \Pi^\downarrow_{a,b}(z)
 =P_a(H_{B_-}<\tau_D,X_{H_{B_-}}=z)\frac{H_D(z,b)}{H_D(a,b)},
 \qquad
 \Pi^\uparrow_{z,b}(w)=H_A(z,w)\frac{H_D(w,b)}{H_D(z,b)}.
\]
The first formula and the preceding estimates prove the downward ratios, and also
$\Pi^\downarrow_{a,b}(\dagger)=1/3+O(e^{-c\mathfrak h})$.
In the second formula $z$ is held fixed, so $H_A(z,w)$ cancels when exit labels are compared.
This proves the upward ratios. Only the physical radius ratios and one-step rounding errors
were used; the entire range \eqref{chVI:eq:R17.3} is covered.
\end{proof}

\begin{lemma}[A rare-event likelihood bound]\label{chVI:lem:rare-likelihood}
Let $\mu,\nu$ be equivalent probability laws on the same measurable space,
$K$ a nonnegative random variable, $m\ge0$, $s>0$, $A\ge1$, and $\delta\ge0$. If
$|\log(d\nu/d\mu)|\le\delta(m+K)$ and $E_\mu e^{sK}\le A^m$,
then every event $E$ with $p=\mu(E)>0$ satisfies
\begin{equation}\label{chVI:eq:rare-likelihood}
 \nu(E)\ge p\exp\left\{-\delta\left(m+
                 \frac{m\log A+\log(1/p)}s\right)\right\}.
\end{equation}
\end{lemma}
\begin{proof}
Conditional Jensen gives
$E_\mu[K\mid E]\le s^{-1}\log E_\mu[e^{sK}\mid E]
\le s^{-1}(m\log A+\log(1/p))$.
Apply Jensen once more to
$\nu(E)/p\ge E_\mu[e^{-\delta(m+K)}\mid E]$.
\end{proof}

For each stopped path retain its downward endpoints and termination symbol,
and every complete upward bridge from $B_-$ to exit from $A$, including
holding times. Discard downward interiors, which cannot visit $\mathsf K$.
This record determines all core occupations and has common support for
every prescribed entrance/exit pair; the prescribed labels are parameters,
not random record coordinates. Given an upward bridge's endpoints, the
Doob factors telescope, so its law is independent of the eventual exit.

For $m$ paths let $K$ be the total number of upward pieces. Their likelihood
has $m+K$ downward factors and $K$ upward endpoint factors. The bridge
kernels cancel. Lemma \ref{chVI:lem:annular-harmonic-kernels}, for
$\varepsilon_q\le1/2$, therefore gives
\begin{equation}\label{chVI:eq:trace-likelihood}
 \left|\log\frac{dP_{\boldsymbol w}}{dP_{\boldsymbol w'}}\right|
 \le4\varepsilon_q(m+K),\qquad
 E_{\boldsymbol w}e^{sK}\le2^m,\qquad s=\log(7/6).
\end{equation}
The exponential moment follows because each piece count is dominated by
the number of failures before a geometric success of probability $1/4$,
whose moment at $s$ is $(1/4)/(1-(3/4)e^s)=2$.
Both bounds persist under arbitrary mixtures of words of length $m$:
the pointwise likelihood bounds are uniform in both words, and the
moment bound integrates directly.

\subsubsection{Transfer to every prescribed critical word}
Let $\mu$ be the record law of the first $m_c$ wired paths, conditional
only on $M'\ge m_c$. Lemma \ref{chVI:lem:wired-thinning} makes this a mixture
of product word laws. Nonnegative occupation and \eqref{chVI:eq:R17.7} give
\[
 \mu(S)\ge\frac{p_q-P(M'<m_c)}{P(M'\ge m_c)}\ge p_q/2.
\]
For any prescribed target word $\nu=P_{m_c}$, apply
\eqref{chVI:eq:rare-likelihood} and \eqref{chVI:eq:trace-likelihood}. Since
$\log q=o(\mathfrak h)$, the logarithmic loss is bounded by
\begin{equation}\label{chVI:eq:seed-likelihood-transfer}
 C\varepsilon_q\bigl(m_c+\log(1/\mu(S))\bigr)
 \le Ce^{-c\mathfrak h}(q^2/\mathfrak h+\mathfrak h+\log q)=o(1).
\end{equation}
Hence $P_{m_c}(S)\ge(1-o(1))p_q/2\ge e^{-C\mathfrak h}$,
which is \eqref{chVI:eq:R17.2}, uniformly in \eqref{chVI:eq:R17.3}.

\paragraph{Small radial perturbations.}
If the three logarithmic radii change by $\delta_K,\delta_A,\delta_D$ with
$\max|\delta_\bullet|\le Ce^{-c\mathfrak h}$, the condenser proof gives, precisely,
\[
 \operatorname{cap}_{\rm LS}(A,D)^{-1}
 =\mathfrak h+\delta_D-\delta_A+O(e^{-cq})
 =\mathfrak h+O(e^{-c'\mathfrak h}).
\]
At $t'=O(q^2)$ the corresponding Poisson-mean change is
$O(q^2\mathfrak h^{-2}e^{-c'\mathfrak h})=o(1)$.
The Gaussian seed exponent changes by $O(e^{-c'\mathfrak h})$, its domain parameters stay
in the same fixed class, and Lemma \ref{chVI:lem:annular-harmonic-kernels} has unchanged positive
radial gaps. Thus \eqref{chVI:eq:R17.2} persists with the same uniform type of constants.

\subsection{A denominator for all block counts}
For the rate-$g$ exact-exit law $P^{a,b}$ with
$a\in\partial_iA$ and $b\in\partial_oD$, the potential kernel and
exit Poisson-kernel ratio imply, uniformly in $x\in\mathsf K$ and feasible endpoints,
\begin{equation}
 P^{a,b}(X_s\ne x\text{ for every }0\le s<\tau_D)
 \ge \exp\{-{\mathfrak h}/q-C{\mathfrak h}^2/q^2\}.
\label{chVI:eq:R17.11}
\end{equation}
If the first $m_c$ paths create a soft point, select one deterministically and require every remaining path to avoid it.
For $m\ge m_c$ this gives
\begin{equation}
P_m(S)\ge
\exp\{2q-W^2/q-C{\mathfrak h}(1+W^2/q^2)\}.
\label{chVI:eq:R17.12}
\end{equation}
For $m<m_c$, append independent paths to reach $m_c$; the soft set only decreases, so
$P_m(S)\ge e^{-C\mathfrak h}$.
For every integer $m\ge0$, with $s=W-\sqrt2q$, these two cases give
\begin{equation}
P_m(S)\ge
\exp\{-2\sqrt2s_+-s_+^2/q-C{\mathfrak h}(1+W^2/q^2)\}.
\label{chVI:eq:R17.13}
\end{equation}
Rounding is absorbed in $C\mathfrak h$.

\subsection{One-point, two-point, and harmonic bounds}
Let $Y:=|\{x\in\mathsf K:\ell(x)\le u\}|$, and define
$H_2:=\#\{(x,x')\in\mathsf K^2:\ell(x),\ell(x')\le u,
\ |x-x'|_2\ge e^{q-C_*\mathfrak h}\}$, with $C_*>0$ fixed.
The exact-exit Laplace calculation below gives
\begin{equation}
\begin{aligned}
E_m Y&\le
\exp\{2q-W^2/q+C{\mathfrak h}W^2/q^2+CW/q+C\},\\
E_m H_2&\le
\exp\{4q-2W^2/q+C{\mathfrak h}W^2/q^2+CW/q+C\}.
\end{aligned}
\label{chVI:eq:R17.14}
\end{equation}
The first bound applies to every count; the second is used at near-critical or higher counts.
For a single vertex, the prescribed-exit block hits it with probability at least
\[
\mathfrak h/q-C\mathfrak h^2/q^2.
\]
For $F=\{x\}$ or $F=\{x,x'\}$ with
$|x-x'|\ge e^{q-C_*\mathfrak h}$, put $j=|F|\in\{1,2\}$.
In the prescribed-exit kernel with entrance $a$ and exit $b$, set
\[
 h_b(z)=P_z(X_{\tau_D}=b),\qquad
 \widehat G_b(z,z')=G_D^{\rm LS}(z,z')\frac{h_b(z')}{h_b(z)},
 \qquad G_F=(\widehat G_b(z,z'))_{z,z'\in F}.
\]
The potential and Poisson-kernel estimates give, uniformly,
\[
 \widehat G_b(a,z)=\mathfrak h+O(\mathfrak he^{-c\mathfrak h}),\quad
 G_F(z,z)=q+O(\mathfrak h),\quad
 0\le G_F(x,x')\le C\mathfrak h.
\]
For one or two sites the inverse $M_F=G_F^{-1}$ has nonpositive
off-diagonal entries, and its row sums satisfy
\[
 e_F:=\min_{z\in F}(M_F\mathbf1)(z)
 \ge q^{-1}-C\mathfrak h q^{-2}>0,\qquad
 e_F\le\lambda:=q^{-1}+C\mathfrak h q^{-2}.
\]
For two sites this follows directly from the inverse of the
$2\times2$ matrix: its determinant is $q^2+O(q\mathfrak h)$
and both inverse-row numerators are $q+O(\mathfrak h)$.
Write $\ell_1(z)$ for the LS occupation at $z$ accumulated during
this one complete record. Let $v_\theta(z)$ be the Laplace transform of occupation in $F$
for this exact-exit kernel started at $z$. The finite killed
resolvent identity gives
\[
 (I+\theta G_F)v_\theta=\mathbf1,\qquad
 E^{a,b}e^{-\theta\sum_{z\in F}\ell_1(z)}
 =1-\theta\sum_{z\in F}\widehat G_b(a,z)v_\theta(z).
\]
At a minimizing coordinate, the equation
$(M_F+\theta I)v_\theta=M_F\mathbf1$ and the sign of the
off-diagonal entries imply $v_\theta(z)\ge e_F/(e_F+\theta)$
for every $z\in F$. Hence, after increasing the fixed $C$,
\[
 \begin{gathered}
 p_j:=j\mathfrak h/q-C\mathfrak h^2/q^2>0,\qquad
 E^{a,b}e^{-\theta\sum_{z\in F}\ell_1(z)}
 \le1-\frac{p_j\theta}{\lambda+\theta}
 \le\exp\!\left\{-\frac{p_j\theta}{\lambda+\theta}\right\},
 \quad\theta\ge0.
 \end{gathered}
\]
Multiply over the $m$ independent records and apply exponential
Markov at level $ju$:
\[
 \begin{aligned}
 P_m\left(\sum_{z\in F}\ell(z)\le ju\right)
 &\le\inf_{\theta\ge0}
       \exp\!\left\{\theta ju-\frac{mp_j\theta}{\lambda+\theta}\right\}\\
 &=\exp\{-(\sqrt{mp_j}-\sqrt{ju\lambda})_+^2\}.
 \end{aligned}
\]
The simultaneous soft event is contained in this sum event.
Sum over at most $Ce^{2q}$ vertices or $Ce^{4q}$ ordered pairs
and use $W^2=\mathfrak h m$. Expansion of the last exponent
proves \eqref{chVI:eq:R17.14}.
When the positive part vanishes, $W=O_u(1)$ and the added constant in that estimate
bounds the trivial number of candidates. Thus all count ranges are covered before
conditioning on $S$.

For fixed $\nu>\gamma_+$,
\begin{equation}
\sup_{W\ge\sqrt2q+q^\nu}
E_m(H_2\mid S)\le e^{-cq^\nu}.
\label{chVI:eq:R17.15}
\end{equation}
Retain the negative quadratic term in \eqref{chVI:eq:R17.14} before dividing by \eqref{chVI:eq:R17.12},
using $s=W-\sqrt2q\ge q^\nu$:
\[
\log E_m(H_2\mid S)
\le-2\sqrt2s-\frac{s^2}{q}
+C\mathfrak h(1+W^2/q^2)+CW/q+C
\le-cq^\nu.
\]
Here $\mathfrak h=o(q^\nu)$ and $\mathfrak h=o(q)$; the negative quadratic term absorbs all higher-count errors.
Let $\tau_{q,y}$ be the first exit from the Euclidean lattice disk
$\{x:d_2(x,y)<e^{q+1}\}$, in the chosen chart, and put
\[
 \omega_{q,y}(\xi):=P_y(X_{\tau_{q,y}}=\xi),\qquad
 \Sigma_{q,y}:=\{\xi:\omega_{q,y}(\xi)>0\}.
\]
Thus $\Sigma_{q,y}$ lies on the one-step outer boundary of this disk and
$\sum_\xi\omega_{q,y}(\xi)=1$. For this exit harmonic measure, set
$H=(\sum_\xi\omega_{q,y}(\xi)\ell(\xi))^{1/2}$.
Chapter~IV's endpointwise Kac and Bernstein estimates give
\begin{equation}
P_m(|H-W|>\Delta)\le 2e^{-c\Delta^2/{\mathfrak h}},
\qquad CWe^{-c{\mathfrak h}}\le\Delta\le W,
\label{chVI:eq:R17.16}
\end{equation}
and the one-sided form
\begin{equation*}
W\le v,\quad \Delta\ge0,\quad \Delta\ge CWe^{-c\mathfrak h}
\quad\Longrightarrow\quad
P_m(H>v+\Delta)\le2e^{-c\Delta^2/\mathfrak h}
\end{equation*}
for every $v\ge W$ and $\Delta\ge0$ satisfying the displayed lower bound.
To prove it, use
$|E_mH^2-W^2|\le CW^2e^{-c\mathfrak h}$. After adjusting
the constants, the event implies
$H^2-E_mH^2\ge c(2v\Delta+\Delta^2)$. In the Bernstein bound
\eqref{chVI:eq:bernstein}, both
\[
 \frac{(2v\Delta+\Delta^2)^2}{W^2\mathfrak h}
 \quad\hbox{and}\quad
 \frac{2v\Delta+\Delta^2}{\mathfrak h}
\]
are at least $c\Delta^2/\mathfrak h$ when $m>0$. For $m=0$,
$H=0$ and the assertion is immediate. This is the form used
when the required upper deviation exceeds $W$.

\subsection{Four local exceptions for sharp clustering}
Fix $\gamma=1/8192$, $\mathfrak h=q^\gamma$, $\eta_j=(j+1)/100$, and $\nu_j=1/2-\eta_j$,
$j=0,\ldots,4$. These parameters are local to this sharp-clustering argument.
In particular,
\begin{equation*}
\nu_0>\cdots>\nu_4>\gamma,\qquad
2\nu_0-\gamma>\nu_1,\qquad2\nu_1-\gamma>\nu_1.
\end{equation*}
Under every exact-count kernel conditioned on $S$, the following implications hold:
\begin{align*}
W>\sqrt2q+q^{\nu_1}
&\ \Longrightarrow\ P_m(H_2>0\mid S)\le e^{-cq^{\nu_1}},\\
W\le\sqrt2q+q^{\nu_1}
&\ \Longrightarrow\ P_m(H>\sqrt2q+q^{\nu_0}\mid S)
\le e^{-cq^{2\nu_0-\gamma}+Cq^{\nu_1}},\\
|W-\sqrt2q|\le q^{\nu_1}
&\ \Longrightarrow\ P_m(Y>e^{q^{\nu_0}}\mid S)
\le e^{-q^{\nu_0}+Cq^{\nu_1}+C\mathfrak h},\\
W<\sqrt2q-q^{\nu_1}
&\ \Longrightarrow\ P_m(H\ge\sqrt2q-q^{\nu_2}\mid S)
\le e^{-cq^{2\nu_1-\gamma}+C\mathfrak h}.
\end{align*}
The first uses \eqref{chVI:eq:R17.15}; the next two use \eqref{chVI:eq:R17.13} with
\eqref{chVI:eq:R17.16} or \eqref{chVI:eq:R17.14}; the last uses the one-sided harmonic estimate and denominator $e^{-C\mathfrak h}$.
All four are at most $e^{-cq^{\nu_1}}$ for large $q$.
Their complement still permits $H<\sqrt2q-q^{\nu_2}$.
For that remaining regime, Appendix~U uses the Gaussian bound
\eqref{chVI:eq:R16.8}, Lemmas~\ref{chVI:lem:energy-probes}
and~\ref{chVI:lem:probe-coarea}, and the probability estimate
\eqref{chVI:eq:R18.16}; its scale sum is bounded by
Corollary~\ref{chVI:cor:low-owner-transfer}.
More precisely, put
\[
 \begin{aligned}
 \mathsf E_1&=\{W>\sqrt2q+q^{\nu_1},\ H_2>0\},\\
 \mathsf E_2&=\{W\le\sqrt2q+q^{\nu_1},\ H>\sqrt2q+q^{\nu_0}\},\\
 \mathsf E_3&=\{|W-\sqrt2q|\le q^{\nu_1},\ Y>e^{q^{\nu_0}}\},\\
 \mathsf E_4&=\{W<\sqrt2q-q^{\nu_1},\ H\ge\sqrt2q-q^{\nu_2}\}.
 \end{aligned}
\]
Then
\[
 \{H_2>0\}\cap
 \bigl(\{H>\sqrt2q+q^{\nu_0}\}\cup\{Y>e^{q^{\nu_0}}\}\bigr)
 \subset\bigcup_{j=1}^4\mathsf E_j
       \cup\{H<\sqrt2q-q^{\nu_2}\}.
\]
Indeed the complement of $\mathsf E_1$ first excludes the upper
count range; $\mathsf E_2$ controls the upper harmonic event,
$\mathsf E_3$ controls the middle-count cardinality, and the
remaining lower count forces the last event unless $\mathsf E_4$
occurs. The intrinsic diameter bound in Appendix~U supplies
$H_2>0$ on the original cloud.

\subsection{Common radial shifts for repelled profiles}
The estimates remain valid when both entrance and exit logarithmic radii are shifted by the same $d\in[0,1]$:
\begin{equation*}
q+\mathfrak h+d+O(e^{-c\mathfrak h}),\qquad
q+2\mathfrak h+d+O(e^{-c\mathfrak h}).
\end{equation*}
Replace $a$ by $\sqrt2(q+2\mathfrak h+d)$; the Gaussian seed exponent changes by $O(1)$.
The capacity denominator and mean block harmonic occupation remain
$\mathfrak h+O(e^{-c\mathfrak h})$, while the annular comparison gaps remain comparable to $\mathfrak h$.
All other changes are absorbed by $C\mathfrak h$.
This concerns a common shift; an arbitrary change of width to $\mathfrak h+O(1)$
would change the critical count center at $m\asymp q^2/\mathfrak h$.

\appendixchapter{U}{Proof of Proposition~\ref{chVI:prop:omitted}: contact bound}
We prove the Gaussian contact estimate \eqref{chVI:eq:R18.3} used
in the omitted-activity argument of Section~6.2.

\subsection{Statement and probability laws}
Let $\ell_N=L^{1/8}$ and $D_A=Ne^{-2\ell_N}$. Throughout this appendix
use the local exponent $\gamma=1/8192$ of Appendix~T. Fix $u_0,v>0$ and $C<\infty$;
the omitted-activity application uses $u_0=1$ and $v=1/4$.
Let $\phi$ have the centered mean-zero Gaussian law of Chapter~II,
with covariance $G_N$, and let $h:=\phi-\phi(r)$ be independent
of the path. Write $P_h,P_\phi$ for these laws and $E_h,E_\phi$
for their expectations, holding the path fixed. Define
\begin{equation}
S(t)=\{x:\ell_t(x)\le u_0\},\quad
G_t(v)=\{x:2\ell_t(x)+h(x)^2\le v\},
\label{chVI:eq:R18.1}
\end{equation}
Use the cover in Section~6.2--6.3 and set $C_z^\circ=S(t)\cap P_z$, $C_z=S(t)\cap P_z^{\rm buf}$.
Same-color extended buffers are $3D_A$-separated, core-to-buffer margins are at least $3D_A$,
buffer diameters are at most $C_{\rm at}Ne^{-\ell_N}$, and overlap multiplicities are bounded.
Here $C_{\rm at}$ bounds the ratio of a buffer's Euclidean diameter
to $Ne^{-\ell_N}$. For the atlas of Section~6.2 and sufficiently large $N$,
\[
 \diam_2 P_z^{\rm buf}\le6\sqrt2\,s
 \le\frac{6\sqrt2}{1000}Ne^{-\ell_N}<Ne^{-\ell_N},
\]
so we take $C_{\rm at}=1$. Put
\[
 \mathbb W_N^{\rm mac}(t):=\{z:C_z\ne\varnothing\},\qquad
 \mathrm{RF}:=\{d_1(x,r)\ge3D_A\text{ for every }x\in S(t)\}.
\]
For $u_0=1$, $\mathbb W_N^{\rm mac}$ is the family
$\mathcal O^{\rm active}$ used in the body. Intrinsic owners below are
those with $C_z^\circ\ne\varnothing$.

For later image tests, use the fixed-root shell law on $(\phi,\sigma)$,
where $\phi$ has the mean-zero Gaussian law $P_\phi$ of Chapter~II and
$\sigma\in\{-1,1\}$. Set $\mathcal{Q}_N=N^{-2}\|\phi\|_2^2$,
$B_s=2s/N^2$, and $a_s=(B_s-\mathcal{Q}_N)^{1/2}$ on $\mathcal{Q}_N<B_s$.
With $b_{r,s}^\sigma=0$ off this set, put
\[
 \begin{split}
 b_{r,s}^\sigma(\phi)&=\frac{[\phi(r)+\sigma a_s]_+}{N^2a_s}
       \quad(\mathcal{Q}_N<B_s),\\
 \varrho_{N,r}(s)&=E_\phi\sum_{\sigma=\pm1}b_{r,s}^\sigma(\phi),\\
 \nu_{N,s}^r(d\phi,\{\sigma\})
   &=\varrho_{N,r}(s)^{-1}b_{r,s}^\sigma(\phi)P_\phi(d\phi).
 \end{split}
\]
These are the branch densities $b_{r,s}^\sigma$ of Chapter~III.
On the strict shell support, the square image is
$W_s^\sigma(\phi)=(\phi+\sigma a_s(\phi)\mathbf1)^2$;
it is set to zero off that support. The law $\nu_{N,s}^r$ is used on the central-shell range of
Section~3.3, where $0<\varrho_{N,r}(s)<\infty$ and
$n\varrho_{N,r}(s)=1+o(1)$ uniformly. The enlarged ranges produced by
the fixed energy cutoff also remain central: for fixed $D$,
$DnL/2=(D/(2g))b_N$, so $I+[0,DnL/2]$ lies in a fixed enlargement of
the central window. For a square-field event $E$ we write
$\nu_{N,s}^r(W_s^\sigma(\phi)\in E)$. For a nonnegative square field $W$, define
\[
 \mathcal R_{q,y}^{\rm sq}(W):=
 \left(\frac1{2g}\sum_{\xi\in\Sigma_{q,y}}\omega_{q,y}(\xi)W(\xi)\right)^{1/2},
\]
using the exit weights and support defined in Appendix~T.

Use the chart domains $U_\alpha$, grids $\mathcal G_{N,q}(a)$,
outer disks $D_{q,y}$ and available-locator sets
$\mathscr A_{N,q}(a)$ defined in
\hyperref[chVI:def:intrinsic]{Section~6.2}. Here
$\gamma=\gamma_{\rm loc}$ and the origin $a$ is fixed (take $a=0$
for Proposition~\ref{chVI:prop:omitted}). For this threshold $u_0$,
define
\[
 q_z:=\inf\{q\in\ZZ:q\ge K_{\rm loc},\
 \exists(\alpha,y)\in\mathscr A_{N,q}(a),\ C_z\subset B(y,e^q)\},
 \qquad \inf\varnothing:=+\infty,
\]
for the owners with $C_z^\circ\ne\varnothing$. For finite $q_z$,
use the fixed deterministic order $\prec$ to bind
\[
 (\alpha_z,y_z):=\min_{\prec}
 \{(\alpha,y)\in\mathscr A_{N,q_z}(a):C_z\subset B(y,e^{q_z})\}.
\]
The local kernel at this intrinsic locator has entrance disk $B(y,e^{q+q^\gamma})$ and exit domain
$D_{q,y}=B(y,e^{q+2q^\gamma-4e^{-q^\gamma}})$.
Its logarithmic perturbation $-4e^{-q^\gamma}$ is covered by
Appendix~T. Root exclusion and count eligibility are imposed
separately below. Set $k_{\rm mid}:=\lceil L^{1/8}+10\rceil$.

\begin{lemma}[Uniform intrinsic-scale geometry]\label{chVI:lem:intrinsic-geometry}
There are $c_*>0$ and $K_{\rm loc}<\infty$ such that, for all sufficiently large $N$,
\begin{equation}\label{chVI:eq:R18.2}
 K_{\rm loc}<q_z<k_{\rm mid}
 \ \Longrightarrow\
 C_z\subset B(y_z,e^{q_z}),\quad D_{q_z,y_z}\subset P_z^{\rm buf},\quad
 \diam_2 C_z\ge c_*e^{q_z}.
\end{equation}
If $q_z\ge k_{\rm mid}$, including $q_z=+\infty$, then
$\diam_2 C_z\ge c_*e^{k_{\rm mid}}$.
On the event that every soft point has $d_1$-distance at least $3D_A$ from the root,
the fine intrinsic outer domains are root-free.
\end{lemma}
\begin{proof}
Use the deep-cover property of the fixed $U_\alpha$ and $\delta_0$
from Section~6.2.
For every $K_{\rm loc}\le q<k_{\rm mid}$, a nearest preceding-scale grid point to any $x$
is at distance at most
\[
 2^{-1/2}\lfloor e^{q-1}\rfloor+2\le\tfrac34e^{q-1}
\]
after increasing $K_{\rm loc}$. This point remains available in a chart
in which $x$ has depth at least $\delta_0N$: its depth is at least
$\delta_0N/2$, and both the preceding outer radius and its distance
to $x$ are $o(N)$, uniformly over these scales.
For an owner with $C_z^\circ\ne\varnothing$, choose $x\in C_z^\circ$. If
$\diam_2 C_z<e^{q-1}/8$, the nearest preceding locator $y'$ satisfies
\[
 \sup_{w\in C_z}d_2(w,y')<\tfrac78e^{q-1}\le\lfloor e^{q-1}\rfloor.
\]
Its core contains $C_z$, contradicting minimality when $q=q_z>K_{\rm loc}$.
Thus $c_*=1/(8e)$ is sufficient. If $q_z\ge k_{\rm mid}$, apply the same argument directly
at the unavailable containment scale $k_{\rm mid}-1$; no finite $q_z$ is required.

For the chosen fine locator, $d_2(y_z,x)\le e^{q_z}$ and
\[
 e^{q_z}+e^{q_z+2q_z^\gamma}+2=o(D_A).
\]
The core-to-buffer margin is at least $3D_A$, so $D_{q_z,y_z}$ and its boundary lie inside
the buffer. The same estimate and $d_2(x,r)\ge3D_A/\sqrt2$ imply root exclusion.
Chart availability, buffer containment, root exclusion and radius rounding have thus been
checked separately, with constants uniform over the full fine-scale range.
\end{proof}

We prove the physical-window contact estimate
\begin{equation}\label{chVI:eq:R18.3}
 \lim_{K\to\infty}\limsup_{N\to\infty}\sup_r
 (\mathsf Q_{N,r}^I\otimes P_h)
 \left(G_t(v)\cap\bigcup_{z:q_z(t)\ge K}C_z^\circ(t)
                 \ne\varnothing\right)=0.
\end{equation}
The rooted window law is
\begin{equation}\label{chVI:eq:R18.4}
 \mathsf Q_{N,r}^I(F)=
 \frac{\int_I p_t(r,r)E_t^{r,r}F\,dt}{\int_I p_t(r,r)\,dt}.
\end{equation}
The field is sampled independently after the path. The energy cutoff
below is retained when applying coarea; independence is only used
under this original product law.

\subsection{Root exclusion, active-count tightness, and a growing-scale guard}
\subsubsection{Soft points near the root}
For fixed $\theta>0$, let $A_{\theta,x}=-\mathcal L_N+
\theta1_{\{x\}}$, where the second term denotes multiplication by
the indicated function. Write $\lambda_x$ for its principal
eigenvalue, equal to $[n(gL+O(1))]^{-1}$. Set
$\Delta:=C_0n\log L$, choosing a fixed sufficiently large $C_0$,
and leave a mixing segment of length $\Delta$
at both kernel ends; domination and the operator norm give
\begin{equation}
e^{-tA_{\theta,x}}(a,b)\le\frac Cn e^{-\lambda_x(t-2\Delta)}
\le Cn^{-2}L^{3/2}
\label{chVI:eq:R18.5}
\end{equation}
uniformly in central $t$ and $a,b,x$.
Exponential Markov and division by $p_t(r,r)\asymp n^{-1}$ give the loop soft bound
$Cn^{-1}L^{3/2}$.
Thus the probability of a soft point within $3D_A$ of $r$ is
$O(L^{3/2}e^{-4\ell_N})=o(1)$.

\subsubsection{Soft packing and active macro-count tightness}
For $R>0$ and $A\subset V_N$, define
\[
 \operatorname{Pack}_R(A):=\max\{|F|:F\subset A,\
 d_2(x,y)\ge R\text{ for all distinct }x,y\in F\}.
\]
For later use first record the exact packing statement needed in Sections~6.2 and~6.6.
\begin{lemma}[Central rooted-loop soft packing]\label{chVI:lem:loop-soft-packing}
For every fixed central physical-window range,
\begin{equation}\label{chVI:eq:loop-soft-packing}
 \lim_{B\to\infty}\limsup_{N\to\infty}\sup_r
 \mathsf Q_{N,r}^I
 \left(\operatorname{Pack}_{D_A}(S(t))>B\sqrt L\right)=0,
\end{equation}
uniformly over physical windows $I$ of length comparable to $b_N$ in that range.
\end{lemma}
\begin{proof}
This is the physical-window packing input of Chapter~IV, recalled
in Section~6.1. It precedes every local raw-count estimate used here.
\end{proof}

The endpoint estimates used above are, for root-far $D_A$-separated points,
\begin{equation}
P_h(|h(x)|\le1)\asymp L^{-1/2},\qquad
P_h(|h(x)|\le1,|h(x')|\le1)\le C/L.
\label{chVI:eq:R18.6}
\end{equation}
They follow from the one-point variance $2gL+O(\ell_N)$ and the two covariance eigenvalues
$gL+O(\ell_N)$ and $3gL+O(\ell_N)$.

For the active macro count, in each fixed macro color choose one soft point from every active
buffer. Same-color buffers are separated by more than $2s>3D_A$, so these points form a
$D_A$-packing. Let $C_{\rm mac}$ be a fixed upper bound on the number of macro colors. Therefore
$\#\mathbb W_N^{\rm mac}\le C_{\rm mac}\operatorname{Pack}_{D_A}(S(t))$, and
\eqref{chVI:eq:loop-soft-packing} yields
\begin{equation}
\lim_{B\to\infty}\limsup_{N\to\infty}\sup_r
\mathsf Q_{N,r}^I(\#\mathbb W_N^{\rm mac}>B\sqrt L)=0.
\label{chVI:eq:R18.7}
\end{equation}
This bound follows directly from the soft-packing estimate
\eqref{chVI:eq:loop-soft-packing}.

\subsubsection{A guard over all fine locators}
Set $M_N=\lceil e^{L/8}\rceil$, $d_N=e^{L/2}$. For a locator
$(q,y)$, put $\mathfrak h=q^\gamma$ and let $m_{q,y}(t)$ count
the completed records from $B(y,e^{q+\mathfrak h})$ to exit from
$D_{q,y}$, using the alternating entrance/exit convention of Section~6.3.
The \emph{guard event} requires $m_{q,y}(t)\le M_N$ and duration
at most $d_N$ for each record, simultaneously at all available
locators with $K_{\rm loc}\le q<k_{\rm mid}$.
Under the ordinary rate-one walk law $\mathbb P_a$, with expectation
$\mathbb E_a$, each block started before $t$ hits a fixed core point
with probability at least $c\mathfrak h/q$;
at most one block straddles $t$. Strong Markov and the visit bound $C(1+t/n+L)$ give
\begin{equation}
\sup_a \mathbb E_a m_{q,y}(t)\le C(q/{\mathfrak h})(1+t/n+L)
\le C L^{17/8}.
\label{chVI:eq:R18.8}
\end{equation}
With the same $C$ as in the final bound, iteration at exits yields
$\sup_a\mathbb E_a\binom{m_{q,y}(t)}j\le(CL^{17/8})^j$, hence
\begin{equation*}
\sup_a\mathbb P_a(m_{q,y}(t)>M_N)\le2e^{-M_N/(4CL^{17/8})}.
\end{equation*}
If $R_D$ is the radius of $D_{q,y}$, then $R_D^2=e^{O(\ell_N)}$
and the unconditional exit tail is $2e^{-cd_N/R_D^2}$.
Sum over the first $M_N+1$ blocks, divide by the loop return probability, and sum at most
$CnL^{1/8}$ locators. The guard failure is at most $e^{-e^{L/9}}$.
For a fixed count $m$, let $A_{d_N}$ be the event that all $m$
record durations are at most $d_N$. Conditional on core softness
$S=\{\min_{B(y,e^q)}L^{\rm loc}\le u_0\}$, the duration estimate
of Chapter~III gives
$P_m(A_{d_N}^c\mid S)\le e^{-e^{0.49L}}$, uniformly for
$m\le M_N$ and all feasible words.

\subsection{Certified selection and local-exception owner counts}
Use the sharp-clustering exponents fixed in the final local-exception subsection of Appendix~T.
For an owner $z$ with finite intrinsic locator $(q_z,y_z)$, set
\[
 \mathsf H_z:=\left(\frac1g\sum_\xi
 \omega_{q_z,y_z}(\xi)\ell_t(\xi)\right)^{1/2}.
\]
At a fixed deterministic scale $q$, the up and big owner events for
owners with $q_z=q$ are
\begin{equation}\label{chVI:eq:R18.9}
 \mathsf H_z>\sqrt2q+q^{\nu_0},\qquad |C_z|>e^{q^{\nu_0}}.
\end{equation}
\phantomsection\label{chVI:def:host-selection}
Fix a deterministic $q$, a macro color and one of the four Appendix~T
count branches. Let $\kappa_q$ be a deterministic coloring of the
outer locator disks and their one-step boundaries with at most
$e^{Cq^\gamma}$ colors, disjoint within each color. For a fixed
locator color $\beta$ and an owner $z$ of the chosen macro color, set
\[
 I(q,z):=\{y\in\mathcal G_{N,q}(a):
 D_{q,y}\cup\partial_oD_{q,y}\subset P_z^{\rm buf},\
 r\notin D_{q,y},\ \kappa_q(y)=\beta\}.
\]
The owner of this host is the specified $z$. The candidate list is
deterministic given the root and grid origin. Simultaneously delete
all complete blocks for this list, as in the local product identity
of Chapter~III. Write $W_{q,y}=\sqrt{q^\gamma m_{q,y}}$ and define
the four real intervals corresponding to the branches of Appendix~T:
\[
 \begin{aligned}
 J_1(q)&=(\sqrt2q+q^{\nu_1},\infty),&
 J_2(q)&=(-\infty,\sqrt2q+q^{\nu_1}],\\
 J_3(q)&=[\sqrt2q-q^{\nu_1},\sqrt2q+q^{\nu_1}],&
 J_4(q)&=(-\infty,\sqrt2q-q^{\nu_1}).
 \end{aligned}
\]
For the currently fixed branch $b\in\{1,2,3,4\}$ define
\[
 I^{\rm elig}_b(q,z):=\{y\in I(q,z):m_{q,y}\le M_N,
 \ W_{q,y}\in J_b(q)\}.
\]
This eligibility set is measurable from the collapsed exterior of
Lemma~\ref{chIII:lem:multihost}, because that coordinate retains host
labels, endpoint words and complete-record counts. We write
$I^{\rm elig}(q,z)$ for the branch currently under consideration.
If $V_j$ denotes the duration of the $j$th local record, put
\[
 A^{\rm dur}_{q,y}:=\bigcap_{j=1}^{m_{q,y}}\{V_j\le d_N\},
\]
with the empty intersection equal to the whole local record space.
The sole local flag is
\[
 A_{q,y}=\{\min_{B(y,e^q)}L^{\rm loc}\le u_0\}\cap A^{\rm dur}_{q,y}.
\]
Here is the precise selection rule, also used in Appendices~W--X.
Fix $B>0$; the \emph{active cap} is the event
$\{\#\mathbb W_N^{\rm mac}(t)\le B\sqrt L\}$, with $B$ held fixed
as $N\to\infty$. Fix deterministic orders of the hosts and owners. For a candidate index $i=(z,y)$, its host has entrance
$B(y,e^{q+q^\gamma})$, exit domain $D_{q,y}$, complete-record
coordinate $\xi_i$, record count $m_i=m_{q,y}$ and endpoint word
$w_i=((a_{i,j},b_{i,j}))_{j=1}^{m_i}$. Its local kernel is
\[
 P_{m_i,w_i}=\bigotimes_{j=1}^{m_i}
 \mathsf B_{D_{q,y}}^{a_{i,j},b_{i,j}},
\]
and $A_i=A_{q,y}$ is its local event. Put
$s_i=\one_{\{\xi_i\in A_i\}}$. After the exterior eligibility sets
$I^{\rm elig}(q,z)$ have been fixed, define the selector as a function
of the flag vector $s$ by
\[
 \begin{aligned}
 I^+(q,z;s)&:=\{y\in I^{\rm elig}(q,z):s_{(z,y)}=1\},\\
 \mathscr S_0(s)&:=\{(z,\min_{\prec}I^+(q,z;s)):
                         I^+(q,z;s)\ne\varnothing\},\\
 \operatorname{Sel}_{M_B}(s)&:=\operatorname{First}_{M_B}(\mathscr S_0(s)),
 \qquad M_B:=\lceil C_{\rm adj}B\sqrt L\rceil.
 \end{aligned}
\]
Here $\operatorname{First}_{M_B}$ uses the fixed candidate order;
the actual selected set is $\operatorname{Sel}_{M_B}(s)$ at the
realized flag vector. There is at most one slot per active owner, so
$C_{\rm adj}=1$ suffices. The fixed bound $C_{\rm adj}\ge1$ in $M_B$
is also used in Appendices~W--X.
For a proposed selected set $C$, let $\mathbf y$ denote the collapsed
exterior together with all unselected records, and let
$s_{C^c}(\mathbf y)$ be their flag vector. With exterior eligibility
held fixed, put
\[
 \begin{aligned}
 \Lambda_C(\mathbf y)&:=\left\{
 \operatorname{Sel}_{M_B}(\mathbf1_C,s_{C^c}(\mathbf y))=C\right\},\\
 \{\operatorname{Sel}_{M_B}(s)=C\}
 &=\Lambda_C(\mathbf y)\cap\bigcap_{i\in C}\{\xi_i\in A_i\}.
 \end{aligned}
\]
On every positive-mass stratum Lemma~\ref{chIII:lem:selection} gives
\[
 \mu_{C,\mathbf y}
 =\bigotimes_{i\in C}P_{m_i,w_i}(\,\cdot\mid A_i),
\]
and the exterior measure retains
\[
 \one_{\Lambda_C(\mathbf y)}
 \prod_{i\in C}P_{m_i,w_i}(A_i).
\]

Define the owner exception sets before applying the selection estimate.
Evaluate $\mathsf E_1,\ldots,\mathsf E_4$ of Appendix~T on the
intrinsic host, using its parameters, count and endpoint word. Assign
the owner to the first index for which the event holds, in the order
$1<2<3<4$. For $\ell\in\{1,2,3,4\}$, let $\mathscr E_{q,\ell}$
be the set of owners in all macro and locator colors such that
$q_z=q<\infty$, the intrinsic locator $(\alpha_z,y_z)$ is defined
and eligible for branch $\ell$, its realized coordinate
$\xi_{(z,y_z)}$ belongs to $A_{(z,y_z)}$, and its first exception
index is $\ell$. Owners failing these allocation conditions belong
to none of these sets. Put
\[
 \mathscr E_q:=\bigcup_{\ell=1}^4\mathscr E_{q,\ell}.
\]
Each owner has its unique macro color and the locator color $\kappa_q(y_z)$.
For the original-path counting step, Lemma~\ref{chVI:lem:intrinsic-geometry}
gives $\varnothing\ne C_z\subset B(y_z,e^q)$ for $q_z=q$.
Every other same-color candidate core in its buffer is disjoint from
this core and contains no soft point. On the guard event, the intrinsic
host is its owner's unique positive flag; on the active cap the bound
$M_B\ge B\sqrt L$ prevents truncation.
For the currently fixed macro color and locator color $\beta$,
apply the selector for branch $b=\ell$. On these original-sample events,
\[
 \left|\{z\in\mathscr E_{q,\ell}:z\text{ has the fixed macro color},\
                            \kappa_q(y_z)=\beta\}\right|
 \le\sum_{i\in\operatorname{Sel}_{M_B}(s)}
        \one_{\{\xi_i\in\mathsf E_\ell\}}.
\]
Here $\mathsf E_\ell$ on the right is evaluated with host $i$'s
parameters and endpoint word. The assignment is injective in the owner.
The conditions $q_z=q$ and $C_z\subset B(y_z,e^q)$ restrict the
original-sample inequality; the selected law retains only the local
flag $A_{q,y}$ defined above.

For the two conditional resamples, selected total durations lie in
\begin{equation}\label{chVI:eq:R18.10}
 [0,M_BM_Nd_N],\qquad M_BM_Nd_N=e^{5L/8+o(L)}=o(b_N).
\end{equation}
Duration reconditioning costs $(1-e^{-e^{0.49L}})^{-1}$ by the previous subsection.
\phantomsection\label{chVI:def:weighted-window}
Here is the weighted window estimate used also in Appendices~W--X.
On a fixed selected stratum put $d=M_BM_Nd_N$ and let
$V(\xi)$ be its total selected duration. The duration flags give
$0\le V(\xi)\le d$ under $\mu_{C,\mathbf y}$. For every nonnegative
selected-record test $F$ and the window functional of
\eqref{chIII:eq:window-functional}, two independent copies satisfy
\[
 \begin{aligned}
 \int F(\xi)J_I(\xi,\mathbf y)\,\mu_{C,\mathbf y}(d\xi)
 &\le\iint F(\xi)J_{I^d}(\xi',\mathbf y)
                \,\mu_{C,\mathbf y}(d\xi)\mu_{C,\mathbf y}(d\xi')\\
 &=\mu_{C,\mathbf y}(F)\int J_{I^d}(\xi',\mathbf y)\,\mu_{C,\mathbf y}(d\xi'),
 \qquad I^d=[l_N-d,r_N+d].
 \end{aligned}
\]
The inequality is exactly \eqref{chIII:eq:short-window}; Tonelli
permits unbounded $F$. In particular, if $F=\sum_{i\in C}F_i$ and
$P_{m_i,w_i}(F_i\mid A_i)\le a_q$ on every stratum, then summing
the exterior measures, endpoint words and selected sets gives
\[
 E_{\mathsf Q_{N,r}^I}\sum_{i\in\operatorname{Sel}_{M_B}(s)}F_i
 \le M_Ba_q\,
 \frac{\int_{I^d}p_t(r,r)\,dt}{\int_Ip_t(r,r)\,dt}
 =M_Ba_q(1+o_N(1)).
\]
The last equality is uniform in the root and admissible windows,
since $d=o(b_N)$ and $|I|\ge\delta b_N$. Window transfer is applied
at the product-kernel level with the factor $J_I$ retained inside the
integral.

Summing over the four exception indices and all macro and locator
colors gives $|\mathscr E_q|$. This random variable is determined
from the original sample. The original-sample counting inequality
and weighted conditional estimate, with the cap, guard and
$\mathrm{RF}$ restrictions retained, give
\begin{equation}\label{chVI:eq:R18.11}
 \begin{aligned}
 &E_{\mathsf Q_{N,r}^I}\Bigl[|\mathscr E_q|;\
     \{\#\mathbb W_N^{\rm mac}(t)\le B\sqrt L\}\cap\mathrm{RF}\\
 &\qquad\cap\!\bigcap_{\substack{K_{\rm loc}\le q'<k_{\rm mid}\\
                        (\alpha,y)\in\mathscr A_{N,q'}(a)}}
       \bigl(\{m_{q',y}(t)\le M_N\}\cap A^{\rm dur}_{q',y}\bigr)\Bigr]
 \\
 &\qquad\le CB\sqrt L\,e^{Cq^\gamma-cq^{\nu_1}}.
 \end{aligned}
\end{equation}
The first exception uses the separated pair supplied by \eqref{chVI:eq:R18.2};
the third uses $C_z=S\cap B(y_z,e^q)$ on the original containment event.
The charging inequality is evaluated on these original-sample restrictions before
the conditional local estimate is integrated. Markov at threshold $e^{-q^{\nu_4}}\sqrt L$ and a scale sum give
\[
 CB\sum_{q\ge K}e^{Cq^\gamma+q^{\nu_4}-cq^{\nu_1}}+o_N(1)\longrightarrow0
\]
in the order $N$, then $K$, because $\gamma<\nu_4<\nu_1$.

\subsection{Low-harmonic owners and exact energy-truncated amplification}
Define
\begin{equation}
 D_q=\{z:q_z=q,\ \mathsf H_z<\sqrt2q-q^{\nu_2}\},\qquad
 E_q^{\rm low}=\{|D_q|>e^{-q^{\nu_4}}\sqrt L\}.
 \label{chVI:eq:R18.12}
\end{equation}
On $E_q^{\rm low}$ and the common root-far event, choose owners in one macro color, with their number satisfying
\[m\ge c\max\{1,e^{-q^{\nu_4}}\sqrt L\}.\]
 Choose one core soft point $x_j$ from each. Use fixed orders throughout; the choices are Borel functions of the path alone.
Their points are $cD_A$-separated. For a fixed $0<\chi<\nu_2$, use the simultaneous probes
\begin{equation}
 B_j=\left\{|h(x_j)|\le1,\quad
 \frac1{2g}\sum_\xi\omega_{q,y_j}(\xi)h(\xi)^2\le q^{1+\chi},\quad
 \|h\|^2\le DnL\right\}.
 \label{chVI:eq:R18.13}
\end{equation}
Lemma \ref{chVI:lem:energy-probes} gives a fixed $D$, independent of $q,N,m$, such that
$P_h(B_j\mid\text{path})\ge c/\sqrt L$ and
$P_h(B_i\cap B_j\mid\text{path})\le C/L$. Consequently,
\begin{equation}
 P_h\left(\bigcup_{j=1}^mB_j\,\middle|\,\text{path}\right)
 \ge c\frac{m/\sqrt L}{1+m/\sqrt L}\ge c'e^{-q^{\nu_4}}.
 \label{chVI:eq:R18.14}
\end{equation}
When the nominal number is below $1$, its rounded value $1$ makes this inequality stronger.
On $B_j$, the square field $W=2\ell_t+h^2$ has $W(x_j)\le2u_0+1$ and
\begin{equation}
 \mathcal R_{q,y_j}^{\rm sq}(W)^2
 =\mathsf H_{z_j}^2+\frac1{2g}\sum_\xi\omega_{q,y_j}(\xi)h(\xi)^2
 \le(\sqrt2q-\tfrac12q^{\nu_2})^2.
 \label{chVI:eq:R18.15}
\end{equation}
Here $q^{1+\chi}=o(q^{1+\nu_2})$. Define the square-field event
\[
 \mathcal A_q:=\{W\in[0,\infty)^{V_N}:\exists(\alpha,y)\in
 \mathscr A_{N,q}(a),\
 \mathcal R_{q,y}^{\rm sq}(W)\le\sqrt2q-\tfrac12q^{\nu_2}\}.
\]
It depends only on the coarea image and the deterministic locator
family; $s$ enters through its testing law $\nu_{N,s}^r$.
By Appendix~Y, after replacing the factor $1/2$ by an exponent strictly between
$\nu_3$ and $\nu_2$,
\[
 \sup_{\substack{r\in V_N\\s\in I_D}}\nu_{N,s}^r(W_s^\sigma(\phi)\in\mathcal A_q)
 \le C_D(e^{-cq^{\nu_3}}+e^{-cLq^a}),\qquad I_D=I+[0,DnL/2].
\]
Applying Lemma \ref{chVI:lem:probe-coarea} to $E_q^{\rm low}$ and the union of the $B_j$ gives
\begin{equation}
 \mathsf Q_{N,r}^I(E_q^{\rm low}\cap\mathrm{RF})
 \le C_D e^{q^{\nu_4}}(e^{-cq^{\nu_3}}+e^{-cLq^a}).
 \label{chVI:eq:R18.16}
\end{equation}
The union bound separates the common root-far complement:
\begin{equation}\label{chVI:eq:low-owner-union}
 \begin{split}
 \mathsf Q_{N,r}^I\left(\mathrm{RF}^c\cup
             \bigcup_{K\le q<k_{\rm mid}}E_q^{\rm low}\right)
 &\le\mathsf Q_{N,r}^I(\mathrm{RF}^c)
    +\sum_{K\le q<k_{\rm mid}}
       \mathsf Q_{N,r}^I(E_q^{\rm low}\cap\mathrm{RF}).
 \end{split}
\end{equation}
Since $\nu_4<\nu_3$, Corollary~\ref{chVI:cor:low-owner-transfer}
makes the sum of the bounds in \eqref{chVI:eq:R18.16} vanish as
$N\to\infty$, then $K\to\infty$.
For a fixed path, put
\[
 \mathcal B_q:=\{z:q_z=q,\ |C_z|>e^{q^{\nu_0}}\text{ or }
                         \mathsf H_z>\sqrt2q+q^{\nu_0}\}.
\]
Together with the four count regimes of Appendix~T and
\eqref{chVI:eq:R18.11}, this yields, outside the active-count exception and a vanishing probability,
\begin{equation}
 |\mathcal B_q|
 \le C e^{-q^{\nu_4}}\sqrt L
 \quad\text{simultaneously for }K\le q<k_{\rm mid}.
 \label{chVI:eq:R18.17}
\end{equation}
Paley--Zygmund is applied to the energy-truncated probes under the
source path--Gaussian law. Exact coarea transfers their image-event
inclusion, giving \eqref{chVI:eq:R18.16}.

\subsection{From exceptional owners to cloud contact exclusion}
The whole-ball estimate \eqref{chVI:eq:R16.9}, summed over
\eqref{chVI:eq:R18.17}, gives
\begin{equation}
P_h\bigl(G_t(v)\cap\!\bigcup_{q\ge K}\bigcup_{z\in\mathcal B_q}C_z
          \ne\varnothing\mid\text{path}\bigr)
\le C_v\sum_{q\ge K}(q+1)e^{-q^{\nu_4}}\to0.
\label{chVI:eq:R18.18}
\end{equation}
For the remaining owners,
$|C_z|\le e^{q^{\nu_0}}$ and $\mathsf H_z\le\sqrt2q+q^{\nu_0}$.
The full repulsion estimate of Appendix~Y with $\delta=\eta_0/2$ gives, outside a joint all-scale exception tending to zero,
\begin{equation}
\mathcal R_{q,y}^{\rm sq}(2\ell_t+h^2)
\ge\sqrt2q+q^{1/2-\eta_0/2}.
\label{chVI:eq:R18.19}
\end{equation}
at a contact point.
On $\|h\|^2\le DnL$, a failure together with a contact implies a square-image
repulsion failure. The exact push-forward inequality used in Lemma \ref{chVI:lem:probe-coarea}
transfers its probability to the fixed-root shell.
The complementary source Gaussian probability is at most $C/D$; let $D\to\infty$ last.
Comparison of squares forces
\begin{equation}
\left(\frac1{2g}\sum_\xi\omega_{q,y}(\xi)h(\xi)^2\right)^{1/2}
\ge c q^{3/4-\eta_0/4}.
\label{chVI:eq:R18.20}
\end{equation}
The same-endpoint estimate \eqref{chVI:eq:R16.11} bounds the joint event of contact and this large radial norm by
$C_vL^{-1/2}e^{-cq^{1/2-\eta_0/2}}$.
There are at most $e^{q^{\nu_0}}$ vertices per cloud and $B\sqrt L$ owners, hence the remaining contact probability is at most
\begin{equation}
CB\sum_{q\ge K}
\exp\{q^{\nu_0}-cq^{1/2-\eta_0/2}\}\to0.
\label{chVI:eq:R18.21}
\end{equation}
Combine \eqref{chVI:eq:R18.7} and \eqref{chVI:eq:R18.17}--\eqref{chVI:eq:R18.21}, taking $N$, then $K$, then $B,D$.
This proves fine-scale exclusion uniformly in the root and central physical window.

\subsection{Large intrinsic scales and a whole-ball closing estimate}
If $q_z\ge k_{\rm mid}$, minimality gives two cloud points separated by at least $ce^{\ell_N+9}$.
Every core point is at distance greater than $e^{\ell_N}/2$ from at least one,
and at distance less than $C_{\rm at}Ne^{-\ell_N}$.
Let $\widehat D$ be the set of points in $S(t)$ with a witness in $S(t)$
at a $d_2$-distance in $(e^{\ell_N}/2,Ne^{-\ell_N})$. Since $d_2=d$,
\[
 \bigcup_{z:q_z\ge k_{\rm mid}}C_z^\circ
 \subset\widehat D\subset D_N^{\rm wide}(S(t)),
\]
by \eqref{chIV:eq:wide-def}. Theorem~\ref{chIV:thm:defects} gives
$|\widehat D|=o_P(\sqrt L)$ uniformly at ordinary central times.
For loops, delete a terminal $C_0n\log L=o(b_N)$ mixing segment:
the prefix density is $1+o(1)$ relative to an ordinary path, and defects decrease with added occupation.
At a root-far defect point, conditional on the path,
\[
 P_h(2\ell_t(x)+h(x)^2\le v\mid\text{path})
 \le P_h(|h(x)|\le\sqrt v)\le C_v/\sqrt L,
\]
by the same one-dimensional Gaussian density estimate underlying \eqref{chVI:eq:R18.6}.
A union bound over the root-far points of $\widehat D$ therefore gives $1\wedge(C_v|\widehat D|/\sqrt L)$;
its expectation tends to zero.

\appendixchapter{V}{Proof of Proposition~\ref{chVI:prop:common}}
\renewcommand{\theequation}{AV.\arabic{equation}}
\renewcommand{\theHequation}{AV.\arabic{equation}}

This appendix gives the lattice mixture used in
Proposition~\ref{chVI:prop:common}. The unnormalized comparison is
useful because the anchor event may have very small probability.
We keep the body's zero-based indexing.

\subsection{The reference measure and the terminal factor}
Use the parameters $\gamma,\eta,\eta'$ of \eqref{chVI:eq:exponents},
put $\mathfrak h=k^\gamma$ and $T=\lfloor(k^2-k)/\mathfrak h\rfloor$, and center all
disks at $0\in\ZZ^2$. Thus $B_i^\pm=B_i^\pm(0)$ are the planar
versions of the disks in Section~6.3. Occupation fields in this
appendix are in rate-one units.
Recall the endpoint laws $\pi_k^\pm$, finite word set
$\mathcal W_k$, weights $p_k(w),\alpha(w)$ and field laws
$\mathcal K_w$ from \hyperref[chVI:def:common-measure]{Section~6.3}.
In particular, the conditioning domain in $\pi_k^-$ is exactly
$B_1^+=\mathbb B(0,k+2\mathfrak h-4e^{-\mathfrak h})$. Their finite mixture on
$\mathscr X_k=[0,\infty)^{Q(k)}$ is
\begin{equation}\label{chVI:eq:reference}
 \mathcal R_k=
 \sum_{w\in\mathcal W_k}p_k(w)\alpha(w)e^{2\sqrt2(\alpha(w)-1)}\mathcal K_w.
\end{equation}
The sum is finite: the corridor bounds $m$ by $9k^4$ for large
$k$, and both endpoint alphabets are finite. Zero-weight or
infeasible words make zero contribution.

For a terminal word $z$ of length $m(z)$ at layer $T-1$, use
$b_z=\sqrt{\mathfrak h m(z)}-\sqrt2q_{T-1}$ as defined in
\hyperref[chVI:def:terminal-measure]{Section~6.3}. Define
\begin{equation}\label{chVI:eq:leading}
 \mathcal L_k(a,b)=\sqrt{2/\pi}\frac{ab}{Tk^{3/2}}
 \exp\{-2(T-1)\mathfrak h+2\sqrt2a-2\sqrt2b-b^2/(T\mathfrak h)\}.
\end{equation}
Its first-variable dependence is
$\mathcal L_k(a,b)=ae^{2\sqrt2(a-1)}\mathcal L_k(1,b)$.
This identity explains the weight in \eqref{chVI:eq:reference}.

\subsection{The closed annular corridor estimate}
For a terminal-record product coordinate, define the measurable
\emph{first-word map} $\operatorname{FW}$ as follows.
In terminal-record index order, and
within each terminal record in chronological order, extract every
complete layer-$0$ record from first entrance into $B_0^-$ to the
subsequent exit from $B_0^+$. Record its entrance--exit pair and
concatenate these pairs in that order. Incomplete layer-$0$ fragments
are assigned to the complementary path pieces and are not included in
the word. The map takes values in the space of finite endpoint
words; on $\mathscr C_{\rm int}$ its value lies in the finite set
$\mathcal W_k$.
Because every omitted path piece stays outside the first entrance disk,
it contributes no occupation on $Q(k)$. Conditional on a first word
$w$, Proposition~\ref{chIII:prop:record-product} gives the field kernel
$\mathcal K_w$ used below.

The local count calculation has the following form. Uniformly
over feasible terminal words $z$ and first-layer words $w$, with
$\alpha(w)$ in the first narrow corridor and $b_z$ in the terminal
corridor. Under the fixed terminal-word product law $P_z$ from Section~6.3,
\begin{equation}\label{chVI:eq:wordfactor}
 P_z(\{\operatorname{FW}=w\}\cap\mathscr C_{\rm int})
 =p_k(w)\mathcal L_k(\alpha(w),b_z)(1+o_k(1)).
\end{equation}
The error is relative and uniform. The sufficient exponent
conditions used here are
\begin{equation}\label{chVI:eq:exponentconditions}
 0\le\eta'<\eta<1/6,\quad 0<\gamma,\quad
 7\gamma/2<1,\quad3\gamma/2<1/2-\eta',\quad
 \gamma<1/2-\eta.
\end{equation}
They hold for \eqref{chVI:eq:exponents}.

We spell out the comparison behind \eqref{chVI:eq:wordfactor} to keep
its discrete endpoints and rare normalization visible. The
continuous occupation corridor calculation is the sharp ballot
calculation for the annular chain: its endpoint mass is
$\mathcal L_k(a,b)(1+o(1))$. The continuous occupation input is
\cite[Proposition 6.3 and Section 6.4]{chVI:LS26}.
In its square-root
coordinates the reflected Gaussian kernel of duration $t$ is
\[
 (\pi t)^{-1/2}
 \{e^{-(a-b)^2/t}-e^{-(a+b)^2/t}\}.
\]
Reflection gives its leading positive-bridge factor; the radial
change of density and the endpoint factors give
\eqref{chVI:eq:leading}. The upper endpoint powers produce an error
$O(k^{-1/2+3\eta'})$, which tends to zero for $\eta'<1/6$.
The lower and upper power barriers may be moved to any fixed
intermediate exponents between $\eta'$ and $\eta$, at vanishing
relative cost. Normalizing the Gaussian bridge tail by the positive-bridge mass
gives the required relative error. The radial-kernel error tends to zero when
$7\gamma/2<1$ and $3\gamma/2<1/2-\eta'$.

For the closed-count comparison corresponding to
\cite[Theorem 5.3]{chVI:LS26}, choose
\[
 \eta'<\eta_-<\eta<\eta_+<\zeta<\min(1/6,1/2-\gamma).
\]
Couple the occupation and count coordinates by the exact annular
occupation kernel. Outside their joint deviation event at
square-root distance $k^{1/2-\zeta}$, the occupation corridor
with exponent $\eta_-$ lies inside the closed count corridor;
the latter lies inside the occupation corridor with exponent
$\eta_+$. The endpoints are kept at their original integer
counts. The joint deviation mass is
$o(\mathcal L_k(a,b))$, uniformly in the two endpoint windows:
its exponential tail has power $1-2\zeta-\gamma>\gamma$, whereas the
endpoint leading factor is retained in the estimate. More explicitly,
after division by $\mathcal L_k(a,b)$ the smooth and jump parts,
and the low- and high-count cutoff parts, are bounded by a constant
times
\begin{equation}\label{chVI:eq:couplingerrors}
\begin{split}
 &k^{15/2}e^{8k^\gamma-k^{1-2\zeta-\gamma}/2}
 +k^{41/2}e^{40k^\gamma-k^{1-2\zeta-\gamma}/128}\\
 &\hspace{8mm}+k^{19/2}e^{40k^\gamma-k/4}
 +k^{15/2}e^{6k^2-k^4/4}.
\end{split}
\end{equation}
The first two terms follow by conditioning at the first bad
coordinate and applying the exponential occupation bound to the
exact count/occupation kernel. The endpoint completion kernel
is kept in that calculation, so the resulting estimate is
relative to the endpoint mass. For the last two terms, sum the
unrestricted transition kernel over paths entering the low or
high cutoff; its exponential tail gives the displayed bounds.
This partition exhausts all integer words of the exact
count/occupation kernel together with its endpoint-completion kernel.
Each term vanishes under $\gamma<1/2-\zeta$.
The two occupation asymptotics squeeze the count mass to the
same expression \eqref{chVI:eq:leading}, including both endpoints
of the closed integer count intervals.

The reference entrance law $\pi_k^-$ in Section~6.3 conditions on
$H_{B_0^-}<\tau_{B_1^+}$. Let $\pi_k^{-,\infty}$ denote harmonic
measure from infinity on $\partial_iB_0^-$. The strong Markov
subtraction and harmonic-measure estimates in the proof of
Lemma~\ref{chVI:lem:annular-harmonic-kernels}, applied to the radii
$e^k$, $e^{k+\mathfrak h}$ and $e^{k+2\mathfrak h}$ with the prescribed
logarithmic perturbation and rounding, give
\[
 \sup_{a\in\partial_iB_0^-}
 \left|\frac{\pi_k^-(a)}{\pi_k^{-,\infty}(a)}-1\right|
 \le Ce^{-c\mathfrak h}.
\]
Indeed the unnormalized entrance mass is
$\pi_k^{-,\infty}(a)\bigl[\mathbb P_{x_k}^{\rm pl}
(H_{B_0^-}<\tau_{B_1^+})+O(e^{-c\mathfrak h})\bigr]$;
the entrance probability is bounded away from zero, so division
preserves the relative error. The exit-law comparison uses the
same proof's harmonic gradient estimate. Thus the exact annular
entrance/exit kernel differs from the reference endpoint kernel
by a factor between
$e^{-C\varepsilon_{\mathfrak h}}$ and $e^{C\varepsilon_{\mathfrak h}}$ per recorded
transition, where $\varepsilon_{\mathfrak h}\le Ce^{-c\mathfrak h}$. There are at most
$9k^4$ paths per layer and at most $k^2$ layers on the corridor.
Multiplying the comparisons therefore costs
$\exp\{O(k^6e^{-c\mathfrak h})\}=1+o(1)$. Strong Markov decomposes the
ordered records without changing their completed-path counts;
all-return records and physical terminal paths describe the
same event on the exact-exit product kernel. This proves
\eqref{chVI:eq:wordfactor} for the original rounded domains.

\subsection{Unnormalized word comparison and field mixtures}
For each feasible terminal word $z$, write
\[
 r_z(w):=P_z(\operatorname{FW}=w,\ \mathscr C_{\rm int}).
\]
Exact record disintegration, as in the common conditional-law
construction in \cite[the proof of Lemma 5.14]{chVI:LS26}, gives, for every Borel
$B\subset\mathscr X_k$,
\[
 \mathcal M_z(B)=\sum_{w\in\mathcal W_k}r_z(w)\mathcal K_w(B).
\]
All core occupation is carried by these first-layer records;
the path pieces omitted in this disintegration stay outside the
first entrance disk. The conditional field kernel $\mathcal K_w$
therefore has no dependence on $z$. Equation~\eqref{chVI:eq:wordfactor}
means that there is one $e_k\downarrow0$ (after replacing it by
 a decreasing upper envelope if necessary) such that
 \[
 \begin{gathered}
  |r_z(w)-p_k(w)\mathcal L_k(\alpha(w),b_z)|
  \le e_kp_k(w)\mathcal L_k(\alpha(w),b_z),\\
  z,w\text{ in their prescribed corridors}.
 \end{gathered}
 \]
Zero reference weights have zero physical weight by the common
endpoint support. Since $\mathcal K_w(B)\ge0$, summation gives
\[
 \begin{aligned}
 |\mathcal M_z(B)-g_z\mathcal R_k(B)|
 &\le\sum_w |r_z(w)-g_zp_k(w)\alpha(w)e^{2\sqrt2(\alpha(w)-1)}|
                    \mathcal K_w(B)\\
 &\le e_kg_z\mathcal R_k(B).
 \end{aligned}
\]
Thus
\[
 (1-e_k)g_z\mathcal R_k\le\mathcal M_z\le(1+e_k)g_z\mathcal R_k.
\]
Here $g_z=\mathcal L_k(1,b_z)$ is defined in
\hyperref[chVI:def:terminal-measure]{Section~6.3}.
This is \eqref{chVI:eq:relative}; $g_z>0$ on the narrow terminal
corridor and cancels after conditioning.
The common mark law is the anchor-centered image of
$\mathcal R_k(\cdot\mid E_0)$, with the empty-mark convention if
$\mathcal R_k(E_0)=0$.

The soft threshold and local units deserve an explicit check.
The interior edge rate in \cite{chVI:LS26} is $1/(2\pi)$, hence total
rate $g$. Multiplication by $g$ converts its occupation field
to the present rate-one field. Thus our threshold $1$ is the
source threshold $1/g$. The cell trim is exactly $Q(k-\epsilon)$;
the annular endpoint powers are $q_i+\mathfrak h$, not $q_i$. Both match
the source's equations (5.13)--(5.16) after the single index
shift. The numerical range used here is verified directly by the sufficient
inequalities \eqref{chVI:eq:exponentconditions} in the annular
comparison.

Original raw centers are classified by whole-owner-cloud containment.
The restriction $\mathcal M_z(B)=P_z(F\in B,\mathscr C_{\rm int})$
and its comparison with $g_z\mathcal R_k$ give the common mark law
after normalization on $E_0$ and application of $\Theta$.
Selected factors also impose the short-record flag; exterior vacancy
enters $\Lambda_C$ in \eqref{chVI:eq:certificate}.
The false-center estimate bounds interior extra soft points on the
original sample. The weighted-window inequality of
\hyperref[chVI:def:weighted-window]{Appendix~U.3} is applied after
the product estimate and duration bound, with $J_I$ retained.

\appendixchapter{W}{Proof of Proposition~\ref{chVI:prop:false}}
The layer index in this appendix is one-based: its $i$
is the body's $i-1$, and its terminal layer $T$ is the body's
$T-1$. Its layer width is denoted by $\mathfrak h$; thus
$q_i=k+(i-1)\mathfrak h$, with the same closed count intervals
as in \eqref{chVI:eq:corridor}. Put $u_0=1$ for the application
to Proposition~\ref{chVI:prop:false} in Section~6.4. A \emph{host}
below is a nested core, entrance disk and exit disk, together with
the completed records between its entrance and exit boundaries.
The local ratio is a direct consequence of Appendix~T.
We verify its host parameters, then use containment, chronological
injection and certified selection to control the full false multiplicity.
The common-profile law is stated separately in Section~6.3--6.4.

\subsection{Host parameters and the all-count ratio}
Use $\gamma_f=1/20$, $\eta=1/12$, $\eta'=1/24$ and $\nu=1/3$.
The lower raw corridors imply the weaker displacement bound with exponent $\nu$.
These parameters are independent of those in
Sections~III--IV. Set
\begin{equation}\label{chVI:eq:R13.1}
 \mathfrak h=k^{\gamma_f},\qquad
 T=\left\lfloor\frac{k^2-k}{\mathfrak h}\right\rfloor,\qquad
 q_i=k+(i-1)\mathfrak h.
\end{equation}
For $3\le i\le T+3$ one has $k\le q_i\le k^2+3\mathfrak h$.
Define the open auxiliary host disks by
\begin{equation}\label{chVI:eq:R13.2}
 A_y=B(y,e^{q_i+e^{-\mathfrak h}}),\quad
 D_y=B(y,e^{q_i+\mathfrak h-5e^{-\mathfrak h}}),\quad
 \mathsf K_y=B(y,e^{q_i-\mathfrak h}).
\end{equation}
For the product of $m$ variable-duration first-exit paths with arbitrary
feasible words, write $q=q_i$, $w=\sqrt{\mathfrak h m}$ and
$S=\{\min_{\mathsf K_y}L^{\rm loc}\le u_0\}$, where $u_0>0$
becomes $u_0/g$ in LS units. Let
\begin{equation}\label{chVI:eq:R13.8}
 H=\#\{(x,x')\in\mathsf K_y^2:L^{\rm loc}(x),L^{\rm loc}(x')\le u_0,
                                  \ |x-x'|_2\ge e^{q-6\mathfrak h}\}.
\end{equation}
Put $p=q-\mathfrak h$ and choose
$0<\gamma_-<\gamma_f/2<\gamma_f<\gamma_+<\nu<1/2$.
Uniformly over the host scales, $p^{\gamma_-}\le\mathfrak h\le p^{\gamma_+}$
for large $k$. The entrance and exit logarithms are
$p+\mathfrak h+e^{-\mathfrak h}$ and
$p+2\mathfrak h-5e^{-\mathfrak h}$, within the radial perturbations
allowed in Appendix~T. Moreover,
$e^{q-6\mathfrak h}=e^{p-5\mathfrak h}$ and
\[
 w\ge\sqrt2q+(q+\mathfrak h)^\nu
 \quad\Longrightarrow\quad w-\sqrt2p\ge p^\nu.
\]
Thus \eqref{chVI:eq:R17.15}, with core scale $p$ and separation constant $5$,
gives
\begin{equation}\label{chVI:eq:R13.13}
 \sup_{\substack{m,\ \mathrm{feasible\ endpoint\ words}\\
       \sqrt{\mathfrak h m}\ge\sqrt2q+(q+\mathfrak h)^\nu}}
 E_m(H\mid S)\le e^{-cq^\nu}.
\end{equation}
Here $p/q\to1$ uniformly, so only the constant changes; the estimate
is uniform over all recentered counts in the displayed range.

\subsection{Host containment and anchorwise multiplicity}
Return to the LS layers
\begin{equation}
B^-[x;i]=\mathbb B(x,q_i),\qquad
B^+[x;i]=\mathbb B(x,q_i+{\mathfrak h}-4e^{-{\mathfrak h}}).
\label{chVI:eq:R13.14}
\end{equation}
For a fixed disjoint subfamily of raw-center domains, remove only
the complete layer-$T$ records.
The exact narrow terminal word is exterior; layers $T+1,\ldots,T+5$ are exterior guards.
Layer $1$ is narrow and layers $2,\ldots,T-1$ form the interior wide corridor.
For an owner $z$, put
$\mathcal R_z^{\rm phys}=\{x\in\mathcal R^{\rm phys}:z_*(x)=z\}$
and $C_z=S_t\cap P_z^{\rm buf}$. Consider owners with
$\mathcal R_z^{\rm phys}\ne\varnothing$.
Exterior vacancy gives $C_z\subset B^+[x;T]$ for each raw center $x$.
The anchor of the raw center lies in $x+Q(k)$ at
boundary distance at least $\vartheta e^k$.
For fixed $\epsilon>0$, choose $\vartheta=(1-e^{-\epsilon})/4$. Then
\begin{equation}
-\log\vartheta=o({\mathfrak h}).
\label{chVI:eq:R13.15}
\end{equation}
\phantomsection\label{chVI:def:first-fitting-host}
Let $i$ be the first $j\ge3$ for which some
$x\in\mathcal R_z^{\rm phys}$ satisfies $C_z\subset B^-[x;j-2]$;
then $i\le T+3$. Choose the first such $x$ in the fixed order,
then the nearest $y\in\mathcal G_{N,q_i-\mathfrak h}(o)$, breaking
ties in that order. The first fitting host is
$(\mathsf K_y,A_y,D_y)$ from \eqref{chVI:eq:R13.2}.
Its layer index exceeds the body's index by one.
Since $|x-y|\le0.72e^{q_i-\mathfrak h}$ for large $k$,
\begin{equation}
C_z\subset\mathbb B(x,q_i-2{\mathfrak h})\subset {\mathsf K}_y,
\quad
B^-[x;i]\subset A_y,\quad
D_y\subset B^+[x;i].
\label{chVI:eq:R13.16}
\end{equation}
The three margins are respectively $0.72+e^{-\mathfrak h}<1$,
$e^{q_i}(e^{e^{-\mathfrak h}}-1)>e^{q_i-\mathfrak h}$, and a quantity comparable to $e^{q_i}$;
they absorb lattice rounding.
Use the raw crossing count $m_{i-1}(x)$ of Section~6.3 and
write $N^{\rm in}[y;i]$ for the completed crossing count from $A_y$ to
exit from $D_y$. For paths starting and ending outside both outer
domains,
\begin{equation}
N^{\rm in}[y;i]\ge m_{i-1}(x).
\label{chVI:eq:R13.17}
\end{equation}
Associate to each raw crossing the host crossing whose entrance into $A_y$ is the last such
entrance before the raw entrance time and whose subsequent exit from $D_y$ occurs before the
raw exit from $B^+[x;i]$. The containments in \eqref{chVI:eq:R13.16} guarantee existence.
Between distinct raw crossings the path is outside $B^+[x;i]\supset D_y$, so these associated
host crossings are distinct. This chronological injection proves
\eqref{chVI:eq:R13.17}.
Root-free loops meet the endpoint condition.
The raw lower corridor implies the host count hypothesis of \eqref{chVI:eq:R13.13}; an upper count comparison is unnecessary.

If $i\ge4$, minimality gives a cloud point outside $\mathbb B(x,q_i-3\mathfrak h)$;
the raw anchor is within $(1/\sqrt2+o(1))e^k$ of $x$, yielding a separated pair.
If $i=3$ and a false candidate exists, its anchor and a soft point outside its cell are at distance
at least $\vartheta e^k\ge e^{q_i-4\mathfrak h}$ for large $k$.
For either pair $x_1,x_2$, each raw anchor in $C_z\subset\mathsf K_y$ is at distance at least
$|x_1-x_2|/2$ from one endpoint.
Distinct cells have distinct anchors, hence distinct ordered pairs.
With the separation threshold in \eqref{chVI:eq:R13.8},
\begin{equation}
|\{x\in\mathcal F_N^{\rm false}:z_*(x)=z\}|\le H_{i,y},
\label{chVI:eq:R13.18}
\end{equation}
where $H_{i,y}$ is the pair count \eqref{chVI:eq:R13.8} for this
host. On the original
containment event the completed records account for all occupation
in $\mathsf K_y$; hence $L^{\rm loc}=\ell_t$ there.

\subsection{Certified host selection and window transfer}
Use the block and duration caps $M_N=\lceil e^{L/8}\rceil$,
$d_N=e^{L/2}$ from Appendix~U. For each deterministic layer $i$
and macro color, the candidate set for an owner $z$ is
\[
 \{y\in\mathcal G_{N,q_i-\mathfrak h}(o):
 D_y\cup\partial_oD_y\subset P_z^{\rm buf},\ r\notin D_y\}.
\]
A host is eligible when its exact count is in the lower range of
\eqref{chVI:eq:R13.13} and at most $M_N$. Color the domains and their one-step boundaries
with at most $e^{C\mathfrak h}$ colors so each color is disjoint.
All eligibility conditions are measurable in the simultaneous collapsed exterior.
Write $S_y:=\{\min_{\mathsf K_y}L^{\rm loc}\le u_0\}$ and
let $A_{d,y}$ be the event that each host-record duration is at
most $d_N$. The path guard $\mathsf G_{N,k}$ requires, for every
candidate host in every layer $3\le i\le T+3$, that $A_{d,y}$ holds
and that $N^{\rm in}[y;i]\le M_N$ whenever $S_y$ holds.
The count and exit-time argument of Appendix~U applies to these
$A_y,D_y$: at fixed $k$ there are $O_k(n)$ candidates, with disk
radii $O_k(1)$. It gives
$\mathsf Q_{N,r}^I(\mathsf G_{N,k}^c)=o_N(1)$ uniformly in $r,o$
and the windows \eqref{chVI:eq:window}.
In each locator color use the local flag
$S_y\cap A_{d,y}$ and the selector
$\operatorname{Sel}_{M_B}$ defined in
\hyperref[chVI:def:host-selection]{Appendix~U}, with these candidates
and this eligibility test. The certificate in
\hyperref[chVI:def:host-selection]{Appendix~U.3} has exactly these
local flags, so Lemma~\ref{chIII:lem:selection} applies on every
positive-mass exterior stratum.

The host constructed in \eqref{chVI:eq:R13.16} contains the entire original buffer cloud in
its core. On $\mathsf G_{N,k}$ and the active cap of Appendix~U,
it is the unique positive same-color flag for that owner. The
first-positive selector therefore selects it. The first fitting layer and
\eqref{chVI:eq:R13.18} are used only for this charging step; neither is conditioned on locally.
Every physical raw center is admissible by Section~6.3, so its
protective disk already gives buffer containment and root exclusion
for this host and its one-step boundary. On the active cap, the
selector bound $M_B\ge B\sqrt L$ prevents truncation of a designated host.

Chapter~III bounds the conditional duration failure by $e^{-e^{0.49L}}$ for the allowed words and
block caps. Thus the conditional mean of $H_{i,y}$ in the selected kernel is at most
$e^{-cq_i^\nu}/(1-e^{-e^{0.49L}})$.
Apply the nonnegative-test inequality of
\hyperref[chVI:def:weighted-window]{Appendix~U.3} with
$F=\sum_{i\in C}H_i$ and
$a_q=e^{-cq_i^\nu}/(1-e^{-e^{0.49L}})$. The duration displacement
is $M_BM_Nd_N=o(b_N)$. Summing the at most $e^{C\mathfrak h}$
locator colors and the fixed macro colors gives the following bound
for each fixed layer $i$. Let $C_{i,c}$ be its selected set in color $c$:
\begin{equation}
 E_{\mathsf Q_N^I}\sum_c\sum_{(z,y)\in C_{i,c}}H_{i,y}
 \le CB\sqrt L(1+o_N(1))e^{C\mathfrak h-cq_i^\nu}.
 \label{chVI:eq:R13.19}
\end{equation}
The original charging relation is the pointwise inequality
\[
 |\mathcal F_N^{\rm false}|
 \one_{\{\#\mathbb W_N^{\rm mac}\le B\sqrt L\}\cap\mathsf G_{N,k}}
 \le\sum_{i=3}^{T+3}\sum_c\sum_{(z,y)\in C_{i,c}}H_{i,y}.
\]
The right-hand side depends only on the selected-host pair counts
$H_{i,y}$; the preceding deterministic domination has eliminated the
first-fitting selector. Applying the selected-kernel bound gives
\begin{equation}
 \begin{aligned}
 &E_{\mathsf Q_N^I}
 [|\mathcal F_N^{\rm false}|;\ \#\mathbb W_N^{\rm mac}(t)\le B\sqrt L,\ \mathsf G_{N,k}]\\
 &\qquad\le CB\sqrt L(1+o_N(1))
       \sum_{i=3}^{T+3}e^{C\mathfrak h-cq_i^\nu}.
 \end{aligned}
 \label{chVI:eq:R13.20}
\end{equation}
There are $O(k^{2-\gamma_f})$ layers and $\nu>\gamma_f$, so the right side divided by
$\sqrt L$ vanishes as $k\to\infty$ at fixed $B$.
First choose $B$ by active-count tightness \eqref{chVI:eq:R18.7}.
For this fixed $B$, choose $k$ to make the normalized bound small,
then choose $N$ for the guard and fixed-$k$ estimates.
Markov's inequality gives the normalized false-multiplicity bound.

\subsection{Scope of the conclusion}
The local numerator and all-count denominator are both proved in
Appendix~T; this section uses only their ratio, one-sided
host containment, and weighted false control.

Repelled coverage, common-field conditioning, and final identification are proved separately in Appendix~X, Section~6.3--6.4, and Section~6.6.

\appendixchapter{X}{Proof of Proposition~\ref{chVI:prop:coverage}}
The layer index in this appendix is one-based: its $i$
is the body's $i-1$, and its terminal layer $T$ is the body's
$T-1$. Write $\mathfrak h$ for its layer width, so that
$q_i=k+(i-1)\mathfrak h$ and all closed count intervals are those
in \eqref{chVI:eq:corridor}. Set $D_A:=R_N$, $u_0:=1$, and let
$h:=\phi-\phi(r)$ be the independent difference-pinned Gaussian
field, whose law is denoted by $P_h$.
Use Appendix~T's all-count denominator, Chapter~IV's Kac/Bernstein bounds, and Appendix~Y's low-endpoint repulsion and shell domination.
We apply \cite[Proposition 5.8]{chVI:LS26} to central polynomial
second-layer counts within its bound $e^{ck^{\gamma_f}}$.

\subsection{Target profile and a common independent grid translation}
Fix $K$ and use $\gamma_f=1/20$, $\eta=1/12$ and $\eta'=1/24$.
Set
\begin{equation*}
{\mathfrak h}=k^{\gamma_f},\quad q_i=k+(i-1){\mathfrak h},\quad
T=\lfloor(k^2-k)/{\mathfrak h}\rfloor.
\end{equation*}
Use the Appendix~W/Section~6.3--6.4 layers
\begin{equation}
B_i^-(x)=\mathbb B(x,q_i),\qquad
B_i^+(x)=\mathbb B(x,q_i+{\mathfrak h}-4e^{-\mathfrak h}),\qquad
W_i(x)=\sqrt{{\mathfrak h}\,m_i(x)}.
\label{chVI:eq:R21.1}
\end{equation}
Here $m_i(x)$ denotes the completed crossing count for these
one-based disks; it is $m_{i-1}(x)$ in Section~6.3. This local
indexing convention applies through layer $T+5$.
The first and $T$th layers use $\eta'$; other layers through $T+5$ use $\eta$.
The exact terminal layer is $T$, with layers $T+1,\ldots,T+5$ exterior guards.
Let $A\sim\operatorname{Unif}(V_N)$ be the origin $o$ of
$\Qaug$, independent of path and field, and use
\begin{equation*}
\mathcal G_{N,k}(A)\quad\text{for raw cells},\qquad
\mathcal G_{N,p}(A)\quad\text{for each integer-scale Gaussian grid }p.
\end{equation*}
Use the fixed lifted atlas $(U_\alpha)$ of Section~6.2; the random
origin $A$ translates only the locator grids. Chapter~II's root-uniform
estimates are therefore applied on one deterministic atlas, uniformly
in the grid phase.

For a retained parent $\mathcal C\in\mathcal P_{N,K}^{\rm ret}$,
with canonical representative $a_{\mathcal C}$, consider the cells
with centers $x\in\mathcal X_{N,k}^{\rm adm}(r,A)$ whose local
field $F_x$ satisfies $E_0$ in \eqref{chVI:eq:localflag} and
\[
 x+b(F_x)=a_{\mathcal C}.
\]
If this set is nonempty, call $\mathcal C$ geometrically designated
and let $x(\mathcal C)$ be its first cell center in the fixed cell
order. This choice uses only the path and grid origin; count
corridors and record-duration bounds are not part of the choice.
In the following profile estimates fix such a parent and write
$x=x(\mathcal C)$. Since $a_{\mathcal C}\in x+Q(k)$, admissibility
gives $z_*(x)=z_*(a_{\mathcal C})$, and hence
\[
 S_t\cap P_{z_*(x)}^{\rm buf}=\mathcal C
\]
by \eqref{chVI:eq:retained}. The admissibility strips of width
$w_k=2(R_k^{\rm ext}+1)+2$ in Section~6.3 keep the later hosts
and their one-step boundaries inside that buffer and away from
the root. Here $B_{T+5}^+$ is the body's outermost disk
$B_{T+4}^+$. Parent-count tightness is
Theorem~\ref{chIV:thm:tight}. Define
\begin{equation}
\theta=\eta'/4,\qquad \theta_1=\theta/2,\qquad \delta=\theta/4.
\label{chVI:eq:R21.2}
\end{equation}
We first prove, except for $o_P(\sqrt L)$ geometrically designated
parents, simultaneously for $2\le i\le T+5$,
\begin{equation}
\sqrt2q_i+q_i^{1/2-\theta}<W_i(x)
<\sqrt2q_i+q_i^{1/2+\theta}.
\label{chVI:eq:R21.3}
\end{equation}
then recover the first layer separately.

\subsection{Admissible Gaussian locators}
For layer $2$ take $p=k$ and use the original second-layer domains
$B_2^-(y),B_2^+(y)$ and count $m_2(y)$. At this layer both
directional implications use this same kernel, from entrance into
$B_2^-(y)$ to first exit from $B_2^+(y)$; both the in-host count
and the out-host count mean $m_2(y)$.
For the geometrically designated parent fixed above, the designated
locator is $y=x(\mathcal C)=x\in\mathcal G_{N,k}(A)$.
Its core $B(y,e^k)$ contains the cell, the harmonic circle at $e^{k+1}$
lies inside $B_2^-(y)$, and the second-layer records reconstruct the
occupation on this circle.
For $i\ge3$, define
\begin{equation}
\begin{gathered}
p_i=\lfloor q_i-\mathfrak h\rfloor,\\
y=\min_{\rm lex}\operatorname*{arg\,min}_{z\in\mathcal G_{N,p_i}(A)}d_2(z,a_{\mathcal C}).
\end{gathered}
\label{chVI:eq:R21.4}
\end{equation}
so $q_i-p_i=\mathfrak h+d_i$, $0\le d_i<1$.
Since the raw center is within $0.71e^k$ of the anchor, for large $k$,
\begin{equation}
|y-x|\le0.72e^{p_i}+0.71e^k<0.74e^{q_i-{\mathfrak h}},
\qquad\mathcal C\subset {\mathsf K}_y:=B(y,e^{p_i}).
\label{chVI:eq:R21.5}
\end{equation}
Define
\begin{equation}
\begin{array}{lll}
A_y^{\rm in}=B(y,e^{q_i+e^{-\mathfrak h}}),&
D_y^{\rm in}=B(y,e^{q_i+{\mathfrak h}-5e^{-\mathfrak h}}),\\
A_y^{\rm out}=B(y,e^{q_i-e^{-\mathfrak h}}),&
D_y^{\rm out}=B(y,e^{q_i+{\mathfrak h}-3e^{-\mathfrak h}}).
\end{array}
\label{chVI:eq:R21.6}
\end{equation}
The radial margins imply
\begin{equation}
B_i^-(x)\subset A_y^{\rm in},\quad D_y^{\rm in}\subset B_i^+(x),
\qquad
A_y^{\rm out}\subset B_i^-(x),\quad B_i^+(x)\subset D_y^{\rm out}.
\label{chVI:eq:R21.7}
\end{equation}
The inner margin is $(1+o(1))e^{q_i-\mathfrak h}$ and the outer one is $(1+o(1))e^{q_i}$,
both absorbing center displacement and rounding.
Let $m_y^{\rm in}$ and $m_y^{\rm out}$ count the completed
records from $A_y^{\rm in}$ to exit from $D_y^{\rm in}$ and from
$A_y^{\rm out}$ to exit from $D_y^{\rm out}$, respectively.
Chronological injection of root-free loop crossings gives
\begin{equation}
m_y^{\rm out}\le m_i(x)\le m_y^{\rm in}.
\label{chVI:eq:R21.8}
\end{equation}
Use the single integer-scale radial variable
\begin{equation}
H_y^2=\frac1g\sum_\xi\omega_{p_i,y}(\xi)\ell_t(\xi).
\label{chVI:eq:R21.9}
\end{equation}
Its circle is inside both host entrances.
Use the deterministic one-sided count comparisons with Appendix~T
at core parameter $p_i$ and common entrance/exit shift $d_i$;
the width is $\mathfrak h+O(e^{-\mathfrak h})$.

\subsection{Two directional local implications}
In an exact host kernel, write $W=\sqrt{\mathfrak h m}$ and
$s=W-\sqrt2p$, and condition on physical core softness
$S:=\{\min_{B(y,e^p)}L^{\rm loc}\le u_0\}$. Under the rate-$g$
time change this is exactly $\{\min_{B(y,e^p)}\ell\le u_0/g\}$,
so every invocation of Appendix~T in this subsection uses the fixed
softness parameter $u=u_0/g$. In this kernel
$H_y^2=g^{-1}\sum_\xi\omega_{p,y}(\xi)L^{\rm loc}(\xi)$;
on the original root-free loop it agrees with
\eqref{chVI:eq:R21.9}.
Appendix~T gives
\begin{equation}
P_m(S)\ge
\exp\{-2\sqrt2s_+-s_+^2/p-C{\mathfrak h}(1+W^2/p^2)\}.
\label{chVI:eq:R21.10}
\end{equation}
Uniformly over layers, $p\asymp q_i$, $q_i-p=O(\mathfrak h+1)$, and
$\mathfrak h=o(p^{1/2-\theta})$.
For a downward raw failure $W_i\le\sqrt2q_i+q_i^{1/2-\theta}$, use the out-host
(the original second-layer kernel for layer $2$).
The count upper bound and one-sided Bernstein give
\begin{equation}
P_m\{H_y>\sqrt2p+p^{1/2-\theta_1}\mid S\}
\le e^{-cp^{3/4}}.
\label{chVI:eq:R21.11}
\end{equation}
The required deviation is at least $cp^{1/2-\theta_1}$, with unconditioned exponent
$-cp^{1-2\theta_1}/\mathfrak h$; dividing by \eqref{chVI:eq:R21.10} adds at most
$Cp^{1/2-\theta}+C\mathfrak h$. The chosen exponent margins leave a bound stronger than $-cp^{3/4}$.

For an upward raw failure $W_i\ge\sqrt2q_i+q_i^{1/2+\theta}$,
use the in-host (the original second-layer kernel for layer $2$).
For its actual count,
\begin{equation}
P_m\{H_y<\sqrt2p+p^{1/2+\theta_1}\mid S\}
\le e^{-cp^{3/4}}.
\label{chVI:eq:R21.12}
\end{equation}
uniformly over all higher counts.
Indeed, $s\ge cp^{1/2+\theta}$ and the required downward deviation is at least $cs$;
the conditional log-probability is at most
\begin{equation}
-cs^2/{\mathfrak h}+C(s+s^2/p+{\mathfrak h}+{\mathfrak h}W^2/p^2)\le-c's^2/{\mathfrak h}.
\label{chVI:eq:R21.13}
\end{equation}
Absorb errors using $\mathfrak h=o(p)$ and $s/\mathfrak h\to\infty$.

\subsection{Certified transfer of the two local exceptions}
Fix one deterministic layer $i\ge3$, a host type
$\sigma\in\{\mathrm{in},\mathrm{out}\}$, a macro color and a
 locator color. For an owner $z$ the candidate list is
 \[
 \begin{split}
  I_i^\sigma(z):=\{y\in\mathcal G_{N,p_i}(A):\;&
  D_y^\sigma\cup\partial_oD_y^\sigma\subset P_z^{\rm buf},\ r\notin D_y^\sigma,\\
  &y\text{ has the prescribed locator color}\}.
 \end{split}
 \]
For layer $2$ define separately
\[
 \begin{split}
 I_2(z):=\{y\in\mathcal G_{N,k}(A):\;&
 B_2^+(y)\cup\partial_oB_2^+(y)\subset P_z^{\rm buf},\quad
 r\notin B_2^+(y),\\
 &y\text{ has the prescribed locator color}\}.
 \end{split}
\]
Its complete-record coordinate uses the original kernel from entrance
into $B_2^-(y)$ to first exit from $B_2^+(y)$, with count $m_2(y)$.
The exterior eligibility is the current second-layer count range and
$m_2(y)\le M_N$; the local flag is physical core softness intersected
with the duration event. This is the same multi-host/certified selector
of Lemma~\ref{chIII:lem:multihost} and Appendix~U.
Use the block cap $M_N=\lceil e^{L/8}\rceil$ and individual
duration cap $d_N=e^{L/2}$ from Appendix~U. The count range and
block cap are functions of the simultaneous collapsed exterior;
the local flag is $S$ intersected with the event that every
record duration is at most $d_N$. The original path guard imposes Appendix~W's guard requirements
on all these candidate lists: both host types for $i\ge3$ and
the original second-layer hosts used in both directions for $i=2$.
The same fixed-$k$ union bound gives failure probability $o_N(1)$.
Apply $\operatorname{Sel}_{M_B}$ from
\hyperref[chVI:def:host-selection]{Appendix~U} to these candidates,
exterior count eligibility tests and the flags
$S\cap\bigcap_j\{V_j\le d_N\}$.
Intrinsic-scale, parent, failing-layer and full-cloud conditions
are used in the original-path charging argument below.
The certificate in \hyperref[chVI:def:host-selection]{Appendix~U.3}
and Lemma~\ref{chIII:lem:selection} supply its conditional product law.

For a geometrically designated parent on the original event being
charged, the designated host contains its entire buffer soft cloud.
The other same-color cores in the buffer are therefore soft-empty.
The first-positive selector captures that host on the original active
cap and guard. Distinct retained parents have distinct stable owners,
so the parent assignment is injective and charges each actual
implication failure to one selected local bad event.

Let $\mathscr E_{K,k}$ be the set of actual parents whose charged
local implication fails, and let $G_{\rm cap}$ denote the intersection
of the original active cap and guard. Write $Z_{\rm sel}$ for the
sum, over the layers and colors used here, of the local failure
indicators on the selected hosts. The injective charging gives
the pointwise inequality
\begin{equation*}
 |\mathscr E_{K,k}|\one_{G_{\rm cap}}\le Z_{\rm sel}.
\end{equation*}
Taking expectations yields
\begin{equation*}
 E[|\mathscr E_{K,k}|;G_{\rm cap}]\le E Z_{\rm sel}.
\end{equation*}
Apply \eqref{chVI:eq:R21.11} or \eqref{chVI:eq:R21.12} to each
selected kernel. Duration conditioning multiplies the bound by at most
$(1-e^{-e^{0.49L}})^{-1}$, as in Appendix~U.3. Apply
\hyperref[chVI:def:weighted-window]{Appendix~U.3} to the sum of
local failure indicators. Summing over the at most $O(k^{2-\gamma_f})$ layers and
$e^{C\mathfrak h}$ colors gives
\begin{equation}
 \begin{aligned}
 \frac1{\sqrt L}E_{\mathsf Q_{N,r}^I}[|\mathscr E_{K,k}|;G_{\rm cap}]
 &\le CB k^{2-\gamma_f}e^{C\mathfrak h-ck^{3/4}}+o_N(1).
 \end{aligned}
 \label{chVI:eq:R21.14}
\end{equation}
For fixed $B$, the first term on the right tends to zero as $k\to\infty$.
The expectation is under $\mathsf Q_{N,r}^I$ with the origin
fixed; averaging over $A$ gives the augmented-law version.

\subsection{A failing layer is assigned before the independent Gaussian probes}
Remove the local-implication failures counted in \eqref{chVI:eq:R21.14}. If more than
$\varepsilon\sqrt L$ remaining geometrically designated parents
fail \eqref{chVI:eq:R21.3}, select
$m=\lceil c\varepsilon\sqrt L\rceil$ of them, one failing layer and its prescribed locator
per parent, and one soft point in each parent. These choices use only the path. The two
possible failures imply respectively
$H_y\le\sqrt2p+p^{1/2-\theta_1}$ or
$H_y\ge\sqrt2p+p^{1/2+\theta_1}$.
The selected points are $D_A$-separated and root-far, and their harmonic circles lie
within $Ce^{p+1}$ of the corresponding point.

Choose $0<\chi<1/2-\delta$ and the fixed $D$ in Lemma \ref{chVI:lem:energy-probes}. For each
selected point require
\begin{equation}
 |h(x)|\le1,\qquad
 \frac1{2g}\sum_\xi\omega_{p,y}(\xi)h(\xi)^2\le p^{1+\chi},\qquad
 \|h\|^2\le DnL.
 \label{chVI:eq:R21.15}
\end{equation}
The probability, conditional on the path, of at least one success is at least
$c\varepsilon/(1+\varepsilon)$. On success the square field $W=2\ell_t+h^2$ has a point
in its locator core with $W(x)\le2u_0+1$ and radial value outside
\[
 \sqrt2p+(p^{1/2-\delta},p^{1/2+\delta}),\qquad\delta=\theta_1/2.
\]
For a downward failure the added quadratic form changes the radial square root by at
most $Cp^\chi=o(p^{1/2-\delta})$; for an upward failure it is nonnegative.
All $p$ are integers at least $k$ on the same translated grids. At fixed $k$ their
finite range lies below $L^{1/8}$ for all sufficiently large $N$.

Define the finite set of scales
$\mathcal P_k:=\{k\}\cup\{\lfloor q_i-\mathfrak h\rfloor:3\le i\le T+5\}$.
For fixed $k$ and all sufficiently large $N$,
\[
 y\in\mathcal G_{N,p}(A)\quad\Longrightarrow\quad
 \exists\alpha:\ (\alpha,y)\in\mathscr A_{N,p}(A),
 \qquad p\in\mathcal P_k.
\]
Indeed the fixed deep-cover property gives a chart with depth at least
$\delta_0N$, while
\[
 \max_{p\in\mathcal P_k}
 (e^{p+2p^{\gamma_{\rm loc}}}+1)<\delta_0N/2
\]
for fixed $k$ and large $N$. The chart index is retained as the
existential witness in the Gaussian image event; the path-side
candidate continues to be indexed by its owner and locator. Define the
square-field event
\[
 \begin{split}
 \mathcal A_{N,k}:=\{W\in[0,\infty)^{V_N}:{}&
 \exists p\in\mathcal P_k,\ (\alpha,y)\in\mathscr A_{N,p}(A),\
 x\in B(y,e^p),\\
 &W(x)\le2u_0+1,\quad
 \mathcal R_{p,y}^{\rm sq}(W)\notin
 \sqrt2p+(p^{1/2-\delta},p^{1/2+\delta})\}.
 \end{split}
\]
This image event is defined by the square field and the existence
of witnesses $(p,\alpha,y,x)$ on the fixed-origin locator grids.
Appendix~Y's all-scale repulsion and its shell transfer give, for
each fixed energy cutoff $D$,
\[
 \lim_{k\to\infty}\limsup_N\sup_{\substack{r\in V_N\\s\in I+[0,DnL/2]}}
 \nu_{N,s}^r(W_s^\sigma(\phi)\in\mathcal A_{N,k})=0.
\]
Lemma \ref{chVI:lem:probe-coarea}, with success parameter
$c\varepsilon/(1+\varepsilon)$, therefore bounds the original bad-parent event by
$C_{D,\varepsilon}\sup_{r\in V_N,\ s\in I+[0,DnL/2]}
\nu_{N,s}^r(W_s^\sigma(\phi)\in\mathcal A_{N,k})$, plus the already paid
local-implication, active-cap and root exceptions. It follows that
\begin{equation}
 \lim_{k\to\infty}\limsup_{N\to\infty}\sup_r \mathsf Q_{N,r}^I
 \{\#\{\mathcal C:\exists i\in\{2,\ldots,T+5\},\ \mathcal C\text{ fails (\eqref{chVI:eq:R21.3}) at }i\}
                    >\varepsilon\sqrt L\}=0.
 \label{chVI:eq:R21.16}
\end{equation}
The joint all-scale estimate applies directly to this one
existential image event $\mathcal A_{N,k}$.

\subsection{First-layer control from a central second layer}
For a second layer satisfying \eqref{chVI:eq:R21.3}, the count is $O(k^2/\mathfrak h)<e^{c\mathfrak h}$,
within the range of \cite[Proposition 5.8]{chVI:LS26}. Set
\begin{equation*}
\Delta_k=k^{1/2-\eta'/2}.
\end{equation*}
Conditioning on the exact second-layer word, its two tails give
$P(|W_1-W_2|>\Delta_k)\le4e^{-c\Delta_k^2/\mathfrak h}$.
Appendix~T's additional softness conditioning costs at most $e^{Ck^{1/2+\theta}}$, hence
\begin{equation}
P(|W_1-W_2|>\Delta_k\mid S,\text{exact second-layer data})
\le e^{-ck^{3/4}}.
\label{chVI:eq:R21.17}
\end{equation}
Use the same first-positive charging rule on the fixed second-layer
host list and the weighted estimate of
\hyperref[chVI:def:weighted-window]{Appendix~U.3}. The expected
number of mismatches on the original active cap and guard is at most
\[
 CB\sqrt L\,e^{C\mathfrak h-ck^{3/4}}+o_N(\sqrt L).
\]
Markov's inequality and then the cap and guard bounds give
$o_P(\sqrt L)$ mismatches in the order $N\to\infty$, $k\to\infty$.
Higher second-layer counts were already excluded by \eqref{chVI:eq:R21.16}.
Since
\begin{equation*}
{\mathfrak h}=o(k^{1/2-\eta'}),\qquad
k^{1/2-\eta'/2}=o(k^{1/2-\eta'/4}),\qquad
k^{1/2+\eta'/4}=o(k^{1/2+\eta'}),
\end{equation*}
the second-layer corridor and $|W_1-W_2|\le\Delta_k$ imply
\[
 \tfrac12 k^{1/2-\theta}<W_1-\sqrt2k<2k^{1/2+\theta}
 \quad(\theta=\eta'/4)
\]
for all sufficiently large $k$. Since $\mathfrak h=o(k)$ and $\theta<\eta'$,
this gives the exact source interval
\begin{equation}
 W_1-\sqrt2k\in[(k+\mathfrak h)^{1/2-\eta'},(k+\mathfrak h)^{1/2+\eta'}].
\label{chVI:eq:R21.18}
\end{equation}
For every other layer, $q_i\ge k$ and the stronger exponent $\theta<\zeta$,
where $\zeta=\eta'$ or $\eta$ as appropriate, give uniformly
\[
 (q_i+\mathfrak h)^{1/2-\zeta}<q_i^{1/2-\theta},\qquad
 q_i^{1/2+\theta}<(q_i+\mathfrak h)^{1/2+\zeta}.
\]
Thus \eqref{chVI:eq:R21.3} implies the closed corridors
\eqref{chVI:eq:corridor}. In integer-count coordinates they are
\[
 \begin{gathered}
 \left\lceil\frac{(\sqrt2q_i+
             (q_i+\mathfrak h)^{1/2-\zeta})^2}{\mathfrak h}\right\rceil
 \le m_i(x),\\
 m_i(x)\le
 \left\lfloor\frac{(\sqrt2q_i+
             (q_i+\mathfrak h)^{1/2+\zeta})^2}{\mathfrak h}\right\rfloor,
 \end{gathered}
\]
with both integer endpoints included.

\subsection{Physical-window coverage and trimming}
For a geometrically designated parent $\mathcal C$, the full
profile is the intersection of the closed count corridors at
its designated cell $x(\mathcal C)$ for the one-based layers
$i=1,\ldots,T+5$, with exponent
\[
 \zeta_i=\begin{cases}
 \eta',&i\in\{1,T\},\\
 \eta,&2\le i\le T+5,\ i\ne T.
 \end{cases}
\]
Here layer $T$ is terminal and $T+1,\ldots,T+5$ are the outer guards.
The layer-$1$ condition is the recovered closed corridor
\eqref{chVI:eq:R21.18}; every integer corridor includes both endpoints
as displayed above. Define the profile-failure count directly by
\[
 \mathrm{NR}_{N,K,k}:=
 \sum_{\substack{\mathcal C\in\mathcal P_{N,K}^{\rm ret}\\
                  \mathcal C\text{ geometrically designated}}}
 \one_{\{\text{at least one of these }T+5\text{ corridors fails at }x(\mathcal C)\}}.
\]
The local field flag, admissibility and anchor equality determine
which parents are geometrically designated. Their geometric losses
are charged in the subsequent trimming step, while duration loss
remains in the later passage from the physical raw family to the
actual raw family. With this definition,
\begin{equation}
\lim_{k\to\infty}\limsup_{N\to\infty}\sup_r
\mathsf Q_{N,r}^I
(\mathrm{NR}_{N,K,k}>\varepsilon\sqrt L)=0
\quad(K\text{ fixed}).
\label{chVI:eq:R21.19}
\end{equation}

Use the parent-count cap $B_{\rm par}$ selected from
Theorem~\ref{chIV:thm:tight} in the proof of
Proposition~\ref{chVI:prop:coverage}; it is distinct from the
active-owner cap $B$ used in the selector. On
$\{|\mathcal P_N(U_1(t))|\le B_{\rm par}\sqrt L\}$, averaging over the
independent origin gives expected trimming loss at most
$C(\epsilon+D_Ke^{-k})B_{\rm par}\sqrt L$. A fixed vertex falls in a
coordinate trim strip with probability at most
$C\epsilon+CD_Ke^{-k}+O(e^k/N)$, by counting residues in the
fundamental square. Sum this bound over the capped parent anchors.
The $O(e^k/N)$ term vanishes with $N$ at fixed $k$.
Macro and seam boundary strips and the root-near region are paid
by the one-point soft estimate. First make the cap failure small,
then choose $\epsilon$, then $k$, and finally $N$.

The center obtained this way satisfies the field flag
\eqref{chVI:eq:localflag}; the entire owner cloud is its retained
parent, so it also satisfies exterior vacancy and is not false.
All the preceding profile estimates concern the physical raw
family. At positive closed-loop times, chronological decomposition
identifies its completed paths with the terminal records used
in the actual raw family. The uniform fixed-$k$ duration tail in
Proposition~\ref{chVI:prop:common}, or its unconditional guard version,
makes the duration exception negligible. For every remaining
parent $\mathcal C$, the designated cell center $x(\mathcal C)$
is therefore a true raw center, and
$S_t\cap P_{z_*(x(\mathcal C))}^{\rm buf}=\mathcal C$.
Thus it is assigned to $\mathcal C$ in the sense of Section~6.5.
This proves coverage by the original family $\mathcal R$ and
completes the proof of Proposition~\ref{chVI:prop:coverage}.

\appendixchapter{Y}{Gaussian estimates for Propositions~\ref{chVI:prop:omitted} and~\ref{chVI:prop:coverage}}
These are the Gaussian inputs to the omitted-activity and retained-coverage
proofs in Appendices~U and~X.
Throughout this appendix, $\phi$ has the mean-zero Gaussian law
$P_\phi$ with covariance $G_N$, and $h:=\phi-\phi(r)$ has law
$P_h$; their expectations are $E_\phi,E_h$. The path is held
fixed when $P_h$ is used. Put $D_A:=R_N$;
``root-far'' means at $d_2$-distance at least a fixed positive
multiple of $D_A$. Norms $\|h\|^2=\sum_{x\in V_N}h(x)^2$
use counting measure. The weights $\omega_{q,y}$ and their support
$\Sigma_{q,y}$ are the exit harmonic measure from Section~6.1.
The external inputs are \cite[Propositions 3.15--3.16]{chIV:LS24}. The last input applies on its prescribed deep-interior $e^q$ locator grid and, for
$0<\eta<\eta'<1/2$, gives a bound $C\exp\{-q^{1/2-\eta'}\}$.
Only those grids and finitely many charts are used.

\subsection{Fixed-root shell domination}
Use the normalized fixed-root energy-shell law $\nu_{N,s}^r$ and the
square-image radial functional $\mathcal R_{q,y}^{\rm sq}$ defined in Appendix~U.1.
Let $\mathcal{Q}_N=N^{-2}\|\phi_N\|_2^2$, $v_N=E\mathcal{Q}_N$, $B_s=2s/N^2\asymp L^2$, and
$a_s=\sqrt{B_s-\mathcal{Q}_N}$ on $\mathcal{Q}_N<B_s$.
Under the fixed-root shell law $\nu_{N,s}^r$ of Chapters~II--III, write
\[
 W_s^\sigma(\phi)=(\phi+\sigma a_s\mathbf1)^2
 \quad(\mathcal{Q}_N<B_s).
\]
For fixed $r$, coarea gives
\begin{equation}
b_{r,s}^\sigma(\phi)
=\frac{[\phi(r)+\sigma a_s]_+}{N^2 a_s},
\qquad \sigma=\pm1.
\label{chVI:eq:R16.1}
\end{equation}
The normalizer satisfies $\varrho_{N,r}(s)\ge c/N^2$ uniformly in $r$ and central $s$:
$\mathcal{Q}_N=O_P(L)$, $|\phi(r)|=O_P(\sqrt L)$, and $a_s\asymp L$, so the positive-branch bracket
is at least $a_s/2$ with probability tending to one.
The branch sum satisfies
\begin{equation}
N^2(b_{r,s}^++b_{r,s}^-)\le 1+\frac{2|\phi(r)|}{a_s}.
\label{chVI:eq:R16.2}
\end{equation}
On $\mathcal{Q}_N\le B_s/2$ its Gaussian second moment is uniformly bounded since $a_s\ge cL$.
Cauchy--Schwarz gives, for every event $E$ in the $(\phi,\sigma)$ space,
\begin{equation}
\nu_{N,s}^r(E,\ \mathcal{Q}_N\le B_s/2)
\le C\left(\sum_{\sigma=\pm1}
P_\phi(E_\sigma)\right)^{1/2},
\label{chVI:eq:R16.3}
\end{equation}
where $E_\sigma$ is any ordinary Gaussian event containing the corresponding branch event.
It may depend on the fixed root.

Let $a_s^0=\sqrt{B_s-v_N}$, $K\le q\le L^{1/8}+O(1)$, and $0<a<1$ fixed.
For sufficiently large $K,N$,
\begin{equation*}
|a_s-a_s^0|>q^a\quad\Longrightarrow\quad \mathcal{Q}_N-v_N>cLq^a.
\end{equation*}
Indeed, the possible negative deviation $\mathcal{Q}_N-v_N\ge-v_N=-O(L)$ changes the amplitude by only $O(1)$.
Chapters~II--III's density bound and $|\phi(r)|\le N\sqrt{\mathcal{Q}_N}$ give
\begin{equation}
\nu_{N,s}^r(|a_s-a_s^0|>q^a)
\le C N L^2e^{-cLq^a}
\le C e^{-c'Lq^a},
\label{chVI:eq:R16.4}
\end{equation}
The singularity $(B_s-\mathcal{Q}_N)^{-1/2}$ is integrable as in the energy-shell density calculation of Chapter~II, Appendix~H.

\subsection{Quantitative downward radial tails}
Let $\Phi_N^{\rm LS}=\phi_N/\sqrt{2g}$ and
$m_L=\sqrt2L-\frac3{4\sqrt2}\log L$. For an admissible locator $(q,y)$, write
\begin{equation}
\mathcal R_{q,y}(f)
=\left(\sum_{\xi\in\Sigma_{q,y}}\omega_{q,y}(\xi)f(\xi)^2\right)^{1/2}.
\label{chVI:eq:R16.5}
\end{equation}
Thus $\mathcal R_{q,y}^{\rm sq}(W)=\mathcal R_{q,y}(\sqrt{W/(2g)})$
for a nonnegative square field $W$; in particular
$\mathcal R_{q,y}^{\rm sq}(F^2)=\mathcal R_{q,y}(F/\sqrt{2g})$.
For fixed $0<\nu_3<\nu_2<1/2$, there exist $c,C,K_0>0$ such that, uniformly for
$K_0\le q\le L^{1/8}+O(1)$,
\begin{equation}
P_\phi\left(
\exists y:\ \mathcal R_{q,y}({\Phi_N^{\rm LS}}+m_L)
\le\sqrt2q-q^{\nu_2}\right)
\le C e^{-c q^{\nu_3}}.
\label{chVI:eq:R16.6}
\end{equation}
Here locators range over $\mathscr A_{N,q}(o)$ in the fixed
finite deep atlas of Section~6.2, with the prescribed origin $o\in V_N$.
In a chart domain $D$, use the mean-zero Gaussian Markov
decomposition of Chapter~II: $\bar\Phi_D$ is the Dirichlet field
conditioned to have zero sum on $D$, and $H_D$ is the independent
remainder, with $\Phi_N^{\rm LS}=\bar\Phi_D+H_D$ on $D$.
If $G_D^{\rm LS}/2$ is the unconditioned Dirichlet covariance, set
\[
 \mathfrak q_D(x):=
 \frac{\sum_{z\in D}G_D^{\rm LS}(x,z)}
 {\sqrt{2\sum_{z,z'\in D}G_D^{\rm LS}(z,z')}}.
\]
Then $\Phi_D:=\bar\Phi_D+\mathfrak q_DY$, with an independent
$Y\sim N(0,1)$, has the unconditioned Dirichlet law.
Fix a deep core $D'\subset D$ containing all the locator supports
$\Sigma_{q,y}$ used in this chart. Then $\sup_{x\in D'}|H_D(x)|$
has a uniform Gaussian tail and $\|\mathfrak q_D\|_\infty\le C$.
Choose $\nu_3<\nu_*<\nu_2$ and $\nu_3/2<a<\nu_2$.
On $\sup_{x\in D'}|H_D(x)|\le q^a$, $|Y|\le1$, Minkowski implies the Dirichlet event
\begin{equation}
\mathcal R_{q,y}(\Phi_D+m_L)
\le \sqrt2q-q^{\nu_*},
\label{chVI:eq:R16.7}
\end{equation}
for large $q$. \cite[Proposition 3.16]{chIV:LS24} bounds this by $Ce^{-q^{\nu_3}}$;
the coarse-field exception is $Ce^{-cq^{2a}}$.
Divide by the fixed positive probability $P(|Y|\le1)$, using its independence of the torus event.
The $O(1)$ change of $m_L$ under a fixed chart rescaling is absorbed by the exponent margins.
The argument is uniform under an extra shift $|b|\le c_0q^a$;
negative branches follow by Gaussian symmetry.
Since $a_s^0/\sqrt{2g}=m_L+O(1)$, \eqref{chVI:eq:R16.3}--\eqref{chVI:eq:R16.4} imply the fixed-root shell estimate
\begin{equation}
\nu_{N,s}^r\left(
\exists y:\ \mathcal R_{q,y}^{\rm sq}(W_s^\sigma(\phi))
\le\sqrt2q-q^{\nu_2}\right)
\le C e^{-c q^{\nu_3}}+C e^{-c'Lq^a}.
\label{chVI:eq:R16.8}
\end{equation}
The square root in \eqref{chVI:eq:R16.3} changes only the exponent constant.
This estimate is summable after any thinning cost $e^{q^\nu}$ with $\nu<\nu_3$.

\subsection{Qualitative full repulsion at low endpoints}\label{chVI:sec:radialrepulsion}
For fixed $V>0$ and $0<\delta<1/2$, a low endpoint for a square field $W$ at
a locator $(q,y)$ is a point $x\in B(y,e^q)$ with $W(x)\le V$.
The repulsion estimate \cite[Proposition 3.15]{chIV:LS24} excludes
locators having such an endpoint while $\mathcal R_{q,y}^{\rm sq}(W)$ lies outside
\begin{equation*}
\sqrt2q+(q^{1/2-\delta},q^{1/2+\delta}),\qquad
K\le q\le L^{1/8}+O(1),
\end{equation*}
jointly over admissible locators, in the order $N\to\infty$ then $K\to\infty$.
Precisely, for a fixed $C_0<\infty$ and origin $o$, let
\[
 \begin{split}
 \mathscr B_{N,K}(V,\delta;o):=\{W:\;&\exists q\in\ZZ\cap[K,L^{1/8}+C_0],\
 (\alpha,y)\in\mathscr A_{N,q}(o),\ x\in B(y,e^q),\\
 &W(x)\le V,\quad \mathcal R_{q,y}^{\rm sq}(W)\notin
 \sqrt2q+(q^{1/2-\delta},q^{1/2+\delta})\}.
 \end{split}
\]
For each fixed central shell range $|s-s_N|\le W_0b_N$,
the transferred statement is
\[
 \lim_{K\to\infty}\limsup_{N\to\infty}
 \sup_{r,o}\sup_{|s-s_N|\le W_0b_N}
 \nu_{N,s}^r\bigl(W_s^\sigma(\phi)\in
                     \mathscr B_{N,K}(V,\delta;o)\bigr)=0.
\]
For the periodic version, truncate coarse fields at a fixed $M$.
Both the low-endpoint window and the radial norm change by $O(M)$.
Choose a narrower corridor in the imported theorem, let $K\to\infty$, and then $M\to\infty$.
Shell transfer uses \eqref{chVI:eq:R16.3}; amplitude exceptions use \eqref{chVI:eq:R16.4} with a fixed large cutoff.
This is a joint qualitative exclusion, distinct from the quantitative per-scale bound \eqref{chVI:eq:R16.8}.

\subsection{Contact of an entire small ball}
Let $\mathcal C\subset B(y,e^q)$, $q\le L^{1/8}+O(1)$, with the whole ball at distance at least $D_A$ from $r$.
For fixed $u>0$ and an independent pinned field $h=\phi-\phi(r)$,
\begin{equation}
P_h(\exists x\in\mathcal C:\ |h(x)|\le u)
\le \frac{C_u(q+1)}{\sqrt L}.
\label{chVI:eq:R16.9}
\end{equation}
Enlarge $\mathcal C$ to the ball and use Chapter~II, Appendix~D's root-conditioned equilibrium decomposition.
The common-coordinate variance is at least $cL$, so its density is bounded by $C/\sqrt L$ for every root value.
Write $s$ for the standardized root value and $m(s)$ for the
conditional mean of $h$ on $\mathcal C$. In that decomposition $\eta_{\mathcal C}$
is the residual field after removal of the common Gaussian
coordinate. Put $\operatorname{osc}_{\mathcal C} m(s):=
\max_{x\in\mathcal C}m_x(s)-\min_{x\in\mathcal C}m_x(s)$. Given the residual,
the contact interval union has length at most
\begin{equation}
2u+\max_{x\in\mathcal C}\eta_{\mathcal C}(x)-\min_{x\in\mathcal C}\eta_{\mathcal C}(x)
+\operatorname{osc}_{\mathcal C} m(s).
\label{chVI:eq:R16.10}
\end{equation}
Since $\max_x\Var\eta_{\mathcal C}(x)\le C(q+1)$ and $|\mathcal C|\le Ce^{2q}$,
the expected Gaussian maximum is at most $C(q+1)$.
The within-block root-regression oscillation is at most
$C|s|e^q/(D_A\sqrt L)$ and remains negligible after integration in standard normal $s$.
Integration of \eqref{chVI:eq:R16.10} proves \eqref{chVI:eq:R16.9}.

\subsection{A radial square-norm tail at the same endpoint}
For root-far $x\in B(y,e^q)$, Gaussian regression yields
\begin{equation}
P_h\left(
|h(x)|\le u,\
\left(\frac1{2g}\sum_\xi\omega_{q,y}(\xi)h(\xi)^2\right)^{1/2}>t
\right)
\le \frac{C_u}{\sqrt L}e^{-ct^2/q},
\qquad t\ge C_u\sqrt q.
\label{chVI:eq:R16.11}
\end{equation}
Conditional on $h(x)=z_0\in[-u,u]$, the weighted boundary-increment covariance has trace and
operator norm at most $Cq$, while the conditional mean squared norm is bounded.
The Gaussian quadratic-form tail supplies the exponential factor.
Integrate over $z_0\in[-u,u]$ using the density bound $C/\sqrt L$.
These conclusions concern Gaussian kernels or shells;
transfer to a deterministic walk time requires the path comparison below.

\subsection{Energy-truncated probes and exact transfer}
\begin{lemma}\label{chVI:lem:energy-probes}
Fix $\chi>0$. Given a path, choose points $x_1,\ldots,x_m$
that are root-far and mutually separated by at least $c_0R_N$.
Let $K_0\le q_i\le L^{1/8}+O(1)$ and let $\omega_i$ be
probability weights supported within distance $C_0e^{q_i+1}$
of $x_i$. For suitable fixed $D,K_0$, the events
\[
 B_i=\left\{|h(x_i)|\le1,\quad
 \frac1{2g}\sum_\xi\omega_i(\xi)h(\xi)^2\le q_i^{1+\chi},
 \quad\|h\|^2\le DnL\right\}
\]
satisfy, conditionally on the path,
\begin{equation}\label{chVI:eq:probe-pz}
 P_h(B_i)\ge cL^{-1/2},\quad P_h(B_i\cap B_j)\le C/L,
 \quad P_h\left(\bigcup_iB_i\right)
       \ge c\frac{m/\sqrt L}{1+m/\sqrt L}.
\end{equation}
\end{lemma}
\begin{proof}
The Green estimates give $\Var h(x_i)=2gL+O(L^{1/8})$ and
the pair covariance matrix
$gL\left(\begin{smallmatrix}2&1\\1&2\end{smallmatrix}\right)
+O(L^{1/8})$. Thus the one-point density on $[-1,1]$ is at
least $c/\sqrt L$ and the pair density is at most $C/L$.
Condition on $h(x_i)=u$, $|u|\le1$. Regression and the Green
difference bound give
\[
 E\left(\frac1{2g}\sum_\xi\omega_i(\xi)h(\xi)^2
                   \,\middle|\,h(x_i)=u\right)\le C(q_i+1).
\]
The pinned covariance matrix $C_r(x,y):=E_h[h(x)h(y)]$ satisfies
\[
 \operatorname{tr}C_r=2nv_N,\qquad
 \|C_r(\cdot,x_i)\|_2^2\le\Var h(x_i)\operatorname{tr}C_r.
\]
The conditional mean of $\|h\|^2$ is therefore at most $CnL$.
Choose $D$ large, then $K_0$ large, so each conditional Markov
failure costs at most $1/4$. Integration over $|u|\le1$ proves
the first bound. Dropping both radial and energy restrictions
gives the pair bound. For $Y=\sum_i\one_{B_i}$,
$EY\ge cm/\sqrt L$ and $EY^2\le C(m/\sqrt L+m^2/L)$.
Cauchy--Schwarz applied to these joint first and second moments
proves the final assertion; each $B_i$ retains the common energy event.
\end{proof}

\begin{lemma}\label{chVI:lem:probe-coarea}
Suppose a path event $E$ in a window $I$ satisfying
\eqref{chVI:eq:window} has a Gaussian probe with conditional
probability at least $p>0$, and every success
satisfies $\|h\|^2\le DnL$ and implies an image event
$\mathcal A(s,W)$ for $s=t+\|h\|^2/2$, $W=2\ell_t+h^2$.
Then
\begin{equation}\label{chVI:eq:probe-transfer}
 \mathsf Q_{N,r}^I(E)\le
 \frac{C_{I,D}}p
 \sup_{s\in I+[0,DnL/2]}\nu_{N,s}^r\bigl(\mathcal A(s,W_s^\sigma(\phi))\bigr).
\end{equation}
\end{lemma}
\begin{proof}
Write $B(t,\omega,h)$ for the probe-success event and
$I_D:=I+[0,DnL/2]$. The hypotheses give the pointwise inclusion
\[
 \one_{\{t\in I\}}\one_E\one_B
 \le\one_{\{s\in I_D\}}\one_{\mathcal A(s,2\ell_t+h^2)},
 \qquad s=t+\tfrac12\|h\|^2.
\]
Multiplication by $\one_{\{X_t=r\}}$, integration over the
original independent path--field law, and
Theorem~\ref{chIII:thm:coarea} yield
\[
 \begin{aligned}
 p\,\mathsf Q_{N,r}^I(E)\int_Ip_t(r,r)\,dt
 &\le\int_I E_r[\one_{\{X_t=r\}}\one_E P_h(B\mid t,\omega)]\,dt\\
 &\le\int_{I_D}\varrho_{N,r}(s)
       \nu_{N,s}^r(\mathcal A(s,W_s^\sigma(\phi)))\,ds.
 \end{aligned}
\]
Here $\varrho_{N,r}(s)\le C/n$ on the enlarged central range,
whereas $\int_Ip_t(r,r)\,dt\ge |I|/(2n)$ eventually. Also
\[
 \frac{|I_D|}{|I|}
 =1+\frac{DnL}{2|I|}\le1+\frac{D}{2g\delta}.
\]
Division gives \eqref{chVI:eq:probe-transfer}. The only test
passed through coarea is the displayed image majorant; the path
event $E$ and its selected labels remain on the source side.
\end{proof}

\begin{corollary}\label{chVI:cor:low-owner-transfer}
If a scale-$q$ event supplies at least
$c\max\{1,e^{-q^{\nu_4}}\sqrt L\}$ separated probes, each
success implying an image event with shell probability at most
$C(e^{-cq^{\nu_3}}+e^{-cLq^a})$, where
$0<\nu_4<\nu_3$, its root-far probability is at most
\[
 C_De^{q^{\nu_4}}(e^{-cq^{\nu_3}}+e^{-cLq^a}).
\]
The sum over $K\le q\le L^{1/8}+O(1)$ tends to zero as
$N\to\infty$, then $K\to\infty$.
\end{corollary}
\begin{proof}
Apply \eqref{chVI:eq:probe-pz} and \eqref{chVI:eq:probe-transfer}.
For sufficiently large $q$, $q^{\nu_4}$ is absorbed by
$cq^{\nu_3}$; the second term is still smaller on this range.
For a common root-far restriction, use the union bound as in
\eqref{chVI:eq:low-owner-union}: its complement contributes once,
outside the scale sum.
\end{proof}


\par\endgroup

\begingroup
\mergedchaptersetup{VII}{The number of unvisited components}
\let\mathscr\mathcal
\appendixchapter{Z}{Proof of the packing estimate~\eqref{chVII:eq:packingtail}}
This appendix proves \eqref{chVII:eq:packingtail}. Its probability
input from Chapter~V is only the scalar convergence
\eqref{chVII:eq:activitylimit}. We first retain a positive mass of
activity in a short loop window, then select many bounded-count
cores. A fixed positive chance of a zero in each selected
core yields a macroscopic packing. All auxiliary letters and host
exponents introduced in this appendix have local scope.

\subsection{Activity survives a short loop extension}
Recall the retained family $\mathcal P_{N,K}^{\rm ret}(t)$ and
omitted activity $O_{N,K}(t)$ from Section~6.2, using the same stable
macro owners and retention conditions. Define
\[
 A_{N,K}^{\rm ret}(t):=\frac1{\sqrt L}
 \sum_{D\in\mathcal P_{N,K}^{\rm ret}(t)}a_D^G(\ell_t;r,1/4).
\]
Thus $A_N(t)=A_{N,K}^{\rm ret}(t)+O_{N,K}(t)$.

\begin{lemma}\label{chVII:lem:terminallower}
For every $\eta>0$ there are $a>0$ and $\delta\in(0,1]$
such that
\begin{equation}\label{chVII:eq:terminalpositive}
 \mathsf Q_N^\delta(A_N(t)<a)\le\eta
\end{equation}
for all sufficiently large $N$. There are also $a'>0$ and
a fixed integer $K\ge1$ such that the retained activity of
Chapter~VI satisfies
\begin{equation}\label{chVII:eq:retainedpositive}
 \mathsf Q_N^\delta(A_{N,K}^{\rm ret}(t)<2a')\le\eta
\end{equation}
eventually, after possibly changing $a,\delta$.
\end{lemma}
\begin{proof}
It suffices to treat $0<\eta<1$. Since $\Gamma Z>0$ a.s.,
choose $a_0>0$ so small that the central activity is at least
$2a_0$ except with probability at most $\eta/4$, eventually.
This follows from \eqref{chVII:eq:activitylimit} by choosing a
continuity point of the limiting distribution near zero.
Use \eqref{chVII:eq:incA}, proved in Appendix~AA from activity
tightness, defect control and \eqref{chVII:eq:bridgehit}, to choose
a fixed $\delta>0$ such that
\[
 \mathsf Q_N^\delta(|A_N(t)-A_N(s_N)|>a_0)\le\eta/4.
\]
Prefix comparison, with error at most $\eta/4$, shows that
$\mathsf Q_N^\delta(A_N(s_N)<2a_0)\le\eta/2$.
If neither of these events occurs, $A_N(t)\ge a_0$.
This proves \eqref{chVII:eq:terminalpositive}, with room in the budget.

For \eqref{chVII:eq:retainedpositive}, first obtain a terminal
threshold $a_1>0$ with failure probability $\eta/2$.
The omitted-activity estimate of Chapter~VI permits a fixed
$K$ for which
$\mathsf Q_N^\delta(O_{N,K}(t)>a_1/3)\le\eta/2$.
Since $A_N(t)=A_{N,K}^{\rm ret}(t)+O_{N,K}(t)$, retained activity is
at least $2a_1/3$ off these two exceptions. Set $a'=a_1/3$.
All parameters are fixed before the limit in $N$.
\end{proof}

\subsection{Bounded-count hosts}
We spell out the selection estimates needed here, to distinguish
them from the finer common-profile comparison of Chapter~VI.
Use the stable macro owners and the translated locator grid
of that chapter. Within this appendix, at a fixed integer scale $q$ let
\begin{equation}\label{chVII:eq:hostscales}
 \gamma_0=1/8192,\quad h=q^{\gamma_0},\quad
 r_q=e^{q+h},\quad d_q=e^{q+2h-4e^{-h}}.
\end{equation}
For the host regions use the open Euclidean lattice disks
$B(y,R)=\{x\in V_N:d(x,y)<R\}$, as in Section~6.1.
A host centered at $y$ has core $B(y,e^q)$, entrance region
$B(y,r_q)$ and exit region $B(y,d_q)$. The small shrinkage
in $d_q$ is retained in every endpoint kernel. The count
$m_y$ is the number of complete entrance-to-exit records,
including all returns within each record. For a feasible endpoint
word $((a_i,b_i))_{i=1}^m$, write
\[
 P_m=\bigotimes_{i=1}^m\mathsf B_{B(y,d_q)}^{a_i,b_i}
\]
for its local product kernel, built from the exact first-exit laws
of Section~3.5, and $E_m$ for its expectation. Under this product
law let $T_i$ and $\ell^{(i)}$
be the duration and rate-one local-time profile of record $i$,
and put $\ell=\sum_{i=1}^m\ell^{(i)}$. The endpoint word is
suppressed in $P_m$; every estimate below is uniform over feasible
words. For $u\ge0$, define the softness and duration events
\[
 S_u:=\{\min_{x\in B(y,e^q)}\ell(x)\le u\},\qquad
 E_{N,m}:=\bigcap_{j=1}^m\{T_j\le e^{L/2}\}.
\]
For $m=0$, $\ell=0$ and both events are the whole one-point
record space. Put
\begin{equation}\label{chVII:eq:countcap}
 M(q)=\left\lceil
 \frac{\bigl(\sqrt2(q+h)+(q+h)^{3/2}\bigr)^2}{h}
 \right\rceil.
\end{equation}
It is a fixed integer once $q$ is fixed, however large it may be.

Let $\mathcal M_N$ be the macro centers of Section~6.2,
with buffers $P_z^{\rm buf}$, and put
\[
 s_q=\lfloor e^q\rfloor,\qquad
 \mathcal G_{N,q}^o
 =\{o+s_q(i,j)\pmod N:0\le i,j\le\lfloor N/s_q\rfloor\}.
\]
Fix macro and locator color maps $\chi_N^{\rm mac}$ and
$\chi_{N,q,o}^{\rm loc}$ into finite sets $\mathfrak A$ and
$\mathfrak B_q$, respectively, whose sizes are independent of $N$.
Choose the macro coloring as in Section~6.2, so that distinct
centers of the same color have Euclidean distance at least $20s$,
where $s$ is the macro mesh of Section~6.2. Choose the
locator coloring so that like-colored disks $B(y,d_q+1)$ with
distinct $y\in\mathcal G_{N,q}^o$ are disjoint. A greedy coloring
of the disk-intersection graph gives $|\mathfrak B_q|\le C e^{4h}$:
each grid point has at most $C(d_q/s_q+1)^2$ neighbors. Write
$\mathfrak C_q=\mathfrak A\times\mathfrak B_q$ and
$J(q)=|\mathfrak C_q|$.

For $\mathfrak c=(a,b)\in\mathfrak C_q$, the available host family is
\[
 \begin{split}
 \mathcal H_{q,\mathfrak c}(r,o):=\{(z,y)\in
 \mathcal M_N\times\mathcal G_{N,q}^o:\;&
 \chi_N^{\rm mac}(z)=a,\ \chi_{N,q,o}^{\rm loc}(y)=b,\\
 &B(y,d_q+1)\subset P_z^{\rm buf},\quad
 d(r,P_z^{\rm buf})\ge R_N\}.
 \end{split}
\]
For $i=(z,y)$ put $\operatorname{owner}(i):=z$ and $y_i:=y$.
The extra unit of radius includes the exit boundary. These
conditions imply root exclusion and disjoint exit domains.
For distinct owners of one color, points in their cores have
graph distance at least $R_N$ for large $N$. A retained parent
has a containing available core when $q$ is sufficiently large
relative to its fixed retention parameter $K$. Indeed, for its
stable owner $z$, $a_D\in P_z$ and retention imply
\[
 d(r,P_z^{\rm buf})
 \ge d(r,a_D)-\sup_{x\in P_z^{\rm buf}}d(x,a_D)
 >D_N^{\rm mac}/\sqrt2-5\sqrt2s>R_N
\]
eventually. Thus the entire buffer is root-far.
Choose $q$ so that $0.72e^q+\sqrt2D_K<e^q$. A nearest grid
point $y$ to $a_D$ then has $D\subset B(y,e^q)$, and, for
sufficiently large $N$,
\[
 d(y,a_D)+d_q+1<s
 \quad\Longrightarrow\quad B(y,d_q+1)\subset P_z^{\rm buf}.
\]
These deterministic containments and the exit-domain coloring
specify the available family used by the selectors below.

We use the crossing and occupation estimates of Chapter~IV
and the exact-kernel estimates of Appendix~T. They
concern the same exact
endpoint words as Chapter~III. Let $\mu_{q,y}$ be the exit
distribution from $\{x:d(x,y)<e^{q+1}\}$ for a walk started at $y$.
In local occupation units define
\[
 \lambda(x)=\ell(x)/g,\qquad
 Y=\sum_x\mu_{q,y}(x)\lambda(x),\qquad w=\sqrt{mh}.
\]
Thus $Y$ is an average of the complete recorded occupation, with
the harmonic weights of Section~4.1. Then
\begin{equation}\label{chVII:eq:hostoccupation}
 P_m(\sqrt Y<w-\Delta)
 \le 2\exp\{-c\Delta^2/h\},
 \qquad C e^{-ch}w\le\Delta\le w.
\end{equation}
The mean satisfies $E_mY=mh(1+O(e^{-ch}))$, uniformly in the
count and feasible entrance and exit words. On the original path put
\[
 Y_{q,y}(t):=\frac1g\sum_x\mu_{q,y}(x)\ell_t(x).
\]
For large $q$, the support of $\mu_{q,y}$ lies inside the entrance
region. On a closed path rooted outside the exit region,
$Y_{q,y}(t)$ therefore equals the same average of the complete-record
occupation. For $1/2<p\le1$, define the following subset of the
macro centers, with the root $r$ and origin $o$ fixed:
\begin{equation}\label{chVII:eq:highowners}
 \begin{split}
 \mathcal M_N^{\rm high}(q,p;t):=\{z\in\mathcal M_N:\;&
 d(r,P_z^{\rm buf})\ge R_N,\ \exists x\in U_1(t)\cap P_z^{\rm buf},\\
 &\exists y\in\mathcal G_{N,q}^o:\ x\in B(y,e^q),\
 \sqrt{Y_{q,y}(t)}>\sqrt2q+q^p\}.
 \end{split}
\end{equation}
The contact and radial estimates of Appendix~Y give, for every
fixed $a>0$,
\begin{equation}\label{chVII:eq:highownertail}
 \lim_{q\to\infty}\limsup_{N\to\infty}\sup_{r,o}
 \mathsf Q_{N,r}^{I}
 \bigl(|\mathcal M_N^{\rm high}(q,p;t)|\ge a\sqrt L\bigr)=0,
\end{equation}
uniformly over the admissible windows \eqref{chVII:eq:admissible},
with $W,\delta$ fixed. The active macro family
$\mathbb W_N^{\rm mac}(t)$ has $u_0=1$ here; its normalized
cardinality is tight by \eqref{chVI:eq:R18.7}.

To verify the normalization in \eqref{chVII:eq:highownertail},
fix an admissible $I=[l_N,r_N]$ and let
$E=\{|\mathcal M_N^{\rm high}(q,p;t)|\ge a\sqrt L\}$.
The choice $I=I_N^\delta$ has $W=2$. One macro color contains at
least $a\sqrt L/|\mathfrak A|$ of the centers in $E$. Choosing
one of their soft points gives an $R_N$-separated set at distance
at least $R_N$ from the root.

Define the set of nonnegative fields
\[
 \mathcal A:=\{f\in[0,\infty)^{V_N}:\exists y\in\mathcal G_{N,q}^o,\
 \exists x\in B(y,e^q):\ f(x)\le3,\
 \mathcal R_{q,y}^{\rm sq}(f)>\sqrt2q+q^p\}.
\]
The energy-truncated Gaussian probe of Appendix~Y gives
$p_{\rm det},D>0$, fixed before $q,N$. A success satisfies
$\|h\|^2\le DnL$ and, for one of the selected soft points,
$|h(x)|\le1$. Thus, with $s=t+\|h\|^2/2$ and $f=2\ell_t+h^2$,
\[
 f(x)\le3,\qquad
 \mathcal R_{q,y}^{\rm sq}(f)^2
 =Y_{q,y}(t)+\frac1{2g}\sum_x\mu_{q,y}(x)h(x)^2
 \ge Y_{q,y}(t),
\]
so $f\in\mathcal A$. The coarea identity gives
\begin{equation}\label{chVII:eq:highowner-coarea}
 \begin{split}
 p_{\rm det}\,\mathsf W_N(I;1)\mathsf Q_{N,r}^{I}(E)
 &\le\int_{l_N}^{r_N+DnL/2}\varrho_{N,r}(s)\,
 \nu_{N,s}^r\bigl(W_s^\sigma(\phi)\in\mathcal A\bigr)\,ds\\
 &\le\frac{r_N-l_N+DnL/2}{n}\,\xi.
 \end{split}
\end{equation}
Here the second bound holds for any prescribed $\xi>0$, by the
upper radial estimate and the uniform bound on $n\varrho_{N,r}(s)$,
first for large $q$ and then for large $N$.
Translation invariance gives
$\mathsf W_N(I;1)=\int_Ip_t(r,r)\,dt\ge\delta b_N/(2n)$.
Also $r_N-l_N+DnL/2\le(2W+D/(2g)+1)b_N$. Hence
\[
 \mathsf Q_{N,r}^{I}(E)
 \le\frac{2(2W+D/(2g)+1)}{p_{\rm det}\delta}\,\xi.
\]
Choose $\xi$ after $a,p_{\rm det},D,\delta$ have been fixed.
All bounds are uniform in $r,o$ and admissible $I$.

For the original loop sampled at time $t$, let $m_i=m_i(t)$
be the number of complete records of host $i$, and $T_{i,j}$
the duration of its $j$th record. Put $m_*(N)=\lceil e^{L/8}\rceil$.
For each color define
\[
 \begin{split}
 B_{\mathfrak c}^{\rm cnt}
 &:=\{\exists i\in\mathcal H_{q,\mathfrak c}(r,o):
        m_i\ge m_*(N),\ \min_{x\in B(y_i,e^q)}\ell_t(x)\le1\},\\
 B_{\mathfrak c}^{\rm dur}
 &:=\{\exists i\in\mathcal H_{q,\mathfrak c}(r,o):
        m_i\le m_*(N),\ \exists j\le m_i:\ T_{i,j}>e^{L/2}\},\\
 B^{\rm guard}&:=\bigcup_{\mathfrak c\in\mathfrak C_q}
       (B_{\mathfrak c}^{\rm cnt}\cup B_{\mathfrak c}^{\rm dur}).
 \end{split}
\]
The count exception concerns soft cores; the duration exception
concerns hosts below the count cap. Root exclusion ensures that
on a closed loop the complete records account for all occupation
in each host core.

We use the following selection rule for an exterior-measurable
eligibility set $\mathcal H^{\rm elig}\subset\mathcal H_{q,\mathfrak c}(r,o)$. Fix a total order
$\preceq$ on each $\mathcal H_{q,\mathfrak c}(r,o)$, and let
$\xi_i$ be the complete record vector of host $i$ at the sampled
terminal time. Its count is $m_i$, and $S_{u,i},E_{N,i}$ denote
the preceding events on its record space. In the following three sets the fixed eligibility set is suppressed.
For $B>0$, set
\[
 \begin{aligned}
 \mathcal A_{\mathfrak c}
   &:=\{i\in\mathcal H^{\rm elig}:\ \xi_i\in S_{1,i}\cap E_{N,i}\},\\
 \mathcal C_{\mathfrak c}^{\rm first}
   &:=\{i\in\mathcal A_{\mathfrak c}:
       i\preceq j\text{ for every }j\in\mathcal A_{\mathfrak c}
       \text{ with }\operatorname{owner}(j)=\operatorname{owner}(i)\},\\
 \mathcal C_{\mathfrak c}^{B}
   &:=\{i\in\mathcal C_{\mathfrak c}^{\rm first}:
       |\{j\in\mathcal C_{\mathfrak c}^{\rm first}:j\preceq i\}|
       \le\lceil B\sqrt L\rceil\}.
 \end{aligned}
\]
The selected cores are $B(y_i,e^q)$ for
$i\in\mathcal C_{\mathfrak c}^{B}$. This takes the first eligible
positive flag per owner, followed by a deterministic cap.
The count and availability conditions use exterior data; the
interior flag is exactly $S_{1,i}\cap E_{N,i}$.

When both choices occur, write $\mathcal C_{\mathfrak c}^{B}
(\mathcal H^{\rm elig})$ for the last set. We use
\[
 \begin{aligned}
 \mathcal H_{\mathfrak c}^{\rm high}
   &:=\{i\in\mathcal H_{q,\mathfrak c}(r,o):M(q)<m_i<m_*(N)\},\\
 \mathcal H_{\mathfrak c}^{\rm low}
   &:=\{i\in\mathcal H_{q,\mathfrak c}(r,o):m_i\le M(q)\}.
 \end{aligned}
\]

For later counting define
\[
 \begin{split}
 \mathcal P_{N,K}^{\rm exc}(t):=\{D\in\mathcal P_{N,K}^{\rm ret}(t):\ &
 \exists\mathfrak c\in\mathfrak C_q,\ i\in\mathcal H_{q,\mathfrak c}(r,o),\\
 &\operatorname{owner}(i)=z_*(a_D),\ D\subset B(y_i,e^q),\ m_i>M(q)\},
 \end{split}
\]
and
\[
 F_q(t):=\sum_{\mathfrak c\in\mathfrak C_q}
 \sum_{i\in\mathcal C_{\mathfrak c}^{B}(\mathcal H_{\mathfrak c}^{\rm high})}
 \one_{\{\sqrt{Y_{q,y_i}(t)}<\sqrt2q+q\}}.
\]
The parameters $q,K,B,r,o$ are fixed in these definitions.
\begin{lemma}[Excess-count exclusion]\label{chVII:lem:count}
Fix $K$, a central physical window as in
\eqref{chVII:eq:admissible}, and $d,\varepsilon>0$. For every
sufficiently large fixed $q$, the number of retained parents
having an available containing host with $m_y>M(q)$ is
at most $d\sqrt L$, except with probability at most
$\varepsilon$ for all sufficiently large $N$. The bound
is uniform in the grid origin. The common guard exception
$B^{\rm guard}$ also has arbitrarily small probability in this
order of limits.
\end{lemma}
\begin{proof}
Here is the count reduction with its probability budget.
If $m>M(q)$, then
\[
 w>\sqrt2q+q^{3/2}.
\]
Write $s=w-\sqrt2q$ and apply
\eqref{chVII:eq:hostoccupation} with $\Delta=s/2$. Its lower
constraint holds for large $q$, since the relative mean
error is exponentially small in $h$. Its upper constraint
is automatic. Outside its failure event,
\[
 \sqrt Y\ge w-s/2=\sqrt2q+s/2
                      \ge\sqrt2q+q,
\]
eventually. Thus an excess-count parent either has a host
whose occupation has an exceptional downward deviation,
or belongs to an owner of height above
$\sqrt2q+q^{3/4}$.

Under the soft condition the deviation probability is divided
by $P_m(S_1)$. The all-count estimate of Chapter~VI, at exponent
$\gamma_0$, gives
\[
 -\log P_m(S_1)
 \le C\{s+s^2/q+h(1+w^2/q^2)\},\qquad s\ge q^{3/2}.
\]
Equation~\eqref{chVII:eq:hostoccupation} gives
$P_m(\sqrt Y<w-s/2)\le2e^{-cs^2/h}$.
The soft-duration estimate below gives
$P_m(E_{N,m}\mid S_1)\ge1/2$ eventually for $m<m_*(N)$.
Consequently
\[
 \begin{split}
 P_m(\sqrt Y<w-s/2\mid S_1\cap E_{N,m})
 &\le\frac{P_m(\sqrt Y<w-s/2)}
 {P_m(S_1)P_m(E_{N,m}\mid S_1)}\\
 &\le4\exp\{-cs^2/h+C(s+s^2/q+h(1+w^2/q^2))\}\\
 &\le C e^{-c's^2/h}
 \end{split}
\]
for large $q$. The last bound is uniform in $m$, since
$(s+s^2/q+h(1+w^2/q^2))/(s^2/h)\to0$ uniformly for
$s\ge q^{3/2}$. We write $c$ for the smaller constant $c'>0$.

First restrict to
$|\mathbb W_N^{\rm mac}(t)|\le B\sqrt L$ and
$(\omega,t)\notin B^{\rm guard}$, where the latter event includes the
sampled time. For each deterministic color use
$\mathcal C_{\mathfrak c}^{B}(\mathcal H_{\mathfrak c}^{\rm high})$.
Since $\sqrt2q+q<w-s/2$ for large $q$, the preceding
conditional estimate bounds each indicator in $F_q(t)$.
The weighted inequality in
\hyperref[chVI:def:weighted-window]{Appendix~U.3}, followed by
Markov's inequality and the sum over $J(q)\le Ce^{4q^{\gamma_0}}$
colors, gives
\[
 (\mathsf Q_N^I\otimes\operatorname{Unif}(V_N))
 \bigl(F_q(t)>d\sqrt L/2\bigr)
 \le\frac{CB}{d}e^{4q^{\gamma_0}-cq^{3-\gamma_0}}+o_N(1).
\]
Here $q$ is fixed in $o_N(1)$; the bound tends to zero first as
$N\to\infty$ and then as $q\to\infty$.

For $D\in\mathcal P_{N,K}^{\rm exc}(t)$ choose a containing host
$i$ as in its definition. Retention gives
$U_1(t)\cap P_{z_*(a_D)}^{\rm buf}=D$.
Every other core with the same owner and color is disjoint from
$B(y_i,e^q)$ and hence contains no soft point. For
$(\omega,t)\notin B^{\rm guard}$,
$i$ has a positive softness-duration flag and $m_i<m_*(N)$;
the active cap prevents its removal by the selector. Thus
$i\in\mathcal C_{\mathfrak c}^{B}(\mathcal H_{\mathfrak c}^{\rm high})$.
If its indicator in $F_q(t)$ is zero, its owner belongs to
$\mathcal M_N^{\rm high}(q,3/4;t)$, using a soft point of $D$.
Distinct retained parents have distinct stable owners. Therefore,
when $|\mathbb W_N^{\rm mac}(t)|\le B\sqrt L$ and
$(\omega,t)\notin B^{\rm guard}$,
\begin{equation}\label{chVII:eq:excess-partition}
 |\mathcal P_{N,K}^{\rm exc}(t)|
 \le F_q(t)+|\mathcal M_N^{\rm high}(q,3/4;t)|.
\end{equation}

Choose $B$ by \eqref{chVI:eq:R18.7}.
Equation~\eqref{chVII:eq:excess-partition} gives the event inclusion
\[
 \begin{split}
 \{|\mathcal P_{N,K}^{\rm exc}(t)|>d\sqrt L\}
 \subset{}&\{|\mathbb W_N^{\rm mac}(t)|>B\sqrt L\}
 \cup B^{\rm guard}\\
 &\cup\{F_q(t)>d\sqrt L/2\}
 \cup\{|\mathcal M_N^{\rm high}(q,3/4;t)|>d\sqrt L/2\}.
 \end{split}
\]
Assign probability budget $\varepsilon/4$ to each event.
Equation~\eqref{chVII:eq:highownertail} bounds the last event;
the conditional and weighted estimates above bound the third.
The original-path count and duration bounds
of Chapters~III--IV give probability at most $N^{-17}$ for each
of $B_{\mathfrak c}^{\rm cnt}$ and $B_{\mathfrak c}^{\rm dur}$ under the
normalized root-window law. Summation over colors gives
$2J(q)N^{-17}$. For the duration bound, condition on the realized endpoint word and
integrate over its capped-length distribution; the word probabilities
sum to one, yielding the stated error. Root exclusion and disjointness are
deterministic properties of the available host family.
Choose $q$ large, and only then $N$ large.
All bounds are uniform in the original origin, so integration
against an independent origin law leaves them unchanged.
\end{proof}

For the soft-duration estimate at fixed $q$, put
$d(q)=h/q+2(h/q)^2$. For an exact product of $m$ exit
paths, endpoint-conditioned point avoidance gives
\begin{equation}\label{chVII:eq:coarseden}
 P_m(S_u)\ge e^{-md(q)},\qquad u\ge0.
\end{equation}
Indeed, let $\mathsf D=B(y,d_q)$ and
$H_{\mathsf D}(v,b)=P_v(X_{\tau_{\mathsf D}}=b)$.
Strong Markov at the first visit to $y$, followed by conditioning
on the prescribed exit $b$, gives
\[
 P^{a,b}(H_y<\tau_{\mathsf D})
 =\frac{G_{\mathsf D}^{\rm LS}(a,y)}{G_{\mathsf D}^{\rm LS}(y,y)}
   \frac{H_{\mathsf D}(y,b)}{H_{\mathsf D}(a,b)}
 =\frac{h+O(he^{-ch})}{q+2h+O(1)}\le\frac hq.
\]
Here the potential-kernel calculation and the Poisson-kernel
ratio used in \eqref{chVI:eq:R17.11} apply with the radii
\eqref{chVII:eq:hostscales}; their logarithmic shrinkage is
$O(e^{-h})$, and Lemma~\ref{chVI:lem:annular-harmonic-kernels}
gives $H_{\mathsf D}(y,b)/H_{\mathsf D}(a,b)=1+O(e^{-ch})$. The last inequality
holds uniformly for sufficiently large $q$, since
$qe^{-ch}=o(h)$. Thus
\[
 P^{a,b}(H_y\ge\tau_{\mathsf D})\ge1-h/q
 \ge\exp\{-h/q-2(h/q)^2\}=e^{-d(q)},
\]
using $h/q\le1/2$. The product avoidance event implies $S_0$.
For a constant $C_q\ge1$ depending on $q$, the exit-duration tail gives
\[
 P_m(\max_i T_i>D)\le mC_q e^{-D/C_q}.
\]
Dividing by \eqref{chVII:eq:coarseden}, with $D=e^{L/2}$,
shows that
\begin{equation}\label{chVII:eq:duration}
 P_m(\max_iT_i>e^{L/2}\mid S_u)
 \le mC_q\exp\{-e^{L/2}/C_q+md(q)\}
 \le e^{-e^{L/4}}
\end{equation}
eventually, uniformly for $m\le\lceil e^{L/8}\rceil$.
This includes the bounded-count hosts selected below.
The broader guard and high-count transfer in
Lemma~\ref{chVII:lem:count} use the count truncation before
applying this duration estimate, as in Chapter~III.

From now on selected cores use the low-count eligibility set
$\mathcal H_{\mathfrak c}^{\rm low}$; write
$\mathcal C_{\mathfrak c}^{B}$ for this selection.

\begin{lemma}[Many selected cores]\label{chVII:lem:selection}
For every $\varepsilon>0$ and every prescribed lower bound
on $q$, there are a fixed larger integer $q$, constants
$\alpha,B>0$, and $0<\delta\le1$ such that the following
holds. In the window $I_N^\delta$, except with probability
at most $\varepsilon$, some deterministic color has at
least $\alpha\sqrt L$ selected cores. In each color at
most $\lceil B\sqrt L\rceil$ cores are selected; each
selected host has count at most $M(q)$ and each of its
records has duration at most $e^{L/2}$.
\end{lemma}
\begin{proof}
Use Lemma~\ref{chVII:lem:terminallower} with failure probability
$\varepsilon/3$ to obtain retained activity at least $2a$
with fixed $a,\delta,K$. Since a retained parent satisfies $D\subset a_D+Q_K$, the
interval-union bound gives
$a_D^G(\ell_t;r,1/4)\le B_*^{-1}|Q_K|$.
Put $b_0=aB_*/|Q_K|$, and let $\mathcal P_{N,K}'(t)$ be the first
$\lfloor b_0\sqrt L\rfloor$ members of
$\mathcal P_{N,K}^{\rm ret}(t)$ in its fixed order.
They exist whenever $A_{N,K}^{\rm ret}(t)\ge2a$.

Apply Lemma~\ref{chVII:lem:count} with threshold $b_0/8$ and
failure probability $\varepsilon/3$; assign another $\varepsilon/3$
to $B^{\rm guard}$. For large $N$, on the complement of these
exceptions put
\[
 \mathcal P_{N,K}''(t):=
 \mathcal P_{N,K}'(t)\setminus\mathcal P_{N,K}^{\rm exc}(t),
 \qquad |\mathcal P_{N,K}''(t)|\ge b_0\sqrt L/4.
\]
For each $D\in\mathcal P_{N,K}''(t)$ choose the first containing
available host in a fixed order on the finite union of all colors;
write $i(D)$ for this host and $\mathfrak c(D)$ for its color.
Then $m_{i(D)}\le M(q)$ and
$\xi_{i(D)}\in S_{1,i(D)}\cap E_{N,i(D)}$ on the original path.
Some color $\mathfrak c$ satisfies
$|\{D\in\mathcal P_{N,K}''(t):\mathfrak c(D)=\mathfrak c\}|
\ge\alpha\sqrt L$, where
\[
 \alpha=\frac{b_0}{4J(q)},\qquad B=1+b_0.
\]

Apply the selection rule
$\mathcal C_{\mathfrak c}^{B}(\mathcal H_{\mathfrak c}^{\rm low})$.
For an assigned $D$, retention gives
$U_1(t)\cap P_{z_*(a_D)}^{\rm buf}=D\subset B(y_{i(D)},e^q)$;
disjointness makes $i(D)$ its owner's only positive flag of that
color. Distinct assigned parents have distinct owners. Hence
\[
 \begin{split}
 |\mathcal C_{\mathfrak c}^{B}(\mathcal H_{\mathfrak c}^{\rm low})|
 &\ge\min\bigl\{\lceil B\sqrt L\rceil,
 |\{D\in\mathcal P_{N,K}''(t):\mathfrak c(D)=\mathfrak c\}|\bigr\}\\
 &\ge\alpha\sqrt L,
 \end{split}
\]
since $B=1+b_0>\alpha$. This proves the assertion on the
complement of the three exceptions.
\end{proof}

\subsection{A strict zero under the selected kernel}
\begin{lemma}[Local clearing]\label{chVII:lem:clearing}
For sufficiently large fixed $q$, put
\begin{equation}\label{chVII:eq:clearprob}
 p_q=\tfrac12\exp\{-M(q)d(q)\}>0.
\end{equation}
Uniformly in a feasible exact endpoint word of count
$m\le M(q)$, the local kernel conditioned on softness
and the short-duration flag has probability at least
$p_q$ of having a point of zero occupation in the core,
for all sufficiently large $N$.
\end{lemma}
\begin{proof}
For the fixed endpoint word, put
$S_0=\{\min_{x\in B(y,e^q)}\ell(x)=0\}\subset S_1$.
The point-avoidance argument proving \eqref{chVII:eq:coarseden} gives
$P_m(S_0)\ge2p_q$, hence
$P_m(S_0\mid S_1)\ge P_m(S_0)\ge2p_q$.
By \eqref{chVII:eq:duration},
$P_m(E_{N,m}^c\mid S_1)\le p_q$ eventually. Consequently
\[
 \begin{aligned}
 P_m(S_0\mid S_1\cap E_{N,m})
 &=\frac{P_m(S_0\cap E_{N,m}\mid S_1)}
         {P_m(E_{N,m}\mid S_1)}\\
 &\ge P_m(S_0\mid S_1)-P_m(E_{N,m}^c\mid S_1)\ge p_q.
 \end{aligned}
\]
The conditioning event has positive probability, as the
same bounds show. For $m=0$ the empty local path family
leaves the whole core unvisited and the assertion also holds.
The local success event is
$S_0=\{\exists x\in B(y,e^q):\ell(x)=0\}$.
\end{proof}

For a realized selected list $\mathcal C$ and its feasible
endpoint words $\mathbf w=(w_i)_{i\in\mathcal C}$, define
\[
 \Pi_{\mathcal C,\mathbf w}
 :=\bigotimes_{i\in\mathcal C}
     P_{m_i,w_i}(\,\cdot\mid S_{1,i}\cap E_{N,i}),
\]
where $P_{m_i,w_i}$ is the exact product kernel for host $i$.
Lemma~\ref{chIII:lem:selection}, using the first-positive
certificate of \hyperref[chVI:def:host-selection]{Appendix~U.3}, gives this product on
each positive-mass selection layer before the loop window is
imposed. The full window retains the duration weight $J_I$.
The selector reads only the flags and exterior data, so
Lemma~\ref{chVII:lem:clearing} applies to every factor of
$\Pi_{\mathcal C,\mathbf w}$.

\begin{lemma}[Many successful cores]\label{chVII:lem:success}
With the parameters of Lemma~\ref{chVII:lem:selection}, the
probability that some color has at least $\alpha\sqrt L$
selected cores but fewer than $p_q\alpha\sqrt L/2$
successful cores tends to zero under the original
normalized loop-window law, with an independent grid origin.
\end{lemma}
\begin{proof}
Fix a color and a realized selection layer with $m$ hosts.
Write $P=\Pi_{\mathcal C,\mathbf w}$ and $E$ for its expectation,
enumerate $\mathcal C$ in its fixed order, and put
$I_i=\one_{S_{0,i}}$ for $i=1,\ldots,m$. They are independent and
$P(I_i=1)\ge p_q$. The exponential Markov inequality gives
\[
 P\left(\sum_{i=1}^m I_i<p_qm/2\right)
 \le e^{p_qm/2}\prod_{i=1}^mE(e^{-I_i})
 \le \exp\{p_qm/2-p_q(1-e^{-1})m\}
 \le e^{-p_qm/8}.
\]
For $m\ge\alpha\sqrt L$ this is bounded by
$e^{-p_q\alpha\sqrt L/8}$. The bound is uniform over
words and positive-probability selection layers.

Pass back to the physical window by
\eqref{chIII:eq:short-window}, integrated as in
\hyperref[chVI:def:weighted-window]{Appendix~U.3}. The total
selected duration in one copy is at most
\[
 d_N^{\rm sel}:=\lceil B\sqrt L\rceil M(q)e^{L/2}=o(b_N).
\]
The duration cap holds for both copies because it is part
of the flag. For each fixed color,
\[
 \begin{split}
 &(\mathsf Q_N^\delta\otimes\operatorname{Unif}(V_N))
 \left(|\mathcal C_{\mathfrak c}^{B}|\ge\alpha\sqrt L,\
 \sum_{i\in\mathcal C_{\mathfrak c}^{B}}\one_{\{\xi_i\in S_{0,i}\}}
                  <p_q\alpha\sqrt L/2\right)\\
 &\quad\le
 \frac{\mathsf W_N((I_N^\delta)^{d_N^{\rm sel}};1)}
      {\mathsf W_N(I_N^\delta;1)}e^{-p_q\alpha\sqrt L/8}
 =(1+o_N(1))e^{-p_q\alpha\sqrt L/8}.
 \end{split}
\]
Apply this bound to each deterministic color and then take the
union over the fixed $J(q)$ colors. The result is at most
$2J(q)e^{-p_q\alpha\sqrt L/8}$ eventually, which tends to zero.

On the original closed path, the root lies outside each selected
exit domain. Every core visit occurs between an entrance and its
following exit. Evaluating the record profiles on $\xi_i$ gives
\[
 \ell_t(x)=\sum_{j=1}^{m_i(t)}\ell^{(j)}(x),
 \qquad x\in B(y_i,e^q).
\]
Consequently
\[
 \xi_i\in S_{0,i}
 \quad\Longleftrightarrow\quad
 B(y_i,e^q)\cap U_0(t)\ne\varnothing.
\]
At the sampled positive closed-loop time these zeros are unvisited
vertices, almost surely, as in \eqref{chVII:eq:zerosvisited}.
\end{proof}

\subsection{From successful cores to central packing}
Fix $\eta>0$ and apply Lemmas~\ref{chVII:lem:selection} and
\ref{chVII:lem:success} with failure budgets $\eta/4$ each.
Let $b=p_q\alpha/2>0$. On the complementary event a
single color has at least $b\sqrt L$ successful cores.
Choose one zero from each. The same-color geometry gives
pairwise distances at least $R_N$. Hence
\begin{equation}\label{chVII:eq:terminalpacking}
 \mathsf Q_N^\delta
 (\Pack_N(U_0(t))<b\sqrt L)\le\eta/2
\end{equation}
eventually. Although the construction used an independent
origin, the terminal packing event does not; integrating
it out gives \eqref{chVII:eq:terminalpacking} on the path window.

For $t\ge s_N$, the inclusion $U_0(t)\subset U_0(s_N)$ gives
$\Pack_N(U_0(t))\le\Pack_N(U_0(s_N))$.
Thus a bad packing at $s_N$ implies a bad packing at $t$.
The prefix density from the proof of Lemma~\ref{chVII:lem:prefix} is
eventually bounded below by $1/2$. It follows that
\[
 \PP_{\pi_N}(\Pack_N(U_0(s_N))<b\sqrt L)
 \le2\mathsf Q_N^\delta(\Pack_N(U_0(s_N))<b\sqrt L)
 \le\eta.
\]
This proves \eqref{chVII:eq:packingtail}.

\appendixchapter{AA}{Proof of Lemma~\ref{chVII:lem:increments}}
Tightness, the vanishing wide-defect count and the endpoint-weighted
short-bridge hitting estimate of Chapters~I and IV give the
hit-and-defect proof of both increments in
Lemma~\ref{chVII:lem:increments}. Its activity case also supplies
the short-window input to Appendix~Z.

\subsection{Pathwise hit and defect budgets}
For $u\ge0$ and $v>0$, write
\[
 M_{u,v}(s)=\frac1{\sqrt L}
   \sum_{D\in\mathcal P_N(U_u(s))}a_D^G(\ell_s;r,v).
\]
Thus $M_{1,1/4}(s)=A_N(s)$. Let $z_D$ be the lexicographically first
vertex of $D$, using the representatives $\{0,\ldots,N-1\}^2$.
This choice depends only on the occupation profile at time $s$.
For $R\ge0$ put
\[
 B_{\rm gr}(z,R):=\{x\in V_N:d_1(x,z)\le R\}.
\]
Define
\begin{align}
 H_{u,v}(s,t)&=\frac1{\sqrt L}
   \sum_{D\in\mathcal P_N(U_u(s))}a_D^G(\ell_s;r,v)
   \one_{\{X_{(s,t]}\cap B_{\rm gr}(z_D,2\rho_N)\ne\varnothing\}},
                                                    \label{chVII:eq:hitbudget}\\
 D_{u,v}^{\rm can}(s)&=\frac{2c_g\sqrt v}{\sqrt L}
                         |D_N^{\rm can}(U_u(s))|.
                                                    \label{chVII:eq:defbudget}
\end{align}
Every weight in the hit budget is evaluated at time $s$.
The budgets are nonnegative. The matching inequality
\eqref{chIV:eq:matching} gives, for $s\le t$,
\begin{align}
 |A_N(t)-A_N(s)|&\le H_{1,1/4}(s,t)+D_{1,1/4}^{\rm can}(s),
                                                    \label{chVII:eq:changeA}\\
 2c_g|C_N(t)-C_N(s)|&\le H_{0,1}(s,t)+D_{0,1}^{\rm can}(s).
                                                    \label{chVII:eq:changeC}
\end{align}
The relevant deterministic matching is the following. Put
$U^-=U_u(s)$ and $U^+=U_u(t)\subset U^-$. Then
\[
 (U^-)^\circ\cap U^+\subset(U^+)^\circ,\qquad
 (U^+)^\circ\setminus(U^-)^\circ\subset D_N^{\rm can}(U^-).
\]
For each early parent $D$ with $D\cap U^+\ne\varnothing$,
all its surviving points lie in one late component, since their
pairwise distances are at most $\rho_N$. Two early parents
cannot map to the same late component, whose diameter is at most
$\rho_N<R_N$. Interval-union subadditivity gives the bound for
additional clean vertices:
\[
 \begin{split}
 \sum_{x\in(U^+)^\circ\setminus(U^-)^\circ}
 a_{\{x\}}^G(\ell_t;r,v)
 &\le2c_g\sqrt v\,|(U^+)^\circ\setminus(U^-)^\circ|\\
 &\le2c_g\sqrt v\,|D_N^{\rm can}(U^-)|.
 \end{split}
\]
On a surviving early parent, unchanged occupation gives unchanged
activity. If an early parent disappears or its occupation changes,
the path visits that parent and therefore its ball in
\eqref{chVII:eq:hitbudget}. Monotonicity of the interval-union
activity as occupation increases, together with the preceding
defect bound, gives \eqref{chVII:eq:changeA}.

For the zero count the same assignment is made between
canonical zero components. A disappearing early component
must be hit. Each nonempty early zero component has
activity at least $2c_g$ at $v=1$, because one of its
zero vertices contributes an interval of length two to
the Gaussian union. The unmatched late components satisfy
\[
 \bigl|\{D\in\mathcal P_N(U^+):
                  D\cap(U^-)^\circ=\varnothing\}\bigr|
 \le |D_N^{\rm can}(U^-)|,
\]
since each such nonempty component consists of early defect
vertices and different components are disjoint. Multiplying the
count difference by $2c_g$ gives \eqref{chVII:eq:changeC}. This
argument allows new clean
components to appear when defect witnesses are removed.

\subsection{Tightness at the zero threshold}
The input from Chapter~IV directly bounds $M_{u,v}$ for
positive $u$, including $(u,v)=(1,1/4)$ and $(1,1)$.
For zero threshold the deterministic subset comparison is
\begin{equation}\label{chVII:eq:zerotight}
 M_{0,1}(s)\le M_{1,1}(s)+D_{1,1}^{\rm can}(s),\qquad
 D_{0,1}^{\rm can}(s)\le D_{1,1}^{\rm can}(s).
\end{equation}
For the first inequality, intersect each clean zero component
with the clean soft set and assign its nonempty intersection to
the containing soft parent. Distinct zero components cannot use
the same soft parent, since its diameter is at most $\rho_N$.
Their activity on the clean soft set is therefore bounded by
$M_{1,1}(s)$. For the remaining vertices, subadditivity and the
singleton bound give
\[
 \begin{split}
 \frac1{\sqrt L}\sum_{x\in U_0(s)\cap D_N^{\rm can}(U_1(s))}
             a_{\{x\}}^G(\ell_s;r,1)
 &\le\frac{2c_g}{\sqrt L}
             |U_0(s)\cap D_N^{\rm can}(U_1(s))|\\
 &\le D_{1,1}^{\rm can}(s).
 \end{split}
\]
These two bounds give the first part of \eqref{chVII:eq:zerotight}.
For the second inequality, a zero-set defect is a soft-set
defect, since its witnessing pair is still present in
$U_1(s)$. Since $D_N^{\rm can}(U_1(s))\subset D_N^{\rm wide}(U_1(s))$,
the wide-defect estimate makes $D_{1,1}^{\rm can}(s)$ small in probability, and therefore \eqref{chVII:eq:zerotight}
proves tightness of $M_{0,1}(s_N)$.

\subsection{Weights far from both endpoints}
Consider either pair $(u,v)=(1,1/4)$ or $(0,1)$, and put
$s=s_N$. Given a cap $B>0$, retain the normalized parent
weight only on $\{M_{u,v}(s)\le B\}$ and only for parents
whose representative is at Euclidean distance at least $R_N/4$
from both $X_s$ and $r$. Explicitly, for $z=z_D$ with
$D\in\mathcal P_N(U_u(s))$, put
\[
 w_{z_D}=L^{-1/2}a_D^G(\ell_s;r,v)\,
 \one_{\{M_{u,v}(s)\le B,\ d(z_D,X_s)\ge R_N/4,\ d(z_D,r)\ge R_N/4\}},
\]
and put $w_z=0$ for other vertices. These weights are
$\mathcal F_s$-measurable, nonnegative, and satisfy
\begin{equation}\label{chVII:eq:weights}
 \sum_z w_z\le B,\qquad
 w_z>0\ \Longrightarrow\ d(z,X_s),d(z,r)\ge R_N/4.
\end{equation}
Set
\[
 F(s,t)=\sum_z w_z
       \one_{\{X_{(s,t]}\cap B_{\rm gr}(z,2\rho_N)\ne\varnothing\}}.
\]
For each fixed terminal time this is a measurable random
variable. The hit event can be written using rational
times and the endpoint on the almost sure set of positive
holding times and nonexplosion.

Lemma~\ref{chIV:lem:envelope}
implies, for either fixed $v\in\{1/4,1\}$ and uniformly over
$U\subset V_N$, $\ell\in[0,\infty)^{V_N}$ and $r,y\in V_N$,
\begin{equation}\label{chVII:eq:near}
 \frac1{\sqrt L}\sum_{\substack{D\in\mathcal P_N(U)\\
                       d(z_D,y)<R_N/4}}
                   a_D^G(\ell;r,v)=o(1).
\end{equation}
Indeed separation permits at most one such representative;
the Gaussian envelope for a component of diameter
$\rho_N$ bounds its unnormalized activity by $CL^{1/6}$.
Its normalized activity is at most $CL^{-1/3}$, which tends
to zero. It follows, by discarding the two near-endpoint
parts of \eqref{chVII:eq:hitbudget}, that for every $h>0$,
\begin{equation}\label{chVII:eq:farcomparison}
 M_{u,v}(s)\le B\quad\Longrightarrow\quad
 H_{u,v}(s,t)<h+F(s,t)
\end{equation}
almost surely for all sufficiently large $N$. The bound
is uniform in $t\ge s$; the two near parts are estimated
before any future path is exposed.

The bridge input is \eqref{chIV:eq:endpoint-hit}. There is $C>0$ such that,
for each fixed $0<\delta\le1$ and all sufficiently large
$N$, uniformly for $\delta b_N\le v'\le2\delta b_N$,
\begin{equation}\label{chVII:eq:bridgehit}
 \PP_x\bigl(X_{v'}=y,\ X_{(0,v']}\cap
             B_{\rm gr}(z,2\rho_N)\ne\varnothing\bigr)
 \le\frac{C(\delta+\eta_N)}n,
 \qquad \eta_N=\frac{L^{1/6}+1}{L},
\end{equation}
provided $d(z,x),d(z,y)\ge R_N/4$.
Conditioning on the prefix in \eqref{chVII:eq:weights}, applying
\eqref{chVII:eq:bridgehit} with $x=X_s,y=r$, and summing yields
\begin{equation}\label{chVII:eq:weightedhit}
 \EE_{\pi_N}[\one_{\{X_t=r\}}F(s,t)]
 \le\frac{CB(\delta+\eta_N)}n,
 \qquad t\in I_N^\delta.
\end{equation}
The past measurability of both the cap and the weights is
essential for this conditioning step.

\subsection{A general increment estimate}
\begin{lemma}\label{chVII:lem:genericincrement}
Let $X_N(t)$ be either of the two observables in
\eqref{chVII:eq:changeA}--\eqref{chVII:eq:changeC}. Suppose, with
the corresponding $(u,v)$ and a fixed $c>0$, that for $0\le s\le t$,
\[
 c|X_N(t)-X_N(s)|\le H_{u,v}(s,t)+D_{u,v}^{\rm can}(s).
\]
Assume $M_{u,v}(s_N)$ is tight and $D_{u,v}^{\rm can}(s_N)\to0$
in probability under the ordinary stationary law.
For every $\tau,\varepsilon>0$ there is
$\delta_0\in(0,1]$ such that for each fixed
$0<\delta\le\delta_0$,
\begin{equation}\label{chVII:eq:genericincrement}
 \mathsf Q_N^\delta(|X_N(t)-X_N(s_N)|>\tau)
 \le\varepsilon
\end{equation}
eventually.
\end{lemma}
\begin{proof}
Set $s=s_N$ and $h=c\tau/4$. Choose $B>0$ so large that the
ordinary-path probability of $M_{u,v}(s)>B$ is at most
$\varepsilon/16$ eventually. The probability of
$D_{u,v}^{\rm can}(s)>h$ is also eventually at most
$\varepsilon/16$. Their union $E$ is $\mathcal F_s$-measurable
and has probability at most $\varepsilon/8$.
Take
\[
 \delta_0=\min\left(1,\frac{\varepsilon h}{8CB}\right),
\]
where $C$ is from \eqref{chVII:eq:bridgehit}. Fix
$0<\delta\le\delta_0$ and take $N$ large enough that
$\eta_N\le\delta$ and \eqref{chVII:eq:farcomparison} holds.
Off $E$, an increment larger than $\tau$ forces
\[
 4h<c|X_N(t)-X_N(s)|\le H_{u,v}(s,t)+D_{u,v}^{\rm can}(s)
 <2h+F(s,t),
\]
so in particular $F(s,t)\ge h$.

At each fixed $t\in I_N^\delta$, the Markov property
and $p_{t-s}(x,y)\le2/n$ imply
\[
 \EE_{\pi_N}[\one_{\{X_t=r\}}\one_E]
 \le\frac2n\PP_{\pi_N}(E).
\]
Markov's inequality, with the endpoint indicator retained, gives
\[
 \begin{split}
 \EE_{\pi_N}\bigl[\one_{\{X_t=r\}}
    \one_{\{|X_N(t)-X_N(s)|>\tau\}}\bigr]
 &\le\frac2n\PP_{\pi_N}(E)
     +\frac1h\EE_{\pi_N}[\one_{\{X_t=r\}}F(s,t)]\\
 &\le\frac2n\frac{\varepsilon}{8}
     +\frac{CB(\delta+\eta_N)}{nh}\\
 &\le\frac{\varepsilon}{4n}+\frac{2CB\delta}{nh}
 \le\frac{\varepsilon}{2n}.
 \end{split}
\]
Its time-section normalizer is
$n^{-1}\sum_r p_t(r,r)\ge1/(2n)$, so each section has exceptional
probability at most $\varepsilon$. The actual increment event
is jointly measurable in path and time, by measurability of local
times and the finite canonical construction. Integrate its
unnormalized bound over $I_N^\delta$ and divide by $\mathsf W_N(I_N^\delta;1)$
to obtain \eqref{chVII:eq:genericincrement}.
\end{proof}

Apply Lemma~\ref{chVII:lem:genericincrement} first with
$X_N=A_N$, $(u,v)=(1,1/4)$ and $c=1$, and then with
$X_N=C_N$, $(u,v)=(0,1)$ and $c=2c_g$.
The positive-threshold tightness and defect estimates
come from Chapter~IV; \eqref{chVII:eq:zerotight} supplies the
zero-threshold versions. Each application is valid for
every smaller positive fixed $\delta$, so the minimum
of the two resulting $\delta_0$ values gives both
\eqref{chVII:eq:incA} and \eqref{chVII:eq:incC}. This proves
Lemma~\ref{chVII:lem:increments}.


\par\endgroup

\begingroup
\mergedchaptersetup{VIII}{The limiting cover-time distribution}
\let\mathscr\mathcal
\appendixchapter{AB}{Proof of Proposition~\ref{chVIII:prop:survival}}
We apply the capacity and killed-eigenfunction estimates
(Propositions~\ref{chI:prop:capacity} and \ref{chI:prop:killed})
and mixing (Lemma~\ref{chI:lem:mixing}) at the canonical
separation scale to prove Proposition~\ref{chVIII:prop:survival}.

\subsection{Capacity and stationary survival}
Fix an integer $k\ge1$ and $W\ge0$, and let $A\subset V_N$
have $k$ points with $d(a,b)\ge R_N/2$ whenever $a\ne b$ in $A$. For large $N$ it is a proper subset of
$V_N$. Let
\[
 G_A=(G_N(a,b))_{a,b\in A},\quad
 q_A=G_A^{-1}\one_A,\quad c_A=\sum_{a\in A}q_A(a).
\]
Proposition~\ref{chI:prop:capacity} gives $q_A\ge0$, $c_A>0$;
the Green estimates of Section~1.3 give
\[
 G_N(a,a)=gL+O(1),\qquad
 |G_N(a,b)|\le C\left(1+\log\frac{N}{d(a,b)\vee1}\right).
\]
For $a\ne b$ in $A$, \eqref{chVIII:eq:scales} makes the
last bound at most $CL^{1/8}$. Sum the equations
$G_Aq_A=\one_A$ and use the nonnegativity of the charges:
\[
 |gLc_A-k|\le C_kL^{1/8}c_A.
\]
For $N$ large enough that $C_kL^{1/8}<gL/2$, this is the
two-sided bound
\[
 \frac{k}{gL+C_kL^{1/8}}\le c_A\le
 \frac{k}{gL-C_kL^{1/8}}.
\]
Expanding the two denominators yields
\begin{equation}\label{chVIII:eq:capacity}
 c_A=\frac{k}{gL}(1+O_k(L^{-7/8})).
\end{equation}
All constants are uniform over the permitted sets $A$.

Let $\gamma_N=\sin^2(\pi/N)$ be the spectral gap,
$\lambda_A$ the principal killed eigenvalue, and
$v_A=\pi_N(\varphi_A^2)-1$ for its nonnegative
eigenfunction, extended by zero on $A$ and normalized by
$\pi_N(\varphi_A)=1$, as in Section~1.4.
With $\alpha_A=c_A/(n\gamma_N)=O_k(L^{-1})$,
Proposition~\ref{chI:prop:killed} gives
\[
 (1-\alpha_A)c_A/n\le\lambda_A\le c_A/n,
 \qquad 0\le v_A\le\frac{\alpha_A}{1-\alpha_A}.
\]
Consequently
\begin{equation}\label{chVIII:eq:eigen}
 b_N\lambda_A=k+O_k(L^{-7/8}),\qquad v_A=O_k(L^{-1}).
\end{equation}
The stationary survival expansion is
\begin{equation}\label{chVIII:eq:spectral}
 S_A(t):=\PP_{\pi_N}(H_A>t)
 =\frac{e^{-\lambda_At}}{1+v_A}+E_A(t),\qquad
 0\le E_A(t)\le v_Ae^{-\gamma_Nt}.
\end{equation}

Put $h_N=4n\log L$. For fixed $|z|\le W$ and
$j\in\{0,1,2\}$, the time $\bar u_N(z)-jh_N$ is positive
for all sufficiently large $N$. Since
$h_N/b_N=(4/g)(\log L)/L$, \eqref{chVIII:eq:eigen} gives
\[
 \lambda_A(\bar u_N(z)-jh_N)
 =\tfrac{k}{2}\log L+kz+O_{k,W}(\omega_N).
\]
Moreover $\gamma_N(\bar u_N(z)-2h_N)\ge cL\log L$ eventually.
The remainder in \eqref{chVIII:eq:spectral} is therefore
negligible on the relative scale: uniformly in $|z|\le W$,
\[
 L^{k/2}e^{kz}E_A(\bar u_N(z)-jh_N)
 \le C_{k,W}L^{k/2-1}e^{-cL\log L}=o(\omega_N).
\]
Also $(1+v_A)^{-1}=1+O_k(L^{-1})$. Thus
\begin{equation}\label{chVIII:eq:stationary}
 S_A(\bar u_N(z)-jh_N)
 =L^{-k/2}e^{-kz}(1+O_{k,W}(\omega_N)),
 \qquad j=0,1,2.
\end{equation}
The fixed-$k$ restriction is used here as well as in
\eqref{chVIII:eq:capacity}.

\subsection{Early hits and multiplicative comparison}
Fix $|z|\le W$ and assume $d(y,A)\ge R_N/8$. For a point $a\in A$, the
truncated Green identity of Chapter~I reads
\[
 \EE_y\ell_t(a)=t/n+G_N(y,a)
                 -\int_t^\infty(p_s(y,a)-n^{-1})\,ds.
\]
The Green bound at distance $R_N/8$, and the uniform
boundedness of the tail for $t\ge h_N$, give
\[
 \EE_y\ell_{2h_N}(a)\le C(L^{1/8}+\log L),\qquad
 \EE_a\ell_{h_N}(a)\ge cL.
\]
If $H_a\le h_N$, at least $h_N$ units of time remain
before $2h_N$. The strong Markov property implies
$\EE_y\ell_{2h_N}(a)\ge
\PP_y(H_a\le h_N)\EE_a\ell_{h_N}(a)$.
Sum the resulting point bounds to obtain
\begin{equation}\label{chVIII:eq:early}
 \PP_y(H_A\le h_N)
 \le C_k\frac{L^{1/8}+\log L}{L}=O_k(\omega_N).
\end{equation}

Set $T=\bar u_N(z)$. Lemma~\ref{chI:lem:mixing} gives
\[
 \sup_{x,w}|np_{h_N}(x,w)-1|\le L^{-8}.
\]
Let $f(x)=\PP_x(H_A>T-h_N)$. The Markov property gives
\[
 \PP_y(H_A>T)
 =\EE_y f(X_{h_N})-\EE_y[\one_{\{H_A\le h_N\}}f(X_{h_N})].
\]
The first term lies between $(1-L^{-8})S_A(T-h_N)$
and $(1+L^{-8})S_A(T-h_N)$. Dropping the no-hit
restriction for one more mixing interval gives
$\sup_x f(x)\le(1+L^{-8})S_A(T-2h_N)$. Hence
\begin{equation}\label{chVIII:eq:comparison}
 \begin{split}
 (1-L^{-8})S_A(T-h_N)
 -(1+L^{-8})\PP_y(H_A\le h_N)S_A(T-2h_N)
 &\le\PP_y(H_A>T)\\
 &\le(1+L^{-8})S_A(T-h_N).
 \end{split}
\end{equation}
The upper bound applies to every starting point and,
with \eqref{chVIII:eq:stationary}, proves \eqref{chVIII:eq:upper}.
For the lower bound under the distance assumption,
substitute \eqref{chVIII:eq:early} and both stationary
estimates from \eqref{chVIII:eq:stationary}. The subtracted
term is $O_k(\omega_N)$ times the desired survival
scale. This proves \eqref{chVIII:eq:relative} with relative
error $O_{k,W}(\omega_N)$ and completes the appendix.


\par\endgroup
\par\addvspace{.75\baselineskip}
\phantomsection\label{sec:availability-disclosure}
\mergedcontentsline{toc}{section}{Code availability and AI-use disclosure}
\noindent\textbf{Code availability.}\enspace
A Lean~4 formalization of the main cover-time distribution theorem
(Theorem~\ref{chVIII:thm:main}) is available at
\href{https://github.com/zyc111-234/torus-cover-time-lean}
{github.com/zyc111-234/torus-cover-time-lean}.

\par\smallskip\noindent\textbf{AI-use disclosure.}\enspace
AI-assisted tools, specifically OpenAI's ChatGPT and Codex, provided
substantive assistance with preliminary literature screening, proof
exploration, Lean~4 code generation and refactoring, consistency checks between
the formal development and the manuscript, selected stages of proof auditing,
manuscript organization, and language editing. Their outputs were treated as
suggestions rather than as authoritative sources. The author determined the
mathematical statements and scope, reviewed the arguments and their formal
counterparts, and assumes full responsibility for every mathematical claim,
proof, citation, and statement of priority in this paper.

\par\medskip
\phantomsection
\mergedcontentsline{toc}{section}{References}

\end{document}